\documentclass[b5paper,10pt]{book}

\usepackage[utf8]{inputenc}
\usepackage[english]{babel}

\usepackage{amsmath}
\usepackage{amsfonts}
\usepackage{amssymb}
\usepackage{amsthm}
\usepackage{mathtools}

\usepackage{graphicx}

\usepackage{caption}
\usepackage{subcaption}
\usepackage{wrapfig}

\usepackage[margin=2cm]{geometry}

\usepackage{paralist}

\usepackage{blkarray} 
\usepackage{fdsymbol} 

\usepackage{pdfpages}

\usepackage{csquotes,etoolbox,keyval,ifthen,url,babel}
\usepackage[backend=biber, defernumbers=true, sorting=anyt, style=alphabetic, maxbibnames=4, maxcitenames=2]{biblatex}
\renewbibmacro{in:}{%
  \ifentrytype{article}{}{\printtext{\bibstring{in}\intitlepunct}}}
\AtBeginDocument{
                 }

\usepackage{hyperref}

\hypersetup{
	colorlinks=true,
	linkcolor=red,
	filecolor=magenta,
	urlcolor=cyan,
	citecolor=red,
	citebordercolor=blue
}

\usepackage{url}

\usepackage{tikz}
\newcommand*\circled[1]{\tikz[baseline=(char.base)]{
            \node[shape=circle,draw,inner sep=2pt] (char) {#1};}}

\usepackage{amsthm}

\usepackage{chngcntr}

\allowdisplaybreaks

\usepackage{cmll} 

\usepackage{lastpage}
\usepackage{fancyhdr}
\DeclareMathOperator{\re}{Re}
\DeclareMathOperator{\im}{Im}
\DeclareMathOperator{\tg}{tg}

\DeclareMathOperator{\arctg}{arctg}
\DeclareMathOperator{\arcctg}{arcctg}

\DeclareMathOperator{\card}{card}

\DeclareMathOperator{\adj}{adj}
\DeclareMathOperator{\cof}{cof}
\DeclareMathOperator{\rang}{rang}
\DeclareMathOperator{\proj}{proj}

\DeclareMathOperator{\km}{\,KM}
\DeclareMathOperator{\lit}{\,l}
\DeclareMathOperator{\kg}{\,kg}
\DeclareMathOperator{\m}{\,m}
\DeclareMathOperator{\g}{\,g}

\DeclareMathOperator{\nzd}{NZD}
\DeclareMathOperator{\nzs}{NZS}

\graphicspath{{Slike/}}

\addto\captionsenglish{

}

\DefineBibliographyStrings{english}{%
  bibliography = {Literatura i reference},
}

\usepackage{calc} 
\usepackage{makeidx} 
\makeindex 

\usepackage{empheq}
\usepackage[most]{tcolorbox}

\tcbset{ mojstil1/.style={%
                          enhanced,
                          colframe=red!35,
                          colback=black!3,
                          arc=3pt,
                          boxrule=.251pt,
                         }
        }
\newtcbox{\mymath}[1][]{%
                         nobeforeafter,
                         math upper,
                         tcbox raise base,
                         enhanced,
                         colframe=red!50!,
                         colback=black!5,
                         boxrule=.25pt,
                         arc=3pt,
                         #1
                         }

\usepackage{amsthm}
\usepackage{thmtools}
\usepackage[framemethod=TikZ]{mdframed}

\declaretheoremstyle[  bodyfont=\it,
                       headfont=\bfseries\sc,
                       postheadspace=2em,
                       mdframed={ backgroundcolor=black!3,
                                  linecolor=black!35,
                                  innertopmargin=6pt,
                                  linewidth=.5pt,
                                  roundcorner=5pt,
                                  innerbottommargin=6pt,
                                  skipabove=\parsep,
                                  skipbelow=\parsep},
                                  postheadspace=\newline,
                                  postheadhook={\textcolor{red!25}{\rule[.136em]{\linewidth}{0.001em}}\\}
                     ]{mystyle}

\declaretheoremstyle[  
                       headfont={\bfseries\sc},
                       headfont=\sc,
                       postheadspace=2em,
                       mdframed={ backgroundcolor=black!3,
                                  linecolor=black!35,
                                  innertopmargin=6pt,
                                  linewidth=.5pt,
                                  roundcorner=5pt,
                                  innerbottommargin=6pt,
                                  skipabove=\parsep,
                                  skipbelow=\parsep },
                                  postheadspace=\newline,
                                  postheadhook={\textcolor{red!25}{\rule[.136em]{\linewidth}{0.001em}}\\}
                       ]{mystyle1}

\declaretheoremstyle[  
                       headfont=\bfseries\sc,
                       mdframed={ backgroundcolor=black!3,
                                  linecolor=black!25,
                                  innertopmargin=6pt,
                                  roundcorner=5pt,
                                  innerbottommargin=6pt,
                                  skipabove=\parsep,
                                  skipbelow=\parsep },
                                  postheadspace=\newline,
                                  postheadhook={\textcolor{green!45}{\rule[.136em]{\linewidth}{0.001em}}\\}
                       ]{mystyle2}

\declaretheorem[style=mystyle, name=Teorema]{theorem}

\declaretheorem[style=mystyle1, name=Definicija]{definition} 
\declaretheorem[style=mystyle1, name=Posljedica]{corollary}

\declaretheorem[style=mystyle2, name=Napomena]{remark}
\declaretheorem[style=mystyle2, name=Primjer]{example}
\declaretheorem[style=mystyle2, name=Propozicija]{proposition}

\usepackage{showexpl}
\mdfdefinestyle{frame1}{%
    linecolor=red,
    outerlinewidth=.5pt,
    roundcorner=5pt,
    innertopmargin=\baselineskip,
    innerbottommargin=\baselineskip,
    innerrightmargin=20pt,
    innerleftmargin=20pt,
    backgroundcolor=gray!3!white}

\counterwithin{theorem}{chapter}
\counterwithin{definition}{chapter}
\counterwithin{example}{chapter}
\counterwithin{remark}{chapter}
\counterwithin{lemma}{chapter}
\counterwithin{corollary}{chapter}

 \makeindex

\allowdisplaybreaks
\usepackage{marvosym}

\usepackage{changepage} 

\nocite{*}

\begin{document}
  \includepdf{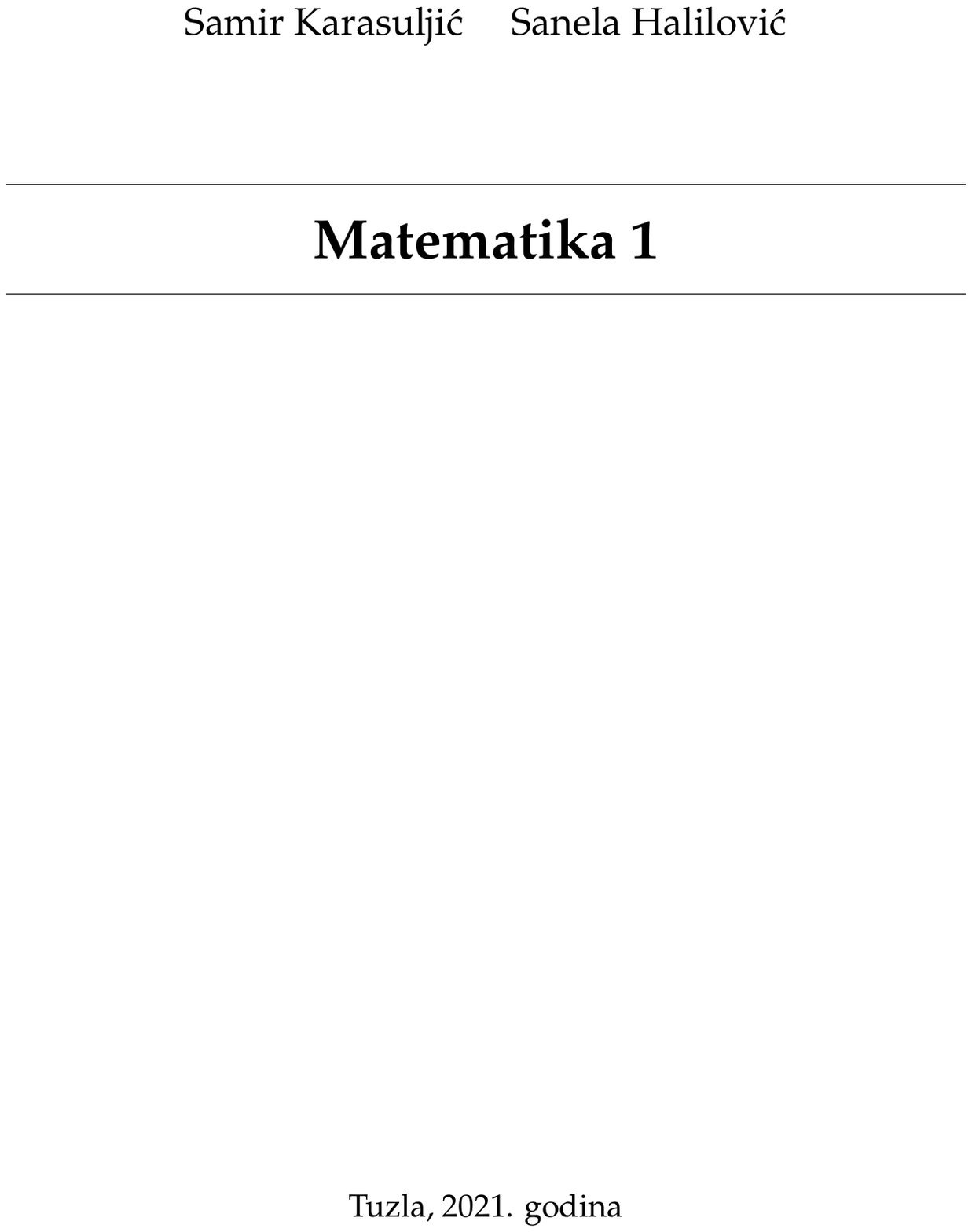}
  \let\cleardoublepage\clearpage
  \frontmatter
  \pagestyle{plain}
  \tableofcontents
  \let\cleardoublepage\clearpage


\chapter*{Predgovor}
\thispagestyle{empty}

Svrha pisanja ove  knjige je da studentima Odsjeka hemija na Prirodno--matemati\v ckom fakultetu, zatim  Odsjeka HIT, PT i IZO na Tehnolo\v skom fakultetu, te studentima na Odjecima pred\v skolski odgoj i razredna nastava Filozofskog fakulteta, kao i studentima sa drugih odsjeka/fakulteta, olak\v sa pripremu i polaganje ispita iz Matematike 1, Osnova matematike 1 kao i drugih predmeta koji obuhvataju izlo\v zeno gradivo.

Knjiga je napisana na osnovu materijala koje smo izlagali na nekoliko fakulteta, bilo kao predavanja, bilo u sklopu auditornih vje\v zbi. Svjesni da u okru\v zenju postoji vi\v se knjiga i ud\v zbenika u kojima su obra\dj ene teme koje su obuhva\'cene ovom knjigom, nastojali smo, na osnovu pomenutog iskustva,  napisati knjigu orjentisanu ka studentima. Knjiga se sastoji od $12$ poglavlja u kojima se uz teoretski obra\dj enu tematsku cjelinu nalaze i rije\v seni reprezentativni primjeri. Na kraju svake cjeline ostavljeno je niz zadataka za samostalni rad kako bi studenti mogli da provjere stepen savladanosti gradiva i adekvatno se pripreme za ispit.

U pisanju materijala za predavanja i auditorne vje\v zbe, a potom i ove knjige, slu\v zili smo se sljede\'com literaturom: U\v s\' cumli\' c i Mili\v ci\' c \cite{uscumlic2002elementi}, Mitrinovi\' c, Mihailovi\' c i Vasi\' c \cite{mitrinovic1973linearna}, Dedagi\' c \cite{dedagic1997uvod, dedagic2005matematicka},  Crnjac, Juki\' c i Scitovski \cite{crnjac1994matematika}, Juki\' c i Scitovski \cite{jukic1998matematika1, jukic2017matematika1}, Vugdali\'c \cite{vugdalic2009matematika, vugdalic2014matematika1}.

Primjedbe, prijedlozi, sugestije, uo\v cene gre\v ske  \Frowny{}, a i pohvale  \Smiley{} su dobrodo\v sle i mo\v zete ih poslati na na\v se email adrese:\\
\href{mailto: samir.karasuljic@unitz.ba}{\textbf{samir.karasuljic@unitz.ba}},
\href{mailto: sanela.halilovic@unitz.ba}{\textbf{sanela.halilovic@unitz.ba}},\\
\href{mailto: samir.karasuljic@gmail.com}{\textbf{samir.karasuljic@gmail.com}},
\href{mailto: halilovicsanela8@gmail.com}{\textbf{halilovicsanela8@gmail.com}}.%
\\

{\small
	\begin{center}
		\Smiley
	\end{center}
	\begin{center}
		\Smiley\quad\quad\Smiley
	\end{center}
}

\noindent Zahvaljujemo se recenzentima:\\
$\ast$ dr.sc. Fehimu Dedagi\'cu, profesor emeritus Univerziteta u Sarajevu,\\
$\ast$ dr.sc. Ramizu Vugdali\'cu, redovni profesor Univerziteta u Tuzli,\\
$\ast$ dr.sc. Nacimi Memi\'c, vanredna profesorica Univerziteta u Sarajevu,\\
koji su svojim zapa\v zanjima i korisnim sugestijama doprinijeli u kona\v cnoj verziji teksta.\\\\

\noindent Tuzla, juli 2021. godine        \hspace{7cm}Autori    
   \addcontentsline{toc}{section}{Predgovor}
   
   \mainmatter
   
   \pagestyle{fancy}
   \pagenumbering{arabic}
   \setcounter{page}{1}

\setcounter{chapter}{0}
\chapter{Matemati\v cka logika}
\pagestyle{fancy}

U savremenoj matematici bitno mjesto zauzima jezik, odnosno kori\v stenje jezika matema-\\ti\v cke logike. Me\dj utim, veliku ulogu matemati\v cka logika ima i van same matematike. Ima zna\v cajnu primjenu kako u prirodnim naukama, tako i u tehni\v ckim disciplinama, na primjer u radu ra\v cunara, u teorijskoj kibernetici itd.

\subsection*{Algebra izkaza}
Polazni objekti teorije iskaza nazivaju se elementarnim iskazima\index{iskaz} ili elementarnim sudovima. Uobi\v cajeno se obilje\v zavaju malim slovima $a,\,b,\,c,\,\ldots$   Pretpostavlja se da elementarni iskaz mora biti ili istinit (ta\v can) ili neistinit (neta\v can). Tako\dj e, smatramo da postoji mogu\' cnost  da se provjeri da li je dati iskaz ta\v can ili neta\v can. U algebri iskaza obi\v cno se ne razmatra sam sadr\v zaj iskaza, nego operacije sa iskazima bez obzira  da li je iskaz ta\v can ili neta\v can. Dakle, pod iskazom podrazumijevamo smi\v sljeno tvr\dj enje, koje ima svojstvo da mo\v ze biti samo ili ta\v cno ili neta\v cno.  Istinitost iskaza $p$ ozna\v cava se, po definiciji, sa gr\v ckim slovom $\tau $ i vrijedi

\begin{empheq}[box=\mymath]{equation*}
	\tau(p)=\left\{
	\begin{array}{l}
		\top, \text{ako je iskaz $p$ ta\v can};\\
		\bot, \text{ako je iskaz $p$ neta\v can}.
	\end{array}
	\right.
	\label{iskaz1}
\end{empheq}

$\top$ se \v cita "te", a $\bot$ "ne te". U literaturi se \v cesto umjesto oznaka $\top$ i $\bot$ koriste oznake $1$ i $0,$ respektivno. \\

\begin{example} Posmatrajmo sljede\' ce re\v cenice:
 \begin{enumerate}
   \item Broj 36 je djeljiv sa 9.
   \item Broj 14 je djeljiv sa 6.
   \item Broj 4 je manji od 10.
   \item Lisabon je glavni grad Ma\dj arske.
 \end{enumerate}
Poja\v snjenje:\\\\
Iskazi 1. i 3. su ta\v cni, dok 2. i 4. nisu.  Sljede\' ce re\v cenice nisu iskazi:
\begin{enumerate}
  \item $2x+3>0.$ (jer nije data vrijednost promjenljive $x$)
  \item Opet samo vrata ciklus. (re\v cenica nema smisla)
  \item Februar ima 28 dana. (nije precizirana godina, te se ne mo\v ze utvrditi ta\v cnost)
  \item Ananas je najukusnije vo\'ce. (zavisi od li\v cnog mi\v sljenja).
\end{enumerate}
\end{example}

\section[L\lowercase{ogi\v cke operacije}]{Logi\v cke operacije}

Od elementarnih (prostih) iskaza  pomo\' cu logi\v ckih operacija formiraju se slo\v zeni iskazi. Rezultati primjene ovakvih operacija mogu se prikazati u tablicama istinitosti. Operacije koje samo jednom iskazu daju odre\dj ene vrijednosti nazivaju se unarne, ima ih ukupno 4. Ako se operacija odnosi na dva iskaza onda je to binarna operacija, a njih ima $2^{2^2}=16.$ Operacije nad $n$ iskaza, tj. $n$--arne operacije ne uvode se direktno, nego preko unarnih i binarnih operacija. Po\v ce\' cemo sa negacijom.\\

\index{logi\v cke operacije!negacija}
\begin{definition}[Negacija]
Negacija iskaza $p$ je iskaz $\neg p$ (\v cita se "ne $p$")  koji je istinit ako i samo ako je iskaz $p$ neistinit.
\end{definition}
\begin{table}[h]\centering
\begin{tabular}{c|c}
$\tau(p)$ & $\tau(\neg p)$ \\\hline
$\top$ &  $\bot$ \\
$\bot$ &  $\top$ \\
\end{tabular}
\caption{Tablica negacije}
\end{table}
\ \\
\begin{example}
 Negirati iskaze:\\
 $p:$ 3 je ve\' ce od 1;\\
 $q:$ 5 je manje od 4.\\\\
\noindent Rje\v senje:\\ \ \\
$\neg p:$ 3 nije ve\' ce od 1;\\
$\neg q:$ 5 nije manje od 4.\\
Vidimo da je $\tau(p)=\top$ ali $\tau(\neg p)=\bot,$ dok je $\tau (q)= \bot$ i $\tau(\neg q)=\top.$
\end{example}

\newpage

\index{logi\v cke operacije!konjukcija}
\begin{definition}[Konjukcija]
  Konjukcija iskaza $p$ i $q$ je slo\v zeni iskaz $p\wedge q$ (\v cita se "$p$ i $q$"), koji je istinit ako i samo ako su oba iskaza $p$ i $q$ istiniti.
\end{definition}
\begin{table}[h]\centering
\begin{tabular}{c|c|c}
$\tau(p)$ & $\tau(q)$ & $\tau(p\wedge q)$ \\\hline
$\top$ & $\top$ & $\top$ \\
$\top$ & $\bot$ & $\bot$ \\
$\bot$ & $\top$ & $\bot$ \\
$\bot$ & $\bot$ & $\bot$ \\
\end{tabular}
\caption{Tablica konjukcije}
\end{table}
\begin{example}
  Formirati konjukciju od zadanih iskaza:\\
  $p:$ Broj 24 je sadr\v zilac broja 6;\\
  $q:$ Broj 24 je djelilac broja 48. \\

\noindent   Rje\v senje:\\  \ \\
$p\wedge q:$ Broj 24 je sadr\v zilac broja 6 i djelilac broja 48. Ova konjukcija je ta\v cna jer su oba iskaza $p$ i $q$ ta\v cni.
\end{example}
\begin{example}
 Formirati konjukciju od sljede\' cih iskaza:\\
 $p:$ \v Cetvorougao ima 2 dijagonale;\\
 $q:$ Petougao ima 4 dijagonale.\\

 \noindent Rje\v senje:\\ \ \\
 $p \wedge q:$ \v Cetvorougao ima 2 dijagonale i petougao ima 4 dijagonale. Konjukcija nije ta\v cna jer petougao ima 5 dijagonala (broj dijagonala je $D_n=\frac{(n-3)n}{2}$).
  \end{example}

\index{logi\v cke operacije!disjunkcija}
\begin{definition}[Disjunkcija]
  Disjunkcija iskaza $p$ i $q$ je slo\v zeni iskaz $p\vee q$ (\v cita se "$p$ ili $q$"),  koji je istinit ako i samo ako je bar jedan od iskaza $p$ ili $q$ istinit.
\end{definition}

\begin{example}
Formirati disjunkciju od iskaza:\\
$p:$ Broj 10 je paran;\\
$q:$ Broj 10 je prost.\\

\noindent Rje\v senje: \\ \ \\
$p\vee q:$ Broj 10 je paran ili broj 10 je prost. Disjunkcija je ta\v cna jer $\tau(p)=\top,$ \v sto je dovoljno da disjunkcija bude ta\v cna.
\end{example}

\begin{table}[h]\centering
	\begin{tabular}{c|c|c}
		$\tau(p)$ & $\tau(q)$ & $\tau(p\vee q)$ \\\hline
		$\top$ & $\top$ & $\top$ \\
		$\top$ & $\bot$ & $\top$ \\
		$\bot$ & $\top$ & $\top$ \\
		$\bot$ & $\bot$ & $\bot$ \\
	\end{tabular}
	\caption{Tablica disjunkcije}
\end{table}

\index{logi\v cke operacije!ekskluzivna disjunkcija}
\begin{definition}[Ekskluzivna disjunkcija]
  Ekskluzivna disjunkcija iskaza $p$ i $q$ je slo\v zeni iskaz $p\veebar q$ (\v cita se "ili $p$ ili $q$") koji je istinit ako i samo ako je istinit  samo jedan od iskaza $p$ i $q.$
\end{definition}
\begin{table}[!h]\centering
\begin{tabular}{c|c|c}
$\tau(p)$ & $\tau(q)$ & $\tau(p\veebar q)$ \\\hline
$\top$ & $\top$ & $\bot$ \\
$\top$ & $\bot$ & $\top$ \\
$\bot$ & $\top$ & $\top$ \\
$\bot$ & $\bot$ & $\bot$ \\
\end{tabular}
\caption{Tablica ekskluzivne disjunkcije}
\end{table}
\begin{example}
Dati su iskazi:\\
$p:$ Sutra u 5 sati poslije podne bi\' cu ku\' ci;\\
$q:$ Sutra u 5 sati poslije podne bi\' cu van ku\' ce.\\

\noindent Rje\v senje:\\ \ \\
Sutra u 5 sati poslije podne ili \' cu biti ku\' ci ili \' cu biti van ku\' ce. Samo jedno od toga mo\v ze biti ta\v cno.
\end{example}

\index{logi\v cke operacije!implikacija}
\begin{definition}[Implikacija]
  Implikacija iskaza $p$ i $q$ je slo\v zeni iskaz $p\Rightarrow q$ (\v cita se "$p$ implicira $q$") koji je neistinit ako i samo ako je istinit iskaz $p,$  a iskaz  $q$ je neistinit.
\end{definition}
\begin{table}[!h]\centering
\begin{tabular}{c|c|c}
$\tau(p)$ & $\tau(q)$ & $\tau(p\Rightarrow  q)$ \\\hline
$\top$ & $\top$ & $\top$ \\
$\top$ & $\bot$ & $\bot$ \\
$\bot$ & $\top$ & $\top$ \\
$\bot$ & $\bot$ & $\top$ \\
\end{tabular}
\caption{Tablica implikacije}
\label{implikacija}
\end{table}

U implikaciji $p\Rightarrow q$ iskaz $p$ se naziva premisa ili pretpostavka,\index{logi\v cke operacije!premisa} a iskaz $q$ posljedica\index{logi\v cke operacije!posljedica} ili konsekvenca implikacije. Iskaz $p\Rightarrow q$ \v cita se jo\v s "iz $p$ proizilazi $q$"; ili "$p$ je dovoljan uslov\index{logi\v cke operacije!dovoljan uslov} za $q$"; ili "$q$ je potreban uslov\index{logi\v cke operacije!potreban uslov} za $p$".

Kada se ka\v ze "$p$ je dovoljan uslov za $q$", onda to zna\v ci da je iskaz $q$ istinit ako je iskaz $p$ istinit. Na primjer, da bi neki broj bio djeljiv sa 2 dovoljno je da bude djeljiv sa 4, Ali ako broj nije djeljiv sa 2, onda nije djeljiv ni sa 4. U ovom slu\v caju djeljiv sa 4 je dovoljan uslov za iskaz djeljiv sa 2, dok je djeljiv sa 2 potreban uslov za djeljiv sa 4. Kada se ka\v ze "$q$ je potreban uslov za $p$", to zna\v ci da iskaz $p$ ne mo\v ze biti istinit ako iskaz $q$ nije istinit.\\

\begin{example}[Poja\v snjenje za tre\' cu i \v cetvrtu vrstu iz Tabele \ref{implikacija} --Tablica implikacije]
Dva dje\v caka sabiraju brojeve 10 i 20. Prvi dje\v cak kada je sabrao ka\v ze da je zbir 40, drugi ka\v ze pogrije\v sio si. Posmatrajmo implikaciju:
Ako je zbir brojeva 10 i 20 jednak 40, onda je prvi dje\v cak pogrije\v sio. Ovdje je \\
$p:\: 10+20=40;$ \\
$q:$ Prvi dje\v cak je pogrije\v sio.\\
Implikacija $p\Rightarrow q$, upravo iskazana, je istinita (iako premisa $p$ nije ta\v cna), i ovakve primjere implikacija susre\' cemo \v cesto u \v zivotu. \\

\noindent A ako je drugi dje\v cak rekao prvom:\\
Ako je $10+20=40,$ onda sam ja visok 5 metara. Ovdje je \\
$p:\:10+20=40,$\\
$q:$ Drugi dje\v cak je visok 5 metara. Oba iskaza su neta\v cna, ali ova implikacija je istinita.
\end{example}

\index{logi\v cke operacije!ekvivalencija}
\begin{definition}[Ekvivalencija]
  Ekvivalencija iskaza $p$ i $q$ je slo\v zeni iskaz $p\Leftrightarrow q$  (\v cita se "$p$ ekvivalentno sa $q$") koji je istinit ako i samo ako su vrijednosti iskaza $p$ i $q$ iste.
\end{definition}
\begin{table}[!h]\centering
\begin{tabular}{c|c|c}
$\tau(p)$ & $\tau(q)$ & $\tau(p\Leftrightarrow  q)$ \\\hline
$\top$ & $\top$ & $\top$ \\
$\top$ & $\bot$ & $\bot$ \\
$\bot$ & $\top$ & $\bot$ \\
$\bot$ & $\bot$ & $\top$ \\
\end{tabular}
\caption{Tablica ekvivalencije}
\end{table}
\newpage
Sljede\'ci primjer nam ilustruje vezu implikacije i ekvivalencije. \\
\begin{example}
 Iskaze, koji su implikacije:
  \begin{enumerate}
  	\item [$p:$] Ako je broj djeljiv sa 6, onda je djeljiv i sa 2 i sa 3;
  	\item [$q:$]  Ako je broj djeljiv i sa 2 i sa 3, onda je on djeljiv i sa 6,
  \end{enumerate}	
mo\v zemo spojiti u jedan slo\v zeni iskaz (ekvivalenciju). To radimo na sljede\'ci na\v cin:
 \begin{enumerate}
 	\item [$p\Leftrightarrow q:$] Da bi prirodan broj bio djeljiv sa 6, potrebno je i dovoljno da bude djeljiv i sa 2 i sa 3, ili
 	\item [$p\Leftrightarrow q:$] Prirodan broj je djeljiv sa 6, ako i samo ako je djeljiv i sa 2 i sa 3, ili
 	\item [$p\Leftrightarrow q:$] Prirodan broj je djeljiv sa 6, onda i samo onda ako je djeljiv i sa 2 i sa 3.
 \end{enumerate}	
\end{example}

\index{logi\v cke operacije!iskazna formula}
\begin{definition}[Iskazna formula]
  Iskazna formula je kona\v can niz iskaza sastavljen pomo\' cu logi\v ckih operacija, konstanti $\top$ i $\bot$ (ili 1 i 0) i glavnih promjenljivih iskaza $p,\,q,\,r,\ldots$
\end{definition}

\begin{example}[Iskazne formule]
  Iskazne formule sa dva i tri iskazna slova.
  \begin{enumerate}[$(a)$]
     \item $(p\Rightarrow\neg r)\wedge(p\veebar r);$
     \item $(p\wedge q)\vee r$.
  \end{enumerate}
\end{example}
\index{logi\v cke operacije!tautologija}
\begin{definition}[Tautologija]
  Iskazna formula koja je ta\v cna bez obzira na istinitosnu vrijednost iskaza u njoj, naziva se \textbf{tautologija}. Iskazna formula koja je uvijek neta\v cna (za sve istinitosne vrijednosti
  prostih iskaza koji u njoj u\v cestvuju), naziva se \textbf{kontradikcija}. \textbf{Alternativna fo-\\rmula} je za neke vrijednosti prostih iskaza ta\v cna, a za preostale vrijednosti prostih iskaza
  neta\v cna.
\end{definition}

\begin{example}
 Ispitati da li je iskazna formula tautologija
    \begin{enumerate}[$(a)$]
       \item $(p\wedge \neg p)\Rightarrow q$;\label{zadatak1}
       \item $p\Rightarrow (q\vee r).$\label{zadatak2}
    \end{enumerate}
\ \\

\noindent Rje\v senje:
   \begin{enumerate}[$(a)$]
   	\item Tabela \ref{tabela1};
   	\item Tabela \ref{tabela2}.
   \end{enumerate}
\end{example}
\begin{table}[!h]\centering\small
   \begin{tabular}{c|c|c|c|c}
   $\tau(p)$ & $\tau(q)$ & $\tau(\neg p)$ & $\tau(p\wedge\neg p)$ & $\tau((p\wedge\neg p)\Rightarrow q)$ \\ \hline
   $\top$ & $\top$ & $\bot$ & $\bot$ & $\top$ \\
   $\top$ & $\bot$ & $\bot$ & $\bot$ & $\top$ \\
   $\bot$ & $\top$ & $\top$ & $\bot$ & $\top$ \\
   $\bot$ & $\bot$ & $\top$ & $\bot$ & $\top$ \\
   \end{tabular}
   \caption{Rje\v senje za iskaznu formulu datu pod \eqref{zadatak1} }
   \label{tabela1}
\end{table}
\begin{table}[!h]\centering\small
    \begin{tabular}{c|c|c|c|c}
      $\tau(p)$&$\tau(q)$&$\tau(r)$&$\tau(q\vee r)$& $\tau(p\Rightarrow(q\vee r))$\\\hline
      $\top$&$\top$&$\top$&$\top$& $\top$\\
      $\top$&$\top$&$\bot$&$\top$& $\top$\\
      $\top$&$\bot$&$\top$&$\top$& $\top$\\
      $\top$&$\bot$&$\bot$&$\bot$& $\bot$\\
      $\bot$&$\top$&$\top$&$\top$& $\top$\\
      $\bot$&$\top$&$\bot$&$\top$& $\top$\\
      $\bot$&$\bot$&$\top$&$\top$& $\top$\\
      $\bot$&$\bot$&$\bot$&$\bot$& $\top$\\
    \end{tabular}
    \caption{Rje\v senje za iskaznu formulu datu pod \eqref{zadatak2} }
    \label{tabela2}
\end{table}
Iz Tabele \ref{tabela1} zaklju\v cujemo da je formula data u \eqref{zadatak1} tautologija, dok na osnovu Tabele \ref{tabela2} vidimo da formula data u \eqref{zadatak2} nije tautologija.
Neke poznate tautologije imaju i svoje nazive, kao \v sto su: De Morganove\footnote{Augustus De Morgan (27.juni 1806. -- 18.mart 1871. godine) bio je britanski matemati\v car i logi\v car. Formulisao je De Morganove zakone, uveo pojam matemati\v cke indukcije} formule, zakon kontrapozicije (kontradikcije), zakon isklju\v cenja tre\' ceg i tako dalje (vidjeti 1.zadatak u sekciji 1.3).

\section[P\lowercase{redikati i kvantifikatori}]{Predikati i kvantifikatori} Mogu\' cnosti primjene algebre iskaza, posebno u matematici, dosta su ograni\v cene jer se ne mogu iskoristiti u iskazivanju osobina i odnosa izme\dj u matemati\v ckih objekata. Ovaj nedostatak otklanja se uvo\dj enjem novih pojmova i operacija. Na primjer, ako su $x$ i $y$ realni brojevi, onda tvrdnja $"x$ je paran broj" ili "$x$ je manji od $y$" mogu biti kako istiniti tako i neistiniti, zavisno od vrijednosti realnih brojeva $x$ i $y.$ Ako je $x=1$ onda je prva tvrdnja neta\v cna, a ako je $x=2$ onda je ta tvrdnja ta\v cna. Za $x=3$ i $y=5$ druga  tvrdnja je  ta\v cna (istinita), dok za $x=6$ i $y=4$ neistinita je. Tvrdnje ovakve vrste smatramo neodre\dj enim.    U ovakvim slu\v cajevima elementarni iskazi tretiraju se kao promjenljive koje uzimaju jednu od vrijednosti $\top$ i $\bot.$ U matematici  susre\' cemo se sa iskazima koji se odnose na izvjesne objekte, ili elemente.

Neka je dat skup $M$ objekata (elemenata) $x,y,z,\ldots$  Svojstva koja se odnose na te objekte obilje\v zavaju se sa
\[P(x),\:Q(y),\:R(x,y),\ldots\]
Ako je $M$ skup prirodnih brojeva, u tom slu\v caju  prethodna svojstva mogu, npr., biti: "$x$ je prost broj", "$y$ je paran broj", "$x$ je manje od $y$", respektivno. Sada ovi iskazi o svojstvu mogu biti kako istiniti, tako i neistiniti. U nastavku samo \' cemo razmatrati iskaze sa ovakvim osobinama.

Neka je $x$ proizvoljni element iz nekog skupa objekata $M$.  U tom slu\v caju  neko svojstvo $P(x)$ potpuno je odre\dj eno kada je $x$ objekat iz $M.$ Isto tako neko svojstvo $R(x,y)$ potpuno je odre\dj eno ako su $x$ i $y$ odre\dj eni objekti iz skupa $M.$ Ovakvi iskazi postaju funkcije od tih objekata, tj. promjenljivih iz skupa $M.$

\begin{definition}[Predikat]
  Neodre\dj eni iskazi koji zavise od jedne ili vi\v se promjenljivih nazivaju se logi\v cke funkcije ili predikati.
\end{definition}
Skup $M$ naziva se predmetna oblast, dok se elementi skupa $M$ nazivaju predmetne promjenljive ili promjenljivi objekti, a ako su neki od tih objekata fiksirani onda se oni nazivaju individualni objekti ili predmetne konstante. Iskazi ili iskazne promjenljive i predikati koji zavise od individualnih objekata nazivaju se, u logici predikata, elementarni iskazi. Iskazne promjenljive i predikati kako od individualnih objekata tako i od predmetnih promjenljivih nazivaju se elementarne formule logike predikata.
\begin{definition}[Formule logike predikata]
  Pod formulom logike predikata podrazumijeva se svaki kona\v can izraz formiran od elementarnih iskaza i elementarnih formula logike predikata, primjenom kona\v cnog broja logi\v ckih operacija $\wedge,\,\vee, \Rightarrow,\Leftrightarrow,\,\neg.$
\end{definition}
Mo\v zemo primjetiti da logika  predikata predstavlja poop\v stenje algebre iskaza, po\v sto sadr\v zi u sebi cijelu algebru iskaza kao i sve njene operacije i formule.

Osim operacija algebre iskaza u logici predikata upotrebljavaju se jo\v s dvije nove\\ operacije, koje su povezane sa osobinama logike ili ra\v cuna predikata.
\begin{enumerate}[$1^{0}$]
  \item Univerzalni kvantifikator. Neka je $P(x)$ potpuno odre\dj en predikat koji ima jednu od dvije vrijednosti $\top$ ili $\bot$ za svaki element skupa $M.$ Tada se
       pod izrazom \[(\forall x)P(x),\]
       koji se \v cita "za svako $x$ $P(x)$ je istinito",  podrazumijeva istinit iskaz, kada je $P(x)$ istinito za svako $x$ iz skupa $M$, a neistinit u suprotnom slu\v caju.
       Sam simbol $\forall$ naziva se univerzalni kvantifikator.
  \item Egzistencijalni kvantifikator. Neka je $R(x)$ neki predikat. Formula \[(\exists x)R(x),\]  koja se \v cita "postoji $x$ za koje je $R(x)$ istinito", defini\v se se tako  da poprima vrijednost $\top,$ ako postoji element $x$ iz skupa $M,$ za koji je $R(x)$ istinito, a vrijednost $\bot$ u suprotnom slu\v caju. Znak $\exists$ naziva se egzistencijalni kvantifikator.
\end{enumerate}

\begin{example}
  \begin{enumerate}
    \item $(\forall x)\:x=x,$ je istinit iskaz;
    \item $(\forall x)(\exists y)\:x-y=2,$ \v cita se "za svako $x$ postoji $y,$ tako da je $x-y=2$  istinit iskaz";
    \item $(\forall x)(\forall y)\:(x=y)\Rightarrow (x^3=y^3),$ \v cita se "za svako $x$ i svako $y,$ iz $x=y$ slijedi $x^3=y^3$";
    \item $(\exists x)(\exists y)\:x^2+y^2=4,$ \v cita se "postoji $x$ i postoji $y,$ takvi da je $x^2+y^2=4$ istinit iskaz";
    \item $(\exists x \in \mathbb{R})\:x^2+1=0,$ ovaj iskaz nije ta\v can.
  \end{enumerate}
\end{example}

\index{kvantifikatori}

\section[Z{\lowercase{adaci}}]{Zadaci}  \index{Zadaci za vje\v zbu!logika}

\begin{enumerate}
	\item Ispitati da li su sljede\' ce formule tautologije:\\
	\begin{inparaenum}
		\item $p\vee \neg p;\quad$ (zakon isklju\v cenja tre\' ceg); \\
		\item $(p\Rightarrow q)\Leftrightarrow (\neg q\Rightarrow \neg p)\,$ (zakon kontrapozicije);\\
		\item $\neg(p\wedge q)\Leftrightarrow (\neg p\vee \neg q)\,$ (De Morganova formula);\\
		\item $(p\vee q)\Rightarrow \neg q;\quad$
		\item $(p\wedge q)\vee r;\quad$
		\item $(p\vee \neg q)\Rightarrow (q\wedge r);\:$
		\item $(p\wedge \neg q)\Rightarrow (p\vee\neg r).$
	\end{inparaenum}
	\item Pomo\' cu kvantifikatora izraziti re\v cenicu: Za svaki realan broj $a$ je $a^3=a\cdot a\cdot a,$ (oznaka za skup realnih brojeva je $\mathbb{R}$).
	\item Izraziti pomo\' cu kvantifikatora re\v cenicu: Proizvod tri realna broja $a,\,b$ i $c$ je nula onda i samo onda ako je $a=0$ ili $b=0$ ili $c=0$.
	\item Napisati pomo\' cu kvantifikatora: Za svaki realan broj $x$ je $(x-1)(x+1)=x^2-1.$
	\item Napisati pomo\' cu kvantifikatora: Za svako $a$ je $a+a+a=3a.$
	\item Napisati pomo\' cu kvantifikatora: Postoji bar jedno $a$ za koje je $a^2=4.$
	\item Napisati pomo\' cu kvantifikatora:
	  \begin{enumerate}
	    \item  Za svako $x$ postoji $y$ takvo da je $x+3y=4;$
        \item  Postoje $x$ i $y$ takvi da je $x+y=7.$
      \end{enumerate}
	\item Napisati pomo\' cu kvantifikatora: Za svako $x$ iz skupa prirodnih brojeva postoji $y$ iz skupa realnih brojeva takav da je $y:x=5,$
	(oznaka za skup prirodnih brojeva je $\mathbb{N}$).
	\item Napisati rije\v cima iskaze date formulama:\\
	\begin{inparaenum}
		\item $(\exists y\in\mathbb{N})\:y+1\geqslant 0;\quad$
		\item $(\forall x,y\in\mathbb{R})\:x^2+y^2\geqslant 0;\quad$
		\item $(\forall x\in\mathbb{N})(\exists y\in\mathbb{Z})x>y,$ (oznaka za skup cijelih brojeva je $\mathbb{Z}$).
	\end{inparaenum}
	\item Utvrditi istinitnu vrijednost iskaza:\\
	\begin{inparaenum}
		\item $(\forall x \in \mathbb{R})(\exists y\in \mathbb{R})\:x+y>0;\quad$
		\item $(\exists y\in\mathbb{R })(\forall x\in \mathbb{R })\:x+y>0.$
	\end{inparaenum}
\end{enumerate}


\setcounter{chapter}{1}
\chapter{Skupovi}
\pagestyle{fancy}

Skup, odnosno pojam skupa, igra zna\v cajnu ulogu u savremenoj matematici, kako zbog \v cinjenice da je teorija skupova postala jedna detaljna i sadr\v zajna matemati\v cka disciplina, tako i zbog uticaja na druge oblasti matematike. U dana\v snjoj se matematici mo\v ze uo\v citi  skupovno--teorijski pristup. Bitno je pomenuti da se jezik teorije skupova sve vi\v se koristi i u drugoj nau\v cnoj i tehni\v ckoj literaturi. Primjena ovog apstraktnog jezika omogu\' cava, izme\dj u ostalog, lak\v se povezivanje i ispitivanje raznih procesa i pojava koje se prili\v cno razlikuju po raznim kriterijumima.

Osnove teorije skupova postavio je krajem devetnaestog  vijeka njema\v cki matemati\v car Cantor. \footnote{George Cantor, 3.3.1845.--6.1.1918. godine, bio je njema\v cki matemati\v car. Osniva\v  c je teorije skupova.}

\section[U\lowercase{vodni pojmovi}]{Uvodni pojmovi, definicije i notacija}

Skup i elemente skupa smatramo intuitivno jasnim i zbog toga ih prihvatamo bez definisanja, tj. smatramo ih osnovnim pojmovima. Primjeri skupova su: skup studentica/ studenata u jednoj grupi na nekom fakultetu, skup stanovnika u nekoj zgradi, skup svih stanovnika jednog mjesta i dr. U  matematici su posebno bitni skupovi sa matemati\v ckim objektima, a to su, na primjer, skup prirodnih brojeva, skup realnih brojeva, skupovi koje sa\v cinjavaju razni geometrijski objekti i dr. Neki skupovi koji se \v ce\v s\' ce upotrebljavaju imaju svoje standardne oznake, a to su:
\begin{enumerate}[$\bullet$]
  \item $\mathbb{N}$  skup svih prirodnih brojeva;
  \item $\mathbb{Z}$  skup svih cijelih brojeva;
  \item $\mathbb{Q}$  skup svih racionalnih brojeva;
  \item $\mathbb{R}$  skup svih realnih brojeva;
  \item $\mathbb{R}^{+}$  skup svih realnih pozitivnih brojeva;
  \item $\mathbb{I}$  skup svih iracionalnih brojeva;
  \item $\mathbb{C}$  skup svih kompleksnih brojeva;
  \item $(a,b)$  otvoren interval u $\mathbb{R}$ ili kra\' ce interval;
  \item $[a,b]$  zatvoren interval u $\mathbb{R}$ ili kra\' ce segment.
\end{enumerate}
\index{skup!prirodnih brojeva}
\index{skup!cijelih brojeva}
\index{skup!racionalnih brojeva}
\index{skup!realnih brojeva}
\index{skup!realnih pozitivnih brojeva}
\index{skup!iracionalnih brojeva}
\index{skup!kompleksnih  brojeva}
\index{skup!otvoreni interval}
\index{skup!zatvoreni interval ili segment}

Elementi mogu pripadati ili ne pripadi nekom skupu. Tvrdnje "element $x$ pripada skupu $A$" ili "$x$ je element skupa $A$" ili "$x$ je ta\v cka iz skupa $A$," imaju isto zna\v cenje (smisao) i ovo pi\v semo simboli\v cki $x\in A$ (\v citamo $x$ pripada skupu $A$), a oznaka $\in$ zna\v ci pripada. U slu\v caju da element $x$ ne pripada skupu $A,$ onda pi\v semo $x\notin A,$ (\v citamo $x$ ne pripada skupu $A$), $\notin$ je oznaka da neki element ne pripada skupu.

U slu\v caju da elementi $a,b,c$ pripadaju skupu $A$ i da su to svi elementi skupa $A,$ pi\v semo

\begin{empheq}[box=\mymath]{equation*}
A=\{ a,b,c\}.
\label{skup1}
\end{empheq}

Skup ne zavisi od rasporeda (poretka) kojim su dati njegovi elementi. Tako su npr. skupovi $\{1,2,3\}$ i  $\{2,1,3\}$ jednaki. Tako\dj e su jednaki i skupovi $\{1,2,3\},$ $\{1,1,2,3\},$ $\{1,2,3,3,3,3\},$ itd., vidjeti definiciju jednakosti dva skupa  \eqref{jednakostSkupova}.

Ako je $n$ prirodan broj, skup $S=\{x_1,x_2,\ldots,x_n\}$ od $n$ elemenata $x_1,x_2,\ldots,x_n$ je kona\v can. Skup je beskona\v can ako broj njegovih elemenata nije kona\v can. Broj elemenata nekog skupa $A$ ozna\v cavamo sa $\card A$ ili $|A|$ i \v citamo kardinalni broj\footnote{Kod skupova sa beskona\v cnim brojem elemenata kardinalni broj predstavlja mo\' c skupa.} skupa $A.$

Elementi skupa imaju ponekad neku zajedni\v cku osobinu, na primjer $P(x),$ pa mo\v zemo pisati
\begin{empheq}[box=\mymath]{equation*}
B=\{x:P(x)\},
\label{skup2}
\end{empheq}
a ovo zna\v ci "$B$ je skup svih elemenata koji imaju osobinu $P(x)$".
\begin{example}
 Napisati simboli\v {c}ki skup realnih brojeva koji su ve\' {c}i od 2 a manji od 10.\\

\noindent Rje\v senje:
\[B=\{x\in\mathbb{R}:2<x<10\},\]
u ovom slu\v{c}aju $P(x)$ je $2<x<10.$
\end{example}

\begin{figure}[!h]\centering
\includegraphics[scale=1]{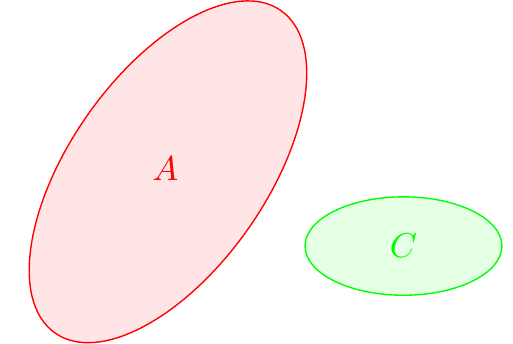}
\includegraphics[scale=.85]{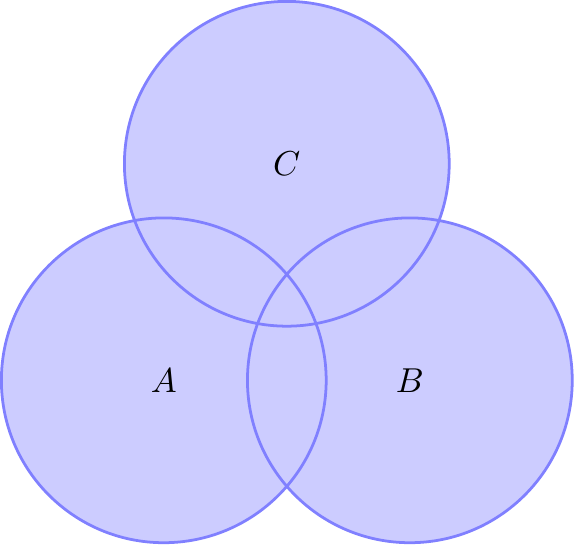}
\caption{Predstavljanje skupova Euler--Vennovim dijagramima}
\label{vennovi1}
\end{figure}

Skup mo\v {z}e sadr\v zavati kona\v can broj elemenata, beskona\v can broj elemenata, a mo\v ze biti i bez elemenata. U slu\v caju da je broj elemenata $0$ (bez elemenata) ka\v zemo da je to prazan skup. Oznaka za prazan skup je $\emptyset .$ 
Skupove mo\v zemo predstavljati i grafi\v cki, a jedan od na\v cina je pomo\' cu Euler\footnote{Leonhard Euler (15.april 1707.--18.septembar 1783. godine) bio je \v svajcarski matemati\v car, fizi\v car, astronom, logi\v car, in\v zinjer. Uspostavio je temelje teorije grafova i topologije, te dao veliki doprinos u mnogim oblastima matematike. }--Vennovih\footnote{John Venn (4.avgust 1834.--4.april 1923. godine) bio je engleski matemati\v car, logi\v car i filozof. Poznat po uvodjenju Euler--Vennovih dijagrama.} dijagrama (Slika \ref{vennovi1}).

Ako su svi elementi jednog skupa A sadr\v zani u drugom skupu B, simboli\v cki zapisano $(\forall x)(x\in A\Rightarrow x\in B),$ onda ka\v zemo da je skup $A$ sadr\v zan u skupu $B,$ ili da je $A$ podskup od $B$ i pi\v semo

\begin{empheq}[box=\mymath]{equation*}
A\subset B.
\end{empheq}

\newpage

\ \\ \ \\

\begin{figure}[!h]\centering
  \includegraphics[scale=1]{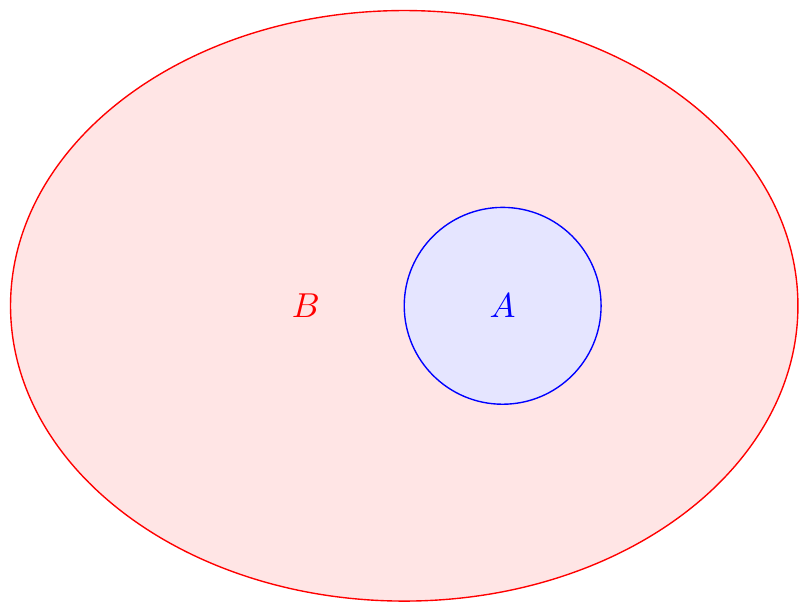}
  \caption{Skup $A$ je podskup skupa $B$}
  \label{podskup}
\end{figure}

\begin{definition}[Podskup]
  Svaki skup koji je sadr\v zan u nekom skupu $B$ naziva se podskup skupa $B.$
\end{definition}

Ako vrijedi da su elementi skupa $A$ sadr\v zani u skupu $B,$ a i elementi skupa $B$ sadr\v zani u skupu $A,$ onda skupovi $A$ i $B$ imaju iste elemente.
Drugim rije\v cima $A$ i $B$ su jednaki i vrijedi

\begin{empheq}[box=\mymath]{equation}\label{jednakostSkupova}
A=B.
\end{empheq}

Prazan skup je podskup svakog skupa.\\
U literaturi susre\'{c}emo i oznaku $\subseteq $ za podskup, npr. $A \subseteq B,$ a tada $A\subset B$ zna\v ci da je $A$ strogi poskup od $B$
(tj. $(A \subset B)\Leftrightarrow (A\subseteq B \wedge A\neq B)$).

\begin{definition}[Partitivni skup]
 Skup svih podskupova skupa $A$ naziva se partitivni skup skupa $A$ i obilje\v zava se sa $P(A).$
\end{definition}

\begin{example}
  Dat je skup $A=\{a,b,c\}.$ Odrediti partitivni skup skupa $A$.\\

\noindent Rje\v senje:
\[P(A)=\{\emptyset,\{a\},\{b\},\{c\},\{a,b\},\{a,c\},\{b,c\},\{a,b,c\}\}.\]
\end{example}
Ako je $n$ broj elemenata skupa $A,$ onda je broj elemenata partitivnog skupa $P(A)$ jednak $2^{n},$ tj. ako je $\card A =n,$ onda je $\card P(A)=2^{n}.$

\v Cesto se uka\v  ze potreba za pojmom univerzalnog skupa. U mnogim matemati\v ckim oblastima javlja se potreba za nekim op\v stim skupom, koji \' ce sadr\v zavati sve skupove koji se razmatraju u toj oblasti. Uobi\v cajeno ovakav skup se naziva univerzalni skup (ili univerzum). Univerzalni skup se ne defini\v se  nego se jednostavno podrazumijeva da znamo o kakvom se skupu radi. Po pravilu univerzalni skup je relativan pojam, razlikuje se od oblasti do oblasti, na primjer on mo\v ze biti skup svih ta\v caka ravni ili skup svih ta\v caka prave ili skup svih racionalnih brojeva i tako dalje. Dakle, univerzalni skup zavisi od potrebe i oblasti gdje se koristi.

\begin{figure}[!h]\centering
     \begin{subfigure}[b]{.45\textwidth}\centering
        \includegraphics[scale=.85]{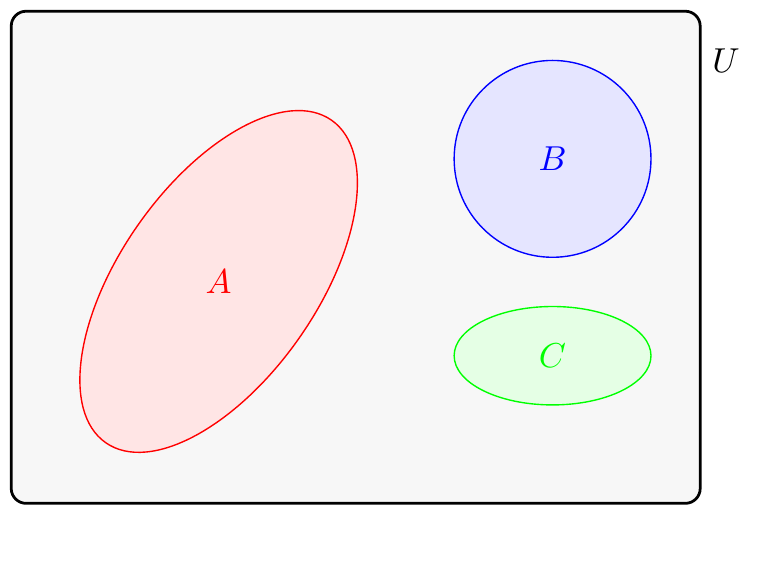}
         \caption{Univerzalni skup $U$}
          \label{univerzalni}
   \end{subfigure}\hspace{.5cm}
   \begin{subfigure}[b]{.45\textwidth}\centering
        \includegraphics[scale=.85]{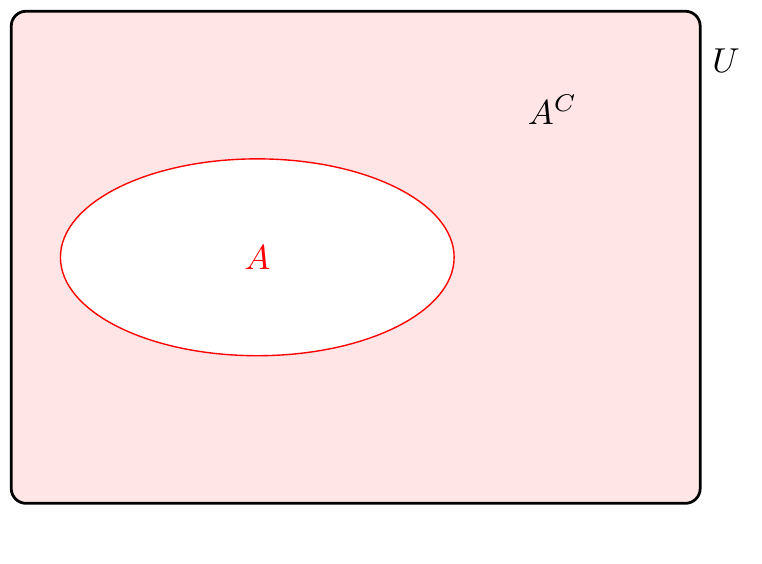}
         \caption{Komplement $A^C$}
      \label{komplement}
     \end{subfigure}
        \caption{Univerzalni skup $U$ i komplement $A^C,$ skupa $A$ u odnosu na univerzalni skup $U$}
     \label{univwrzalni1}
  \end{figure}

U odnosu na univerzalni skup mo\v zemo uvesti pojam komplementa skupa.
\begin{definition}[Komplement skupa]
  Ako je $A\subset U,$ onda se pod komplementom skupa $A$ u odnosu na skup $U$ podrazumijeva skup $A^C$  za koji vrijedi
  \[A^{C}=\{x:x\in U\wedge x\notin A\}.\]
\end{definition}
Na Slici \ref{komplement} prikazan je svijetlocrvenom bojom komplement skupa $A$ u odnosu na skup $U.$

\section[O\lowercase{peracije sa skupovima}]{Operacije sa skupovima}

Za neki dati univerzalni skup $U,$ mo\v zemo od elemenata njegovog partitivnog skupa $P(U)$ formirati nove skupove primjenom  raznih operacija.
\paragraph{Unija skupova}
\begin{definition}[Unija skupova]\index{skup!unija skupova}
  Unija skupova $A$ i $B$ (u oznaci $A\cup B$) je skup svih elemenata koji se nalaze bar u jednom od skupova $A$ ili $B.$
\end{definition}
Simboli\v cki zapisana unija skupova $A$ i $B$ izgleda ovako
\begin{empheq}[box=\mymath]{equation*}
  A\cup B=\{ x:x\in A\vee x\in B\}.
\end{empheq}

\begin{figure}[!h]\centering
 \includegraphics[scale=1]{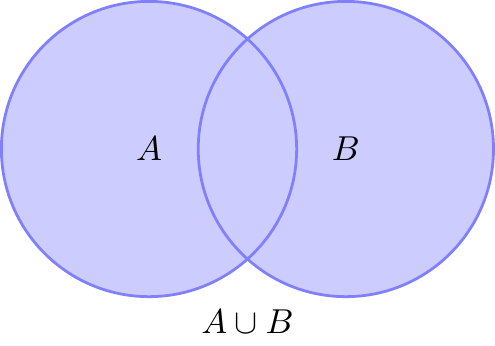}
 \caption{Unija skupova $A$ i $B$}
\end{figure}

\begin{example}
  Ako je $A=\{1,2,3,4,5,6,7\},\:B=\{3,4,9,10\},$ odrediti $A\cup B.$ \\\\
\noindent Rje\v senje: \[A\cup B=\{1,2,3,4,5,6,7,9,10\}.\]
\end{example}

\paragraph{Presjek skupova}
\begin{definition}[Presjek skupova]\index{skup!presjek skupova}
  Presjek skupova $A$ i $B$ (u oznaci $A\cap B$) je skup svih elemenata koji istovremeno pripadaju i skupu  $A$ i skupu $B.$
\end{definition}
Simboli\v cki zapisan presjek skupova $A$ i $B$ izgleda ovako
\begin{empheq}[box=\mymath]{equation*}
   A\cap B=\{ x:x\in A\wedge x\in B\}.
\end{empheq}
Ako dva skupa nemaju zajedni\v ckih elemenata, za njih imamo naziv dat u sljede\' coj definiciji.

\begin{definition}[Disjunktni skupovi]\index{skup!disjunktni skupovi}
Za dva skupa ka\v zemo da su disjunktni ako nemaju zajedni\v ckih elemenata, tj. \\ $A\cap B=\emptyset.$
\end{definition}

\begin{figure}[!h]\centering
  \includegraphics[scale=1]{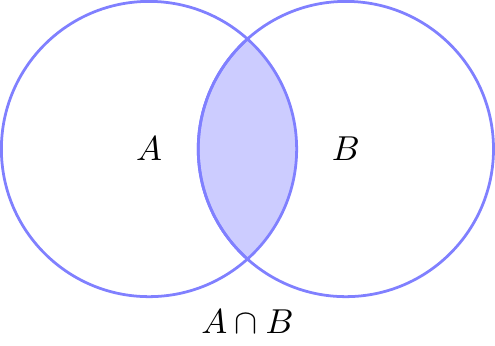}
  \caption{Presjek skupova $A$ i $B$}
\end{figure}
\begin{example}
  Dati su skupovi $A=\{1,2,3,5\},\,B=\{1,3,7,8\},\,C=\{6,7,9\}.$ Odrediti $A\cap B,\,B\cap C$ i $A\cap C.$\\\\
  \noindent Rje\v senje: \[A\cap B=\{1,3\},\,B\cap C=\{7\},\,A\cap C=\emptyset.\]
\end{example}

\paragraph{Razlika i simetri\v cna razlika skupova}
\begin{definition}[Razlika skupova]\index{skup!razlika skupova}
  Razlika (diferencija) skupova $A$ i $B$ (u oznaci $A\setminus B$) je skup svih elemenata skupa $A$ koji ne pripadaju skupu $B.$
\end{definition}
Simboli\v cki zapisana razlika skupova $A$ i $B$ izgleda ovako
\begin{empheq}[box=\mymath]{equation*}
  A\setminus B=\{ x:x\in A\wedge x\notin B\}.
\end{empheq}

\newpage

\begin{figure}[!h]\centering
  \includegraphics[scale=1]{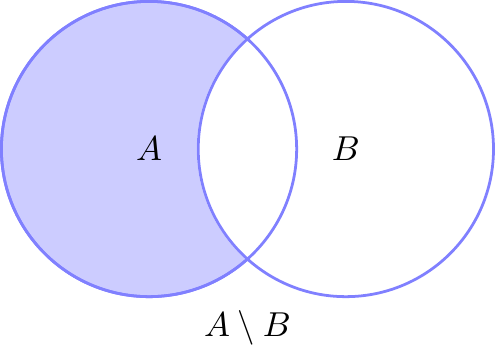}\hspace{1cm}
  \includegraphics[scale=1]{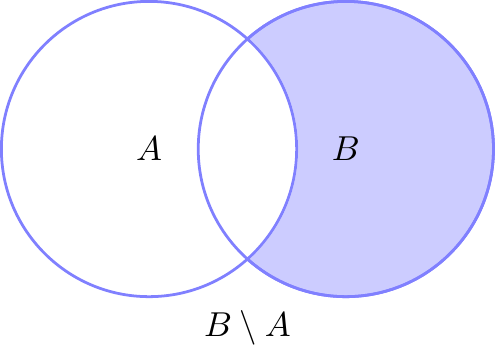}
  \caption{Razlika skupova $A\setminus B$ i $B\setminus A.$}
\end{figure}

\begin{definition}[Simetri\v cna razlika skupova]
  Skup $(A\setminus B)\cup (B\setminus A)$ je simetri\v cna razlika skupova $A$ i $B$ i obilje\v zava se sa $A\triangle B.$
\end{definition}
\begin{figure}[!h]\centering
  \includegraphics[scale=1]{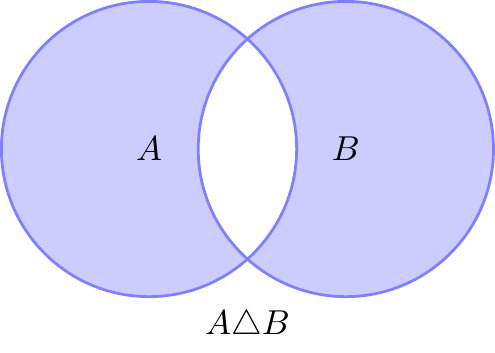}
  \caption{Simetri\v cna razlika skupova}
\end{figure}
Mo\v zemo primijetiti da vrijedi:
$$A\triangle B= (A\setminus B)\cup (B\setminus A)=(A\cup B)\setminus (A\cap B).$$
\begin{example}
   Dati su skupovi $A=\{a,b,c,d,e,f\},\,B=\{d,e,f,g,h\}.$ Odrediti\\
        $A\cup B,\,A\cap B,\,A\setminus B,\,B\setminus A,\,A\triangle B.$\\\\
\noindent Rje\v senje:
  \begin{align*}
     A\cup B&=\{a,b,c,d,e,f,g,h \};\:
     A\cap B=\{ d,e,f\};\:
     A\setminus B=\{ a,b,c\};\\
     B\setminus A&=\{ g,h\},\:
     A\triangle B=\{ a,b,c,g,h\}.
  \end{align*}
\end{example}

\section[D\lowercase{escartesov} (D\lowercase{ekartov}) \lowercase{proizvod skupova}]{Descartesov\footnote{Ren\'e Descartes (31.mart 1596.--11.februar 1650.godine) bio je francuski matemati\v car, filozof i generalno nau\v cnik. Osnovao analiti\v cku geometriju, povezuju\'ci prethodno razdvojene oblasti u geometriji i algebri.         } (Dekartov) proizvod skupova}
 Simboli $\{a,b\}$ i $\{b,a\}$ ozna\v cavaju isti skup od dva elementa $a$ i $b.$ U ovakvim skupovima redoslijed elemenata nije bitan. Me\dj utim, nekada je potrebno da znamo i redoslijed elemenata, na primjer koordinate ta\v cke u pravouglom koordinatnom sistemu u ravni \v cine dva broja i bitno je da znamo \v sta je prva koordinata ($x$--koordinata), a \v sta druga koordinata ($y$--koordinata). Da bi rije\v sili ovaj problem, uve\v s\' cemo pojam ure\dj enog para \v cija je prva komponenta (projekcija) $a,$ a druga komponenta (projekcija) $b.$
\begin{definition}[Ure\dj eni par]\index{ure\dj eni par}
Ure\dj eni par elemenata $a$  i $b$, u oznaci $(a,b)$,  je
\[(a,b)=\{\{a\},\{a,b\}\}.\]
\end{definition}
Ovaj ure\dj eni par ozna\v cavamo sa $(a,b)$\footnote{Ne mije\v sati sa intervalom $(a,b)$} ili $\langle a,b\rangle.$ Smatramo da je $(a,b)$ razli\v cito od $(b,a)$ osim u slu\v caju $a=b.$  Ure\dj eni parovi $(a,b)$ i $(c,d)$ su jednaki ako i samo ako je $a=c$ i $b=d.$

\begin{definition}[Ure\dj ena trojka]\index{ure\dj ena trojka}
  Ure\dj ena trojka $(a,b,c)$ elemenata $a,b,c$ defini\v se se pomo\' cu jednakosti
  \[(a,b,c)=((a,b),c).\]
\end{definition}

Na sli\v can se na\v cin defini\v se i ure\dj ena $n$--torka $(a_1,a_2,\ldots,a_n)$.

\begin{definition}[Descartesov proizvod]\index{Descartesov proizvod}
  Dekartov proizvod dva neprazna skupa $A$ i $B$ je skup $C$ \v ciji su elementi svi ure\dj eni parovi sa prvom komponentom iz skupa $A,$ a drugom komponentom iz skupa $B,$ tj.
  \[C=A\times B=\{(a,b):a\in A \wedge b\in B\}.\]
\end{definition}
Descartesov proizvod poznat je i pod nazivima Cartesiev ili direktni proizvod.
Ako skup $A$ ima $m$ elemenata, a skup $B$ ima $n$ elemenata, onda skupovi $A\times B$ i $B\times A$ imaju po $m\cdot n$ elemenata
(ako je $\card A=m$ i $\card B=n,$ onda je $\card A\times B=\card A \cdot \card B= m\cdot n$).
\begin{example}
  Dati su skupovi $A=\{a,b,c\},\,B=\{1,2\}.$ Odrediti $A\times B$ i $B\times A.$\\\\
\noindent Rje\v senje:
 \begin{align*}
   A\times B&=\{(a,1),\,(a,2),\,(b,1),\,(b,2),\,(c,1),\,(c,2)\};\\
   B\times A&=\{ (1,a),\,(1,b),\,(1,c),\,(2,a),\,(2,b),\,(2,c)\}.
 \end{align*}
Vidimo da je $\card A=3, \, \card B=2,$ a $\card A\times B= \card B\times A=6.$
\end{example}

\begin{figure}[!h]\centering
  \includegraphics[scale=1]{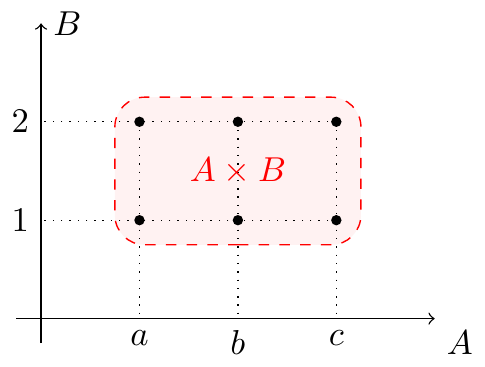}\hspace{1cm}
  \includegraphics[scale=1]{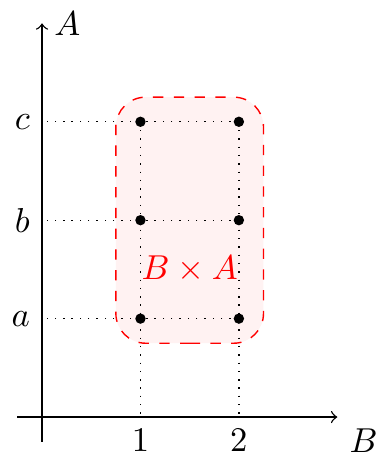}
  \caption{ Descartesovi (Dekartovi) proizvodi $A\times B$ i $B\times A$}
\end{figure}


\section[Z\lowercase{adaci}]{Zadaci}

 \begin{example}
	Dati su skupovi $A=\{x\in\mathbb{N}:1\leqslant x<4\},\:B=\{1,4\}.$ Odrediti\\
	\begin{inparaenum}
		\item $(A\cup B)\times(A\cap B);\:$
		\item $(A\cup B)\times(A\setminus B).$
	\end{inparaenum}\\ \\
\noindent 	Rje\v senje:\\ \ \\
	\begin{inparaenum}
		\item Vrijedi $A=\{1,2,3\},\:A\cup B=\{1,2,3,4\},\:A\cap B=\{1\},$ sada je\\
		$(A\cup B)\times(A\cap B)=\{1,2,3,4\}\times\{1\}=\{(1,1),(2,1),(3,1),(4,1)\}.$\\
		\item Vrijedi $A=\{1,2,3\},\:A\cup B=\{1,2,3,4\},\:(A\setminus B)=\{ 2,3 \},$ pa je sada
		\begin{align*}
		  (A\cup B)\times(A\setminus B)&=\{1,2,3,4\}\times\{2,3\} \\
		                               & =\{(1,2),(1,3),(2,2),(2,3),(3,2),(3,3),(4,2),(4,3)\}.
        \end{align*}		
	\end{inparaenum}
\end{example}

\begin{figure}[!h]\centering
  \includegraphics[scale=.95]{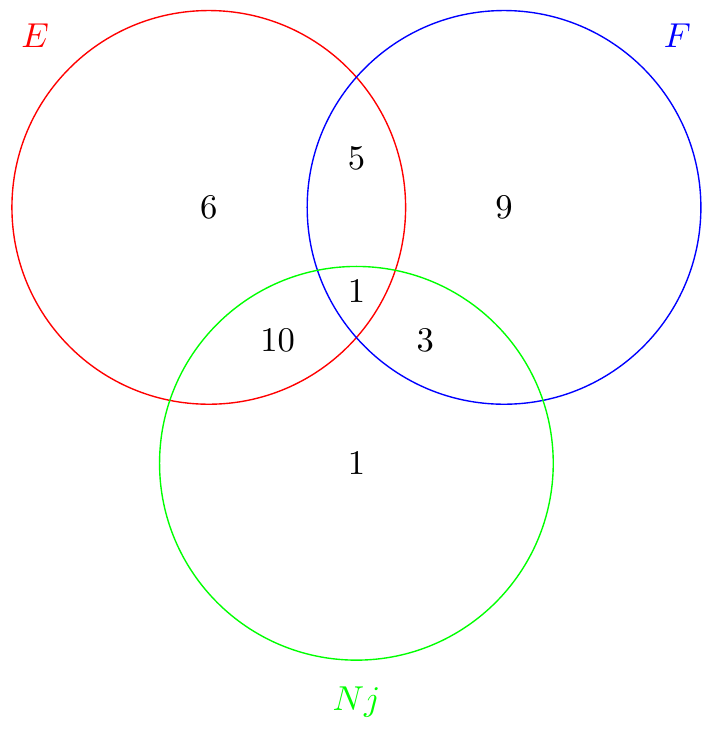}
  \caption{Polaznici kursa stranih jezika predstavljeni Euler--Vennovim dijagramom}
  \label{jezici}
\end{figure}

\begin{example}
  Na jednom kursu stranih jezika, svaki od 35 polaznika u\v ci bar jedan od tri strana jezika (francuski, njema\v cki i engleski) i to: 18 polaznika u\v ci francuski,   22 engleski, 6 polaznika u\v ci engleski i francuski, 11 engleski i njema\v cki, 4 francuski i njema\v cki   i jedan polaznik u\v ci sva tri jezika. Koliko polaznika u\v ci njema\v cki jezik?
  \newpage

\noindent Rje\v senje:\\

Ozna\v cimo sa $F,Nj$ i $E$ skupove polaznika koji u\v ce francuski, njema\v cki i engleski jezik, respektivno.\\
Ukupan broj polaznika je $\card(F\cup Nj\cup E)=35,$ dok jedan polaznik u\v ci sva tri jezika pa je  $\card( F \cap Nj\cap E)=1.$ Broj 1 pi\v semo u presjeku sva tri skupa (vidjeti Sliku \ref{jezici}), 6 polaznika u\v ci francuski i engleski ali jedan u\v ci i njema\v cki pa
$\card(E\cap F)-\card(E\cap F\cap Nj)=5$ pi\v semo u presjeku skupova $E$ i $F,$ 11 polaznika  u\v ci engleski i njema\v cki pa je
$\card(Nj\cap E)-\card(E\cap F\cap NJ)=10$ pi\v semo u presjeku skupova $E$ i $Nj,$ 4 polaznika u\v ce francuski i njema\v cki, pa sada $\card(F\cap Nj)-\card(E\cap F\cap Nj)=3$ i to pi\v semo u presjeku skupova $E$ i $Nj.$  Po\v sto ukupno 18 polaznika u\v ci francuski (samo francuski ili u kombinaciji sa ostalim jezicima), to je broj polaznika koji u\v ce samo francuski $\card F-\card(F\cap E)-\card(F\cap Nj)-\card(E\cap F\cap Nj)=9.$ Na isti na\v cin postupamo sa engleskim jezikom, 22 polaznika u\v ci engleski (ponovo sve mogu\' cnosti), sada je broj polaznika koji u\v ce samo engleski $\card E-\card(F\cap E)-\card(E\cap Nj)-\card(E\cap F\cap Nj)=6.$

Ukupan broj polaznika koji u\v ce engleski je $6+5+1+10=22$ (crvena kru\v znica), 9 polaznika u\v ci samo francuski (plava kru\v znica) i 3 polaznika u\v ce francuski i njema\v cki. Broj ovih polaznika je $22+9+3=34$, preostali polaznici su oni koji slu\v saju samo njema\v cki i njih je $35-34=1.$ I na kraju samo saberemo broj polaznika u zelenoj kru\v znici, $10+1+3+1=15,$ dakle ukupno je 15 polaznika koji slu\v saju njema\v cki (u bilo kojoj varijanti).

\end{example}


\subsection*{Zadaci za vje\v zbu}\index{Zadaci za vje\v zbu!skupovi}

\begin{enumerate}
   \item Odrediti uniju skupova $A$ i $B$ ako je \\
      \begin{inparaenum}
        \item $A=\{x\in \mathbb{Z}: x\leqslant 5 \wedge x>1 \},\:B=\{x\in\mathbb{Z}: -2\leqslant x\leqslant 3 \}; $
        \item $A=\{x\in \mathbb{R}: x\geqslant -1 \},\:B=\{x\in\mathbb{R}: x<1 \}; $
        \item $A=\{x\in \mathbb{Z}: x<3 \},\:B=\{x\in\mathbb{Z}: x>10 \}. $
      \end{inparaenum}
    \item Odrediti presjek kupova $A$ i $B$ ako je \\
      \begin{inparaenum}
        \item $A=\{x\in \mathbb{Z}: -4 \leqslant x<5 \},\:B=\{x\in\mathbb{N}: 3<x\leqslant 12 \}; $
        \item $A=\{x\in \mathbb{N}: 3\leqslant <10 \},\:B=\{x\in\mathbb{N}: 5\leqslant x<15 \}; $
        \item $A=\{x\in \mathbb{R}:6<x\leqslant 8 \},\:B=\{x\in\mathbb{R}: 3\leqslant y\leqslant 7 \}.$
      \end{inparaenum}
   \item Odrediti razlike $A\setminus B$  i $B\setminus A$ ako je \\
      \begin{inparaenum}
        \item $A=\{x\in \mathbb{Z}: 4<x\leqslant 9 \},\:B=\{x\in\mathbb{Z}: -2< x\leqslant 5 \}; $
        \item $A=\{x\in \mathbb{R}: x\leqslant 6 \},\:B=\{x\in\mathbb{R}: x>-2 \}.$
      \end{inparaenum}
   \item Dati su skupovi $A=\{1,2,3,4,5 \},\:B=\{1,3,5,7 \},\: C=\{2,4,6 \}.$ Odrediti\\
       \begin{inparaenum}
         \item $(A\cup B)\setminus C;\:$
         \item $(A\cap B)\setminus C;\:$
         \item $(A  \cup B)\setminus(B\cap C);\: $
         \item $(A\cup C)\setminus (B\setminus C),;\:$
         \item $(A\setminus B)\cup (A\setminus C).$
       \end{inparaenum}
   \item Odrediti komplement skupa $B$ u odnosu na skup $E=\{1,2,3,\ldots,15\}$ ako je \\
      \begin{inparaenum}
         \item $B=\{1,2,3\},\:$
         \item $B=\{x\in \mathbb{N }:x<8 \};\:$
         \item $B=\{x \in \mathbb{N} :x\leqslant 15 \};\:$
         \item $B=\{ 3,4,5,6,7,8,9\}.$
      \end{inparaenum}

   \item U jednoj porodici bilo je mnogo djece. Sedmoro njih voljelo je kupus, \v sestoro mrkvu, petoro krompir, \v cetvoro je voljelo kupus  i mrkvu, troje kupus i krompir, dvoje mrkvu i krompir i samo jedno dijete voljelo je sve mrkvu, krompir i kupus. Koliko je ukupno djece bilo u toj porodici?
   \item U jednom razredu 10 u\v cenika bavi se fudbalom, 12 \v sahom, 18 ko\v sarkom, 5 u\v cenika fudbalom i \v sahom, 7 fudbalom i ko\v sarkom, 8 \v sahom i ko\v sarkom i 3 u\v cenika \v sahom, fudbalom i ko\v sarkom. Koliko je sportista u razredu?
\item Odrediti partitivni  skup za date skupove\\
     \begin{inparaenum}
       \item $A=\{a,b\};\:$
       \item $A=\{1,2,3\};\:$
       \item $A=\{a,b,c,d\}.$
     \end{inparaenum}
\item Izra\v cunati\\
     \begin{inparaenum}
      \item $A\times B,\, B\times A,\,A\times A,\, B\times B,$ ako je $A=\{1,3,4\},\,B=\{a,c,g\};$

      \item $A\times B,\, B\times A,\,A\times A,\, B\times B,$ ako je $A=\{a,2,k\},\,B=\{1,2,b\}.$
    \end{inparaenum}
\item
    \begin{enumerate}
    	\item  Dati su skupovi $A=\{x: 3\leqslant x\leqslant 9\wedge x\in\mathbb{N}\},\: B=\{x:-2\leqslant x\leqslant 7\wedge x\in \mathbb{N}\}$ i $C=\{ x: x\text{ dijeli broj 12} \wedge x\in\mathbb{N}\}.$ Odrediti skupove $(A\cup C)\setminus(B\cap C)$ i $B\cap A.$
    	\item Dati su skupovi $A=\{x: 3\leqslant x\leqslant 9\wedge x\in\mathbb{N}\},\: B=\{x:-2\leqslant x\leqslant 7\wedge x\in \mathbb{N}\}$ i $C=\{ x: x\text{ dijeli broj 12} \wedge x\in\mathbb{N}\}.$ Odrediti skupove $(A\cup C)\setminus(B\cap C)$ i $B\cap A.$
    	\item  Dati su skupovi $A=\{a,b,c,d\},\,B=\{c,d\},\,C=\{1,2,3,4\}$ i $D=\{1,3,4,5,6,7\}.$\\ Odrediti
    	$(A\setminus B)\times(C\cap D).$
    \end{enumerate}

\end{enumerate}

\newpage

\chapter{Binarne relacije}
\pagestyle{fancy}

Razli\v cite veze izme\dj u elemenata skupova u matematici prou\v cavaju se u teoriji relacija. Pojam relacije spada me\dj u najop\v stije pojmove matematike.

\section[O\lowercase{snovni pojmovi}]{Osnovni pojmovi}

Izme\dj u elemenata jednog ili dva i vi\v se skupova, mogu postojati razli\v citi odnosi. U skupu realnih brojeva znamo da brojevi mogu da budu ve\' ci, manji ili jednaki; prave u geometriji mogu da budu paralelne, okomite,\ldots ; trouglovi mogu da budu podudarni, sli\v cni i tako dalje.

\begin{definition}[Binarna relacija]\index{relacija!binarna}
 Neka su $A$ i $B$ dva neprazna skupa. Svaki podskup Dekartovog proizvoda $A\times B$ naziva se binarna ili dvo\v clana relacija na skupu $A\times B$.
\end{definition}
Ako je $A=B$,  relacija na skupu $A\times A$ zove se i relacija u skupu $A$. Binarnu relaciju naj\v ce\v s\' ce ozna\v cavamo sa $\varrho.$ Po\v sto je binarna relacija podskup Dekartovog proizvoda $A\times B$ mo\v zemo pisati $\varrho\subset A\times B.$ Ako je $\varrho\subset A\times B$ i ako $(a,b)\in\varrho,$ gdje je $a\in A$ i $b\in B,$ onda ka\v zemo da su elementi $a$ i $b$ u relaciji $\varrho$ i pi\v semo $a\varrho b.$ Ako $a$ i $b$ nisu u relaciji $\varrho,$ onda pi\v semo  $\neg(a\varrho b),$ (mo\v ze i $a\neg\varrho b$)  ili $(a,b)\notin\varrho.$ Binarnu relaciju $\varrho$  simboli\v cki zapisujemo
\begin{empheq}[box=\mymath]{equation*}
\varrho\subset \left\{ (a,b)\in A\times B:a\in A\wedge b\in B\right\}.
\end{empheq}

 Relacije $=,\,<,\,>,\,\leqslant,\,\geqslant$ su binarne relacije u skupu realnih brojeva $\mathbb{R}.$ U skupu pravih imamo npr. relaciju paralelnosti, relaciju okomitosti i tako dalje.
\begin{remark}
  Relacije mo\v zemo prikazati kao skup ure\dj enih parova, osim toga i preko tabele, strelastog dijagrama--\v seme, koordinatnog sistema u ravni--mre\v ze (kao skup ta\v caka) ili preko Venovih dijagrama.
\end{remark}

\begin{figure}[!h]\centering
  \includegraphics[scale=.8]{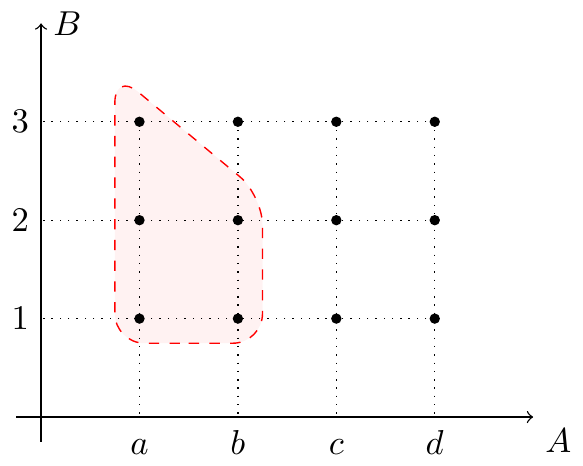}\hspace{.25cm}
  \includegraphics[scale=.8]{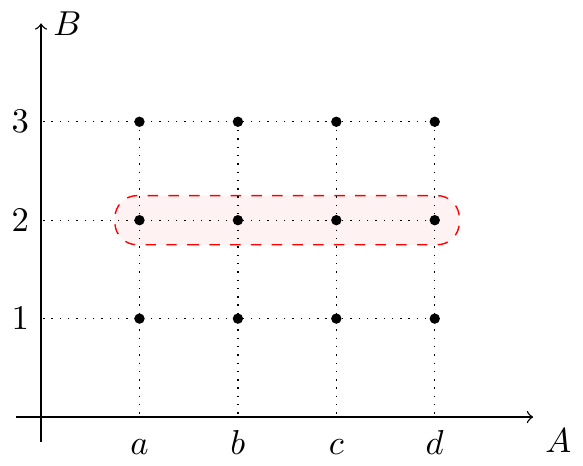}\hspace{.25cm}
  \includegraphics[scale=.8]{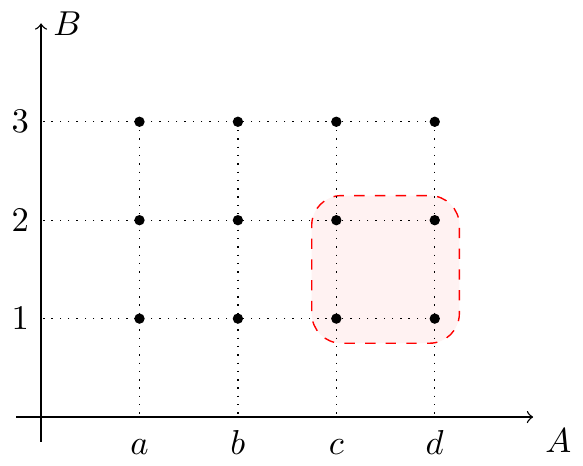}
  \caption{Primjeri binarnih relacija}
  \label{relacije1}
\end{figure}
\newpage
\begin{figure}[!h]\centering
   \includegraphics[scale=.8]{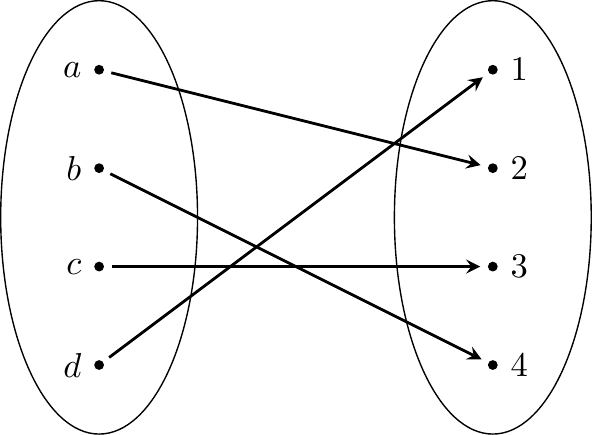}\hspace{1cm}
   \includegraphics[scale=1]{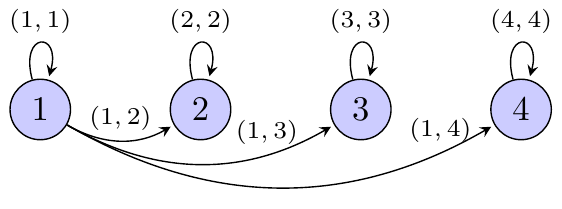}
   \caption{Primjeri predstavljanja binarnih relacija}
   \label{relacije2}
\end{figure}\index{relacija!predstavljanje}

Na Slici \ref{relacije1} predstavljeni su razni podskupovi Dekartovog proizvoda $A\times B.$ Svi ovi podskupovi su neke binarne relacije na $A\times B.$ Slike \ref{relacije1} i \ref{relacije2} predstavljaju razne na\v cine predstavljanja binarnih relacija.


Skup svih prvih koordinata relacije $\varrho\subset A\times B$ naziva se skup definisanosti ili domen relacije $\varrho$ i obilje\v zava se sa $D(\varrho)$  ili $D_\varrho,$ dok se skup svih drugih koordinata  ove relacije naziva skup vrijednosti relacije $\varrho$ ili kodomen i ozna\v cava $R(\varrho).$

\begin{example}
   Dati su skupovi $A=\{x\in \mathbb{N}:2\leqslant x\leqslant 5\},\,B=\{x\in\mathbb{N}:4\leqslant x\leqslant 6\}.$ Odrediti binarnu relaciju $\varrho=\{(a,b)\in A\times B:b=a+2\},$ zatim odrediti $D_\varrho$ i $R(\varrho).$\\\\
\noindent Rje\v senje: \\

\noindent Vrijedi $A=\{2,3,4,5\},\,B=\{4,5,6\},$ te je $\varrho=\{(2,4),(3,5),(4,6)\},$\\
\noindent i $D_\varrho=\{2,3,4\},\:R(\varrho)=\{4,5,6\}.$ Vidjeti sliku 3.3.
\end{example}

\begin{figure}[!h]\centering
 \includegraphics[scale=1]{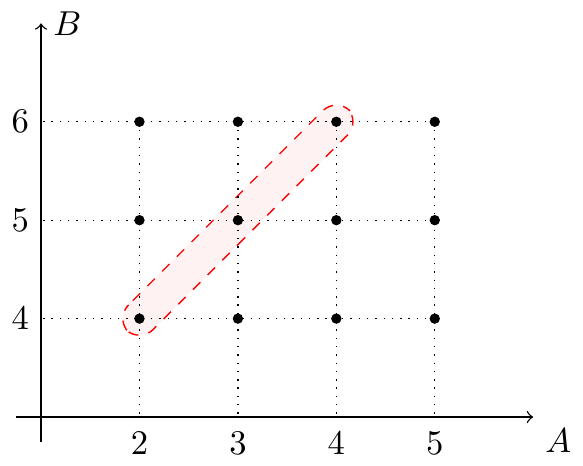}
 \caption{Relacija $\varrho$ predstavljena grafi\v cki}
 \label{relacije3}
\end{figure}

\begin{definition}[Presjek skupa $\varrho$ elementom $a$]
   Pod presjekom skupa $\varrho \subset A\times B$ elementom $a\in A$ podrazumijeva se skup svih elemenata $b\in B$ takvih da je $(a,b)\in\varrho,$ u oznaci $\varrho_a.$
\end{definition}
Vrijedi \[\varrho_a=\{b:b\in B\wedge (a,b)\in\varrho\}.\]
\begin{definition}[Skup svih presjeka]
   Skup svih presjeka relacije $\varrho$ naziva se faktor--skup skupa $B$ po relaciji $\varrho,$ u oznaci $\left.B\right|_\varrho.$
\end{definition}

\begin{example}
   Neka je $A=\{a_1,a_2,a_3,a_4,a_5\},\,B=\{b_1,b_2,b_3,b_4\}$ i\\
   $\varrho=\left\{ (a_1,b_1),(a_1,b_3),(a_2,b_1),(a_2,b_3),(a_2,b_4),(a_3,b_1),(a_3,b_2),(a_3,b_4),(a_4,b_3),\right.$\\
   $\left. (a_5,b_2),(a_5,b_4) \right\}.$ \\
Vrijedi   $\varrho_{a_1}=\{b_1,b_3\},$ $\varrho_{a_2}=\{b_1,b_3,b_4\},$ $\varrho_{a_3}=\{b_1,b_2,b_4\},$ $\varrho_{a_4}=\{b_3\}$ i $\varrho_{a_5}=\{b_2,b_4\}.$ Faktor--skup je \\
$\left.B\right|_{\varrho}=\{\varrho_{a_1}, \varrho_{a_2},\varrho_{a_3}, \varrho_{a_4}, \varrho_{a_5}\}=
\{\{b_1,b_3,\},  \{b_1,b_3,b_4\}, \{b_1,b_2,b_4 \}, \{b_3\}, \{b_2,b_4\}  \}.$
\end{example}

\begin{definition}[Inverzna relacija]
   Ako je $\varrho\subset A\times B,$ onda se pod inverznom relacijom podrazumijeva skup $\varrho^{-1}\subset B\times A$ pri \v cemu je
     \[\varrho^{-1}=\{(b,a):(a,b)\in\varrho \}.\]
\end{definition}\index{relacija!inverzna}
Drugim rije\v cima $b\varrho^{-1} a\Leftrightarrow a\varrho b.$

\begin{example}
    Dati su skupovi $A=\{-1,1,3,4\}$ i $B=\{0,-2,2\}.$ Odrediti $\varrho$ i $\varrho^{-1},$ ako je $\varrho=\{(a,b)\in A\times B:a+b=1\}.$\\\\
Rje\v senje: \\\\
Vrijedi $\varrho=\{(-1,2),(1,0),(3,-2) \}$ i $\varrho^{-1}=\{ (2,-1), (0,1),(-2,3)\}.$
\end{example}

\begin{definition}[Proizvod ili kompozicija relacija]
    Ako je $\varrho\subset A\times B$ i $r\subset B\times C,$ onda se pod proizvodom (ili kompozicijom) relacija $\varrho$ i $r$ podrazumijeva relacija definisana u skupu $A\times C,$  u oznaci $r \varrho$ takva da je
    \[a (r \varrho )c\Leftrightarrow (\exists x\in B)(a\varrho x\wedge xrc),\]
    gdje je $a\in A,x\in B$ i $c\in C.$
\end{definition}\index{relacija!kompozicija}

\begin{example}
   Neka je $A=\{1,2,3,4\},\,B=\{5,6\}$ i $C=\{1,2,4,7,8\},$ te neka vrijedi $\varrho=\{(1,5)\}$ i $r=\{(5,7),(6,8)\},$ onda je kompozicija $r \varrho =\{(1,7)\}.$
\end{example}

\begin{definition}[Relacija jednakosti]
 Pod relacijom jednakosti u skupu $A$ podrazumijeva se dijagonala skupa $A\times A$ i obilje\v zava se sa $\triangle.$
\end{definition}

\begin{figure}[!h]\centering
	\begin{subfigure}[b]{.45\textwidth}\centering
		\includegraphics[scale=.7]{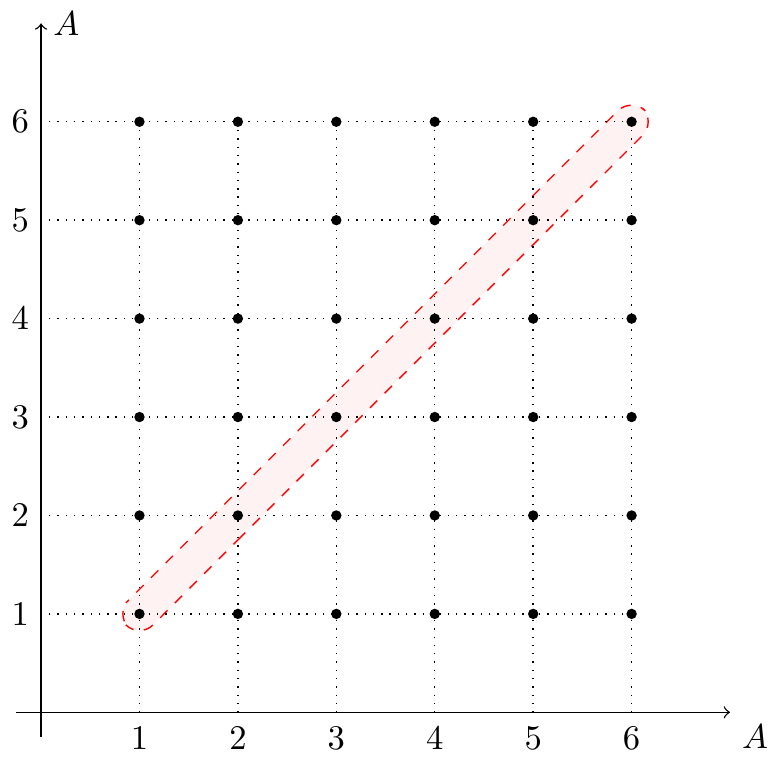}
		\caption{Mre\v za }
		\label{relacija4a}
	\end{subfigure}\hspace{.5cm}
	\begin{subfigure}[b]{.45\textwidth}\centering
		\includegraphics[scale=.7]{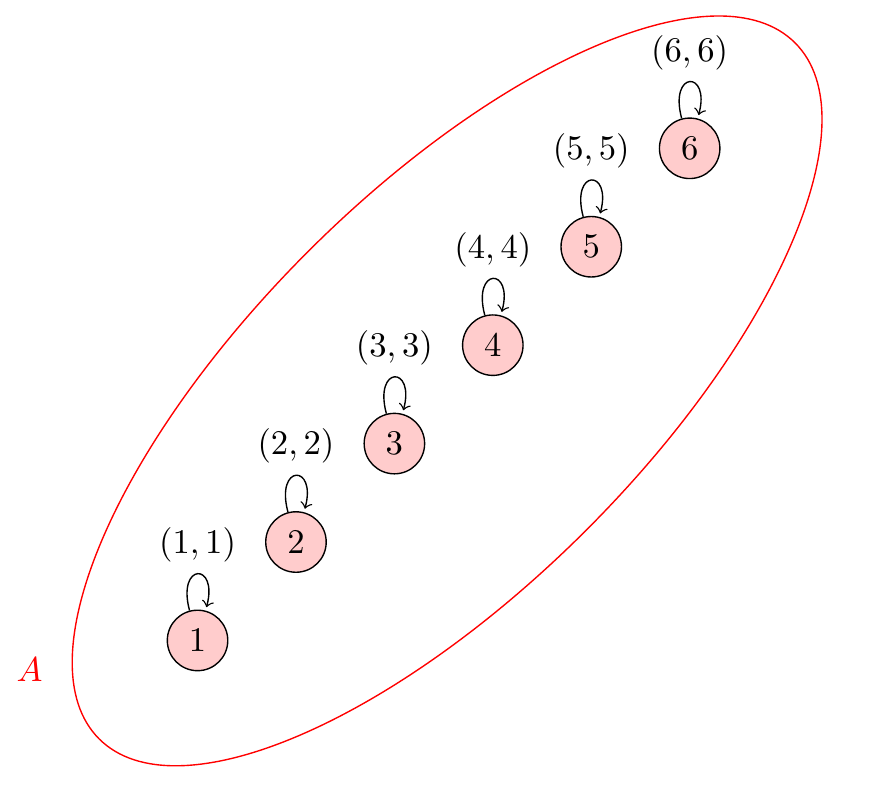}
		\caption{\v Sema}
		\label{relacija4b}
	\end{subfigure}
	\caption{Relacija jednakosti predstavljena mre\v zom i \v semom}
	\label{relacija4}
\end{figure}

\newpage
Dakle vrijedi \[\triangle=\{(a,a):a\in A\}.\]

\begin{example}
  Neka je $A=\{1,2,3,4,5,6\},$ vrijedi $\triangle=\{(1,1),(2,2),(3,3),(4,4),(5,5),(6,6)\},$ vidjeti Sliku \ref{relacija4}.
\end{example}

\section[O\lowercase{sobine binarnih relacija}]{Osobine binarnih relacija u jednom skupu}

 Neka je binarna relacija $\varrho$ definisana u skupu $A.$ Nave\v s\' cemo koje osobine mo\v ze imati relacija $\varrho $.
    \begin{enumerate}
         \item Relacija $\varrho$ je refleksivna ako je $(\forall a\in A)\:a\varrho a.$
         \item Relacija $\varrho$ je antirefleksivna ako je $(\forall a\in A)\: (a,a)\notin \varrho$.
         \item Relacija $\varrho$ je simetri\v cna ako je
                $(\forall a,b\in A)\:a\varrho b\Rightarrow b\varrho a.$
         \item Relacija $\varrho$ je antisimetri\v cna ako je
                $(\forall a,b\in A)\:(a\varrho b\wedge b\varrho a)\Rightarrow a=b.$
         \item Relacija $\varrho$ je tranzitivna ako je
                $(\forall a,b,c\in A)\:(a\varrho b\wedge b\varrho c)\Rightarrow a\varrho c.$
    \end{enumerate}
\index{relacija!osobine}

\begin{example}
   \begin{enumerate}
       \item Dat je skup $A=\{1,2,3,4\}$ i relacija $\varrho=\{(1,1),(2,2),(3,3),(4,3),(4,4)\}.$ Ova relacija je refleksivna.
       \item Dat je skup $A=\{1,4,5,7\}$ i relacija $\varrho=\{(1,1),(4,7),(7,4),(5,5)\}.$ Ova relacija nije refleksivna.
       \item Dat je skup $A=\{1,3,4,6\}$ i relacija $\varrho=\{(1,3),(3,4),(4,3),(3,1),(4,4),(6,6)\}.$ Ova relacija je simetri\v cna.
       \item Dat je skup $A=\{1,3,4\}$ i relacija $\varrho=\{(1,3),(3,4),(3,3),(4,4)\}.$ Ova relacija nije simetri\v cna.
       \item Dat je skup $A=\{1,3,5,7\}.$ \\
           \begin{inparaenum}
              \item $\varrho=\{ (1,5),(3,3),(5,5),(5,7),(7,3)\}$ je antisimetri\v cna.\\
              \item $\varrho=\{(1,1),(3,3),(5,5),(7,7) \}$ je i simetri\v cna i antisimetri\v cna.\\
              \item $\varrho=\{(1,5),(3,5),(5,1),(5,5) \}$ nije ni antisimetri\v cna ni simetri\v cna.
           \end{inparaenum}
       \item Dat je skup $A=\{2,3,4,5\}$ i relacija $\varrho=\{(2,3),(3,4),(2,4),(4,5),(3,5),(2,5)\}.$ Ova relacija je tranzitivna.
  \end{enumerate}
\end{example}

\subsection{Relacija ekvivalencije} \index{relacija!ekvivalencije}
Relacija ekvivalencije ima veliki zna\v caj u matematici.

\begin{definition}[Relacija ekvivalencije]
   Binarna relacija $\varrho\subset A\times A$ naziva se relacija ekvivalencije ako je refleksivna, simetri\v cna i tranzitivna.
\end{definition}
Relacija ekvivalencije ozna\v cava se sa $\sim,$ pa umjesto $a\varrho b$ pi\v semo $a\sim b$ i \v citamo "$a$ je ekvivalentno $b"$. Jednakost brojeva je jedna relacija ekvivalencije, zatim paralelnost pravih, podudarnost trouglova i tako dalje.

\begin{definition}[Klase  ekvivalencije]
  Skup svih elemenata skupa $ A$ koji su u relaciji ekvivalencije sa elementom $a$ naziva se klasa ekvivalencije i obilje\v zava sa $C_a,$  simboli\v cki
  \[C_a=\{ x\in A:x\sim a\}. \]
\end{definition}
Element $a$ nazivamo predstavnik klase $C_a.$
\begin{definition}[Faktorski ili kvocijentni (koli\v cni\v cki) skup]
   Skup svih klasa ekvivalencije nekog skupa $A$ zove se faktorski skup i obilje\v zava se sa $ \left. A\right|_\sim .$
   \[ \left. A\right|_\sim= \{C_a : a\in A \}. \]
\end{definition}

\begin{example}
  U skupu $A=\{16,17,18,19,20,21,22,23,24,25,25,26,27,28\}$ zadana je relacija $\varrho :$ imati jednak zbir cifara. Ispitati da li je $\varrho$ relacija ekvivalencije. Ako jeste  odrediti klase i faktorski skup. \\\\
  Rje\v senje.\\\\
Vrijedi $\varrho=\left\{ (16,16),\ldots,(28,28),(16,25),(25,16),(17,26),(26,17),(18,27),\right.$\\
$\left.(27,18),(19,28),(28,19)\right\}.$\\ Vidimo da je relacija $\varrho$ refleksivna, simetri\v cna i tranzitivna, dakle $\varrho$ je relacija ekvivalencije. Klase su $C_{16}=\{16,25\},\,C_{17}=\{17,26\},\, C_{18}=\{18,27\},\, C_{19}=\{19,28\},\, C_{20}=\{ 20\},\, C_{21}=\{ 21\},\,  C_{22}=\{ 22\},\,C_{23}=\{ 23\},\,C_{24}=\{ 24\}.$\\ Koli\v cni\v cki skup je\\
$\left. A\right|_{\sim}=\{\{16,25\},\, \{ 17,26\},\, \{ 18,27\},\,\{19,28 \} ,\,\{20 \},\, \{ 21\},\,\{ 22\},\,\{ 23\},\{ 24\}\}. $
\end{example}

\begin{figure}[!h]\centering
     \begin{subfigure}[b]{.45\textwidth}\centering
        \includegraphics[scale=.75]{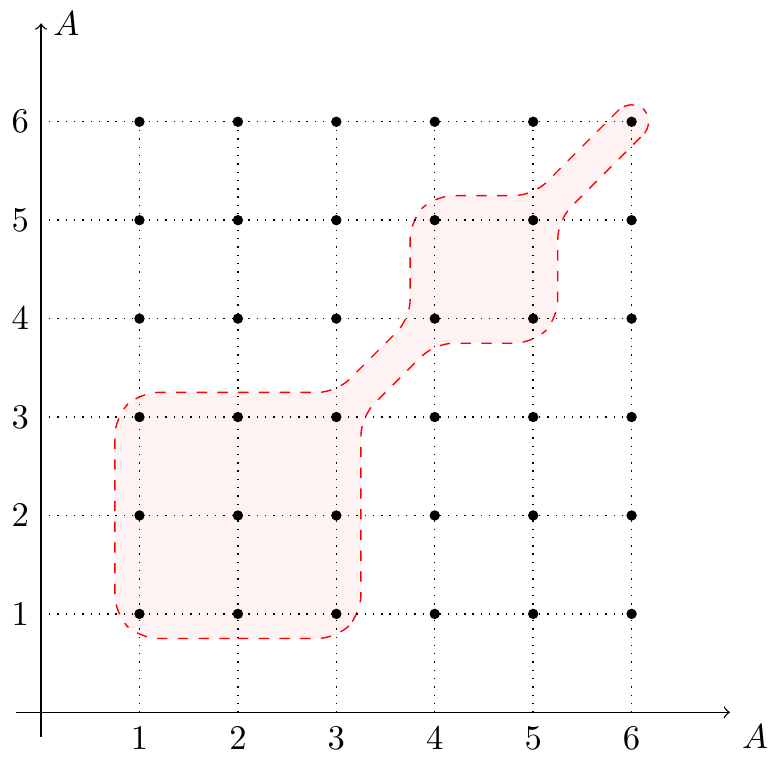}
         \caption{Mre\v za }
          \label{relacija8a}
   \end{subfigure}\hspace{.5cm}
   \begin{subfigure}[b]{.45\textwidth}\hspace{-.5cm}
        \includegraphics[scale=.75]{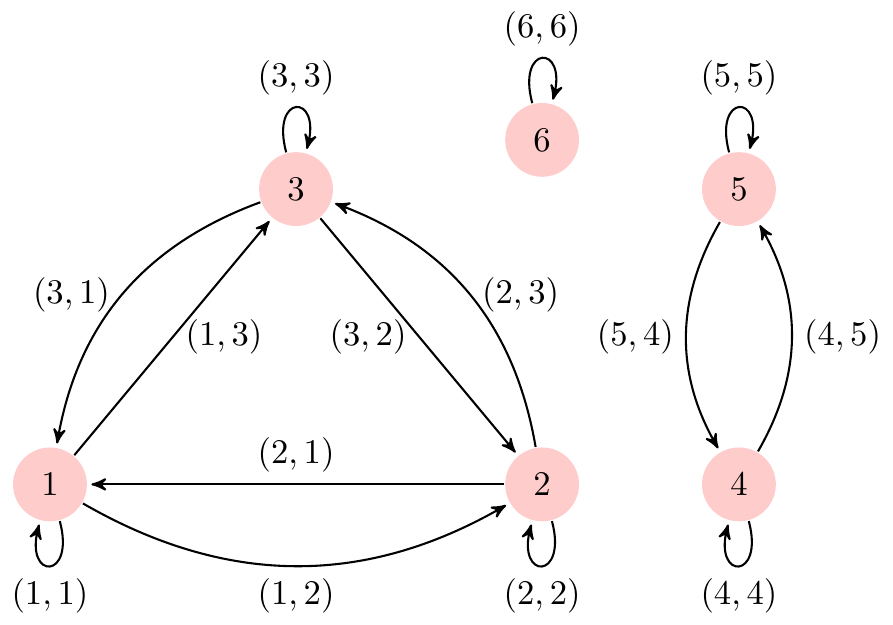}
         \caption{\v Sema}
      \label{relacija8b}
     \end{subfigure}
        \caption{Relacija ekvivalencije predstavljena mre\v zom i \v semom}
     \label{relacija8}
  \end{figure}

\begin{example}
  Dat je skup $A=\{1,2,3,4,5,6\}$ i binarna relacija \\$\varrho=\left\{(1,1),(2,2),(3,3),(4,4),(5,5),(6,6),(1,2),(2,1),(1,3),(3,1),(2,3),(3,2),\right.$\\
  $\left.(5,4),(4,5) \right\}.$ Ispitati da li je relacija $\varrho$ relacija ekvivalencije, ako jeste odrediti klase, faktorski skup, nacrtati mre\v zu i \v semu.\\\\
Rje\v senje.\\\\ Relacija $\varrho$ ima osobine refleksivnosti, simetri\v cnosti i tranzitivnosti. Klase su
$C_1=\{1,2,3\},\,C_4=\{4,5\}$ i $C_6=\{6\}.$\\  Koli\v cni\v cki skup je $\left.A\right|_\varrho=\{C_1,C_4,C_6\}=\{ \{1,2,3\},\,\{4,5\},\,\{6\}   \}.$ Vidjeti Sliku \ref{relacija8}.
\end{example}

\subsection{Relacija poretka}\index{relacija!poretka}
Osim relacija ekvivalencije, tako\dj e i relacija poretka ima bitnu ulogu u matematici.
\begin{definition}[Relacija poretka]
  Binarna relacija $\varrho\subset A\times A$ naziva se relacija poretka ako je refleksivna, antisimetri\v cna i tranzitivna.
\end{definition}

Relaciju poretka obilje\v zavamo sa $\leqslant$ i \v citamo "manje ili jednako". Ako je $a\leqslant b$ ili $b\leqslant a$ onda se ka\v ze da su elementi $a$ i $b$ uporedivi. Ako na skupu $A$ postoji relacija poretka $\leqslant$ onda se ka\v ze da je skup $A$ ure\dj en relacijom poretka $\leqslant$ i pi\v semo $(A,\leqslant).$

Ako su svi elementi skupa $A$ uporedivi, onda se on naziva potpuno ure\dj en ili lanac ili linearno ure\dj en. U slu\v caju da svi elementi skupa $A$  nisu uporedivi onda se ka\v ze da je $A$ djelimi\v cno ure\dj en skup.

\begin{definition}[Relacija strogog poretka]
  Binarna relacija koja posjeduje osobine tranzitivnosti i antirefleksivnosti, naziva se relacija strogog poretka.
\end{definition}
Relacija $<$ na skupu realnih brojeva $\mathbb{R}$ je relacija strogog poretka.
Ako je $a<b$ ka\v ze se i da je $a$ ispred $b$ ili $b$ je iza $a.$ Ako je $a<b$ i $b<c$ onda ka\v zemo da je $b$ izme\dj u $a$ i $c.$
Kod skupova koristimo relaciju podskup $\subseteq $ (relacija poretka) i relaciju strogi poskup $ \subset $ (strogi poredak).

\begin{figure}[!h]\centering
     \begin{subfigure}[b]{.45\textwidth}\centering
        \includegraphics[scale=.9]{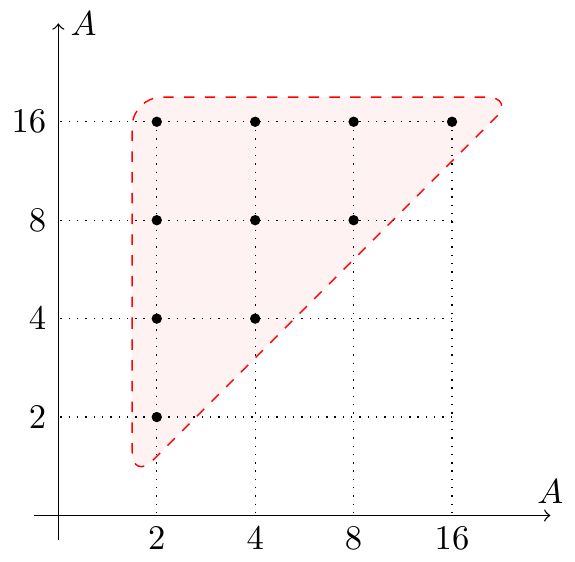}
         \caption{Mre\v za }
          \label{relacija10a}
   \end{subfigure}
   \begin{subfigure}[b]{.45\textwidth}\centering
        \includegraphics[scale=.9]{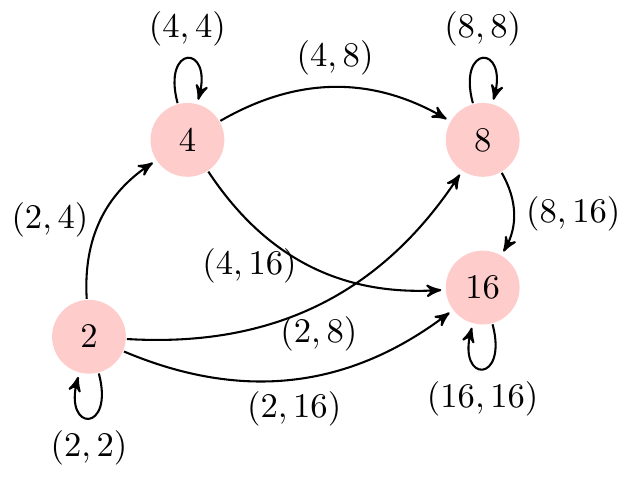}
         \caption{\v Sema}
      \label{relacija10b}
     \end{subfigure}
        \caption{Relacija poretka predstavljena mre\v zom i \v semom}
     \label{relacija10}
  \end{figure}

\begin{example}
 U skupu $A=\{2,4,8,16\}$ definisana je operacija $\varrho$ sa $a\varrho b\Leftrightarrow a\mid b $ ($a$ dijeli $b$).
  \begin{enumerate}
     \item Napisati relaciju $\varrho$ kao skup ure\dj enih parova, nacrtati mre\v zu i \v semu.
     \item Ispitati da li je $\varrho$ relacija poretka
  \end{enumerate}
Rje\v senje:
 \begin{enumerate}
   \item Vrijedi\\
        $\varrho=\{(2,2),(2,4), (2,8),(2,16),(4,4),(4,8),(4,16),(8,8),(8,16),(16,16)\}.$ Vi-\\djeti Sliku \ref{relacija10}.
   \item Operacija $\varrho$ je refleksivna, antisimetri\v cna i tranzitivna, dakle $\varrho$ je relacija poretka.
  \end{enumerate}
\end{example}

\section{Z\lowercase{adaci}}\index{Zadaci za vje\v zbu!relacije}

\begin{enumerate}
   \item U skupu $A=\{1,2,3,4,5,6,7 \}$ zadana je relacija sa $\varrho=\{(a,b)\in A\times A: b\geqslant a+1 \}.$ \\
     \begin{inparaenum}
        \item Odrediti relaciju $\varrho$ kao skup ure\dj enih parova. \\
        \item Predstaviti relaciju $\varrho$ grafi\v cki (u koordinatnom sistemu--mre\v zi i u strelastom dijagramu--\v semi). \\
        \item Odrediti domen $D_\varrho$ i kodomen relacije $R_\varrho.$\\
        \item Odrediti presjeke relacije sa svakim elementom i uniju tih presjeka.
    \end{inparaenum}

   \item Date su relacije $\varrho=\{(1,a),(2,b),(2,c),(3,e)\}$ i $r=\{ (a,a),(a,b),(b,c),(c,c)\}.$  Odrediti
         kompoziciju relacija $\varrho $ i $r.$
   \item Dati su skupovi $A=\{x\in\mathbb{Z}:-2\leqslant x\leqslant 2\},\:B=\{x\in\mathbb{Z}:-1\leqslant x\leqslant \},\:
         C=\{x\in\mathbb{N}:x<5\}.$ Odrediti\\
         \begin{inparaenum}
           \item Relacije $\varrho=\{(a,b)\in A\times B: ab>0\},$ $r=\{ (b,c)\in B\times C:b+c<4\}.$\\
           \item Inverzne relacije $\varrho^{-1}$ i $r^{-1}.$\\
           \item Kompozicije relacija $\varrho r$ i $(\varrho r)^{-1}.$
\end{inparaenum}
   \item U skupu $A=\{-2,-1,0,1,2,3 \}$ definisana je relacija $\varrho$ sa $x\varrho y\Leftrightarrow x+y=2.$ \\
       \begin{inparaenum}
          \item Napisati relaciju $\varrho$ (kao skup ure\dj enih parova).\\
          \item Predstaviti relaciju $\varrho$ grafi\v cki (u koordinatnom sistemu--mre\v zi i u strelastom dijagramu--\v semi).\\
          \item Koje osobine ima relacija $\varrho?$
       \end{inparaenum}
    \item U skupu $A=\{-2,-1,0,1,2\}$ definisana je relacija $\varrho$ sa $x\varrho y\Leftrightarrow x^2=y^2.$ \\
       \begin{inparaenum}
          \item Napisati relaciju $\varrho$ (kao skup ure\dj enih parova).\\
          \item Predstaviti relaciju $\varrho$ grafi\v cki (u koordinatnom sistemu--mre\v zi i u strelastom dijagramu--\v semi).\\
          \item Koje osobine ima relacija $\varrho?$ \\
          \item Ako je relacija $\varrho$ relacija ekvivalencije odrediti klase i koli\v cni\v cki skup.
       \end{inparaenum}
     \item U skupu $A=\{1,2,3,4,5,6,7\}$ za relaciju $\varrho$ date su njene klase $\{1,5,6\}\:\{2\},\:\{3,7\},\:\{4\}.$ \\
       \begin{inparaenum}
          \item Napisati relaciju $\varrho$ (kao skup ure\dj enih parova).\\
          \item Predstaviti relaciju $\varrho$ grafi\v cki (u koordinatnom sistemu--mre\v zi i u strelastom dijagramu--\v semi).\\
          \item Koje osobine ima relacija $\varrho?$\\
          \item \v Sta je koli\v cni\v cki skup $ \left.A\right|_{\varrho} ?$
       \end{inparaenum}
     \item U skupu $A=\{15,16,17,21,24,25,33,34,37\}$ definisana je relacija $\varrho$ sa
       $x\varrho y\Leftrightarrow x \text{ i }y \text{ imaju isti zbir cifara.}$  \\
        \begin{inparaenum}
           \item Napisati relaciju $\varrho$ kao skup ure\dj enih parova.\\
           \item Predstaviti relaciju $\varrho$ grafi\v cki (u koordinatnom sistemu--mre\v zi i u strelastom dijagramu--\v semi). \\
           \item Navesti osobine relacije $\varrho.$
        \end{inparaenum}

    \item U skupu $A=\{3,9,27,81\}$ definisana je operacija $\varrho$ sa $a\varrho b\Leftrightarrow a\mid b$. \\
    \begin{inparaenum}[($i$)]
    	\item Napisati relaciju $\varrho$ kao skup ure\dj enih parova.\\
    	\item Predstaviti relaciju $\varrho$ grafi\v cki (u koordinatnom sistemu--mre\v zi i u strelastom dijagramu--\v semi).
    \end{inparaenum}
      \item U skupu $a=\{5,10,20,30,40\}$ zadana je relacija $\varrho$ sa $x\varrho y\Leftrightarrow x\mid y.$ \\
         \begin{inparaenum}
           \item Pokazati da je $\varrho$ relacija poretka.\\
           \item Predstaviti relaciju $\varrho$ grafi\v cki (u koordinatnom sistemu--mre\v zi i u strelastom dijagramu--\v semi).
        \end{inparaenum}

    \item  Dati su skupovi $A=\{-3,-2,-1,0,1,2\},\:B=\{-1,0,1,2\}$ i $C=\{1,2,3,4\}$ i relacije $\varrho=\{ (x,y)\subseteq A\times B:x\varrho y\Leftrightarrow xy>0\}$ i $\psi=\{(x,y)\subseteq B\times C:x \psi y\Leftrightarrow x+y<4\}.$ Napisati sve elemente relacija $\varrho$ i $\psi,$ te odrediti relaciju $(\varrho \circ \psi)^{-1}.$

\end{enumerate}

\chapter{Binarne operacije}
\pagestyle{fancy}

\begin{definition}\index{binarna operacija}
  Neka je $S$ neprazan skup. Binarna operacija u skupu $S$ je postupak kojim se elementima skupa $S\times S$ jednozna\v cno pridru\v zuje jedan element iz skupa $S$.
\end{definition}
Binarna operacija zove se tako\dj e i binarna kompozicija, apstraktna operacija ili samo operacija.  Znak koji pokazuje da na elemente $a$ i $b$ skupa $S$ treba primijeniti odre\dj eni postupak da bi se dobio odre\dj eni element $c$ iz istog skupa $S$ zove se operator i \v cesto se ozna\v cava sa $\circ.$  Rezultat operacije  $\circ$ koja je izvr\v sena sa elementima $a$ i $b$ ozna\v cava se $a\circ b=c.$

Primjeri poznatih oznaka operatora i binarnih operacija su:
 \begin{enumerate}[\null]
 \item $+$ za sabiranje (brojeva, polinoma, vektora,\ldots);
 \item $-$ za oduzimanje (brojeva, polinoma, vektora,\ldots);
 \item $\times$ ili $\cdot$ za mno\v zenje (brojeva, polinoma, vektora,\ldots);
 \item $:$ ili $/$ za dijeljenje (brojeva, polinoma, \ldots);
 \item $\cup$ za uniju skupova;
 \item $\cap$ za presjek skupova;
 \item $\setminus$ za razliku skupova;\ldots
\end{enumerate}

\section[O\lowercase{sobine binarnih operacija}]{Osobine binarnih operacija}
\index{binarna operacija!osobine}
\begin{definition}
Ka\v zemo da je binarna operacija $\circ $ u skupu $S$ zatvorena ili unutra\v snja ako je za bilo koja dva elementa $a,b \in S$ rezultat operacije element iz $S,$ odnosno
 \[(\forall a,b \in S) \, a \circ b \in S.\]
\end{definition}

\begin{definition}\index{binarna operacija!grupoid} \label{definicijaGrupoid}
 Skup $S$ i zatvorena binarna operacija $\circ$ u skupu $S$ \v cine algebarsku strukturu koja se zove grupoid. Grupoid se ozna\v cava sa $(S,\circ).$
\end{definition}
 Jasno je da je u skupu prirodnih brojeva $\mathbb{N}=\{1,2,3, \ldots \}$ operacija sabiranja zatvorena, dakle $\left(\mathbb{N}, + \right)$ je grupoid.
 Sa druge strane, operacija oduzimanja nije zatvorena, jer razlika dva prirodna broja nije uvijek prirodan broj. To zna\v ci da $\left(\mathbb{N}, - \right)$
 nije grupoid.
Za binarnu operaciju mo\v zemo definisati jo\v s neke osobine.


\begin{definition}[Asocijativnost]
  Neka je dat grupoid $(S,\circ)$. Ako za svaka tri elementa $a,b,c\in S$ va\v zi
   \[(a\circ b)\circ c=a\circ(b\circ c),\]
   ka\v ze se da je operacija $\circ$ asocijativna na skupu $S$.
\end{definition}

\begin{definition}[Komutativnost]
 Neka je dat grupoid $(S,\circ)$. Ako za svaka dva elementa $a,b\in S$ va\v zi jednakost
 \[a\circ b=b\circ a,\]
 ka\v zemo da je operacija $\circ$ komutativna na $S$.
\end{definition}

\begin{definition}[Neutralni element]
  Ako u grupoidu $(S,\circ)$ postoji takav element $e\in S$ za koji je
  \[(\forall a\in S)\:a\circ e=e\circ a=a\]
  on se zove neutralni ili jedini\v cni element.
\end{definition}

\begin{proposition}
 Ako u grupoidu $(S,\circ)$ postoji neutralni element $e$ za binarnu operaciju $\circ,$ onda je on jedinstven.
\end{proposition}
\begin{proof}[Dokaz]
Neka su $e, e' \in \mathbb{S}$ neutralni elementi za binarnu operaciju $\circ .$ Kako je $e$ neutralni element, vrijedi:
\[e\circ e'=e'.\]
Sa druge strane, s obzirom da je $e'$ neutralni element, imamo:
\[e\circ e'=e. \]
Prema tome, zaklju\v cujemo da je $e=e',$ tj. neutralni element je jedinstven.
\end{proof}

\begin{definition}[Inverzni element]
 Neka grupoid $(S,\circ)$ ima neutralni element $e$. Ka\v zemo da element $a\in S$ ima sebi inverzni element (inverz), ili suprotan element $a'\in S$, ako vrijedi jednakost
   \[a'\circ a=a\circ a'=e.\]
\end{definition}

\begin{proposition}
 Neka je $(S,\circ)$ grupoid u kojem je binarna operacija $\circ $ asocijativna i ima neutralni element $e.$ Ako je $a \in S$ invertibilan element,
 onda je njegov inverz jedinstven.
\end{proposition}
\begin{proof}[Dokaz]
Neka su $y$ i $y'$ inverzi od $a.$ Tada vrijedi:
\[ y\circ a=a \circ y=e, \quad y' \circ a=a \circ y'=e.  \]
Zbog asocijativnosti operacije $\circ$ imamo
\[ y\circ a \circ y'=(y \circ a) \circ y'=e \circ y'=y' ; \]
\[ y\circ a \circ y'=y \circ (a \circ y')=y \circ e =y.  \]
Iz ove dvije jednakosti slijedi da je $y=y'.$ Dakle, pokazali smo da je invertibilni element jedinstven.
\end{proof}

\begin{example}
   U skupu $\mathbb{Z}$ definisane su binarne operacije $\star$ i $\circ,$ na sljede\' ci na\v cin $x\star y=3x-2y+1$ i $x\circ y=2x-3y+1.$ Odrediti \[((-1\star 2)\circ 5)\star(2\circ 2).\]  \\
   \noindent Rje\v senje. Skup $\mathbb{Z}$ je skup cijelih brojeva, tj. $\mathbb{Z}=\{\ldots,-3,-2,-1,0,1,2,3,\ldots\}.$ Vrijedi
   \begin{align*}
   -1\star 2&=3\cdot(-1)-2\cdot 2+1=-6\\
   2\circ 2 &=2\cdot 2-3\cdot 2+1=-1 \\
   (-1\star 2)\circ 5&=(-6)\circ 5=2\cdot (-6)-3\cdot 5+1=-26,
   \end{align*}
   pa je sada
   \begin{align*}
   ((-1\star 2)\circ 5)\star(2\circ 2)&=-26\star (-1)=3\cdot(-26)-2\cdot (-1)+1=-75.
   \end{align*}
\end{example}
Za pregledniji prikaz rezultata binarne operacije definisane u nekom skupu $S,$ mo\v ze se formirati tabela koja nosi naziv Cayleyeva\footnote{Arthur Cayley (16.avgust 1821.--26.januar 1895. godine) bio je britanski matemati\v car koji je uglavnom radio na algebri.}  tablica (tabela),
\v sto \' cemo prikazati u nare-\\dnom primjeru.
\begin{example}
	Dati su skup $S=\{-1,1,-i,i\}$ i operacija mno\v zenja $\bullet.$ Formirati Cayleyjevu tabelu za operaciju $\bullet.$ \\
	\noindent Rje\v senje. Po\v sto je
	\begin{align*}
	 &(-1)\bullet (-1)=1,\: (-1)\bullet (-i)=(-i)\bullet (-1)=i,\: (-1)\bullet i=i\bullet (-1)=-i,\\
	 &1\cdot (-i)=(-i)\bullet 1=-i,\: 1\bullet i=i\bullet 1=i,\\
	 & (-i)\bullet i=i\bullet (-i)=1,
	\end{align*}
	dobijamo Tabelu \ref{kejli1}.
\end{example}

\begin{center}
	\begin{table}[!h]\centering
		\begin{tabular}{r|rrrr} 
			$\bullet$  & $-1$ & $1$  & $i$  &$-i$ \\    \hline
			$-1$  & $1$  & $-1$ & $-i$ & $i$ \\  
			$1$   & $-1$ & $1$  & $i$  & $-i$ \\ 
			$i$   & $-i$ & $i$  & $-1$ & $1$ \\
			$-i$  & $i$  & $-i$ & $1$  & $-1$
		\end{tabular}
		\caption{Cayleyjeva tabela}
		\label{kejli1}
	\end{table}
\end{center}

\section[A\lowercase{lgebarske strukture}]{Algebarske strukture sa jednom i dvije operacije}
\index{algebarske strukture!grupa}

Vidjeli smo da je u slu\v caju osnovne algebarske strukture grupoid $(S, \circ),$ binarna operacija $\circ $ zatvorena. Ako ta operacija ima jo\v s
neke osobine, onda dobijamo bogatije algebarske strukture.
\begin{definition}[Polugrupa]
 Grupoid $(S,\circ)$ u kome je operacija $\circ$ asocijativna naziva se polugrupa (ili semigrupa).
\end{definition}
 \begin{definition}[Grupa]
    Ka\v ze se da grupoid $(S,\circ)$ ima strukturu grupe ili da \v cini grupu, ako operacija $\circ$ ima i ove tri osobine:
     \begin{enumerate}[$1^0$]
       \item Operacija je asocijativna \[(\forall a,b,c\in S)\:(a\circ b)\circ c=a\circ(b\circ c);\]
       \item Operacija ima neutralni element $e$
         \[(\exists e\in S)(\forall a\in S)\:a\circ e=e\circ a=a;\]
       \item  Svaki element $a$ ima simetri\v can element $a'$ za operaciju $\circ$
         \[(\forall a\in S)(\exists a'\in S)\: a\circ a'=a'\circ a=e.\]
     \end{enumerate}
Ako je $(S,\circ)$ grupa, ka\v ze se  da elementi skupa $S$ obrazuju grupu.
 \end{definition}
\index{algebarske strukture!komutativna ili Abelova grupa}Ako je operacija $\circ$ jo\v s i komutativna dobijamo algebarsku stukturu nazvanu po Abelu\footnote{Niels Henrik Abel (5.avgust 1802.--6.april 1829. godine) bio je norve\v ski matemati\v car. Dao pionirski doprinos u teoriji polja.}.
\begin{definition}[Komutativna ili Abelova grupa]
   Ako je $(S,\circ)$ grupa i ako je operacija $\circ$ komutativna
   \[( \forall a,b\in S)\:a\circ b=b\circ a,\]
   ka\v ze se da je grupa komutativna ili Abelova.
\end{definition}

Operacija sabiranja u skupu prirodnih brojeva je asocijativna, pa je $(\mathbb{N},+)$ ne samo grupoid nego i polugrupa, a kako sabiranje ima i osobinu komutativnosti, ka\v zemo da je to komutativna polugrupa. Ne postoji neutralni element za sabiranje u ovom skupu, pa ni inverzni elementi.
U skupu cijelih brojeva $\mathbb{Z}$ postoji neutralni element za sabiranje, to je broj $0,$ a za svaki cijeli broj $a$ postoji inverzni (suprotni) element $-a$, pa je $(\mathbb{Z},+)$ komutativna ili Abelova grupa.

\begin{example}
	Dat je skup $A=\{0,1,2,3,4,5\}$ i u njemu operacija $+_6$ definisana sa: \\
	$a+_6 b=c,$ gdje je $c$ je ostatak pri dijeljenju broja $a+b$ sa $6.$
	\begin{enumerate}[$(a)$]
		\item Formirati Cayleyjevu tabelu.
		\item Koju algebarsku strukturu predstavlja $(A,+_6)?$
	\end{enumerate}
	Rje\v senje:
	\begin{inparaenum}[$(a)$]
		\item Vidjeti Tabelu \ref{kejli2}.
		\item $(A,+_6)$ je komutativna grupa.
	\end{inparaenum}
\end{example}

\begin{table}[!h]\centering
	\begin{tabular}{c|cccccc}\\
		$+_6$& $0$ & $1$ &$2$ & $3$ &$4$ &$5$\\     \hline
		$0$ &  $0$ & $1$ & $2$ & $3$ & $4$ & $5$ \\
		$1$ &  $1$ & $2$ & $3$ & $4$ & $5$ & $0$ \\
		$2$ &  $2$ & $3$ & $4$ & $5$ & $0$ & $1$ \\
		$3$ &  $3$ & $4$ & $5$ & $0$ & $1$ & $2$ \\
		$4$ &  $4$ & $5$ & $0$ & $1$ & $2$ & $3$ \\
		$5$ &  $5$ & $0$ & $1$ & $2$ & $3$ & $4$ \\
	\end{tabular}
	\caption{Cayleyjeva tabela za $(A,+_6)$}
	\label{kejli2}
\end{table}

\begin{example}
	U skupu $\mathbb{Z}$ definisana je binarna operacija $\star$ sa $x\star y=x+y-4.$ Dokazati da je $(\mathbb{Z},\star)$ komutativna grupa.\\
	\noindent Rje\v senje.
	\begin{enumerate}
		\item\label{grupa1} Poka\v zimo prvo da vrijedi $\forall x,y\in\mathbb{Z}\Rightarrow x\star y\in\mathbb{Z}.$ Zbir dva cijela broja je cio broj tj. $\forall x,y\in\mathbb{Z} \Rightarrow x+y\in\mathbb{Z},$ takodje je i $4\in\mathbb{Z},$ pa je $x+y-4\in\mathbb{Z},$ drugim rije\v cima $\forall x,y\in\mathbb{Z}\Rightarrow x\star y\in\mathbb{Z}.$ Prethodno pokazano zna\v ci, na osnovu Definicije  \ref{definicijaGrupoid},  da je $(\mathbb{Z},\star)$ grupoid.
		\item\label{grupa2} Poka\v zimo da vrijedi asocijativni zakon za operaciju $\star,$ $\forall x,y,z\in\mathbb{Z}\Rightarrow (x\star y)\star z=x\star(y\star z).$  Vrijedi
		\begin{align*}
		(x\star y)\star z&=(x+y-4)\star z=(x+y-4)+z-4=x+y+z-8,\\
		x\star (y\star z)&=x\star(y+z-4)=x+(y+z-4)-4=x+y+z-8.
		\end{align*}
		Dakle vrijedi asocijativni zakon.
		\item\label{grupa3}  Odredimo neutralni element $e$ za operaciju $\star.$ Vrijedi
		\begin{align*}
		x\star e=x\Leftrightarrow x+e-4=x\Leftrightarrow e=4,\\
		e\star x=x\Leftrightarrow e+x-4=x\Leftrightarrow e=4.
		\end{align*}
		Neutralni element je $e=4.$
		\item\label{grupa4} Odredimo sada inverzni $x'$  (suprotni, simetri\v cni) element, vrijedi 	
		\begin{align*}
		x\star x'&=e\Leftrightarrow x+x'-4=4\Leftrightarrow x'=-x+8,\\
		x'\star x&=e\Leftrightarrow x'+x-4=4\Leftrightarrow x'=-x+8,
		\end{align*}
		pa je inverzni element $x'=-x+8.$
		Sada na osnovu \ref{grupa1}.--\ref{grupa4}. zaklju\v cujemo da $(\mathbb{Z},\star)$ ima strukturu grupe.
		\item Poka\v zimo da je $\star$ komutativna operacija. Vrijedi
		\begin{align*}
		x\star y=x+y-4=y+x-4=y\star x.
		\end{align*}
		Dakle $(\mathbb{Z},\star)$ ima strukturu komutativne ili Abelove grupe.
	\end{enumerate}
\end{example}

Za dvije binarne operacije definisane u istom skupu mo\v zemo definisati osobinu\\ distributivnosti.
\index{algebarske strukture!distributivnost}
\begin{definition}[Distributivnost]
   Neka su na jednom skupu $S$ definisane dvije binarne operacije $\circ_1$ i $\circ_2.$ Ako va\v zi
   \[(\forall a,b,c\in S)\quad a\circ_2(b\circ_1 c)=(a\circ_2 b)\circ_1(a\circ_2 c)\text{ i }
        (a \circ_1 b)\circ_2 c=(a\circ_2 c)\circ_1(b\circ_2 c),\]
        ka\v zemo da je $\circ_2$ distributivna u odnosu na $\circ_1.$
\end{definition}
Ako dvije zatvorene binarne operacije u nekom skupu zadovoljavaju odre\dj ene osobine, dobijamo nove algebarske strukture sa dvije binarne operacije: prsten, tijelo i polje.
\index{algebarske strukture!prsten}
\begin{definition}[Prsten]
  Neka su na skupu $S$ definisane dvije binarne operacije redom ozna\v cene sa $+$ i $\cdot$. Ka\v zemo da $S$ \v cini prsten u odnosu na te dvije operacije ako su ispunjeni uslovi:
   \begin{enumerate}[$1^0$]
    \item Skup $S$    \v cini komutativnu grupu u odnosu na operaciju $+$;
    \item Operacija $\cdot$ je asocijativna, tj va\v zi
        \[(\forall a,b,c\in S)\:a\cdot(b\cdot c)=(a\cdot b)\cdot c;\]
     \item Operacija $\cdot$ je distributivna prema operaciji $+$, tj.
        \[(\forall a,b,c\in S)\: (a+b)\cdot c=a\cdot c+b\cdot c \text { i }
                  a\cdot(b+ c)=a\cdot b+a\cdot c.  \]
    \end{enumerate}
Prsten koji \v cine skup $S$   i njegove operacije $+$ i $\cdot$ ozna\v cava se $(S,+,\cdot).$
\end{definition}

\index{algebarske strukture!tijelo}
\begin{definition}[Tijelo]
  Prsten $(S,+,\cdot)$ zove se tijelo, ako skup $S\setminus\{ 0\}$ \v cini grupu u odnosu na operaciju $\cdot.$
\end{definition}

\index{algebarske strukture!polje}
\begin{definition}[Polje]
  Prsten $(S,+,\cdot)$ zove se polje, ako je skup $S\setminus\{0\}$ \v cini komutativnu grupu u odnosu na operaciju $\cdot.$
\end{definition}

Ako posmatramo operacije sabiranja i mno\v zenja u skupu cijelih brojeva nalazimo da je $(\mathbb{Z},+, \cdot)$ prsten, a u skupu realnih brojeva vidimo da $(\mathbb{R},+, \cdot)$ ima algebarsku strukturu polja.
\begin{example}
 Na skupu $\mathbb{Q}$ definisane su dvije binarne operacije $\oplus$ i $\odot$ sa $x\oplus y=x+y+1$ i $x\odot y=xy+x+y.$ Pokazati da $(\mathbb{Q},\oplus,\odot)$  ima strukturu polja.\\
 \noindent Rje\v senje. Skup $\mathbb{Q}$ je skup racionalnih brojeva, $\mathbb{Q}=\{\frac{m}{n}:m\in \mathbb{Z} \wedge n\in\mathbb{N }\}.$
 \begin{enumerate}[($\oplus$1)]
 	\item\label{grupa5} Vrijedi $\forall x,y\in\mathbb{Q}\Rightarrow x+y+1\in\mathbb{Q},$ pa je $x\oplus y\in\mathbb{Q},$  $(\mathbb{Q},\oplus)$  ima strukturu grupoida.
 	\item\label{grupa6} Poka\v zimo da vrijedi asocijativni zakon za binarnu operaciju $\oplus,$
 	\begin{align*}
 	(x\oplus y)\oplus z&=(x+y+1)\oplus z=(x+y+1)+z+1=x+y+z+2,\\
 	x\oplus(y\oplus z)&=x\oplus (y+z+1)=x+(y+z+1)+1=x+y+z+2,
 	\end{align*}
 	pa je $\forall x,y,z\in\mathbb{Q}\:\: (x\oplus y)\oplus z=x\oplus(y\oplus z).$
 	\item\label{grupa7} Neutralni element za binarnu operaciju $\oplus,$ na osnovu
 	    \begin{align*}
 	    x\oplus e=x\Leftrightarrow x+e+1=x\Leftrightarrow e=-1,\\
 	    e\oplus x=x\Leftrightarrow e+x+1=x\Leftrightarrow e=-1,
 	    \end{align*}
 	    je $e=-1.$
 	\item\label{grupa8} Simetri\v cni element za binarnu operaciju $\oplus,$     na osnovu
 	    \begin{align*}
 	    x\oplus x'=e\Leftrightarrow x+x'+1=-1\Leftrightarrow x'=-x-2,\\
 	    x'\oplus x=e\Leftrightarrow x'+x+1=-1\Leftrightarrow x'=-x-2,
 	    \end{align*}
 	    je $x'=-x-2.$
 	\item\label{grupa9}  Operacija $\oplus,$ na osnovu
 	  \begin{align*}
 	  x\oplus y=x+y+1=y+x+1=y\oplus x
 	  \end{align*}
 	  je komutativna.
 \end{enumerate}
Dakle, zbog ($\oplus$\ref{grupa5})--($\oplus$\ref{grupa9}), $(\mathbb{Q},\oplus)$ ima strukturu komutativne grupe.
\begin{enumerate}[$(\oplus \odot)$]
	\item\label{grupa10} Poka\v zimo da je $\odot$ distributivna prema $\oplus.$ Vrijedi
	  \begin{align*}
	  x\odot (y\oplus z)&=x\odot(y+z+1)=xy+xz+x+x+y+z+1\\
	                    &=xy+xz+2x+y+z+1=(xy+x+y)+(xz+x+z)+1\\
	                    &=(x\odot y)+(x\odot z)+1\\
	                    &=(x\odot y)\oplus (x\odot z),\\
      (x\oplus y)\odot z&=(x+y+1)\odot z=xz+yz+z+x+y+1+z\\
                        &=xz+yz+x+y+2z+1=(xz+x+z)+(yz+y+z)+1\\
                        &=(x\odot z)+(y\odot z)+1=    (x\odot z)\oplus(y\odot z).
   	  \end{align*}
\end{enumerate}
Ispitajmo sada kakvu ima strukturu $(\mathbb{Q}\setminus \{-1\},\odot),$ odnosno $(\mathbb{Q},\oplus, \odot).$
\begin{enumerate}[($\odot$1)]
	\item\label{grupa11} Po\v sto vrijedi $\forall x,y\in\mathbb{Q}\setminus \{-1\}\Rightarrow xy+x+y\in\mathbb{Q}\setminus\{-1\},$ onda je $x\odot y\in\mathbb{Q}\setminus\{-1\},$ tj. $(\mathbb{Q}\setminus \{-1\},\odot)$ sada ima strukturu grupoida.
	\item\label{grupa12} Poka\v zimo da je binarna operacija $\odot$ asocijativna. Vrijedi
	    \begin{align*}
	    (x\odot y)\odot z&= (xy+x+y)\odot z=xyz+xz+yz+xy+x+y+z\\
	                     &=x(yz+y+z)+(yz+y+z)+x  \\
	                     &=x\odot(yz+y+z)=x\odot(y\odot z).
	    \end{align*}
	    Do sada je pokazano da $(\mathbb{Q},\oplus, \odot)$ ima strukturu prstena.
	 \item\label{grupa13} Odredimo neutralni element binarne operacije $\odot,$
	     \begin{align*}
	     x\odot e=x\Leftrightarrow xe+x+e=x\Leftrightarrow e(x+1)=0\Leftrightarrow e=0,\\
	     e\odot x=x\Leftrightarrow ex+e+x=x\Leftrightarrow e(x+1)=0\Leftrightarrow e=0.
	     \end{align*}
	     Pa je neutralni element $e=0.$
	 \item\label{grupa14} Simetri\v cni (inverzni) element, na osnovu
	     \begin{align*}
	      x\odot x'=e\Leftrightarrow xx'+x+x'=0\Leftrightarrow x'=-\frac{x}{x+1},\\
	      x'\odot x=e\Leftrightarrow x'x+x'+x=0\Leftrightarrow x'=-\frac{x}{x+1},
	     \end{align*}
	     je $x'=-\frac{x}{x+1}.$
	     Na osnovu prethodnog, $(\mathbb{Q},\oplus, \odot)$ ima strukturu tijela.
	  \item\label{grupa15}  Na kraju poka\v zimo da je $\odot$ komutativna operacija. Vrijedi
	  \begin{align*}
	  x\odot y&=xy+x+y=yx+y+x=y\odot x,
	  \end{align*}
	  pa na osnovu ($\odot$\ref{grupa11})--($\odot$\ref{grupa15}),  $(\mathbb{Q}\setminus\{-1\})$ ima strukturu komutativna grupe, dok na osnovu  ($\oplus$\ref{grupa5})--($\odot$\ref{grupa15}),$ \\	   (\mathbb{Q},\oplus, \odot)$ ima strukturu polja.
\end{enumerate}

\end{example}

\section[Z\lowercase{adaci}]{Zadaci}\index{Zadaci za vje\v zbu!operacije}

 \begin{enumerate}
 	\item U skupu $\mathbb{Z}$ definisane su operacije $\star$ i $\circ$ sa: $a\star b=-2a+b-4$ i $a\circ b=a-2b+6.$ Rije\v siti jedna\v cinu  $(3x)\circ [(-1\circ 2)\star(4\circ -3)]=[(-4\star 1)\circ (2\star -3)].$
    \item Dat je skup $\mathbb{Z}$ i u njemu operacija definisana sa $(\forall a,b\in\mathbb{Z})\:a\star b=a+b+1.$ Ispitati da li vrijede
          komutativni i asocijativni zakon za operaciju $\star.$
    \item   Dat je skup $\mathbb{Z}$ i u njemu operacija definisana sa $(\forall a,b\in\mathbb{Z})\:a\star b=2a+b+1.$ Ispitati da li vrijede
          komutativni i asocijativni zakon za operaciju $\star.$
    \item  Dat je skup $\mathbb{Z}$ i u njemu operacije definisana sa
          \begin{enumerate}
             \item $(\forall a,b\in\mathbb{Z})\:a\star b=ab+a+b;$
             \item $(\forall a,b\in\mathbb{Z})\:a\circ b=a+ab+b.$
          \end{enumerate} Ispitati da li vrijede komutativni i asocijativni zakon za operacije $\star$ i $\circ.$
    \item Neka je na skupu $S=\{1,2,3,4,6,9\}$ definisana operacija $x\star y=NZD(x,y),$ ($NZD(x,y)$  je najve\' ci zajedni\v cki djelilac). Sastaviti tabelu operacije.  Da li je operacija $\star$ komutativna operacija?
    \item U skupu $\mathbb{R}$ definisane su operacije\\
       \begin{inparaenum}
          \item $x\star y=\frac{x}{y},\:$
          \item $x\star y=x(x+y),\:$
          \item $x\star y=\frac{xy}{x+y},\:$
          \item $x\star y=x-y,\:$
          \item $x\star y=x^2+xy+y^2,\:$
          \item $x\star y=x^2+xy+2y^2.$
       \end{inparaenum}          \\
       Koja je od ovih operacija komutativna?
    \item U skupu $A=\{1,2,3,4,5,6\}$ definisane su operacije $x\star y=\max\{x,y\}$ i  $x\circ y=\min\{x,y\}.$ \\
       \begin{inparaenum}
          \item Sastaviti tabele za operacije $\star$ i $\circ.$
          \item Pokazati da su operacije $\star$ i $\circ$ komutativne.
      \end{inparaenum}
    \item   U skupu $\mathbb{Z}$ definisane su operacije \\
        \begin{inparaenum}
            \item $(\forall x,y\in\mathbb{Z})\:x\star y=x+y+6;\:$
            \item $(\forall x,y\in\mathbb{Z})\:x\circ y=x+y+1;\:$
            \item $(\forall x,y\in\mathbb{Z})\:x\bullet y=x+y-2.$
         \end{inparaenum}
       Ispitati da li su algebarske strukture $(\mathbb{Z},\star),\,(\mathbb{Z},\circ),\,(\mathbb{Z},\bullet)$ Abelove grupe.
   \item Dat je skup $A=\{0,1,2,3,4,5,6\}$ i njemu je definisana operacija $+_7,$ (sabiranje po modulu $7$). \\
         \begin{inparaenum}
            \item Sastaviti tabelu za operaciju $+_7.$
            \item Koju algebarsku strukturu predstavlja $(A,+_7).$
         \end{inparaenum}
   \item Na skupu $\mathbb{Q}$ definisane su operacije $\oplus$ i $\odot$ sa $x\oplus y=x+y-2$ i $x\odot y=xy+x+y.$  Ispitati \v sta je  struktura  $(\mathbb{Q},\oplus,\odot).$
 \end{enumerate}

\chapter{Preslikavanja}
\pagestyle{fancy}

U ovom poglavlju bi\' ce uvedeni pojmovi koji se koriste u \v citavoj matematici, a posebno u nekim oblastima kao \v sto je matemati\v cka analiza. Ove pojmove najlak\v se je matemati\v cki precizno uvesti koriste\' ci teoriju skupova, odnosno odgovaraju\' ce pojmove iz teorije skupova. Ovakav pristup omogu\' cava nam uvo\dj enje jedinstvene terminologije koja se onda koristi u razli\v citim oblastima kako matematike tako i drugih nauka, a kori\v stenje jedinstvene terminologije dovelo je do lak\v seg povezivanja razli\v citih nauka.

\section[O\lowercase{snovni pojmovi}]{Osnovni pojmovi}
Pojam preslikavanja mo\v zemo uvesti preko pojma binarne relacije, drugim rije\v cima, pre-\\slikavanje je jedna posebna vrsta binarne relacije, tj. binarna relacija sa nekim dodatnim osobinama koje su navedene u sljede\' coj definiciji.

    \begin{figure}[!h]\centering
     \begin{subfigure}[b]{.33\textwidth}\centering
        \includegraphics[scale=.45]{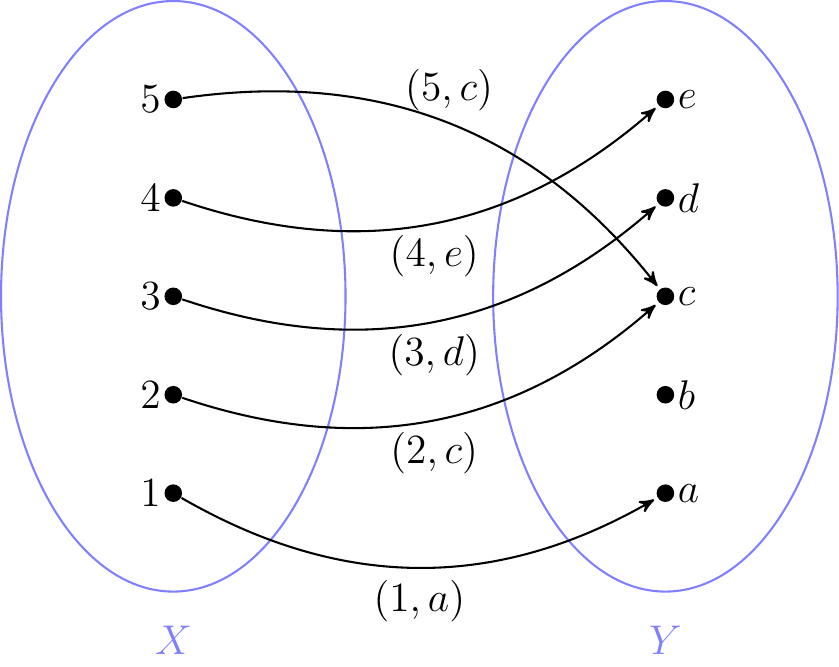}
         \caption{}
          \label{funkcija1a}
   \end{subfigure}
   \begin{subfigure}[b]{.33\textwidth}\centering
        \includegraphics[scale=.45]{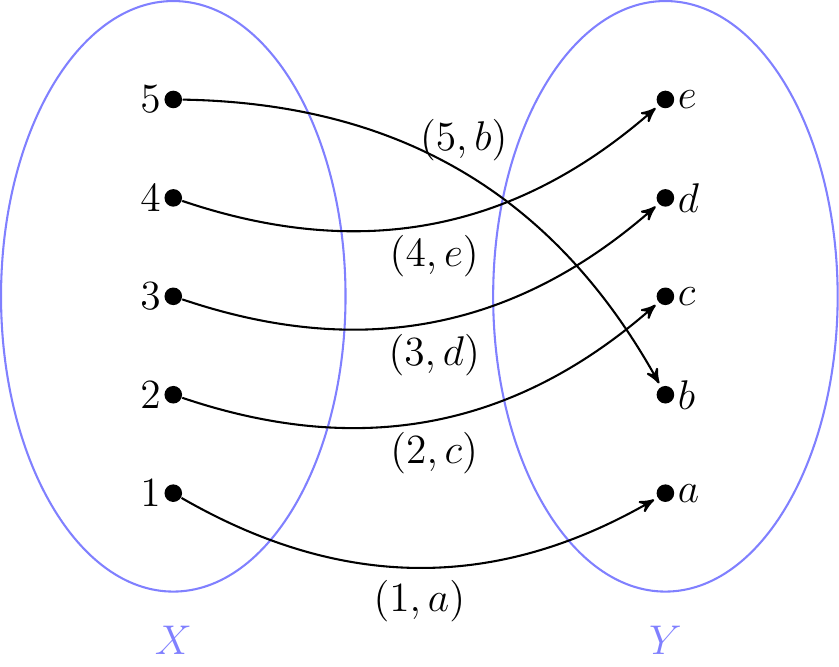}
         \caption{}
      \label{funkcija1b}
   \end{subfigure}
      \begin{subfigure}[b]{.33\textwidth}\centering
        \includegraphics[scale=.45]{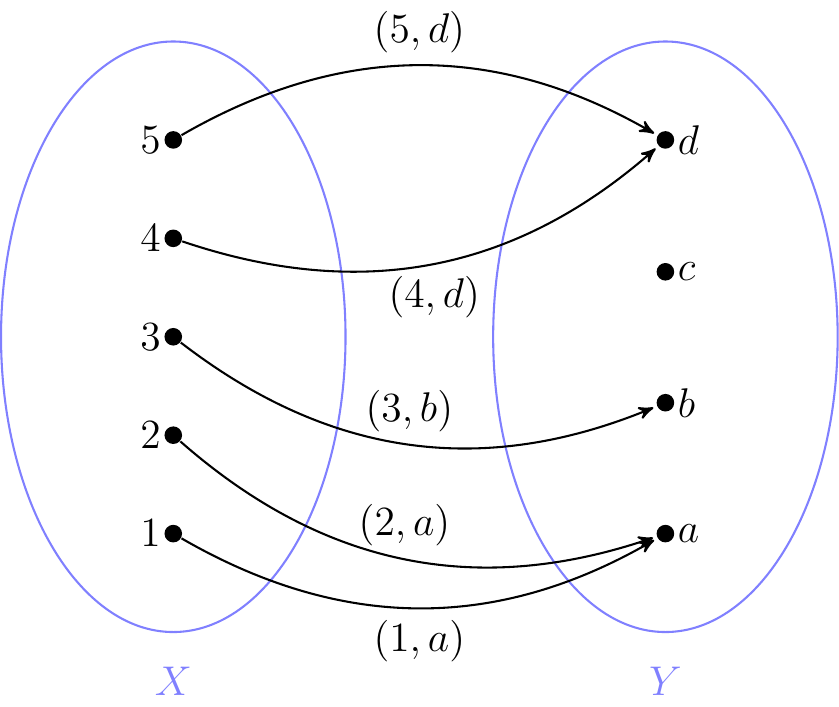}
         \caption{}
      \label{funkcija1c}
   \end{subfigure} \\
   \begin{subfigure}[b]{.33\textwidth}\centering
        \includegraphics[scale=.45]{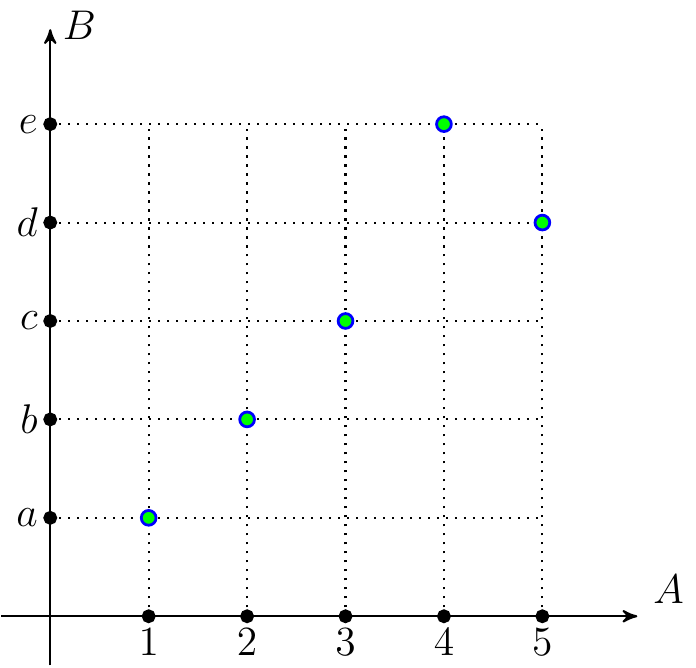}
         \caption{}
          \label{funkcija1a1}
   \end{subfigure}
   \begin{subfigure}[b]{.33\textwidth}\centering
        \includegraphics[scale=.45]{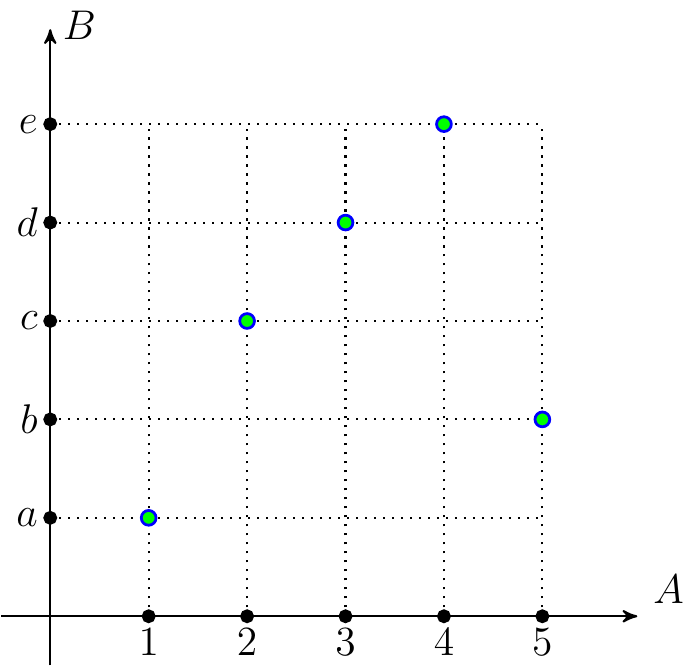}
         \caption{}
      \label{funkcija1b1}
   \end{subfigure}
      \begin{subfigure}[b]{.33\textwidth}\centering
        \includegraphics[scale=.45]{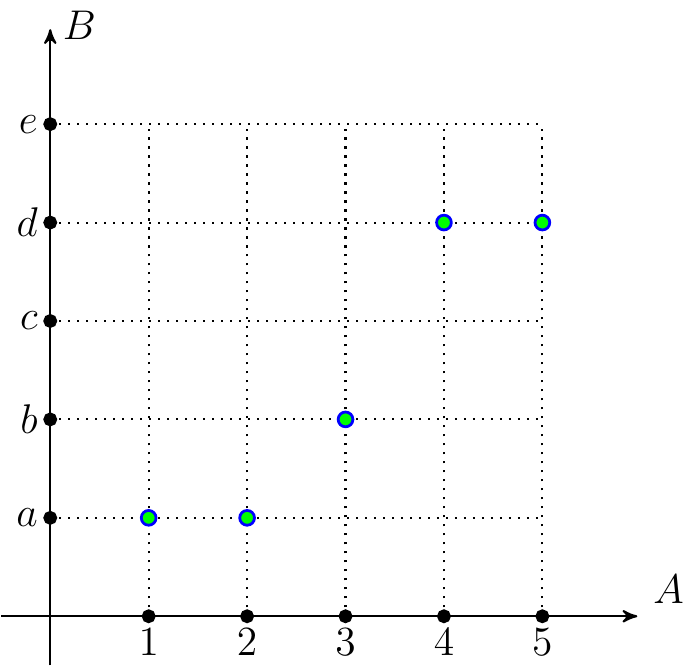}
         \caption{}
      \label{funkcija1c1}
   \end{subfigure}
  \caption{Binarne relacije koje su preslikavanja}
  \label{funkcija1}
  \end{figure}

\begin{figure}[!h]\centering
     \begin{subfigure}[b]{.33\textwidth}\centering
        \includegraphics[scale=.45]{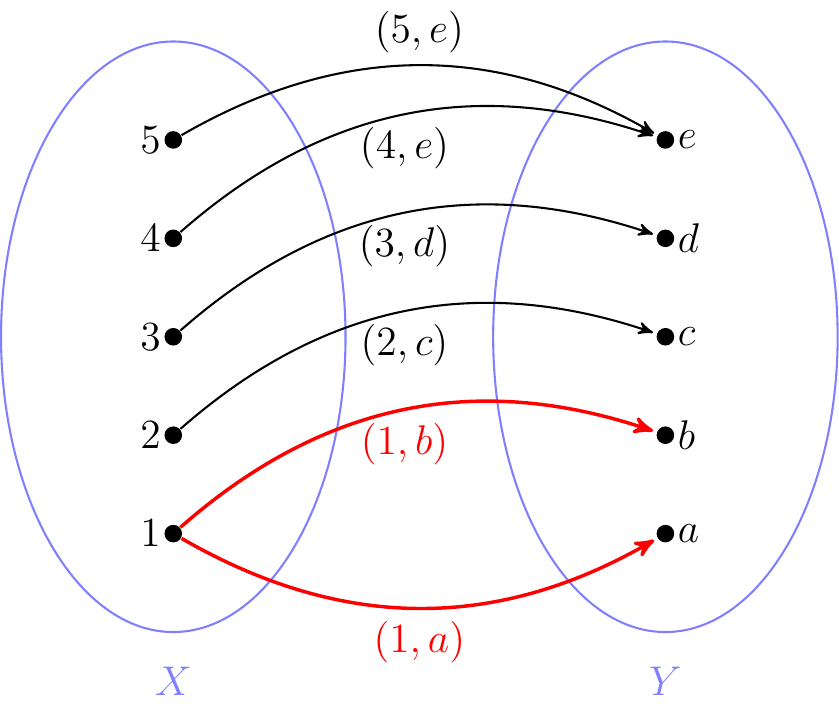}
         \caption{}
          \label{funkcija2a}
   \end{subfigure}
   \begin{subfigure}[b]{.33\textwidth}\centering
        \includegraphics[scale=.45]{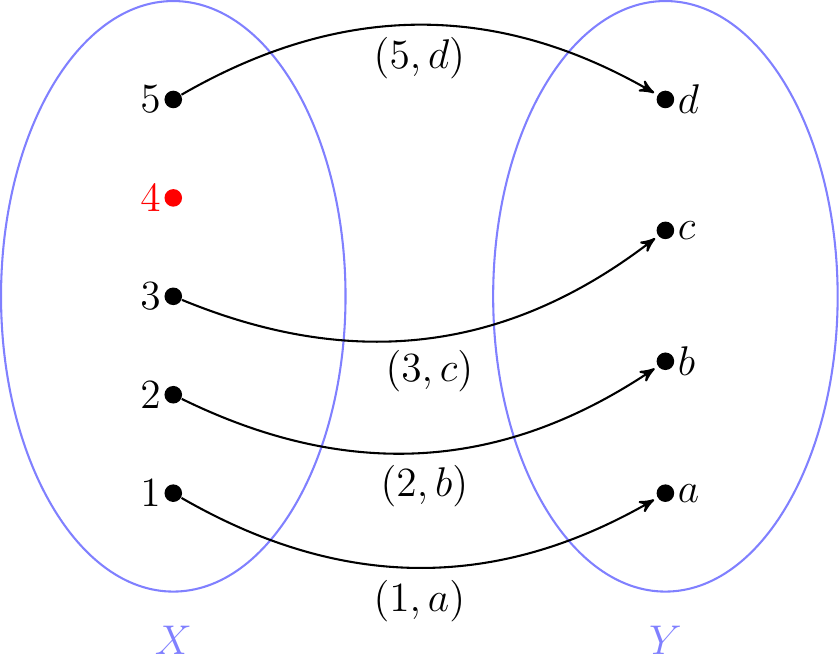}
         \caption{}
      \label{funkcija2b}
   \end{subfigure}
   \begin{subfigure}[b]{.33\textwidth}\centering
        \includegraphics[scale=.45]{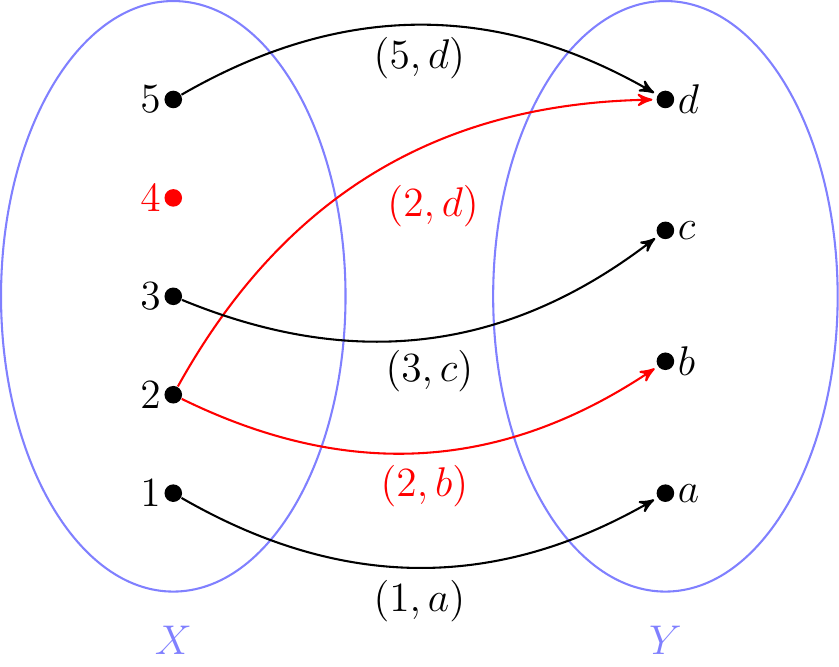}
         \caption{}
     \label{funkcija2c}
  \end{subfigure}
       \begin{subfigure}[b]{.33\textwidth}\centering
        \includegraphics[scale=.45]{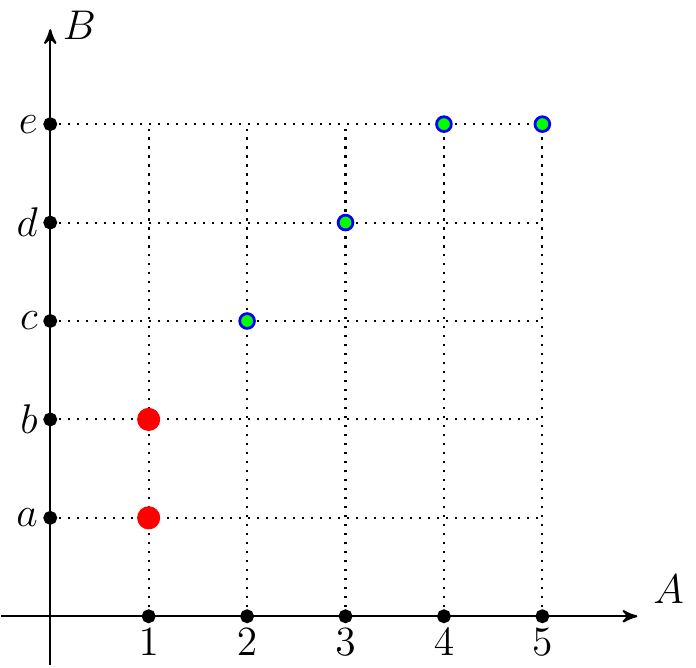}
         \caption{}
          \label{funkcija2a}
   \end{subfigure}
   \begin{subfigure}[b]{.33\textwidth}\centering
        \includegraphics[scale=.45]{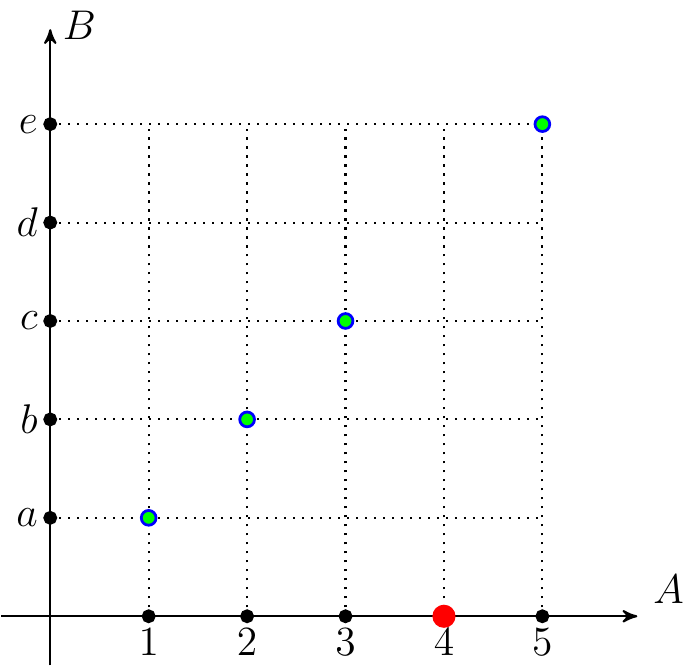}
         \caption{}
      \label{funkcija2b}
   \end{subfigure}
   \begin{subfigure}[b]{.33\textwidth}\centering
        \includegraphics[scale=.45]{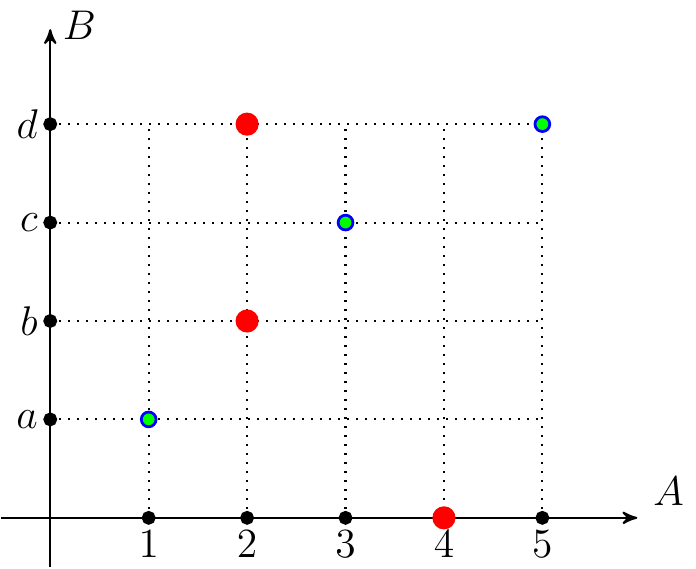}
         \caption{}
     \label{funkcija2c}
  \end{subfigure}
 \caption{Binarne relacije koje nisu preslikavanja}
  \label{funkcija2}
  \end{figure}

\begin{definition}[Preslikavanje--funkcija]
  Neka su $X$ i $Y$ neprazni skupovi. Binarnu relaciju $f\subset X\times Y$ nazivamo preslikavanje ili funkcija, sa skupa  $X$ u skup $Y$ ako vrijedi
  \begin{enumerate}
      \item $(\forall x\in X)(\exists y\in Y)(x,y)\in f;$
      \item $(x,y)\in f\wedge (x,y')\in f\Rightarrow y=y'.$
 \end{enumerate}
\end{definition}
Drugim rije\v cima, svakom elementu $x\in X$ takvom da je $(x,y)\in f$ odgovara jedan i samo jedan element $y\in Y.$ Na Slici \ref{funkcija1} predstavljene su relacije koje su preslikavanja. Dakle, svi elementi $x\in X$ moraju se preslikati u neki element iz skupa $Y$ i ne smiju se preslikati u vi\v se od jednog elementa. Na Slici \ref{funkcija2} dati su primjeri relacija koje nisu preslikavanja; crvenom bojom su ozna\v ceni ili elementi ili ure\dj eni parovi zbog kojih ove relacije nisu preslikavanja sa skupa $X$ u skup $Y.$
\begin{remark}
  \v Ce\v s\'ce se koristi termin funkcija kada su $X$ i $Y$ skupovi brojeva, dok se u ostalim slu\v cajevima koristi termin preslikavanje.
\end{remark}	

Skup $X$ zovemo domen (ili definiciono podru\v cje ili oblast definisanosti) preslikavanja $f$ i ozna\v cavamo $D(f)$ ili $D_f,$ a skup $f(X)=\{ f(x):x\in X\}$ zovemo kodomen (ili skup vrijednosti ili rang) preslikavanja $f$. Oznake koje se koriste za kodomen su $V(f),\,V_f,\,R(f),\,R_f,\,f(X)$ i druge. Elemente skupa $D_f$ zovemo originali (ili argumenti ili promjenljive), a elemente skupa $R_f$ slike. Relacija $f$ je preslikavanje ako  svakom ori-\\ginalu odgovara ta\v cno jedna slika.

Da preslikavanje $f$ preslikava skup $X$ u skup $Y$ ozna\v cavamo sa
\begin{empheq}[box=\mymath]{equation*}
f:X\mapsto Y,
\label{fun1}
\end{empheq}
ako je $(x,y)\in f$ onda pi\v semo

\begin{empheq}[box=\mymath]{equation*}
y=f(x) \text{ ili } f:x\mapsto y.
\label{fun2}
\end{empheq}

\begin{remark}
\v Cesto se u literaturi sre\'{c}emo sa zadavanjem preslikavanja iz skupa $X$ u skup $Y,$ pri \v cemu nije svaki element iz $X$ original. Tada
se nalazi definiciono podru\v cje $D(f)$ koje je strogi podskup od $X.$ Dakle, generalno vrijedi $D(f)\subseteq X.$
\end{remark}	
Preslikavanje mo\v zemo zadati na vi\v se na\v cina: analiti\v cki, tabelarno, grafi\v cki ili strelastim dijagramom.
Analiti\v cki oblici zadavanja preslikavanja su sljede\' ci:
\begin{enumerate}
\item  eksplicitni oblik $y=f(x);$
\item implicitni oblik $F(x,y)=0;$
\item parametarski oblik $x=x(t), \, y=y(t)$ ($t$ je parametar);
\item polarni oblik $\rho =\rho (\varphi)$  ($\rho$-modul, $\varphi$-argument).
\end{enumerate}
Primjeri strelastog dijagrama su na Slici \ref{funkcija1a}, \ref{funkcija1b} i \ref{funkcija1c}. Tabelarni prikaz preslikavanja mo\v zemo  vidjeti u Primjeru \ref{primjer51}
kao Tabela \ref{fun3}. Preslikavanja odnosno funkciju mo\v zemo prikazati i grafi\v cki u Dekartovom koordinatnom sistemu kao skup ta\v caka (ili linija). Preciznije je
dato u sljede\' coj definiciji.

\begin{figure}[!h]\centering
	\begin{subfigure}[b]{.45\textwidth}\centering
		\includegraphics[scale=.8]{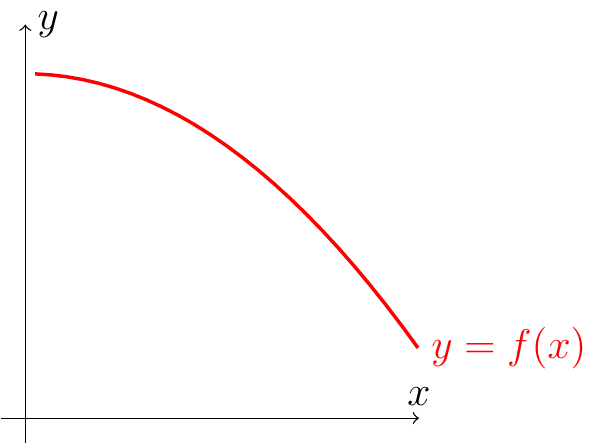}
		\caption{}
		\label{funkcija2a}
	\end{subfigure}
	\begin{subfigure}[b]{.45\textwidth}\centering
		\includegraphics[scale=.8]{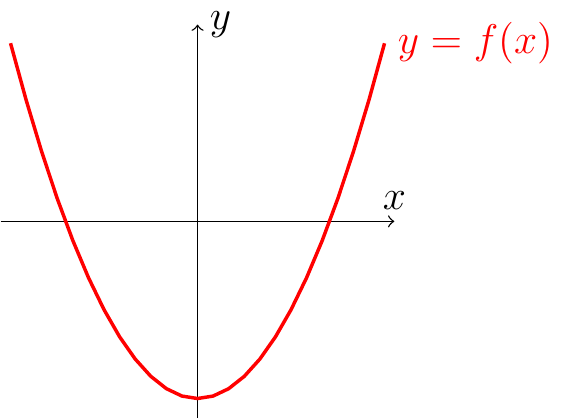}
		\caption{}
		\label{funkcija2b}
	\end{subfigure}
	\caption{Grafici funkcija}
	\label{grafik3}
\end{figure}

\begin{definition}[Grafik funkcije]
  Ako je $D_f\subset\mathbb{R}$ i $R_f\subset\mathbb{R}$ onda skup ta\v caka $\Gamma(f)=\{(x,f(x)):x\in D_f\} \subseteq\mathbb{R }\times \mathbb{R }$ zovemo grafik funkcije $f.$
\end{definition}
Primjeri dva grafika funkcija su na Slici \ref{grafik3}.

\section[V\lowercase{rste preslikavanja}]{Vrste preslikavanja}

Neka je $f:X\mapsto Y$ preslikavanje, dakle svakom elementu $x\in X$ preslikavanje $f$ pridru\v zuje neki element $y\in Y.$  Ne znamo ni\v sta o tome da li su svi elementi skupa $Y$ slika nekog elementa (originala) iz skupa $X,$ kao ni to da li je neki element $y\in Y$ slika jednog ili vi\v se originala (vidjeti Sliku \ref{funkcija1}). U svrhu daljeg analiziranja preslikavanja, razmatraju se neke dodatne osobine.

Slijede definicije preslikavanja koje imaju neke od tih dodatnih osobina.

\begin{figure}[!h]\centering
   \begin{subfigure}[b]{.33\textwidth}\centering
        \includegraphics[scale=.45]{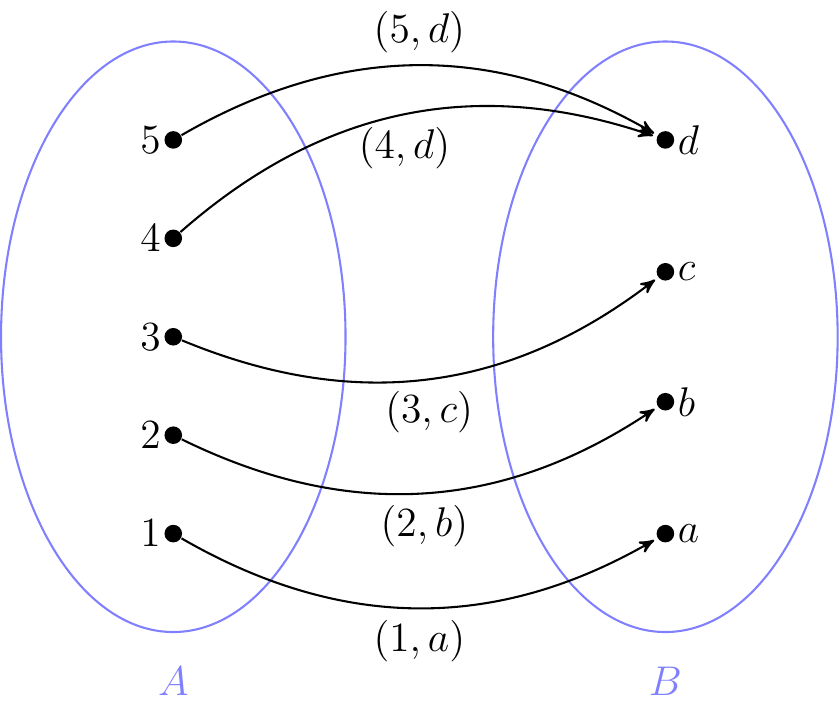}
         \caption{Sirjekcija}
      \label{funkcija3b}
   \end{subfigure}
  \begin{subfigure}[b]{.33\textwidth}\centering
        \includegraphics[scale=.45]{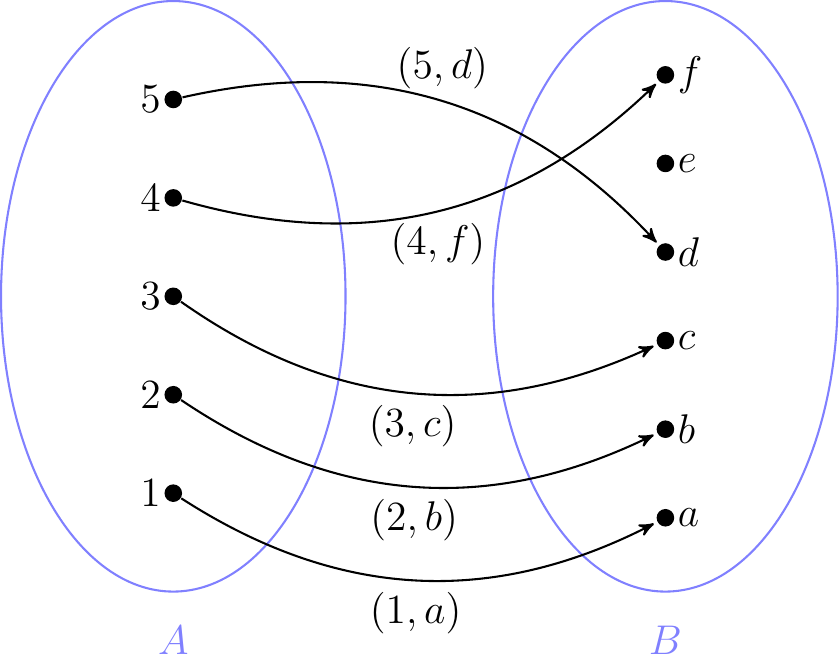}
         \caption{Injekcija}
          \label{funkcija3a}
   \end{subfigure}
   \begin{subfigure}[b]{.33\textwidth}\centering
        \includegraphics[scale=.45]{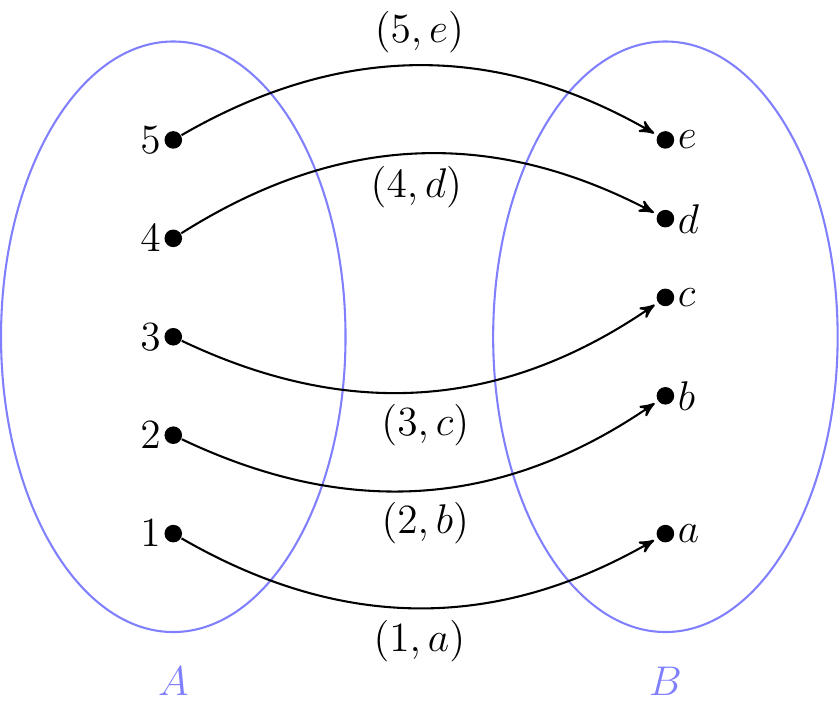}
         \caption{Bijekcija}
     \label{funkcija3c}
   \end{subfigure}

   \begin{subfigure}[b]{.33\textwidth}\centering
        \includegraphics[scale=.45]{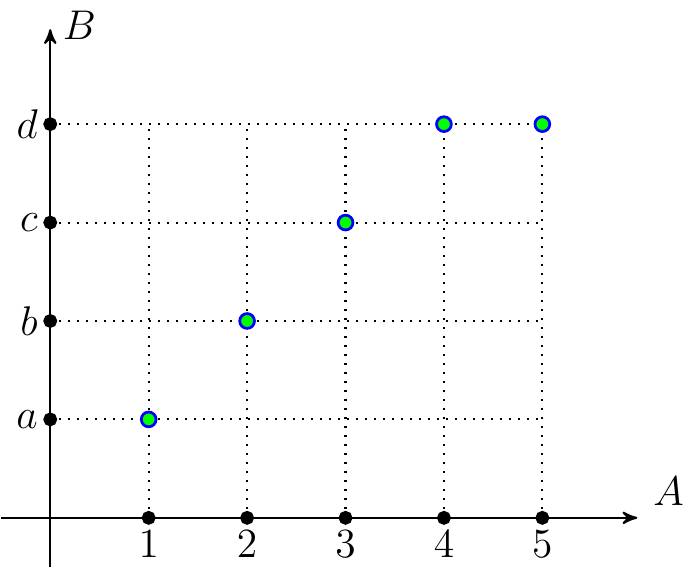}
         \caption{Sirjekcija}
      \label{funkcija3bk}
   \end{subfigure}
   \begin{subfigure}[b]{.33\textwidth}\centering
        \includegraphics[scale=.45]{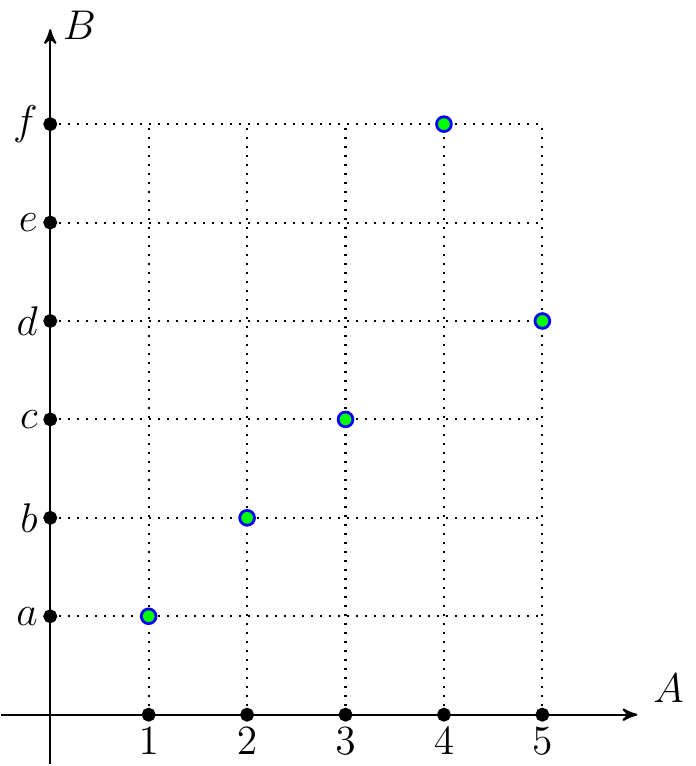}
         \caption{Injekcija}
          \label{funkcija3ak}
   \end{subfigure}
   \begin{subfigure}[b]{.33\textwidth}\centering
        \includegraphics[scale=.45]{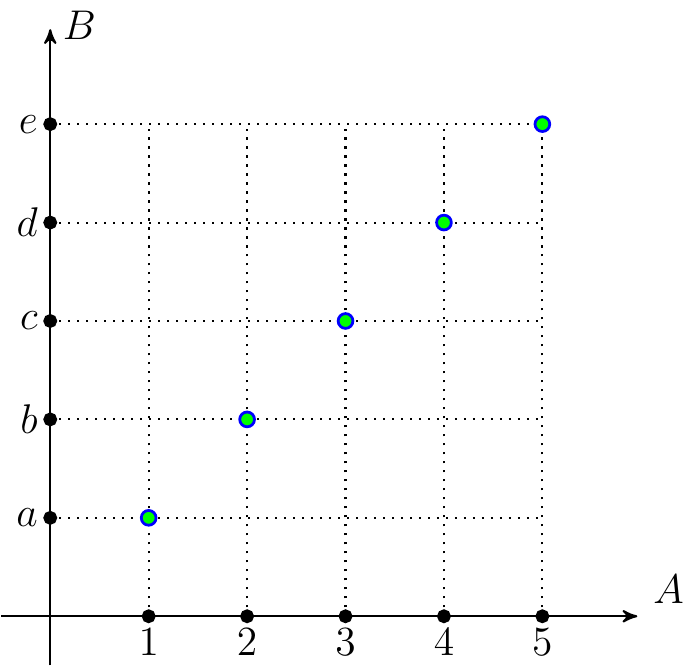}
         \caption{Bijekcija}
     \label{funkcija3ck}
  \end{subfigure}
\caption{Vrste preslikavanja: Sirjekcija, injekcija i bijekcija}
  \label{funkcija3}
  \end{figure}

\begin{definition}[Sirjekcija]
      Za preslikavanje $f:X\mapsto Y$ ka\v ze se da je sirjekcija ili sirjektivno preslikavanje sa skupa $X$ na skup $Y$ ako je $R_f=Y,$ tj. ako je svaki $y\in Y$
      slika bar jednog elementa $x\in X.$
\end{definition}
Sirjekcija se naziva i "preslikavanje na".
Na Slici \ref{funkcija3b} i \ref{funkcija3bk} je prikazan primjer sirje-\\ktivnog preslikavanja.

\begin{definition}[Injekcija]
  Za preslikavanje $f:X\mapsto Y$ ka\v ze se da je injekcija ili injektivno preslikavanje sa skupa $X$ u skup $Y$, ako vrijedi
 \[(\forall x, x'\in X) \, x\neq x'\Rightarrow f(x) \neq f(x'), (f(x), f(x') \in Y).   \]
\end{definition}
Drugim rije\v cima, kod injektivnog preslikavanja razli\v citim originalima odgovaraju razli\v cite slike, vidjeti Sliku \ref{funkcija3a} i \ref{funkcija3ak}.
Injektivno preslikavanje se naziva jo\v s i "$1-1$ preslikavanje".
Na osnovu zakona kontrapozicije vrijedi da je $f:X\mapsto Y$ injekcija ako va\v zi
 \[(\forall x, x') \, f(x)=f(x')\Rightarrow x=x'.\]
Ako je neko preslikavanje istovremeno i sirjektivno i injektivno dobijamo novo preslikavanje dato u sljede\' coj definiciji.
\begin{definition}[Bijekcija]
  Preslikavanje $f:X\mapsto Y$ naziva se bijekcija ili bijektivno preslikavanje, ako je ona istovremeno i sirjekcija i injekcija.
\end{definition}
Bijekcija se naziva jo\v s i "uzajamno jednozna\v cno preslikavanje" ili "obostrano jednozna\v cno preslikavanje", vidjeti Sliku \ref{funkcija3c} ili \ref{funkcija3ck}.

\begin{definition}[Kompozicija]
  Ako su $f:X\mapsto Y$ i $g:Y\mapsto Z$ preslikavanja, tada preslikavanje $h:X\mapsto Z$ zadano sa \[(\forall x \in X) \, h(x)=g(f(x))\]
zovemo slo\v zeno preslikavanje  (ili kompozicija ili slaganje preslikavanja), preslikavanja $f$ i $g$ i pi\v semo $h=g\circ f.$
\end{definition}
Dakle vrijedi $h(x)=(g \circ f)(x)=g(f(x)).$

Asocijativni zakon vrijedi za kompoziciju preslikavanja. Za preslikavanja $f:X\mapsto Y,$ $g:Y\mapsto Z,$ i $h:Z\mapsto W$ vrijedi \[(h\circ g)\circ f=h\circ (g\circ f).\]

Svakom  preslikavanju $f:X\mapsto Y,$ kao specijalnom slu\v caju relacije iz $X$ u $Y$ odgovara inverzna  relacija $f^{-1}\subset Y\times X$ iz $Y$ u $X.$ U op\v stem slu\v caju relacija $f^{-1}$ nije preslikavanje. U slu\v caju kada je $f^{-1}$ preslikavanje, onda se ono naziva inverzno preslikavanje preslikavanja $f.$
\begin{definition}
Neka je $f:X\mapsto Y$ bijektivno preslikavanje. Inverzno preslikavanje preslikavanja $f$ je preslikavanje koje zadovoljava uslov
\[ (\forall x \in X) \, f^{-1}(f(x))=x.  \]
\end{definition}

Uslov kada je relacija $f^{-1}$ inverzno preslikavanje preslikavanja $f,$ dat je u sljede\' coj teoremi.

\begin{theorem}
   Inverzna relacija $f^{-1}$ preslikavanja $f:X\mapsto Y$ bi\' ce preslikavanje  $f^{-1}:Y\mapsto X$ ako i samo ako  je $f$ bijekcija.
\end{theorem}

 \begin{figure}[!h]\centering
   \begin{subfigure}[b]{.45\textwidth}\centering
        \includegraphics[scale=.6]{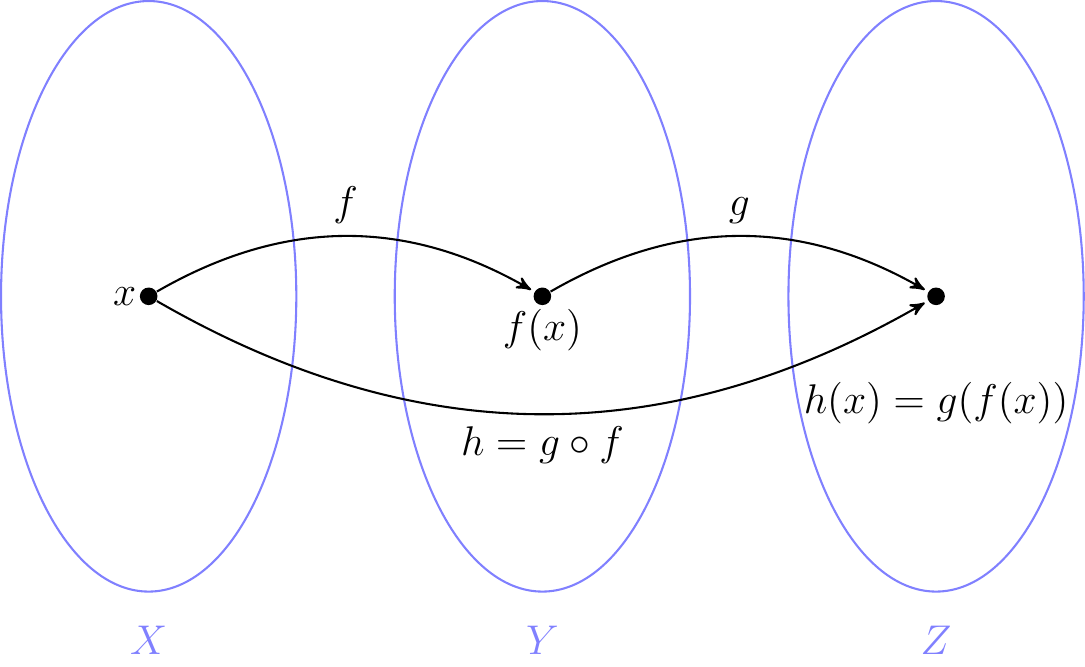}
         \caption{Kompozicija funkcija}
      \label{kompozicija}
   \end{subfigure} \hspace{.5cm}
  \begin{subfigure}[b]{.45\textwidth}\centering
        \includegraphics[scale=.6]{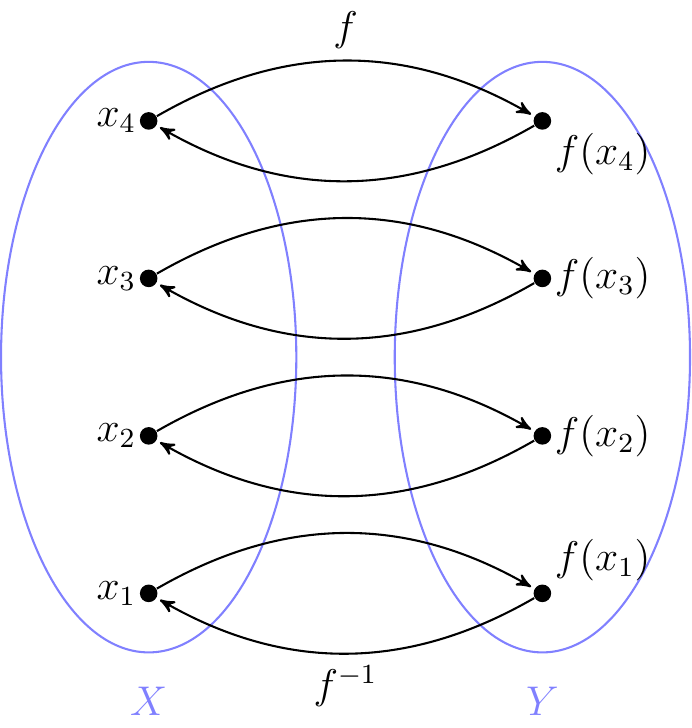}
         \caption{Inverzna funkcija $f^{-1}$ funkcije $f$}
          \label{Inverzna}
   \end{subfigure}
   \caption{Kompozicija preslikavanja (funkcija) i inverzno preslikavanje (funkcija)}
  \label{komp}
  \end{figure}

\begin{figure}[!h]\centering
   \begin{subfigure}[b]{.45\textwidth}\centering
        \includegraphics[scale=.6]{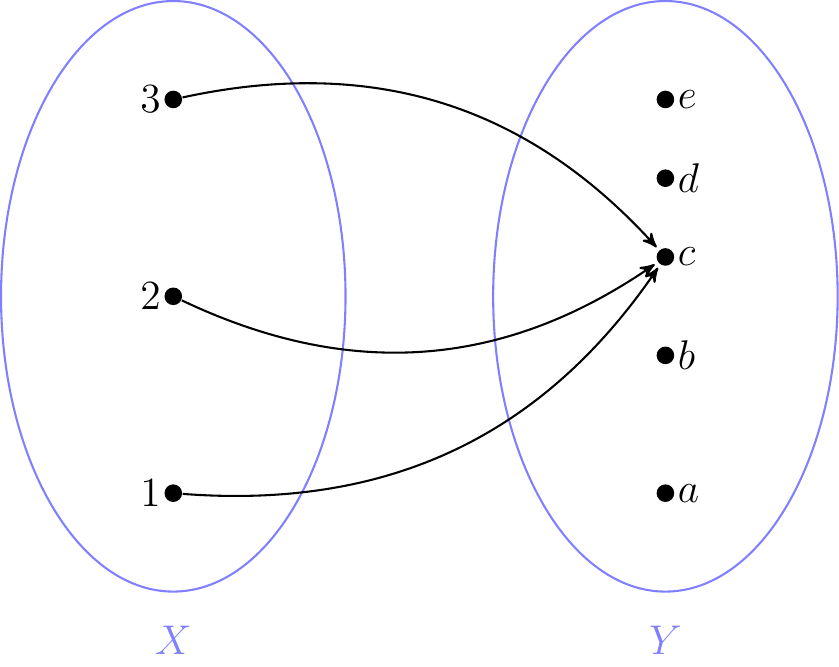}
         \caption{Svi elementi skupa $X$ se slikaju u element $c$ }
      \label{konstantno1}
   \end{subfigure} \hspace{.5cm}
  \begin{subfigure}[b]{.4\textwidth}\centering
        \includegraphics[scale=1]{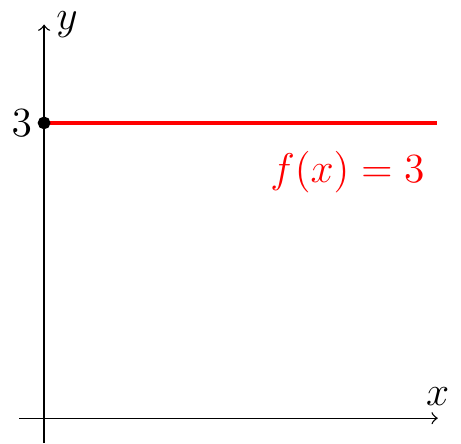}
         \caption{Grafik konstantne funkcije $f(x)=3,$ $x\geq 0$}
          \label{konstantno2}
   \end{subfigure}
   \caption{Primjeri konstantnih funkcija}
  \label{konstantno3}
  \end{figure}

\begin{figure}[!h]\centering
   \begin{subfigure}[b]{.45\textwidth}\centering
        \includegraphics[scale=.6]{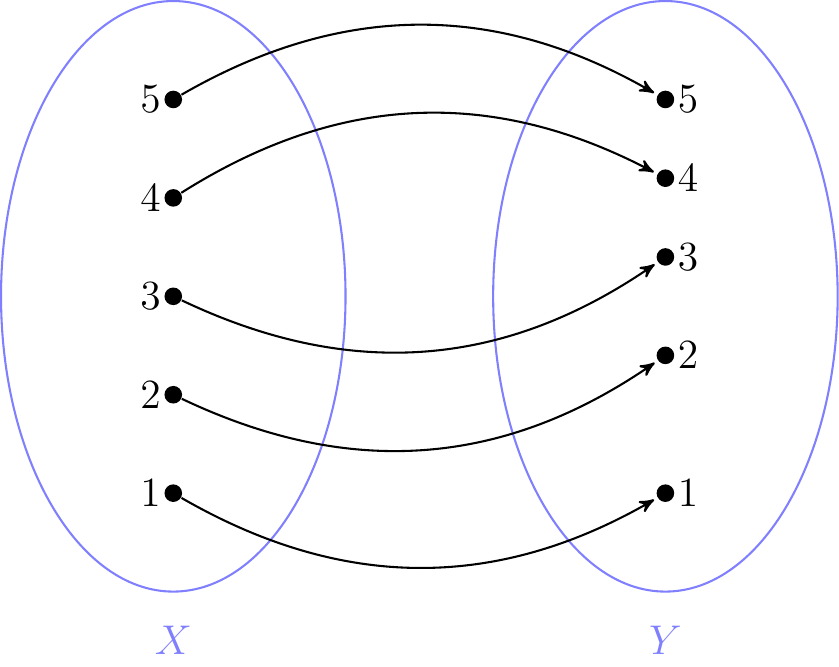}
         \caption{Svi elementi skupa $X$ se slikaju u iste elemente}
      \label{identicko1}
   \end{subfigure} \hspace{.5cm}
  \begin{subfigure}[b]{.4\textwidth}\centering
        \includegraphics[scale=1]{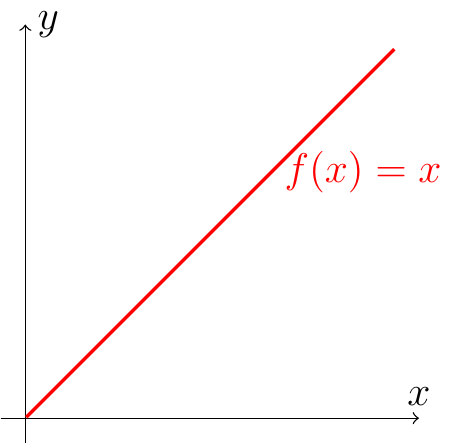}
         \caption{Grafik identi\v cke funkcije $f(x)=x,$ $x\geq 0$}
          \label{identicko2}
   \end{subfigure}
   \caption{Primjeri identi\v ckih funkcija}
  \label{identicko3}
  \end{figure}

Osim navedenih vrsta preslikavanja/funkcija, \v cesto susre\' cemo i konstantnu funkciju Slika \ref{konstantno3} i identi\v cku funkciju Slika \ref{identicko3}.

\vspace{2cm}

\newpage

\section[Z\lowercase{adaci}]{Zadaci}

\begin{figure}[!h]\centering
   \begin{subfigure}[b]{.45\textwidth}\centering
        \includegraphics[scale=.55]{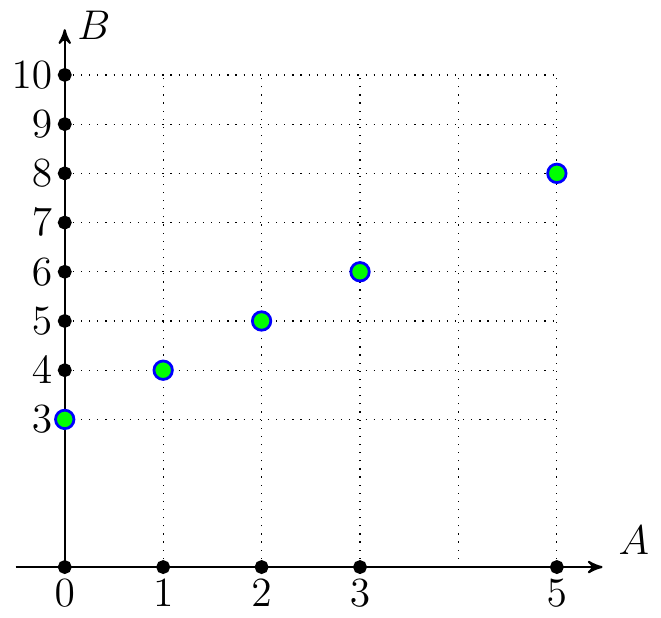}
         \caption{Funkcija $f(x)=x+3$ predstavljena u koordinatnom sistemu}
      \label{funkcija9a}
   \end{subfigure} \hspace{.5cm}
  \begin{subfigure}[b]{.45\textwidth}\centering
        \includegraphics[scale=.55]{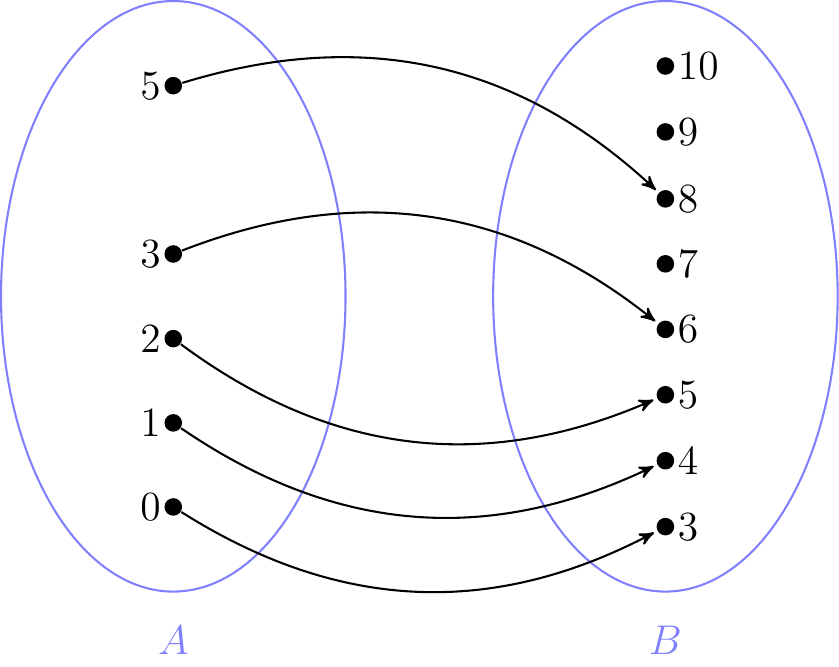}
         \caption{Funkcija $f(x)=x+3$ predstavljena u strelastom dijagramu}
          \label{funkcija9b}
   \end{subfigure}
   \caption{Funkcija $f(x)=x+3$}
  \label{funkcija9c}
  \end{figure}

\begin{example}\label{primjer51}
  Dati su skupovi $A=\{0,1,2,3,5\},\,B=\{3,4,5,6,8,9,10\}$ i funkcija $f(x)=x+3.$ Predstaviti funkciju pomo\' cu tabele, ure\dj enih parova, koordinatnog sistema i strelastog dijagrama.\\\\
  Rje\v senje:\\\\
Vidi Tabelu \ref{fun3} za tabelarno predstavljenja fukcije. Predstavljenje funkcije ure\dj enim parovima  $f=\{(0,3),(1,4),(2,5),(3,6),(5,8)\}.$ Vidjeti Sliku \ref{funkcija9c} za koordinatni sistem i strelasti dijagram.
\end{example}

\begin{table}[!h]\centering
   \begin{tabular}{c|ccccc}
            $x$ & $0$ & $1$ & $2$ & $3$ & $5$\\\hline
            $f(x)$ & $3$ & $4$ & $5$ & $6$ & $8$
  \end{tabular}
  \caption{Predstavljanje funkcije $f(x)=x+3$ tabelarno}
  \label{fun3}
\end{table}

\begin{figure}[!h]\centering
	\includegraphics[scale=.65]{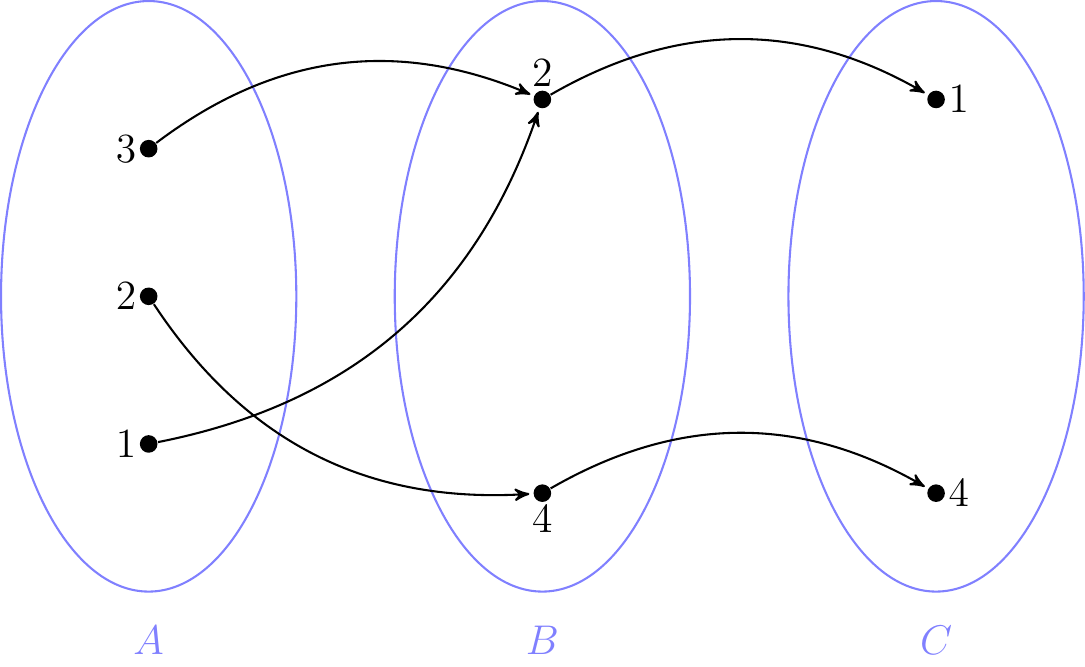}
	\caption{Kompozicija preslikavanja}
	\label{funkcija11a}
\end{figure}

\begin{example}
   Dati su skupovi $A=\{ 1,2,3\},\,B=\{2,4\}$ i $C=\{1,4\},$ funkcija $f:A\mapsto B$ zadana sa $f(1)=2,\,f(2)=4,\,f(3)=2,$  te funkcija $g:B\mapsto C$ zadana sa $g(2)=1,\,g(4)=4.$ Odrediti funkciju $h=g\circ f,$ tj. $h(x)=g(f(x)).$ \\ \ \

\noindent Rje\v senje: \\\\
Vrijedi $h(1)=1,\: h(3)=1,\:h(2)=4,$ vidjeti Sliku \ref{funkcija11a}.
\end{example}

\newpage

\begin{example}
   Date su funkcije $f(x)=x+1,\,g(x)=2x-1,\,k(x)=3x+1.$ Odrediti kompozicije funkcija $f\circ g,\,g\circ f,\, f\circ(g\circ k),\,(f\circ g)\circ k,\,f\circ(f\circ g).$\\\\
Rje\v senje:

\begin{align*}
 &\text{Izra\v cunajmo prvo $f\circ g$ i $g\circ f,$ vrijedi}\\
 &  (f\circ g)(x)=f(g(x))=g(x)+1=2x-1+1=2x,\\
 &  (g\circ f)(x)=g(f(x))=2(f(x))-1=2(x+1)-1=2x+1.\\
 &\text{Odredimo sada } g\circ k\\
 &  (g\circ k)(x))=2(3x+1)-1=2(3x+1)-1=6x+1.\\
 &\text{Sada izra\v cunajmo $(f\circ (g\circ k)), (f\circ g)\circ k$ i $(f\circ (f\circ g))$}\\
 & (f\circ (g\circ k)))=g(k(x))+1=6x+1+1=6x+2,\\
 &(f\circ g)\circ k=2k(x)=2 (3x+1)=6x+2,\\
 & (f\circ(f\circ g))(x)=f(g(x))+1=2x+1.
\end{align*}
Vidimo da u op\v stem slu\v caju kompozicija  preslikavanja (funkcija) nije komutativna operacija ($(f\circ g)(x)\neq (g\circ f)(x)$).
\end{example}

\begin{example}
   Date su funkcije $f(x)=2x,\,g(x)=x^2-5.$ Odrediti $f\circ g$ i $g\circ f.$\\\\
Rje\v senje:
\begin{align*}
  &  (f\circ g)(x)=f(g(x))=2(g(x))=2(x^2-5)=2x^2-10,\\
  &  (g\circ f)(x)=g(f(x))=(f(x))^2-5=(2x)^2-5=4x^2-5.
\end{align*}
\end{example}

\begin{example}
   Date su funkcije $f(x+1)=3x+3,\,g(2x-1)=x,\,k(x-3)=x^2+x.$ Odrediti
      \begin{enumerate}[$(a)$]
         \item $f(x),\,g(x),\,k(x);$
         \item $f\circ g,\, g\circ k,\,f\circ (g\circ k).$
      \end{enumerate}

\noindent Rje\v senje:
\begin{enumerate}[$(a)$]
    \item Odredimo prvo $f(x),$ uvodimo smjenu $x+1=t,$ pa je $x=t-1.$ Uvrstimo prethodne izraze u $f(x+1)=3x+3, $ vrijedi
         $f(t)=3(t-1)+3=3t,$  dakle vrijedi $f(t)=3t$ pa je i $f(x)=3x.$ Sljede\' ca smjena je $2x-1=t,$ pa je $x=\frac{t+1}{2},$
         pa je  $g(t)=\frac{t+1}{2}$  te vrijedi $g(x)= \frac{x+1}{2}.$  I na kraju smjena je $x-3=t$ i $x=t+3,$
         pa je $k(t)=(t+3)^2+t+3=t^2+7t+12$ i $k(x)=x^2+7x+12.$ Odredimo sada   $f\circ g,\, g\circ k,\,f\circ (g\circ k).$
    \item Vrijedi
\begin{align*}
  &  (f\circ g)(x)= 3 g(x)=3\frac{x+1}{2}=\frac{3(x+1)}{2},\\
  &  (g\circ k)(x)=\frac{k(x)+1}{2}=\frac{x^2+7x+12+1}{2}=\frac{x^2+7x+13}{2},\\
  &  (f\circ (g\circ k))=3g(k(x))=3\frac{x^2+7x+13}{2}=\frac{3(x^2+7x+13)}{2}.
\end{align*}
\end{enumerate}
\end{example}

\begin{figure}[!h]\centering
	\includegraphics[scale=.85]{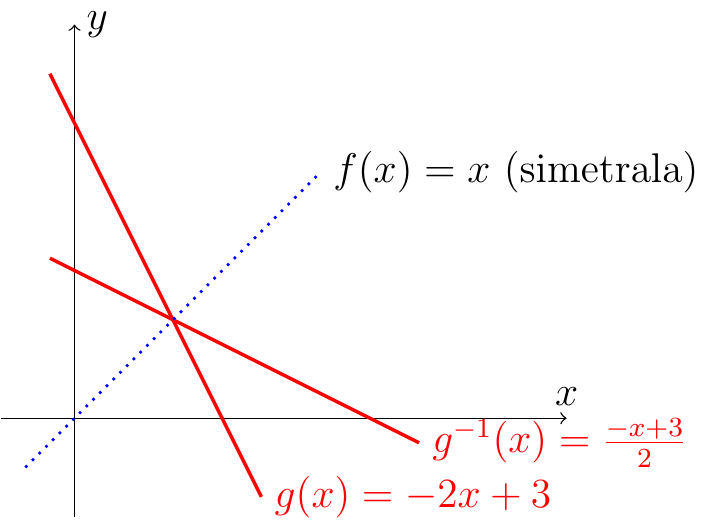}
	\caption{Kompozicija preslikavanja}
	\label{funkcija12a}
\end{figure}

\begin{example}
  Data je bijekcija $g:\mathbb{R}\mapsto\mathbb{R}$ sa $g(x)=3-2x.$
    \begin{enumerate}[$(a)$]
       \item Odrediti $g^{-1},\,g^{-1}\circ g,\, g\circ g^{-1};$
       \item Nacrtati u istom koordinatnom sistemu grafike funkcija $g$ i $g^{-1}.$
    \end{enumerate}
Rje\v senje:
\begin{enumerate}[$(a)$]
    \item Odredimo prvo $g^{-1}.$ Izrazimo $x$ iz $g(x)=3-2x,$ vrijedi $x=\frac{3-g(x)}{2},$ sada zamijenimo $x$ sa $g^{-1}(x)$
             a $g(x)$ sa $x,$ te dobijamo $g^{-1}(x)=\frac{3-x}{2}.$  Dalje je
        \begin{align*}
            & (g^{-1}\circ g)(x)=g^{-1}(g(x)))=\frac{3-g(x)}{2}=\frac{3-(3-2x)}{2}=x,\\
            & (g\circ g^{-1})(x)= g(g^{-1}(x))=3-2g^{-1}(x)=3-2\frac{3-x}{2}=3-(3-x)=x.
        \end{align*}
    \item Vidjeti Sliku   \ref{funkcija12a}.
 \end{enumerate}
\end{example}

\begin{example}
  Date su funkcije $f\left( \frac{x+3}{2}\right)=x,\,g(2x-1)=x+4.$ Odrediti
    \begin{enumerate}[$(a)$]
       \item $f(x),\,g(x);$
       \item $f^{-1}(x),\,g^{-1}(x);$
       \item $f^{-1}\circ g^{-1},\,(f^{-1}\circ g)\circ f.$
    \end{enumerate}
\noindent Rje\v senje:
\begin{enumerate}[$(a)$]
    \item Uvodimo smjenu $\frac{x+3}{2}=t$ i $x=2t-3,$ sada je  $f(t)=2t-3$ i $f(x)=2x-3.$ Sljede\' ca smjena je $2x-1=t$ i $x=\frac{t+1}{2}$
          i vrijedi $g(t)=\frac{t+1}{2}+4 =\frac{t+9}{2}$ i $g(x)=\frac{x+9}{2 }.$
    \item Iz $f(x)=2x-3$ dobijamo $x=\frac{f(x)+3}{2},$ pa je $f^{-1}(x)=\frac{x+3}{2}.$ Na isti na\v cin iz $g(x)=\frac{x+9}{2}$ je
          $x=2g(x)-9,$ i $g^{-1}(x)=2x-9.$
    \item  Vrijedi
      \begin{align*}
       &  (f^{-1}\circ g^{-1})(x)=f^{-1}(g^{-1}(x))=\frac{g^{-1}(x)+3}{2}=\frac{2x-9+3}{2}=\frac{2x-6}{2}=x-3,\\
       & (f^{-1}\circ g)(x)=f^{-1}(g(x))=\frac{g(x)+3}{2}=\frac{\frac{x+9}{2}+3}{2}=\frac{x+15}{4},\\
       &  ((f^{-1}\circ g)\circ f)(x)=\frac{f(x)+15}{4}=\frac{2x-3+15}{4}=\frac{2x+12}{4}=\frac{x+6}{2}.
      \end{align*}
\end{enumerate}
\end{example}

\subsection*{Zadaci za vje\v zbu}\index{Zadaci za vje\v zbu!funkcije}

\begin{enumerate}
   \item Date su funkcije $f(x)=\frac{x}{3},\,g(x)=\frac{x+2}{2},\,h(x)=x-6.$ Odrediti \\
      \begin{inparaenum}
         \item $f\circ g;\:$
         \item $g\circ f;\:$
         \item $(f\circ g)\circ h.$
      \end{inparaenum}
   \item Date su funkcije $f:A\mapsto B$ i $g:B\mapsto C,$ gdje su
       $A=\{x\in\mathbb{R}:0\leqslant x\leqslant 5 \},\: B=\{x\in\mathbb{R}: 1\leqslant y\leqslant 11 \} ,\:C=\{ x\in\mathbb{R}:2\leqslant z\leqslant 32\},$ i $f(x)=2x+1,\:g(x)=3x-1.$ \\
      \begin{inparaenum}
         \item Odrediti $f\circ g$ i $g\circ f;\:$
         \item Da li se kompozicijom $g\circ f$ skup $A$ preslikava u skup $C?$
     \end{inparaenum}
  \item Date su funkcije $f(x)=\frac{x-3}{x},\: g(x)=\frac{2}{x+1}.$ \\
     \begin{inparaenum}
         \item Odrediti definiciono podru\v cje funkcija $f$ i $g.$\:
         \item Odrediti $f\circ g$ i $g\circ f$ i domene dobijenih funkcija.
         \item Izra\v cunati $(f\circ g)(4)$ i  $(g\circ f)(6).$
    \end{inparaenum}
 \item Date su funkcije $f(x)=\frac{3x}{2x+5}$ i $g(x)=\frac{x-1}{2x+3}.$ \\
      \begin{inparaenum}
         \item Odrediti definiciono podru\v cje funkcija.\:
         \item Odrediti $f^{-1}$ i $g^{-1}.$
     \end{inparaenum}
  \item Za funkciju $f(x)=\frac{3x-5}{7}.$ Izra\v cunati \\
        \begin{inparaenum}
          \item $f^{-1}(1);\:$
          \item $(f^{-1}\circ f)(3);\:$
          \item $f^{-1}(-4).$
       \end{inparaenum}
  \item   Date su funkcije $f(3x+2)=x+3,\:g(x-1)=4x+1;\: k(x)=x^2-2x.$ Odrediti\\
         \begin{inparaenum}
             \item $f(x);\:g(x);\:k(x).$
             \item $f\circ g;\: f\circ k;\:g\circ (f\circ k).$
         \end{inparaenum}
  \item Date su $f(x)=\frac{3x-1}{x+7};\: g(x)=\frac{x+6}{x+9}.$ Odrediti\\
        \begin{inparaenum}
            \item Definiciono podru\v cje datih funkcija.
            \item $f^{-1},\:g^{-1}.$
       \end{inparaenum}
  \item Date su $f(x)=\frac{2x+1}{x-3};\: g\left(\frac{4x+5}{x+6}\right)=x.$ Odrediti\\
      \begin{inparaenum}
       	\item $g(x); $
   	     \item Definiciono podru\v cje datih funkcija.
   	     \item $f^{-1},\:g^{-1}\,, f\circ g^{-1}.$
      \end{inparaenum}
  \item Za funkciju $f\left( \frac{3x+1}{x+2}\right)=x+6$ odrediti \\
        \begin{inparaenum}
            \item $f(x);\:$
            \item $f^{-1}(x);\:$
            \item $f^{-1}\circ f.$
        \end{inparaenum}

  \item Za funkciju $f\left( \frac{-x+5}{4x+5}\right)=\frac{x}{7x+3}$ odrediti \\
        \begin{inparaenum}
        	\item $f(x);\:$
        	\item $f^{-1}(x);\:$
	       \item $f^{-1}\circ f.$
       \end{inparaenum}
\end{enumerate}

\chapter{Skup realnih brojeva}\label{poglavljeRealni}

Skup realnih brojeva $\mathbb{R}$ se sastoji od skupa racionalnih brojeva $\mathbb{Q}$ i skupa iracionalnih brojeva $\mathbb{I},$ pri \v cemu je $\mathbb{R}=\mathbb{Q}\cup\mathbb{I}$ i $\mathbb{Q}\cap\mathbb{I}=\emptyset.$
Tako\dj e vrijedi lanac inkluzija $\mathbb{N}\subset \mathbb{Z}\subset\mathbb{Q}\subset\mathbb{R}.$ Sa $\mathbb{N}$ je ozna\v cen skup prirodnih brojeva\footnote{prirodni brojevi ili pozitivni cijeli brojevi; nenegativni cijeli brojevi su $0,1,2,\ldots,$ tj. prirodni brojevi i jo\v s nula}, a sa $\mathbb{Z}$ skup cijelih brojeva. Sve ove skupove \' cemo detaljnije razmotriti u narednim sekcijama ovog poglavlja.

\section[P\lowercase{rirodni brojevi}]{Prirodni brojevi i skup prirodnih brojeva}

Istorijski gledano, od svih brojeva, prirodni brojevi su se najranije pojavili i kori\v steni su za brojanje raznih objekata. Skup svih prirodnih brojeva ozna\v cavamo sa $\mathbb{N}=\{1,2,3, \ldots\}.$  Prirodne brojeve mo\v zemo sabirati i mno\v ziti i rezultat ovih operacija je ponovo prirodan broj, pa ka\v zemo da je skup prirodnih brojeva zatvoren  u odnosu na ove dvije operacije. Problem se ve\' c javlja kada ho\' cemo da primijenimo operaciju oduzimanja u skupu priro-\\dnih brojeva. Na primjer rezultat oduzimanja  $1-2$ nije prirodan broj, tj. rezultat je broj koji ne postoji u skupu prirodnih brojeva. Dakle ve\' c za operaciju oduzimanja skup prirodnih brojeva je "mali" ili "preuzak", pa je bilo potrebno "pro\v siriti ovaj skup".

Skup prirodnih brojeva mo\v zemo uvesti ili definisati na vi\v se na\v cina. Jedan od na\v cina definisanja je preko Peanovih\footnote{Giuseppe Peano (27.avgust 1858.--20.april 1932. godine), bio je italijanski matemati\v car, ve\' ci dio karijere predavao je matematiku na Univerzitetu u Torinu. Jedan od osniva\v ca matemati\v cke logike i teorije skupova.}  aksioma. Peanovi aksiomi glase
\begin{enumerate}[$({A}1)$]
  \item\label{peano1} $1$ je prirodan broj, tj. $1\in \mathbb{N}.$
  \item\label{peano2}  Svaki prirodan broj $n$ ima svog sljedbenika $n'=n+1,$ koji je tako\dj e prirodan broj, tj. vrijedi
         $n\in\mathbb{N}\Rightarrow n+1\in\mathbb{N}.$
  \item $1$ nije sljedbenik ni jednog prirodnog broja.
  \item Ako je $m'=n',$ tada je $m=n,$ tj. svaki prirodan broj sljedbenik je najvi\v se jednog prirodnog broja.
  \item Ako je $M\subset\mathbb{N}$ i ako u skupu $M$ va\v ze \hyperref[peano1]{(A1)} i \hyperref[peano1]{(A2)} tada je $M=\mathbb{N}.$
\end{enumerate}

Ovi aksiomi defini\v su skup prirodnih brojeva $\mathbb{N}.$ Posljednji, peti aksiom je i osnova za jednu metodu dokaza poznatu u matematici pod nazivom "princip matemati\v cke indukcije" i o njemu \' ce biti vi\v se rije\v ci kasnije.

Osim kori\v stenja prirodnih brojeva za brojanje, potrebno je koristiti i ra\v cunske operacije sabiranja, mno\v zenja, oduzimanja i dijeljenja prirodnih brojeva.

\paragraph{Sabiranje prirodnih brojeva.}  Broj $a+b,$ je zbir ili suma brojeva $a$ i $b,$ brojevi $a$ i $b$ su sabirci. Uobi\v cajena oznaka za operaciju sabiranja je znak $"+".$ Operacija sabiranja je zatvorena u skupu prirodnih brojeva, tj. $(\forall a,b\in\mathbb{N})\:a+b\in\mathbb{N}.$

\paragraph{Oduzimanje prirodnih brojeva.} Broj $a-b$ je razlika brojeva $a$ i $b,$ broj $a$ je umanjenik a broj $b$ je umanjilac. Oznaka za operaciju oduzimanja je $"-".$ Operacija oduzimanja nije zatvorena u skupu prirodnih brojeva, naime, razlika dva prirodna broja nije uvijek prirodan broj, npr. $1-2\notin\mathbb{N}.$

\paragraph{Mno\v zenje  prirodnih brojeva.} Broj $ab$ ili $a\cdot b$ je proizvod brojeva $a$ i $b.$ Brojevi $a$ i $b$ su faktori, oznaka za operaciju mno\v zenja je $"\cdot".$ Operacija mno\v zenja je zatvorena u skupu $\mathbb{N},$ tj.   $(\forall a,b\in\mathbb{N})\:a\cdot b\in\mathbb{N}.$

\paragraph{Dijeljenje prirodnih brojeva.} Broj $a:b$ je koli\v cnik brojeva $a$ i $b,$ broj $a$ je dijeljenik, a  broj $b$ je djelilac. Operacija dijeljenja nije zatvorena u skupu prirodnih brojeva $\mathbb{N},$ npr. $1:2\notin\mathbb{N}.$ Oznaka za dijeljenje je $":"$.\\

\v Cesto se javlja potreba da pro\v sirimo skup prirodnih brojeva sa brojem $0,$ pa tada koristimo oznaku $\mathbb{N}_0.$ Dakle $\mathbb{N}_0 =\mathbb{N} \cup \{0\}.$

\subsection{Reprezentacija prirodnih brojeva}

Uobi\v cajeno predstavljanje brojeva je u decimalnoj notaciji. Zapisujemo brojeve koriste\' ci cifre  $\{0\,,1\,,2\,,3\,,4\,,5\,,6\,,7\,,8\,,9\}$ da predstavimo umno\v ske stepena baze $10.$ Na primjer broj $234597$ mo\v zemo zapisati ovako:
\[2\cdot 10^5+3\cdot 10^4+4\cdot 10^3+5\cdot 10^2+9\cdot 10+7\cdot 10^0.\]

Nema posebnog razloga za\v sto koristimo broj $10$ kao bazu na\v se notacije, mo\v zda \v sto imamo $10$ prstiju na rukama. Babilonci su koristili bazu $60$, Maje bazu $12,$ dok dana\v sni kompjuteri koriste bazu $2.$

Sljede\' ca teorema ka\v ze da svaki prirodan broj ve\' ci od $1$ mo\v zemo koristiti kao bazu.
\begin{theorem}
   Neka je $b$ prirodan broj ve\' ci od $1.$ Tada svaki prirodan broj $n$ mo\v zemo na jedinstven na\v cin predstaviti u obliku
   \begin{equation}a_k b^k+a_{k-1}b^{k-1}+\ldots+a_1b^1+a_0,
      \label{ekspanzija1}
   \end{equation}
gdje je $k$ prirodan broj, $a_j$ je broj iz $\mathbb{N}_0=\mathbb{N }\cup\{0\}$ za koji je $0\leqslant a_j\leqslant b-1$ za $j=0,1,\ldots,k$ i $a_k\neq 0.$
\end{theorem}

Broj $b$ iz razvoja \eqref{ekspanzija1} nazivamo baza razvoja, ekspanzije ili predstavljanja. Razvoj sa bazom $10$ (na\v sa uobi\v cajena baza), zovemo  decimalna notacija, predstavljanje ili razvoj. U slu\v caju  baze $2,$  ka\v zemo da je to binarno predstavljanje, za bazu $8$ oktalno i za bazu $16$ heksadecimalno predstavljanje. Koeficijenti $a_j$ su cifre razvoja.

Da bi razlikovali predstavljanje brojeva u razli\v citim bazama, koristimo specijalnu notaciju. Pi\v semo
\[(a_ka_{k-1}\ldots a_1a_0)_b\]
da bi predstavili broj
\[a_kb^k+a_{k-1}b^{k-1}+\ldots+a_1b+a_0.\]

\begin{example}
Predstaviti brojeve $(236)_7$ i $(10010011)_2$ u decimalnoj notaciji.\\\\
Rje\v senje:\\\\
Vrijedi $(236)_7=2\cdot 7^2+3\cdot 7+6=2\cdot 49+21+6=125$ i $(10010011)_2\\=1\cdot 2^7+1\cdot 2^4+1\cdot 2+1=128+16+2+1=147.$
\end{example}

\begin{example}
 Predstaviti brojeve $125 $ u bazi $7$ i $147$ u bazi $2.$\\\\
 Rje\v senje:\\\\
 Postupamo na sljede\' ci na\v cin
 \begin{align*}
        125&=17\cdot 7+{\color{magenta}6}\\
        17&=2\cdot 7+{\color{magenta}3}\\
        2&=0\cdot 7+{\color{magenta}2},
 \end{align*}
te je $125=(236)_7,$  i
 \begin{align*}
     147&=73\cdot 2+{\color{magenta}1}\\
     73&=36\cdot 2+{\color{magenta}1}\\
     36&=18\cdot 2+{\color{magenta}0}\\
     18&=9\cdot 2+{\color{magenta}0}\\
     9&=4\cdot 2+{\color{magenta}1}\\
     4&=2\cdot 2+{\color{magenta}0}\\
     2&=1\cdot 2+{\color{magenta}0}\\
     1&=0\cdot 2+{\color{magenta}1}
 \end{align*}
$147=(10010011)_2.$
\end{example}

\begin{example}
  Rije\v siti jedna\v cine
      \begin{enumerate}[$(a)$]
          \item $2012_3\cdot x_2=2183.$
          \item $5\cdot x_4=11001_2\cdot 2222_3.$
          \item $x_6:321_4=432_5.$
      \end{enumerate}
Rje\v senje:
\begin{enumerate}[$(a)$]
    \item Prvo pretvorimo $2012_3$ u decimalnu notaciju, vrijedi $2012=2\cdot 3^3+1\cdot 3+2=59,$ zatim $\frac{2183}{59}=37.$ Na kraju $37$ pretvorimo u binarni broj, tj. predstavimo ga u binarnoj notaciji
    \begin{align*}
         37&=18\cdot 2+1\\
         18&=9\cdot 2+0\\
         9&=4\cdot 2+1\\
         4&=2\cdot 2+0\\
         2&=1\cdot 2+0\\
         1&=0\cdot 2+1,
    \end{align*}
   pa je $x_2=100101.$
\item  Pretvorimo brojeve u decimalnu notaciju, vrijedi $11001_2=1\cdot 2^4+1\cdot 2^3+1=25,$ $2222_3=2\cdot 3^3+2\cdot 3^2+2\cdot 3+2=80,$ dalje je $5x=25\cdot 80\Leftrightarrow x=\frac{25\cdot 80}{5}=400.$ Sada  pretvorimo $400$ iz decimalne u notaciju sa bazom $b=4$
  \begin{align*}
         400&=100\cdot 4+0\\
         100&=25\cdot 4+0\\
         25&=6\cdot 4+1\\
         6&=1\cdot 4+2\\
         1&=0\cdot 4+1,
   \end{align*}
pa je  $x_4=12100_4.$
\item Vrijedi $321_4=3\cdot 4^2+2\cdot 4+1=57,$ $ 432_5=4\cdot 5^2+3\cdot 5+2=117,$ pa je $x=57\cdot 117=6669,$ te je na kraju
  \begin{align*}
      6669&=1111\cdot 6+3\\
      1111&=185\cdot 6+1\\
      185&=30\cdot 6+5\\
      30&=5\cdot 6+0\\
      5&=0\cdot 6+5,
  \end{align*}
  $x_6=50513_6.$
\end{enumerate}
\end{example}

\subsection{Matemati\v cka indukcija}\index{matemati\v cka indukcija}

Princip matemati\v cke indukcije koristan je alat za dokazivanje raznih rezultata sa priro-\\dnim brojevima. Princip matemati\v cke indukcije zasnovan je na petom Peanovom aksi-\\omu.
\paragraph{Peti Peanov aksiom:}   Neka je $M\subset\mathbb{N},$   ako vrijedi $1\in M$ i $n\in M\Rightarrow n+1\in M,$  tada je $M=\mathbb{N}.$

Da bi pokazali matemati\v ckom indukcijom da neki rezultat vrijedi za svaki prirodan broj,  potrebno je uraditi sljede\' ca tri koraka. Prvi korak je provjera, tj. treba provjeriti da li taj rezultat vrijedi za prvi broj iz tvrdnje (naj\v ce\v s\'ce je to broj $1$).
 Ovaj korak zovemo baza indukcije. Drugi korak je postavljanje pretpostavke da tvrdnja vrijedi za neki prirodan broj $n.$ Ova pretpostavka se naziva i induktivna hipoteza. Tre\' ci korak je dokazivanje da iz induktivne hipoteze
 slijedi da je tvrdnja ta\v cna za $n+1.$ Ovaj korak zovemo korak indukcije.
Nakon ovih uspje\v sno obavljenih koraka, zaklju\v cujemo na osnovu principa matemati\v cke indukcije da je posmatrana tvrdnja ta\v cna za sve prirodne brojeve.
Ovo je osnovna \v sema dokaza, ali postoje i brojne druge varijante primjene metode matemati\v cke indukcije.
\begin{example}
 Dokazati matemati\v ckom indukcijom $(\forall n\in\mathbb{N })\:1+2+\ldots+n=\frac{n(n+1)}{2}.$\\ \ \\
 Rje\v senje:
 \begin{enumerate}
 	\item Za $n=1$ vrijedi $1=\frac{1(1+1)}{2}\Leftrightarrow 1=1.$ \\
\item Pretpostavimo sada da je formula ta\v cna za neko $n,$ tj. da vrijedi
\[1+2+\ldots+n=\frac{n(n+1)}{2}.\]
\item Ako sada dodamo i lijevoj strani prethodne jednakosti $n+1,$  dobi\'cemo
\[1+2+\ldots+n+(n+1)=\frac{n(n+1)}{2}+n+1. \]
Naravno ako je pretpostavka ta\v cna onda je i ova jednakost ta\v cna. Lijeva strana predstavlja zbir prvih $n+1$ prirodnih brojeva i tu ne treba ni\v sta mijenjati. Na desnoj strani je izraz $\frac{n(n+1)}{2}+n+1$ i trebamo pokazati da je on jednak $\frac{(n+1)(n+2)}{2}.$ Zaista vrijedi:
\begin{align*}
 \frac{n(n+1)}{2}+n+1=\frac{n(n+1)+2(n+1)}{2}=\frac{(n+1)(n+2)}{2},
\end{align*}
\v sto zna\v ci da smo uspjeli dokazati da formula vrijedi i za $n+1$.
\end{enumerate}
Na osnovu principa matemati\v cke indukcije, dokazali smo da je jednakost ta\v cna za svako $n\in\mathbb{N }.$
\end{example}

\begin{example}
	Dokazati matemati\v ckom indukcijom $(\forall n\in\mathbb{N})\:1+3+6+\ldots+\frac{n(n+1)}{2}=\frac{n(n+1)(n+2)}{6}.$\\ \ \\
\noindent 	Rje\v senje.
\begin{enumerate}
	\item  Neka je $n=1,$ sada je $1=\frac{1\cdot(1+1)(1+2)}{6}    \Leftrightarrow 1=1.$ Dakle prvi korak (baza ili osnova indukcije) je ispunjen. \\
\item Pretpostavimo da je formula ta\v cna za neko $n \in \mathbb{N} :$
\begin{equation*}
   1+3+\ldots+\frac{n(n+1)}{2}=\frac{n(n+1)(n+2)}{6}.
   \label{indukcija2}
\end{equation*}
\item Trebamo sada pokazati da iz pretpostavke \eqref{indukcija2} slijedi ta\v cnost formule za\\ $n+1 \in \mathbb{N},$ odnosno da vrijedi:
\[ 1+3+\ldots+\frac{n(n+1)}{2}+\frac{(n+1)(n+2)}{2}=\frac{(n+1)(n+2)(n+3)}{6}.\]
Dodajmo na obje strane jednakosti \eqref{indukcija2} izraz $\frac{(n+1)(n+2)}{2},$ pa je
\[1+3+\ldots+\frac{n(n+1)}{2}+\frac{(n+1)(n+2)}{2}=\frac{n(n+1)(n+2)}{6}+\frac{(n+1)(n+2)}{2}.\]
Lijeva strana predstavlja zbir $n+1$ razlomaka oblika $\frac{n(n+1)}{2}$ i tu se nema \v sta dokazivati, ali treba pokazati da je desna strana jednaka
$\frac{(n+1)(n+2)(n+3)}{6}.$ Dalje je
\begin{multline*}
   \frac{n(n+1)(n+2)}{6}+\frac{(n+1)(n+2)}{2}\\= \frac{n(n+1)(n+2)}{6}+\frac{(n+1)(n+2)\cdot 3}{2\cdot3} =\frac{(n+1)(n+2)(n+3)}{6},
\end{multline*}
\v cime smo pokazali ta\v cnost formule i za $n+1.$
\end{enumerate}
Na osnovu principa matemati\v cke indukcije, jednakost je ta\v cna za svako $n\in\mathbb{N}.$
\end{example}

\begin{example}
     Dokazati matemati\v ckom indukcijom $3\mid 5^n+2^{n+1},\,(\forall n\in\mathbb{N}).$\\  \ \\
Rje\v senje.
\begin{enumerate}
	\item Za $n=1$ vrijedi, $5^1+2^{1+1}=5+4=9,\:9:3=3.$ Prvi korak je pokazan.
	\item Induktivna hipoteza: izraz $5^n+2^{n+1}$ je djeljiv sa $3.$
	\item Doka\v zimo da iz pretpostavke (induktivne hipoteze) slijedi djeljivost sa $3$ izraza $5^{n+1}+2^{n+2}.$ Ako je $3\mid 5^n+2^{n+1},$ onda postoji prirodan broj $a$ takav da je $5^{n}+2^{n+1}=3\cdot a.$ Sada je
\begin{align*}
 5^{n+1}+2^{n+2}=&5\cdot 5^n+2\cdot 2^{n+1}=5\underbrace{\left( 5^n+2^{n+1}\right)}_{=3a}-3\cdot 2^{n+1}=5\cdot 3a-3\cdot 2^{n+1}\\
                =&3\left(5 a-2^{n+1}\right)=3\cdot b,\:b=5a-2^{n+1}\in\mathbb{Z}.
\end{align*}
Ovim smo pokazali da iz $3\mid 5^n+2^{n+1}$ slijedi $3\mid 5^{n+1}+2^{n+2}.$
\end{enumerate}
Sada na osnovu principa matemati\v cke indukcije, zaklju\v cujemo da je tra\v zeni izraz djeljiv sa $3$ za svako $n\in\mathbb{N }.$
\end{example}

\begin{example}
Dokazati   matemati\v ckom indukcijom $6\mid n^3+11n, \, (\forall n\in\mathbb{N}).$\\ \ \\
Rje\v senje.
\begin{enumerate}
	\item Za $n=1$ vrijedi $1^3+11\cdot 1=12,\:12:6=2.$ Prvi korak je pokazan.
	\item Pretpostavimo da za neko $n\in\mathbb{N}$ vrijedi $6\mid n^3+11n.$
	\item Iz $6\mid n^3+11n$ slijedi da postoji prirodan broj $a$ takav da je $n^3+11n=6\cdot a.$ Dalje je
         \begin{align*}
             (n+1)^3+11(n+1)&=n^3+3n^2+3n+1+11n+11\\
                            &=n^3+11n+3n^2+3n+12\\
                            &=\underbrace{n^3+11n}_{=6\cdot a}+3n(n+1)+12.
          \end{align*}
Treba da poka\v zemo da je $3n(n+1)$ djeljivo sa $6.$ Ovaj izraz je ve\'c djeljiv sa $3,$ a brojevi $n$ i $n+1$ su uzastopni prirodni brojevi, pa je jedan od njih paran a drugi neparan, tako da je proizvod $n(n+1)$ paran, a samim tim i djeljiv sa $2.$ S obzirom da je izraz $3n(n+1)$ istovremeno djeljiv i sa $2$ i sa $3,$ to zna\v ci da je djeljiv i sa $6.$ Na kraju je
\[(n+1)^3+11(n+1)=6\cdot a+6\cdot b+12=6\cdot c,\, (b,c\in\mathbb{N}).\]
pa je dokazana djeljivost.
\end{enumerate}
Na osnovu principa matemati\v cke indukcije, dokazali smo da je tra\v zeni izraz djeljiv sa $6$ za svako $n\in\mathbb{N }.$
\end{example}
Jedna od zna\v cajnijih nejednakosti koja se \v cesto koristi u matematici je tzv. Bernoullijeva (Bernulijeva) \footnote{Daniel Bernoulli (8.februar 1700.--27.mart 1782. godine) bio je \v svajcarski matemati\v car i fizi\v car i jedan od najistaknutijih matemati\v cara u porodici Bernoulli. Posebno \'ce biti zapam\'cen po primjeni matematike u mehanici fluida, te pionirskom doprinosu u razvoju vjerovatno\'ce i statistike. } nejednakost. I ona se mo\v ze dokazati matemati\v ckom indukcijom.
\ \\
\begin{example}
Dokazati Bernoullijevu nejednakost: \\
Za svaki prirodan $n\geq 2$ vrijedi $(1+h)^n > 1+nh,$ gdje je $h\neq 0$ i $h>-1.$ \\  \ \\
Rje\v senje.\\\\
Ovdje \' cemo provjeru po\v ceti od broja $n=2,$ jer je to prvi od prirodnih brojeva za koje se tvrdi ta\v cnost nejednakosti.
\begin{enumerate}
	\item Za $n=2$ je $(1+h)^2=1+2h+h^2>1+2h$ jer je $h^2>0$.
	\item Pretpostavka: nejednakost je ta\v cna za neko $n=k>2,$ tj. $(1+h)^k > 1+kh.$
	\item Doka\v zimo da iz pretpostavke slijedi da je nejednakost je ta\v cna za $n=k+1.$ Pomno\v zimo li obje strane prethodne nejednakosti (pretpostavke)
sa $(1+h)>0,$ dobi\' cemo \\
 $(1+h)^{k+1}>(1+kh)(1+h)=1+(1+k)h+kh^2>1+(1+k)h,$
 \v sto je i trebalo pokazati.
\end{enumerate}
Na osnovu principa matemati\v cke indukcije dokazali smo tra\v zenu nejednakost.
\end{example}
Dakle, u koraku provjere kre\' cemo od prvog broja za koji treba pokazati da je tvrdnja ta\v cna.\\
\begin{remark}
Princip matemati\v cke indukcije mo\v ze i da se uop\v sti tako da glasi: \\
Neka je dat niz tvr\dj enja $A_s, A_{s+1}, A_{s+2}, ...$ gdje je $s$ neki prirodan broj i ako vrijedi:
\begin{enumerate}[a)]
	\item $A_s$ je ta\v cno
	\item za svaki broj $r \geq s$ iz istinitosti tvr\dj enja $A_r$ slijedi istinitost tvr\dj enja $A_{r+1},$
\end{enumerate}	
tada su sva tvr\dj enja $A_s, A_{s+1}, A_{s+2}, ...$ ta\v cna. Drugim rije\v cima $A_s$ je ta\v cno za svako $n\geq s.$
\end{remark}

\subsection{Binomni obrazac}\index{binomni obrazac}

Binomni obrazac nam daje formulu po kojoj mo\v zemo izraz $(a+b)^n$ (stepen binoma) predstaviti kao sumu sabiraka odre\dj enog oblika.
Prije nego \v sto damo ovu formulu, neophodno je uvesti neke pojmove.\\
\begin{definition}
Funkcija faktorijel (u oznaci $!$) svakom broju iz skupa $\mathbb{N} \cup \{0\}$ pridru\v zuje prirodan broj, tako da je
$0!=1$ i $n!=1\cdot 2 \cdot ... \cdot n \,(\forall n \in \mathbb{N}).$
\end{definition} \index{binomni obrazac!faktorijel, faktorijelna funkcija}

\newpage

\index{binomni obrazac!binomni koeficijent}
\begin{definition}
Funkcija koja svakom uredjenom paru $(n,k) \in \mathbb{N}_0\times \mathbb{N}_0$ pridru\v zuje nenegativan cijeli broj ozna\v cen sa
$\left(\begin{array}{c}
n \\
k
\end{array} \right)$
i definisan sa:\\
\begin{equation*}
\left(\begin{array}{c}
0 \\
0
\end{array} \right)= \left(\begin{array}{c}
n \\
0
\end{array} \right)=1, \textrm{za sve} \quad n \in \mathbb{N} \quad \textrm{i}
\end{equation*}
\begin{equation*}
 \left(\begin{array}{c}
n \\
k
\end{array} \right)= \frac{n(n-1)\cdot\ldots\cdot   (n-k+1)}{k!} \, \textrm{za sve} \, n \in \mathbb{N}_0, \, k \in \mathbb{N},
\end{equation*}
naziva se binomni koeficijent.

\end{definition}
\begin{example}
\begin{equation*}
5!=1\cdot 2 \cdot 3 \cdot 4 \cdot 5=120, \quad \left(\begin{array}{c}
7 \\
2
\end{array} \right)= \frac{7 \cdot 6 }{2!}=21.
\end{equation*}
\end{example}
Faktorijel i binomni koeficijent imaju primjene u kombinatorici. Neka imamo $n$ razli\v citih elemenata $k_1,k_2,\ldots , k_n.$ Ovih $n$ elemenata
mogu\'{c}e je poredati na vi\v se razli\v citih na\v cina. Svaka takva ure\dj ena $n$--torka se zove permutacija bez ponavljanja.
Broj svih mogu\' cih permutacija bez ponavljanja elemenata $n$--\v clanog skupa je $n !.$ Ako od $n$ razli\v citih elemenata pravimo slogove
od po $k$ razli\v citih elemenata, a pri tom ne pridajemo va\v znost poretku unutar sloga, dobijamo kombinacije $n$--tog reda i $k$--tog razreda
bez ponavljanja. Broj ovakvih kombinacija jednak je upravo binomnom koeficijentu $\binom{n}{k}.$
Ako razlomak
$\binom{n}{k}$
pro\v sirimo sa $(n-k)!$ (kad je $n>k$), onda dobijamo
\begin{equation*}
\left(\begin{array}{c}
n \\
k
\end{array} \right)
 =\frac{n(n-1)\cdot\ldots\cdot (n-k+1)\cdot (n-k)(n-k-1) \cdot\ldots\cdot 1}{k!\cdot (n-k)!},
\end{equation*}
 odnosno
\begin{equation*}
\left(\begin{array}{c}
n \\
k
\end{array} \right) =\frac{n !}{k!(n-k)!}.
\end{equation*}

Za binomne koeficijente kad je $n\geq k,$ vrijede sljede\' ce osobine:
\begin{enumerate}
\item  \begin{equation*}
\left(\begin{array}{c}
n \\
k
\end{array} \right)=
\left(\begin{array}{c}
n \\
n-k
\end{array} \right)
\end{equation*}
\item Pascalova formula \begin{equation*}
\left(\begin{array}{c}
n \\
k
\end{array} \right) +\left(\begin{array}{c}
n \\
k+1
\end{array} \right)=\left(\begin{array}{c}
n+1 \\
k+1
\end{array} \right)
\end{equation*}
\item \begin{equation*}
\left(\begin{array}{c}
n \\
k
\end{array} \right)=
\frac{n}{k}\left(\begin{array}{c}
n-1 \\
k-1
\end{array} \right)
\end{equation*}

\item \begin{equation*}
k\left(\begin{array}{c}
n \\
k
\end{array} \right)=
n\left(\begin{array}{c}
n-1 \\
k-1
\end{array} \right).
\end{equation*}
\end{enumerate}
\newpage
Prva osobina zna\v ci simetri\v cnost binomnih koeficijenata.
Binomne koeficijente mo\v zemo poredati u obliku beskona\v cnog trougla,
kojeg nazivamo Pascalov trougao, u \v cast Blaise Pascala\footnote{Blaise Pascal (19.juni 1623.--19.avgust 1662. godine) bio je francuski matemati\v car, fizi\v car, izumitelj, filozof, pisac i katoli\v cki teolog. Dao veliki doprinos u matematici, fizici te drugim prirodnim i primijenjenim naukama.} koji ga je opisao
u 17.vijeku (vidjeti Tabela \ref{tablePascalov1}). Ali ovo je bilo poznato i mnogim matemati\v carima prije njega. \index{binomni obrazac!Pascalov trougao}

\begin{figure}[!h]
	\begin{adjustwidth}{0cm}{}

	
	\begin{tabular}{cc cc cc cc cc cc cc }
		$n=0$ &&&&&&&$\binom{0}{0}$&&&&&&\\
\\
		$n=1$ &&&&&&$\binom{1}{0}$&&$\binom{1}{1}$&&&&&\\
\\
		$n=2$ &&&&&$\binom{2}{0}$&&$\binom{2}{1}$&&$\binom{2}{2}$&&&&\\
\\
		$n=3$ &&&&$\binom{3}{0}$&&$\binom{3}{1}$&&$\binom{3}{2}$&&$\binom{3}{3}$&&&\\
\\
		$n=4$ &&&$\binom{4}{0}$&&$\binom{4}{1}$&&$\binom{4}{2}$&&$\binom{4}{3}$&&$\binom{4}{4}$&&\\
\\
		$n=5$ &&$\binom{5}{0}$&&$\binom{5}{1}$&&$\binom{5}{2}$&&$\binom{5}{3}$&&$\binom{5}{4}$&&$\binom{5}{5}$&\\
\\
		$n=6$ &$\binom{6}{0}$&&$\binom{6}{1}$&&$\binom{6}{2}$&&$\binom{6}{3}$&&$\binom{6}{4}$&&$\binom{6}{5}$&&$\binom{6}{6}$\\
		$\vdots$ &&&&&&&$\vdots$&&&&&&\\\\

		$n=0$ &&&&&&&1&&&&&&\\
\\
		$n=1$ &&&&&&1&&1&&&&&\\
\\
		$n=2$ &&&&&1&&2&&1&&&&\\
\\
		$n=3$ &&&&1&&3&&3&&1&&&\\
\\
		$n=4$ &&&1&&4&&6&&4&&1&&\\
\\
		$n=5$ &&1&&5&&10&&10&&5&&1&\\
\\
		$n=6$ &1&&6&&15&&20&&15&&6&&1\\
		$\vdots$ &&&&&&&$\vdots$&&&&&&
		
	\end{tabular}
	\caption{Pascalov trougao}
	\label{tablePascalov1}

 \end{adjustwidth}
\end{figure}

Ovdje dolazi do izra\v zaja Pascalova formula, tj. da je zbir svaka dva susjedna binomna koeficijenta
u istom redu i $k$ i $(k+1)$--oj koloni jednak binomnom koeficijentu u sljede\' cem redu i $(k+1)$--oj koloni.
Sada mo\v zemo navesti i dokazati binomni obrazac.
\newpage

\begin{theorem}
Za bilo koja dva realna broja $a$ i $b$ i bilo koji prirodan broj $n$ vrijedi:
\begin{equation}\label{bob}
\left(a+b \right)^n=\displaystyle\sum_{k=0}^{n}\left(\begin{array}{c}
n \\
k
\end{array} \right) a^{n-k}\cdot b^k.
\end{equation}
\end{theorem}

\begin{proof}[Dokaz]
Teoremu \' cemo dokazati metodom matemati\v cke indukcije. \\
Formulu (\ref{bob}) mo\v zemo napisati u obliku:
\begin{equation*}
\left(a+b \right)^n=\left(\begin{array}{c}
n \\
0
\end{array} \right)a^n+\left(\begin{array}{c}
n \\
1
\end{array} \right)a^{n-1}b+\ldots+\left(\begin{array}{c}
n \\
n-1
\end{array} \right)ab^{n-1}+\left(\begin{array}{c}
n \\
n
\end{array} \right)b^n.
\end{equation*}
\begin{equation*}
=a^n+\left(\begin{array}{c}
n \\
1
\end{array} \right)a^{n-1}b+\ldots +\left(\begin{array}{c}
n \\
n-1
\end{array} \right)ab^{n-1}+b^n.
\end{equation*}
\begin{enumerate}
	\item Provjeravamo da li je jednakost ta\v cna za $n=1$. Lijeva strana je $\left(a+b \right)^1,$ a desna je $1 \cdot a^1+1\cdot b^1,$
pa je data jednakost ta\v cna.

\item Pretpostavimo da je tvrdnja ta\v cna za $n=k.$ To zna\v ci da je
\begin{equation}
\left(a+b \right)^k=\displaystyle\sum_{j=0}^{k}\left(\begin{array}{c}
k \\
j
\end{array} \right) a^{k-j}\cdot b^j.
\end{equation}

\item Doka\v zimo da iz pretpostavke slijedi da je tvrdnja teorema ta\v cna za $n=k+1.$
Krenimo od lijeve strane jednakosti koju trebamo dokazati,
\begin{align*}
  \left(a+b \right)^{k+1}
     &= (a+b)\cdot (a+b)^k= (a+b) \cdot \displaystyle\sum_{j=0}^{k}\left(\begin{array}{c}k \\ j \end{array} \right) a^{k-j}\cdot b^j \\
     &=\displaystyle\sum_{j=0}^{k}\left(\begin{array}{c} k \\j \end{array} \right) a^{k+1-j}\cdot b^j
            + \displaystyle\sum_{j=0}^{k}\left(\begin{array}{c} k \\ j \end{array} \right) a^{k-j}\cdot b^{j+1}.
\end{align*}


Ako u lijevoj sumi izdvojimo prvi sabirak (tj. kad je $j=0$), a u desnoj sumi izdvojimo zadnji sabirak (tj. kad je $j=k$),
dobijamo
\begin{align*}
(a+b)^{k+1}
    &=a^{k+1}+\displaystyle\sum_{j=1}^{k}\left(\begin{array}{c}k \\j\end{array} \right) a^{k+1-j}\cdot b^j
        + \displaystyle\sum_{j=0}^{k-1}\left(\begin{array}{c}k \\ j \end{array} \right) a^{k-j}\cdot b^{j+1}+b^{k+1} \\
    & =a^{k+1}+\displaystyle\sum_{j=1}^{k}\left(\begin{array}{c}k \\ j\end{array} \right) a^{k+1-j}\cdot b^j
           + \displaystyle\sum_{j=1}^{k}\left(\begin{array}{c} k \\ j-1 \end{array} \right) a^{k-(j-1)}\cdot b^{j}+b^{k+1} \\
    &=a^{k+1}+\displaystyle\sum_{j=1}^{k}\left[\left(\begin{array}{c} k \\ j \end{array} \right)
          +\left(\begin{array}{c} k \\ j-1 \end{array} \right) \right] a^{k+1-j}\cdot b^j +b^{k+1}.
\end{align*}


Na osnovu Pascalove formule (2. osobine binomnih koeficijenata) je
\begin{equation*}
\left(\begin{array}{c}
k \\
j
\end{array} \right) +\left(\begin{array}{c}
k \\
j-1
\end{array} \right)=\left(\begin{array}{c}
k+1 \\
j
\end{array} \right),
\end{equation*}  \index{binomni obrazac!binomni obrazac}
\newpage
pa kona\v cno dobijamo
\begin{displaymath}
(a+b)^k=a^{k+1}+\displaystyle\sum_{j=1}^{k}\left(\begin{array}{c}
k+1 \\
j
\end{array} \right) a^{k+1-j}\cdot b^j +b^{k+1}=\displaystyle\sum_{j=0}^{k+1}\left(\begin{array}{c}
k+1 \\
j
\end{array} \right) a^{k+1-j}\cdot b^j .
\end{displaymath}
Dakle, tvrdnja teoreme vrijedi i za $n=k+1$.
\end{enumerate}
Na osnovu principa matemati\v cke indukcije, tvrdnja teoreme vrijedi za svaki prirodan broj $n$.
\end{proof}
Poop\v stenje binomnog obrasca na racionalne eksponente dao je Isaac Newton \footnote{Isaac Newton (25.decembar 1642.--20.mart 1726/27. godine) bio je engleski matemati\v car, fizi\v car, astronom, teolog i autor (opisan u njegovo vrijeme kao "prirodni filozof"). \v Siroko je rasprostranjeno mi\v sljenje da je jedan od najve\'cih matemati\v cara i jedan od najuticajnijih nau\v cnika svih vremena.}.
Primijetimo da je u razvoju $n$--tog stepena binoma (\ref{bob}) ukupno $n+1$ sabiraka, a $(k+1)$--vi \v clan iznosi
\begin{equation}\label{kClanBinoma}
T_{k+1}=\left(\begin{array}{c}
n \\
k
\end{array} \right) a^{n-k}\cdot b^k.
\end{equation}

\begin{example}
  Razviti po binomnom obrascu $(3+2x)^4.$ 	\\ \ \\
  \noindent Rje\v senje.\\\\
   Za dati binom je $a=3,$ $b=2x$ i $n=4.$ U binomni obrazac $(a+b)^n=\binom{n}{0}a^n+\binom{n}{1}ab^{n-1}+\ldots+ \binom{n}{n-1}ab^{n-1}+\binom{n}{n}b^n$ uvrstimo date podatke,  dobijamo
  \begin{align*}
  (3+2x)^4&=3^4+\binom{4}{1}3^3\cdot(2x)+\binom{4}{2}3^2\cdot(2x)^2+\binom{4}{3}3\cdot (2x)^3+\binom{4}{4}(2x)^4\\
          &= 81+216x+216x^2+96x^3+16x^4.
  \end{align*}
\end{example}

\begin{example}
Odrediti peti  \v clan u razvijenom obliku binoma $\left(x^{\frac{1}{2}}+x^{\frac{2}{3}} \right)^{12}.$ \\\\
\noindent Rje\v senje. \\
Binomni obrazac mo\v zemo pisati u obliku $\displaystyle(a+b)^n=\sum_{k=0}^{n}\binom{n}{k}a^{n-k}b^{k},$ pa prvom \v clanu  odgovara $k=0,$ drugom \v clanu $k=1,$ tj. kao \v sto je u \eqref{kClanBinoma} napisano, $(k+1)$--vi \v clan mo\v zemo pisati u obliku $T_{k+1}=\binom{n}{k}a^{n-k}b^k.$ Kako treba da izra\v cunamo $5$--ti \v clan, to je $k=4,$ pa je
\begin{align*}
T_{5}&=\binom{12}{4}\left(x^{\frac{1}{2}}\right)^{12-4} \cdot \left( x^{\frac{2}{3}}\right)^{4}=\frac{12!}{4!(12-4)!}x^{4}\cdot x^{\frac{8}{3}}=495x^{\frac{20}{3}}.
\end{align*}
\end{example}
\begin{example}
Odrediti \v clan koji ne sadr\v zi $x$ u razvijenom obliku binoma $(x+x^{-2})^{12}.$\\ \\
\noindent Rje\v senje. \\\\
Koristimo ponovo formulu  $T_{k+1}=\binom{n}{k}a^{n-k}b^k.$ Potrebno je $k$ izra\v cunati, to radimo iz uslova
\begin{align*}
 T_{k+1}&=\binom{n}{k}a^{n-k}b^k=\binom{12}{k}\left(x\right)^{12-k}\cdot \left(x^{-2}\right)^{k}=\binom{12}{k}x^{12-3k}.
\end{align*}
\v Clan koji ne sadr\v zi $x$ zbog toga \v sto je eksponent od $x$ jednak $0,$ dakle $12-3k=0\Leftrightarrow k=4,$ tj. peti \v clan ne sadr\v zi $x$
\[T_5=\binom{12}{4}x^{12-4}\cdot \left (x^{-2} \right)^4=\frac{12!}{4!\cdot 8!}x^{8}\cdot x^{-8}=495.\]  	
\end{example}

\begin{example}
	 Dat je binom $\displaystyle\left(\sqrt{2^x}+\frac{1}{\sqrt{2^{x-1}}}\right)^{n}.$
	 \begin{enumerate}[$(a)$]
	 	\item Odrediti $n$ tako da je zbir binomnih koeficijenta posljednja tri \v clana u razvijenom obliku binoma jednak $22.$
	 	\item Odrediti onu vrijednost od $x$ za koju je zbir tre\'ceg i petog \v clana u razvijenom obliku binoma jednak $105.$
	 \end{enumerate}
\noindent Rje\v senje.
\begin{enumerate}[$(a)$]
	\item Vrijedi
	   \begin{align*}
	    \binom{n}{n-2}+\binom{n}{n-1}+\binom{n}{n}=22\Leftrightarrow \frac{n(n-1)}{2}+n+1=22\Leftrightarrow n^2+n-42=0,
	   \end{align*}
	   pa su rje\v senja ove kvadratne jedna\v cine $n_1=-7\vee n_2=6.$  Po\v sto $n$ mora biti prirodan broj, vrijedi $n=6.$
	\item   Sada je
	\begin{align*}
	T_3+T_5=135
	     &\Leftrightarrow \binom{6}{2}\left(\sqrt{2^x}\right)^4\cdot \left(\frac{1}{\sqrt{2^{x-1}}} \right)^2
	         +\binom{6}{4}\left(\sqrt{2^x}\right)^2\cdot \left(\frac{1}{\sqrt{2^{x-1}}} \right)^4=135 \\
	     &\Leftrightarrow \frac{6!}{2!4!}2^{2x}\cdot 2^{-(x-1)}+\frac{6!}{4!2!}2^x\cdot 2^{-2(x-1)}=135 \\
	     &\Leftrightarrow  2\cdot 2^x+4\cdot 2^{-x}=9\\
	     &\Leftrightarrow  2\cdot 2^{2x}-9\cdot 2^{x}+4=0.
    \end{align*}
    Uvedimo smjenu $2^x=t,$ dobijamo kvadratnu jedna\v cinu $2t^2-9t+4=0,$ \v cija su rje\v senja $t_1=\frac{1}{2}\vee t_2=4,$ odnosno
    $2^x=\frac{1}{2}\Leftrightarrow x_1=-1$ i $2^x=4\Leftrightarrow x_2=2.$
\end{enumerate}
	
\end{example}

\section[C\lowercase{ijeli brojevi}]{Cijeli brojevi i skup cijelih brojeva}
Pro\v sirivanjem skupa prirodnih brojeva dobijamo skup cijelih brojeva $\mathbb{Z}.$ Skup cijelih brojeva \v cine, osim prirodnih brojeva, jo\v s i broj $0$ i negativni cijeli brojevi $\ldots,-3,-2,-1,$ pa vrijedi $\mathbb{Z}=\{\ldots,-2,-1,0,1,2,\ldots\}$ i $\mathbb{N}\subset\mathbb{Z}.$ Pri formiranju skupa cijelih brojeva vo\dj eno je ra\v cuna da se osobine skupa prirodnih brojeva prenesu i na skup cijelih brojeva, preciznije skup cijelih brojeva je zatvoren u odnosu na operacije sabiranja i mno\v zenja, ali zatvoren je i u odnosu na operaciju oduzimanja.

Za operacije sabiranje i mno\v zenja vrijede osobine
\begin{enumerate}
  \item Zatvorenost: $a+b$ i $a\cdot b$ su cijeli brojevi, za sve cijele brojeve $a$ i $b.$
  \item Komutativni zakon: $ a+b=b+a$  i $a\cdot b=b\cdot a,$ za sve cijele brojeve $a$ i $b.$
  \item Asocijativni zakon: $(a+b)+c=a+(b+c)$ i $(a\cdot b)\cdot c=a\cdot(b\cdot c),$ za sve cijele brojeve $a,\:b$ i $c.$
  \item Distributivni zakon: $(a+b)\cdot c=a\cdot c+b\cdot c,$ za sve cijele brojeve $a,$  $b$ i $c.$
  \item Neutralni elementi: $0$ za sabiranje $a+0=a$ i $1$ za mno\v zenje   $a\cdot 1=a,$ za sve cijele brojeve $a.$
  \item Suprotni element: za svaki cijeli broj $a$ postoji broj $x$ takav da je $a+x=0,$ tj. postoji u skupu cijelih brojeva
        rje\v senje jedna\v cine   $a+x=0.$ Broj $x$ zovemo suprotni element (ili inverzni element) broja $a$ u odnosu na sabiranje
        i ozna\v cavamo ga sa $-a.$
  \item Zakon kancelacije (skra\'civanja): Ako su $a,\:b$ i $c$ cijeli brojevi i $a\cdot c=b\cdot c,\: c\neq 0,$ tada je $a=b.$
\end{enumerate}

\paragraph{Ure\dj enje na skupu cijelih brojeva.} Ure\dj enje  na skupu cijelih brojeva $\mathbb{Z}$ mo\v ze se definisati koriste\' ci prirodne brojeve.
\begin{definition}
   Ako su $a$ i $b$ cijeli brojevi, tada je $a<b,$ ako je $b-a$  prirodan broj. Ako je $a<b$ mo\v zemo pisati i $b>a.$
\end{definition}

Vrijede fundamentalne osobine ure\dj enosti cijelih brojeva:
\begin{enumerate}[$\circ$]
     \item Zatvorenost za prirodne brojeve: $a+b$ i $a\cdot b$ su prirodni brojevi, kad god su $a$ i $b$ prirodni brojevi.
     \item Zakon trihonomije: Za svaki cijeli broj $a$ vrijedi jedno od pravila:
         \begin{inparaenum}
             \item $a>0;\quad$\item $a=0;\quad$ \item $ a<0.$
         \end{inparaenum}
\end{enumerate}

Ka\v zemo da je skup cijelih brojeva ure\dj en skup jer sadr\v zi podskup  koji je zatvoren u odnosu na sabiranje i mno\v zenje, te zakon trihotomije vrijedi za svaki cio broj.

\subsection{D\lowercase{jeljivost cijelih brojeva}}

Teorija brojeva je grana matematike koja se ponajprije bavi prou\v cavanjem osobina skupa prirodnih brojeva $\mathbb{N}$. Osim toga, prou\v cavaju se i osobine skupa cijelih brojeva $\mathbb{Z},$ te racionalnih brojeva $\mathbb{Q}$. Pojam djeljivosti je jedan od najjednostavnijih, ali i najva\v znijih pojmova u teoriji brojeva. Kada cio broj podijelimo drugim cijelim brojem razli\v citim od nule, rezultat mo\v ze a i ne mora biti cio broj. Na primjer, $32:4=8$ je cio broj, dok $12:5=2.4$ nije cio broj. Ovi primjeri vode nas do sljede\' ce definicije.\\
\begin{definition}
	Ako su $a$ i $b$ cijeli brojevi i $a\neq 0,$ ka\v zemo da $a$ dijeli $b$ ili da je $b$ djeljivo sa $a,$ ako postoji cio broj $c,$ takav da je $b=a\cdot c.$ Ako $a$ dijeli $b,$ ka\v zemo $a$ je djelilac ili faktor od broja $b.$
\end{definition}

Ako $a$ dijeli $b$ pi\v semo $a\mid b,$ a ako ne dijeli onda $a\nmid b.$\\
\begin{example}
	Sljede\' ci primjeri ilustruju koncept djeljivosti cijelih brojeva: $7\mid 35,\:13\mid 65,\:11\mid 55,\:-5\mid 45,\:14\mid 0,\:3\nmid 10,\:4\nmid 17.$
\end{example}

\begin{example}
	Djelioci broja $6$ su $\pm 1,\:\pm 2,\:\pm 3,\:\pm 6.$ Djelioci broja $19$ su $\pm 1$ i $\pm 19.$ Djelioci broja $100$ su\\ $\pm 1,\,\pm 2,\, \pm 4,\,\pm 5,\,\pm 10,\,\pm 20,\,\pm 25,\, \pm 50$ i $\pm 100.$
\end{example}

Slijede neke osnovne osobine djeljivosti.\\

\begin{theorem}\label{djeljivost1}
	Ako su $a,\,b$ i $c$ cijeli brojevi takvi da je $a\mid b$ i $b\mid c,$ tada je $a\mid c.$
\end{theorem}

\begin{example}
	Ako je $7\mid 35$ i $35\mid 105,$ tada po Teoremi  \ref{djeljivost1} slijedi da je $7\mid 105.$
\end{example}

\begin{theorem}\label{djeljivost2}
	Ako je su $a,\,b,\,c,\,m$ i $n$ cijeli brojevi takvi da je $c\mid a,\,c\mid b,$ tada je $c\mid (ma+nb).$
\end{theorem}
\newpage
\begin{example}
	Kako je $3\mid 12$ i $3\mid 21,$ tada po Teoremi \ref{djeljivost2}, $3$ dijeli i $39$ jer je
	\[5\cdot 12-1\cdot 21=39.\]
\end{example}

\begin{theorem}[Teorema o dijeljenju sa ostatkom]
	Ako su $a$ i $b,\,b>0$ cijeli brojevi, tada postoje jedinstveni   cijeli brojevi $q$ i $r$ takvi da je
	\[a=b\cdot q+r\text{  i  } 0\leqslant r<b.\]
\end{theorem}

U prethodnoj teoremi broj $a$ je dijeljenik, $b$ je djelilac, $q$ je koli\v cnik i  $r$ je ostatak. Ova teorema se u dijelu literature naziva i algoritam dijeljenja. Iz $a=bq+r$ mo\v zemo primijetiti da je $a$ djeljivo sa $b$ jedino ako je $r=0.$\\

\begin{example}
	Ako je $a=133$ i $b=21$, tada je $q=6$  i $r=7,$ jer je $133=21\cdot 6+7.$
\end{example}

Za dati pozitivan cio broj $d,$ mo\v zemo cijele brojeve svrstati u razli\v cite podskupove zavisno od toga koliki je ostatak pri dijeljenju nekog cjelog broja sa $d$. Ako je broj $d=2$ onda su mogu\' ci ostaci $0$ i $1.$ Ovo nas vodi do sljede\' ce definicije.\\

\begin{definition}
	Ako je ostatak dijeljenja cijelog pozitivnog broja $n$ sa $2$ jednak $0,$ tada je $n=2k$ za neki pozitivan cio broj $k$ i ka\v zemo da je $n$ paran broj. U suprotnom, ako je ostatak dijeljenja $n$ sa $2$ jednak $1,$ onda je $n=2k+1,$ za neki pozitivan cio broj $k$ i ka\v zemo da je $n$ neparan broj.
\end{definition}

\paragraph{Prosti i slo\v zeni brojevi.} Broj $1$ ima samo jednog pozitivnog djelioca i to je $1.$ Svaki drugi pozitivan cio broj ima najmanje dva pozitivna djelioca i to $1$ i samog sebe. Pozitivni cijeli brojevi koji imaju samo dva djelioca jako su bitni.\\

\begin{definition}[Prosti brojevi]
	Prost broj je pozitivan cio broj koji je djeljiv jedino sa $1$  i samim sobom.
\end{definition}

\begin{example}
	Prosti brojevi su $2,\,3,\,5,\,7,\,11,\, 13, \ldots$
\end{example}

\newpage

\begin{definition}[Slo\v zeni brojevi]
	Pozitivni cijeli brojevi ve\' ci od $1$ koji nisu prosti, nazivaju se slo\v zeni brojevi.
\end{definition}

\begin{example}
	Slo\v zeni brojevi su $12=2\cdot 2\cdot 3,\,15=3\cdot 5,\,60=2\cdot 2\cdot  3\cdot 5,\ldots$
\end{example}

Vrijedi teorema.
\begin{theorem}[Euklidov teorem]
	Prostih brojeva ima beskona\v cno mnogo.
\end{theorem}
\begin{proof}[Dokaz]
	Pretpostavimo da prostih brojeva ima kona\v cno mnogo i da su to $p_1,\ldots,p_{k-1},p_k,$\\$p_{k+1},\ldots,p_n.$ Konstrui\v simo sljede\' ci broj $q=p_1\cdot\ldots \cdot p_{k-1}\cdot p_k\cdot p_{k+1}\cdot\ldots\cdot p_n+1,$ ovaj broj je po na\v soj pretpostavci slo\v zen, jer su samo $p_1,\ldots,p_{k-1},p_k,p_{k+1},\ldots, p_n$ prosti brojevi.  Po\v sto je $q$ slo\v zen broj, mora postojati njegov djelioc razli\v cit od $1$ i njega samog. Podijelimo broj $q$ bilo kojim prostim brojem, na primjer sa $p_k,$ vrijedi
	\[\frac{q}{p_k}=\frac{p_1\cdot\ldots \cdot p_{k-1}\cdot p_k\cdot p_{k+1}\cdot\ldots\cdot p_n+1}{p_k}
	=p_1\cdot\ldots \cdot p_{k-1}\cdot p_{k+1}\cdot\ldots\cdot p_n+\frac{1}{p_k}.\] Dobijeni broj nije cio (racionalan broj bi dobili u slu\v caju dijeljenja sa bilo kojim od brojeva $p_1,\ldots,p_n$), pa $q$ nije djeljiv ni sa jednim prostim brojem, tj. $q$ je i sam prost broj. Proces konstruisanja brojeva na na\v cin kako je konstruisan broj $q$ mo\v zemo nastaviti proizvoljan broj puta. Dakle pretpostavka da prostih brojeva ima kona\v cno mnogo je pogre\v sna.
\end{proof}

Fundametalna teorema aritmetike je va\v zan rezultat, koji ka\v ze da su prosti brojevi fundamentalni blokovi za izgradnju slo\v zenih brojeva.
\begin{theorem}[Fundamentalni teorem aritmetike]\label{fundamentalnith}
	Svaki slo\v zeni broj mo\v ze biti predstavljen na jedinstven na\v cin kao proizvod prostih brojeva, pri \v cemu su prosti brojevi u proizvodu zapisani u neopadaju\' cem redoslijedu.
\end{theorem}
\begin{example}
	Brojeve $12,\,30,\,33,\,48,$ zapisujemo na sljede\' ci na\v cin:
	\[12=2^2\cdot 3,\:30=2\cdot3\cdot 5,\:33=3\cdot 11,\:48=2^4\cdot 3.\]
\end{example}

\paragraph{Najve\' ci zajedni\v cki djelilac i najmanji zajedni\v cki  sadr\v zilac.} Neka su $a$ i $b$ cijeli brojevi, razli\v citi od nule, tada je skup zajedni\v ckih djelilaca brojeva $a$ i $b$ neki kona\v can skup cijelih brojeva koji uvijek sadr\v zi $-1$ i $1.$ Zanima nas koji je najve\' ci djelilac ova dva broja $a$ i $b.$

\begin{definition}[Najve\' ci zajedni\v cki djelilac]
	Najve\' ci zajedni\v cki djelilac cijelih brojeva $a$ i $b,$ razli\v citih od nule, je najve\' ci djelilac koji dijeli istovremeno i $a$ i $b.$
\end{definition}
Najve\' ci zajedni\v cki djelilac brojeva $a$ i $b$ ozna\v cavamo sa $\nzd(a,b)$ (ili samo $(a,b)$). Broj $0$ ima beskona\v cno mnogo djelilaca.

\begin{example}
	Zajedni\v cki djelioci brojeva $24$ i $84$ su: $\pm 1,\,\pm 2,\,\pm 3,\, \pm 4,\,\pm 6$ i $\pm 12.$ Dakle $\nzd(24,84)=12.$ Vrijedi i $\nzd(15,81)=3,\,\nzd(100,5)=5,\,\nzd(17,25)=1,\,\nzd(0,44)=44,\,\nzd(-15,-6)=3$ i $\nzd(-17,289)=17.$
\end{example}
Ako je najve\' ci zajedni\v cki djelilac dva cijela broja $a$ i $b$ jednak $1$ tada taj par ima posebno ime.
\begin{definition}[Relativno prosti brojevi]
	Cijeli brojevi $a$ i $b$ su relativno prosti ako je najve\' ci zajedni\v cki djelilac ova dva broja $1,$ tj. $\nzd(a,b)=1.$
\end{definition}

\begin{example}
	Kako je $\nzd(25,42)=1,$ to su $25$ i $42$ relativno prosti.
\end{example}

\paragraph{Euklidov algoritam.}
Euklidov algoritam predstavlja sistematski metod za ra\v cunanje najve\' ceg zajedni\v ckog djelioca dva cijela broja $a$ i $b$ bez direktnog ra\v cunanja zajedni\v ckih djelilaca. Ime je dobio po gr\v ckom matemati\v caru Euklidu\footnote{Euclid (ro\dj en sredinom 4. vijeka pne, umro sredinom 3. vijeka pne), napisao Elemente - djelo od 13 knjiga, jedno od najuticajnijih djela u razvoju matematike}.\\

\begin{theorem}[Euklidov algoritam]
	Neka su $r_0=a$ i $r_1=b$ cijeli brojevi takvi da je $a\geqslant b>0.$ Ako algoritam dijeljenja sukscesivno primijenimo da dobijemo $r_j=r_{j+1}q_{j+1}+r_{j+2}$ gdje je $0<r_{j+2}<r_{r+1}$ za $j=0,1,\ldots,n-2$ i $r_{n+1}=0$ tada je $\nzd(a,b)=r_n,$ tj. zadnji  ostatak razli\v cit od nule je najve\' ci zajedni\v cki djelilac brojeva $a$ i $b.$
\end{theorem}

\begin{example}
	Odrediti \\
	\begin{inparaenum}
		\item $\nzd(252,198).\:$\\
		\item $\nzd(102,222).\:$
	\end{inparaenum}\\\\
	Rje\v senje:
	\begin{enumerate}
		\item    \begin{align*}
		252=&1\cdot 198+54\\
		198=&3\cdot 54+36\\
		54=&1\cdot 36+{\color{magenta}18}\\
		36=&2\cdot 18,
		\end{align*}
		pa je $\nzd(252,198)=18.$
		\item
		\begin{align*}
		222=&2\cdot 102+18\\
		102=& 5\cdot 18+12\\
		18=&1\cdot 12+{\color{magenta}6}\\
		12=&2\cdot 6,
		\end{align*}
		pa je $\nzd(222,102)=6.$
	\end{enumerate}
\end{example}

\paragraph{Najmanji zajedni\v cki sadr\v zilac.}

Faktorizacija na proste brojeve mo\v ze biti iskori\v stena za ra\v cunanje najmanjeg cijelog pozitivnog broja koji sadr\v zi dva pozitivna cijela broja. Na ovaj problem nailazimo, na primjer, prilikom sabiranja razlomaka.
\begin{definition}[Najmanji zajedni\v cki sadr\v zilac]
	Najmanji zajedni\v cki sadr\v zilac pozitivnih cijelih brojeva $a$ i $b$ je najmanji pozitivan broj koji je djeljiv i sa $a$ i sa $b.$
\end{definition}
Najmanji zajedni\v cki sadr\v zilac brojeva $a$ i $b$ ozna\v cavamo sa $\nzs(a,b)$ (ili sa $[a,b]$).

\begin{example}
	Koriste\' ci faktorizaciju na proste brojeve izra\v cunati $\nzd(90,24)$ i $ \nzs(90,24).$\\\\
	Rje\v senje.\\\\
	Izvr\v simo prvo faktorizaciju brojeva $90$ i $24.$ Dijelimo brojeve $90$ i $ 24$ prostim brojevima po\v cev od $2,\,3,\,\ldots,$ ako su brojevi $90$ i $24$ djeljivi pi\v semo rezultat. Vrijedi
	\begin{align*}
	90&=45\cdot 2 & 24&=12\cdot 2\\
	45&=15\cdot 3 & 12&=6\cdot 2\\
	15&=5\cdot 3  & 6&=3\cdot 2\\
	5&=1\cdot 5   & 3&=1\cdot 3,
	\end{align*}
	pa je
	\begin{align*}
	90&=2\cdot 3\cdot 3\cdot 5\\
	24&=2\cdot 2\cdot 2\cdot 3.
	\end{align*}
	
	Iz oba razvoja biramo zajedni\v cke faktore, to su $2$ i $3.$ Njihov proizvod je najve\' ci zajedni\v cki djelilac brojeva $90$ i $24,$ tj. $\nzd(90,24)=2\cdot 3=6.$
	
	Najmanji zajedni\v cki sadr\v zilac dobijamo tako \v sto uzimamo sve faktore koji se pojavljuju u oba razvoja, ali bez ponavljanja. Iz drugog razvoja uzimamo proizvod $2\cdot 2\cdot 2,$ broj $2$ iz razvoja broja $90$ ne\' cemo uzimati jer ve\' c imamo $2$ u proizvodu $2\cdot 2\cdot 2,$ zatim uzimamo proizvod $3\cdot 3$ iz prvog razvoja, opet ne\' cemo uzeti $3$ iz drugog razvoja jer smo ve\' c uzeli  $3\cdot 3$ i na kraju uzimamo $5$ u proizvod, pa je $\nzs(90,24)=2\cdot 2\cdot 2\cdot 3\cdot 3\cdot 5=360.$
\end{example}



\section[S\lowercase{kup racionalnih brojeva}]{Skup racionalnih brojeva} Skup racionalnih brojeva ozna\v cavamo sa $\mathbb{Q}$ i defini\v semo ga sa $\mathbb{Q}=\left\{ \frac{m}{n}:m\in\mathbb{Z},\,n\in\mathbb{N}\right\}.$ Cijele brojeve mo\v zemo predstaviti u obliku $m=\frac{m}{1},$ pa vrijedi $\mathbb{N}\subset\mathbb{Z}\subset\mathbb{Q}.$ Skup racionalnih brojeva je zatvoren za sabiranje, mno\v zenje, oduzimanje i dijeljenje (osim nulom), te u skupu racionalnih brojeva mo\v zemo rje\v savati i jedna\v cine oblika $a\cdot x=b,\,a\neq 0$. Znamo da je kod racionalnog broja $\frac{m}{n},$ $m$ brojnik, a $n$ je  nazivnik ili imenilac, te da se dva racionalna broja sabiraju, oduzimaju, mno\v ze  i dijele po pravilima:
\begin{gather*}
 \frac{a_1}{b_1}\pm\frac{a_2}{b_2}=\frac{a_1b_2\pm a_2b_1}{b_1b_2},\:\frac{a_1}{b_1}\cdot\frac{a_2}{b_2}=\frac{a_1a_2}{b_1b_2},b_1,b_2\neq 0,
\end{gather*}
\begin{gather*}
 \frac{a_1}{b_1}:\frac{a_2}{b_2}=\frac{a_1b_2}{b_1a_2},\:b_1,b_2,a_2\neq 0.
\end{gather*}
Me\dj utim, ni u skupu racionalnih brojeva $\mathbb{Q}$ ne mo\v zemo rije\v siti neke relativno jednostavne jedna\v cine, kao na primjer $x^2=2.$ Rje\v senje ove jedna\v cine je $x_{1/2}=\pm\sqrt{2}.$ Poka\v zimo da $\sqrt{2}$ nije racionalan broj. Pretpostavimo suprotno, tj. da postoje brojevi $\m\in\mathbb{Z}$ i $n\in\mathbb{N},$ za koje vrijedi $\sqrt{2}=\frac{m}{n}$ i da su $m$ i $n$ relativno prosti ($\nzd(m,n)=1$). Kvadriranjem dobijamo \[2=\frac{m^2}{n^2}\Leftrightarrow 2n^2=m^2,\] sada zaklju\v cujemo da su $m^2$ i $m$ parni brojevi, pa je $m=2k,\,k\in\mathbb{Z}.$ Dalje je
\[2n^2=m^2\Leftrightarrow 2n^2=(2k)^2\Leftrightarrow n^2=2k^2,\] \v sto zna\v ci da je i $n^2,$ odnosno $n,$ paran broj. Sada smo dobili da su i $m$ i $n$ parni brojevi, \v sto je u suprotnosti sa na\v som pretpostavkom da su $m$ i $n$ uzajamno prosti brojevi. Mo\v zemo zaklju\v citi da je gre\v ska u pretpostavci da je $\sqrt{2}$ racionalan broj, dakle $\sqrt{2}$ nije racionalan broj.

\section[S\lowercase{kup iracionalnih brojeva}]{Skup iracionalnih brojeva} Prethodni primjer doveo je do formiranja skupa iracionalnih brojeva $\mathbb{I}.$ Ovaj skup \v cine iracionalni brojevi, tj. brojevi koje ne mo\v zemo napisati u obliku $\frac{m}{n},\,m\in\mathbb{Z},\,n\in\mathbb{N}.$ Poznato je da ovi brojevi imaju beskona\v can decimalni zapis, \v cije se cifre ne ponavljaju periodi\v cki (kod racionalnih brojeva u decimalnom zapisu cifre se ponavljaju periodi\v cki, zato je $\mathbb{Q}\cap\mathbb{I}=\emptyset$). U iracionalne brojeve, osim $\sqrt{2}$ spadaju  $\sqrt{3},\sqrt{5},\sqrt[3]{2}, \pi=3.14159\ldots,\,e=2.718281\ldots$ i mnogi drugi. Iracionalnih brojeva
ima zapravo beskona\v cno mnogo. Mo\v zemo ih podijeliti u dvije grupe: algebarski i transcendentni iracionalni brojevi. Algebarski brojevi su brojevi koji se mogu dobiti kao rje\v senja algebarskih
jedna\v cina $a_n x^n+a_{n-1} x^{n-1}+ \ldots +a_1x+a_0=0,$ gdje je $n$ prirodan broj, a koeficijenti $a_0, a_1, \ldots, a_n$ su cijeli brojevi. Pojam algebarskog broja je prirodna generalizacija racionalnog broja. Naime, svaki racionalan broj je ujedno i algebarski broj, ali ima i brojeva koji su algebarski a nisu racionalni. Na primjer $\sqrt{2}$ nije racionalan broj ali jeste algebarski broj, jer se mo\v ze dobiti kao rje\v senje jedna\v cine $x^2-2=0.$
Joseph Liouville\footnote{Joseph Liouville (24.mart 1809.--8.septembar 1882. godine) bio je francuski matemati\v car i in\v zinjer. }  je 1844. godine dokazao da postoje brojevi koji nisu algebarski. Oni su nazvani transcedentni brojevi, a takvi su npr. brojevi $e$ i $\pi.$ Dokaz da je broj $2^{\sqrt{2}}$ transcedentan dali su nezavisno jedan od drugog, C.L.Siegel\footnote{Carl Ludwig Siegel (31.decembar 1896.--4.april 1981. godine) bio je njema\v cki matemati\v car koji se bavio analiti\v ckom teorijom brojeva. Va\v zio je za jednog od najva\v znijih matemati\v cara 20.vijeka.} i A.Gelfond \footnote{Alexander Osipovich Gelfond (24.oktobar 1906.--7.novembar 1968. godine) bio je sovjetski matemati\v car. Objavio je va\v zne rezultate u nekoliko matemati\v ckih oblasti: teorija brojeva, analiti\v cke funkcije, integralne jedna\v cine i istorija matematike. Gelfondova teorema nosi ime po njemu. }. To je bio jedan od problema koje je David Hilbert\footnote{David Hilbert (23.januar 1862.--14.februar 1943. godine) bio je njema\v cki matemati\v car, jedan od najuticajnijih matemati\v cara u 19. i 20. vijeku.} predstavio u svom \v cuvenom govoru na Internacionalnom kongresu matemati\v cara u Parizu $1900.$ godine. Vrijedi i op\v stije:
svaki broj oblika $a^b,$ gdje je $a$ algebarski broj razli\v cit od $0$ i $1,$ a $b$ bilo koji iracionalan algebarski broj, je transcedentan.

\section[S\lowercase{kup realnih brojeva}]{Skup realnih brojeva} Racionalni i iracionalni brojevi sa\v cinjavaju zajedno skup realnih brojeva $\mathbb{R}$ i zajedno ih nazivamo realni brojevi, tj. vrijedi $\mathbb{R}=\mathbb{Q}\cup\mathbb{I}.$

U skupu realnih brojeva definisane su operacije sabiranja $+$ i operacija mno\v zenja $\cdot.$ Za svaka dva realna broja $x$ i $y,$ jedozna\v cno su odre\dj eni realni brojevi $x+y$ i $x\cdot y$. Ove operacije imaju sljede\' ce osobine:

\begin{enumerate}[({A}$1$)]
  \item Za svako $x,y,z\in\mathbb{R}$ vrijedi $(x+y)+z=x+(y+z),$ (asocijativni zakon). \label{aksiom1}
  \item Postoji samo jedan realan broj $0\in\mathbb{R}$ takav da je svaki $x\in\mathbb{R}$ vrijedi $x+0=0+x=x,$ ($0$ je neutralni element za sabiranje).
  \item Za svaki $x\in\mathbb{R}$ postoji samo jedan $-x\in\mathbb{R}$ takav da je $x+(-x)=(-x)+x=0,$ ($-x$ je inverzni element elementa $x$ u odnosu na sabiranje).
  \item Za sve $x,y\in\mathbb{R}$ vrijedi $x+y=y+x$ (komutativni zakon za sabiranje).
  \item Za sve $x,y,z\in\mathbb{R}$ vrijedi $(xy)z=x(yz)$ (asocijativni zakon za mno\v zenje).
  \item Postoji samo jedan realan broj $1\in\mathbb{R}$ takav da za svaki $x\in\mathbb{R}$ vrijedi $x\cdot 1=1\cdot x=x$ ($1$ je neutralni element za mno\v zenje).
  \item Za svaki $x\in\mathbb{R},\,x\neq 0,$ postoji jedinstven element $x^{-1}=\frac{1}{x}\in\mathbb{R}$ takav da je $x\cdot x^{-1}=x^{-1}\cdot x=1$ ($x^{-1}$ je inverzni  element  elementa $x$ za mno\v zenje).
  \item Za sve $x,y\in\mathbb{R}$ vrijedi $xy=yx$ (komutativni zakon za mno\v zenje).
  \item Za sve $x,y,z\in\mathbb{R}$ vrijedi $x(y+z)=xy+xz$ (distributivni zakon mno\v zenja prema sabiranju).  \label{aksiom9}
\end{enumerate}

\begin{remark}
   Svaku ure\dj enu trojku $(S,+,\cdot)$ koju \v cine proizvoljan neprazan skup, te binarne ope-\\racije $+$ i $\cdot$ za koje va\v ze osobine (aksiomi)
   \hyperref[aksiom1]{(A1)}--\hyperref[aksiom9]{(A9)} nazivamo polje. Dakle polja su $(\mathbb{R},+,\cdot)$ i $(\mathbb{Q},+,\cdot),$ dok $(\mathbb{N},+,\cdot)$ i $(\mathbb{Z},+,\cdot)$
    nisu polja.
\end{remark}

\begin{figure}[!h]\centering
  \includegraphics[scale=1]{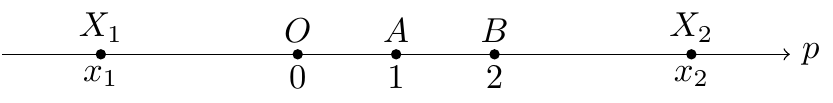}
  \caption{Brojna prava}
  \label{slika1}
\end{figure}

\paragraph{Brojna osa.} Neka je data prava $p,$ vidi Sliku \ref{slika1}. Na pravoj $p$ uo\v cimo bilo koju ta\v cku, na primjer O.  Na ovaj na\v cin dijelimo pravu $p$ na tri dijela, prvi dio su ta\v cke koje su lijevo od ta\v cke O, samu ta\v cku O i ta\v cke koje su desno od ta\v cke O. Sada ho\' cemo da svakoj ta\v cki prave $p$ pridru\v zimo ta\v cno jedan realan broj. Samoj ta\v cki O mo\v zemo pridru\v ziti nulu $0,$  ta\v cki $A$ broj $1,$ ta\v cki $B$ broj $2$ itd. Ta\v ckama koje su desno od ta\v cke O pridru\v  zujemo pozitivne realne brojeve. Na primjer, ta\v cki $X_2$ dodjeljujemo broj $x_2=\overline{OX_2},$ tj. du\v zinu du\v zi $OX_2,$ dok ta\v ckama lijevo od ta\v cke O dodjeljujemo negativne broje, na primjer  $x_1=-\overline{OX_1}.$ Na ovaj na\v cin uspostavljena je bijekcija izme\dj u ta\v caka brojne prave i realnih brojeva (svakoj ta\v cki brojne prave odgovara samo jedan realan broj i obrnuto). Pravu $p$ zovemo brojna osa (ili prava).

\paragraph{Ure\dj enje na skupu realnih brojeva.} Vidjeli smo da svakom realnom broju odgovara jedna ta\v cka na brojnoj osi. Ako se ta\v cka koja odgovara broju $x$ nalazi lijevo od ta\v cke koja odgovara broju $y,$ onda je $x<y$ ili $y>x,$ tj. $x$ je manje od $y$ ili $y$ je ve\' ce od $x$ (ili $y$ je desno od $x$ pa je $y$ ve\' ce od $x$).  Dakle, bilo koja dva realna broja $x$ i $y$ mo\v zemo uporediti. Mogu\' ca su tri slu\v caja $x<y,$ $x=y$ i $x>y.$ Ako koristimo relaciju ure\dj enja (ili poretka)  $\leqslant,$ onda $x\leqslant y$ zna\v ci ili je $x<y$ ili $x=y.$ Relacija ure\dj enja ima osobine:
\begin{enumerate}[({A}$1$)]\setcounter{enumi}{9}
\item Za bilo koja dva realna broja $x,y\in\mathbb{R}$ vrijedi $x\leqslant y$ ili $y\leqslant x.$
\item ($\forall x \in \mathbb{R}$) $x\leq x$ (refleksivnost).
\item $x\leqslant y$ i $y\leqslant x$ ako i samo ako je $x=y$ (antisimetri\v cnost).
\item Ako je $x\leqslant y$ i $y\leqslant z$ onda je $x\leqslant z$ (tranzitivnost).
\end{enumerate}
Relacija ure\dj enja $\leqslant$ je kompatibilna sa sabiranjem i mno\v zenjem  i vrijedi
\begin{enumerate}[({A}$1$)]\setcounter{enumi}{12}
   \item Ako je $x\leqslant y$ onda za svaki $z\in\mathbb{R}$ vrijedi $x+z\leqslant y+z.$
   \item Iz $(0\leqslant x \wedge 0\leqslant y)\Rightarrow 0\leqslant xy.$\label{aksiom14}
\end{enumerate}
\begin{remark}
 Svaku ure\dj enu trojku $(S,+,\cdot)$ koju \v cine neprazan skup $S$ i binarne operacije $+$ i $\cdot$ za koje va\v ze osobine \hyperref[aksiom1]{(A1)}--\hyperref[aksiom14]{(A14)} zovemo ure\dj eno polje. Ure\dj eno polje je i $(\mathbb{R},+,\cdot)$ i $(\mathbb{Q},+,\cdot),$ dok  $(\mathbb{N},+,\cdot)$ i $(\mathbb{Z},+,\cdot)$ nisu polja.
\end{remark}

\paragraph{Intervali i segmenti.} Osim podskupova $\mathbb{N},\mathbb{Z},\mathbb{Q},\mathbb{I}$ skupa realnih brojeva $\mathbb{R}$ \v cesto se koriste i intervali odnosno segmenti. Uobi\v cajena naziv je otvoreni interval ili samo interval, dok zatvoreni interval nazivamo i segment. Otvoreni interval (ili samo interval) u oznaci $(a,b)$ je skup elemenata $x\in\mathbb{R}$ za koje vrijedi $a<x<b,$ gdje su $a$ i $b$ realni brojevi i $a<b,$ tj. \[(a,b)=\{x\in\mathbb{R}:a<x<b\}.\]
Zatvoreni interval (ili samo segment) u oznaci $[a,b],$ je skup elemenata $x\in\mathbb{R}$ za koje vrijedi $a\leqslant x\leqslant b,$ i ovdje su $a$ i $b$ realni brojevi i $a<b,$ tj. \[[a,b]=\{x\in\mathbb{R}:a\leqslant x\leqslant b\}.\]

Pored intervala i segmenta mogu se definisati i poluotvoreni intervali
\[(a,b]=\{x\in\mathbb{R}:a<x\leqslant b\},\quad [a,b)=\{x\in\mathbb{R}:a\leqslant x<b\},\]
i beskona\v cni intervali
\[(-\infty,a)=\{x\in\mathbb{R}:x<a\},\quad (-\infty,a]=\{x\in\mathbb{R}:x\leqslant a\},\]
\[(a,+\infty)=\{x\in\mathbb{R}:x>a\},\quad [a,+\infty)=\{x\in\mathbb{R}:x\geqslant a\}.\]


\paragraph{Supremum i infimum.}
\begin{definition}
  Ka\v zemo da je skup $S\subseteq\mathbb{R}$ odozgo ome\dj en ili ograni\v cen, ako postoji realan broj $M$ takav da je $x\leqslant M$ za svaki $x\in S.$ Svaki broj $M$ sa navedenim svojstvom nazivamo majoranta skupa $S.$ Ako skup $S$ nije odozgo ome\dj en, ka\v zemo da je  odozgo neome\dj en.

  Ka\v zemo da je skup $S\subseteq\mathbb{R}$ odozdo ome\dj en ili ograni\v cen, ako postoji realan broj $m$ takav da je $x\geqslant m$ za svaki $x\in S.$ Svaki broj $m$ sa navedenim svojstvom nazivamo minoranta skupa $S.$ Ako skup $S$ nije odozdo ome\dj en, ka\v zemo da je odozdo neome\dj en. Skup $S\subseteq \mathbb{R}$ je ome\dj en, ako je i odozgo i odozdo ome\dj en.  U protivnom se    se ka\v ze da je $S$ neome\dj en.
\end{definition}

\begin{example}
  Skup prirodnih brojeva ograni\v cen je odozdo, ali nije ograni\v cen odozgo, za minorantu dovoljno je uzeti bilo koji broj manji ili jednak 1, dok sa druge strane za bilo koji prirodan broj $n$ uvijek postoji prirodan broj ve\' ci od $n,$ na primjer $n+1,$ pa je $\mathbb{N}$ neome\dj en odozgo.
\end{example}

\begin{example}
  Svi skupovi $(a,b),\,[a,b],\,(a,b],[a,b)$ su ograni\v ceni, dok su $(-\infty,a),\,(-\infty,a]$ ograni\v ceni odozgo, a skupovi $(a,+\infty),\,[a,+\infty)$ su ograni\v ceni odozdo.
\end{example}

Svaki odozgo ograni\v cen skup $S$ ima vi\v se majoranti, isto tako ograni\v cen skup odozdo ima vi\v se minoranti. Vidjeli smo da su minorante skupa prirodnih brojeva svi brojevi manji ili jednaki $1$. Od svih majoranti najzanimljivija je najmanja majoranta, pa je pitanje kada ona postoji. Isto tako  postavlja se pitanje egzistencije najve\' ce minorante.

\begin{definition}
  Najmanju majorantu skupa $S$ nazivamo supremum i ozna\v cavamo sa $\sup S.$ Ako je $\sup S\in S,$ nazivamo ga maksimalnim elementom i ozna\v cavamo ga sa $\max S.$

  Najve\' cu minorantu skupa $S$ nazivamo infimum i ozna\v cavamo sa $\inf S.$ Ako je $\inf S\in S,$ nazivamo ga minimalnim elementom skupa $S$ i ozna\v cavamo ga sa $\min S.$
\end{definition}

U skupu $\mathbb{R}$ vrijedi va\v zna osobina:
\paragraph{Aksiom o supremumu.}
Svaki neprazan skup $S\subset\mathbb{R}$ ograni\v cen sa gornje strane ima $\sup S$ koji pripada $\mathbb{R}.$
Sli\v cna osobina vrijedi i za infimum, pa mo\v zemo objediniti:
\begin{enumerate}[({A}$1$)]\setcounter{enumi}{14}
 \item Svaki odozgo ograni\v cen skup $S\subset\mathbb{R}$ ima supremum, a svaki odozdo ograni\v cen skup $S\subset\mathbb{R}$ ima infimum.
\end{enumerate}

\section[A\lowercase{psolutna vrijednost realnog broja}]{Apsolutna vrijednost realnog broja}

\begin{definition}[Apsolutna vrijednost realnog broja]\index{apsolutna vrijednost realnog broja}
   Apsolutna vrijednost (ili modul ili norma) realnog broja $x,$ u oznaci $|x|$,  je pre-\\slikavanje $|\cdot|:\mathbb{R}\mapsto\mathbb{R}^{+}\cup\{0\}$, definisano sa
   \begin{equation*}
       |x|=\left\{
                    \begin{array}{cl}
                      x,&\text{ ako je } x>0,\\
                      0,& \text{ ako je } x=0,\\
                      -x,&\text{ ako je } x<0.
                    \end{array}
           \right.
   \end{equation*}
   \label{apsolutni}
\end{definition}

\begin{remark}
  Mo\v zemo koristiti i sljede\' ce definicije koje su ekvivalentne sa Definicijom {\rm\ref{apsolutni}}
     \begin{align*}
       |x|&\overset{\mathrm{def}}{=}\max\{x,-x\},\\
       |x|&\overset{\mathrm{def}}{=}\sqrt{x^2}.
     \end{align*}
\end{remark}

\begin{example}Vrijedi
  \begin{align*}
      |3|=3,\quad |-3|&=3,\quad |0|=0.
  \end{align*}
\end{example}
Vrijede sljede\' ce osobine apsolutne vrijednosti realnog broja

\begin{inparaenum}
   \item $|xy|=|x||y|,\:$
   \item $\left| \frac{x}{y}\right|=\frac{|x|}{|y|},\:$
   \item $|x+y|\leqslant |x|+|y|$ (nejednakost trougla),\,

   \item $|x-y|\leqslant |x|+|y|,\:$
   \item $|x|\leqslant a\Leftrightarrow -a\leqslant x\leqslant a,\, a>0,\:$
   \item $|x|\geq a\Leftrightarrow (x\leq -a) \vee (x\geq a),\:$
   \item $||x|-|y||\leqslant |x-y|.$
\end{inparaenum}

\begin{figure}[!h]\centering
 \includegraphics[scale=.8]{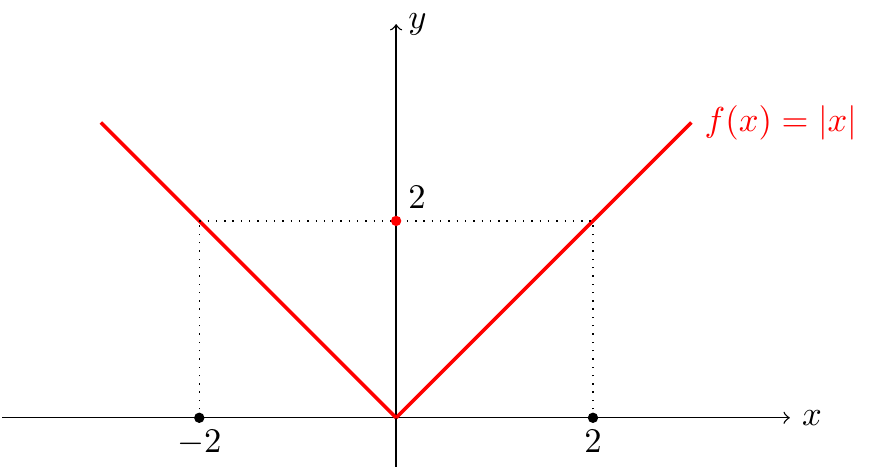}
 \caption{Grafik funkcije $f(x)=|x|$ i udaljenost ta\v caka $2$ i $-2$ od koordinatnog po\v cetka}
 \label{apsolutna1}
\end{figure}

\begin{remark}
   Geometrijski gledano apsolutna vrijednost realnog broja $x$ predstavlja udaljenost tog realnog broja od koordinatnog po\v cetka, vidjeti Sliku \ref{apsolutna1}.
\end{remark}

\begin{figure}[!h]
     \begin{subfigure}[b]{.5\textwidth}\centering
        \includegraphics[scale=.8]{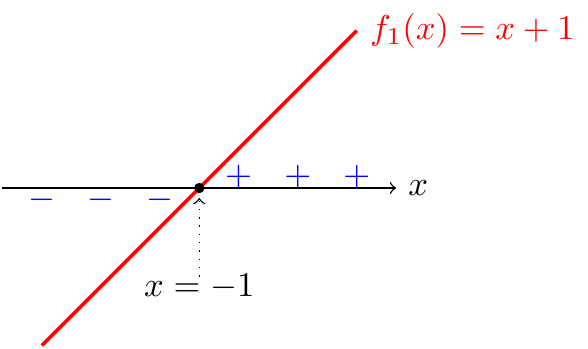}
         \caption{Grafik funkcije $f_1(x)=x+1$}
   \end{subfigure}
   \begin{subfigure}[b]{.5\textwidth}\centering
        \includegraphics[scale=.8]{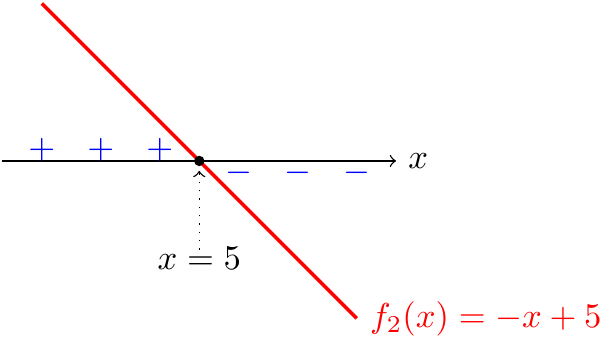}
         \caption{Grafik funkcije $f_2(x)=-x+5$}
     \end{subfigure}
     \caption{Grafici linearnih funkcija, koji \' ce poslu\v ziti za osloba\dj anje od zagrada za apsolutne vrijednosti}
     \label{slikalin}
   \end{figure}

\begin{example}
  Rije\v siti jedna\v cinu $3x-2|x+1|-|5-x|=3.$
\end{example}\ \\
\noindent Rje\v senje:\\\\
   Odredimo prvo znak izraza $x+1$ i $5-x,$ iskoristimo grafike linearnih funkcija $f_1(x)=x+1$ i $f_2(x)=-x+5,$ datih na Slici \ref{slikalin}.

Sada na osnovu grafika dobijamo Tabelu \ref{tabelalin1}, te na osnovu znaka funkcija $f_1(x)=x+1$ i $f_2(x)=-x+5,$ dijelimo domen na tri intervala.

\begin{table}[h]\centering
  \begin{tabular}{cccccccc}\\\\
  &$\infty$&&-1&&5&&$+\infty$\\\\\hline\\
  $|x+1|$&&$-(x+1)$&$0$&$x+1$&&$x+1$&\\\\
   $|5-x|$&&$5-x$&&$5-x$&$0$&$-(5-x)$&\\\\\hline\\
       &&I&&II&&III&
  \end{tabular}
  \caption{Znak funkcija $f_1(x)=x+1$ i $f_2(x)=-x+5$}
   \label{tabelalin1}
\end{table}

  \noindent Interval I,   $x\in (-\infty,-1]$ vrijedi
    \[ 3x-2[-(x+1)]-(5-x)=3\Leftrightarrow 3x+2x+2-5+x=3 \Leftrightarrow 6x=6
         \Leftrightarrow x=1\notin (-\infty,-1].\]
   Interval II, $x\in(-1,5]$ vrijedi
     \[3x-2(x+1)-5+x=3\Leftrightarrow 2x=10\Leftrightarrow x=5\in(-1,5].\]
  Interval III,  $x\in(5,\infty)$
    \[ 3x-2(x+1)-[-(5-x)]=3\Leftrightarrow 3=3,\]
    ovo zna\v ci da su rje\v senja jedna\v cine na intervalu III svi $x\in (5,\infty).$ Pa rje\v senje jedna\v cine unija svih dobijenih rje\v senja na sva tri intervala $x\in\{5\}\cup(5,\infty)=[5,\infty).$

\begin{example}
 Rije\v siti jedna\v cinu $\left|x^2+x-2\right|+\left|-x^2+x+2\right|=2.$
 \end{example}\ \\
\noindent Rje\v senje:\\\\
Da bi odredili znak kvadratnih trinoma unutar zagrada za apsolutne vrijednosti, po-\\smatrajmo kvadratne funkcije $f_1(x)=x^2+x-2$ i $f_2(x)=-x^2+x+2.$
   \begin{figure}[!h]
     \begin{subfigure}[b]{.5\textwidth}\centering
        \includegraphics[scale=.8]{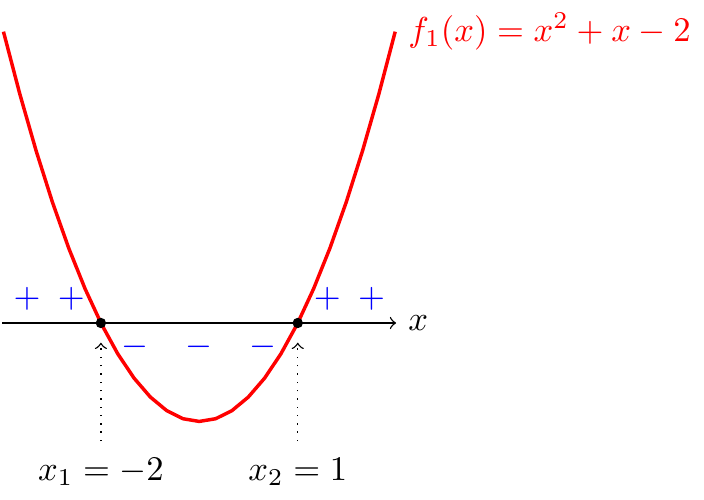}
         \caption{Grafik funkcije $f_1(x)=x^2+x-2$}
   \end{subfigure}
   \begin{subfigure}[b]{.5\textwidth}\centering
        \includegraphics[scale=.8]{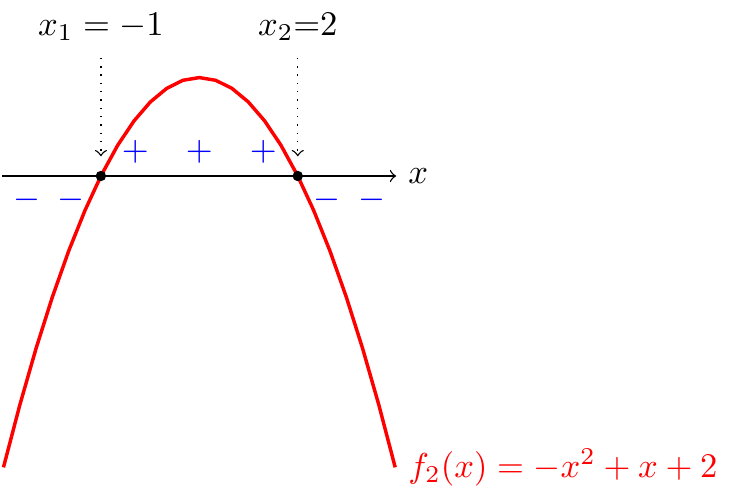}
         \caption{Grafik funkcije $f_2(x)=-x^2+x+2$}
     \end{subfigure}
     \caption{Grafici kvadratnih funkcija, koji \' ce poslu\v ziti za osloba\dj anje od zagrada za apsolutne vrijednosti}
     \label{slikakv}
   \end{figure}

Nule funkcija $f_1$ i $f_2$ odre\dj ujemo rje\v savaju\' ci odgovaraju\' ce kvadratne jedna\v cine po formuli $x_{1/2}=\frac{-b\pm\sqrt{b^2-4ac}}{2a},$ pa tako nalazimo 
  \[x^2+x-2=0\Leftrightarrow (x_1=-2\vee x_2=1);\quad -x^2+x+2=0\Leftrightarrow (x_1=-1\vee x_2=2).\]

Znak izraza pod apsolutnim zagradama odre\dj ujemo koriste\' ci prethodne grafike date na Slici \ref{slikakv}.
\begin{table}[!h]\centering\tiny
  \begin{tabular}{cccccccccccc}\\
  &$\infty$&&-2&&-1&&1&&2&&$+\infty$\\\\\hline\\
  $|x^2+x-2|$&&$x^2+x-2$&&$-(x^2+x-2)$&&$-(x^2+x-2)$&&$x^2+x-2$&&$x^2+x-2$ &\\\\
  $|-x^2+x+2|$&&$-(-x^2+x+2)$&&$-(-x^2+x+2)$&&$-x^2+x+2$&&$-x^2+x+2$ &&$-(-x^2+x+2)$ &\\\\\hline\\
       &&I&&II&&III&&IV&&V&
       \label{tabelakv1}
  \end{tabular}
    \caption{Znak funkcija $f_1(x)=x^2+x-2$ i $f_2(x)=-x^2+x+2$}
\end{table}

\noindent Interval I, $x\in (\infty, -2]$, vrijedi
\[x^2+x-2-(-x^2+x+2)=2\Leftrightarrow x^2=3\Leftrightarrow x_{1/2}=\pm\sqrt{3},\:x_{1/2}\notin(-\infty,-2].\]
\noindent Interval II, $x\in (-2,-1]$, vrijedi
\[-(x^2+x-2)-(-x^2+x+2)=2\Leftrightarrow x=-1,\:x=-1\in(-2,-1].\]
\noindent Interval III, $x\in (-1,1]$, vrijedi
\[-(x^2+x-2)-x^2+x+2=2\Leftrightarrow x^2=1\Leftrightarrow x_{1/2}=\pm 1,\:x_2=1\in(-1,1].\]
\noindent Interval IV, $x\in (1,2]$, vrijedi
\[x^2+x-2-x^2+x+2=2\Leftrightarrow x=1,\:x=1\notin(1,2].\]
\noindent Interval V, $x\in (2,\infty)$, vrijedi
\[x^2+x-2-(-x^2+x+2)=2\Leftrightarrow x^2=3\Leftrightarrow x_{1/2}=\pm\sqrt{3},\:x_{1/2}\notin(2,\infty).\]
\noindent Pa je kona\v cno rje\v senje jedna\v cine $x\in\{-1,1\}.$
\ \\

\begin{example}\label{primjernejed}
Rije\v siti nejedna\v cinu $|2-x|>|x+1|-3.$
\end{example}
\noindent Rje\v senje:\\\\

Osloba\dj anju od zagrada apsolutne vrijednosti pristupamo na isti na\v cin kao i slu\v caju jedna\v cina, koriste\' ci sljede\' ce grafike.

   \begin{figure}[!h]
     \begin{subfigure}[b]{.5\textwidth}\centering
        \includegraphics[scale=.8]{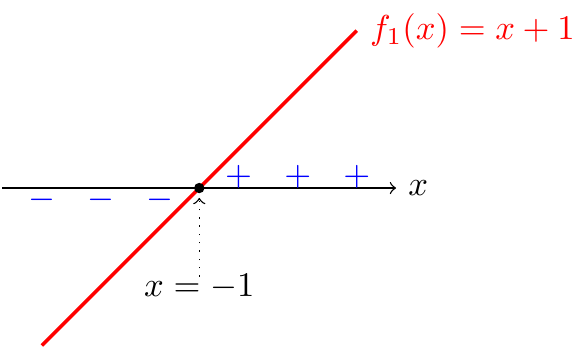}
         \caption{Grafik funkcije $f_1(x)=x+1$}
   \end{subfigure}
   \begin{subfigure}[b]{.5\textwidth}\centering
        \includegraphics[scale=.8]{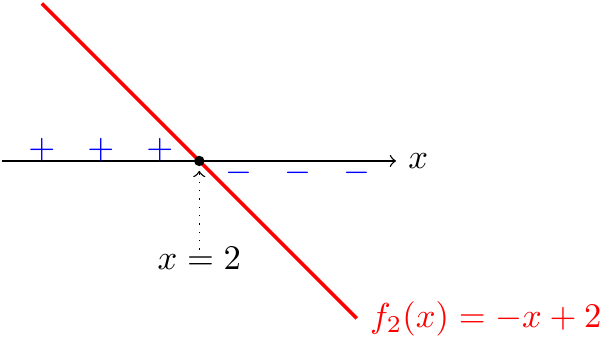}
         \caption{Grafik funkcije $f_2(x)=-x+2$}
     \end{subfigure}
     \caption{Grafici linearnih funkcija, koji \' ce poslu\v ziti za osloba\dj anje od apsolutnih vrije-\\dnosti}
     \label{slikalin2}
   \end{figure}
Podatke iz grafika prenesemo u tabelu.
\begin{table}[h]\centering
  \begin{tabular}{cccccccc}\\\\
  &$\infty$&&-1&&2&&$+\infty$\\\\\hline\\
  $|x+1|$&&$-(x+1)$&$0$&$x+1$&&$x+1$&\\\\
   $|2-x|$&&$2-x$&&$2-x$&$0$&$-(2-x)$&\\\\\hline\\
       &&I&&II&&III&
  \end{tabular}
  \label{tabelalin2}
  \caption{Znak funkcija $f_1(x)=x+1$ i $f_2(x)=-x+2$}
\end{table}

\noindent Interval I, $x\in (-\infty,-1],$  vrijedi
\[2-x>-(x+1)-3\Leftrightarrow 2-x>-x-1-3\Leftrightarrow 2>-4.\]
Nejednakost $2>-4$ je ta\v cna bez obzira na vrijednost $x,$ tj. $(\forall x\in\mathbb{R})\:2>-4,$ te je rje\v senje nejedna\v cine na Intervalu I,
\[x\in \mathbb{R}\cap(-\infty,-1]=(-\infty,-1].\]
\noindent Interval II, $x\in (-1,2],$  vrijedi
\[ 2-x>x+1-3\Leftrightarrow -2x>-4\Leftrightarrow x<2,\:x\in(-\infty,-2)\]
pa je rje\v senje nejedna\v cine
\[x\in(-\infty,2)\cap(-1,2]=(-1,2).\]
\noindent Interval III, $x\in (2,\infty),$  vrijedi
\[-( 2-x)>x+1-3\Leftrightarrow -2>-2.\]
Ovo je neta\v cna brojna nejednakost pa rje\v senje nejedna\v cine na $x\in (2,\infty),$  prazan skup, tj. $\emptyset$.

Rje\v senje nejedna\v cine je unija rje\v senja na svakom posmatranom intervalu
\[x\in(-\infty,-1]\cup (-1,2)=(-\infty,2).\]

\begin{example}
 Rije\v siti nejedna\v cinu $\left|\dfrac{2x-1}{x+1}\right|\leqslant 2.$
\end{example}

 \noindent Rje\v senje:\\\\
Razlomak koji se pojavljuje u datoj nejedna\v cini nije definisan za $x=-1,$ pa to treba obavezno uzeti u obzir kod rje\v savanja ove nejedna\v cine, drugim rije\v cima rje\v senje tra\v zimo na skupu $\mathbb{R}\setminus\{-1\}.$

 Ovu nejedna\v cinu mo\v zemo rije\v siti kao nejedna\v cinu u Primjeru  \ref{primjernejed}, me\dj utim mo\v zemo iskoristiti sljede\' cu osobinu
 \[|x|<a\Leftrightarrow -a<x<a\Leftrightarrow (-a<x\wedge x<a),\,a\in\mathbb{R}^{+}.\]
Koriste\' ci prethodnu osobinu realnih brojeva, vrijedi
\[\left|\frac{2x-1}{x+1}\right|\leqslant 2\Leftrightarrow -2\leqslant \frac{2x-1}{x+1}\leqslant 2\Leftrightarrow
 -2\leqslant \frac{2x-1}{x+1}\wedge \frac{2x-1}{x+1}\leqslant 2.\]
Sada je \[-2\leqslant \frac{2x-1}{x+1}\Leftrightarrow \frac{4x+1}{x+1}\geqslant 0,\]
zatim odredimo znak brojnika i nazivnika, koriste\'ci grafike linernih funkcija, te unesemo te podatke u sljede\'cu tabelu.


\begin{table}[h]\centering
 	\begin{tabular}{cc cc cc cc}\\
     		 & $-\infty$ & & $-1$ & & $-\frac{1}{4}$& & $+\infty$	\\ \hline
 $4x+1$      && $-$ && $-$ && $+$ \\
 $x+1$       && $-$ && $+$ && $+$ \\\hline
 $\dfrac{4x+1}{x+1}$ && $+$ && $-$ && $+$
 	\end{tabular}
\caption{ Znak razlomka  $\frac{4x+1}{x+1}$}
\label{tabelaLin3}
 \end{table}
\noindent Iz Tabele \ref{tabelaLin3} je  $ \frac{4x+1}{x+1}\geqslant 0$ za
\begin{equation}\label{rjesenje1}
	x\in(-\infty,-1)\cup [-\tfrac{1}{4},+\infty).
\end{equation}
Dalje je
\[ \frac{2x-1}{x+1}\leqslant 2\Leftrightarrow -\frac{3}{x+1}\leqslant 0,\]
odgovaraju\'ca tabela je

 \begin{table}[!h]\centering
	\begin{tabular}{cc cc cc }\\
		& $-\infty$ & & $-1$ & &   $+\infty$	\\ \hline
		$-3$      && $-$ && $-$   \\
		$x+1$       && $-$ && $+$   \\\hline
		$-\dfrac{3}{x+1}$ && $+$ && $-$
	\end{tabular}
	\caption{ Znak razlomka  $-\frac{3}{x+1}$}
\end{table}
\noindent Sada iz tabele dobijamo da je $-\frac{3}{x+1}\leqslant 0$ za
\begin{equation}\label{rjesenje2}
 x\in(-1,+\infty).
 \end{equation}
Rje\v senje polazne nejedna\v cine je presjek  dobijenih rje\v senja (rje\v senja datih u \eqref{rjesenje1} i \eqref{rjesenje2})),  u ovom slu\v caju
 \[x\in\left[-\tfrac{1}{4},\infty\right).\]


\begin{example}
  Rije\v siti nejedna\v cinu $\left|\frac{2x-1}{x+1}\right|\geqslant 2.$
\end{example}
\noindent Rje\v senje:\\\\
U slu\v caju nejedna\v cine $\left|\frac{2x-1}{x+1}\right|\geqslant 2,$ koristimo se osobinom
  \[|x|>a\Leftrightarrow (x<-a\vee x>a),\,a\in\mathbb{R}^{+},\]
pa je
\[\left|\frac{2x-1}{x+1}\right|\geqslant 2\Leftrightarrow \frac{2x-1}{x+1}\geqslant 2\vee \frac{2x-1}{x+1}\leqslant -2.\]
Rije\v se se obje nejedna\v cine, ali rje\v senje polazne nejedna\v cine je sada unija dobijenih rje\v senja, pa je na kraju rje\v senje
\[x\in(-\infty,-1)\cup\left(-1,-\tfrac{1}{4}\right].\]

\begin{remark}
I kod rje\v savanja ove nejedna\v cine vodimo ra\v cuna o definicionom podru\v cju, tj. rje\v senje tra\v zimo na skupu $\mathbb{R}\setminus\{-1\}.$ 	
	\end{remark}

\section[Z\lowercase{adaci}]{Zadaci}
\index{Zadaci za vje\v zbu!apsolutna vrijednost realnog broja}
\index{Zadaci za vje\v zbu!matemati\v cka indukcija}
\index{Zadaci za vje\v zbu!binomni obrazac}
   \begin{enumerate}
   \item Matemati\v ckom indukcijom pokazati ta\v cnost sljede\' cih jednakosti
    \begin{enumerate}
       \item $(\forall n\in\mathbb{N})\: 1+3+\ldots+(2n-1)=n^2;$
       \item $(\forall n\in\mathbb{N})\: 1\cdot 2+2\cdot 3+\ldots+n(n+1)=\frac{n(n+1)(n+2)}{3};$
       \item $(\forall n\in\mathbb{N})\: 1^2+2^2+\ldots+n^2=\frac{n(n+1)(2n+1)}{6};$
       \item $(\forall n\in\mathbb{N})\: \frac{1}{1\cdot 2}+\frac{1}{2\cdot 3}+\cdot\frac{1}{n(n+1)}=\frac{n}{n+1};$
       \item $(\forall n\in\mathbb{N})\: \frac{1}{1\cdot 3}+\frac{1}{3\cdot 5}+\ldots+\frac{1}{(2n-1)(2n+1)}=\frac{n}{2n+1}.$
    \end{enumerate}
    \item Matemati\v ckom indukcijom pokazati djeljivost \\
     \begin{inparaenum}
       \item  $(\forall n\in\mathbb{N})\: 6\mid 7^n-1;\:$
       \item  $(\forall n\in\mathbb{N})\: 3\mid n^3-n;\:$ \\
       \item  $(\forall n\in\mathbb{N})\: 6\mid n^3+5n;\:$
       \item  $(\forall n\in\mathbb{N})\: 17\mid 7\cdot 5^{2n-1} +2^{3n+1}.$
     \end{inparaenum}
   \item
   \begin{enumerate}
      \item Koriste\'ci binomni obrazac razviti binom $(3+2x)^5.$
      \item Odrediti tre\'ci \v clan u razvijenom obliku binoma $(2x\sqrt{x}-\sqrt[3]{x})^{8}.$
      \item U  razvijenom obliku binoma $\left(x^2+\frac{2}{x}\right)^n$ koeficijenti \v cetvrtog i trinaestog \v clana  su jednaki. Odrediti \v clan koji ne sadr\v zi $x.$
       \item \v Cetvrti \v clan u razvijenom obliku binoma $\left(10^{\log\sqrt{x}}+10^{-\frac{1}{\log x}}\right)^7 $  je $3500000.$\\ Odrediti $x.$
    \end{enumerate}


   	\item Primjenjuju\'ci algoritam dijeljenja, odrediti brojeve $q$ i $r$ tako da je $a=b\cdot q+r,$ ako je    \\
   	\begin{inparaenum}
   		\item $a=18,$ $b=4;$
   		\item $a=65,$ $b=9;$
   		\item $a=763,$ $b=24;$
   		\item $a=1000,$ $b=143.$
   	\end{inparaenum}
   	
   	\item   Odrediti najve\'ci zajedni\v cki djelilac $\nzd$ \\
   	\begin{inparaenum}
   		\item $\nzd(12,15);$
   		\item $\nzd(9,21);$
   		\item $\nzd(24, 60);$
   		\item $\nzd(52,76);$
   		
   		\item $\nzd(13,25);$
   		\item $\nzd(225, 102);$
   		\item $\nzd(524, 88);$
   		\item $\nzd(1331,1001).$          	
   	\end{inparaenum}
   	
   	\item Odrediti najmanji zajedni\v cki sadr\v zilac $\nzs$ \\
   	\begin{inparaenum}
   		\item  $\nzs(6,15 );$
   		\item  $\nzs(7,12 );$
   		\item  $\nzs(9,24 );$
   		\item  $\nzs(10, 48 );$
   		
   		\item  $\nzs(15,24 );$
   		\item  $\nzs(32, 92 );$
   		\item  $\nzs(54, 216 );$
   		\item  $\nzs(42, 200 ).$
   	\end{inparaenum}


   \item Rije\v siti jedna\v cine  \\
    \begin{inparaenum}
       \item $|2x-1|-2|1-x|=1;\:$
       \item $|x-3|+|1-4x|=2|x+2|;\:$
       \item $|2x-1|-3(3x+1)=|x-1|-2(x+3)-9;\:\vspace{.2cm}$\\
       \item $|2x+1|-3|x-3|=|x-1|+x+2;\:$
       \item $4(1-x)+|2x-1|=3x-|x-2|-1;\vspace{.2cm}$\\
       \item $|2x-1|-3(3x+1)=|x-1|-2(x+3)-9;\:\vspace{.2cm}$\\
       \item $|2x+1|-3|x-3|=|x-1|+x+2;\:$
       \item $4(1-x)+|2x-1|=3x-|x-2|-1.$
   \end{inparaenum}
   \item Rije\v siti jedna\v cine  \\
      \begin{inparaenum}
       \item $|x^2-4x|+3=x^2+|x-5|;\:$
       \item $|x-1|+|x^2+3x-4|=5 ;\:$
       \item $|x-5|+|x^2-2x-8|=7 ;\:\vspace{.2cm}$\\
       \item $2x-|5-|x-2||=1.$
    \end{inparaenum}

  \item Rije\v siti nejedna\v cine  \\
    \begin{inparaenum}
       \item $2-3|1-x|-2x\leqslant 1-4x-2|2x+3|;$
       \item $|x+1|+|3x-1|>2 ;\:\vspace{.2cm}$\\
       \item $|2x-1|-3(3x+1)<|x-1|-2(x+3)-9;\vspace{.2cm}$\\
       \item $|2x+1|-3|x-3|\leqslant|x-1|+x+2;\:$
       \item $4(1-x)+|2x-1|\geqslant3x-|x-2|-1;$\\

       \item $\left| \dfrac{2x-1}{x+1}\right|\leqslant 2;\:$
       \item $\left|\dfrac{x-1}{2x+1} \right|<1;$
       \item $\left|\dfrac{x+1}{2x-3} \right|\geqslant\dfrac{1}{2}.$
   \end{inparaenum}

    \item Rije\v siti nejedna\v cine  \\
     \begin{inparaenum}
         \item $|x^2-x|-|x|<1;$
         \item $|x^2-3x+2|-1>|x-3|\;\:$
         \item $ |x^2-7x+10|-|x-3|<6\;\:\vspace{.2cm}$\\
         \item $|x-1|+|x^2+3x-4|\geqslant5 ;\:$
         \item $|x^2-2x-3|+2-2x\geqslant|x-4|+x^2;\:$
         \item $|x^2+2x-3|<3x+3;\:\vspace{.2cm}$\\
         \item $|x^2-4x|+3\geqslant x^2+|x-5| ;\:$
         \item $\left| \dfrac{x^2+2x}{x^2-4x+3}\right|<1.$
   \end{inparaenum}
\end{enumerate}

   \chapter{Kompleksni brojevi}


\index{skup!kompleksnih brojeva}
\index{kompleksni brojevi}

U Poglavlju \ref{poglavljeRealni} smo vidjeli da se zbog potrebe za pro\v sirenjem pojma broja do\v slo do polja realnih brojeva $(\mathbb{R},\,+,\,\cdot).$ U polju realnih brojeva sve jedna\v cine oblika
\begin{align*}
  a+x&=b,\: c\cdot x=d,\\  x^n&=p,\:  x^{2n+1}=q,
\end{align*}
gdje su $a,b, d \in \mathbb{R},\,c\in \mathbb{R}\setminus\{0\},\, p\in \mathbb{R}^{+}\cup\{0\},q\in \mathbb{R},$ a $n\in\mathbb{N},$ mo\v zemo rije\v siti.\\ 
Me\dj utim, ponovo je potrebno izvr\v siti pro\v sirenje, ovog puta skupa realnih brojeva $\mathbb{R},$ jer npr. jednostavna jedna\v cina
\[x^2+1=0,\]
nema rje\v senje u skupu realnih brojeva $\mathbb{R}.$   

Defini\v simo na skupu
\[\mathbb{R}\times\mathbb{R}=\mathbb{R}^2=\{(a,b):a,b\in\mathbb{R}\},\]
operacije sabiranja $+$ i mno\v zenja $\cdot,$ ure\dj enih parova  na sljede\' ci na\v cin \index{kompleksni brojevi!sabiranje}\index{kompleksni brojevi!mno\v zenje}
\begin{empheq}[box=\mymath]{equation*}
 (\forall a,b,c,d\in\mathbb{R})\:(a,b)+(c,d)=(a+c, b+d),
\end{empheq}
 \begin{empheq}[box=\mymath]{equation*}
   (\forall a,b,c,d\in\mathbb{R})\:(a,b)\cdot( c,d)=(ac-bd, ad+bc).
 \end{empheq}

Ure\dj ena trojka $(\mathbb{R}^2,+,\cdot)$ je polje. Ovo polje naziva se polje kompleksnih brojeva i obilje\v zava se sa $\mathbb{C},$ dok njegove elemente nazivamo kompleksni brojevi. Pokazuje se da vrijede svi aksiomi polja.

\begin{remark}
 Za polje kompleksnih brojeva koristi se \v ce\v s\'ce oznaka $\mathbb{C},$ umjesto $(\mathbb{C}, +,\cdot).$ Ista oznaka $\mathbb{C}$ koristi se i za sam skup kompleksnih brojeva. Iz samog konteksta vidi se radi li se o skupu ili polju, a ako se ne vidi onda se naglasi o \v cemu se radi.
\end{remark}

\index{kompleksni brojevi!nula} \index{kompleksni brojevi!jedinica}\index{kompleksni brojevi!imaginarna jedinica}
Kompleksni broj $(0,0)$ zva\' cemo kompleksna nula, a kompleksni broj $(1,0)$ kompleksna jedinica i radi jednostavnijeg zapisa ozna\v cava\' cemo ih sa
\[(0,0)=0 \text{ i } (1,0)=1,\]
respektivno. Kompleksni broj $(0,1)$ ozna\v ci\' cemo sa $i,$ tj.
\[(0,1)=i,\]
te \' cemo ga zvati imaginarna jedinica.  Sada vrijedi,
\begin{align*}
z&=(x,y)=(x,0)+(0,y)=(x,0)+(y,0)\cdot(0,1)=x(1,0)+y(1,0)(0,1)\\
 &=x\cdot1+y\cdot1\cdot i=x+iy.
\end{align*}
Oblik \begin{empheq}[box=\mymath]{equation*}
z=x+iy,
\end{empheq} naziva se algebarski oblik kompleksnog broja. Pri tome je
\[i^2=(0,1)\cdot(0,1)=(-1,0)=-1.\]\index{kompleksni brojevi!algebarski oblik}
Vrijede sljede\'ce jednakosti
\begin{align*}
   i^1&=i\\
   i^2&=-1\\
   i^3&=i\cdot i^2=-i\\
   i^4&=-1\cdot(-1)=1\\
   i^5&=i\cdot i^4=i\\
   &\:\vdots
\end{align*}
te je
\begin{empheq}[box=\mymath]{equation*}
 i^{4n}=1, \: i^{4n+1}=i,\: i^{4n+2}=-1,\:  i^{4n+3}=-i,\,n\in\mathbb{N}.
\end{empheq}
S obzirom da smo definisali da je $(1,0)=1,$ mo\v zemo primijetiti da za svaki kompleksni broj oblika $(a,0)$ dobijamo
$(a,0)=a \cdot 1=a,$ \v sto je u stvari realan broj. Prema tome, jasno je da je zaista $\mathbb{R}\subset \mathbb{C}.$
Sa druge strane broj $(0,b)=b\cdot (0,1)=i\cdot b$ nazivamo imaginarni broj.
\section[A\lowercase{lgebarski oblik}]{Algebarski oblik kompleksnog broja}
Vratimo se ponovo algebarskom obliku kompleksnog broja $ z=x+iy.$
Realni broj $x$ je realni dio kompleksnog broja $z$ i pi\v semo
\begin{empheq}[box=\mymath]{equation*}
x=\re(z),
\end{empheq}
a realan broj $y$ je imaginarni dio kompleksnog broja $z,$ ovdje pi\v semo
\begin{empheq}[box=\mymath]{equation*}
y=\im(z).
\end{empheq}\index{kompleksni brojevi!realni dio}\index{kompleksni brojevi!imaginarni dio}
\newpage
Dva kompleksna broja $z_1$ i $z_2$ su jednaka ako je realni dio prvog kompleksnog broja jednak realnom dijelu drugog kompleksnog broja, te imaginarni dio prvog kompleksnog broja jednak imaginarnom dijelu drugog kompleksnog broja, drugim rije\v cima ako vrijedi:
\begin{equation*}
 \tcbhighmath[mojstil1]{ z_1=z_2\Leftrightarrow\re(z_1)=\re(z_2)\wedge \im(z_1)=\im(z_2).}
\label{kompl1}
\end{equation*}
Ako je sada imaginarni dio kompleksnog broja $z$ jednak nuli $\im(z)=0,$ tada je $z\in\mathbb{R},$ tj. taj kompleksni broj je \v cisto realan. Sa druge strane  ako je $\im(z)\neq 0\wedge \re(z)=0,$ tada je kompleksni broj \v cisto imaginaran.

\begin{remark}
Polje $(\mathbb{C},+,\cdot)$ nije ure\dj eno kao polje $(\mathbb{R},+,\cdot).$ Pretpostavimo suprotno da se mo\v ze definisati relacija poretka $\leqslant$ u $(\mathbb{C},+,\cdot).$ Tada bi za  $\forall z_1,z_2\in\mathbb{C}$ vrijedilo $z_1\leqslant z_2$ ili $z_2\leqslant z_1.$ Prethodno vrijedi i za $0$ i za $i,$ dakle imamo $0\leqslant i$ ili $i\leqslant 0.$  Iz $0\leqslant i\Rightarrow 0\leqslant i^2=-1.$ Kontradikcija. Sli\v cno, iz $i\leqslant 0$ slijedi $0\leqslant -i$ a sada je $0\leqslant i^2=-1,$ \v sto ponovo vodi ka kontradikciji.
\end{remark}

\paragraph{Konjugovano--kompleksni broj.} Za kompleksni broj $z=x+iy,$ kompleksni broj koji ima isti realni dio a imaginarni dio mu je suprotan broj od imaginarnog dijela $\im{z},$ ima posebnu ulogu. Taj broj $\overline{z}=x-iy$ je konjugovano--kompleksni broj kompleksnog broja $z.$
Vrijede sljede\'ce jednakosti
\begin{align*}
   \overline{\overline{z}}&=z,\\
   \frac{z+\overline{z}}{2}&=\frac{x+iy+x-iy}{2}=x=\re(z),\\
   \frac{z-\overline{z}}{2i}&=\frac{x+iy-(x-iy)}{2i}=y=\im(z),\\
   z\cdot\overline{z}&=(x+iy)\cdot(x-iy)=x^2+y^2=\left( \re(z)\right)^2+\left( \im(z)\right)^2.
\end{align*}
\index{kompleksni brojevi!konjugovano kompleksni broj}
\subsection[O\lowercase{peracije sa kompleksnim brojevima u algebarskom obliku}]{Operacije sa kompleksnim brojevima u algebarskom obliku}
Neka su kompleksni brojevi $z_1$ i $z_2$ dati u algebarskom obliku, tj. $z_1=x_1+iy_1,\:z_2=x_2+iy_2.$  Sabiranje, oduzimanje, mno\v zenje i dijeljenje se svodi na pravila data formulama

\begin{empheq}[box=\mymath]{equation*}
z_1+z_2=x_1+iy_1+(x_2+iy_2)=(x_1+ x_2)+i(y_1+ y_2),
\end{empheq}
\begin{empheq}[box=\mymath]{equation*}
 z_1-z_2=x_1+iy_1-(x_2+iy_2)=(x_1- x_2)+i(y_1- y_2),
\end{empheq}
\begin{empheq}[box=\mymath]{equation*}
z_1\cdot z_2=(x_1+iy_1)\cdot(x_2+iy_2)=(x_1x_2-y_1y_2)+i(x_1y_2+x_2y_1)
\end{empheq}
\begin{empheq}[box=\mymath]{equation*}
  \frac{z_2}{z_1}=\frac{x_2+iy_2}{x_1+iy_1}=\frac{x_2+iy_2}{x_1+iy_1}\cdot\frac{x_1-iy_1}{x_1-iy_1}=
      \frac{x_1x_2+y_1y_2}{x_1^2+y_1^2}+i\frac{x_1y_2-x_2y_1}{x_1^2+y_1^2},\:z_1\neq 0.
\end{empheq}
\index{kompleksni brojevi!operacije: sabiranje, oduzimanje,\\ mno\v zenje, dijeljenje}
\newpage
Dakle, dva kompleksna broja sabiramo tako \v sto saberemo realni dio prvog sa realnim dijelom  drugog kompleksnog broja, te isto uradimo i sa imaginarnim dijelovima. Odu-\\zimanje je analogno. Mno\v zenje dva kompleksna broja radimo kao mno\v zenje dva binoma, primenjujemo pravilo "mno\v zenje svaki sa svakim" i koristimo $i^2=-1.$  I na kraju, kod dijeljenja je potrebno ukloniti imaginarni dio iz nazivnika, a ovo radimo koriste\'ci ve\'c navedenu osobinu $z\cdot\overline{z}=(\re{(z)})^2+(\im{(z)})^2.$\\

\begin{example}
 Dati su kompleksni brojevi $z_1=3+4i,\:z_2=2-5i.$\\ Izra\v cunati
  \begin{inparaenum}[$(a)$]
   \item $z_1+z_2;$
   \item $z_1-z_2;$
   \item $z_1\cdot z_2;$
   \item $\frac{z_2}{z_1},$ te odrediti $\re\left( \frac{z_2}{z_2}\right)$
          i $\re\left( \frac{z_2}{z_2}\right);$
   \item $\frac{\overline{z}_1}{z_1+\overline{z}_2}.$
  \end{inparaenum}\ \\\\
\noindent Rje\v senje:\\\\
Vrijedi
 \begin{enumerate}[(a)]
  \item $z_1+z_2=3+4i+2-5i=5-i;$
  \item $z_1-z_2=3+4i-(2-5i)=1+9i;$
  \item $z_1\cdot z_2=(3+4i)\cdot(2-5i)=6-15i+8i-20i^2=6-7i-(-20)=26-7i;$
  \item $\dfrac{z_2}{z_1}=\dfrac{2-5i}{3+4i}\cdot\dfrac{3-4i}{3-4i}=\dfrac{(2-5i)(3-4i)}{9+16}
        =\dfrac{6-8i-15i+20i^2}{25}=\dfrac{-14-23i}{25}=-\dfrac{14}{25}-\dfrac{23}{25}i;$\\ te je
        $\re\left( \frac{z_2}{z_1}\right)=-\frac{14}{25},\;\im\left( \frac{z_2}{z_1}\right)=-\frac{23}{25};$
  \item  $\dfrac{\overline{z}_1}{z_1+\overline{z}_2}=\dfrac{3-4i}{3+4i+2+5i}=\dfrac{3-4i}{5+9i}
         \cdot\dfrac{5-9i}{5-9i}=\dfrac{15-27i-20i+36i^2}{25+81}=\dfrac{-21-47i}{106}.$
 \end{enumerate}
\end{example}
\ \\
\begin{example}
 Odrediti $x$ i $y$ iz jednakosti $z+3x-2\overline{z}=2+3i.$\\\\
 \noindent Rje\v senje:\\\\
 Iz $z=x+yi$ i $\overline{z}=x-yi,$ dobijamo
 \[x+yi+3x-2(x-yi)=2+3i\Leftrightarrow 2x+3yi=2+3i,\] sada je iz prethodne jednakosti $2x=2\Leftrightarrow  x=1$ i $3y=3\Leftrightarrow y=1.$
\end{example}

\section[G\lowercase{eometrijska interpretacija}]{Geometrijska interpretacija kompleksnog broja}

Kompleksni broj $z$ mo\v zemo intrepretirati kao ta\v cku u Dekartovom pravouglom koordinatnom sistemu u ravni. Svakom kompleksnom broju $z$ obostrano--jednozna\v cno \' cemo pridru\v ziti ta\v cku u $xOy$ ravni. Ova ravan je Gaussova\footnote{Johann Carl Friedrich Gau{\ss} (30. april 1777.--23. februar 1855.godine) bio je njema\v cki matemati\v car koji je dao doprinos u mnogim oblastima matematike, kao npr. teorija brojeva, algebra, statistika, analiza, diferencijalna geometrija i dr.} ili kompleksna ravan. Na $x$--osu nanosimo realni dio kompleksnog broja $z$ i zovemo je realna osa, dok na $y$--osu nanosimo imaginarni dio kompleksnog broja i zovemo je imaginarna osa.

 \begin{figure}[!h]\centering
  \includegraphics[scale=.8]{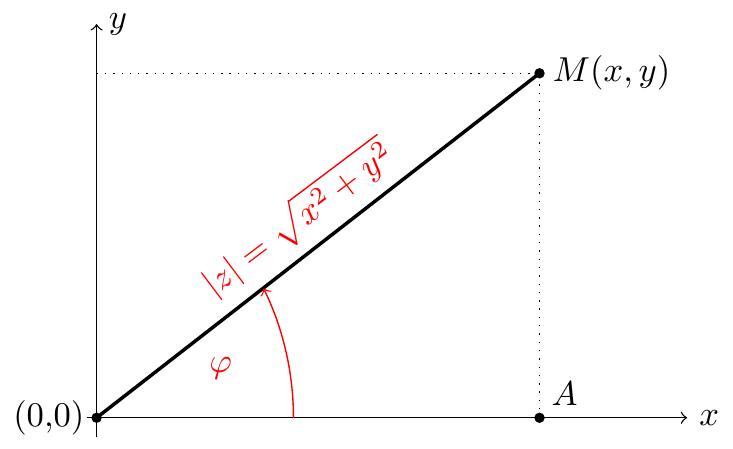}
     \caption{Predstavljanje ta\v cke $M(x,y)$ u Gaussovoj ravni}
     \label{kompl2}
   \end{figure}

Posmatrajmo proizvoljan kompleksni broj $z=x+iy$ predstavljen na Slici \ref{kompl2} ta\v ckom $M(x,y).$
Funkcija $|\cdot|:\mathbb{C}\mapsto\mathbb{R}^{+}\cup\{0\}$ definisana sa
\begin{empheq}[box=\mymath]{equation*}
|z|=\sqrt{x^2+y^2}
\end{empheq}
je modul kompleksnog broja $|z|$ (ili $\rho$). Kako je $z\cdot \overline{z}=x^2+y^2,$ to je  $|z|^2=z\cdot\overline{z}.$ Geometrijski modul $|z|$  kompleksnog broja $z$ predstavlja rastojanje ta\v cke $M(x,y)$ od koordinatnog po\v cetka $(0,0).$

\index{kompleksni brojevi!modul}

Posmatrajmo funkciju $\arg:\mathbb{C}\setminus\{0\}\mapsto(-\pi,\pi]$ definisanu sa
\begin{equation}\label{argumentiz}
 \arg z=\left\{
          \begin{array}{cc}
              \arctg\frac{y}{x},& x>0\\
              \pi+\arctg\frac{y}{x},& x<0\wedge y\geqslant 0\\
              -\pi+\arctg\frac{y}{x},& x<0\wedge y<0\\
              \frac{\pi}{2},&x=0\wedge y>0\\
              -\frac{\pi}{2},& x=0\wedge y<0,
          \end{array}
       \right.
\end{equation}
gdje je $z=x+iy\in\mathbb{C}\setminus\{0\}.$ Broj $\arg z$ (\v cesto se ozna\v cava sa $\varphi$ ili $\theta$) je glavna vrijednost ili glavni argument kompleksnog broja $z.$ \index{kompleksni brojevi!argument}
\begin{remark}
 Argument kompleksnog broja $\arg z$ ($\varphi$ ili $\theta$) predstavlja vrijednost ugla, \v ciji je jedan krak pozitivni dio realne ose a drugi krak je poluprava koja polazi iz koordinatnog po\v cetka i prolazi kroz ta\v cku $M(x,y)$ (Slika \ref{kompl2}).
\end{remark}

 \begin{figure}[!h]
 	\hspace{-1cm}
	\includegraphics[scale=.75]{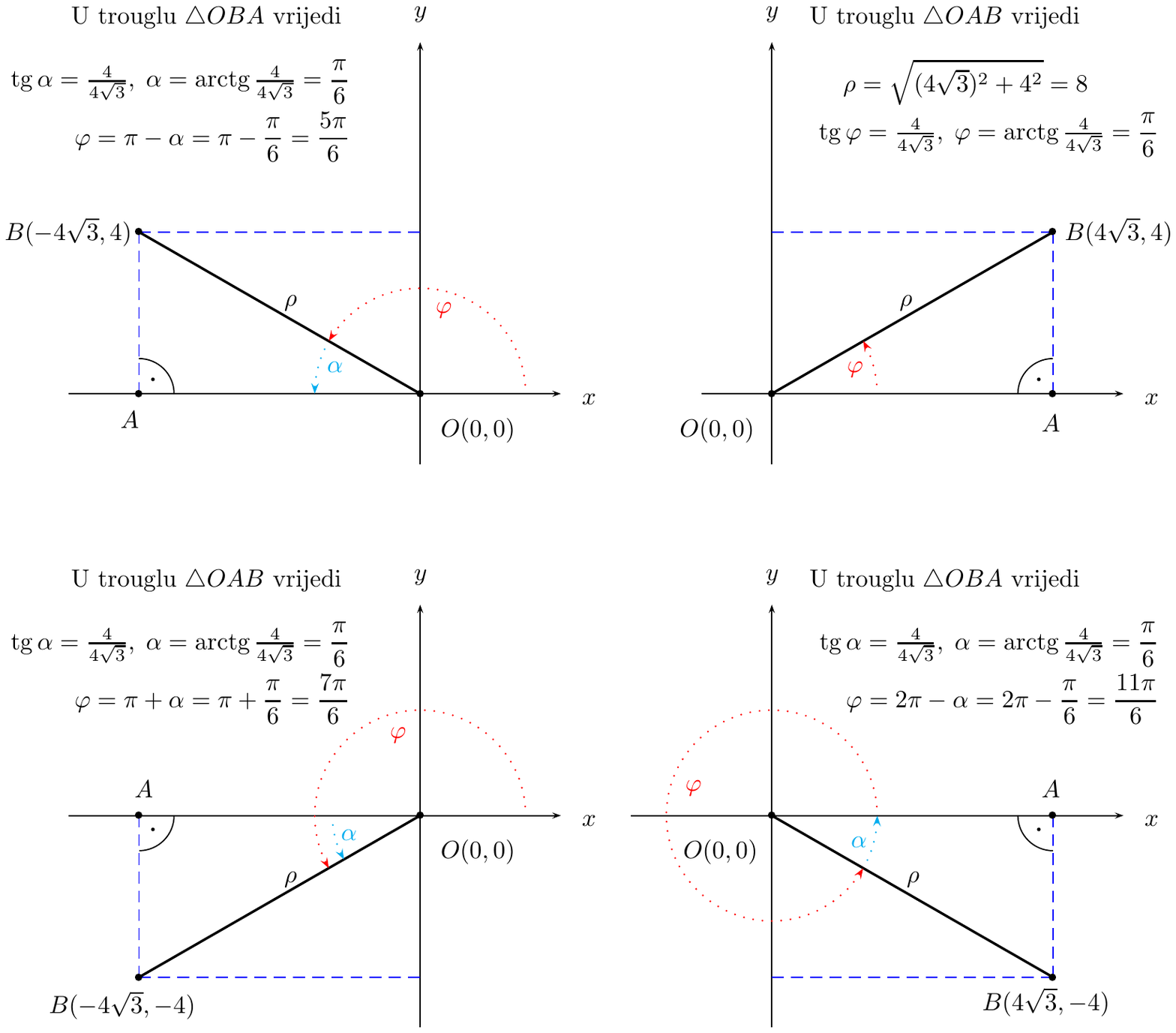}\vspace{-10cm}
	\caption{Odre\dj ivanje ugla}
	\label{kompl3}
\end{figure}

\begin{figure}[!h]\centering
	\includegraphics[scale=.65]{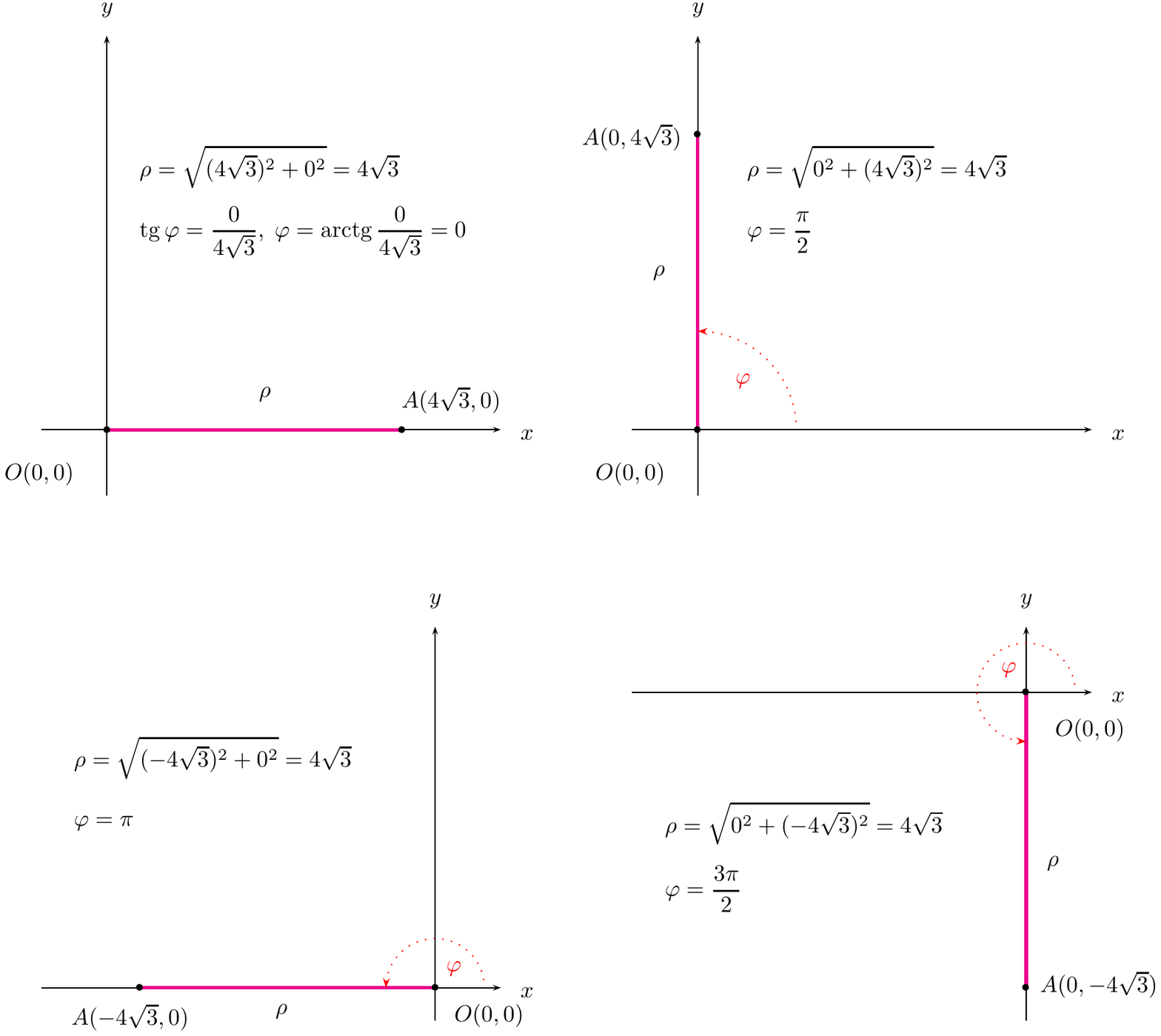}\vspace{-9cm}
	\caption{Odre\dj ivanje ugla za $\varphi=0,\,\pi/2,\, \pi,\, 3\pi/2$}
	\label{kompl3a}
\end{figure}

\begin{example}[Prakti\v cno odre\dj ivanje ugla $\varphi$]\label{primjerUgao1}
 Odrediti vrijednost modula $|z|$ ( ili $\rho$ ) i argumenta $\arg z$ ( ili $\varphi$ ), ako je \\
   \begin{inparaenum}[$(a)$]
      \item $z=4\sqrt{3}+4i;$
      \item $z=-4\sqrt{3}+4i;$
      \item $z=-4\sqrt{3}-4i;$
      \item $z=4\sqrt{3}-4i;$
      \item $z=4\sqrt{3};$
      \item $z=4\sqrt{3}i;$
      \item $z=-4\sqrt{3};$
      \item $z=-4\sqrt{3}i.$
   \end{inparaenum}\\\\
\noindent Rje\v senje: \\\\
Vrijedi (vidi Sliku \ref{kompl3})
\begin{enumerate}[$(a)$]
   \item $x=4\sqrt{3},\,y=4,$ pa je $|z|=\sqrt{(4\sqrt{3})^2+4^2}=8,$ dok je  $\varphi=\arctg\frac{y}{x}=\frac{4}{4\sqrt{3}}=\frac{1}{\sqrt{3}},\:\varphi=\frac{\pi}{6},$

   Slika \ref{kompl3} gore desno.
   \item $x=-4\sqrt{3},\,y=4,$ Slika \ref{kompl3} gore lijevo,  $|z|=\sqrt{(4\sqrt{3})^2+4^2}=8,$ sa slike je
       $\varphi=\pi-\alpha,$ iz trougla $\triangle AOB$ je $\alpha=\arcctg\frac{4}{4\sqrt{3}}=\frac{\pi}{6},$ pa je na kraju $\varphi=\frac{5\pi}{6};$
   \item $x=-4\sqrt{3},y=-4,$   sa Slike \ref{kompl3} dole lijevo je $\varphi=\pi+\alpha,$ dok je
         $\alpha=\arctg\frac{4}{4\sqrt{3}}=\frac{\pi}{6},$    te je $\varphi=\frac{7\pi}{6};$
    \item $x=4\sqrt{3},\,y=-4,$ sa Slike \ref{kompl3} dole desno je          $\varphi=2\pi-\alpha,\,\alpha=\frac{\pi}{6},\:\varphi=\frac{11\pi}{6};$
   \item $x=4\sqrt{3},\,y=0,$ $|z|=4\sqrt{3},\;\varphi=0,$ Slika \ref{kompl3a} gore lijevo;
   \item $x=0,\,y=4\sqrt{3},$ $|z|=4\sqrt{3},\;\varphi=\frac{\pi}{2},$ Slika \ref{kompl3a} gore desno;
   \item $x=-4\sqrt{3},\,y=0,$ $|z|=4\sqrt{3},\;\varphi=\pi,$ Slika \ref{kompl3a} dole lijevo;
   \item $x=0,\,y=-4\sqrt{3},$ $|z|=4\sqrt{3},\;\varphi=\frac{3\pi}{2},$ Slika \ref{kompl3a} dole desno.
\end{enumerate}
\end{example}



\section[T\lowercase{rigonometrijski oblik}]{Trigonometrijski i eksponencijalni oblik \\kompleksnog broja} Trigonometrijski i eksponencijalni oblik kompleksnog broja su prakti\v cniji od algebarskog oblika u slu\v caju kada treba izvr\v siti operacije  mno\v zenja, dijeljenja, zatim stepenovanja i korjenovanja kompleksnih brojeva. Da bi iz algebarskog dobili trigonometrijski oblik kompleksnog broja radimo na sljede\'ci na\v cin.

Za kompleksni broj $z=x+iy,$ $\rho=|z|,\varphi=\arg z$ sa Slike \ref{kompl2} iz trougla $\triangle OAM$ vrijedi
\[x=\rho\cos\varphi,\:y=\rho\sin\varphi,\:\tg\varphi=\frac{y}{x},\,x\neq 0.\]
Vidimo sa Slike \ref{kompl2}, da je ta\v cka $M$ koja predstavlja kompleksni broj $z$ u prvom kvadrantu, pa nalazimo $\varphi=\arctg \frac{y}{x}$.
Ako je $x=0,\,y>0,$ tada je $\varphi=\frac{\pi}{2},$ odnosno za $x=0,\,y<0,$ je $\varphi =-\frac{\pi}{2}.$ U slu\v caju kad je $z$ u nekom drugom kvadrantu, za
nala\v zenje argumenta koristimo \eqref{argumentiz} ili Primjer \ref{primjerUgao1}. Na opisani na\v cin dobijamo trigonometrijski oblik kompleksnog broja $z=x+iy$
\begin{empheq}[box=\mymath]{equation*}
z=\rho(\cos\varphi+i\sin\varphi).
\end{empheq}
\index{kompleksni brojevi!trigonometrijski oblik} \index{kompleksni brojevi!eksponencijalni ili Eulerov oblik}
\index{kompleksni brojevi!operacije: mno\v zenje, dijeljenje }

Ako iskoristimo identitet $\cos \varphi+i\sin\varphi=e^{i\varphi},$ kompleksni broj $z=\rho(\cos\varphi+i\sin\varphi)$ mo\v zemo pisati u obliku
\begin{empheq}[box=\mymath]{equation*}
z=\rho e^{\varphi i}.
\end{empheq}
Ovaj oblik zovemo eksponencijalni ili Eulerov\footnote{Leonhard Euler (15. april 1707.--18. septembar 1783. godine) bio je \v svajcarski matemati\v car, fizi\v car, astronom, logi\v car, in\v zinjer. Dao veliki doprinos u mnogim oblastima matematike.} oblik kompleksnog broja.

\subsection[M\lowercase{no\v zenje i dijeljenje}]{Mno\v zenje i dijeljenje kompleksnog broja u trigonometrijskom obliku} Neka su data dva kompleksna broja u trigonometrijskom obliku $z_1=\rho_1(\cos\theta_1+i\sin\theta_1)$ i $z_2=\rho_2(\cos\theta_2+i\sin\theta_2).$ Koriste\'ci trigonometrijske identitete $\cos(\alpha+\beta)=\cos\alpha\cos\beta-\sin\alpha\sin\beta$ i $\sin(\alpha+\beta)=\sin\alpha\cos\beta+\cos\alpha\sin\beta,$ vrijedi

\begin{align*}
  z_1\cdot z_2&=[\rho_1(\cos\theta_1+i\sin\theta_1)]\cdot[\rho_2(\cos\theta_2+i\sin\theta_2)]\\
              &=\rho_1\rho_2[\cos \theta_1\cos\theta_2 -\sin\theta_1\sin\theta_2 +i(\cos\theta_1\sin\theta_2+\cos\theta_2\sin\theta_1)]\\
              &=\rho_1\rho_2[\cos(\theta_1+\theta_2)+i(\sin(\theta_1+\theta_2))].
\end{align*}
Mno\v zenje dva kompleksna broja u trigonometrijskom obliku vr\v simo po formuli
\begin{empheq}[box=\mymath]{equation*}
z_1\cdot z_2=\rho_1\cdot\rho_2[\cos(\theta_1+\theta_2)+i\sin(\theta_1+\theta_2)].
\end{empheq}
Na sli\v can na\v cin, koriste\'ci poznate trigonometrijske identite $\cos(\alpha-\beta)=\cos\alpha\cos\beta+\sin\alpha\sin\beta,$  $\sin(\alpha-\beta)=\sin\alpha\cos\beta-\cos\alpha\sin\beta$ i $\sin^2\alpha+\cos^2\alpha=1,$  te $z\cdot\overline{z}=(\re z)^2+(\im z)^2,$ dobijamo
\begin{align*}
  \frac{z_1}{z_2}&=\frac{\rho_1}{\rho_2}\frac{\cos\theta_1+i\sin\theta_1}{\cos\theta_2+i\sin\theta_2}\cdot\frac{\cos\theta_2-i\sin\theta_2}{\cos\theta_2-i\sin\theta_2}\\
            &=\frac{\rho_1}{\rho_2}\cdot\frac{\cos\theta_1\cos\theta_2+\sin\theta_1\sin\theta_2+i(\sin\theta_1\cos\theta_2-\cos\theta_1\sin\theta_2)}
              {\cos^2\theta_2+\sin^2\theta_2}\\
            &=\frac{\rho_1}{\rho_2}[\cos(\theta_1-\theta_2)+i\sin(\theta_1-\theta_2)]  .
\end{align*}
Dijeljenje dva kompleksna broja u trigonometrijskom obliku vr\v simo po formuli
\begin{empheq}[box=\mymath]{equation*}
\frac{z_1}{z_2}=\frac{\rho_1}{\rho_2}[\cos(\theta_1-\theta_2)+i\sin(\theta_1-\theta_2)].
\end{empheq}

Dva kompleksna broja data u eksponencijalnom obliku $z_1=\rho_1 e^{\theta_1 i},\, z_2=\rho_2 e^{\theta_2 i},$ mno\v zimo i dijelimo po sljede\'cim formulama

\begin{empheq}[box=\mymath]{equation*}
z_1\cdot z_2=\rho_1e^{\theta_1 i}\cdot\rho_2e^{\theta_2i}=\rho_1\rho_2e^{(\theta_1+\theta_2)i}
\end{empheq}
i
\begin{empheq}[box=\mymath]{equation*}
\frac{z_1}{z_2}=\frac{\rho_1e^{\theta_1 i}}{\rho_2e^{\theta_2 i}}=\frac{\rho_1}{\rho_2}e^{(\theta_1-\theta_2)i}.
\end{empheq}

\begin{figure}[h!]\centering
  \includegraphics[scale=.8]{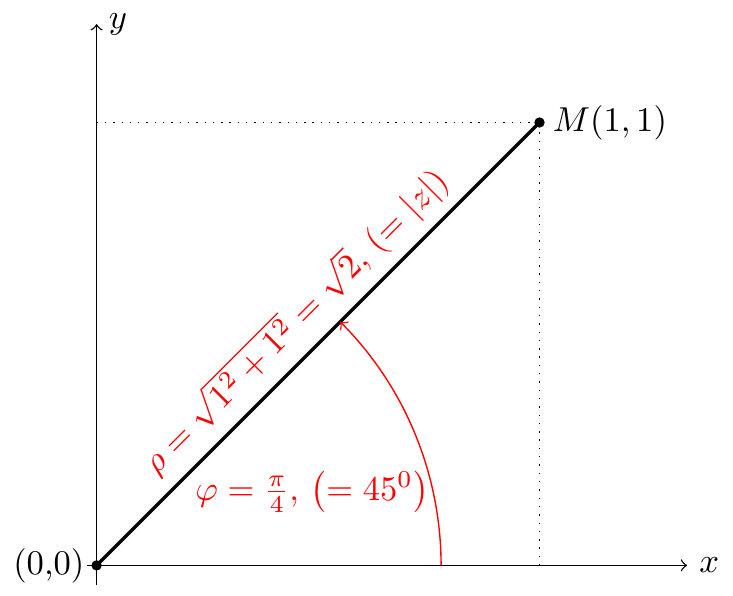}
    \caption{Kompleksni broj $z_1=1+i$ u Gaussovoj ravni}
    \label{kompl4}
\end{figure}

\begin{example}
Dati su kompleksni brojevi $z_1=1+i,\,z_2=8(\cos\frac{5\pi}{6}+i\sin\frac{5\pi}{6}),\,z_3=2(\cos121^0+i\sin 121^0),\,z_4=\cos 14^0+i\sin 14^0.$ Izra\v cunati $\frac{z_1\cdot z_2}{z_3\cdot z_4}.$\\\\
\noindent Rje\v senje:\\\\
Vidi Sliku \ref{kompl4}, prvo prevedemo $z_1$ u trigonometrijski oblik kompleksnog broja, sa Slike \ref{kompl4} je $z_1=\sqrt{2}(\cos 45^0+i\sin45^0),$  pa je sada
  \begin{align*}
     \frac{z_1\cdot z_2}{z_3\cdot z_4}&=\frac{\sqrt{2}(\cos 45^0+i\sin 45^0)\cdot 8(\cos 150^0+i\sin 150^0)}{2(\cos 121^0+i\sin 121^0)\cdot(\cos 14^0+i\sin 14^0)}\\
       &=\frac{8\sqrt{2}(\cos 195^0+i\sin 195^0)}{2(\cos 135^0+i\sin 135^0)}=4\sqrt{2}[\cos(195^0-135^0)+i\sin(195^0-135^0)]\\
       &=4\sqrt{2}(\cos 60^0+i\sin 60^0).
  \end{align*}
  Rezultat mo\v zemo napisati i u algebarskom obliku kompleksnog broja
  \[ \frac{z_1\cdot z_2}{z_3\cdot z_4}=4\sqrt{2}(\cos 60^0+i\sin 60^0)=4\sqrt{2}\left(\frac{1}{2}+\frac{\sqrt{3}}{2}i\right)=2\sqrt{2}+2\sqrt{6} i.\]
\end{example}

\subsection{Stepenovanje kompleksnog broja} Neka je $z_k=\rho_k(\cos\theta_k+i\sin\theta_k),\,k=1,2,\ldots,n.$ Na osnovu prethodno izvedenih formula za mno\v zenje dva kompleksna broja u trigonometrijskom obliku vrijedi
\begin{align*}
z_1\cdot z_2=&\rho_1\cdot\rho_2[\cos(\theta_1+\theta_2)+i\sin(\theta_1+\theta_2)],\\
(z_1\cdot z_2)\cdot z_3=&\rho_1\cdot\rho_2\cdot\rho_3[\cos(\theta_1+\theta_2+\theta_3)+i\sin(\theta_1+\theta_2+\theta_3)].
\end{align*}
Na osnovu matemati\v cke indukcije\footnote{Princip matemati\v cke indukcije: Jedan iskaz $P(n)$ istinit je za svaki prirodan brod $n,$ 1) Ako je istinit za prirodan broj 1; 2) Ako implikacija $P(n)\Rightarrow P(n+1)$  va\v zi za svaki prirodan broj $n$.} se lako pokazuje da vrijedi
\[z_1\cdot z_2\cdot\ldots\cdot z_n=\rho_1\cdot\rho_2\cdot\ldots\cdot\rho_n[\cos(\theta_1+\theta_2+\ldots+\theta_n)+i\sin(\theta_1+\theta_2+\ldots+\theta_n)].\]

Ako je $z_1=z_2=\ldots= z_n,$ onda je $\rho_1=\rho_2=\ldots=\rho_n$ i $\theta_1=\theta_2=\ldots=\theta_n,$ dobijamo formule za stepenovanje dva kompleksna broja u trigonometrijskom
\begin{empheq}[box=\mymath]{equation*}
z^n=\rho^n(\cos n\theta+i\sin n\theta)
\end{empheq}
ili u eksponencijalnom obliku
\begin{empheq}[box=\mymath]{equation*}
z^n=\rho^ne^{n\theta i}.
\end{empheq}
\index{kompleksni brojevi!stepenovanje}
\index{kompleksni brojevi!Moivreov  obrazac ili formula}

Izvode\'ci formulu za stepenovanje dva kompleksna broja u trigonometrijskom obliku, dobija se i Moivreov\footnote{Abraham de Moivre (26.maj 1667.--27.novembar 1754. godine) bio je francuski matemati\v car poznat po Moivreovom obracscu koji povezuje kompleksne brojeve i trigonometriju, te doprinosu u vjerovatno\'ci i statistici. } (Moavrov) obrazac (ili formula)
\begin{empheq}[box=\mymath]{equation*}
(\cos\theta+i\sin\theta)^n=\cos n\theta+i\sin n\theta,
\end{empheq}
koja vrijedi i kada je $n$ ne samo prirodan  nego i negativan cio broj i nula, tj. vrijedi za sve $n\in\mathbb{Z}.$

\subsection{Korjenovanje kompleksnog broja}

\index{kompleksni brojevi!korjenovanje}
Neka je data jedna\v cina \[z^n=u,\] gdje je $n\in\mathbb{N},$ a $u$ je kompleksni broj razli\v cit od nule. Tada pi\v semo \[z=\sqrt[n]{u}.\]  Potrebno je za dati broj $u$ odrediti kompleksni broj $z.$ Neka je $u=r(\cos\varphi+i\sin\varphi),\,z=\rho(\cos\theta+i\sin\theta).$ Tada je, na osnovu Moivreovog obrasca,  $(\rho(\cos\theta+i\sin\theta))^n=\rho^n(\cos n\theta +i\sin n\theta)=r(\cos\varphi+i\sin\varphi),$ pa na osnovu jednakosti kompleksnih brojeva $\rho^n=r,\,n\theta=\varphi+2k\pi,\,k\in\mathbb{Z},$ odakle je
\[\rho=\sqrt[n]{r},\,\theta=\frac{\varphi+2k\pi}{n},\,k\in\mathbb{Z}.\]

Me\dj utim, razli\v cite vrijednosti cijelih brojeva $k$ ne daju uvijek razli\v cite argumente. Ako uzmemo da je $k=n$ bi\' ce $\frac{\varphi+2n\pi}{n}=\frac{\varphi}{n}+2\pi,$
ista se vrijednost dobije i za $k=0.$ Sli\v cno, dobijamo da je $\frac{\varphi+2(n+k)\pi}{n}=\frac{\varphi+2k\pi}{n}+2\pi.$ Dakle, poslije $n$ uzastopnih cijelih brojeva kompleksni brojevi se ponavljaju, pa imamo $n$ razli\v citih vrijednosti $\sqrt[n]{u},$ a to su
\begin{empheq}[box=\mymath]{equation*}
z_k=\sqrt[n]{r}\left(\cos\frac{\varphi+2k\pi}{n}+i\sin\frac{\varphi+2k\pi}{n} \right),\:k=0,\ldots,n-1,
\end{empheq}
ili u eksponencijalnom obliku
\begin{empheq}[box=\mymath]{equation*}
z_k=\sqrt[n]{u}=\sqrt[n]{r}e^{\frac{\varphi+2k\pi}{n}i},\,k=0,\ldots,n-1.
\end{empheq}

Sva ova rje\v senja le\v ze na kru\v znici polupre\v cnika $\sqrt[n]{r}$ i \v cine tjemena pravilnog $n$--tougla, \v cije je jedno tjeme ($z_0$) sa argumentom $\frac{\varphi}{n},$ a za svako sljede\' ce tjeme argument se pove\' cava za $\frac{2\pi}{n}.$

\begin{example}
Dat je kompleksni broj $z=i^{81}+i^{43}+\frac{(-1-i)^{80}}{{2^{40}}}+i^{19}.$ Izra\v cunati $\sqrt[3]{z}.$\\\\
\noindent Rje\v senje:\\\\
Izra\v cunajmo prvo $\frac{(-1-i)^{80}}{2^{40}},$ prevedemo ovaj kompleksni broj u trigonometrijski oblik, sa Slike \ref{slika32} je
  \begin{align*}
      \frac{(-1-i)^{80}}{2^{40}}&=\frac{\left[\sqrt{2}\left( \cos\frac{5\pi}{4}+i\sin\frac{5\pi}{4}\right)\right]^{80}}{2^{40}}\\
      &=\frac{2^{40}\left[\cos\left(\frac{5\pi}{4}\cdot 80\right)+i\sin\left(\frac{5\pi}{4}\cdot 80\right)\right]}{2^{40}}=\cos 100\pi+i\sin 100\pi=1,
  \end{align*}
  pa je sada
  \[z=i^{81}+i^{43}+\frac{(-1-i)^{80}}{{2^{40}}}+i^{19}=i^{4\cdot 20+1}+i^{4\cdot 10+3}+1+i^{4\cdot 4+3}=1-i.\]
Koriste\' ci Sliku \ref{slika33} pretvorimo kompleksni broj iz algebarskog u trigonometrijski oblik, i dobijamo
\[z=\sqrt{2}\left(\cos\frac{7\pi }{4}+i\sin\frac{7\pi}{4}\right).\]
Sada je
\[\sqrt[3]{z}=\sqrt[3]{\sqrt{2}}\left(\cos\frac{\frac{7\pi}{4}+2k\pi}{3}+i \sin\frac{\frac{7\pi}{4}+2k\pi}{3}\right),\]
te rje\v senja $w_0,\,w_1$ i $w_2$ dobijamo uvr\v stavaju\' ci $k=0,1,2$ u prethodnu jednakost. Vrijedi
\begin{align*}
   k&=0\\
   w_0&=\sqrt[6]{2}\left(\cos\frac{7\pi}{12}+i\sin\frac{7\pi}{12}\right),\\
   k&=1\\
   w_1&=\sqrt[6]{2}\left(\cos\frac{15\pi}{12}+i\sin\frac{15\pi}{12}\right),\\
   k&=2\\
   w_2&=\sqrt[6]{2}\left(\cos\frac{23\pi}{12}+i\sin\frac{23\pi}{12}\right).\\
\end{align*}
Rje\v senja su predstavljena na Slici \ref{slika34}.
\end{example}

\begin{figure}[!h]
     \begin{subfigure}[b]{.3\textwidth}\centering
        \includegraphics[scale=.8]{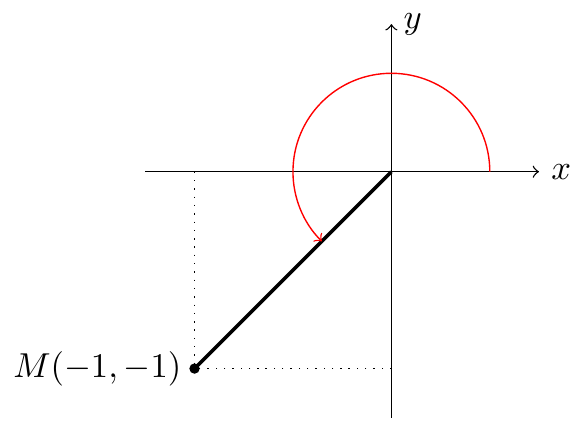}
         \caption{$z=-1-i$}
          \label{slika32}
   \end{subfigure}\hspace{1cm}
   \begin{subfigure}[b]{.3\textwidth}\centering
        \includegraphics[scale=.8]{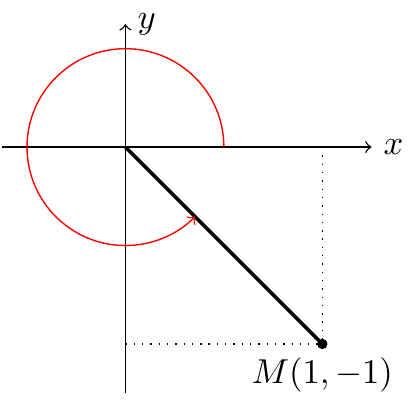}
         \caption{$z=1-i$}
      \label{slika33}
     \end{subfigure}\hspace{-.5cm}
     \begin{subfigure}[b]{.3\textwidth}\centering
        \includegraphics[scale=.8]{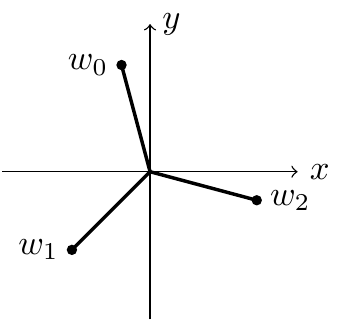}
         \caption{Rje\v senje jedna\v cine $\sqrt[3]{z}$}
      \label{slika34}
      \end{subfigure}
        \caption{Gau\ss ova ravan u kojoj su predstavljeni kompleksni brojevi $-1-i,\,1-i,\,w_0,\,w_1$ i $w_2$ }
     \label{slika35}
  \end{figure}

\section{Z\lowercase{adaci}}\index{Zadaci za vje\v zbu!kompleksni brojevi}

\begin{enumerate}
    \item Dati su kompleksni brojevi $z_1=2-3i,\,z_2=1+2i.$ Izra\v cunati\\
        \begin{inparaenum}
            \item $z_1+z_2;\,$
            \item $z_1-z_2;\,$
            \item $z_1\cdot z_2;\,$
            \item $\dfrac{z_1}{z_2};\,$
            \item $|z_1|;\,$
            \item $\left| \frac{z_1}{z_2}\right| ;\,$
            \item $\re\left( \frac{z_1}{z_2}\right);\,$
            \item $\im\left( \frac{z_1}{z_2} 	\right).$
        \end{inparaenum}
  \item Izra\v cunati \\
       \begin{inparaenum}
          \item $i^{25}+(-i)^{50}+i^{40}+i^{62}+i^{83} ;\,$
          \item $\dfrac{1+3i}{(-1+i)^2}+\dfrac{(-4+i)(-4-i)}{1+i} ;\,$
          \item Za $z=1-3i,\,\dfrac{2z-2z\overline{z}}{z\overline{z}+i}.$
       \end{inparaenum}
  \item  Izra\v cunati\\
      \begin{inparaenum}
        \item $x$ i $y$ iz jedna\v cine $3x+xi-2y=12-yi-i,\,$
        \item $z$ iz uslova $(2+i)z+2z-3=4-6i.$
\end{inparaenum}
  \item Odrediti kompleksni broj $z$ iz uslova $|z-4|=|z-2|\wedge |z-3|=|z-2i|.$
  \item Dati su kompleksni brojevi $z_1= 2-2\sqrt{3}i,\,z_2= -6\sqrt{3}-6i,\,z_3 -5+5i.$ \\
     \begin{enumerate}
         \item Pretvoriti date kompleksne brojeve iz algebarskog u trigonometrijski oblik;
         \item Izra\v cunati $\dfrac{z_1\cdot z_2}{z_3};$
         \item Izra\v cunati $\dfrac{z_3^{40}\cdot z_2^{60}}{z_2^{30}}.$
     \end{enumerate}
  \item Kompleksni broj $z=\frac{1+4i+i^{29}}{2-3i}$ napisati u trigonometrijskom obliku;\:

  \item    Izra\v cunati\\
     \begin{inparaenum}
        \item $ \sqrt[3]{-3-\sqrt{3}i};\,$
        \item $ \sqrt[4]{3-3\sqrt{3}i};\,$
        \item $z$ ako je $z^3=2-2\sqrt{3}i.$
     \end{inparaenum}

   \item Izra\v cunati kompleksni broj $ z$ iz uslova
         $\re\left( \frac{z-\overline {z}}{2}\right)+10i\im\left( \frac{\overline{z}-1}{3+i}\right)=x+1+2i,\:$ a zatim izra\v cunati $\sqrt{z}.$
    \item Rije\v siti jedna\v cinu $z^3+2+2i=0.\:$

  \item  Razni zadaci
    \begin{enumerate}
      \item Rije\v siti jedna\v cinu $z^5-\sqrt{3}+i=0.\:$
      \item Dat je kompleksni broj $z=4\sqrt{3}-12i$, izra\v cunati $\sqrt[4]{z}.\:$
      \item Dat je kompleksni broj $z=-6-2\sqrt{3}i$, izra\v cunati $\sqrt[4]{z}.\:$
      \item Rije\v siti jedna\v cinu $z^5-\sqrt{3}+i=0.\:$
      \item Ako je $z=5-5\sqrt{3}i$, odrediti $z^5$ i $\sqrt[3]{z}.\:$
      \item Rije\v siti jedna\v cinu u skupu racionalnih brojeva $4z^3+7\sqrt{3}i=7.\:$
      \item U skupu kompleksnih brojeva rije\v siti jedna\v cinu $z^4-3\sqrt{3}=-3i.\:$
      \item U skupu kompleksnih brojeva rije\v siti jedna\v cinu $z^4+2i=\sqrt{3}-i.\:$
      \item Rije\v siti jedna\v cinu $z=\sqrt{\dfrac{1+i}{\sqrt{2}}}.\:$
       \item Rije\v siti jedna\v cinu $z=\sqrt[3]{\dfrac{1+i}{\sqrt{2}}}.\:$
      \item U skupu kompleksnih brojeva rije\v siti jedna\v cinu $\sqrt{2}z^3+2i=3i+1.\:$
      \item Izra\v cunati kompleksni broj $z$ iz uslova
         $i\re\left( \frac{\overline{z}+2}{1+i}\right)+\im\left( \frac{2\overline{z}+z}{2}\right)+z=1+3i,\:$
      \item Izra\v cunati kompleksni broj $z$ iz uslova
        $\re\left( \frac{\overline{z}}{z_1}\right)=\frac{2-3\sqrt{3}}{13}$ i $\im\left(\frac{z\cdot z_1}{2} \right)=\frac{-3-2\sqrt{2}}{2},$ gdje je $z_1=2-3i,$ a zatim izra\v cunati $z^{12}.$
      \item Izra\v cunati kompleksni broj $z$ iz uslova $\left|\dfrac{z-12}{z-8i} \right|=\dfrac{5}{3}$  i $\left| \dfrac{z-4}{z-8}\right|=1.$
      \item Izra\v cunati $z$ i $\sqrt[3]{z}$ ako je $|z-2|+\overline{z}\re(2z)+z^2-4z=4i.\:$
      \item Izra\v cunati $z^{50}$ ako je $2\overline{z}-i\re\left( z+\frac{1}{2}\right)=1-i\sqrt{3}-i.\:$
      \item Izra\v cunati $\sqrt{z}$ ako je $\frac{2}{i}\left( z(2-i)+\frac{\sqrt{3\sqrt{3}}}{2}\right)=-3(1+2\sqrt{3}+2i).$
      \item Pokazati da $\re\left(\frac{1}{z}-\frac{2}{3} \right)=0$   predstavlja jedna\v cinu kru\v znice.
    \end{enumerate}
\end{enumerate}

   \chapter{Matrice i determinante}


Pojam matrice i determinante je izuzetno va\v zan kako iz ugla same matematike, tako i njene primjene u drugim oblastima. Vidje\'{c}emo da je izu\v cavanje sistema linearnih algebarskih jedna\v cina znatno olak\v sano kori\v stenjem matrica, odnosno matri\v cnog ra\v cuna, zatim tu je upotreba u vektorskoj algebri, analiti\v ckoj geometriji i dr. U ovom poglavlju su dati osnovni pojmovi o matricama i determinantama, neophodni za njihovu primjenu u oblastima koje se analiziraju u nastavku.

\section[O\lowercase{snovni pojmovi o  matricama}]{Osnovni pojmovi o  matricama}\index{matrica}

Prelazak promjenljivih $x_1,x_2,\ldots,x_n$ na nove promjenljive $y_1,y_2,\ldots,y_n$ pomo\' cu formula
\begin{align}
   y_1&=a_{11}x_1+a_{12}x_2+\ldots+a_{1n}x_n,\nonumber\\
   y_2&=a_{21}x_1+a_{22}x_2+\ldots+a_{2n}x_n,\nonumber\\
   &\hspace{1cm}\vdots\label{matrica3}\\
   y_n&=a_{n1}x_1+a_{n2}x_2+\ldots+a_{nn}x_n,\,a_{ij}\in\mathbb{C},\,i,j=1,\ldots,n\nonumber
 \label{matrice1}
\end{align}
zove se linearna transformacija promjenljivih $x_1,\ldots,x_n$ u promjenljive $y_1,\ldots,y_n.$  Koeficijente promjenljivih $x_1,\ldots,x_n$  mo\v zemo izdvojiti i zapisati ih u sljede\' cu shemu
\begin{equation}
         \left( \begin{array}{cccc}
                a_{11}&a_{12}&\ldots&a_{1n}\\
                a_{21}&a_{22}&\ldots&a_{2n}\\
                &\vdots&&\\
                a_{n1}&a_{n2}&\ldots&a_{nn}
                \end{array}
         \right).
  \label{matrica2}
\end{equation}
Shemu \eqref{matrica2} zovemo matrica (ili matrica linearne transformacije \eqref{matrica3}).
\begin{enumerate}[$\bullet$]
  \item Koeficijenti $a_{ij},\,i,j=1,\ldots,n,$ su elementni matrice.
  \item Elementi $a_{i1},\ldots,a_{in},\,i=1,\ldots,n,$ \v cine $i$--tu vrstu matrice.
  \item Elementi $a_{1j},\ldots,a_{nj},\,j=1,\ldots,n,$ \v cine $j$--tu kolonu matrice.
  \item Elementi $a_{11},\ldots,a_{nn},$ \v cine glavnu dijagonalu matrice.
  \item Elementi $a_{1n},\,a_{2,n-1},\ldots,a_{n1},$ \v cine sporednu dijagonalu matrice.
  \item Matricu kra\' ce zapisujemo $A=(a_{ij}).$
\end{enumerate}
\newpage
Matrice mo\v zemo definisati i na sljede\'{c}i na\v cin:
\begin{definition}
Preslikavanje $A$ iz skupa $S=\{1,2,\ldots , m\}\times \{1,2,\ldots , n\}$ u polje $F,$ koje svakom ure\dj enom paru $(i,j) \in S$ pridru\v zuje element
$A(i,j)=a_{ij} \in F,$ naziva se matrica $A=(a_{ij})$ tipa (formata) $m \times n$ nad poljem $F.$
\end{definition}
Za oznaku matrice mo\v zemo koristiti malu ili srednju zagradu, dakle $(\cdot )$ ili $[ \cdot ].$ Ako \v zelimo naglasiti kojeg je formata matrica, koristimo
oznaku $(a_{ij})_{m \times n}.$
U formuli \eqref{matrica3} ve\'c je navedeno da su elementi matrica $a_{ij}\in\mathbb{C},$ tj. da su u op\v stem slu\v caju kompleksni brojevi. Kako je $\mathbb{R}\subset\mathbb{C},$ to realne brojeve mo\v zemo shvatiti kao specijalni slu\v caj kompleksnih brojeva i to su oni kompleksni brojevi kod kojih je imaginarni dio jednak nuli. Stoga elementi matrica mogu biti i realni brojevi. U  nastavku \'cemo uglavnom raditi sa matricama \v ciji su elementi realni brojevi.

Matrica kod koje je jednak broj vrsta i kolona, kao kod matrice $A$ iz \eqref{matrica2}, zovemo kvadratna matrica. U slu\v caju da broj kolona i vrsta nije jednak dobijamo pravougaonu matricu
\index{matrica!kvadratna}\index{matrica!pravougaona}
\begin{equation*}
         \left( \begin{array}{cccc}
                a_{11}&a_{12}&\ldots&a_{1n}\\
                a_{21}&a_{22}&\ldots&a_{2n}\\
                &\vdots&&\\
                a_{m1}&a_{m2}&\ldots&a_{mn}
                \end{array}
         \right).
  \label{matrica2a}
\end{equation*}

Za matricu koja ima $m$ vrsta i $n$ kolona, ka\v zemo da je tipa ili formata ili dimenzije  $m\times n$ ili jednostavno matrica $m\times n$ i pi\v semo \index{matrica!format ili dimenzija}
\begin{empheq}[box=\mymath]{equation*}
A_{m\times n}.
\end{empheq}
Matrica formata $1\times n$ je tako\dj e matrica formata $m\times n,$ tj. $m=1,$
\[ \left(a_1\:a_2\:\ldots\:a_{n}\right), \]
a zovemo je matrica vrsta. Na isti na\v cin rezonujemo za format $m\times 1,$ sada je $n=1$
\[ \left(
          \begin{matrix}
           	b_1\\b_2\\\vdots\\b_m
          \end{matrix}
   \right)
\]
i zovemo je matrica kolona. U slu\v caju matrice kolone i matrice vrste, obi\v cno pi\v semo jedan indeks elementa. Matrice kolone i matrice vrste nazivaju se i vektori. U slu\v caju $m=n=1$ dobijamo matricu sa samo jednim elementom $(a_1),$ koju zovemo skalar.
Skup svih vektora tipa $n\times 1$ sa realnim odnosno kompleksnim elementima, obilje\v zava se sa $\mathbb{R}^n,$ odnosno sa $\mathbb{C}^n$,  respektivno.

\index{matrica!kolona}\index{matrica!vrsta}

\begin{example}
  Date su matrice
  \[A=\left( \begin{array}{rrrr}4&-1&2&0\\3&1&2&1\end{array}\right),\:
   B=\left( \begin{array}{rrr}-1&2&0\\1&2&1\\0&0&-4\end{array}\right),\:
   C=\left( \begin{array}{rr}4&-1\\2&1\\3&5\end{array}\right).\]
  Matrice $A$ i $C$ su pravougaone, formata $2\times 4,\:3\times 2,$ respektivo,  dok je matrica $B$ kvadratna matrica, formata $3\times 3.$ Prethodno mo\v zemo zapisati: $A_{2\times 4},\,B_{3\times 3},\,C_{3\times 2}.$
\end{example}

\paragraph{Jednakost matrica.}\index{matrica!jednake}
Kao i kod drugih matemati\v ckih objekata, koji mogu biti jednaki, tako postoje i matrice koje su jednake.  Uslovi pod kojima su dvije matrice jednake dati su u sljede\'coj definiciji.
\begin{definition}
Dvije matrice jednake su ako su istog formata i ako su im odgovaraju\' ci elementi jednaki.
\end{definition}
Drugim rije\v cima, matrica $A=(a_{ij})$ i $B=(b_{ij})$ formata $m\times n$ su jednake ako je ispunjeno $m\times n$ uslova
\begin{empheq}[box=\mymath]{equation*}
a_{ij}=b_{ij},\,i=1,\ldots,m,\:j=1,\ldots,n.
\end{empheq}
\begin{example}
  Date su matrice
  \[A=\left( \begin{matrix}4&-1&2&0\\3&1&2&1\end{matrix}\right),\,
   B=\left( \begin{matrix} -1&2&0\\1&2&1\\0&0&-4\end{matrix}\right),\,
   C=\left( \begin{matrix} 4&-1\\2&1\\3&5\end{matrix}\right),\,
   D=\left( \begin{matrix}4&-1&2&0\\3&1&2&1\end{matrix}\right).\]
  Samo su matrice $A$ i $D$ jednake, tj. $A=D.$
\end{example}

\subsection{Operacije sa matricama}

U ovoj podsekciji obra\dj eno je sabiranje i oduzimanje matrica, te mno\v zenje matrice sa brojem (skalarom).

\paragraph{Sabiranje matrica.}\index{matrica!sabiranje} Sabiranje ili oduzimanje matrica svodi se na sabiranje ili odu-\\zimanje elemenata matrica. Sabiranje odnosno oduzimanje dvije matrice mo\v zemo vr\v siti samo ako su istog formata.  Slijedi precizna definicija.
\begin{definition}
  Zbir matrica $A=(a_{ij})$ i $B=(b_{ij}),$ formata $m\times n$, u oznaci $A+B$ je matrica $C=(c_{ij}),$ formata $m\times n,$
  gdje je $c_{ij}=a_{ij}+b_{ij},\,i=1,\ldots,m,\:j=1,\ldots,n.$
\end{definition}

Drugim rije\v cima, matrice sabiramo ako su istog formata i to radimo tako \v sto im sabiramo elemente na istim pozicijama. Analogno se vr\v si i oduzimanje, tj. $C=A-B$ vr\v simo tako \v sto od elementa $a_{ij}$ matrice $A$ oduzimamo  odgovaraju\' ci element $b_{ij}$ matrice $B,$ te je $c_{ij}=a_{ij}-b_{ij}.$
\paragraph{Mno\v zenje matrice brojem.}\index{matrica!mno\v zenje skalarom} Osim sabiranja i oduzimanja matrica, nekada je potrebno i pomno\v ziti matricu sa brojem, ovaj postupak dat je u sljede\'coj definiciji.
\begin{definition}
  Proizvod matrice \[
         \left( \begin{array}{cccc}
                a_{11}&a_{12}&\ldots&a_{1n}\\
                a_{21}&a_{22}&\ldots&a_{2n}\\
                &\vdots&&\\
                a_{m1}&a_{m2}&\ldots&a_{mn}
                \end{array}
         \right)
                   \] i broja $\alpha$ defini\v se se sa
        \[\alpha\cdot \left( \begin{array}{cccc}
                a_{11}&a_{12}&\ldots&a_{1n}\\
                a_{21}&a_{22}&\ldots&a_{2n}\\
                &\vdots&&\\
                a_{m1}&a_{m2}&\ldots&a_{mn}
                \end{array}
         \right)= \left( \begin{array}{cccc}
                \alpha a_{11}&\alpha a_{12}&\ldots&\alpha a_{1n}\\
                \alpha a_{21}&\alpha a_{22}&\ldots&\alpha a_{2n}\\
                &\vdots&&\\
                \alpha a_{m1}&\alpha a_{m2}&\ldots&\alpha a_{mn}
                \end{array}
         \right).\]
\end{definition}

Pomno\v ziti  matricu $A$ sa brojem $\alpha,$ zna\v ci svaki element matrice $a_{ij} $ pomno\v ziti sa brojem $\alpha$.

\begin{example}
  Date su matrice
   \[A=\left( \begin{array}{rrrr}2&1&2&0\\3&1&2&1\end{array}\right),\:
     B=\left( \begin{array}{rrrr}3&-1&-5&0\\0&-1&2&1\end{array}\right),\:
     C=\left( \begin{array}{rrrr}-2&0&2&0\\3&1&2&1\end{array}\right).\]
     Izra\v cunati
     \begin{enumerate}[$(a)$]
       \item $A+B;$
       \item $C-A$;
       \item $2A-3B+4C.$
     \end{enumerate}
\noindent Rje\v senje:\\\\
Vrijedi

 \begin{align*}
	(a) \: \quad\quad \quad \:\:
    A+B    &=\left( \begin{array}{rrrr}2&1&2&0\\3&1&2&1\end{array}\right)+\left( \begin{array}{cccc}3&-1&-5&0\\0&-1&2&1\end{array}\right)   \\
           &=\left( \begin{array}{rrrr}2+3&1-1&2-5&0+0\\3+0&1-1&2+2&1+1\end{array}\right)=\left( \begin{array}{cccc}5&0&-3&0\\3&0&4&2\end{array}\right);\\
    (b) \: \quad\quad \quad \:\:
    C-A    &=\left( \begin{array}{rrrr}-2&0&2&0\\3&1&2&1\end{array}\right)-\left( \begin{array}{cccc}2&1&2&0\\3&1&2&1\end{array}\right) \\
           &=\left( \begin{array}{rrrr}-2-2&0-1&2-2&0-0\\3-3&1-1&2-2&1-1\end{array}\right)=\left( \begin{array}{cccc}-4&-1&0&0\\0&0&0&0\end{array}\right);\\
    (c) \:
  2A-3B+4C &=2\left( \begin{array}{rrrr}2&1&2&0\\3&1&2&1\end{array}\right)-3\left( \begin{array}{cccc}3&-1&-5&0\\0&-1&2&1\end{array}\right) \\
           &\quad +4\left( \begin{array}{rrrr}-2&0&2&0\\3&1&2&1\end{array}\right) \\
           &=\left( \begin{array}{rrrr}2\cdot 2&2\cdot1&2\cdot2&2\cdot0\\2\cdot3&2\cdot1&2\cdot2&2\cdot1\end{array}\right) \\
           &\quad-\left( \begin{array}{rrrr}3\cdot3&3\cdot(-1)&3\cdot(-5)&3\cdot0\\3\cdot0&3\cdot(-1)&3\cdot2&3\cdot1\end{array}\right) \\
           &\quad+\left( \begin{array}{rrrr}4\cdot (-2)&4\cdot0&4\cdot2&4\cdot0\\4\cdot3&4\cdot1&4\cdot2&4\cdot1\end{array}\right) \\
           &=\left( \begin{array}{rrrr}4&2&4&0\\6&2&4&2\end{array}\right)
                -\left( \begin{array}{rrrr}9&-3&-15&0\\0&-3&6&3\end{array}\right)
                +\left( \begin{array}{rrrr}-8&0&8&0\\12&4&8&4\end{array}\right) \\
           &=\left( \begin{array}{rrrr}4-9-8&2+3+0&4+15+8&0-0+0\\6-0+12&2+3+4&4-6+8&2-3+4\end{array}\right) \\
           &=\left( \begin{array}{rrrr}-13&5&27&0\\18&9&6&3\end{array}\right).
	\end{align*}

\end{example}
U sljede\'coj teoremi navedene su neke osobine mno\v zenja matrice brojem. Ove osobine proizilaze iz upravo definisanih operacija sa matricama.
\begin{theorem}\index{teorema!mno\v zenje matrice brojem--osobine}
 Mno\v zenje matrice brojem ima osobine
 \begin{enumerate}[$1.$]
   \item $\alpha(A+B)=\alpha A+\alpha B;$
   \item $(\alpha+\beta)A=\alpha A+\beta A;$
   \item $(\alpha \beta)A=\alpha(\beta A);$
   \item $1\cdot A=A.$
 \end{enumerate}
 gdje su $\alpha,\beta\in\mathbb{C}$ i $A,\,B$ su matrice istog formata.
\end{theorem}

\paragraph{Mno\v zenje matrica.} \index{matrica!proizvod matrica} Osim mno\v zenja matrice brojem, matricu mo\v zemo, pod odre\dj enim uslovima, pomno\v ziti i drugom matricom. I ova se operacija na kraju svodi na mno\v zenje i sabiranje elemenata matrice. Detalji su dati u sljede\'coj definiciji.
\begin{definition}
  Proizvod matrica $A=(a_{ij})_{m\times n}$ i $B=(b_{jk})_{n\times p},$ u oznaci $A\cdot B$ ili $AB$, je matrica
  \[C=(c_{ik})_{m\times p}=\left(\sum_{j=1}^{n}{a_{ij}b_{jk}} \right)_{m\times p}.\]
\end{definition}
Dakle, za dvije matrice $A=(a_{ij})_{m\times {\textbf{\color{magenta}n}}}$ i $B=(b_{jk})_{{\textbf{\color{magenta}n}}\times p}$ mo\v zemo izra\v cunati proizvod $A\cdot B$ ili $AB,$ samo ako je broj kolona matrice $A$ jednak broju vrsta matrice $B,$ (ozna\v ceno sa {\color{magenta}magentom}/{\textbf{bold}}). Sada, element $c_{ik}$ matrice $C$ dobijamo na taj na\v cin \v sto se elementi $i$--te vrste matrice $A,$ \[(a_{i1} \,a_{i2} \,\ldots\,a_{in})\]
pomno\v ze redom sa odgovaraju\' cim elementima $k$--kolone matrice $B$
\[\left( \begin{array}{c}b_{1k}\\b_{2k}\\\vdots\\b_{nk}\end{array} \right)\] i dobijeni proizvodi se saberu.

Mno\v zenje kompleksnih, pa samim tim i realnih brojeva, je komutativna operacija, vrijedi $ab=ba,\:\forall a,b\in\mathbb{C}.$ Medjutim, mno\v zenje matrica u op\v stem slu\v caju nije komutativna operacija. Samo za neke matrice vrijedi ova osobina, za one matrice koje imaju ovu osobinu dajemo posebno ime. \index{matrica!komutativne matrice}
\begin{definition}
 Ako je $AB=BA,$ matrice   $A$ i $B$ zovu se komutativne matrice.
\end{definition}
U sljede\'coj teoremi date su neke osobine mno\v zenja matrica.

\begin{theorem}\index{teorema! proizvod matrica--osobine}
  Mno\v zenje matrica ima osobine
   \begin{enumerate}
     \item $(AB)C=A(BC)$ (asocijativnost);
     \item $(A+B)C=AC+BC$ (distributivnost mno\v zenje prema sabiranju sa desne strane);
     \item $A(B+C)=AB+AC$ (distributivnost mno\v zenja prema sabiranju sa lijeve strane),
   \end{enumerate}
kada su formati matrica takvi da su svi gore navedeni proizvodi definisani.
\end{theorem}

\begin{example}
 Date su matrice
 \[A=\left(\begin{array}{ccr}3&1&-2\\0&2&0 \end{array}\right),\,
   B=\left(\begin{array}{rc}2&1\\0&4\\-1&0 \end{array}\right),\,
   C=\left(\begin{array}{r}-2\\0\\1 \end{array}\right),\,
   D=\left(\begin{array}{ccc}1&1&-2 \end{array}\right).\]
   Izra\v cunati
   \begin{inparaenum}[$(a)$]
      \item $A B;$
      \item $BA;$
      \item $AC;$
      \item $DB.$
   \end{inparaenum}\ \\\\
\noindent Rje\v senje: \\\\
Vrijedi
\begin{align*}
   (a)\quad  AB&=\left(\begin{array}{ccr}3&1&-2\\0&2&0 \end{array}\right)_{{\color{red}2\times3}}
                         \left(\begin{array}{rc}2&1\\0&4\\-1&0 \end{array}\right)_{{\color{red}3\times2}} \\
               &=\left(\begin{array}{cc} 3\cdot2+1\cdot0-2\cdot(-1)& 3\cdot1+1\cdot4-2\cdot0\\
               0\cdot2 +2\cdot0+0\cdot(-1)& 0\cdot1+2\cdot4+0\cdot0      \end{array}\right)_{{\color{red}2\times2}}
      =\left(\begin{array}{cc}8&7\\0&8 \end{array}\right)_{{\color{red}2\times2}}; \\
    (b)\quad   BA&=\left(\begin{array}{rc}2&1\\0&4\\-1&0 \end{array}\right) \left(\begin{array}{ccr}3&1&-2\\0&2&0 \end{array}\right)=
          \left(\begin{array}{rrr} 6&4&-4\\0&8&0\\-3&-1&2 \end{array}\right) ; \\
    (c) \quad AC&=\left(\begin{array}{rrr}3&1&-2\\0&2&0 \end{array}\right)\left(\begin{array}{r}-2\\0\\1 \end{array}\right)
          =\left(\begin{array}{r}-8\\0\end{array}\right); \\
    (d)\quad DB&=\left(\begin{array}{ccc}1&1&-2 \end{array}\right)\left(\begin{array}{rc}2&1\\0&4\\-1&0 \end{array}\right)=
           \left(4 \:\:\:5\right).
\end{align*}
\end{example}

\subsection{Trougaona, dijagonalna, skalarna, jedini\v cna i transponovana matrica}
Neke su matrice zbog specifi\v cnog rasporeda svojih elemenata ili vrijednosti tih elemenata, dobile i posebna imena. Ina\v ce ovakve matrice sre\'cemo \v cesto u drugim oblastima gdje se primjenjuju matrice, odnosno matri\v cni ra\v cun. Definicije nekih od takvih matrica su date u nastavku.\\

\begin{definition}[Trougaone matrice]
 Kvadratne matrice oblika
  { \footnotesize
 	\[\left(
 	    \begin{matrix}
 	    a_{11}&0&0&\ldots&0&0\\a_{21}&a_{22}&0&\ldots&0&0\vspace{.25cm}\\&\vdots&&& \vspace{.25cm}
 	    \\a_{n-1,1}&a_{n-1,2}&a_{n-1,3}&\ldots&a_{n-1,n-1}&0 \\a_{n1}&a_{n2}&a_{n3}&\ldots&a_{n,n-1}&a_{nn}
 	    \end{matrix}
 	  \right)
 	   \text{ i }
 	  \left(
 	   \begin{matrix}
 	   a_{11}&a_{12}&a_{13}&\ldots&a_{1,n-1}&a_{1n}\\
 	   0&a_{22}&a_{23}&\ldots&a_{2,n-1}&a_{2n}\vspace{.25cm}\\&\vdots&&&\vspace{.25cm}
 	   \\0&0&0&\ldots&a_{n-1,n-1}&a_{n-1,n}\\0&0&0&\ldots&0&a_{nn}
 	   \end{matrix}
 	  \right)
 	\]
 }

zovu se donja i gornja trougaona matrica, respektivno.\index{matrica!donja i gornja trougaona}
\end{definition}
Vidjeti Primjer \ref{primjerTrougaona}.\\
\begin{definition}[Dijagonalna matrica]\index{matrica!dijagonalna}
	Kvadratna matrica, \v ciji su svi elementi van glavne dijagonale nule zove se dijagonalna matrica, tj.
	\[\left(\begin{array}{cccccc} a_{11}&0&0&\ldots&0&0\\0&a_{22}&0&\ldots&0&0\vspace{.25cm}\\&\vdots&&& \vspace{.25cm}
	\\0&0&0&\ldots&a_{n-1,n-1}&0 \\0&0&0&\ldots&0&a_{nn} \end{array}\right).\]
\end{definition}
Vidjeti Primjer \ref{primjerDijagonalna}.

\newpage

\begin{example}\label{primjerTrougaona}
\[A=\left(\begin{array}{ccc}3&0&0\\2&2&0\\1&6&2 \end{array}\right),\:
  B=\left(\begin{array}{ccc}3&-2&1\\0&2&-4\\0&0&2 \end{array}\right),\]
matrica $A$ je donja trougaona, a matrica $B$ je gornja trougaona.
\end{example}

\begin{example}\label{primjerDijagonalna}
   Matrica
      \[\left(\begin{array}{crcc} 4&0&0&0\\ 0&-2&0&0\\0&0&3&0\\0&0&0&5\end{array}\right)\]
   je dijagonalna matrica.
\end{example}

\begin{definition}[Skalarna matrica]\index{matrica!skalarna}
 Dijagonalna matrica \v ciji su svi elementi na glavnoj dijagonali jednaki zove se skalarna matrica.
\end{definition}
Vidjeti Primjer \ref{primjerSkalarna}.
\begin{definition}[Jedini\v cna matrica]\index{matrica!jedini\v cna}
  Dijagonalna matrica \v ciji su svi elementi na glavnoj dijagonali jedinice, zove se jedini\v cna matrica i obilje\v zava se sa $I$ (ili $E$).
\end{definition}

\begin{example}\label{primjerSkalarna}
   Matrica
      \[A=\left(\begin{array}{crcc} 4&0&0&0\\ 0&4&0&0\\0&0&4&0\\0&0&0&4\end{array}\right)\]
   je skalarna matrica, dok je matrica
     \[I=\left(\begin{array}{cccc} 1&0&0&0\\ 0&1&0&0\\0&0&1&0\\0&0&0&1\end{array}\right)\]
     jedini\v cna matrica.
\end{example}

\index{matrica!transponovana}
\begin{definition}[Transponovana matrica]
  Transponovana matrica, matrice
  \[    A= \left( \begin{array}{cccc}
                a_{11}&a_{12}&\ldots&a_{1n}\\
                a_{21}&a_{22}&\ldots&a_{2n}\\
                &\vdots&&\\
                a_{m1}&a_{m2}&\ldots&a_{mn}
                \end{array}
         \right)
  \]
je matrica
  \[    A^T= \left( \begin{array}{cccc}
                a_{11}&a_{21}&\ldots&a_{m1}\\
                a_{12}&a_{22}&\ldots&a_{m2}\\
                &\vdots&&\\
                a_{1n}&a_{2n}&\ldots&a_{mn}
                \end{array}
         \right)
  \]
(za $\forall i=1,\ldots,m;$ $i$--ta vrsta u $A$ postaje $i$-ta kolona u $A^T$).
\end{definition}
Neke osobine transponovanja, odnosno operacija sa transponovanim matricama, date su u sljede\'coj teoremi.

\begin{theorem}\index{teorema! trasponovana matrica--osobine}
  Operacija transponovanja ima osobine
   \begin{enumerate}
     \item $\left(A^T\right)^T=A;$
     \item $(A+B)^T=A^T+B^T;$
     \item $(\alpha A)^T=\alpha A^T;$
     \item $(AB)^T=B^TA^T;$
     \item $(A_1A_2\cdots A_n)^T=A_n^T\cdots A_2^TA_1^T.$
   \end{enumerate}
\end{theorem}

\begin{example}
Date su matrice
 \[A=\left(\begin{array}{ccr}3&1&-2\\0&2&0 \end{array}\right),\,
   B=\left(\begin{array}{rc}2&1\\0&4 \end{array}\right).\]
   Odrediti $A^T$ i $B^T$.\\\\
\noindent Rje\v senje: \\\\
$A^T=\left(\begin{array}{rr} 3&0\\1&2\\-2&0 \end{array} \right)$ i
$B^T=\left(\begin{array}{rr}2&0\\1&4 \end{array} \right)$
\end{example}

\subsection{Z\lowercase{adaci}}\index{Zadaci za vje\v zbu!osnovni pojmovi o matricama}

\begin{enumerate}
\item Date su matrice \\
$A=\left(\begin{array}{cc} 3&-2\\5&-4\end{array}\right),\,
 B=\left(\begin{array}{cc} 3&-2\\5&-4\end{array}\right),\,
 C=\left(\begin{array}{ccc} 2&1&1\\3&0&1\end{array}\right),\,
 D=\left(\begin{array}{cc} 3&1\\2&1\\ 1&0\end{array}\right),\,$\\
$G=\left(\begin{array}{ccr} 5&8&-4\\6&9&-5\\4&7&3\end{array}\right)$ i
$F=\left(\begin{array}{ccc} 3&1&5\\4&-1&3\\6&9&5\end{array}\right).\vspace{.3cm}$\\
Izra\v cunati\\
 \begin{inparaenum}
    \item $A+B;\,$
    \item $2A+3B;\,$
    \item $3A-B;\,$
    \item $A\cdot B;\,$
    \item $B\cdot A;\,$
    \item $C\cdot D;\,$
    \item $D\cdot C;\,$
    \item $G\cdot F;\,$
    \item $F\cdot G.$
 \end{inparaenum}
\item Date su matrice\\
 $A=\left(\begin{array}{cccc}-1&2&5&-9\\-2&6&4&0\\1&1&0&4 \end{array} \right),\,
  B=\left(\begin{array}{ccc}1&4&7\\-2&0&10\\3&5&9\\0&0&-1 \end{array} \right),\,
  C=\left(\begin{array}{ccccc}-1&2&5&-9&6\\-3&-2&6&4&0 \end{array} \right),\,\\
  D=\left(\begin{array}{cc}4&7\\-2&10\\5&9\\0&-1\\3&3 \end{array} \right).$\\
Izra\v cunati\\
  \begin{inparaenum}
     \item $A\cdot B;\,$
     \item $B\cdot A;\,$
     \item $C\cdot D;\,$
     \item $D\cdot C.$
  \end{inparaenum}
\item Date su matrice
$A=\left(\begin{array}{cc}2& -1\\-3&7 \end{array} \right),\;
B=\left(\begin{array}{cc}4&0\\-2&1 \end{array} \right)$ i
$C=\left(\begin{array}{cc} 5&-3\\1&1\end{array} \right).$\\
Izra\v cunati \\
\begin{inparaenum}
   \item $2A;\,$
   \item $A+B;\,$
   \item $B-C;\,$
   \item $A^2-3AB^T+4C-2I;\,$
   \item $B^2C^T-2A+B.$
\end{inparaenum}
\item Izra\v cunati $AB-2A+B$ ako je
$A=\left(\begin{array}{rrr} 1&2&2\\2&1&2\\1&2&3\end{array}\right)$ i
$B=\left(\begin{array}{rrr} 4&1&1\\-4&2&0\\1&2&1\end{array}\right).$
\item Izra\v cunati $A^2-3B$ ako je
$A=\left(\begin{array}{rrr} 2&1&0\\1&1&2\\-1&2&1\end{array}\right)$ i
$B=\left(\begin{array}{rrr} 3&1&-2\\3&-2&4\\-3&5&-1\end{array}\right).$
\end{enumerate}

\newpage
\section[D\lowercase{eterminante}]{Determinante}\index{determinante}

Pojam determinante usko je povezan sa kvadratnim  matricama. Determinante su skalarne vrijednosti pridru\v zene kvadratnim matricama, isto tako nezaobilazan su alat u nekim metodama rje\v savanja sistema algebarskih  linearnih jedna\v cina. Zbog tih, a i drugih razloga, va\v zno je nau\v citi \v sta su determinante i kako se ra\v cuna sa njima.

Posmatrajmo sistem od dvije linearne algebarske jedna\v cine
 \[
  \begin{cases}
     a_{11}x+a_{12}y=b_1\\
     a_{21}x+a_{22}y=b_2.
  \end{cases}
  \]
Rje\v senja mo\v zemo napisati u obliku
\begin{align*}
  x&=\frac{a_{22}b_1-a_{12}b_2}{a_{11}a_{22}-a_{21}a_{12}}\\
  y&=\frac{a_{11}b_1-a_{21}b_2}{a_{11}a_{22}-a_{21}a_{12}},
\end{align*}
za $a_{11}a_{22}-a_{21}a_{12}\neq 0.$
Koeficijente koji se nalaze u nazivniku mo\v zemo izdvojiti u kvadratnu matricu drugog reda
\[\left(\begin{array}{cc} a_{11}&a_{12}\\a_{21}&a_{22}\end{array}\right),\]
ovoj matrici mo\v zemo pridru\v ziti broj $ a_{11}a_{22}-a_{21}a_{12}$ ili zapisan u obliku kvadratne sheme
\[ \left|\begin{array}{cc} a_{11}&a_{12}\\a_{21}&a_{22}\end{array}\right|.\]

Generalno, svakoj kvadratnoj matrici $A=(a_{ij})$ reda $n$ mo\v zemo pridru\v ziti broj ili odgovaraju\' ci izraz. U tu svrhu, posmatrajmo najprije sve permutacije skupa $\{1,2, \ldots , n\}.$
Ure\dj eni par $(k_i, k_j)$ jedne permutacije $p,$ gdje je $i<j,$ zovemo inverzija permutacije $p$ ako je $k_i>k_j.$ Ozna\v cimo broj inverzija permutacije $p$ sa $i(p).$ Na primjer, permutacija $p=(312)$ skupa
$\{1,\, 2,\, 3 \}$ ima $2$ inverzije, a to su parovi $(3,1)$ i $(3,2).$
\begin{definition}
 Neka je data kvadratna matrica $A=(a_{ij})$ reda $n,$ \v ciji su elementi $a_{ij},\:i,j=1,2,\ldots,n$ realni ili kompleksni brojevi. Pod determinantom matrice $A$ podrazumijeva se broj
 \[ \det A=\sum_{p \in S_n} (-1)^{i(p)} a_{1k_1}a_{2k_2}\cdots a_{nk_n},\]
 gdje se sumira po svim permutacijama $p=(k_1 \, k_2 \, \ldots , k_n  )$ skupa $\{1,2, \ldots , n\}.$  Taj skup permutacija ozna\v cavamo sa $S_n$  i zovemo simetri\v cna grupa.
\end{definition}
Suma u ovoj definiciji se sastoji od $n!$ sabiraka jer toliko ima svih permutacija $n$--\v clanog skupa (card $S_n=n!$).
Determinanta matrice  $A$ (determinanta $n$--tog reda) mo\v ze se, osim sa $\det A,$ obilje\v ziti i sa $|A|$ ili
\[\left|\begin{array}{cccc}
                   a_{11}&a_{12}&\ldots&a_{1n}\\
                   a_{21}&a_{22}&\ldots&a_{2n}\\
                   &\vdots&&\\
                   a_{n1}&a_{n2}&\ldots&a_{nn}
\end{array}
\right|.\]

\begin{remark}
  U literaturi su dostupne i definicije koje su prakti\v cnije sa stanovi\v sta samog ra\v cunanja vrijednosti determinante. Jedna takva definicija bi\' ce navedena u nastavku ove se-\\kcije.
\end{remark}
\subsection{Osobine determinanti}\index{determinante!osobine}
Neka su date kvadratne matrice $A=(a_{ij})$ i $B=(b_{ij})$ proizvoljnog reda. Za determinante $\det A$ i $\det B$ vrijede sljede\' ce osobine.
\begin{enumerate}
  \item $\det A=\det A^T.$ Npr.,
  \[\left|\begin{array}{cc}1&2\\3&4\end{array}\right|=1\cdot4-2\cdot3=-2;\quad
    \left|\begin{array}{cc}1&3\\2&4\end{array}\right|=1\cdot4-3\cdot2=-2,\text{ vidimo da je } \det A=\det A^T.
  \]
  \item Ako su u determinanti $\det A$ elementi jedne vrste (ili kolone) jednaki sa elementima  druge vrste (ili kolone), tada je $\det A=0.$ Npr.,
  \[\left|\begin{array}{cc}1&2\\1&2\end{array}\right|=1\cdot2-1\cdot2=0.\]
  \item Determinanta $\det A$  pomno\v zi se sa brojem $\lambda$ tako \v sto se svaki element samo {\bfseries \color{red}jedne vrste ili kolone} te determinante pomno\v zi sa tim brojem. Npr.,
  \[-5\cdot\left|\begin{array}{cc}1&2\\3&4\end{array}\right|
    =\left|\begin{array}{rr}-5\cdot1&-5\cdot2\\3&4\end{array}\right|
    =\left|\begin{array}{cc}-5\cdot1&2\\-5\cdot3&4\end{array}\right|=10.\]
  \item Ako su elementi jedne  vrste (ili kolone) u determinanti $\det A$ proporcionalni elementima druge vrste (ili kolone), tada je $\det A=0.$ Npr.,
  \[\left|\begin{array}{cc}1&2\\3&6\end{array}\right|=1\cdot6-2\cdot3=0.\]
  \item Ako je svaki element $k$--te vrste  determinante $\det A$ zbir dva broja
        (tj. $a_{kj}=a_{kj}^{(1)}+a_{kj}^{(2)}$),  onda je determinanta $\det A$ jednaka zbiru dvije determinante, istog reda pri \v cemu su elementi $k$--te vrste u jednoj determinanti prvi sabirci ($a_{kj}^{(1)}$), a elementi $k$--te vrste druge determinante drugi sabirci ($a_{kj}^{(2)}$). Ostali elementi u obje determinante jednaki su odgovaraju\' cim elementima determinante $\det A$ (ista osobina vrijedi i za kolone). Npr.,
        \begin{align*}
        &  \det A=\left|\begin{array}{cc}2+3&2\\1+2&4\end{array}\right|=14\\
        &  \det A^{(1)}=\left|\begin{array}{cc}2&2\\1&4\end{array}\right|=6\\
        &  \det A^{(2)}=\left|\begin{array}{cc}3&2\\2&4\end{array}\right|=8.
        \end{align*}
        Vidimo da je $\det A=\det A^{(1)}+\det A^{(2)}.$
\item Determinanta ne mijenja vrijednost ako se elementima jedne vrste (ili kolone) dodaju odgovaraju\' ci elementi druge vrste (ili kolone)  pomno\v zeni istim brojem. Npr., ako kod determinante drugu kolonu pomno\v zimo sa 3 i dodamo prvoj koloni, dobi\'cemo
 \begin{align*}
         \left|\begin{array}{cc}5&2\\3&4\end{array}\right|&=14\\
         \left|\begin{array}{cc}5+3\cdot2&2\\3+3\cdot4&4\end{array}\right|&=
         \left|\begin{array}{cc}11& 2\\15&4\end{array}\right|=44-30=14.
  \end{align*}
 \item Ako u determinanti $\det A$ dvije vrste (ili kolone) me\dj usobno promijene mjesta determinanta $\det A$ mijenja znak. Npr.,
   \begin{align*}
         \left|\begin{array}{cc}5&2\\3&4\end{array}\right|&=14\\
         \left|\begin{array}{cc}3&4\\5&2\end{array}\right|&=6-20=-14.
  \end{align*}
\item   $\det (AB)=\det A\det B,$ za kvadratne matrice istog reda. Npr.,
      \begin{align*}
        A&= \left(\begin{array}{cc}1&2\\-1&4\end{array}\right)\\
        B&= \left(\begin{array}{cc}3&0\\2&-5\end{array}\right)\\
        AB&= \left(\begin{array}{cc}7&-10\\5&-20\end{array}\right)\\
        \det A&= \left|\begin{array}{cc}1&2\\-1&4\end{array}\right|=6\\
        \det B&= \left|\begin{array}{cc}3&0\\2&-5\end{array}\right|= -15 \\
        \det AB&= \left|\begin{array}{cc}7&-10\\5&-20\end{array}\right|=-90 .
      \end{align*}
      Vidimo da je $\det (AB)=\det A\det B.$
\end{enumerate}

\subsection{Ra\v cunanje determinanti}
\index{determinante!ra\v cunanje vrijednosti}

Determinanta prvog reda se sastoji samo od jednog elementa i ona je jednaka tom elementu.\\
Za ra\v cunanje determinanti drugog i tre\' ceg reda koristimo sljede\' ce formule
\begin{align}
   \left|\begin{array}{cc}a_{11}&a_{12}\\a_{21}&a_{22}\end{array}\right|&=a_{11}a_{22}-a_{12}a_{21},\label{determinanta1}\\\nonumber \\
   \left|\begin{array}{ccc}a_{11}&a_{12}&a_{13}\\a_{21}&a_{22}&a_{23}\\a_{31}&a_{32}&a_{33}\end{array}\right|&=
   a_{11}a_{22}a_{33}+ a_{12}a_{23}a_{31}+a_{13}a_{21}a_{32} \nonumber  \\
   &\quad\quad  -\left( a_{31}a_{22}a_{13}+ a_{32}a_{23}a_{11}+a_{33}a_{21}a_{12}\right)  .\label{determinanta2}
\end{align}
Za lak\v se ra\v cunanje determinante tre\' ceg reda mo\v zemo koristiti tzv. Sarrusovo\footnote{ Pierre Fr\'ed\'eric Sarrus (10.mart 1798--20.novembar 1861. godine), bio je francuski matemati\v car} pravilo.\\
Ovo pravilo se sastoji u tome da se determinanti sa desne strane dopi\v su prve dvije kolone. Nakon toga, mno\v zimo trojke brojeva u pravcu
glavne dijagonale i saberemo ta tri umno\v ska, a onda od toga oduzmemo sumu umno\v zaka trojki brojeva koji se nalaze u pravcu sporedne dijagonale.
\begin{displaymath}
\left|\begin{array}{ccc}
a_{11}&a_{12}&a_{13}\\
a_{21}&a_{22}&a_{23}\\
a_{31}&a_{32}&a_{33}
\end{array}\right| \,
\begin{array}{cc}
a_{11}&a_{12}\\
a_{21}&a_{22}\\
a_{31}&a_{32}
\end{array}=a_{11}a_{22}a_{33}+ a_{12}a_{23}a_{31}+a_{13}a_{21}a_{32}-
\end{displaymath}
\begin{flushright}
$-\left( a_{31}a_{22}a_{13}+ a_{32}a_{23}a_{11}+a_{33}a_{21}a_{12} \right).$
\end{flushright}

\begin{remark}
Umjesto opisanog dopisivanja kolona desno od determinante, mogli smo obaviti dopisivanje prve dvije vrste ispod determinante, a dalji postupak je isti.
\end{remark}
\begin{example}
 Izra\v cunati vrijednost determinati.

  \begin{enumerate}
      \item $\left|\begin{array}{cc}7&-4\\3&4\end{array}\right|=7\cdot4-(-4\cdot3)=28+12=40;$
      \item $\left|\begin{array}{rrr}3&2&-1\\1&2&4\\0&6&-2\end{array}\right|
       =3\cdot2\cdot(-2)+2\cdot4\cdot0+(-1)\cdot 1\cdot 6-(0\cdot2\cdot(-1)+6\cdot4\cdot3+(-2)\cdot 1\cdot 2)\\=-86.$
  \end{enumerate}
\end{example}

\paragraph{Razvijanje determinante po elementima neke vrste ili kolone} Za ra\v cunanje determinanti reda ve\'ceg od tri postoje formule analogne formulama za ra\v cunanje determinati drugog i tre\'ceg reda koje su date u \eqref{determinanta1} i \eqref{determinanta2}. Me\dj utim, te su formule dosta slo\v zene i neprakti\v cne su. Mnogo bolje rje\v senje je razvijanje determinate po vrsti ili koloni. Na ovaj na\v cin determinantu  $n$--tog reda predstavljamo preko $n$ odgovaraju\'cih determinanti reda $n-1,$ odnosno ovim postupkom smanjujemo red determinante za 1. Postupak nastavlja-\\mo dok ne dobijemo determinate odgovaraju\'ceg reda (1., 2. ili 3.).  Navedeni postupak dat je u teoremi koja slijedi, a prije toga definisana su dva pojma koja se koriste u ovom postupku.

Posmatrajmo ponovo kvadratnu matricu $A=(a_{ij})$ reda $n,$ njena determinanta je
\[\det A=\left| \begin{array}{cccc}
                a_{11}&a_{12}&\ldots&a_{1n}\\
                a_{21}&a_{22}&\ldots&a_{2n}\\
                &\vdots&&\\
                a_{n1}&a_{n2}&\ldots&a_{nn}
                 \end{array}
      \right|.
\]
\index{determinante!minor}
\begin{definition}[Minor]
   Determinanta koja se dobije iz $\det A$ odbacivanjem $i$--te vrste i $j$--te kolone, naziva se minor elementa $a_{ij}$ i obilje\v zava se sa $M_{ij}.$
\end{definition}
\begin{definition}[Algebarski kofaktor]\index{determinante!algebarski kofaktor}
  Broj $(-1)^{i+j}M_{ij}$ naziva se algebarski komplement (ili kofaktor) elementa $a_{ij}$ i ozna\v cava se sa $A_{ij}.$
\end{definition}

Sljede\' ca teorema daje pravilo razlaganja determinante po elementima proizvoljne vrste ili kolone.\index{determinante!razlaganje po vrsti ili koloni--Laplaceov razvoj}
\begin{theorem}[Laplaceova$^{1}$]\index{teorema! Laplaceov razvoj determinante}
  Ako je  $\det A$ reda $n,$ onda je
  \begin{enumerate}[$\bullet$]
     \item $\det A=a_{i1}A_{i1}+a_{i2}A_{i2}+\ldots+a_{in}A_{in},$  (Laplaceov razvoj po $i$--toj vrsti);
     \item $\det A=a_{1j}A_{1j}+a_{2j}A_{2j}+\ldots+a_{nj}A_{nj},$ (Laplaceov razvoj po $j$--toj koloni);
  \end{enumerate}
  za $\forall i,\,j\in\{1,2,\ldots,n\}.$
\end{theorem}
\footnotetext[1]{ Pierre--Simon, marquis de Laplace (23.mart 1749.--05.mart 1827. godine) bio je francuski nau\v cnik \v ciji rad je imao zna\v cajan uticaj u razvoju matematike, fizike, filozofije i raznih in\v zinjerskih disciplina.   }
\begin{example}
  Determinatu
  \[\det A=\left|\begin{array}{rrr}3&2&-1\\1&2&4\\0&6&-2\end{array}\right|,\]
  razviti po elementima druge vrste i prve kolone.\\\\
\noindent Rje\v senje: \\\\
Razvoj po elementima druge vrste
 \begin{align*}
   \left|\begin{array}{rrr}3&2&-1\\1&2&4\\0&6&-2\end{array}\right|&
   =1\cdot(-1)^{2+1}\left|\begin{array}{rr}2&-1\\6&-2\end{array}\right|
     +2\cdot(-1)^{2+2}\left|\begin{array}{rr}3&-1\\0&-2\end{array}\right| \\
   &\quad\quad   +4\cdot(-1)^{2+3}\left|\begin{array}{rr}3&2\\0&6\end{array}\right|\\
   &=-1(-4+6)+2(-6)-4\cdot 18=-86.
 \end{align*}
Razvoj po elementima prva kolone
 \begin{align*}
   \left|\begin{array}{rrr}3&2&-1\\1&2&4\\0&6&-2\end{array}\right|&
   =3\cdot(-1)^{1+1}\left|\begin{array}{rr}2&4\\6&-2\end{array}\right|
     +1\cdot(-1)^{2+1}\left|\begin{array}{rr}2&-1\\6&-2\end{array}\right|\\
   &\qquad\quad  +0\cdot(-1)^{3+1}\left|\begin{array}{rr}2&-1\\2&4\end{array}\right|\\
   &=3(-4-24)-(-4+6)=-86.
 \end{align*}
\end{example}
U posljednjem primjeru vidimo da je korisno ako je neki element u vrsti ili koloni, po kojoj vr\v simo razvoj, jednak nuli. Ovo \' ce  biti iskori\v steno u sljede\' cem primjeru.\\
\begin{example}
  Izra\v cunati vrijednost determinante
  \[\det A=\left|\begin{array}{rrrr}
                       2&1&2&1\\
                       2&-3&1&-3\\
                       4&2&2&2\\
                       -2&4&-1&5
                 \end{array}
          \right| \]
\noindent Rje\v senje:\\\\
Da bi pojednostavili razvijanje determinante, od prve kolone oduzmimo tre\' cu kolonu prethodno pomno\v zenu sa 2, te od \v cetvrte kolone oduzmino drugu kolonu, vrijedi
\begin{align*}
\begin{blockarray}{ccccc}
\begin{block}{|rrrr|r}
  2&1&2&1&\text{I kol-$2\times$III kol }\\
  2&-3&1&-3&\\
  4&2&2&2&\\
  -2&4&-1&5&\\
\end{block}
\end{blockarray}&=
\begin{blockarray}{ccccc}
\begin{block}{|rrrr|r}
  -2&1&2&1&\text{IV kol-II kol }\\
  0&-3&1&-3&\\
  0&2&2&2&\\
  0&4&-1&5&\\
\end{block}
\end{blockarray}\\
&=
\begin{blockarray}{ccccc}
\begin{block}{|rrrr|r}
  -2&1&2&0&\\
  0&-3&1&0&\\
  0&2&2&0&\\
  0&4&-1&1&\\
\end{block}
\end{blockarray}
 \end{align*}
Razvijmo sada determinatnu po elementima  \v cetvrte kolone
\[ =(-1)^{4+4}\left|\begin{array}{rrr}
                                    -2&1&2\\0&-3&1\\0&2&2
                                    \end{array}
                              \right|.
\]
Razvijmo posljednju determinantu po elementima prve kolone
\[ =-2\cdot  (-1)^{1+1}\left|\begin{array}{rr}
              -3&1\\2&2
          \end{array}
    \right|=-2\cdot (-6-2)=16.
\]
\end{example}
\newpage

\begin{example}
 Dokazati \[ \left| \begin{array}{ccc}ax&a^2+x^2&1\\ay&a^2+y^2&1\\az&a^2+z^2&1 \end{array}\right|=a(x-y)(x-z)(z-y).\]
\noindent Rje\v senje:
\begin{align*}
& \begin{blockarray}{cccc}
    \begin{block}{|ccc|l}
      ax & a^2+x^2&1&\\
      ay&a^2+y^2&1&\text{ IIv-Iv}\\
      az&a^2+z^2&1&\text{ IIIv-Iv}\\
    \end{block}
 \end{blockarray}\\
 &\quad\quad=
 \begin{blockarray}{cccc}
    \begin{block}{|ccc|l}
      ax & a^2+x^2&1&\text{ izvu\v cemo $a$ iz Ik}\\
      a(y-x)&(y-x)(y+x)&0&\text{ izvu\v cemo   $y-x$ iz IIv}\\
      a(z-x)&(z-x)(z+x)&0&\text{ izvu\v cemo   $z-x$ iz IIIv}\\
    \end{block}
 \end{blockarray} \\
 &\quad\quad=a(y-x)(z-x)
 \begin{blockarray}{cccc}
    \begin{block}{|ccc|l}
      x & a^2+x^2&1&\\
      1&y+x&0&\text{razvijemo determinantu po   IIIk}\\
      1&z+x&0&\\
    \end{block}
 \end{blockarray}\\
 &\quad\quad=
 a(y-x)(z-x)\cdot1\cdot(-1)^{1+3}
            \begin{blockarray}{cc}
                \begin{block}{|cc|}
                         1&y+x\\1&z+x\\
                \end{block}
            \end{blockarray}\\
&\quad\quad=a(y-x)(z-x)[z+x-(y+x)]\\
&\quad\quad=a(y-x)(z-x)(z-y)=a[-(x-y)][-(x-z)](z-y)\\
&\quad\quad=a(x-y)(x-z)(z-y).
\end{align*}
\end{example}

\begin{example}
 Dokazati \[\begin{blockarray}{cccc}
               \begin{block}{|cccc|}
                   a&b&c&d\\
                   a&-b&-c&-d\\
                   a&b&-c&-d\\
                   a&b&c&-d\\
               \end{block}
              \end{blockarray}=-8abcd.
          \]
\noindent Rje\v senje:
   \begin{align*}
&    \begin{blockarray}{ccccl}
               \begin{block}{|rrrr|l}
                   a&b&c&d&\text{izvu\v cemo $a$ iz Ik} \\
                   a&-b&-c&-d&\text{izvu\v cemo $b$ iz IIk}\\
                   a&b&-c&-d&\text{izvu\v cemo $c$ iz IIIk}\\
                   a&b&c&-d&\text{izvu\v cemo $d$ iz IVk}\\
               \end{block}
              \end{blockarray}\\
         &\quad\quad =abcd\begin{blockarray}{ccccl}
               \begin{block}{|rrrr|l}
                   1&1&1&1& \\
                   1&-1&-1&-1&\text{IIv+Iv}\\
                   1&1&-1&-1&\\
                   1&1&1&-1&\\
               \end{block}
              \end{blockarray}  \\
         &\quad\quad=abcd\begin{blockarray}{ccccl}
               \begin{block}{|rrrr|l}
                   1&1&1&1& \\
                   2&0&0&0&\text{ razvijemo determinantu po elementima IIv}\\
                   1&1&-1&-1&\\
                   1&1&1&-1&\\
               \end{block}
              \end{blockarray} \\
         &\quad\quad=abcd\cdot2\cdot(-1)^{2+1}
         \begin{blockarray}{cccl}
               \begin{block}{|rrr|l}
                   1&1&1& \\
                   1&-1&-1&\text{ IIv+Iv}\\
                   1&1&-1&\\
               \end{block}
         \end{blockarray} \\
       &\quad\quad=-2abcd
         \begin{blockarray}{cccl}
               \begin{block}{|rrr|l}
                   1&1&1& \\
                   2&0&0&\text{ razvijemo determinantu po elementima IIv}\\
                   1&1&-1&\\
               \end{block}
         \end{blockarray} \\
      &\quad\quad=-2abcd\cdot2\cdot(-1)^{2+1}
           \begin{blockarray}{ccl}
               \begin{block}{|rr|l}
                   1&1& \\
                   1&-1&\\
               \end{block}
         \end{blockarray} \\
     &\quad\quad=4abcd(-1-1)=-8abcd.
   \end{align*}
\end{example}

\subsection{Z\lowercase{adaci}}\index{Zadaci za vje\v zbu!determinante}

\begin{enumerate}

\item Izra\v cunati vrijednost determinati

\begin{inparaenum}
\item $\left|\begin{array}{cc}2&4\\-3&9 \end{array} \right|;\:$
\item $\left|\begin{array}{cc}2+2i&3-i\\\\-3+4i& -1-i \end{array} \right|;\:$
\item $\left|\begin{array}{cc}\dfrac{x^2+1}{1-x^2}&\dfrac{2x}{1-x^2}\\\dfrac{2x}{1-x^2}&\dfrac{x^2+1}{1-x^2} \end{array} \right|;\:$

\item $\left|\begin{array}{cc}\sin x+\sin y&\cos x+\cos y\\\\\cos y-\cos x&\sin x-\sin y \end{array} \right|;\:$
\item $\left|\begin{array}{cc}1&\log_a b\\\\\log_b a&1 \end{array} \right|;\:$

\item $\left|\begin{array}{ccc}
5&\:6&\:3\\0&\:1&\:2\\7&\:4&\:5\\	
\end{array}\right|;\: $
\item $\left|\begin{array}{ccc}1&2&0\\-3&2&9\\-1&2&1 \end{array} \right|;\:$
\item $\left|\begin{array}{ccc}a&a&a\\-a&a&x\\-a&-a&x \end{array} \right|;\:$
\item $\left|\begin{array}{ccc}1&i&1+i\\-i&1&0\\1-i&0&1 \end{array} \right|;\:$

\item  $\left|\begin{array}{rrrr}
3&1&2&3\\4&-1&2&4\\1&-1&1&1\\4&-1&1&5\\	
\end{array}\right|;\:$
\item $\left|\begin{array}{cccc}
1&\:1&\:3&\:4\\2&\:0&\:0&\:8\\3&\:0&\:0&\:2\\4&\:4&\:7&\:5\\		
\end{array}\right|;\:$
\item $\left|\begin{array}{cccc}
3&\:1&\:2&\:3\\4&\:-1&\:2&\:4\\1&\:-1&\:1&\:1\\4&\:-1&\:1&\:5\\		
\end{array}\right|;\:$

\item $\left|\begin{array}{ccccc}
7&\:2&\:1&\:3&\:4\\1&\:0&\:2&\:0&\:3\\3&\:0&\:4&\:0&\:7 \\6&\:3&\:2&\:4&\:5\\5&\:1&\:2&\:2&\:3		
\end{array}\right|;\;$
\item $\left|\begin{array}{ccccc}
x&\:0&\:-1&\:1&\:0\\1&\:x&\:-1&\:1&\:0\\1&\:0&\:x-1&\:0&\:1 \\0&\:1&\:-1&\:x&\:1\\0&\:1&\:-1&\:0&\:x		
\end{array}\right|.$

\end{inparaenum}

\item Pokazati da je
   \begin{inparaenum}

   \item $\left|\begin{array}{ccc}
         b^2+c^2&\:ab&\:ca\\ ab&\:c^2+a^2&\:bc\\ ca&\:bc&\:a^2+b^2			
        \end{array}\right|=4a^2b^2c^2;\:$
   \item $\left|\begin{array}{ccc}1&\:a&\:a^{2}+a^3\\ 1&\:a^2&\: a^3+a\\ 1&\:a^3&\:a+a^2\\			
         \end{array}\right|=0;\vspace{.2cm}$

   \item $\left|\begin{array}{ccc}1&1&1\\1&a&1\\1&1&b\end{array} \right|=(1-a)(1-b);\:$
   \item $\left|\begin{array}{ccc}1&\:a&\:bc\\ 1&\:b&\:ca\\ 1&\:c&\:ab			
         \end{array}\right|=(b-c)(b-a)(a-c); \vspace{.2cm}$

   \item $\left|\begin{array}{ccc}
          1&\:a&\:a^2\\ 1&\:b&\:b^2\\ 1&\:c&\:c^2			
          \end{array}\right|=(b-a)(c-a)(c-b);\:$

  \item $\left|\begin{array}{ccc}
        ax&\:a^2+x^2&\:1\\ ay&\:a^2+y^2&\:1\\ az&\:a^2+z^2&\:1			
        \end{array}\right|=a(x-y)(x-z)(z-y);\vspace{.2cm}$

  \item $\left|\begin{array}{rrr}
               	1&bc&b+c\\1&ac&a+c\\1&ab&a+b\\
            \end{array}\right|=(a-b)(b-c)(c-a);\vspace{.2cm}$
  \item  $\left|\begin{array}{ccc}a-b & 2a & 2a\\2b & b-a & 2b\\a-b & 2a & a-b \end{array}\right|=(a+b)^3;\vspace{.2cm}$

  \item $\left| \begin{array}{ccc}a+b&-a&-b\\-b&b+c&-c\\-a&-c&c+a \end{array}\right|=0;\:$
  \item $\left|\begin{array}{cccc}
        -x&\:y&\:z&\:1\\x&\:-y&\:z&\:1\\x&\:y&\:-z&\:1\\x&\:y&\:z&\:-1		
        \end{array}\right|=-8xyz;\vspace{.2cm}$

  \item $\left|\begin{array}{cccc}
        a&\:b&\:c&\:d\\a&\:-b&\:-c&\:-d\\a&\:b&\:-c&\:-d\\a&\:b&\:c&\:-d		
        \end{array}\right|=-8abcd;\vspace{.2cm}$
  \item $\left|\begin{array}{cccc}
       a&\:a&\:a&\:a\\a&\:b&\:b&\:b\\a&\:b&\:c&\:c\\a&\:b&\:c&\:d 	
       \end{array}\right|=-a(a-b)(c-b)(d-c).$

  \end{inparaenum}

\item Izra\v cunati \\
     \begin{inparaenum}
         \item $\left|\begin{array}{ccc}
               a^{-4}&\:a^{-3}&\:a^{-2}\\ a^{-1}&\:1&\: a\\ a^2&\:a^3&\:a^4\\			
              \end{array}\right|\:(=0);\:$
         \item  $\left|\begin{array}{cccc}
               a&\:b&\:a&\:1\\b&\:a&\:b&\:1\\a&\:-a&\:b&\:1\\b&\:-b&\:a&\:1  		
              \end{array}\right|\:(=2(a+b)(b-a)^2);\:\vspace{.2cm}$

          \item  $\left|\begin{array}{rrr}
                 1&1&1\\1&z&z^2\\1&z^2&z\\	
                \end{array}\right|$, ako je $z=\cos\dfrac{4\pi}{3}+i\sin\dfrac{4\pi}{3};\:$

          \item  $\left|\begin{array}{rrr}
                 1&1&z\\1&1&z^2\\z^2&z&1\\	
                \end{array}\right|$, ako je $z=\cos\dfrac{2\pi}{3}+i\sin\dfrac{2\pi}{3};\:$

          \item $\left|\begin{array}{ccc}
                \sin 2x&\:-\cos 2x&\:1\\ \sin x&\:-\cos x&\:\cos x\\ \cos x&\:\sin x&\:\sin x\\			
                  \end{array}\right|;\;(=0)$
          \item $\left|\begin{array}{ccc}
                1+\cos x&\:1+\sin x&\:1\\ 1-\sin x&\:1+\cos x&\:1\\ 1&\:1&1\:\\			
                  \end{array}\right|\;(=1).$

     \end{inparaenum}

\item  Rije\v siti jedna\v cine \\
\begin{inparaenum}
\item $\left|\begin{array}{ccc}x&2&3\\3&1&2\\1&3&4 \end{array} \right|=0;\:$
\item $\left|\begin{array}{ccc}x&1&1\\1&x&1\\1&1&x\end{array} \right|=0;\:$

\item $\left|\begin{array}{ccc}
\log_c x & \log_c x-n\\ \log_c x-m & \log_c x			
\end{array}\right|=0,\:(0<c\neq 1);\vspace{.2cm}$

\item $\left|\begin{array}{cccc}
1&\:1&\:2&\:3\\1&\:2-x^2&\:2&\:3\\2&\:3&\:1&\:5\\2&\:3&\:1&\:9-x^2\\	
\end{array}\right|=0;\:$
\item $\left|\begin{array}{rrr}
	x-3&x+2&x-1\\x+2&x-4&x\\x-1&x+4&x-5\\
\end{array}\right|=0;\vspace{.2cm}$

\item $\left|\begin{array}{ccc}
\sin\left(x+\tfrac{\pi}{4}\right)&\sin x&\cos x\\
\sin\left(x+\tfrac{\pi}{4}\right)&\cos x&\sin x\\1&a&1-a\\	
\end{array}\right|=\dfrac{\sqrt{2}-2}{4};\:$
\item $\left|\begin{array}{cccc}1&1&2&3\\1&2-x^2&2&3\\2&3&1&5\\2&3&1&1-x^2
\end{array} \right|=0.$
\end{inparaenum}	

\item Rije\v siti nejedna\v cine\\
  \begin{inparaenum}
    \item $\left|
              \begin{array}{ccc}x&-1&0\\5&-1&-6\\-1&0&x
              \end{array}
            \right|\geqslant 0;$
     \item  $\left|
              \begin{array}{ccc}1&0&x\\0&1&1\\1&x&0
              \end{array}
            \right
             |\leqslant
             \left|
              \begin{array}{ccc}x&3&x\\2&1&3\\1&x&1
              \end{array}
            \right|.
            $
  \end{inparaenum}
\end{enumerate}

\newpage

\section{I\lowercase{nverzna matrica}} \index{matrica! adjungovana}\index{matrica! kofaktora}
Neka je $a\neq 0$ realan broj i neka je $"\cdot"$ uobi\v cajeno mno\v zenje na skupu realnih brojeva. Za broj $a^{-1}$ ili $\frac{1}{a}$ ka\v zemo da je inverzni broj broja $a$ ako je ispunjen uslov $a\cdot a^{-1}=1$ ili $a\cdot\frac{1}{a}=1.$   \v Cesto je potrebno za kvadratnu matricu $A$ odrediti matricu $A^{-1}$ tako da je ispunjen uslov $A\cdot A^{-1}=I$ (ili $A^{-1}\cdot A=I$) gdje je $I$ jedini\v cna matrica, dok je  $"\cdot"$ sada operacija mno\v zenja matrica.  U nastavku \'cemo vidjeti koje uslove mora ispuniti matrica $A$ da bi postojala matrica $A^{-1}$ sa navedenim osobinama, kako se onda ra\v cuna matrica $A^{-1}$ i gdje se primjenjuje.

Prvi korak u odre\dj ivanju tra\v zene matrice je adjungirana matrica, slijedi njena definicija.

\begin{definition}[Adjungirana matrica]
  Neka je $A=(a_{ij})$ kvadratna matrica reda $n$, i $A_{ij}$ je algebarski kofaktor elementa $a_{ij}.$ Tada se matrica
  \[\adj A=\left(\begin{array}{cccc}
                   A_{11}&A_{21}&\ldots&A_{n1}\\
                   A_{12}&A_{22}&\ldots&A_{n2}\\
                   &\vdots&&\\
                   A_{1n}&A_{2n}&\ldots&A_{nn}
\end{array}   \right)\]
naziva adjungirana matrica matrice $A$. Determinanta matrice $\adj A$ zove se adjungirana determinanta determinante matrice $A$.
\end{definition}
Vrijedi sljede\' ca teorema u kojoj su date neke osobine $\adj A$ i $\det(\adj A),$ koje koristimo u nastavku.\\
\begin{theorem}\index{teorema! kvadratna matrica--osobine}\label{toremaAdjungirana1}
Neka je $A$ kvadratna matrica reda $n$, tada vrijedi
 \begin{enumerate}[$1$]
   \item $A\cdot\adj A=(\adj A)\cdot A=(\det A)\cdot I;$
   \item $\det A\cdot \det(\adj A)=(\det A)^n;$
   \item $\det (\adj A)=(\det A)^{n-1}.$
 \end{enumerate}
\end{theorem}

\begin{example}
  Za date matrice odrediti adjungirane matrice
    \begin{enumerate}[$(a)$]
      \item $A=\left(\begin{array}{rr}3&-5\\4&6\end{array}       \right);$
      \item $A=\left(\begin{array}{rrr}1&2&-5\\0&2&1\\1&1&3\end{array}\right).$
     \end{enumerate}\ \\
\noindent Rje\v senje: \\\\
Ra\v cunanje adjungirane matrice podijeli\' cemo u dva koraka.
\begin{enumerate}[I {korak:}]
 \item Izra\v cunamo algebarske kofaktore $A_{ij}$ i formiramo matricu kofaktora $\cof A =(A_{ij});$
 \item Zatim transponujemo matricu kofaktora $\cof A$ i dobijemo adjungiranu matricu $\adj A$ matrice $A,$ tj. $\adj A=(\cof A)^T.$
\end{enumerate}

\begin{enumerate}[$(a)$]
  \item Izra\v cunamo prvo kofaktore
  \begin{align*}
        A_{11}&=(-1)^{1+1}\cdot6=6,&&A_{12}=(-1)^{1+2}\cdot4=-4\\
        A_{21}&=(-1)^{2+1}\cdot(-5)=5,&&A_{22}=(-1)^{2+2}\cdot3=3
  \end{align*}
  formiramo matricu kofaktora $\cof A$ i transponujemo je
  \begin{align*}
     \cof A&=\left(\begin{array}{rr}6&-4\\5&3\end{array}       \right),\\
     \adj A&=(\cof A)^T=\left(\begin{array}{rr}6&-4\\5&3\end{array}\right)^T
                      =\left(\begin{array}{rr}6&5\\-4&3\end{array}\right).
  \end{align*}
\item Ponovo izra\v cunajmo prvo kofaktore
{\scriptsize
\begin{align*}
   A_{11}&=(-1)^{1+1}\left|\begin{array}{rr}2&1\\1&3 \end{array}\right|=5,&&\hspace{-.2cm}
   A_{12}=(-1)^{1+2}\left|\begin{array}{rr}0&1\\1&3 \end{array}\right|=1,&&\hspace{-.2cm}
   A_{13}=(-1)^{1+3}\left|\begin{array}{rr}0&2\\1&1 \end{array}\right|=-2,&\hspace{-.2cm} \\
   A_{21}&=(-1)^{2+1}\left|\begin{array}{rr}2&-5\\1&3 \end{array}\right|=-11,&&\hspace{-.2cm}
   A_{22}=(-1)^{2+2}\left|\begin{array}{rr}1&-5\\1&3 \end{array}\right|=8,&&\hspace{-.2cm}
   A_{23}=(-1)^{2+3}\left|\begin{array}{rr}1&2\\1&1 \end{array}\right|=1,&\hspace{-.2cm} \\
   A_{31}&=(-1)^{3+1}\left|\begin{array}{rr}2&-5\\2&1 \end{array}\right|=12,&&\hspace{-.2cm}
   A_{32}=(-1)^{3+2}\left|\begin{array}{rr}1&-5\\0&1 \end{array}\right|=-1,&&\hspace{-.2cm}
   A_{33}=(-1)^{3+3}\left|\begin{array}{rr}1&2\\0&2 \end{array}\right|=2.&\hspace{-.2cm}
\end{align*}
}
Matrica kofaktora je
\[\cof A=\left(\begin{array}{rrr} 5&1&-2\\-11&8&1\\12&-1&2\end{array} \right),\]
i na kraju adjungirana matrica
\[\adj A=(\cof A)^T=\left(\begin{array}{rrr} 5&-11&12\\1&8&-1\\-2&1&2\end{array} \right).\]
\end{enumerate}
\end{example}

Slijedi definicija inverzne matrice.

\begin{definition}[Inverzna matrica]\index{matrica! inverzna}
  Ako je $A$ kvadratna matrica reda $n$ i ako postoji matrica $X$ takva da je $XA=AX=I,$ tada ka\v zemo da je $X$ inverzna matrica matrice $A.$ Inverznu matricu matrice $A,$ ozna\v cavamo sa $A^{-1}.$
\end{definition}

Neke kvadratne matrice imaju inverznu matricu, dok neke nemaju. Da bi napravili razliku izme\dj u ove dvije klase matrica uvodimo dva pojma sljede\'com definicijom.

\begin{definition}[Regularna i singularna matrica] \index{matrica! regularna}\index{matrica! singularna}
   Kvadratna matrica $A$ je regularna (nesingularna) matrica ako ima inverznu matricu. Ako kvadratna matrica $A$ nema inverznu matricu, ka\v zemo da je $A$ singularna (neregularna) matrica.
\end{definition}

Sada se postavlja pitanje kako napraviti razliku izmedju regularnih i singularnih matrica, tj. za neku datu kvadratnu matricu kako odrediti je li ona regularna ili singularna? Ako smo utvrdili da je ta matrica regularna, onda na osnovu prethodne definicije postoji njena inverzna matrica, kako sada izra\v cunati tu inverznu matricu? Zatim za datu regularnu matricu da li postoji jedna ili vi\v se njenih inverznih matrica i obrnuto, za inverznu matricu da li postoji jedna ili vi\v se regularnih matrica za koje je ona inverzna matrica? Odgovori na ova pitanja dati su u sljede\'cim teoremama.

\begin{theorem}\index{teorema! regularna matrica}
   Kvadratna matrica $A=(a_{ij})$ je regularna ako i samo ako je $\det A\neq 0.$ Ako je $\det A\neq 0,$ inverzna matrica matrice $A$ je
   \[A^{-1}=\frac{1}{\det A}\adj A.\]
\end{theorem}
\begin{proof}[Dokaz]
Pretpostavimo da je matrica $A$ regularna, tj. ima inverznu matricu $A^{-1}.$ Tada je $\det \left(A \cdot A^{-1} \right)=\det I=1.$ Kako je
$\det \left(A \cdot A^{-1} \right)=\det A \cdot \det A^{-1}=1,$ slijedi da je $\det A \neq 0.$  Osim toga vidimo da je
$\det A^{-1}=\frac{1}{\det A}=\left(\det A \right)^{-1}.$ \\
Pretpostavimo sada da je $\det A \neq 0.$ Treba dokazati da postoji matrica $A^{-1}.$\\
 \[    A\cdot \adj A= \left( \begin{array}{cccc}
                a_{11}&a_{12}&\ldots&a_{1n}\\
                a_{21}&a_{22}&\ldots&a_{2n}\\
                &\vdots&&\\
                a_{n1}&a_{n2}&\ldots&a_{nn}
                \end{array}
         \right) \cdot \left( \begin{array}{cccc}
                A_{11}&A_{21}&\ldots&A_{n1}\\
                A_{12}&A_{22}&\ldots&A_{n2}\\
                &\vdots&&\\
                A_{1n}&A_{2n}&\ldots&A_{nn}
                \end{array}
         \right)=
  \]
  \[ =\left(a_{i1}A_{j1}+a_{i2}A_{j2}+ \ldots +a_{in}A_{jn} \right)_{n\times n} . \]
  S obzirom da je
  \[ a_{i1}A_{j1}+a_{i2}A_{j2}+ \ldots +a_{in}A_{jn}=\left\{\begin{array}{cl}
                                                          \det A, & \text{ako je } i=j \\
                                                          0, & \text{ako je } i \neq j
                                                          \end{array}
                                                    \right.
  \]
  dobijamo da je
   \[    A\cdot \adj A= \left( \begin{array}{cccc}
                \det A & 0 &\ldots& 0\\
                0 & \det A &\ldots& 0\\
                &\vdots&&\\
                0 & 0 &\ldots & \det A
                \end{array}
         \right)= \det A \cdot I.
  \]
 Kako je po pretpostavci $\det A \neq 0,$ zaklju\v cujemo da je $A \cdot (\det A)^{-1}\cdot \adj A=I.$ To zna\v ci da je matrica
 $(\det A)^{-1}\cdot \adj A$ desni inverz matrice $A.$ Na sli\v can na\v cin bismo zaklju\v cili da je $\adj A \cdot A=\det A \cdot I,$
 odakle slijedi da je matrica $(\det A)^{-1}\cdot \adj A$ tako\dj er i lijevi inverz matrice $A.$ Dakle, ako je kvadratna matrica $A$ regularna,
 tada postoji njena inverzna matrica koja je jednaka
 \[ A^{-1}=(\det A)^{-1}\cdot \adj A= \frac{1}{\det A}\adj A . \]
  \end{proof}
\begin{theorem}\label{regularna1}\index{teorema! regularna i inverzna matrica}
 Svaka regularna matrica $A$ ima jednu i samo jednu inverznu matricu $A^{-1}.$
\end{theorem}

U sljede\'coj teoremi date su neke osobine inverznih matrica.

\begin{theorem}[Osobine inverznih matrica]\index{teorema! osobine inverznih matrica}
Ako su $A$ i $B$ regularne matrice istog reda, tada je
\begin{enumerate}[$(1)$]
  \item $(AB)^{-1}=B^{-1}A^{-1}$;
  \item $\left(A^{-1}\right)^{-1}=A;$
  \item $\left(A^T\right)^{-1}=\left(A^{-1}\right)^{T};$
  \item $\det A^{-1}=(\det A)^{-1}.$
\end{enumerate}
\end{theorem}

\begin{example}
Date su matrice
  \begin{enumerate}[$(a)$]
     \item $A=\left(\begin{array}{rr}2&-3\\-4&6\end{array}\right);$
     \item $B=\left(\begin{array}{rr}2&-3\\0&1\end{array}\right);$
     \item $C=\left(\begin{array}{rrr}-1&0&-2\\0&2&1\\1&-1&2\end{array}\right);$
     \item $D=\left(\begin{array}{rrr}-2&0&-2\\3&2&1\\1&-1&2\end{array}\right).$
  \end{enumerate}
Ispitati jesu li date matrice regularne. Ako jesu, izra\v cunati njihove inverzne matrice.   \\\\
\noindent Rje\v senje:\\\\
Postupak ra\v cunanja inverzne matrice matrice $X$, je sljede\' ci:
\begin{enumerate}[\quad I --]
   \item Izra\v cunamo $\det X.$ Ako je $\det X\neq 0$ nastavljamo dalje
   \item Ra\v cunamo matricu $\cof X,$ a zatim $\adj X=(\cof X)^T $
   \item Ra\v cunamo inverznu matricu po formuli $X^{-1}=\frac{1}{\det X}\adj X.$
\end{enumerate}
Na kraju mo\v zemo provjeriti rezultat, odnosno da li vrijedi $X^{-1} \cdot X=I.$
 \begin{enumerate}[$(a)$]
   \item Po\v sto je $\det A=\left|\begin{array}{rr}2&-3\\-4&6\end{array}\right|=12-12=0,$ matrica $A$ je singularna, tj. nema inverzne matrice.
   \item Kako je $\det B=\left|\begin{array}{rr}2&-3\\0&1\end{array}\right|=2,$ matrica je regularna pa mo\v zemo ra\v cunati\\ inverznu matricu. Prvo izra\v cunajmo kofaktore i matricu kofaktora
    \begin{align*}
        B_{11}&=(-1)^{1+1}\cdot1=1,&&B_{12}=(-1)^{1+2}\cdot0=0\\
        B_{21}&=(-1)^{2+1}\cdot(-3)=3,&&B_{22}=(-1)^{2+2}\cdot2=2
  \end{align*}
  formiramo matricu kofaktora $\cof B$ i transponujemo je
  \begin{align*}
     \cof B&=\left(\begin{array}{rr}1&0\\3&2\end{array}       \right),\\
     \adj B&=(\cof B)^T=\left(\begin{array}{rr}1&0\\3&2\end{array}\right)^T
                      =\left(\begin{array}{rr}1&3\\0&2\end{array}\right).
  \end{align*}
  Inverzna matrica je \[B^{-1}=\frac{1}{\det B}\adj B
            =\frac{1}{2}\left(\begin{array}{rr}1&3\\0&2\end{array}\right)=
            \left(\begin{array}{rr}\frac{1}{2}&\frac{3}{2}\\0&1\end{array}\right).\]
  Provjera:
 \[B^{-1} \cdot B=\frac{1}{2}\left(\begin{array}{rr}1&3\\0&2\end{array}\right) \cdot \left(\begin{array}{rr}2&-3\\0&1\end{array}\right)=
 \left(\begin{array}{rr}1&0\\0&1\end{array}\right) \]

\item Kako je $\det C=-1,$ matrica je regularna. Izra\v cunajmo prvo kofaktore
{\scriptsize
       \begin{align*}
   C_{11}&=(-1)^{1+1}\left|\begin{array}{rr}2&1\\-1&2 \end{array}\right|=5,\,&&\hspace{-.2cm}
   C_{12}=(-1)^{1+2}\left|\begin{array}{rr}0&1\\1&2 \end{array}\right|=1,\,&&\hspace{-.2cm}
   C_{13}=(-1)^{1+3}\left|\begin{array}{rr}0&2\\1&-1 \end{array}\right|=-2,&\hspace{-.2cm} \\
   C_{21}&=(-1)^{2+1}\left|\begin{array}{rr}0&-2\\-1&2 \end{array}\right|=2,\,&&\hspace{-.2cm}
   C_{22}=(-1)^{2+2}\left|\begin{array}{rr}-1&-2\\1&2 \end{array}\right|=0,\,&&\hspace{-.2cm}
   C_{23}=(-1)^{2+3}\left|\begin{array}{rr}-1&0\\1&-1 \end{array}\right|=-1,&\hspace{-.2cm} \\
   C_{31}&=(-1)^{3+1}\left|\begin{array}{rr}0&-2\\2&1 \end{array}\right|=4\,&&\hspace{-.2cm}
   C_{32}=(-1)^{3+2}\left|\begin{array}{rr}-1&-2\\0&1 \end{array}\right|=1,\,&&\hspace{-.2cm}
   C_{33}=(-1)^{3+3}\left|\begin{array}{rr}-1&0\\0&2 \end{array}\right|=-2.&\hspace{-.2cm}
\end{align*}
}
Matrica kofaktora je
\[\cof C=\left(\begin{array}{rrr} 5&1&-2\\2&0&-1\\4&1&-2\end{array} \right),\]
adjungovana matrica
\[\adj C=(\cof C)^T=\left(\begin{array}{rrr} 5&2&4\\1&0&1\\-2&-1&-2\end{array} \right)\]
i na kraju inverzna matrica je
\[C^{-1}=\frac{1}{\det C}\adj C
            =\frac{1}{-1}\left(\begin{array}{rrr} 5&2&4\\1&0&1\\-2&-1&-2\end{array} \right)=
            \left(\begin{array}{rrr} -5&-2&-4\\-1&0&-1\\2&1&2\end{array} \right).\]
\item Kako je $\det D=0,$ matrica $D$ je singularna.
 \end{enumerate}
\end{example}\index{matrica! matri\v cne jedna\v cine}
\begin{example}
  Rije\v siti matri\v cnu jedna\v cinu $AX+B=3X+I,$ ako su
  \[A=\left(\begin{array}{rr}2&-3\\-4&6\end{array}\right),\:
    B=\left(\begin{array}{rr}-1&0\\2&3\end{array}\right).\]
\noindent Rje\v senje:\\\\
 Datu matri\v cnu jedna\v cinu mo\v zemo napisati u obliku \[AX+B=3X+I\Leftrightarrow AX-3X=I-B\Leftrightarrow(\underbrace{A-3I}_{=C})X=\underbrace{I-B}_{=D}
               \Leftrightarrow CX=D,\]
({\bfseries \color{red}$X$ izvla\v cimo sa desne strane}) gdje je
 \begin{align*}
   C&=A-3I= \left(\begin{array}{rr}2&-3\\-4&6\end{array}\right)-3\left(\begin{array}{rr}1&0\\0&1\end{array}\right)
           =\left(\begin{array}{rr}-1&-3\\-4&3\end{array}\right)\\
   D&=I-B=\left(\begin{array}{rr}1&0\\0&1\end{array}\right)-\left(\begin{array}{rr}-1&0\\2&3\end{array}\right)
          =\left(\begin{array}{rr}2&0\\-2&-2\end{array}\right).
 \end{align*}
Vrijedi
 \[\det C=\left|\begin{array}{cc}-1&-3\\-4&3\end{array}\right|=-3-12=-15\neq 0.\]
Budu\'ci da je $\det C=15,$ matri\v cnu jedna\v cinu $CX=D,$ pomno\v zimo sa lijeve strane, jer mno\v zenje matrica nije komutativno, sa matricom $C^{-1}$ i dobijamo \[X=C^{-1}D,\] tj. da bi izra\v cunali nepoznatu matricu $X$ treba jo\v s izra\v cunati inverznu matricu matrice $C$ i pomno\v ziti matricu $D$ sa lijeve strane sa matricom $C^{-1}.$
Izra\v cunajmo sada kofaktore, vrijedi
\begin{align*}
 C_{11}&=(-1)^{1+1}\cdot2=3&&C_{12}=(-1)^{1+2}\cdot(-4)=4\\
 C_{21}&=(-1)^{2+1}\cdot(-3)=3&&C_{22}=(-1)^{2+2}\cdot(-1)=-1,
\end{align*}
pa je \begin{align*}
   \cof C&=\left(\begin{array}{cc}3&4\\3&-1\end{array}\right),\quad\adj C=(\cof C)^T=\left(\begin{array}{cc}3&3\\4&-1\end{array}\right),\\
   C^{-1}&=\frac{1}{\det C}\adj C=\frac{1}{-15}\left(\begin{array}{cc}3&3\\4&-1\end{array}\right)\\
   X&=C^{-1}D=-\frac{1}{15}\left(\begin{array}{cc}3&3\\4&-1\end{array}\right)\left(\begin{array}{rr}2&0\\-2&-2\end{array}\right)=
        -\frac{1}{15}\left(\begin{array}{rr}0&-6\\10&2\end{array}\right)
        =\left(\begin{array}{rr}0&\frac{2}{5}\\-\frac{2}{3}&-\frac{2}{15}\end{array}\right).
   \end{align*}
\end{example}

\begin{example}
  Rije\v siti matri\v cnu jedna\v cinu $XA-A=2X+I,$ ako je
  \[A=\left(\begin{array}{rrr}0&1&2\\2&3&4\\1&0&1\end{array}\right).\]
\noindent Rje\v senje:\\\\
 Matri\v cnu jedna\v cinu napi\v simo u sljede\' cem obliku
 \[XA-A=2X+I\Leftrightarrow XA-2X=I+A\Leftrightarrow X(\underbrace{A-2I}_{=C})=\underbrace{I+A}_{=D},\]
 ({\bfseries \color{red} $X$ izvla\v cimo sa lijeve strane}), gdje je

\begin{align*}
 C&=A-2I=\left(\begin{array}{rrr}0&1&2\\2&3&4\\1&0&1\end{array}\right)-2\left(\begin{array}{rrr}1&0&0\\0&1&0\\0&0&1\end{array}\right)
    =\left(\begin{array}{rrr}-2&1&2\\2&1&4\\1&0&-1\end{array}\right),\\
 D&=\left(\begin{array}{rrr}1&0&0\\0&1&0\\0&0&1\end{array}\right)+\left(\begin{array}{rrr}0&1&2\\2&3&4\\1&0&1\end{array}\right)
     =\left(\begin{array}{rrr}1&1&2\\2&4&4\\1&0&2\end{array}\right).
 \end{align*}
Vrijednost $\det C$ je
\[\det C=\left|\begin{array}{rrr}-2&1&2\\2&1&4\\1&0&-1\end{array}\right|=6\neq 0,\]
sada matri\v cnu jedna\v cinu $XC=D $ pomno\v zimo sa desne strane inverznom matricom $C^{-1}$ matrice $C$ i dobijamo \[X=DC^{-1}.\]  Matrica $C$ je regularna pa izra\v cunajmo matricu kofaktora i adjungovanu matricu
{\scriptsize
       \begin{align*}
   C_{11}&=(-1)^{1+1}\left|\begin{array}{rr}1&4\\0&-1 \end{array}\right|=-1,\,&&\hspace{-.2cm}
   C_{12}=(-1)^{1+2}\left|\begin{array}{rr}2&4\\1&-1 \end{array}\right|=6,\,&&\hspace{-.2cm}
   C_{13}=(-1)^{1+3}\left|\begin{array}{rr}2&1\\1&0 \end{array}\right|=-1,& \hspace{-.2cm}  \\
   C_{21}&=(-1)^{2+1}\left|\begin{array}{rr}1&2\\0&-1 \end{array}\right|=1,\,&&\hspace{-.2cm}
   C_{22}=(-1)^{2+2}\left|\begin{array}{rr}-2&2\\1&-1 \end{array}\right|=0,\,&&\hspace{-.2cm}
   C_{23}=(-1)^{2+3}\left|\begin{array}{rr}-2&1\\1&0 \end{array}\right|=1,& \hspace{-.2cm}  \\
   C_{31}&=(-1)^{3+1}\left|\begin{array}{rr}1&2\\1&4 \end{array}\right|=2\,&&\hspace{-.2cm}
   C_{32}=(-1)^{3+2}\left|\begin{array}{rr}-2&2\\2&4 \end{array}\right|=12,\,&&\hspace{-.2cm}
   C_{33}=(-1)^{3+3}\left|\begin{array}{rr}-2&1\\2&1 \end{array}\right|=-4.&\hspace{-.2cm}
\end{align*}
}
 Matrica kofaktora je
\[\cof C=\left(\begin{array}{rrr} -1&6&-1\\1&0&1\\2&12&-4\end{array} \right),\]
adjungovana matrica
\[\adj C=(\cof C)^T=\left(\begin{array}{rrr} -1&1&2\\6&0&12\\-1&1&-4\end{array} \right),\]
te je inverzna matrica
\[C^{-1}=\frac{1}{\det C}\adj C
            =\frac{1}{6}\left(\begin{array}{rrr} -1&1&2\\6&0&12\\-1&1&-4\end{array} \right).\]
Sada je
\[X=DC^{-1}=\left(\begin{array}{rrr}1&1&2\\2&4&4\\1&0&2\end{array}\right)\frac{1}{6}\left(\begin{array}{rrr} -1&1&2\\6&0&12\\-1&1&-4\end{array} \right)
       =\frac{1}{6}\left(\begin{array}{rrr} 3&3&6\\18&6&36\\-3&3&-6\end{array} \right). \]
\end{example}

\subsection{Z\lowercase{adaci}}\index{Zadaci za vje\v zbu!inverzna matrica}
\begin{enumerate}
  \item Odrediti inverznu matricu matrice \\
\begin{inparaenum}
\item $\left(\begin{array}{rr}
1&2\\3&4\\	
\end{array}\right)$;\quad
\item $\left(\begin{array}{rrr}
1&2&-3\\0&1&2\\0&0&1\\
\end{array}\right);\;$ \quad
\item $\left(\begin{array}{rrr}
3&-4&5\\2&-3&1\\3&-5&-1\\
\end{array}\right).$
\end{inparaenum}

\item Rije\v siti matri\v cne jedna\v cine\\
\begin{inparaenum}
\item $\left(\begin{array}{rr}
2&1\\1&-1\\	
\end{array}\right)\cdot X=\left(\begin{array}{rr}1&3\\2&0	
\end{array}\right)$; \quad
\item $X\cdot \left(\begin{array}{rr}
3&-2\\5&-4\\	
\end{array}\right)=\left(\begin{array}{rr}-1&2\\-5&6	
\end{array}\right).\vspace{.2cm}$
\end{inparaenum}

\item Rije\v siti matri\v cnu jedna\v cinu $(I-2A)X=B,$ ako je
\[ A=\left(\begin{array}{cc}2&-3\\1&0\end{array}\right)  \text{ i } B=\left(\begin{array}{cc}-3&0\\-6&1\end{array} \right).\]





\item Rije\v siti matri\v cnu jedna\v cinu
  \[ \left(\begin{array}{ccc}1&2&-3\\3&2&-4\\2&-1&0\end{array}\right)\cdot X=\left(\begin{array}{ccc}1&-3&0\\10&2&7\\10&7&8\end{array}\right).\]
\item Rije\v siti matri\v cnu jedna\v cinu
  \[X\cdot\left(\begin{array}{ccc}5&3&1\\1&-3&-2\\-5&2&1\end{array}\right)=\left(\begin{array}{ccc}-8&3&0\\-5&9&0\\-2&15&0\end{array}\right).\]


\end{enumerate}

\section{R\lowercase{ang matrice}}\index{matrica! rang}

Posmatrajmo sljede\' ce jednakosti
\begin{align*}
\alpha_1\left(\begin{array}{c} 1\\0\\0\end{array}\right)+\alpha_2\left(\begin{array}{c} 0\\1\\0\end{array}\right)+\alpha_3\left(\begin{array}{c} 0\\0\\1\end{array}\right)
= \left(\begin{array}{c} \alpha_1\\0\\0\end{array}\right)+\left(\begin{array}{c} 0\\\alpha_2\\0\end{array}\right)+\left(\begin{array}{c} 0\\0\\\alpha_3\end{array}\right)
= \left(\begin{array}{c} \alpha_1\\\alpha_2\\\alpha_3\end{array}\right),
\end{align*}
i
\begin{align*}
\alpha_1\left(\begin{array}{c} 1\\0\\0\end{array}\right)+\alpha_2\left(\begin{array}{c} 2\\0\\0\end{array}\right)+\alpha_3\left(\begin{array}{c} 5\\0\\0\end{array}\right)
= \left(\begin{array}{c} \alpha_1\\0\\0\end{array}\right)+\left(\begin{array}{c} 2\alpha_2\\0\\0\end{array}\right)+\left(\begin{array}{c} 5\alpha_3\\0\\0\end{array}\right)
= \left(\begin{array}{c} \alpha_1+2\alpha_2+5\alpha_3\\ 0\\ 0\end{array}\right).
\end{align*}
Postavlja se pitanje kolika je vrijednost koeficijenata $\alpha_1,\,\alpha_2,\,\alpha_3$ tako da rezultantna matrica kolona (na desnoj strani jednakosti) bude jednaka nuli? U prvom slu\v caju zbir \' ce biti jedino nula ako je $\alpha_1=\alpha_2=\alpha_3=0,$ dok u drugom slu\v caju  vrijednost zbira (rezultantne matrice kolone) mo\v ze biti jednaka nuli i kada je npr. $\alpha_1=1,\,\alpha_2=2,\,\alpha_3=-1.$ Drugim rije\v cima, u drugom primjeru linearna kombinacija matrica kolona i koeficijenata se mo\v ze anulirati za neke koeficijente $\alpha_1,\,\alpha_2,\,\alpha_3,$ koji nisu svi nule. A u prvom primjeru to nije slu\v caj, tj. mogu\' ce je jedino za $\alpha_1=\alpha_2=\alpha_3=0.$ Ovo nam daje opravdanje za uvo\dj enje sljede\' ceg pojma.

\begin{definition}[Linearna nezavisnost]
  Za matrice vrste (ili matrice kolone) $A_1,A_2,\ldots,A_n$ ka\v zemo da su linearno nezavisne ako za realne brojeve $\alpha_1,\alpha_2,\ldots,\alpha_n$ iz jednakosti
  \begin{equation}
   \alpha_1A_1+\alpha_2A_2+\ldots+\alpha_nA_n=\textbf{0}
   \label{nezavisnost1}
  \end{equation}
  slijedi da je jedino rje\v senje $\alpha_1=\alpha_2=\ldots=\alpha_n=0.$ U suprotnom, tj. ako matrice vrste (ili matrice kolone) nisu linearno nezavisne, onda za njih ka\v zemo da su linearno zavisne.
\end{definition}
Simbol $\textbf{0}$ ozna\v cava nula matricu vrstu ili nula matricu kolonu. U nekoj matrici zanima nas koliki je broj maksimalno nezavisnih vrsta i kolona. Sljede\'ca teorema ka\v ze da je svejedno koji \'cemo broj odrediti:  maksimalan broj nezavisnih vrsta ili kolona te matrice.

\begin{theorem}\index{teorema! linearno nezavisne vrste i kolone}
 Maksimalan broj linearno nezavisnih vrsta posmatrane matrice jednak je maksima-\\lnom broju linearno nezavisnih kolona te matrice.
\label{nezavisnost2}
\end{theorem}

Vrijedi sljede\' ca definicija.
\begin{definition}[Rang matrice]
Pod pojmom ranga matrice $A=(a_{ij})$ formata $m\times n$ podrazumijevamo maksimalan broj linearno nezavisnih vrsta (ili kolona) te matrice.
\end{definition}
Rang matrice ozna\v cava\' cemo sa $\rang A.$ Mo\v ze se jo\v s rang  matrice ozna\v cavati i sa $\rang(A)$ ili $r(A).$ Za matricu $A$ formata $m\times n$ jasno je da vrijedi
\[\rang A\leqslant\min\{m,n\}.\]

Prakti\v cno ra\v cunanje ranga matrice prili\v cno je lako izvesti kori\v stenjem elementarnih transformacija. Primjenom elementarnih transformacija rang matrice se ne mijenja.

\begin{definition}[Elementarne trasformacije]\index{matrica! elementarne transformacije}
 Pod elementarnim transformacijama jedne matrice smatraju se operacije
 \begin{enumerate}
  \item Razmjena mjesta dvije vrste (ili dvije kolone);
  \item Dodavanje elementima jedne vrste (ili kolone) elemenata neke druge vrste (ili kolone), po\v sto su prethodno posljednje pomno\v zene nekim brojem razli\v citim od nule;
  \item Mno\v zenje elemenata jedne vrste (ili kolone) nekim brojem razli\v citim od nule.
 \end{enumerate}
\end{definition}
Matrice dobijene elementarnim transformacijama su ekvivalentne (koristimo oznaku "$ \sim $").  Datu matricu, primjenjuju\' ci elementarne trasformacije svodimo na trougaonu matricu ili na trapeznu ili stepenu formu matrice (za matrice u trapeznoj ili stepenoj  formi u literaturi se  koristi  i naziv kvazitrougaone matrice). U datim primjerima bi\' ce jasnije o kakvim se formama radi.
\begin{remark}
  Elementarne transformacije za odre\dj ivanje ranga matrice na\v celno se koriste samo za manipulaciju sa vrstama, zbog primjene elementarnih transformacija pri rje\v savanju sistema linearnih algebarskih jedna\v cina. O ovome \' ce biti vi\v se rije\v ci u narednom poglavlju.
\end{remark}

\begin{example}\label{primjerRang1}
  Odrediti rang matrice $A=\left(\begin{matrix}4&1&1\\1&2&1\\1&1&2\end{matrix}\right).$\\
\noindent Rje\v senje:
 \begin{align*}
    \begin{blockarray}{rrrl}
       \begin{block}{(rrr)l}
         4&1&1& $Ik$\leftrightarrow $IIk$\\
         1&2&1&\\
         1&1&2&\\
       \end{block}
   \end{blockarray}&\sim
   \begin{blockarray}{rrrl}
       \begin{block}{(rrr)l}
         1&4&1& \\
         2&1&1&\text{ IIv-2Iv}\\
         1&1&2&\text{ IIIv-Iv}\\
       \end{block}
   \end{blockarray}\sim
     \begin{blockarray}{rrrl}
       \begin{block}{(rrr)l}
         1&4&1& \\
         0&-7&-1&\\
         0&-3&1&\text{ 7IIIv-3IIv}\\
       \end{block}
   \end{blockarray}\\
     &\sim
      \begin{blockarray}{rrrl}
       \begin{block}{(rrr)l}
         {\color{magenta}\circled{1}}&4&1& \\
         0&{\color{magenta}\circled{-7}}&-1&\\
         0&0&{\color{magenta}\circled{10}}&,\\
       \end{block}
   \end{blockarray}
    \end{align*}
    pa je $\rang A=3.$
\end{example}
\begin{example}\label{primjerRang2}
  Odrediti rang matrice $A=\begin{blockarray}{rrrr}
                              \begin{block}{(rrrr)}
                                      2&3&-1&4\\
                                      5&-3&8&19\\
                                      1&-2&3&5\\
                              \end{block}
                            \end{blockarray}.$\ \\
\noindent Rje\v senje:
  \begin{align*}
  A&=\begin{blockarray}{rrrrl}
        \begin{block}{(rrrr)l}
                 2&3&-1&4&\text{Iv}\leftrightarrow \text{IIIv}\\
                 5&-3&8&19\\
                 1&-2&3&5\\
        \end{block}
     \end{blockarray}\sim
     \begin{blockarray}{rrrrl}
        \begin{block}{(rrrr)l}
                 1&-2&3&5&\\
                 5&-3&8&19&\text{ IIv-5Iv}\\
                 2&3&-1&4&\text{IIIv-2Iv}\\
        \end{block}
     \end{blockarray}\\
     &\sim
     \begin{blockarray}{rrrrl}
        \begin{block}{(rrrr)l}
                 1&-2&3&5&\\
                 0&7&-7&-6&\\
                 0&7&-7&-6&\text{IIIv-IIv}\\
        \end{block}
     \end{blockarray}
     \sim
          \begin{blockarray}{rrrrl}
        \begin{block}{(rrrr)l}
                 {\color{magenta}\circled{1}}&-2&3&5&\\
                 0&{\color{magenta}\circled{7}}&-7&-6&\\
                 0&0&0&0&\\
        \end{block}
     \end{blockarray}.
  \end{align*}
  U ovom slu\v caju je $\rang A=2.$
\end{example}

\begin{example}\label{primjerRang3}
Odrediti rang matrice $A=\begin{blockarray}{rrrrr}
                           \begin{block}{(rrrrr)}
                             3&6&6&9&1\\
                             2&4&1&2&0\\
                             -1&-2&4&5&1\\
                           \end{block}
                         \end{blockarray}.$\ \\
\noindent Rje\v senje:
\begin{align*}
   A&=\begin{blockarray}{rrrrrl}
         \begin{block}{(rrrrr)l}
            3&6&6&9&1&\text{Iv}\leftrightarrow \text{IIIv}\\
            2&4&1&2&0\\
            -1&-2&4&5&1\\
        \end{block}
      \end{blockarray} \sim
      \begin{blockarray}{rrrrrl}
         \begin{block}{(rrrrr)l}
           -1&-2&4&5&1\\
            2&4&1&2&0&\text{IIv+2Iv}\\
            3&6&6&9&1&\text{IIIv+3Iv}\\
         \end{block}
      \end{blockarray} \\
      &\sim
      \begin{blockarray}{rrrrrl}
         \begin{block}{(rrrrr)l}
           -1&-2&4&5&1\\
            0&0&9&12&2\\
            0&0&18&24&4&\text{IIIv-2IIv}\\
         \end{block}
      \end{blockarray} \sim
       \begin{blockarray}{rrrrrl}
         \begin{block}{(rrrrr)l}
           {\color{magenta}\circled{-1}}&-2&4&5&1\\
            0&0&{\color{magenta}\circled{9}}&12&2\\
            0&0&0&0&0&\text{IIIv-2IIv},\\
         \end{block}
      \end{blockarray}
\end{align*}
pa je $\rang A=2.$
\end{example}
\begin{remark}
	Rezultantna matrica u Primjeru \ref{primjerRang1} je u trougaonoj formi, u Primjeru \ref{primjerRang2} to je trapezna forma i u Primjeru \ref{primjerRang3} rezultantna matrica je u stepenoj formi. Iz prezentovanih primjera vidimo da se rang matrice odre\dj uje tako \v sto prebrojimo broj vrsta koje imaju nenulte elemente nakon svo\dj enja matrice na jednu od pomenutih formi.
\end{remark}

\subsection{Z\lowercase{adaci}} \index{Zadaci za vje\v zbu!rang matrice}
\begin{enumerate}

\item Odrediti rang matrica \\
      \begin{inparaenum}
      	  \item $\begin{pmatrix}
      	            1 & 0 & 0 & 0\\ 2 & 0 & 0 & 0\\ -1 & 0 & 0 & 0 \\
             	  \end{pmatrix}; $ \quad
          \item $\begin{pmatrix}
                     1 & 0 & 0 & 0 \\ 2 & 1 & 0 & 0\\ -1 & 0 & 0 & 0 \\
               \end{pmatrix}; $ \quad
          \item $\begin{pmatrix}
                    1 & 0 & 0 & 0\\ 2 & 1 & 2 & 0\\ -1 & 4 & 1 & 0 \\
                \end{pmatrix}. $ \quad   	
      \end{inparaenum}

\item Odrediti rang matrica \\
\begin{inparaenum}
	\item $\begin{pmatrix}
		         1 & 0 & 0  \\ -3 & 0 & 0 \\ 2 & 0 & 0\\ 4 & 0 & 0 \\
	     \end{pmatrix}; $ \quad
	\item $\begin{pmatrix}
	            1 & 0 & 0  \\ -3 & 1 & 0 \\ 2 & 0 & 0\\ 4 & 0 & 0 \\
	      \end{pmatrix}; $ \quad
	\item $\begin{pmatrix}
                1 & 0 & 0  \\ -3 & 1 & 0 \\ 2 & 2 & 0\\ 4 & 0 & 0 \\
           \end{pmatrix}; $ \quad
    	\item $\begin{pmatrix}
                     1 & 0 & 0  \\ -3 & 1 & 0 \\ 2 & 2 & 0\\ 4 & -1 & 2 \\
               \end{pmatrix}. $ \quad
           	
\end{inparaenum}
\item Odrediti rang matrica\\
  \begin{inparaenum}
   \item $\begin{pmatrix}
             1&1&1&1\\2&3&-1&1\\3&4&0&2\\
         \end{pmatrix};\quad$
   \item $\begin{pmatrix}
             2&0&2&0&2\\0&1&0&1&0\\2&1&0&2&1\\0&1&0&1&0
         \end{pmatrix}; \quad $

   \item $\begin{pmatrix}
              1&1&-2\\2&1&-3\\2&-1&-1\\6&-1&-5\\7&-3&-4
          \end{pmatrix};\quad$
   \item $\begin{pmatrix}
	          1&6&7&1&4\\3&5&11&1&6\\12&5&3&1&4\\15&25&10&5&30\\	
         \end{pmatrix}.$
  \end{inparaenum}

\item Odrediti rang matrica  u zavisnost od parametra $a$

     \begin{inparaenum}
     	\item $\begin{pmatrix}
                   1 & a & -1 & 2\\ 2 & -1 & a & 5 \\ 1 & 10 & -6 & 1\\
              \end{pmatrix}; $ \quad
        \item $\begin{pmatrix}
                 3 & 1 & 1 & 4 \\ a & 4 & 10 & 1 \\ 1 & 7 & 17 & 3 \\ 2 & 2 & 4 & 1\\
              \end{pmatrix}. $ \quad
      \end{inparaenum}

\item Odrediti rang matrice  u zavisnost od parametra $a$  i $b$

    \begin{inparaenum}
         \item $\begin{pmatrix}
                      1 & 3 & 1 & -2 \\ 2 & 6 & -3 & -4 \\ a & b & 6 & -2 \\
                \end{pmatrix}. $ \quad
    \end{inparaenum}
\end{enumerate}

   \chapter[SLAJ]{Sistemi linearnih algebarskih jedna\v cina}



Mnogi problemi u matematici i njenim primjenama svode se na sisteme linearnih algebarskih jedna\v cina\footnote{U daljem tekstu umjesto punog naziva ovog poglavlja "Sistemi linearnih algebarskih jedna\v cina", nekada \'cemo koristiti "Sistemi linearnih jedna\v cina", "Sistemi jedna\v cina ", "Sistemi" a nekada samo skra\'cenicu -- SLAJ.}, odnosno na nala\v zenje rje\v senja ovih sistema. Na primjer, u najve\' cem broju slu\v cajeva veoma je te\v sko ili nemogu\'ce ta\v cno izra\v cunati ili rije\v siti neke druge vrste jedna\v cina (nelinearne, diferencijalne, integralne i dr.), pa se pribjegava ra\v cunanju pribli\v znih rje\v senja ovih jedna\v cina. Sastavni dio tih metoda za ra\v cunanje ovih pribli\v znih rje\v senja su upravo SLAJ i njihovo rje\v savanje.  SLAJ koje sre\'cemo u konkretnim pri-\\mjenama uglavnom imaju veliki broj jedna\v cina i nepoznatih, pa je stoga potrebno prvo utvrditi da li dati SLAJ ima rje\v senje (jedno ili vi\v se njih), a zatim je potrebno izra\v cunati  to rje\v senje ili rje\v senja, naravno ako postoje.

U ovom poglavlju bavimo se upravo navedenim problemima. Prvo, kako utvrditi da li SLAJ ima jedno ili vi\v se rje\v senja ili uop\v ste nema rje\v senja, a zatim ako ima, kako to rje\v senje ili rje\v senja izra\v cunati.

\section[O\lowercase{snovni pojmovi.} T\lowercase{eoreme}]{Osnovni pojmovi. Teoreme o egzistenciji\\ i jedinstvenosti rje\v senja}

Analizira\'cemo SLAJ koji imaju $n$ nepoznatih i $m$ jedna\v cina. Brojevi $n$ i $m$ mogu biti jednaki, a i razli\v citi. U  nastavku \'cemo vidjeti kako se radi sa SLAJ u oba slu\v caja, kada je $m=n,$ odnosno kada je $m\neq n.$  Slijedi definicija.

\begin{definition}[SLAJ]\label{sistem}\index{sistemi}
 Skup jedna\v cina
 \begin{align}
    a_{11}x_1+a_{12}x_2+\ldots+a_{1n}x_n&=b_1\nonumber\\
    a_{21}x_1+a_{22}x_2+\ldots+a_{2n}x_n&=b_2\label{sistem1}\\
    \vdots\hspace{2cm}&\nonumber\\
    a_{m1}x_1+a_{m2}x_2+\ldots+a_{mn}x_n&=b_m,\nonumber
 \end{align}
 gdje su $a_{ij}$ i $b_i$ dati brojevi, dok su $x_j$ nepoznate, $i=1,2,\ldots,m;\:j=1,2,\ldots,n;$ zove se sistem od $m$ algebarskih linearnih jedna\v cina sa $n$ nepoznatih. Brojevi $a_{ij}$ nazivaju se koeficijenti a brojevi $b_i$ slobodni \v clanovi.
 Ako je $b_1=b_2=\ldots=b_m=0,$ onda se ka\v ze da je sistem \eqref{sistem1}  homogen, a ako je bar jedan od slobodnih \v clanova $b_i\neq 0,$ onda ka\v zemo da je sistem \eqref{sistem1} nehomogen.
\end{definition}
Rje\v siti sistem zna\v ci izra\v cunati rje\v senje (ili rje\v senja) tog sistema, ako ono postoji (ili ona postoje). Slijedi definicija rje\v senja sistema linearnih algebarskih jedna\v cina.

\begin{definition}[Rje\v senje sistema]
  Brojevi $\alpha_1,\alpha_2,\ldots,\alpha_n$ nazivaju se rje\v senje sistema \eqref{sistem1} ako se zamjenom $x_j=\alpha_j,\:j=1,2,\ldots,n$ u sistem \eqref{sistem1} dobiju ta\v cne brojne jednakosti.
\end{definition}

Za dva sistema linearnih jedna\v cina ka\v ze se da su ekvivalentni  ako je svako rje\v senje jednog sistema ujedno rje\v senje i drugog sistema i obrnuto.

Jedan sistem mo\v ze imati jedno ili vi\v se rje\v senja, a mo\v ze i da nema rje\v senja. Ako sistem nema rje\v senja, ka\v ze se da je nesaglasan (protivrje\v can, nemogu\'c), a ako ima rje\v senje, ka\v ze se da je saglasan. Ako sistem ima ta\v cno jedno rje\v senje, ka\v ze se da ima jedinstveno rje\v senje, a ako ima vi\v se rje\v senja,  onda je taj sistem neodre\dj en.
\index{sistemi!nepoznate}\index{sistemi!koeficijenti}\index{sistemi!homogen}\index{sistemi!nehomogen}\index{sistemi!rje\v senje}\index{sistemi!ekvivalentni}

Posmatrajmo sljede\' ce sisteme linearnih jedna\v cina
\begin{align}
x+y=2&&x+y=2&&x+y=2\nonumber\\
x-y=0&&2x+2y=4&&x+y=3.\label{sistem2}
\end{align}

\begin{figure}[!h]\centering
	\begin{subfigure}[b]{.3\textwidth}\centering
		\includegraphics[scale=.65]{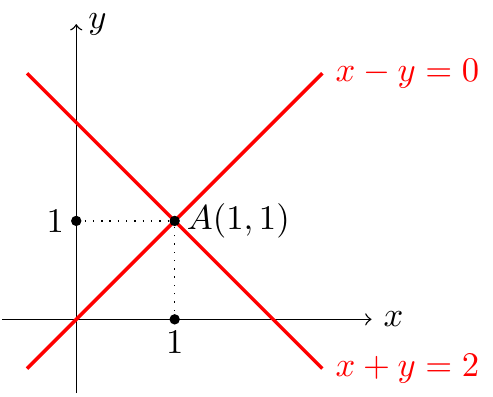}
		\caption{"Lijevi" sistem, prave se sijeku, jedna ta\v cka je zajedni\v cka. Koordinate zajedni\v cke
			ta\v cke predstavljaju rje\v senje sistema i rje\v senje je jedinstveno.}
		\label{slika32}
	\end{subfigure}\hspace{.5cm}
	\begin{subfigure}[b]{.3\textwidth}\centering
		\includegraphics[scale=.65]{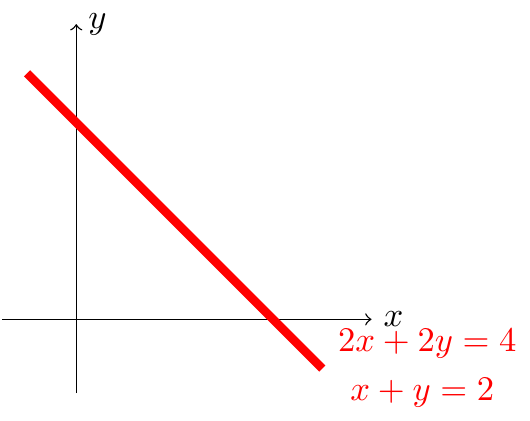}
		\caption{"Srednji" sistem, prave se poklapaju, sve su ta\v cke zajedni\v cke. Ovaj sistem ima beskona\v cno mnogo rje\v senja.}
		\label{slika33}
	\end{subfigure}
	\begin{subfigure}[b]{.3\textwidth}\centering\hspace{.5cm}
		\includegraphics[scale=.65]{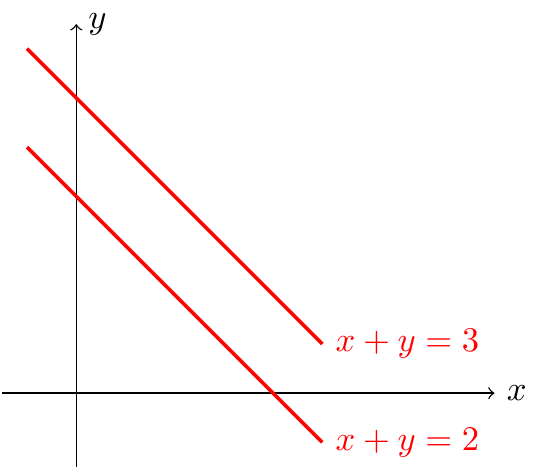}
		\caption{"Desni" sistem, nema zajedni\v ckih ta\v caka, pa prema tome ovaj sistem nema rje\v senja.}
		\label{slika34}
	\end{subfigure}
	\caption{Grafi\v cki prikaz rje\v senja sistema \eqref{sistem2} }
	\label{slika32a}
\end{figure}

Sistemi su jednostavni i ve\' c na prvi pogled vidimo da je rje\v senje "lijevog" sistema $x=1,\,y=1,$ u "srednjem" sistemu drugu jedna\v cinu smo dobili tako \v sto je prva pomno\v zena sa 2, pa ovaj sistem ima beskona\v cno mnogo rje\v senja, jer postoji beskona\v cno mnogo brojeva \v ciji je zbir 2. I na kraju, tre\' ci "desni"  sistem nema rje\v senja, jer ne postoje dva broja \v ciji je zbir istovremeno jednak i 2 i 3. Svaka jedna\v cina iz \eqref{sistem2} predstavlja jedna\v cinu prave u ravni. Na ovaj su na\v cin sistemi \eqref{sistem2} predstavljeni na Slici \ref{slika32a}.

\noindent Rje\v senje sistema su koordinate zajedni\v cke ta\v cke dviju pravih, koje smo nacrtali na osnovu jedna\v cina  koje \v cine taj sistem.  Sa Slike \ref{slika32} vidimo  da prave imaju jednu zajedni\v cku ta\v cku--"lijevi" sistem ima jedno rje\v senje, prave na Slici \ref{slika33} imaju beskona\v cno mnogo zajedni\v ckih ta\v caka--"srednji" sistem ima beskona\v cno mnogo rje\v senja i na Slici \ref{slika34} prave nemaju zajedni\v ckih ta\v caka--"desni" sistem nema rje\v senja.

\index{sistemi!saglasan}\index{sistemi!nesaglasan}\index{sistemi!jedinstveno rje\v senje}\index{sistemi!beskona\v cno mnogo rje\v senja}\index{sistemi!nema rje\v senja}

Postavlja se sada pitanje kako ustanoviti egzistenciju i broj rje\v senja za ne tako jedno-\\stavne sisteme kao u slu\v caju sistema \eqref{sistem2}. Odgovor na ova pitanja (egzistencija--postojanje i broj rje\v senja) dati su u narednim teoremama. Sistem \eqref{sistem1} mo\v zemo zapisati u matri\v cnoj formi
\begin{equation}
\begin{blockarray}{cccc}
       \begin{block}{(cccc)}
          a_{11}&a_{12}&\ldots&a_{1n}\\
          a_{21}&a_{22}&\ldots&a_{2n}\\
          &\vdots&&\\
          a_{m1}&a_{m2}&\ldots&a_{mn} \\
       \end{block}
\end{blockarray}\,
\begin{blockarray}{c}
       \begin{block}{(c)}
        x_1\\x_2\\\vdots\\x_n  \\
       \end{block}
\end{blockarray}
=
\begin{blockarray}{c}
   \begin{block}{(c)}
     b_1\\b_2\\ \vdots\\b_m  \\
   \end{block}
\end{blockarray}
\label{sistem3}
\end{equation}
ili
\[A\,X=B,\]
gdje su
\[A=\begin{blockarray}{cccc}
       \begin{block}{(cccc)}
          a_{11}&a_{12}&\ldots&a_{1n}\\
          a_{21}&a_{22}&\ldots&a_{2n}\\
          &\vdots&&\\
          a_{m1}&a_{m2}&\ldots&a_{mn} \\
       \end{block}
\end{blockarray},\quad
X=\begin{blockarray}{c}
       \begin{block}{(c)}
        x_1\\x_2\\\vdots\\x_n  \\
       \end{block}
\end{blockarray},\quad
B=\begin{blockarray}{c}
   \begin{block}{(c)}
     b_1\\b_2\\ \vdots\\b_m  \\
   \end{block}
\end{blockarray}.
\]
Matrica $A$ je matrica sistema, matrica kolona $X$ je matrica kolona nepoznatih i matrica kolona $B$ je matrica slobodnih \v clanova. Mo\v zemo od matrice $A$ formirati jo\v s jednu matricu, u oznaci $A| _B,$ tako \v sto \' cemo matrici $A$ dodati sa desne strane jo\v s jednu kolonu,  a to \' ce biti elementi matrice--kolone $B.$ Ovo je pro\v sirena matrica sistema i vrijedi
\[A| _B=\begin{blockarray}{ccccc}
       \begin{block}{(cccc|c)}
          a_{11}&a_{12}&\ldots&a_{1n}&b_1\\
          a_{21}&a_{22}&\ldots&a_{2n}&b_2\\
          &\vdots&&\\
          a_{m1}&a_{m2}&\ldots&a_{mn}&b_m \\
       \end{block}
\end{blockarray}.\]
Sada mo\v zemo navesti teoreme.
\begin{theorem}[Kronecker$^8$--Capelli$^9$]\index{sistemi!Kronecker--Capellijeva teorema}\index{teorema! Kronecker--Capellijeva}
\label{thkroneker1}
  Sistem \eqref{sistem1} je saglasan ako i samo ako je rang matrice sistema $A$ jednak rangu pro\v sirene matrice sistema $A| _B,$ tj. \[\rang A=\rang A| _B.\]
\label{sistem4}
\end{theorem}
\footnotetext[8]{Leopold Kronecker (7. decembar 1823.--29. decembar 1891.) bio je njema\v cki matemati\v car koj je radio na teoriji brojeva, algebri i logici.}
\footnotetext[9]{Alfredo Capelli (5. august 1855.--28. januar 1910.) bio je italijanski matemati\v car.}
\begin{proof}[Dokaz]
Ako je $\rang A=\rang A| _B ,$ tada matrica $A$ i pro\v sirena matrica $ A| _B$ imaju isti maksimalan broj linearno nezavisnih kolona. Odatle slijedi da
kolona $B$ slobodnih \v clanova sistema linearno zavisi od kolona matrice $A,$ pa se mo\v ze prikazati kao linearna kombinacija kolona matrice $A$ sa odgovaraju\'{c}im koeficijentima $\alpha _1, \alpha _2, \ldots , \alpha _n.$ Dakle, postoji $n$--torka skalara $(\alpha _1, \alpha _2, \ldots , \alpha _n)$
za koju vrijedi
\begin{equation}
 \begin{blockarray}{c}
   \begin{block}{(c)}
     b_1\\b_2\\ \vdots\\b_m  \\
   \end{block}
\end{blockarray}= \alpha _1 \cdot \begin{blockarray}{c}
   \begin{block}{(c)}
     a_{11}\\a_{21}\\ \vdots\\a_{m1}  \\
   \end{block}
\end{blockarray}+ \alpha _2 \cdot \begin{blockarray}{c}
   \begin{block}{(c)}
     a_{12}\\a_{22}\\ \vdots\\a_{m2}  \\
   \end{block}
\end{blockarray}+ \ldots +\alpha _n \cdot \begin{blockarray}{c}
   \begin{block}{(c)}
     a_{1n}\\a_{2n}\\ \vdots\\a_{mn}  \\
   \end{block}
\end{blockarray}
\label{trokal}
\end{equation}

\begin{equation}
 \begin{blockarray}{c}
   \begin{block}{(c)}
     b_1\\b_2\\ \vdots\\b_m  \\
   \end{block}
\end{blockarray}= \begin{blockarray}{c}
       \begin{block}{(c)}
          a_{11}\alpha _1 +  a_{12} \alpha_2 + \ldots + a_{1n}\alpha_n\\
          a_{21}\alpha _1+ a_{22}\alpha_2 + \ldots + a_{2n}\alpha_n\\
          \vdots \\
          a_{m1}\alpha_1 +a_{m2}\alpha_2 +\ldots +a_{mn}\alpha_n \\
       \end{block}
\end{blockarray}.
\label{krokal}
\end{equation}
Na osnovu jednakosti matrica iz \eqref{krokal} dobijamo: 
\begin{align}
    a_{11}\alpha_1+a_{12}\alpha_2+\ldots+a_{1n}\alpha_n&=b_1\nonumber\\
    a_{21}\alpha_1+a_{22}\alpha_2+\ldots+a_{2n}\alpha_n&=b_2\label{istemip}\\
    \vdots\hspace{2cm}&\nonumber\\
    a_{m1}\alpha_1+a_{m2}\alpha_2+\ldots+a_{mn}\alpha_n&=b_m,\nonumber
 \end{align}

Kad uporedimo \eqref{istemip} sa sistemom \eqref{sistem1}, vidimo da je jedno rje\v senje tog sistema upravo
\[x_1=\alpha _1, \, x_2= \alpha_2, \ldots , x_n=\alpha _n ,\]
pa je sistem saglasan.
Ovim smo pokazali da iz $\rang A=\rang A| _B $ slijedi saglasnost sistema \eqref{sistem1}.\\
Doka\v zimo da vrijedi i obrnuto. Pretpostavimo da je sistem \eqref{sistem1} saglasan. To zna\v ci da ima bar jedno rje\v senje:
\[x_1=\alpha _1, \, x_2= \alpha_2, \ldots , x_n=\alpha _n ,\]
pa koriste\'{c}i zapis sistema u matri\v cnom obliku \eqref{sistem3} dobijamo da vrijedi \eqref{trokal}. Ova relacija je ekvivalentna sa
\[ \begin{blockarray}{c}
   \begin{block}{(c)}
     b_1\\b_2\\ \vdots\\b_m  \\
   \end{block}
\end{blockarray}- \alpha _1 \cdot \begin{blockarray}{c}
   \begin{block}{(c)}
     a_{11}\\a_{21}\\ \vdots\\a_{m1}  \\
   \end{block}
\end{blockarray}- \alpha _2 \cdot \begin{blockarray}{c}
   \begin{block}{(c)}
     a_{12}\\a_{22}\\ \vdots\\a_{m2}  \\
   \end{block}
\end{blockarray}- \ldots - \alpha _n \cdot \begin{blockarray}{c}
   \begin{block}{(c)}
     a_{1n}\\a_{2n}\\ \vdots\\a_{mn}  \\
   \end{block}
\end{blockarray}= \begin{blockarray}{c}
   \begin{block}{(c)}
     0 \\ 0 \\ \vdots\\ 0 \\
   \end{block}
\end{blockarray}. \]
Zato, ako koloni slobodnih \v clanova pro\v sirene matrice $A| _B $ dodamo prvu kolonu pomno\v zenu sa $-\alpha _1,$ drugu kolonu pomno\v zenu
sa $-\alpha _2,\ldots,$ $n$--tu kolonu pomno\v zenu sa $-\alpha _n,$ dobi\'{c}emo da je
\[A| _B \sim \begin{blockarray}{ccccc}
       \begin{block}{(ccccc)}
          a_{11}&a_{12}&\ldots&a_{1n}&0\\
          a_{21}&a_{22}&\ldots&a_{2n}&0\\
          &\vdots&&\\
          a_{m1}&a_{m2}&\ldots&a_{mn}&0 \\
       \end{block}
\end{blockarray}. \]
Ovo zna\v ci da je $\rang A=\rang A| _B ,$ \v sto je i trebalo dokazati.
\end{proof}
Prethodna teorema ka\v ze kako da utvrdimo da li je sistem saglasan ili nije, drugim rije\v cima ima li rje\v senje/rje\v senja ili nema. Me\dj utim, ova teorema ne ka\v ze ni\v sta o tome, ako je sistem saglasan, da li ima jedno ili vi\v se rje\v senja. Odgovor na ovo pitanje daje sljede\'ca teorema.

\begin{theorem}\label{thkroneker2}\index{teorema!o broju rje\v senja sistema lin.alg.jed.}
 Ako je $\rang A=\rang A| _B=r,$ tada sistem
   \begin{enumerate}
      \item ima jedinstveno rje\v senje u slu\v caju $r=n;$
      \item beskona\v cno mnogo rje\v senja u slu\v caju $r<n$ $(n$ je broj nepoznatih u \eqref{sistem1}$)$.
   \end{enumerate}
\end{theorem}
Prethodne dvije teoreme mo\v zemo "testirati" i na sistemima \eqref{sistem2}, za koje smo utvrdili na dva na\v cina da li su saglasni i koliko imaju rje\v senja.

\begin{example}
  Posmatrajmo ponovo sisteme \eqref{sistem2}
  \begin{align*}
    x+y=2&&x+y=2&&x+y=2\nonumber\\
    x-y=0&&2x+2y=4&&x+y=3.\label{sistem03}
  \end{align*}
Odrediti egzistenciju i broj rje\v senja koriste\' ci Teoreme \ref{thkroneker1} i \ref{thkroneker2}.\\\\
\noindent Rje\v senje:\\\\
Zaklju\v cili smo da "lijevi" sistem ima jedinstveno rje\v senje $x=1,\,y=1,$ "srednji"  sistem ima beskona\v cno mnogo rje\v senja, dok "desni" sistem nije saglasan, tj. nema rje\v senja. Zapi\v simo pro\v sirenu matricu sistema "lijevog" sistema i  odredimo njen rang
\begin{equation}
  \begin{blockarray}{rrrc}
     \begin{block}{(rr|r)c}
     1&1&2\\
     1&-1&0&\text{ IIv-Iv}\\
     \end{block}
\end{blockarray}\sim
\begin{blockarray}{rrr}
     \begin{block}{(rr|r)}
     {\color{magenta}\circled{1}}&1&{\mathbf{\color{blue}2}}\\
     0&{\color{magenta}\circled{-2}}&{\mathbf{\color{blue}-2}}\\
     \end{block}
\end{blockarray}.
\label{sistem4}
\end{equation}
\index{sistemi!pro\v sirena matrica sistema}
Pro\v sirena matrica sistema $A|_B$ formirana je od matrice sistema $A$ sa dodatkom matrice kolone slobodnih \v clanova $B.$ Zbog toga iz posljednje matrice (na desnoj strani) u \eqref{sistem4} mo\v zemo odrediti i $\rang A$ i $\rang A| _B.$ Bez kolone sa {\color{blue}plavim}/\textbf{bold} elementima je matrica istog ranga kao matrica sistema $A,$ pa je $\rang A=2$ (zakru\v zeni elementi obojeni u {\color{magenta}magentu}/zaokru\v zeno), ako dodamo kolonu elemenata sa elementima matrice $B$ (elementi obojeni u {\color{blue}plavo}/\textbf{bold}) rang se ne mijenja, tj. $\rang A|_B=2.$ Pa vrijedi
\[\rang A=\rang A|_B=2=r,\]
pa je po Teoremi \ref{thkroneker1} "lijevi" sistem saglasan, a kako je i $n=2,$ (broj nepoznatih), tj.
\[n=r,\]
to po Teoremi \ref{thkroneker2} "lijevi" sistem ima jedinstveno rje\v senje. Posmatrajmo sada "srednji" sistem i njegovu pro\v sirenu matricu sistema, vrijedi
\begin{equation}
\begin{blockarray}{rrrr}
   \begin{block}{(rr|r)r}
    1&1&2&\\
    2&2&4&\text{ IIv-2$\cdot$Iv}\\
   \end{block}
\end{blockarray}
\sim
\begin{blockarray}{rrc}
   \begin{block}{(rr|c)}
    {\color{magenta}\circled{1}}&1&{\color{blue}\textbf{2}}\\
    0&0&{\color{magenta}\circled{  { \color{blue}\textbf{0} }  }}\\
   \end{block}
\end{blockarray}
\end{equation}
Sada je $\rang A=\rang A|_B=1=r$ i $r=1<2=n,$ pa je sistem saglasan, ali ima beskona\v cno mnogo rje\v senja. I na kraju za "desni" sistem vrijedi
\begin{equation}
  \begin{blockarray}{rrrr}
    \begin{block}{(rr|r)r}
     1&1&2&\\
     1&1&3&\text{ IIv-Iv} \\
    \end{block}
  \end{blockarray}
  \sim
 \begin{blockarray}{rrc}
    \begin{block}{(rr|c)}
      {\color{magenta}\circled{1}}&1&{\color{blue}\textbf{2}}\\
      0&0&{\color{magenta}\circled{ { \color{blue} \textbf{1} } }}   \\
     \end{block}
 \end{blockarray}
 \label{sistem5}
\end{equation}
Vrijedi $\rang A=1$ ali $\rang A|_B=2,$ tj. $\rang A\neq \rang A|_B,$ te desni sistem nema rje\v senja, tj. nesaglasan je.
\end{example}
\begin{remark}[Tuma\v cenja za\v sto je ovaj sistem nesaglasan -- posljednji/desni]
Kako smo samo kod ovih sistema primjenjivali elementarne transformacije  na vrste to su u desnoj matrici iz \eqref{sistem5}  elementi vrsta  zapravo koeficijenti uz nepoznate i slobodni \v clanovi, pa tako druga vrsta iz desne matrice predstavlja jedna\v cinu
\[0\cdot x+0\cdot y=1.\]
Sada se postavlja pitanje: Koje brojeve treba uvrstiti umjesto $x$ i $y$ da bi lijeva strana prethodne jednakosti imala vrijednost 1? Znamo da takvi brojevi ne postoje,  jedna\v cina nije saglasna pa nije ni "desni" sistem.
\end{remark}
\section[M\lowercase{etode za rje\v savanje}]{Metode za rje\v savanje sistema linearnih jedna\v cina}
Sada znamo kako utvrditi za dati sistem da li je saglasan ili nesaglasan, te ako je saglasan da li ima jedno ili beskona\v cno mnogo rje\v senja. Sljede\'ci korak, u slu\v caju saglasnog sistema, je izra\v cunati to rje\v senje (ili rje\v senja).  Metode koje \' ce biti ovdje obra\dj ene su Gaussova\footnote{Johann Carl Friedrich Gauss (30. april 1777.--23. februar  1855.) bio je njema\v cki matemati\v car,  smatra se najve\' cim matemati\v carem ili jednim od najve\' cih u istoriji \v covje\v canstva}, Cramerova\footnote{Gabriel Cramer (31. juli 1704.--4. januar 1752.) bio je \v svajcarski matemati\v car} (metoda determinanti) i matri\v cna metoda.
\subsection{Gaussova metoda}\index{sistemi!Gaussova metoda}
Sistem
 \begin{align*}
    a_{11}x_1+a_{12}x_2+\ldots+a_{1n}x_n&=b_1\nonumber\\
    a_{21}x_1+a_{22}x_2+\ldots+a_{2n}x_n&=b_2\nonumber\\
    \vdots\hspace{2cm}&\nonumber\\
    a_{m1}x_1+a_{m2}x_2+\ldots+a_{mn}x_n&=b_m,\nonumber
 \end{align*}
zapi\v simo u obliku pro\v sirene matrice sistema
\[A|_B=
 \begin{blockarray}{rrrrr}
    \begin{block}{(cccc|c)}
      a_{11}&a_{12}&\ldots&a_{1n}&b_1\\
      a_{21}&a_{22}&\ldots&a_{2n}&b_2\\
      a_{31}&a_{32}&\ldots&a_{3n}&b_3\\
      &\vdots&&& \vdots   \\
      a_{m1}&a_{m2}&\ldots&a_{mn}&b_m\\
    \end{block}
 \end{blockarray}.
\]

Cilj nam je matricu $A|_B$ svesti na odgovaraju\' cu trapeznu ili stepenu formu matrice. Da bi to uradili trasformisa\' cemo datu matricu u matricu kod koje je u prvoj koloni samo element $a_{11}$ razli\v cit od nule, dok su svi ostali jednaki nuli. U slu\v caju da je $a_{11}=0,$ zamijeni\'cemo prvu i $i$--tu vrstu, kod je je $a_{i1}\neq 0,$ dakle $i$--ta vrsta postaje prva, a prva postaje $i$--ta vrsta.  Sada nastavljamo postupak u anuliranju elemenata ispod elementa $a_{11}.$ Ovo radimo tako \v sto elemente druge vrste pomno\v zimo sa $a_{11},$ pa od njih oduzmeno elemente prve vrste pomno\v zene sa $a_{21},$ zatim elemente tre\' ce vrste pomno\v zimo sa $a_{11}$ i od njih oduzimamo elemente prve vrste  pomno\v zene sa $a_{31}.$ Postupak nastavljamo dok ne dobijemo matricu kod koje je samo element $a_{11}$ u prvoj koloni razli\v cit od nule. Dobijamo sljede\' cu matricu
\[
 \begin{blockarray}{rrrrrr}
    \begin{block}{(ccccc|c)}
      a_{11}&a_{12}&a_{13}&\ldots&a_{1n}&b_1\\
      0&a^{(1)}_{22}&a^{(1)}_{23}&\ldots&a^{(1)}_{2n}&b^{(1)}_2\\
      0&a^{(1)}_{32}&a^{(1)}_{33}&\ldots&a^{(1)}_{3n}&b^{(1)}_3\\
      &\vdots&&&& \vdots   \\
      0&a^{(1)}_{m2}&a^{(1)}_{m3}&\ldots&a^{(1)}_{mn}&b^{(1)}_m\\
    \end{block}
 \end{blockarray},
\]
sa $a^{(1)}_{ij}$ i $b^{(1)}_i,\,i=2,\ldots,m,\:j=2,\ldots,n$ su ozna\v ceni elementi matrice dobijeni poslije oduzimanja odgovaraju\' cih vrsta.

Zatim prelazimo na drugu kolonu. U drugoj koloni ostavljamo prva dva elementa razli\v cita od nule, "ispod"  elementa $a^{(1)}_{22}$ formiramo nule na ve\'c opisani na\v cin, te dobijamo

\[
 \begin{blockarray}{rrrrrrr}
    \begin{block}{(cccccc|c)}
      a_{11}&a_{12}&a_{13}&a_{14}&\ldots&a_{1n}&b_1\\
      0&a^{(1)}_{22}&a^{(1)}_{24}&a^{(1)}_{23}&\ldots&a^{(1)}_{2n}&b^{(1)}_2\\
      0&0&a^{(2)}_{33}&a^{(2)}_{34}&\ldots&a^{(2)}_{3n}&b^{(2)}_3\\
      &\vdots&&&&& \vdots   \\
      0&0&a^{(2)}_{m3}&a^{(2)}_{m4}&\ldots&a^{(2)}_{mn}&b^{(2)}_m\\
    \end{block}
 \end{blockarray},
\]
u sljede\'cem koraku dobijamo matricu
\[
\begin{blockarray}{rrrrrrr}
\begin{block}{(cccccc|c)}
a_{11}&a_{12}&a_{13}&a_{14}&\ldots&a_{1n}&b_1\\
0&a^{(1)}_{22}&a^{(1)}_{24}&a^{(1)}_{23}&\ldots&a^{(1)}_{2n}&b^{(1)}_2\\
0&0&a^{(2)}_{33}&a^{(2)}_{34}&\ldots&a^{(2)}_{3n}&b^{(2)}_3\\
&\vdots&&&&& \vdots   \\
0&0&0&a^{(3)}_{m4}&\ldots&a^{(3)}_{mn}&b^{(3)}_m\\
\end{block}
\end{blockarray}.
\]
Dakle ispod elemenata $a_{ii},\,i=1,\ldots,n,$ trebamo dobiti nule. Postupak ponavljamo dok ne dobijemo jedan od sljede\'cih oblika matrice $A|_B:$

\[
\begin{blockarray}{rrrrrrr}
\begin{block}{(cccccc|c)}
a_{11}&a_{12}&a_{13}&a_{14}&\ldots&a_{1n}&b_1\\
0&a^{(1)}_{22}&a^{(1)}_{24}&a^{(1)}_{23}&\ldots&a^{(1)}_{2n}&b^{(1)}_2\\
0&0&a^{(2)}_{33}&a^{(2)}_{34}&\ldots&a^{(2)}_{3n}&b^{(2)}_3\\
0&0&0&a^{(3)}_{44}&\ldots&a^{(3)}_{4n}&b^{(3)}_4\\
&\vdots&&&&& \vdots   \\
0&0&0&0&\ldots&a^{(n-1)}_{mn}&b^{(n-1)}_m\\
\end{block}
\end{blockarray},\quad (\text{ovdje je $m=n$})
\]
\[
\begin{blockarray}{rrrrrrr}
\begin{block}{(cccccc|c)}
a_{11}&a_{12}&a_{13}&a_{14}&\ldots&a_{1n}&b_1\\
0&a^{(1)}_{22}&a^{(1)}_{24}&a^{(1)}_{23}&\ldots&a^{(1)}_{2n}&b^{(1)}_2\\
0&0&a^{(2)}_{33}&a^{(2)}_{34}&\ldots&a^{(2)}_{3n}&b^{(2)}_3\\
0&0&0&a^{(3)}_{44}&\ldots&a^{(3)}_{4n}&b^{(3)}_4\\
&\vdots&&&&& \vdots   \\
0&0&0&0&\ldots&a^{(n-1)}_{nn}&b^{(n-1)}_n\\
0&0&0&0&\ldots&0&0\\
&\vdots&&&&& \vdots   \\
0&0&0&0&\ldots&0&0\\
\end{block}
\end{blockarray},
\]

ili

\[
\begin{blockarray}{rrrrrrrrr}
\begin{block}{(cccccccc|c)}
a_{11}&a_{12}&a_{13}&a_{14}&\ldots&a_{1p}&\ldots&a_{1n}&b_1\\
0&a^{(1)}_{22}&a^{(1)}_{24}&a^{(1)}_{23}&\ldots&a^{(1)}_{2p}&\ldots&a^{(1)}_{2n}&b^{(1)}_2\\
0&0&a^{(2)}_{33}&a^{(2)}_{34}&\ldots&a^{(2)}_{3p}&\ldots&a^{(2)}_{3n}&b^{(2)}_3\\
0&0&0&a^{(3)}_{44}&\ldots&a^{(3)}_{4p}&\ldots&a^{(3)}_{4n}&b^{(3)}_4\\
&\vdots&&&   &&&& \vdots   \\
0&0&0&0&\ldots&a^{(r-1)}_{rr}&\ldots&a^{(r-1)}_{rn}&b^{(r-1)}_r\\
0&0&0&0&\ldots&0&\ldots&0&0\\
&\vdots&&&&& && \vdots   \\
0&0&0&0&\ldots&0&\ldots&0&0\\
\end{block}
\end{blockarray}.
\]

Kako je prethodno opisani postupak kori\v sten i za odre\dj ivanje ranga pro\v sirene matrice sistema, odredimo sada rang matrice sistema $A$  i rang pro\v sirene matrice sistema $A|_B$.
U slu\v caju saglasnosti sistema i ako je $n=r,$ iz dobijene matrice rekonstrui\v semo sistem, te iz posljednje jedna\v cine izra\v cunamo posljednju nepoznatu, te njenu vrijednost uvrstimo u prethodnu jedna\v cinu. Ovaj postupak nastavimo dok ne izra\v cunamo i prvu nepoznatu iz prve jedna\v cine.

U slu\v caju  saglasnosti sistema i $r<n,$ na desnu stranu prebacujemo  $n-r$ nepoznatih i ponavljamo prethodno opisani postupak.

Gaussovom metodom mo\v zemo rje\v savati kako pravougaone tako i kvadratne sisteme. Sljede\' cim metodama mo\v zemo rje\v savati samo kvadratne sisteme linearnih jedna\v cina.

\subsection{Metode za rje\v savanje kvadratnih sistema}
Prva metoda za rje\v savanje kvadratnih SLAJ je matri\v cna metoda.  Kvadratni sistem
 \begin{align}
    a_{11}x_1+a_{12}x_2+\ldots+a_{1n}x_n&=b_1\nonumber\\
    a_{21}x_1+a_{22}x_2+\ldots+a_{2n}x_n&=b_2\label{sistem7}\\
    \vdots\hspace{2cm}&\nonumber\\
    a_{n1}x_1+a_{n2}x_2+\ldots+a_{nn}x_n&=b_n,\nonumber
 \end{align}
mo\v zemo zapisati u matri\v cnom obliku
\[AX=B,\]
gdje je
\[A=\begin{blockarray}{cccc}
       \begin{block}{(cccc)}
          a_{11}&a_{12}&\ldots&a_{1n}\\
          a_{21}&a_{22}&\ldots&a_{2n}\\
          &\vdots&&\\
          a_{n1}&a_{n2}&\ldots&a_{nn} \\
       \end{block}
\end{blockarray},\quad
X=\begin{blockarray}{c}
       \begin{block}{(c)}
        x_1\\x_2\\\vdots\\x_n  \\
       \end{block}
\end{blockarray},\quad
B=\begin{blockarray}{c}
   \begin{block}{(c)}
     b_1\\b_2\\ \vdots\\b_n  \\
   \end{block}
\end{blockarray}.
\]
Vrijedi sljede\' ca teorema.
\begin{theorem}\index{teorema! o broju rje\v senja kvadratnog sistema}
 Kvadratni sistem  linearnih jedna\v cina \eqref{sistem7} ima jedinstveno rje\v senje ako je $A$ regularna matrica.
\end{theorem}
\index{sistemi!matri\v cna metoda}

Ako po\dj emo od matri\v cnog zapisa sistema linearnih jedna\v cina $AX=B,$ a  po\v sto je $A$ regularna matrica, to ona ima jedinstvenu inverznu matricu (Teorema \ref{regularna1}). Vrijedi
\begin{empheq}[box=\mymath]{equation*}
AX=B\Leftrightarrow A^{-1}AX=A^{-1}B\Leftrightarrow X=A^{-1}B.
\end{empheq}
\begin{remark}
Dakle, rje\v savanje sistema linearnih jedna\v cina matri\v cnom  metodom sastoji se\\
u sljede\' cem:
\begin{enumerate}
  \item Zapi\v semo kvadratni sistem u matri\v cnom obliku;
  \item Ispitamo da li je $A$ regularna matrica;
  \item Ako jeste izra\v cunamo $A^{-1};$
  \item Rje\v senje sistema ra\v cunamo iz jednakosti $X=A^{-1}B.$
\end{enumerate}
\end{remark}

\index{sistemi!Cramerova metoda}

Sljede\' ca metoda je Cramerova\footnote{Gabriel Cramer (31.juli 1704.--4.januar 1752. godine) bio je \v zenevski matemati\v car}  metoda (metoda determinanti). Ova je metoda za-\\snovana na narednoj teoremi.
\begin{theorem}[Cramerova teorema]\label{thcramer}\index{teorema! Cramerova}
 Ako je $\det A\neq 0,$ sistem \eqref{sistem7} ima jedinstveno rje\v senje $(x_1,x_2,\ldots,x_n)$ dato sa
 \[x_k=\frac{\det A_k}{\det A},\:k=1,2,\ldots,n,\]
 gdje je $\det A_k$ determinanta matrice $A_k$ koja je dobijena zamjenom $k$--te kolone matrice $A$ matricom kolonom $B.$
\end{theorem}
\begin{proof}[Dokaz]
Za $k=1,2, \ldots, n$ je
\[\det A_k =\left| \begin{array}{ccccccc}
                a_{11}&\ldots & a_{1 \ k-1}& b_1 &a_{1\ k+1}\ldots & a_{1n}\\
                a_{21}&\ldots & a_{2 \ k-1}& b_2 &a_{2\ k+1}\ldots & a_{2n}\\
                &\vdots&&\\
                a_{n1}&\ldots & a_{n \ k-1}& b_n &a_{n\ k+1}\ldots & a_{nn}
                 \end{array}
      \right|,
\]
pa je Laplaceov razvoj determinante $\det A_k$ po $k$--toj koloni jednak
\begin{equation}
\det A_k=b_1A_{1k}+b_2A_{2k}+ \ldots +b_n A_{nk}.
\label{ovoDe}
\end{equation}
Ako prvu jedna\v cinu sistema \eqref{sistem7} pomno\v zimo sa $A_{1k},$ drugu jedna\v cinu sistema pomno-\\\v zimo sa $A_{2k},$  itd.,
posljednju jedna\v cinu sistema pomno\v zimo sa $A_{nk},$ a zatim tako dobijene jedna\v cine saberemo, tada dobijamo jedna\v cinu \\
$(a_{11}A_{1k}+a_{21}A_{2k}+ \ldots a_{n1}A_{nk})x_1+(a_{12}A_{1k}+a_{22}A_{2k}+ \ldots a_{n2}A_{nk})x_2+ \ldots +
(a_{1n}A_{1k}+a_{2n}A_{2k}+ \ldots a_{nn}A_{nk})x_n = b_1A_{1k}+b_2A_{2k}+ \ldots +b_n A_{nk}.$ \\
Iz relacije \eqref{ovoDe} vidimo da je desna strana posljednje jednakosti jednaka $\det A_k.$ Na lijevoj strani jednakosti uz nepoznatu
$x_k$ u zagradi stoji $(a_{1k}A_{1k}+a_{2k}A_{2k}+ \ldots a_{nk}A_{nk})$, a to je jednako determinanti $\det A$ (Laplaceov razvoj po determinante
po $k$--toj koloni). Svi ostali koeficijenti uz ostale nepoznate $x_i$ ($i \neq k$) na lijevoj strani su jednaki $0.$ Prema tome dobijamo
da vrijedi
\[ \det A \cdot x_k=\det A_k, \, \left(\forall k\in \{1,2,\ldots, n\} \right).  \]
Kako je $\det A \neq 0,$ to slijedi da je
\[x_k=\frac{\det A_k}{\det A}, \, \left( k\in \{1,2,\ldots, n\} \right).  \]
\end{proof}
\ \\
\begin{remark}

Cramerova metoda rje\v savanja sistema linearnih jedna\v cina sastoji se u sljede\' cem:
\begin{enumerate}
 \item Izra\v cunamo $\det A.$ Ako je $\det A\neq 0,$ nastavljamo dalje.
 \item Izra\v cunamo determinante $\det A_k,\,k=1,\ldots,n,$ po uputstvu iz Teoreme \ref{thcramer}.
 \item Nepoznate ra\v cunamo po formuli $x_k=\dfrac{\det A_k}{\det A},\:k=1,2,\ldots,n.$
\end{enumerate}
\end{remark}

\ \\

\begin{example}
 Dat je sistem \[\left\{\begin{array}{l}x+y+z=3\\2x+3y-z=4\\-x+2y+z=2\\3x+y-3z=1.	
\end{array}\right.\]
\begin{enumerate}[$(a)$]
	\item Ispitati saglasnost Kroneker-Capellijevom teoremom, i u slu\v caju saglasnosti
    \item Rije\v siti sistem Gaussovom, Cramerovom i matri\v cnom metodom.
\end{enumerate}\ \\
\noindent Rje\v senje:\\\\
\begin{enumerate}[$(a)$]
\item Saglasnost
 \begin{align*}
    A|_B&=\begin{blockarray}{rrrrr}
         \begin{block}{(rrrr)r}
           1&1&1&3&\\
           2&3&-1&4&\text{ IIv--2Iv}\\
           -1&2&1&2&\text{ IIIv+Iv}\\
           3&1&-3&1&\text{ IVv--3Iv}\\
         \end{block}
    \end{blockarray}
    \sim
     \begin{blockarray}{rrrrr}
         \begin{block}{(rrrr)r}
           1&1&1&3&\\
           0&1&-3&-2&\\
           0&3&2&5&\text{ IIIv--3IIv}\\
           0&-2&-6&-8&\text{ IVv+2IIv}\\
         \end{block}
    \end{blockarray} \\
    &\sim
    \begin{blockarray}{rrrrr}
       \begin{block}{(rrrr)r}
           1&1&1&3&\\
           0&1&-3&-2&\\
           0&0&11&11&\\
           0&0&-12&-12&\text{ 11IVv+12IIIv}\\
       \end{block}
    \end{blockarray}
    \sim
     \begin{blockarray}{rrrrr}
       \begin{block}{(rrrr)r}
           1&1&1&3&\\
           0&1&-3&-2&\\
           0&0&11&11&\text{IIIv:11}\\
           0&0&0&0&\\
       \end{block}
    \end{blockarray}\\
   & \sim
         \begin{blockarray}{rrrrr}
       \begin{block}{(rrrr)r}
           1&1&1&3&\\
           0&1&-3&-2&\\
           0&0&1&1&\\
           0&0&0&0&.\\
       \end{block}
    \end{blockarray}
 \end{align*}
 Dakle, $\rang A=\rang A|_B=3=r.$ Kako je i $n=3,$   sistem je saglasan i ima jedinstveno rje\v senje.

\item Rje\v savamo sistem dobijen od koeficijenata iz posljednje matrice.
     \begin{enumerate}[(i)]
       \item Gaussova metoda
             \begin{align*}
                   x+y+z&=3\\y-3z&=-2\\z&=1\\
                   -----&----\\
                   x+y+1&=3\\
                   y-3&=-2\\
                   z&=1\\
                   -----&----\\
                   x+y&=2\\
                   y&=1\\
                   z&=1\\
                   -----&----\\
                   x+1&=1\\
                   y&=1\\
                   z&=1\\
                   -----&----\\
                   x&=1\\
                   y&=1\\
                   z&=1.
             \end{align*}

     \item Cramerova metoda
         \begin{align*}
            D&=\begin{blockarray}{rrr}
                 \begin{block}{|rrr|}
                  1&1&1\\0&1&-3\\0&0&1  \\
                 \end{block}
           \end{blockarray}=1\\
            D_x&=\begin{blockarray}{rrr}
                 \begin{block}{|rrr|}
                  3&1&1\\-2&1&-3\\1&0&1 \\
                 \end{block}
           \end{blockarray}=1\\
            D_y&=\begin{blockarray}{rrr}
                 \begin{block}{|rrr|}
                  1&3&1\\0&-2&-3\\0&1&1 \\
                 \end{block}
           \end{blockarray}=1\\
            D_z&=\begin{blockarray}{rrr}
                 \begin{block}{|rrr|}
                  1&1&3\\0&1&-2\\0&0&1  \\
                 \end{block}
           \end{blockarray}=1.
         \end{align*}
     Dakle,
        \begin{align*}
           x&=\frac{D_x}{D}=\frac{1}{1}=1\\
           y&=\frac{D_y}{D}=\frac{1}{1}=1\\
           z&=\frac{D_z}{D}=\frac{1}{1}=1.
        \end{align*}
    \item Matri\v cna metoda. Napi\v simo sistem u matri\v cnoj formi $AX=B,$ tj.
       \[\begin{blockarray}{rrr}
           \begin{block}{(rrr)}
             1&1&1\\0&1&-3\\0&0&1\\
           \end{block}
       \end{blockarray}
       \,
       \begin{blockarray}{r}
            \begin{block}{(r)}
              x\\y\\z\\
            \end{block}
        \end{blockarray}
        =
        \begin{blockarray}{r}
           \begin{block}{(r)}
              3\\-2\\1\\
           \end{block}
        \end{blockarray}
      \]
  Izra\v cunajmo inverznu matricu matrice $A=\begin{blockarray}{rrr}
           \begin{block}{(rrr)}
             1&1&1\\0&1&-3\\0&0&1\\
           \end{block}
       \end{blockarray}.$

       Vrijednost determinante je
       \[\det A=1.\]
       Kofaktori su
 { \scriptsize
  \begin{align*}
    A_{11}&=(-1)^{1+1}\begin{blockarray}{rr}\begin{block}{|rr|}1&-3\\0&1\\\end{block}\end{blockarray}=1,&&
    A_{12}=(-1)^{1+2}\begin{blockarray}{rr}\begin{block}{|rr|}0&-3\\0&1\\\end{block}\end{blockarray}=0,&&
    A_{13}=(-1)^{1+3}\begin{blockarray}{rr}\begin{block}{|rr|}0&1\\0&0\\\end{block}\end{blockarray}=0,\\
    A_{21}&=(-1)^{2+1}\begin{blockarray}{rr}\begin{block}{|rr|}1&1\\0&1\\\end{block}\end{blockarray}=-1,&&
    A_{22}=(-1)^{2+2}\begin{blockarray}{rr}\begin{block}{|rr|}1&1\\0&1\\\end{block}\end{blockarray}=1,&&
    A_{23}=(-1)^{2+3}\begin{blockarray}{rr}\begin{block}{|rr|}1&1\\0&0\\\end{block}\end{blockarray}=0,\\
    A_{31}&=(-1)^{3+1}\begin{blockarray}{rr}\begin{block}{|rr|}1&1\\1&-3\\\end{block}\end{blockarray}=-4,&&
    A_{32}=(-1)^{3+2}\begin{blockarray}{rr}\begin{block}{|rr|}1&1\\0&-3\\\end{block}\end{blockarray}=3,&&
    A_{33}=(-1)^{3+3}\begin{blockarray}{rr}\begin{block}{|rr|}1&1\\0&1\\\end{block}\end{blockarray}=1.\\
  \end{align*}
}
Matrica kofaktora je
       \[ \cof A=\begin{blockarray}{rrr}
                      \begin{block}{(rrr)}1&0&0\\-1&1&0\\-4&3&1\\
                      \end{block}
                 \end{blockarray}.
        \]
Adjungirana matrica je
        \[
           \adj A= \begin{blockarray}{rrr}
                      \begin{block}{(rrr)}1&-1&-4\\0&1&3\\0&0&1\\
                      \end{block}
                 \end{blockarray}.
        \]
Inverzna matrica $A^{-1}$ matrice $A$  je
\[A^{-1}=\frac{1}{\det A}\adj A=\frac{1}{1} \begin{blockarray}{rrr}
                      \begin{block}{(rrr)}1&-1&-4\\0&1&3\\0&0&1\\
                      \end{block}
                 \end{blockarray}= \begin{blockarray}{rrr}
                      \begin{block}{(rrr)}1&-1&-4\\0&1&3\\0&0&1\\
                      \end{block}
                 \end{blockarray}.\]
Pa je na kraju,
\[X=A^{-1}B= \begin{blockarray}{rrr}
                      \begin{block}{(rrr)}1&-1&-4\\0&1&3\\0&0&1\\
                      \end{block}
                 \end{blockarray}\,
                   \begin{blockarray}{r}
           \begin{block}{(r)}
              3\\-2\\1\\
           \end{block}
        \end{blockarray}
        =
          \begin{blockarray}{r}
           \begin{block}{(r)}
              1\\1\\1\\
           \end{block}
        \end{blockarray}.\]
        Dakle rje\v senje je $x=1,\,y=1,\,z=1$ ili $(x,y,z)=(1,1,1)$.
     \end{enumerate}
\end{enumerate}
\end{example}

\begin{example}
Dat je sistem $\left\{\begin{array}{l}x+y+z=3\\2x+3y-z=4\\x+2y-2z=1\\3x+5y-3z=5.	
\end{array}\right.$
\begin{enumerate}[(a)]
	\item Ispitati saglasnost Kroneker-Capellijevom teoremom, i u slu\v caju saglasnosti;
    \item Rije\v siti sistem Gassuovom, Cramerovom i matri\v cnom metodom.
\end{enumerate}\ \\
\noindent Rje\v senje:\\\\
Odredimo, kako i u prethodnim slu\v cajevima, prvo $\rang A$ i $\rang A_{|B},$ da bi ispitali da li je dati sistem saglasan i ukoliko jeste, koliko ima rje\v senja.
\begin{enumerate}[(a)]
 \item Saglasnost
	\[A_{|B}=\left(\begin{array}{rrrr}1&1&1&3\\2&3&-1&4\\1&2&-2&1\\3&5&-3&5	\end{array}\right)
	   \sim \left(\begin{array}{rrrr}1&1&1&3\\0&1&-3&-2\\0&1&-3&-2\\0&2&-6&-4\end{array}\right)
	   \sim \left(\begin{array}{rrrr}1&1&1&3\\0&1&-3&-2\\0&0&0&0\\0&0&0&0\end{array}\right),\] vidimo da je $\rang A=\rang A|_B=2=r,$ a po\v sto je broj nepoznatih $n=3,$
	       to je sistem saglasan i ima beskona\v cno mnogo rje\v senja.
  \item Rije\v simo sada sistem

\[ \left\{ \begin{split}x+y+z=&3\\y-3z=&-2,	\end{split}\right.\] po\v sto je $n-r=1,$ to \' cemo jednu nepoznatu prebaciti na desnu stranu, i dobijamo\\\\
\[\left\{\begin{split}x+y=&3-z\\y=&-2+3z.	
\end{split}\right.\]
Nepoznata $z$ poprima proizvoljne vrijednosti i ozna\v cava\'cemo je sa $t,\,t\in\mathbb{R}.$
\begin{enumerate}[(i)]
	\item  Gaussova metoda
	\begin{align*}
	  x+y=&3-t\\y=&-2+3t,\\
	  ----&----\\
	  x=&-y+3-t\\y=&-2+3t,\\
	  ----&----\\
	  x=&-4t+5\\y=&-2+3t,\\
	  ----&----\\
	  x=&-4t+5\\y=&3t-2\\z=&t,\:t\in\mathbb{R}.	
   \end{align*}

\item Cramerova metoda
\begin{align*}
    D&=\left|\begin{array}{cc}1&1\\0&1\end{array}\right|=1,\\ 	
    D_x&=\left|\begin{array}{cc}3-t&1\\3t-2&1\end{array}\right|=-4t+5,\\ 	
    D_y&=\left|\begin{array}{cc}1&3-t\\0&3t-2\end{array}\right|=3t-2,\\
\end{align*}
sada je
\begin{align*}
   x&=\dfrac{D_x}{D}=\dfrac{-4t+5}{1}=-4t+5,\\
   y&=\dfrac{D_y}{D}=\dfrac{3t-2}{1}=3t-2,\\
   z&=t,\:t\in\mathbb{R}.
\end{align*} \ \\
\item  Matri\v cna metoda. Napi\v simo sistem
\[\left\{\begin{split}x+y=&3-t\\y=&-2+3t,\end{split}\right.\]
u obliku
\[\left(\begin{array}{cc}1&1\\0&1\end{array}\right)
  \left(\begin{array}{c}x\\y\end{array}\right)
  =\left(\begin{array}{c}3-t\\3t-2\end{array}\right),\]

odnosno $AX=B,$ gdje je $A=\left(\begin{array}{cc}1&1\\0&1\end{array}\right),\:
                         X=\left(\begin{array}{c}x\\y\end{array}\right),\:
                         B=\left(\begin{array}{c}3-t\\3t-2 \end{array}\right),$\\\\
dalje je\\\\
$A_{11}=1,\:A_{12}=0,\:A_{21}=-1,\:A_{22}=1,$\\\\
matrica kofaktora je
\[\cof A=\left(\begin{array}{rr}1&0\\-1&1\end{array}\right),\]
adjungirana matrica, matrice $A$ je
\[\adj A=(\cof A)^T=\left(\begin{array}{rr}1&-1\\0&1\end{array}\right),\]
determinanta je $\det A=1,$ pa je
\[A^{-1}=\dfrac{1}{\det A}\adj A=\left(\begin{array}{rr}1&-1\\0&1\end{array}\right),\]
te je rje\v senje sistema
\[X=\left(\begin{array}{cc}1&-1\\0&1	\end{array}\right)
   \left(\begin{array}{c}3-t\\3t-2 \end{array}\right)
   =\left(\begin{array}{c}-4t+5\\3t-2\end{array}\right),\]
   odnosno \[ x=-4t+5,\:y=3t-2,\:z=t,\:t\in\mathbb{R}.\]	
\end{enumerate}
\end{enumerate}
\end{example}

\newpage

\begin{example} Dat je sistem $\left\{\begin{array}{l}x+y+z+w=4\\2x+3y+z-2w=3\\3x+4y+2z-w=7.	
\end{array}\right.$
\begin{enumerate}[(a)]
  \item Ispitati saglasnost Kronecker-Capellijevom teoremom, u slu\v caju saglasnosti;
  \item Rije\v siti sistem Gaussovom, Cramerovom i matri\v cnom metodom.
\end{enumerate}\ \\
\noindent Rje\v senje:\\\\
\begin{enumerate}[(a)]
	\item Saglasnost
	
	$A_{|B}=\left(\begin{matrix}1&1&1&1&4\\2&3&1&-2&3\\3&4&2&-1&7	\end{matrix}\right)
	\sim \left(\begin{matrix}1&1&1&1&4\\0&1&-1&-4&-5\\0&1&-1&-4&-5\end{matrix}\right)
	\sim \left(\begin{matrix}1&1&1&1&4\\0&1&-1&-4&-5\\0&0&0&0&0\end{matrix}\right),$\\\\\\
dakle $\rang A=\rang A_{|B}=2=r<n=4,$ sistem ima beskona\v cno mnogo rje\v senje.
\item  Rije\v simo sada sistem \[\left\{\begin{array}{c}x+y+z+w=4\\y-z-4w=-5,\end{array}\right.\]
po\v sto je $\rang A=\rang A_{|B}=2=r,$  a $n=4,$ pa je $r=2<4=n$ i $n-r=2,$ prebacimo $z$ i $w$ na desne strane jedna\v cina. Sada uzmimo $z=a,\:w=b,\:a,b\in\mathbb{R},$
      te dalje rje\v savajmo sistem
      \[\left\{\begin{array}{c}x+y=-a-b+4\\y=a+4b-5.\end{array}\right.\]
\begin{enumerate}[(i)]
	\item  Gaussova metoda
	\begin{align*}
	     x+y&=-a-b+4\\
	       y&=a+4b-5\\
	       ----&----------\\
	      x&=-(a+4b-5)-a-b+4\\
	      y&=a+4b-5\\
	      ----&---------\\
	      x&=-2a-5b+9\\
	      y&=a+4b-5\\
	      z&=a\\w&=b,\:a,b\in\mathbb{R}.	
	 \end{align*}
\\
\item  Cramerova metoda\\\\
$D=\left|\begin{array}{cc}1&1\\0&1\end{array}\right|=1$\\\\\\
$D_x=\left|\begin{array}{cc}-a-b+4&1\\a+4b-5&1\end{array}\right|=-2a-5b+9$\\\\\\
$D_y=\left|\begin{array}{cc}1&-a-b+4\\0&a+4b-5	\end{array}\right|=a+4b-5,$\\\\\\
pa je $x=\frac{D_x}{D}=-2a-5b+9,\:y=\frac{D_y}{D}=a+4b-5,\:z=a,\:w=b,\:a,b\in\mathbb{R}.$\\
\item  Matri\v cna metoda. Napi\v simo sistem \[\begin{split} x+y&=-a-b+4\\y&=a+4b-5\end{split}\] u matri\v cnom obliku
    \[ \left(\begin{array}{cc}1&1\\0&1\end{array}\right)
    \left(\begin{array}{c}x\\y\end{array}\right)
    =\left(\begin{array}{c} -a-b+4\\a+4b-5\end{array}\right).\]
    Izra\v cunajmo matricu kofaktora, vrijedi $A_{11}=1,\:A_{12}=0,\:A_{21}=-1,\:A_{22}=1,$ pa je
    \[ \cof A=\left(\begin{array}{cc}1&0\\-1&1\end{array}\right),\]
   adjungirana matrica je
   \[ \adj A=\left(\cof A\right)^T=\left(\begin{array}{cc}1&-1\\0&1\end{array}\right),\]
   determinanta matrice sistema $A$ je $\det A=1,$ te je inverzna matrica
   \[ A^{-1}=\left(\begin{array}{cc}1&-1\\0&1\end{array}\right).\]
   Rje\v senje sistema sada je
   \[ X=A^{-1}B=\left(\begin{array}{cc}1&-1\\0&1\end{array}\right)\left(\begin{array}{c} -a-b+4\\a+4b-5\end{array}\right)
        =\left(\begin{array}{c} -2a-5b+9\\a+4b-5\end{array}\right),\]
   odnosno \[ x=-2a-5b+9,\:y=a+4b-5,\:z=a,w=b,\:a,b\in\mathbb{R}.\]
\end{enumerate}
\end{enumerate}
\end{example}

 \newpage

\begin{example}
 Dat je sistem $\left\{\begin{array}{l}x+y+z+w=4\\x-y+z+2w=3\\2x+y+z-w=3\\4x+y+3z+2w=10.\end{array}\right.$
\begin{enumerate}[(a)]
	\item Ispitati saglasnost Kronecer-Capellijevom teoremom,  u slu\v caju saglasnosti;
    \item Rije\v siti sistem Gaussovom metodom.
\end{enumerate}\ \\
\noindent Rje\v senje:\\
\begin{enumerate}[(a)]
	\item  Odredimo $\rang\, A$ i $\rang\,A_{|B}.$ Vrijedi
	   \begin{multline*}
    A_{|B}=\left(\begin{array}{ccccc}1&1&1&1&4\\1&-1&1&2&3\\2&1&1&-1&3\\4&1&3&2&10\end{array}\right)
	    \sim\left(\begin{array}{ccccc}1&1&1&1&4\\0&-2&0&1&-1\\0&-1&-1&-3&-5\\0&-3&-1&-2&-6\end{array}\right)\\
	    \sim\left(\begin{array}{ccccc}1&1&1&1&4\\0&-1&-1&-3&-5\\0&0&2&7&9\\0&0&2&7&9\end{array}\right)
	    \sim\left(\begin{array}{ccccc}1&1&1&1&4\\0&-1&-1&-3&-5\\0&0&2&7&9\\0&0&0&0&0\end{array}\right).
	    \end{multline*}
Dakle $\rang A=\rang A_{|B}=3=r,$ tj. sistem ima rje\v senje, ali po\v sto sistem ima \v cetiri nepoznate $n=4$, $n=4<3=r$, rje\v senja ima beskona\v cno mnogo.
\item Polazni sistem  ekvivalentan je sljede\' cem sistemu
\begin{align*}x+y+z+w&=4\\-y-z-3w&=-5\\2z+7w&=9\\w&=a,\:a\in\mathbb{R}
                      \\--------&---------------\\
              x+y+z&=-a+4\\-y-z&=3a-5\\2z&=-7a+9\\w&=a,\:a\in\mathbb{R}
                      \\--------&---------------\\
              x+y+z&=-a+4\\-y-z&=3a-5\\z&=\dfrac{-7a+9}{2}\\w&=a,\:a\in\mathbb{R} \allowdisplaybreaks
                      \\--------&---------------\\
              x+y+\dfrac{-7a+9}{2}&=-a+4\\-y+\dfrac{7a-9}{2}&=3a-5\\z&=\dfrac{-7a+9}{2}\\w&=a,\:a\in\mathbb{R}
                      \\--------&---------------\\
              x+y&=\dfrac{7a-9}{2}-a+4\\-y&=\dfrac{-7a+9}{2}+3a-5\\z&=\dfrac{-7a+9}{2}\\w&=a,\:a\in\mathbb{R}
                      \\--------&---------------\\
              x&=\dfrac{-a-1}{2}+\dfrac{7a-9}{2}-a+4\\y&=\dfrac{a+1}{2}\\z&=\dfrac{-7a+9}{2}\\w&=a,\:a\in\mathbb{R}
                      \\--------&---------------\\
              x&=\dfrac{4a-2}{2}=2a-1\\y&=\dfrac{a+1}{2}\\z&=\dfrac{-7a+9}{2}\\w&=a,\:a\in\mathbb{R}.
\end{align*}
\end{enumerate}
\end{example}

\subsection{Homogeni sistemi linearnih algebarskih jedna\v cina}\index{sistemi!homogen}
U Definiciji \ref{sistem} ve\' c je uveden pojam homogenog sistema. Da ponovimo, to je sistem algebarskih linearnih jedna\v cina u kojem su svi slobodni \v clanovi $b_i,\,i=1,\ldots,m$ jednaki nuli, tj.
 \begin{align}
    a_{11}x_1+a_{12}x_2+\ldots+a_{1n}x_n&=0\nonumber\\
    a_{21}x_1+a_{22}x_2+\ldots+a_{2n}x_n&=0\label{sistem8}\\
    \vdots\hspace{2cm}&\nonumber\\
    a_{m1}x_1+a_{m2}x_2+\ldots+a_{mn}x_n&=0,\nonumber
 \end{align}
ili u matri\v cnoj formi \[AX=0.\]
O\v cigledno svaki homogeni sistem linearnih jedna\v cina ima rje\v senje $x_1=0,\,x_2=0,\,\ldots,\,x_n=0.$ Uvr\v stavanjem ovih vrijednosti u lijeve strane jedna\v cina sistema \eqref{sistem8}, svaka lijeva strana posmatranih jedna\v cina postaje jednaka nuli. Dobili smo ta\v cne brojne jednakosti, tj. $x_1=x_2=\,\ldots\,x_n=0$ je rje\v senje datog homogenog sistema. Ovo rje\v senje se naziva trivijalno\index{sistemi!trivijalno rje\v senje homogenog sistema} ili nulto rje\v senje. I iz pro\v sirene matrice sistema mo\v zemo lako zaklju\v citi da je sistem uvijek saglasan, tj. da ima rje\v senje
\[A| _B=\begin{blockarray}{ccccc}
       \begin{block}{(cccc|c)}
          a_{11}&a_{12}&\ldots&a_{nn}&0\\
          a_{21}&a_{22}&\ldots&a_{2n}&0\\
          &\vdots&&\\
          a_{m1}&a_{m2}&\ldots&a_{mn}&0 \\
       \end{block}
\end{blockarray},\]
jer posljednja kolona matrice sa slobodnim \v clanovima ne uti\v ce na vrijednost ranga matrice $A|_B,$ po\v sto su svi njeni elementi jednaki nuli. Osnovni problem u rje\v savanju homogenih sistema linearnih jedna\v cina je u ispitivanju da li ovaj sistem ima i netrivijalnih rje\v senja. Rje\v senje ovog problema dato je u sljede\' coj teoremi i njenim posljedicama.
\begin{theorem}\label{sistem9} \index{sistemi!netrivijalno rje\v senje\\ homogenog sistema}\index{teorema! o broju rje\v senja homogenog sistema}
  Homogeni sistem linearnih jedna\v cina \eqref{sistem8} ima netrivijalna rje\v senja ako i samo ako je rang matrice sistema manji od broja nepoznatih, tj. ako je $\rang A=r<n.$
\end{theorem}

\begin{corollary}
  Svaki homogeni sistem linearnih jedna\v cina u kome je broj jedna\v cina manji od broja nepoznatih ima netrivijalna rje\v senja.
\end{corollary}

\begin{corollary}\label{sistem10}
  Homogeni sistem u kome je broj jedna\v cina jednak broju nepoznatih ima netrivijalna rje\v senja ako i samo ako je determimanta sistema jednaka nuli.
\end{corollary}

\begin{example}
  Rije\v siti homogeni sistem $\left\{\begin{array}{c}x+2y-3z=0\\2x+5y+2z=0\\3x-y-4z=0.	\end{array}\right.$\\\\
\noindent Rje\v senje:\\\\ Izra\v cunajmo rang pro\v sirene matrice sistema
\begin{align*}
     A_{|B}= \begin{blockarray}{rrrr}\begin{block}{(rrr|r)}1&2&-3&0\\2&5&2&0\\3&-1&-4&0\\\end{block}\end{blockarray}
  \sim\begin{blockarray}{rrrr}\begin{block}{(rrr|r)}1&2&-3&0\\0&1&8&0\\0&-7&5&0\\\end{block}\end{blockarray}
  \sim\begin{blockarray}{rrrr}\begin{block}{(rrr|r)}1&2&-3&0\\0&1&8&0\\0&0&61&0\\\end{block}\end{blockarray}.
\end{align*}
Dakle vrijedi $\rang A=\rang A|_B=3,$ a kako je broj nepoznatih $n=3,$ sistem ima samo trivijalna rje\v senja. Vidimo da posljednja kolona  ne igra nikakvu ulogu u odre\dj ivanju ranga pro\v sirene matrice sistema kod homogenih sistema, pa se zato i ra\v cuna samo rang matrice sistema, kao \v sto je navedeno u Teoremi \ref{sistem9}. Da sistem ima samo trivijalna rje\v senja mo\v zemo zaklju\v citi koriste\' ci Posljedicu \ref{sistem10} jer broj jedna\v cina jednak broju nepoznatih, pa je
\[\left|\begin{array}{ccc}1&2&-3\\2&5&2\\3&-1&-4\end{array}\right|=61\neq 0.\]
Dakle postoji samo trivijalno rje\v senje $x=y=z=0.$
\end{example}
\begin{example}
 Rije\v siti homogeni sistem
 $\begin{cases} x+2y+z=0\\2x+3y+z=0\\3x+5y+2z=0\\2x+4y+2z=0.\end{cases}$\\\\
\noindent Rje\v senje:\\\\
 \begin{align*}
   A_{|B}= \begin{blockarray}{rrr}
       \begin{block}{(rrr)}1&2&1\\2&3&1\\3&5&2\\2&4&2\\\end{block}
    \end{blockarray}
    \sim
    \begin{blockarray}{rrr}
       \begin{block}{(rrr)}1&2&1\\0&-1&-1\\0&-1&-1\\0&0&0\\\end{block}
    \end{blockarray}
    \sim
    \begin{blockarray}{rrr}
       \begin{block}{(rrr)}1&2&1\\0&-1&-1\\0&0&0\\0&0&0\\\end{block}
    \end{blockarray}
 \end{align*}
 pa je $\rang A=2<3=n.$ Dakle sistem ima i netrivijalna rje\v senja.
 \begin{align*}
   x+2y+z&=0\\-y-z&=0\\-----&-----\\x+2y&=-z\\y&=-z\\z&=a,\,a\in\mathbb{R}\\-----&-----\\x+2y&=-a\\y&=-a\\z&=a\\-----&-----\\x&=a\\y&=-a\\z&=a,\,a\in\mathbb{R}.
 \end{align*}
\end{example}

\section{Z\lowercase{adaci}} \index{Zadaci za vje\v zbu!sistemi linearnih algebarskih jedna\v cina}

\begin{enumerate}
  \setcounter{enumi}{0}
    \item Rije\v siti sisteme  Gaussovom metodom\\
           \begin{inparaenum}
              \item        $
                   \begin{cases}
                       x+y+z=3\\
                       x+2y+3z=6\\
                       3x+3y+z=7;
                   \end{cases}\quad$
              \item $
                         \begin{cases}2x-y+z=2\\2x+2y+4z=8
                                         \\3x-y-2z=0.	
                         \end{cases}
                    \quad$
            \end{inparaenum}
\item Rije\v siti sisteme Cramerovom metodom\\
           \begin{inparaenum}
              \item $\begin{cases}x+3y-z=3\\-3x+2y+z=0
                                \\2x+3y-3z=2;	
                    \end{cases}\quad $
             \item $\begin{cases}x-2y+3z=2\\2x+y+4z=7\\3x-2y-2z=-1.	
                           \end{cases}
                    \quad$
            \end{inparaenum}
\item Rije\v siti sisteme matri\v cnom metodom\\
           \begin{inparaenum}
             \item $\begin{cases}-x-2y+z=1\\2x+5y-z=-2\\3x-y-2z=5;	
                           \end{cases}
                    \quad$
             \item $\begin{cases}3x-2y+3z=4\\x+2y+z=4\\2x-3y-4z=-5.	
                   \end{cases}$
            \end{inparaenum}

\item  Ispitati saglasnost sistema Kronecker-Capellijevim stavom i slu\v caju saglasnosti\\ rije\v siti ga\\
          \begin{inparaenum}
            \item  $
                   \begin{cases}
                       x+y+z=3\\
                       2x+y+z=4\\
                       -x-y+z=-1\\
                       3x-y-z=1;	
                   \end{cases} $
           \item $\begin{cases}x+y+z=3\\2x+y+z=4\\-x-y+z=-1\\3x-y-z=1;	
                   \end{cases}$\quad
           \item $\begin{cases}2x+2y+z=5\\x-y+z=1\\-2x-3y+z=-4\\x-y-z=-1.	
                 \end{cases}$

          \end{inparaenum}
\item  Ispitati saglasnost sistema Kronecker-Capellijevim stavom i slu\v caju saglasnosti\\ rije\v siti dati sistem sa sve tri metode\\
          \begin{inparaenum}
           \item $\begin{cases}x+2y+3z=7\\2x-y-5z=-4\\-x+y+z=1\\3x+y+z=5;	
                  \end{cases}$
           \item $\begin{cases}x+y+z=1\\2x+y+z=4\\-x-y+z=-1\\3x-y-z=1;	
                  \end{cases}$\quad
           \item $\begin{cases}x+4y-3z=2\\-2x-y+6z=3\\x-7y+z=-5\\x-y-z=-1.	
                  \end{cases}$
          \end{inparaenum}
\item Ispitati saglasnost sistema Kronecker-Capellijevim stavom i u slu\v caju saglasnosti rije\v siti ga\\
          \begin{inparaenum}
            \item $\begin{cases}x+y-z+t=2\\-2x-y-z+2t=-3\\2x+2y-3z-3t=4\\-4x-3y+z-t=-7; \end{cases}    \quad$
            \item $\begin{cases}x+2y+z-t=0\\-x-y+2z-2t=-3\\-2x+y-2z+t=-1\\-2x+4y-z-2t=-4. \end{cases}   $
           \end{inparaenum}

\item  Ispitati saglasnost sistema Kronecker-Capellijevim stavom i slu\v caju saglasnosti\\ rije\v siti dati sistem\\
         \begin{inparaenum}
            \item $\begin{cases}2x+7y+3z+u=6\\3x+5y+2z+2u=4\\9x+4y+z+7u=2;	
                   \end{cases}$
            \item $ \begin{cases}x-2y+z+u=1\\x-2y+z-u=-1\\x-2y+z+5u=5.	
                   \end{cases} $
         \end{inparaenum}
   \newpage
\item  Ispitati saglasnost sistema Kronecker-Capellijevim stavom i slu\v caju saglasnosti\\ rije\v siti dati sistem sa sve tri metode\\
         \begin{inparaenum}
            \item $\begin{cases}3x-y-2z+u=5\\-2x-5y+z-2u=2\\x+2y-z+u=-1\\-x+y-3z+5u=-6;	
                   \end{cases}$
            \item $\begin{cases}2x-y+z+u=1\\x+2y-z+4u=2\\x+7y-4z+11u=3.
                  \end{cases}$
         \end{inparaenum}
\item Ispitati ima li homogeni sistem rje\v senja osim trivijalnog, i slu\v caju da ima, odrediti ih\\
    \begin{inparaenum}
       \item $\begin{cases}x+y+2z=0\\-2x+2y+3z=0\\-5x+2y+z=0 ;\end{cases}     \:$
       \item $\begin{cases}x+y+2z=0\\-2x+2y+3z=0\\-x+3y+5z=0 ;\end{cases}     \:$

       \item $\begin{cases}x+2y-2z=0\\-x+3y+z=0\\-x+8y=0\\-x+13y-z=0; \end{cases}        \:$
       \item $\begin{cases}2x+y-2z=0\\-x+2y+z=0\\x+3y-z=0\\2x+6y-2z=0\\3x+4y-3z=0. \end{cases}        \:$
    \end{inparaenum}
\item Odredi realni parametar $\lambda$ tako da homogeni sistem ima netrivijalna rje\v senja i odre-\\diti ih\\
    \begin{inparaenum}
       \item $\begin{cases}x+2y-z=0\\3x-y+2z=0\\4x+y+\lambda z=0; \end{cases}     \:$
       \item $\begin{cases}2x+y+2z=0\\ \lambda x-z=0\\-x+3y+\lambda z=0. \end{cases}     $

    \end{inparaenum}
\end{enumerate}


   \chapter{ Vektori}


\index{vektori}

U matematici, prirodnim i tehni\v ckim naukama, neke veli\v cine su odre\dj ene svojom brojnom vrijedno\v s\' cu. Takve su na primjer, masa, zapremina, vrijeme i dr. Me\dj utim, postoje veli\v cine za koje nije dovoljno samo znati njihove brojne vrijednosti, kao na primjer brzina, sila, ubrzanje, ja\v cina elektri\v cnog polja, ja\v cina magnetnog polja i dr. Veli\v cine sa kojima radimo mogu se podijeliti na:
\begin{enumerate}
  \item Veli\v cine koje su definisane samo brojnom vrijedno\v s\' cu i njih nazivamo  skalarnim veli\v cinama ili kra\' ce skalarima. \index{skalari}
  \item Veli\v cine za koje je potrebno poznavati brojnu vrijednost, pravac i smjer. Ove veli\v cine se nazivaju vektorske veli\v cine ili vektori.\index{vektori}
\end{enumerate}

\section[O\lowercase{snovni pojmovi}]{Osnovni pojmovi}

\subsection{ Pojam vektora} Vektore geometrijski predstavljamo usmjerenim (ili orjentisanim) du\v zima. Neka su date ta\v cke $A,\,B,\,A	\neq B,$ prostora.   Du\v z $\overline{AB}$ kojoj je jedna ta\v cka npr. $A$ progla\v  sena  za po\v cetak, a druga  za kraj (ta\v cka $B$) zva\' cemo vektor, u oznaci $\overrightarrow{AB}$. Vektor  $\overrightarrow{AB}$ ima pravac kao i prava na kojoj le\v zi du\v z $\overline{AB},$ smjer od ta\v cke $A$ ka ta\v cki $B,$ a intenzitet (ili modul ili norma) mu je jednak du\v zini du\v zi $\overline{AB}.$         \index{vektori!osnovni pojmovi}

Vektor je definisan ako znamo njegov intenzitet, pravac i smjer. Vektori se obi\v cno ozna\v cavaju malim slovima latinice sa strelicom iznad, na primjer $\vec{a},\,\vec{b},\ldots,$  a ako su poznate po\v cetna i krajnja ta\v cka onda koristimo oznaku $\overrightarrow{AB},\,\overrightarrow{DF},\ldots $ Osim ovakvog ozna\v cavanja koriste se i mala slova latinice, ali podebljana (boldovana) $\textbf{a},\,\textbf{b},\ldots$ Intenziteti se ozna\v cavaju sa $|\vec{a}|,\,|\vec{b}|,\ldots$ ili $\left|\overrightarrow{AB}\right|,\,\left|\overrightarrow{DF}\right|,\ldots$ ili $|\textbf{a}|,\,|\textbf{b}|\ldots$ Skup svih vektora prostora ozna\v cava\' cemo sa
{\hspace{-4cm}
\begin{empheq}[box=\mymath]{equation*}
V^3.
\end{empheq}
}

 \begin{figure}[!h]\centering
     \begin{subfigure}[b]{.3\textwidth}\centering
        \includegraphics[scale=.85]{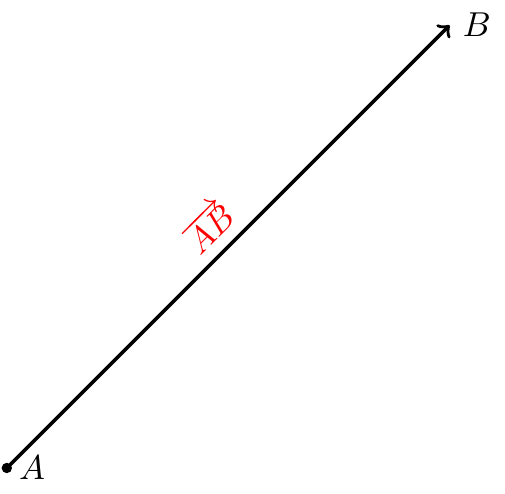}
         \caption{}
          \label{slika35}
   \end{subfigure}\hspace{.5cm}
   \begin{subfigure}[b]{.3\textwidth}\centering
        \includegraphics[scale=.85]{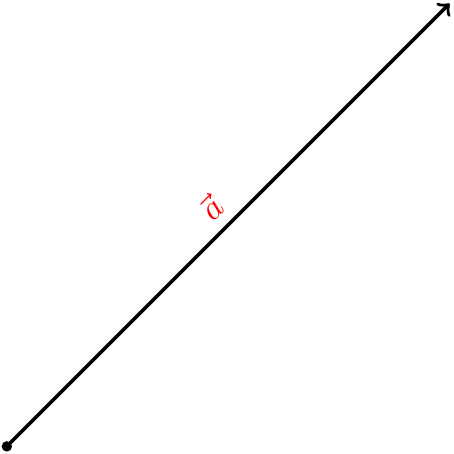}
         \caption{}
      \label{slika36}
     \end{subfigure}
     \begin{subfigure}[b]{.3\textwidth}\centering\hspace{.5cm}
        \includegraphics[scale=.85]{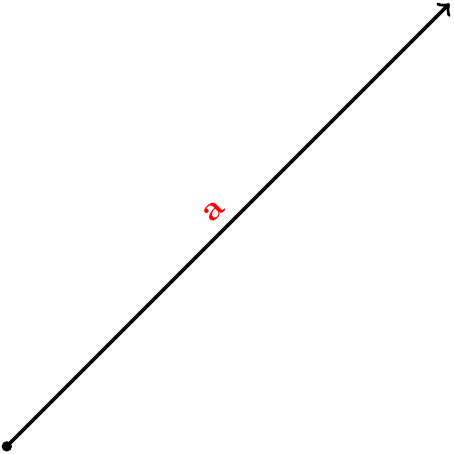}
         \caption{}
      \label{slika37}
      \end{subfigure}
        \caption{Ozna\v cavanje vektora}
     \label{slika35a}
  \end{figure}

Vektor \v ciji je intenzitet jednak nuli naziva se nula--vektor. Po\v cetna i krajnja ta\v cka nula--vektora poklapaju se (tj. iste su), dok pravac i smjer nisu definisani. Nula--vektor ozna\v cavamo $\vec{0}$ ili $\textbf{0}.$  Vektor \v ciji je intenzitet jednak $1$ zovemo jedini\v cni vektor.
\newpage

Uobi\v cajeno, vektori se dijele u tri kategorije:
\begin{enumerate}
 \item slobodni vektori;
 \item vektori vezani za pravu ili nosa\v c;
 \item vektori vezani za ta\v cku.
\end{enumerate}

\textbf{Slobodni vektori.} Za vektore $\vec{a}$ i $\vec{b}$ ka\v zemo da su jednaki, tj. $\vec{a}=\vec{b}$ ako imaju jednake intenzitete, iste smjerove a nosa\v ci su im paralelni ili se poklapaju. Ovakve vektore nazivamo slobodni. Slobodne vektore mo\v zemo paralelno pomjerati du\v z njihovog nosa\v ca ili paraleno nosa\v cu, zadr\v zavaju\' ci smjer i intenzitet nepromijenjenim.

\textbf{Vektori vezani za pravu (nosa\v c).} Vektori $\vec{a}$ i $\vec{b}$ su jednaki ako imaju iste module, smjerove i isti nosa\v c.  Ovo su vektori vezani za pravu. Mo\v zemo ih pomjerati du\v z njihovog nosa\v ca.

\textbf{Vektori vezani za ta\v cku.} Vektori $\vec{a}$ i $\vec{b}$ su jednaki ako imaju iste module, smjerove, intenzitete a krajne ta\v cke im se poklapaju. Ovi vektori su vezani za ta\v cku.\\

U nastavku \' cemo koristiti {\color{magenta}\textbf{slobodne vektore}}.

\subsection{Operacije sa vektorima}

\begin{figure}[!h]\centering
     \begin{subfigure}[b]{.3\textwidth}\centering
        \includegraphics[scale=.85]{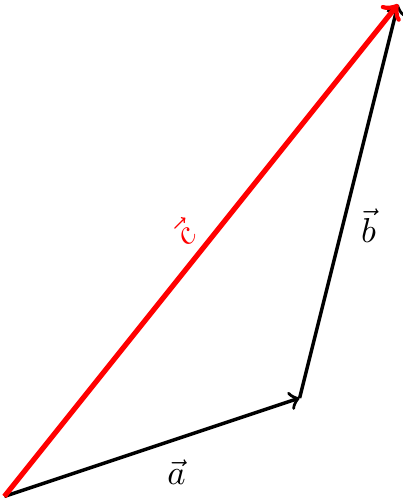}
         \caption{}
          \label{slika38}
   \end{subfigure}\hspace{.5cm}
   \begin{subfigure}[b]{.3\textwidth}\centering
        \includegraphics[scale=.85]{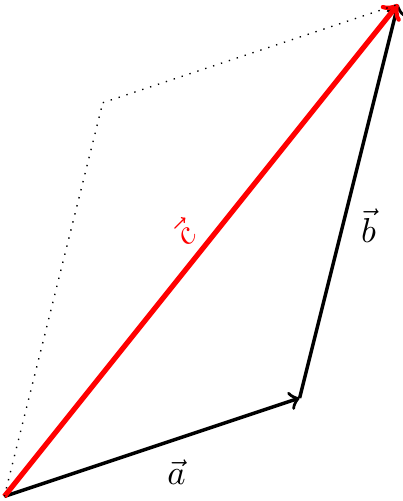}
         \caption{}
      \label{slika39}
     \end{subfigure}
     \begin{subfigure}[b]{.3\textwidth}\centering\hspace{.5cm}
        \includegraphics[scale=.85]{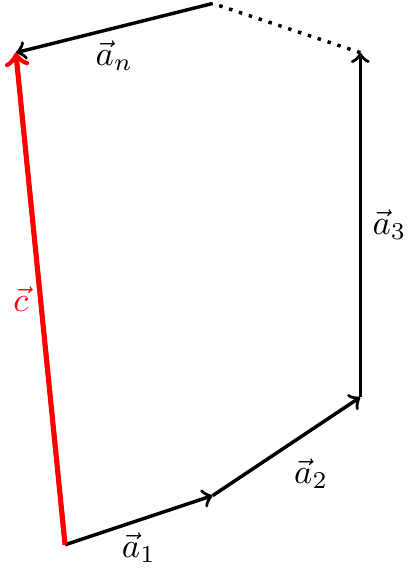}
         \caption{}
      \label{slika40}
      \end{subfigure}
        \caption{Sabiranje vektora}
     \label{slika38a}
  \end{figure}

\textbf{Sabiranje vektora.} Vektore $\vec{a}$ i $\vec{b}$ sabiramo tako \v sto po\v cetak vektora $\vec{b}$ dovedemo pa-\\ralelnim pomjeranjem do kraja vektora $\vec{a}.$ Rezultantni vektor $\vec{c}$,  \v ciji se po\v cetak poklapa sa po\v cetkom vektora $\vec{a}$ a kraj sa krajem vektora $\vec{b},$ je zbir vektora $\vec{a}$ i $\vec{b},$ tj. $\vec{c}=\vec{a}+\vec{b}.$ Ovo je pravilo trougla (Slika \ref{slika38}) ili paralelograma (nadoponjavanjem dobijemo paralelogram, vidjeti Sliku \ref{slika39} -- isprekidane linije). U slu\v caju da trebamo sabrati vektore $\vec{a}_1,\vec{a}_2,\ldots,\vec{a}_{n-1},\vec{a}_n$ postupamo na sljede\' ci na\v cin. Po\v cetak vektora  $\vec{a}_2$ dovedemo na kraj vektora $\vec{a}_1,$ po\v  cetak vektora $\vec{a}_3$ na kraj vektora $\vec{a}_2,$ postupak nastavimo do vektora $\vec{a}_n$ \v ciji po\v cetak  dovodimo do kraja vektora $\vec{a}_{n-1}.$  Sada rezultantni vektor $\vec{c},$ \v ciji se po\v cetak poklapa sa po\v cetkom vektora $\vec{a}_1$ a kraj sa krajem vektora $\vec{a}_n,$ predstavlja tra\v zeni zbir, tj. $\vec{c}=\vec{a}_1+\ldots+\vec{a}_n$ (Slika \ref{slika40}).  Izlo\v zeni postupak \v cesto se zove i pravilo (ili princip) nadovezanih vektora. \index{vektori!sabiranje}

\textbf{Oduzimanje vektora.}  Suprotan vektor vektora $\vec{b}$ je vektor koji ima isti intenzitet, isti ili paralelan nosa\v c, ali suprotan smjer od vektora $\vec{b}.$ Suprotan vektor, vektora $\vec{b}$ ozna\v cavamo sa $-\vec{b}.$ \index{vektori!oduzimanje}

Oduzimanje vektora $\vec{b}$ od vektora $\vec{a}$ svodimo na sabiranje vektora $\vec{a}$ sa suprotnim vektorom, vektora $\vec{b},$ tj. $\vec{a}-\vec{b}=\vec{a}+(-\vec{b}).$

 \begin{figure}[!h]\centering
        \includegraphics[scale=.85]{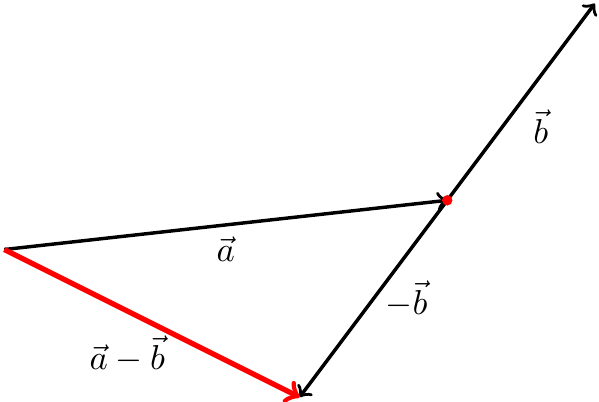}
        \caption{Oduzimanje vektora}
     \label{slika41}
  \end{figure}

\begin{figure}[!h]\centering
     \begin{subfigure}[b]{.35\textwidth}\centering
        \includegraphics[scale=.85]{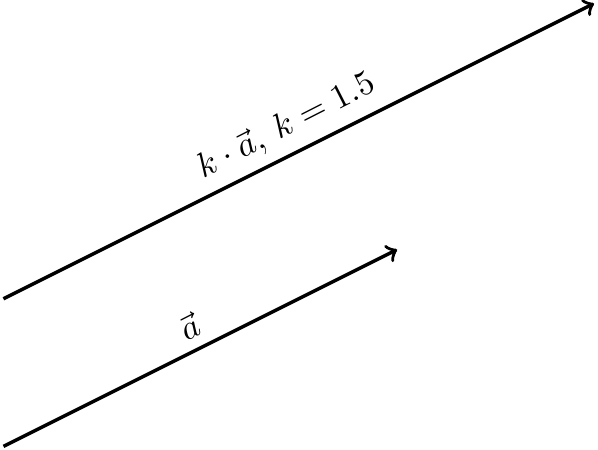}
         \caption{}
          \label{slika42}
   \end{subfigure}\hspace{.5cm}
   \begin{subfigure}[b]{.35\textwidth}\centering
        \includegraphics[scale=.85]{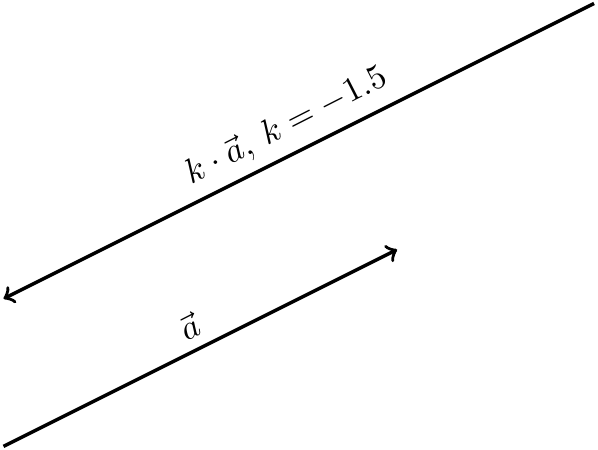}
         \caption{}
      \label{slika43}
     \end{subfigure}\\
     \begin{subfigure}[b]{.35\textwidth}\centering\hspace{.5cm}
        \includegraphics[scale=.85]{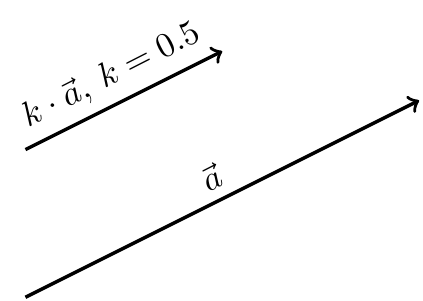}
         \caption{}
      \label{slika44}
      \end{subfigure}
    \begin{subfigure}[b]{.35\textwidth}\centering\hspace{.5cm}
        \includegraphics[scale=.85]{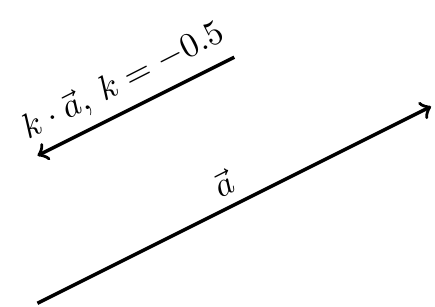}
         \caption{}
      \label{slika45}
    \end{subfigure}
        \caption{Mno\v zenje vektora skalarom}
     \label{slika38a}
  \end{figure}

\textbf{Mno\v zenje vektora skalarom.} Proizvod vektora $\vec{a}$ i skalara $k\in\mathbb{R}$ je vektor, \v ciji se pravac poklapa sa pravcem vektora $\vec{a},$ intenzitet mu je jednak $|k|\cdot|\vec{a}|,$ a smjer mu je isti kao i smjer vektora $\vec{a}$ ako je $k>0,$ a suprotan od smjera $\vec{a}$ ako je $k<0.$   Na Slici \ref{slika42} je vektor $\vec{a}$ i vektor $1.5\vec{a}.$ Vidimo da rezultantni vektor $1.5\vec{a}$ ima isti smjer i pravac kao vektor $\vec{a},$ ali mu je du\v zina (intenzitet) 1.5 puta ve\' ca od du\v zine vektora $\vec{a}.$ Na Slici \ref{slika43} je $k=-1.5<0,$ pa rezultantni vektor ima isti pravac kao i $\vec{a}$, intenzitet mu je 1.5 puta ve\' ci nego kod $\vec{a},$ ali ima suprotan smjer. Na Slikama \ref{slika44} i \ref{slika45} su slu\v cajevi kada je $|k|<1.$ Rezultantni vektor ima intenzitet manji od intenziteta vektora $\vec{a}.$ Za $k=-1$ dolazi samo do promjene smjera vektora. \index{vektori!mno\v zenje vektora skalarom}

Upravo obra\dj ene operacije uvedene su koriste\'ci vektore kao potpuno geometrijske\\ obje~kte, nisu kori\v stene nikakve formule da bismo predstavili ove vektore (osim oznaka za vektore $\vec{a},\ldots$). Ovakav pristup znatno ote\v zava i su\v zava primjenu vektora. Puno prakti\v cnije bi bilo kada bismo vektore mogli predstaviti analiti\v cki, dakle putem neke formule, a samim tim dosta bi olak\v sali kalkulacije sa vektorima i pro\v sirili bi njihovu pri-\\mjenu.  U nastavku bi\'ce  razmatrano kako predstaviti vektor nekom formulom, tj. kako ga predstaviti preko drugih vektora, koje uslove treba da ispune vektori, pa da bi druge vektore mogli preko njih izraziti.

\subsection{Linearna kombinacija vektora. Baza prostora $V^3$}

Neka je dato $n$ vektora $\vec{a}_1,\vec{a}_2,\ldots,\vec{a}_n$ i $n$ realnih skalara $\alpha_1,\alpha_2,\ldots,\alpha_n.$  Linearnom kombinacijom vektora $\vec{a}_1,\vec{a}_2,\ldots,\vec{a}_n$ nazivamo sljede\' ci zbir

\begin{empheq}[box=\mymath]{equation*}
\alpha_1\vec{a}_1+\alpha_2\vec{a}_2+\cdots+\alpha_n\vec{a}_n.
\end{empheq} \index{vektori!linearna kombinacija}
Realni skalari $\alpha_1,\ldots, \alpha_n\in\mathbb{R }$ su koeficijenti linearne kombinacije.

Linearna kombinacija nam je nezaobilazni element u predstavljanju jednog vektora preko drugog ili drugih vektora. Da bismo predstavili vektor nekom linearnom kombinacijom potrebni su nam pojmovi dati u sljede\'coj definiciji.

\begin{figure}[!h]\centering\hspace{-2cm}
	\begin{subfigure}[b]{.35\textwidth}\centering
		\includegraphics[scale=.8]{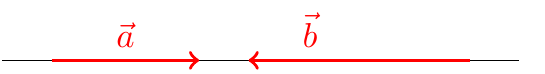}
		\caption{$\vec{a}=-1.5\vec{b}$ ili $\vec{b}=-\frac{1}{1.5}\vec{b}$}
		\label{kolinearnostSlika}
	\end{subfigure}\hspace{.5cm}
	\begin{subfigure}[b]{.35\textwidth}\centering
		\includegraphics[scale=.8]{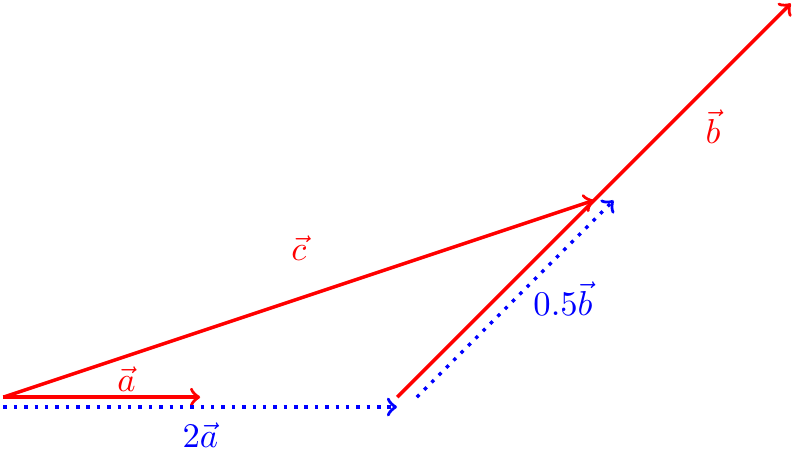}
		\caption{$\vec{c}=2\vec{a}+0.5\vec{b}$}
		\label{komplanarnotsSlika}
	\end{subfigure}\\

	\caption{Predstavljanje vektora}
	\label{zavisnostSlika}
\end{figure}

\begin{definition}[Linearno zavisni--nezavisni vektori]
Vektori $\vec{a}_1,\vec{a}_2,\ldots,\vec{a}_n$ su linearno zavisni ako postoje skalari  $\alpha_1,\alpha_2,\ldots,\alpha_n,$ od kojih je bar jedan razli\v cit od nule,
tako da je
\begin{equation}
\alpha_1\vec{a}_1+\alpha_2\vec{a}_2+\cdots+\alpha_n\vec{a}_n=\vec{0}.
\label{linzavisnost}
\end{equation}
Ako je jednakost \eqref{linzavisnost} zadovoljena jedino kada je $\alpha_1=\alpha_2=\ldots=\alpha_n=0,$ onda se ka\v ze da su vektori $\vec{a}_1,\vec{a}_2,\ldots,\vec{a}_n$ linearno nezavisni.
\end{definition}\index{vektori!zavisni--nezavisni}
Pojmovi zavisni i nezavisni vektori sre\'cu se indirektno u ranijim stadijima obrazovanja. Pa tako za dva vektora ili vi\v se vektora ka\v zemo da su {\bfseries\color{magenta}kolinearni} ako le\v ze na istim ili paralelnim nosa\v cima (radimo sa slobodnim vektorima). Recimo da dva vektora $\vec{a},\,\vec{b},$ le\v ze na istoj pravoj (nosa\v cu), koriste\'ci mno\v zenje vektora skalarom, mo\v zemo predstaviti (izraziti) jedan vektor preko drugog $\vec{a}=k\vec{b}$ ili $\vec{b}=\frac{1}{k}\vec{a},\:k\neq 0,$ Slika \ref{kolinearnostSlika}. Isto tako za tri ili vi\v se vektora ka\v zemo da su {\bfseries\color{magenta}komplanarni} ako le\v ze u jednoj ravni ili su sa njom paralelelni. Neka su sada data tri komplanarna vektora $\vec{a},\, \vec{b},\, \vec{c},$ kao na Slici \ref{komplanarnotsSlika}. Vidimo na koji na\v cin mo\v zemo vektor $\vec{c}$ izraziti preko vektora $\vec{a}$  i $\vec{b}.$  U ovom slu\v caju vektor $\vec{a}$ \'cemo produ\v ziti, pomno\v zi\'cemo ga sa $2,$ dok \'cemo intenzitet vektora $\vec{b}$ smanjiti (skrati\'cemo du\v zinu vektora $\vec{b}$). Na kraju, vektor $\vec{c}$ jednak je linearno kombinaciji $2\vec{a}+0.5\vec{b},$ tj. $\vec{c}=2\vec{a}+0.5\vec{b}.$ Vektori $\vec{a}$ i $\vec{b}$ (Slika \ref{kolinearnostSlika}) su zavisni i mo\v zemo ih izraziti jedan preko drugog. Isto tako, vektori dati na Slici \ref{komplanarnotsSlika}, $\vec{a},\, \vec{b},\, \vec{c}$ su zavisni jer mo\v zemo jedan vektor izraziti preko druga dva.

\index{vektori!kolinearni}\index{vektori!komplanarni}

Vidjeli smo primjer da na pravoj, koja predstavlja jednodimenzionalni prostor, mo\v zemo  izraziti vektor preko samo jednog vektora. Zatim, u ravni, koja predstavlja dvodimenzionalni prostor, vektor mo\v zemo izraziti preko druga dva vektora. Vratimo se sada na prostor $V^3,$ koji je naravno trodimenzionalni.  Vrijedi sljede\' ca teorema.

\begin{theorem}\label{vektori1}\index{teorema! o reprezentacije vektora iz prostora $V^3$}
  Neka su su $\vec{a}_1,\, \vec{a}_2,\,\vec{a}_3\in V^3$  bilo koja tri linearno nezavisna vektora. Tada se svaki
  vektor $\vec{a}\in V^3$ mo\v ze na jedinstven na\v cin prikazati kao njihova linearna kombinacija.
\end{theorem}
Vrijedi i sljede\'ca definicija.

\begin{definition}\label{vektori2}
  Ure\dj enu trojku $(\vec{a}_1,\,\vec{a}_2,\,\vec{a}_3)$ tri linearno nezavisna vektora iz $V^3$ nazivamo
      baza prostora $V^3.$
\end{definition}

Drugim rije\v cima bilo koji vektor $\vec{a}\in V^3$ mo\v zemo izraziti preko nezavisnih vektora $\vec{a}_1,\,\vec{a}_2,\,\vec{a}_3,$ tj. predstaviti  vektor $\vec{a}$ kao linearnu kombinaciju vektora $\vec{a}_1,\,\vec{a}_2,\,\vec{a}_3$
\begin{empheq}[box=\mymath]{equation*}
\vec{a}=\alpha_1\vec{a}_1+\alpha_2\vec{a}_2+\alpha_3\vec{a}_3.
\end{empheq}
Skalare $\alpha_1,\alpha_2,\alpha_3$ nazivamo koeficijenti ili koordinate vektora $\vec{a}$ u bazi $(\vec{a}_1,\,\vec{a}_2,\,\vec{a}_3).$
\index{vektori!koeficijenti ili koordinate vektora}


\subsection{Prostorni koordinatni sistem. Ortogonalna baza}
 Dakle, vektor $\vec{a}\in V^3$ mo\v zemo predstaviti preko bilo koje baze prostora $V^3.$ Na osnovu prethodno izlo\v zenog, svaka tri nezavisna vektora iz $V^3$ \v cine bazu tog prostora i mo\v zemo ih iskoristiti za predstavljanje (razlaganje) nekog vektora iz $V^3.$  Cilj nam je odabrati bazu tako da to bude \v sto jednostavnije, kao i da kalkulacije sa takvim vektorima budu isto tako \v sto  jednostavnije. Da bi odredili jednu takvu bazu, potrebno je definisati nekoliko novih pojmova.

\begin{definition}[Orjentisana prava ili osa]
 Orjentisanom pravom ili osom naziva se prava za \v cije je dvije proizvoljne ta\v cke utvr\dj eno, koja se od njih smatra prethodnom a koja sljede\' com. Orjentisana prava ili osa mo\v ze biti okarakterisana jedini\v cnim vektorom $\vec{u}.$
\end{definition}

\index{vektori!orjentisana prava ili osa}

\begin{figure}[!h]\centering
        \includegraphics[scale=1]{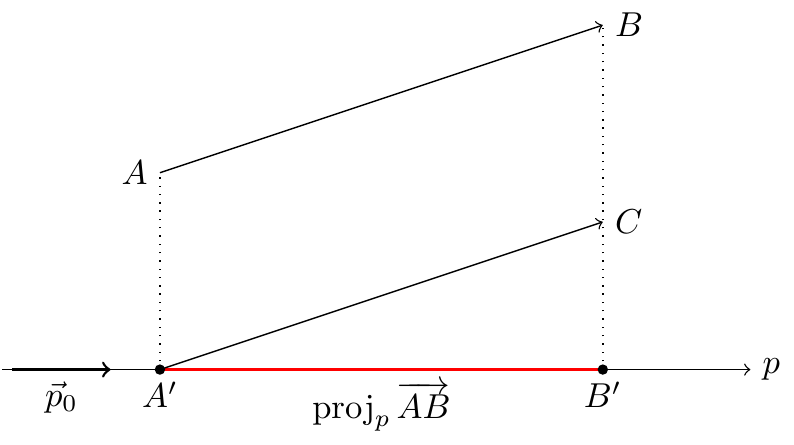}
        \caption{Projekcija vektora  na osu }
     \label{slika48}
  \end{figure}

\begin{definition}[Ortogonalna projekcija, Slika \ref{slika48}]
Pod ortogonalnom projekcijom vektora $\overrightarrow{AB}$ na osu $p$ podrazumijeva se du\v zina du\v zi $A'B',$ tj.  $ \left|A'B'\right|,$ pri \v cemu su ta\v cke $A'$ i $B'$ ortogonalne projekcije (respektivno) ta\v caka $A$ i $B$ na osu $p,$ uzeta sa predznakom $+$ ako je vektor $\overrightarrow{AB}$ orjentisan na istu stranu kao i osa $p$ u odnosu na ravan koja prolazi kroz ta\v cke $A$ i $A'$ i normalna je na pravu $p,$ a sa znakom $-$ u suprotnom slu\v caju.
\end{definition}
\index{vektori!ortogonalna projekcija}
Prava $p$ odnosno osa $p$ je orjentisana prava i to je uradjeno vektorom $\vec{p}_0,$ Slika \ref{slika48}.  Na istoj slici vidimo i ortogonalnu projekciju vektora $\overrightarrow{ AB}$ na osu $p.$ Ova projekcija je uradjena veoma jednostavno, odredi  se ortogonalna projekcija ta\v cke $A$ (spustimo normalu iz ta\v cke $A$ na osu $p$ i odredimo presje\v cnu ta\v cku normale i ose) i to je ta\v cka $A',$ zatim isto uradimo i sa ta\v ckom $B,$ dobijamo ta\v cku $B'$. Sada je du\v zina du\v zi $A'B',$ sa odgovaraju\'cim predznakom, ortogonalna projekcija vektora $\overrightarrow{ AB}$ na osu $p.$

\begin{figure}[!h]\centering
      \begin{subfigure}[b]{.45\textwidth}\centering
        \includegraphics[scale=.65]{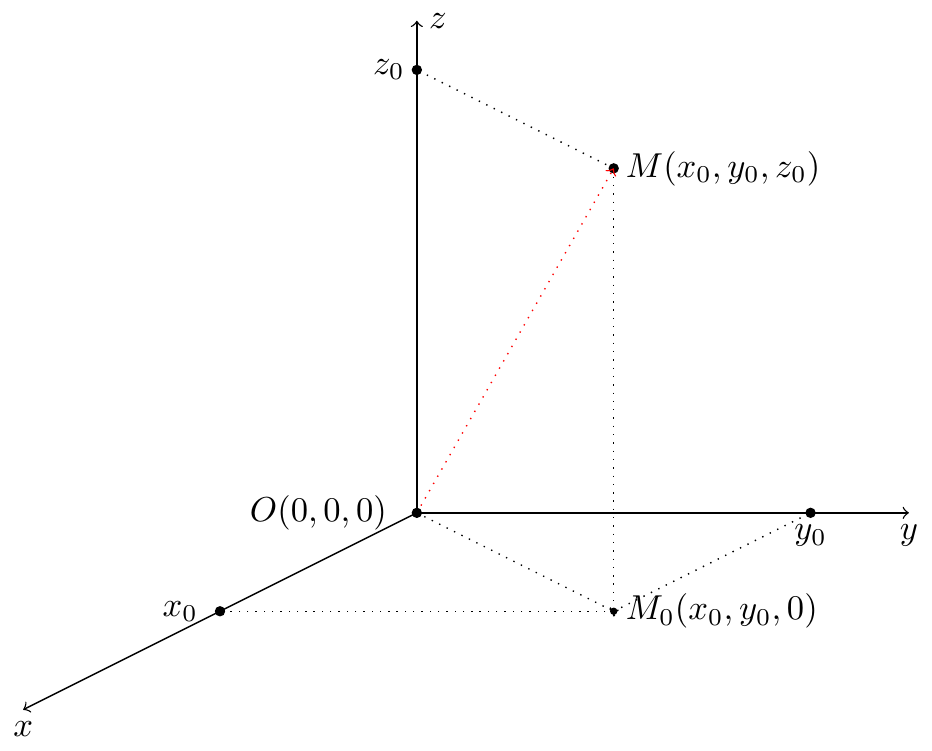}
         \caption{Pravougli koordinatni sistem}
          \label{slika46a}
   \end{subfigure}\hspace{.5cm}
     \begin{subfigure}[b]{.45\textwidth}\centering
        \includegraphics[scale=.65]{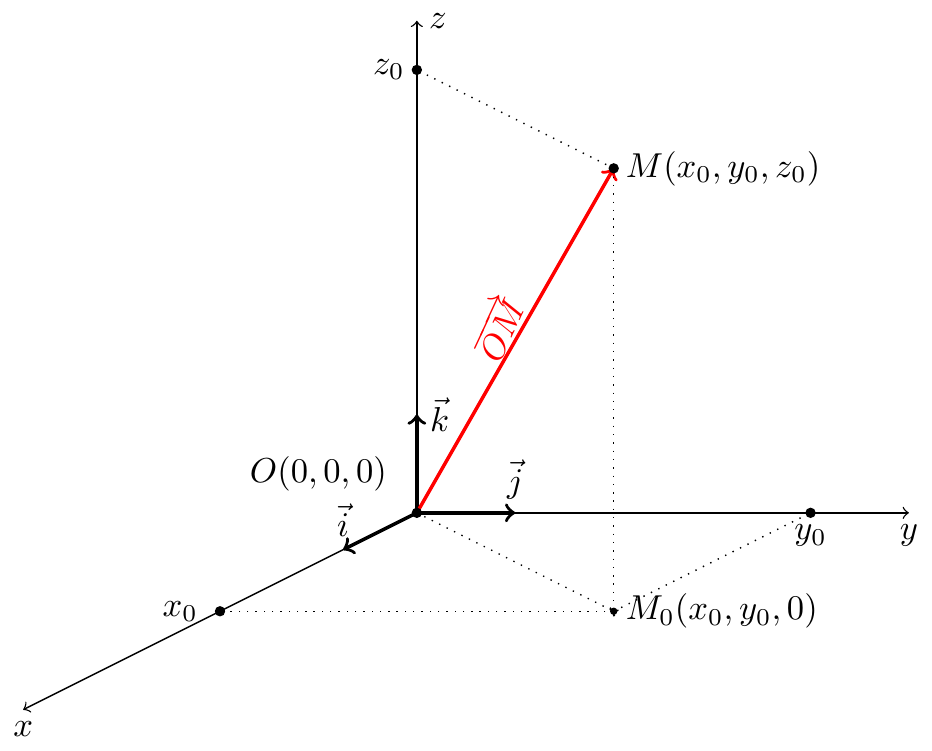}
         \caption{Vektor $\overrightarrow{OM}$}
          \label{slika46b}
   \end{subfigure}
        \caption{Pravougli koordinatni sistem u prostoru}
     \label{slika46}
  \end{figure}

\begin{definition}[Pravougli koordinatni sistem, Slika \ref{slika46a}]
 Ure\dj eni skup tri ose koje prolaze kroz utvr\dj enu ta\v cku $O$ (pol ili koordinatni po\v cetak) i koje su uzajamno okomite, obrazuju Descartesov$^1$ pravougli koordinatni sistem
\end{definition}
\footnotetext[1]{Ren\' e Descartes (31.mart 1596.--11.februar 1650.) bio je francuski filozof, matemati\v car i nau\v cnik.}
\index{vektori!pravougli koordinatni sistem}\index{Descartesov pravougli koordinatni sistem}

Dakle tri uzajamno ortogonalne orjentisane prave, koje prolaze kroz ta\v cku $O,$ obrazuju Descartesov pravougli koordinatni sistem. Prave koje obrazuju pravougli koordinatni sistem nazivaju se koordinatne ose, a njihova zajedni\v cka ta\v cka $O$ je koordinatni po\v cetak.

Ako se koordinatne ose obilje\v ze sa $Ox,\,Oy,\,Oz$ i orjenti\v su kao na Slici \ref{slika46a} ka\v ze se da one obrazuju desni prostorni koordinatni sistem i nazivaju se respektivno apscisna osa, ordinatna osa i aplikativna osa, odnosno $x$--osa, $y$--osa i $z$--osa. Jedini\v cne vektore (ort--vektori ili ortovi) koordinatnih osa $Ox,\,Oy,\,Oz,$ (Slika \ref{slika46b}) ozna\v cava\' cemo sa $\vec{i},\,\vec{j}$ i  $\vec{k}.$ Koordinatne ose $x,\,y$ i $z,$ odre\dj uju koordinatne ravni, pa je tako $xOy$ koordinatna ravan odre\dj ena koordinatnim osama $x$ i $y,$ analogno su odre\dj ene i koordinatne ravni $xOz$ i $yOz.$

\index{vektori!apscisa, ordinata, aplikata}
\index{vektori! ort vektori}

\begin{figure}[!h]\centering
        \includegraphics[scale=.85]{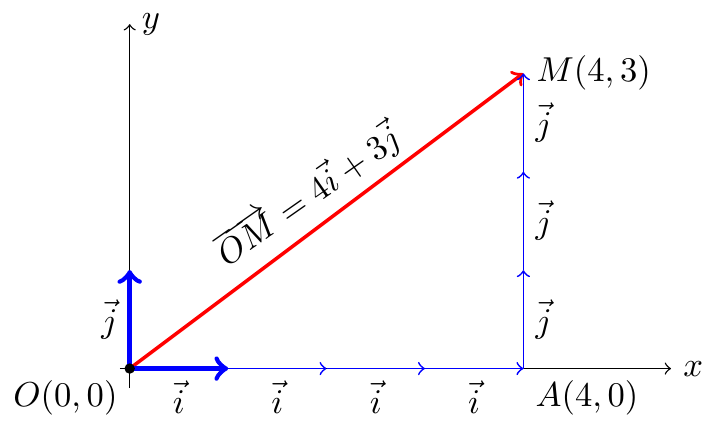}
        \caption{Predstavljanje vektora u ravni }
     \label{slika47}
  \end{figure}

Poka\v zimo kako se ta\v cka a zatim i vektor predstavljaju, prvo u ravni a zatim i u prostoru. Posmatrajmo Sliku \ref{slika47}, tj. slu\v caj u ravni ($z=0$). Ta\v cka $M$ ima koordinate 4 i 3, 4 je njena $x$--koordinata i to je udaljenost od $y$--ose, dok je 3 njena $y$--koordinata i to je udaljenost od $x$--ose. Na ovaj na\v cin jednozna\v cno je odre\dj ena ta\v cka $M(4,3)$ u ravni. Vektor $\overrightarrow{OM}$ mo\v zemo predstaviti u bazi $(\vec{i},\vec{j})$  pomo\' cu ovih projekcija ta\v cke. Sa Slike \ref{slika47} vidimo da je
\begin{empheq}[box=\mymath]{equation*}
\overrightarrow{OM}=4\vec{i}+3\vec{j}.
\end{empheq}
Vidimo da je vektor $\overrightarrow{OM}$ zbir dva vektora i to $\overrightarrow{OA}=4\vec{i}$ i $\overrightarrow{AM}=3\vec{j},$ odnosno linearna kombinacija vektora $\vec{i}$ i $\vec{j}.$

Na sli\v can na\v cin mo\v zemo predstaviti ta\v cku u prostoru, odnosno njen odgovaraju\'ci vektor. U ovom slu\v caju imamo jednu dodatnu koordinatu $z$.  Posmatrajmo Sliku \ref{slika46a} ili \ref{slika46b}. Koordinate ta\v cke $M$ su upravo njena odstojanja od koordinatnih ravni i to $x_0$ je odstojanje od $yOz$ ravni, $y_0$ od $xOz$ ravni i $z_0$ je odstojanje od  $xOy$ koordinatne ravni. Odredili smo polo\v zaj ta\v cke u pravouglom koordinatnom sistemu, Slika \ref{slika46a} ili \ref{slika46b}.

Odredimo sada vektor $\overrightarrow{OM}.$ Ta\v cka $M_0$ je projekcija ta\v cke $M$ na koordinatnu ravan $xOy$ (spustimo normalu kroz ta\v cku $M$ na $xOy$ i presje\v cna ta\v cka normale i koordinatne ravni je tra\v zena ta\v cka). Odredimo sada vektor $\overrightarrow{OM_0}.$ Njegova po\v cetna ta\v cka je  $O,$ a krajnja $M_0(x_0,y_0,0),$ pa je $\overrightarrow{OM_0}=x_0\vec{i}+y_0\vec{j}.$ Prava koja prolazi kroz ta\v cke $M$ i $M_0$ paralelna je sa $z$ osom, pa je vektor kojem je po\v cetna ta\v cka  $M_0$ a krajnja $M,$ $\overrightarrow{M_0M}=z_0\vec{k}.$ Na kraju, vektor koji spaja koordinatni po\v cetak i ta\v cku $M$ ra\v cunamo kao zbir vektora $\overrightarrow{OM_0}$ i $\overrightarrow{M_0M},$ tj. $\overrightarrow{OM}=\overrightarrow{OM_0}+\overrightarrow{M_0M}=x_0\vec{i}+y_0\vec{j}+z_0\vec{k}.$    Drugim rije\v cima,  da bi odredili vektor kojem je po\v cetna ta\v cka koordinatni po\v cetak $O,$ a krajnja ta\v cka $M,$ dovoljno je da znamo koordinate te krajnje ta\v cke $M,$ tj. vrijedi
\begin{empheq}[box=\mymath]{equation*}
\overrightarrow{OM}=x_0\vec{i}+y_0\vec{j}+z_0\vec{k}.
\end{empheq}
Dakle i ovdje smo vektor $\overrightarrow{OM}$ predstavili kao linearnu kombinaciju vektora, ali sada tri vektora $\vec{i},\, \vec{j}$ i $\vec{k}.$
\index{vektori!predstavljanje vektora}
Osim notacije $\overrightarrow{OM}=x_0\vec{i}+y_0\vec{j}+z_0\vec{k}$ koristimo i sljede\' cu notaciju
\begin{empheq}[box=\mymath]{equation*}
\overrightarrow{OM}=(x_0,y_0,z_0).
\end{empheq}
Ina\v ce vektor $\overrightarrow{OM}$ se naziva i vektor polo\v zaja ta\v cke $M.$

\begin{figure}[!h]\centering
        \includegraphics[scale=1]{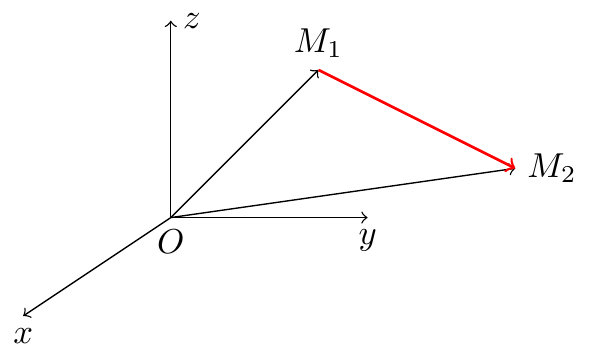}
        \caption{Predstavljenje vektora  $\vec{M_{1} M_{2}}$}
     \label{slika49}
  \end{figure}

Neka su date dvije ta\v cke  $M_1(x_1,y_1,z_1)$ i $M_2(x_2,y_2,z_2)$ kao na Slici \ref{slika49}. Kako izraziti vektor $\overrightarrow{M_1M_2}?$ Vidimo da je $\overrightarrow{OM_2}=\overrightarrow{OM_1}+\overrightarrow{M_1M_2}$.
\noindent Sada je
\begin{align*}
\overrightarrow{M_1M_2}&=\overrightarrow{OM_2}-\overrightarrow{OM_1}
            =(x_2\vec{i}+y_2\vec{j}+z_2\vec{k})-(x_1\vec{i}+y_1\vec{j}+z_1\vec{k})\\
            &=(x_2-x_1)\vec{i}+(y_2-y_1)\vec{j}+(z_2-z_1)\vec{k},
\end{align*}
tj. vrijedi
\begin{empheq}[box=\mymath]{equation}\label{vektorDvijeTacke}
\overrightarrow{M_1M_2}=(x_2-x_1)\vec{i}+(y_2-y_1)\vec{j}+(z_2-z_1)\vec{k}.
\end{empheq}
Dakle, ako su date koordinate po\v cetne i krajnje ta\v cke vektora, njegove koeficijente ra\v cu-\\namo tako \v sto od koordinata krajnje ta\v cke oduzimamo koordinate  po\v cetne ta\v cke.

\begin{example}
 Ispitati da li su vektori $\vec{a}=(4,-6,10)$ i $\vec{b}=(-6,9,-15)$ kolinearni.\\\\
\noindent Rje\v senje:\\\\
  Znamo da su dva vektora kolinearna ako le\v ze na istom nosa\v cu  (pravoj) ili se mogu paralelnim pomjeranjem dovesti na isti nosa\v c. Ako le\v ze na istom nosa\v cu onda je
  \begin{align*}
    \vec{a}&=k\vec{b}\Leftrightarrow (a_1,a_2,a_3)=k(b_1,b_2,b_3),
  \end{align*}
dva vektora su jednaka ako su im odgovaraju\' ci koeficijenti jednaki pa je
\begin{align*}
   a_1&=kb_1\Leftrightarrow\frac{a_1}{b_1}=k,\\
   a_2&=kb_2\Leftrightarrow\frac{a_2}{b_2}=k,\\
   a_3&=kb_3\Leftrightarrow\frac{a_3}{b_3}=k,
\end{align*}
pa dobijamo produ\v zenu   proporciju
\begin{equation*}
 \frac{a_1}{b_1}=\frac{a_2}{b_2}=\frac{a_3}{b_3}=k.
\label{vektor50}
\end{equation*}
Sada samo treba ispitati da li koeficijenti vektora zadovoljavaju prethodnu proporciju. Vrijedi
\[\frac{4}{-6}=\frac{-6}{9}=\frac{10}{-15}=-\frac{2}{3},\]
dakle vektori su kolinearni.
\end{example}

\begin{example}
 Odrediti parametre $\alpha$ i $\beta$ tako da vektori $\vec{a}=(\alpha,1,-4)$ i $\vec{b}=(3,4,\beta)$ budu kolinearni.\\\\
\noindent Rje\v senje:\\\\
Koristimo ponovo uslov \eqref{vektor50}	
\begin{align*}
  \frac{\alpha}{3}&=\frac{1}{4}=\frac{-4}{\beta},\\
  \frac{\alpha}{3}&=\frac{1}{4}\Leftrightarrow\alpha=\frac{3}{4},\\
  \frac{1}{4}&=\frac{-4}{\beta}\Leftrightarrow\beta=-16.
\end{align*}
\end{example}

\begin{example}
Razlo\v ziti vektor $\vec{a}$ u pravcu vektora $\vec{b}$ i $\vec{c},$ ako je $\vec{a}=3\vec{p}-2\vec{q},\vec{b}=-2\vec{p}+\vec{q}$ i $\vec{c}=7\vec{p}-4\vec{q}.$\\\\
\noindent Rje\v senje:\\\\
  Vrijedi
  \begin{align*}
    \vec{a}&=\alpha\vec{b}+\beta\vec{c}\\
   3\vec{p}-2\vec{q}&=\alpha(-2\vec{p}+\vec{q})+\beta(7\vec{p}-4\vec{q})\\
   3\vec{p}-2\vec{q}&=(-2\alpha+7\beta)\vec{p}+(\alpha-4\beta)\vec{q},
  \end{align*}
izjedna\v cavaju\' ci koeficijente dobijamo sistem
\begin{align*}
  3&=-2\alpha+7\beta\\
  -2&=\alpha-4\beta,
\end{align*}
\v cijim rje\v savanjem dobijamo $\alpha=2,\,\beta=1,$ pa je
\[\vec{a}=2\vec{b}+\vec{c}.\]
\end{example}

\begin{example}
Ispitati linearnu zavisnost vektora $\vec{a}=(1,2,0),\,\vec{b}=(2,-4,1),\,\vec{c}=(1,-1,-1).$\\\\
\noindent Rje\v senje:\\\\
\noindent \textbf{I na\v cin}--rang matrice\\
Formirajmo matricu \v ciji su elementi jedne vrste upravo koeficijenti jednog od ve-\\ktora, druge vrste koeficijenti drugog vektora i tre\'{c}e vrste
koeficijenti tre\'{c}eg vektora. Zatim odredimo rang te matrice. Vrijedi
\begin{align*}
 A=\begin{blockarray}{rrr}\begin{block}{(rrr)}1&2&0\\2&-4&1\\1&-1&-1\\\end{block}\end{blockarray}\sim
 \begin{blockarray}{rrr}\begin{block}{(rrr)}1&2&0\\0&-8&1\\0&-3&-1\\ \end{block}\end{blockarray}\sim
 \begin{blockarray}{rrr}\begin{block}{(rrr)}1&2&0\\0&-8&1\\0&0&11\\ \end{block}\end{blockarray},
\end{align*}
dakle $\rang A=3,$ tj. imamo tri nezavisne vrste ili kolone, a po\v sto su u matrici elementi vrsta koeficijenti od vektora to su na\v si vektori nezavisni.

\noindent\textbf{II na\v cin}--sistem homogenih jedna\v cina\\
Vektori $\vec{a},\vec{b},\vec{c}$ su linearno nezavisni ako jedna\v cina $\alpha\vec{a}+\beta\vec{b}+\gamma\vec{c}=\vec{0}$ ima samo trivijalno rje\v senje, tj. $\alpha=\beta=\gamma=0 .$ U suprotnom, ako ima i drugih rje\v senja, vektori su zavisni. Vrijedi
\begin{align*}
&    \alpha\vec{a}+\beta\vec{b}+\gamma\vec{c}=\vec{0}\Leftrightarrow\\
&   \alpha(\vec{i}+2\vec{j})+\beta(2\vec{i}-4\vec{j}+\vec{k})+\gamma(\vec{i}-\vec{j}-\vec{k})=\vec{0}  \Leftrightarrow\\
&(\alpha+2\beta+\gamma)\vec{i}+(2\alpha-4\beta-\gamma)\vec{j}+(\beta-\gamma)\vec{k}=\vec{0},\text{ ($\vec{0}$ predstavlja nula--vektor)}.
\end{align*}
Iz posljednje jednakosti dobijamo homogeni sistem
\[
  \begin{cases}
    \alpha+2\beta+\gamma=0\\
    2\alpha-4\beta-\gamma=0\\
    \beta-\gamma=0,
   \end{cases}
\]
a vrijednost determinante je
\[\begin{blockarray}{rrr}\begin{block}{|rrr|}1&2&1\\2&-4&-1\\0&1&-1\\\end{block}\end{blockarray}=11\neq 0. \]
Dakle sistem ima samo trivijalno rje\v senje $\alpha=\beta=\gamma=0.$ Zaklju\v cujemo da su posmatrani vektori linearno nezavisni.
\end{example}

Predstavljaju\'ci vektore analiti\v cki, tj. nekom formulom, znatno je pro\v sirena pri-\\mjena vektora, a samim tim i kalkulacije sa vektorima su u\v cinjene znatno jednostavnijim. Na po\v cetku ovog poglavlja vidjeli smo kako se ra\v cuna proizvod skalara i vektora. Sada mo\v zemo, imaju\'ci u vidu kako se predstavljaju vektori analiti\v cki, pomno\v ziti me\dj usobno dva ili vi\v se vektora. Prvo \'cemo vidjeti kako se ra\v cuna proizvod dva vektora. Rezultat ovih operacija mo\v ze biti skalar ili vektor.

\section{S\lowercase{kalarni proizvod}}

Neka su data dva vektora $\vec{x}$ i $\vec{y}.$  Kada je rezultat operacije mno\v zenja ova dva vektora skalar, onda se radi o skalarnom proizvodu. Slijedi definicija.

\begin{definition}[Skalarni proizvod]\label{definicijaSkalarni}
  Broj, odnosno skalar \[|\vec{x}||\vec{y}|\cos\measuredangle(\vec{x},\vec{y})\] naziva se skalarni proizvod vektora $\vec{x}$ i $\vec{y}$ i obilje\v zava se sa \[\vec{x}\cdot\vec{y}.\]
\end{definition}\index{vektori! skalarni proizvod}
Dakle, vrijedi
\[\vec{x}\cdot\vec{y}=|\vec{x}||\vec{y}|\cos\measuredangle(\vec{x},\vec{y}),\]
gdje je $\measuredangle(\vec{x},\vec{y})$ ugao izme\dj u vektora $\vec{x}$ i $\vec{y}.$

Skalarni proizvod je binarna operacija, koja nije zatvorena, po\v sto rezultat proizvoda $\vec{x}\cdot\vec{y}$  nije vektor nego skalar.

Primjetimo da je $|\vec{y}|\cos\measuredangle(\vec{x},\vec{y})$ ortogonalna projekcija vektora $\vec{y}$ na osu vektora $\vec{x}$ (ili preciznije na pravu ili nosa\v c vektora $\vec{x}$). Sada je
\[\vec{x}\cdot\vec{y}=|\vec{x}|\proj_{\vec{x}}\vec{y},\] vrijedi i
\[\vec{x}\cdot\vec{y}=|\vec{y}|\proj_{\vec{y}}\vec{x}.\]

 \begin{figure}[!h]\centering
        \includegraphics[scale=.9]{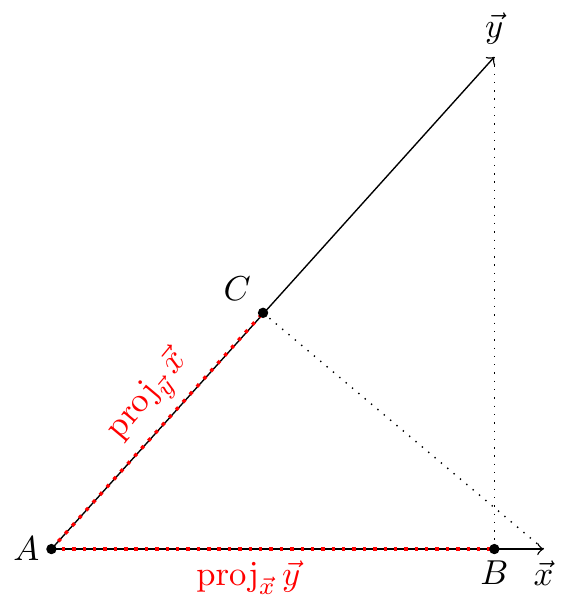}
         \caption{Projekcije vektora}
          \label{slika50}
  \end{figure}

Osobine skalarnog proizvoda date su u sljede\' coj teoremi.

\index{vektori! osobine skalarnog proizvoda}

\begin{theorem}[Osobine skalarnog proizvoda]\index{teorema! o osobinama skalarnog proizvoda}\label{teoremaSkalarni1}
 Neka su $\vec{x},\,\vec{y},\,\vec{z}$ proizvoljni vektori, a $\alpha$ proizvoljan skalar. Vrijedi
 \begin{enumerate}[$(1)$]
  \item $\vec{x}\cdot\vec{y}=\vec{y}\cdot\vec{x};$
  \item $\vec{x}\cdot(\vec{y}+\vec{z})=\vec{x}\cdot\vec{y}+\vec{x}\cdot\vec{z};$
  \item $(\alpha\vec{x})\cdot\vec{y}=\alpha(\vec{x}\cdot\vec{y});$
  \item $\vec{x}\cdot\vec{x}=\left|\vec{x}\right|^2;$\label{skalarni1}
  \item $\left|\vec{x}\cdot\vec{y}\right|\leqslant\left|\vec{x}\right|\cdot\left|\vec{y}\right|$ (Cauchy$^{1}$--Schwarzova$^{2}$ nejednakost).
 \end{enumerate}
\end{theorem}
\footnotetext[1]{Augustin-Louis Cauchy (21.avgust 1789.--23.maj 1857. godine) bio je francuski matemati\v car, in\v zinjer i fizi\v car koji je napravio pionirske doprinose u mnogim oblastima matematike. Bio je jedan od prvih matemati\v cara koji koristi  rigorozne dokaze u matemati\v ckoj analizi.}
\footnotetext[2]{Karl Hermann Amandus Schwarz  (25.januar 1843. -- 30.novembar 1921. godine) bio je njema\v cki matemati\v car najvi\v se poznat po doprinosu u kompleksnoj analizi. }
\begin{proof}[Dokaz.]
Doka\v zimo samo osobinu \eqref{skalarni1}. Na osnovu definicije skalarnog proizvoda je
\[\vec{x}\cdot\vec{x}=\left|\vec{x}\right|\cdot\left|\vec{x}\right|\cos\measuredangle(\vec{x},\vec{x})
           =\left|\vec{x}\right|\cdot\left|\vec{x}\right|\cos 0=\left|\vec{x}\right|\cdot\left|\vec{x}\right|=\left|\vec{x}\right|^2.\]
\end{proof}

Imaju\' ci u vidu da je $\cos 0=1,\,\cos 90^0=\cos\frac{\pi}{2}=0,$ izra\v cunajmo skalarne proizvode sa ortovima $\vec{i},\,\vec{j},\,\vec{k}.$ Vrijedi
\begin{align*}
  \vec{i}\cdot\vec{i}&=|\vec{i}|\cdot|\vec{i}|\cos 0=1\cdot 1\cdot 1=1,\\
  \vec{j}\cdot\vec{j}&=|\vec{j}|\cdot|\vec{j}|\cos 0=1,\\
  \vec{k}\cdot\vec{k}&=|\vec{k}|\cdot|\vec{k}|\cos 0=1,\\
  \vec{i}\cdot\vec{j}&=\vec{j}\cdot\vec{i}=|\vec{i}|\cdot|\vec{j}|\cos \frac{\pi}{2}=1\cdot 1\cdot 0=0,\\
  \vec{i}\cdot\vec{k}&=\vec{k}\cdot\vec{i}=|\vec{i}|\cdot|\vec{k}|\cos \frac{\pi}{2}=0,\\
  \vec{j}\cdot\vec{k}&=\vec{k}\cdot\vec{j}=|\vec{j}|\cdot|\vec{k}|\cos \frac{\pi}{2}=0.
\end{align*}
Dakle 
\begin{equation}
 \tcbhighmath[mojstil1]{\vec{i}\cdot\vec{i}=\vec{j}\cdot\vec{j}=\vec{k}\cdot\vec{k}=1}
 \label{skalarni2}
\end{equation}
i
\begin{equation}
 \tcbhighmath[mojstil1]{\vec{i}\cdot\vec{j}=\vec{j}\cdot\vec{i}=\vec{i}\cdot\vec{k}=\vec{k}\cdot\vec{i}=\vec{j}\cdot\vec{k}=\vec{k}\cdot\vec{j}=0.}
 \label{skalarni3}
\end{equation}

U sljede\' coj teoremi date su neke geometrijske osobine vektora koje mo\v zemo iskazati preko skalarnog proizvoda.

\begin{theorem}\index{teorema! neke geom. osobine vektora izra\v zene preko skalarnog proizvoda}
  Neka su $\vec{x},\,\vec{y}$  proizvoljni vektori. Vrijedi
  \begin{enumerate}[$(1)$]
     \item $|\vec{x}|=(\vec{x}\cdot\vec{x})^{\frac{1}{2}}$ (du\v zina, intenzitet, modul ili norma vektora);
     \item $\cos\measuredangle(\vec{x},\vec{y})=\dfrac{\vec{x}\cdot\vec{y}}{|\vec{x}||\vec{y}|}$ (ugao izme\dj u vektora);
     \item $\proj_{\vec{y}}\vec{x}=\dfrac{\vec{x}\cdot\vec{y}}{|\vec{y}|},\:\vec{y}\neq \vec{0}$ (projekcija vektora $\vec{x}$ na vektor $\vec{y}$);
     \item $\vec{x}\perp\vec{y}\Leftrightarrow(\vec{x}\cdot\vec{y}=0,\:\vec{x}\neq \vec{0}\wedge\vec{y}\neq \vec{0})$ (ortogonalnost vektora).
  \end{enumerate}
\end{theorem}

Poka\v zimo kako se ra\v cunaju vrijednosti iskazane u prethodnoj teoremi u ortonormiranoj bazi, tj. bazi koju \v cine vektori $\vec{i},\,\vec{j},\,\vec{k}.$ 
Neka su $\vec{x}=x_1\vec{i}+x_2\vec{j}+x_3\vec{k}$ i $\vec{y}=y_1\vec{i}+y_2\vec{j}+y_3\vec{k}.$ Tada je
\[\vec{x}\cdot\vec{y}=(x_1\vec{i}+x_2\vec{j}+x_3\vec{k})\cdot(y_1\vec{i}+y_2\vec{j}+y_3\vec{k}),\]
a kako vrijedi \eqref{skalarni2} i \eqref{skalarni3} to je
\begin{empheq}[box=\mymath]{equation*}
\vec{x}\cdot\vec{y}=x_1y_1+x_2y_2+x_3y_3.
\end{empheq}
Posljednja formula nam ka\v ze kako se ra\v cuna skalarni proizvod dva vektora koji su zadani u ortonormiranoj bazi $(\vec{i},\, \vec{j},\, \vec{k})$.
Koriste\'ci prethodnu formulu, vrijedi
\begin{empheq}[box=\mymath]{equation*}
\vec{x}\cdot\vec{x}=x^2_1+x^2_2+x^2_3,
\end{empheq}
a kako je $\vec{x}\cdot\vec{x}=|\vec{x}|^2,$ to dobijamo
\begin{empheq}[box=\mymath]{equation*}
|\vec{x}|=\sqrt{x^2_1+x^2_2+x^2_3}.
\end{empheq}
\newpage
Ovo je formula za ra\v cunanje du\v zine (ili norme ili intenziteta) vektora koji je zadan u ortonormiranoj bazi.
Dalje vrijedi
\[\tcbhighmath[mojstil1]{\cos\measuredangle(\vec{x},\vec{y})=\frac{x_1y_1+x_2y_2+x_3y_3}{\sqrt{x^2_1+x^2_2+x^2_3}\sqrt{y^2_1+y^2_2+y^2_3}},}\]

\[\tcbhighmath[mojstil1]{\proj_{\vec{y}}\vec{x}=\frac{x_1y_1+x_2y_2+x_3y_3}{\sqrt{y^2_1+y^2_2+y^2_3}},}\]
te
\[\tcbhighmath[mojstil1]{\vec{x}\perp\vec{y}\Leftrightarrow (x_1y_1+x_2y_2+x_3y_3=0,\:\vec{x}\neq \vec{0}\wedge\vec{y}\neq \vec{0}).}\]
Prva od posljednje tri formule je formula za ra\v cunanje ugla izmedju dva vektora $\vec{x}$ i $\vec{y},$ sljede\'ca formula je za ra\v cunanje projekcije vektora $\vec{x}$ na osu odredjenu vektorom $\vec{y},$ a tre\'ca formula slu\v zi za ispitivanje ortogonalnosti dva vektora.\\

\index{vektori! ra\v cunanje skalarnog proizvoda}
\index{vektori! ra\v cunanje intenziteta vektora}
\index{vektori! ra\v cunanje ugla izme\dj u vektora}
\index{vektori! projekcija vektora}
\index{vektori! uslov ortogonalnosti}

\begin{example}
Dati su vektori $\vec{a}=(-1,2,1)$ i $\vec{b}=(1,-3,2).$\\ Izra\v cunati
\begin{inparaenum}[$(a)$]
  \item $|\vec{a}|;\,$
  \item $|\vec{b}|;\,$
  \item $\vec{a}\cdot\vec{b};\,$
  \item $\cos\measuredangle(\vec{a},\vec{b});\,$
  \item $\proj_{\vec{a}}\vec{b}.$
\end{inparaenum}\ \\
\noindent Rje\v senje:\\\\
Vrijedi

   \begin{enumerate}[$(a)$]
    \item $|\vec{a}|=\sqrt{a^2_1+a^2_2+a^2_3}=\sqrt{1+4+1}=\sqrt{6};$
    \item $|\vec{b}|=\sqrt{b^2_1+b^2_2+b^2_3}=\sqrt{1+9+4}=\sqrt{14};$
    \item $\vec{a}\cdot\vec{b}=a_1b_1+a_2b_2+a_3b_3=(-1,2,1)\cdot(1,-3,2)=-1-6+2=-5;$
    \item $\cos\measuredangle(\vec{a},\vec{b})=\dfrac{a_1b_1+a_2b_2+a_3b_3}{\sqrt{a^2_1+a^2_2+a^2_3}\sqrt{b^2_1+b^2_2+b^2_3}}
                =\dfrac{(-1,2,1)\cdot(1,-3,2)}{\sqrt{1+4+1}\sqrt{1+9+4}}=\dfrac{-1-6+2}{\sqrt{6}\sqrt{14}}=\dfrac{-5}{\sqrt{84}};$
    \item $\proj_{\vec{a}}\vec{b}=\dfrac{\vec{a}\cdot\vec{b}}{|\vec{a}|}=\dfrac{-5}{\sqrt{6}}.$
   \end{enumerate}
\end{example}
\ \\

\begin{figure}[!h]\centering
   \includegraphics[scale=1.]{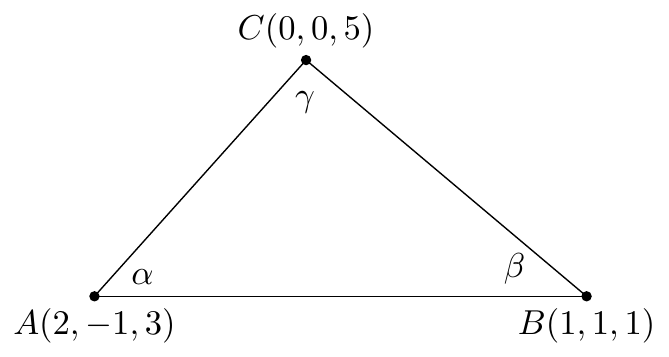}
         \caption{Trougao $\triangle ABC$}
          \label{slika51}
  \end{figure}

\begin{example}
Izra\v cunati unutra\v snje uglove i obim trougla trougla $\triangle ABC$ \v cija su tjemena\\ $A(2,-1,3),\,B(1,1,1),\,C(0,0,5).$\\\\
\noindent Rje\v senje: \\\\
Vidjeti Sliku \ref{slika51}, vrijedi
\begin{align*}
 \overrightarrow{ AB}\cdot\overrightarrow{ AC}&= \left|\overrightarrow{ AB}\right|\cdot\left|\overrightarrow{ AC}\right|\cos\alpha
     \Leftrightarrow\cos\alpha=\frac{ \overrightarrow{ AB}\cdot\overrightarrow{ AC}}{\left|\overrightarrow{ AB}\right|\cdot\left|\overrightarrow{ AC}\right|}\\
 \overrightarrow{ BA}\cdot\overrightarrow{ BC}&= \left|\overrightarrow{ BA}\right|\cdot\left|\overrightarrow{ BC}\right|\cos\beta
     \Leftrightarrow\cos\beta=\frac{ \overrightarrow{ BA}\cdot\overrightarrow{ BC}}{\left|\overrightarrow{ BA}\right|\cdot\left|\overrightarrow{ BC}\right|}.
\end{align*}
Uglove $\alpha$  i $\beta$ ra\v cunamo koriste\' ci prethodne formule, dok tre\' ci ugao mo\v zemo dobiti iz $\alpha+\beta+\gamma=\pi.$ Dalje je

\begin{align*}
  \overrightarrow{AB}&=(-1,2,-2),\,\overrightarrow{BA}=-\overrightarrow{AB}=(1,-2,2)\\
  \overrightarrow{AC}&=(-2,1,2)\\
  \overrightarrow{BC}&=(-1,-1,4)
\end{align*}
pa je

\begin{align*}
   \cos\alpha&=\frac{(-1,2,-2)\cdot(-2,1,2)}{\sqrt{1+4+4}\sqrt{4+1+4}}=0\Leftrightarrow\alpha=\frac{\pi}{2}\\
   \cos\beta&=\frac{(1,-2,2)\cdot(-1,-1,4)}{\sqrt{1+4+4}\sqrt{1+1+16}}=\frac{\sqrt{2}}{2}\Leftrightarrow\beta=\frac{\pi}{4}
\end{align*}
i \[\gamma=\pi-\alpha-\beta=\frac{\pi}{4}.\]
Obim trougla je
\begin{align*}
O&=\left|\overrightarrow{AB}\right|+\left|\overrightarrow{BC}\right|+\left|\overrightarrow{CA}\right|=
  \sqrt{1+4+4}+\sqrt{1+1+16}+\sqrt{4+1+4}\\
 &=6+\sqrt{18}=6+3\sqrt{2}.
\end{align*}
\end{example}

\newpage

\begin{figure}[!h]\centering
        \includegraphics[scale=1.]{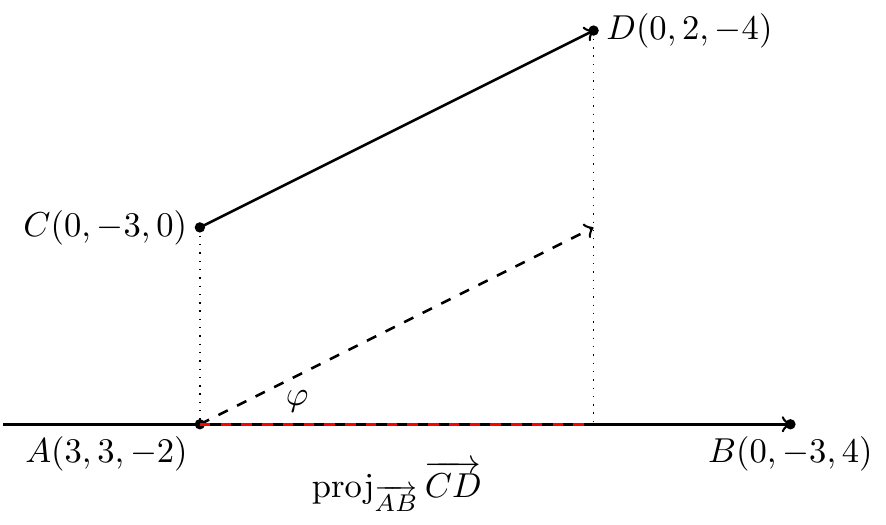}
         \caption{Projekcija vektora}
          \label{slika56}
  \end{figure}

\begin{example}
 Date su ta\v cke $A(3,3,-2),\,B(0,-3,4),\,C(0,-3,0),\,D(0,2,-4).$ Izra\v cunati $\proj_{\overrightarrow{AB}}\overrightarrow{CD}.$\\\\
\noindent Rje\v senje:\\

\noindent Sa Slike \ref{slika56} je
$\cos\varphi=\frac{\proj_{\overrightarrow{AB}}\overrightarrow{CD}}{|\overrightarrow{CD}|}\text{ pa je }\proj_{\overrightarrow{AB}}\overrightarrow{CD}=|\overrightarrow{CD}|\cos\varphi,$
znamo da je $\cos\varphi=\frac{\overrightarrow{CD}\cdot\overrightarrow{AB}}{|\overrightarrow{CD}||\overrightarrow{AB}|},$
te vrijedi
\[\proj_{\overrightarrow{AB}}\overrightarrow{CD}=|\overrightarrow{CD}|\frac{\overrightarrow{CD}\cdot\overrightarrow{AB}}{|\overrightarrow{CD}||\overrightarrow{AB}|}
     =\frac{\overrightarrow{CD}\cdot\overrightarrow{AB}}{|\overrightarrow{AB}|}.\]
     Odredimo sada $\overrightarrow{AB},\,\overrightarrow{CD},\,|\overrightarrow{AB}|$
     \begin{align*}
       \overrightarrow{CD}&=(0-0)\vec{i}+(2-(-3))\vec{j}+(-4-0)\vec{k}=(0,5,-4),\\
       \overrightarrow{AB}&=(0-3)\vec{i}+(-3-3)\vec{j}+(4-(-2))\vec{k}=(-3,-6,6),\\
       |\overrightarrow{AB}|&=\sqrt{9+36+36}=9,
     \end{align*}
pa je na kraju
\[\proj_{\overrightarrow{AB}}\overrightarrow{CD}=\frac{(0,5,-4)\cdot(-3,-6,6)}{9}=\frac{0-30-24}{9}=-6.\]
\end{example}
\begin{remark}
Vidimo da na Slici \ref{slika56} vektori nisu dobro orjentisani, zbog negativne projekcije.
\end{remark}

\section{V\lowercase{ektorski proizvod}}

Kao \v sto je prethodno pomenuto, rezultat operacije mno\v zenja dva vektora mo\v ze biti i vektor, pa stoga imamo potrebu uvesti i vektorski proizvod dva vektora.

Iz Definicije skalarnog proizvoda \ref{definicijaSkalarni} lako je vidjeti da je operacija skalarnog proizvoda, a to je navedeno i u Teoremi o osobinama skalarnog proizvoda \ref{teoremaSkalarni1}, komutativna operacija, tj. da vrijedi $\vec{x}\cdot\vec{y}= \vec{y}\cdot \vec{x}.$  Vidje\'cemo da vektorski proizvod nema ovu osobinu. Zbog toga \'{c}e nam trebati preciznija informacija o rasporedu vektora u prostoru, pa najprije defini\v semo pojam desnog, odnosno lijevog, triedra tri vektora.

\begin{figure}[!h]\centering
      \begin{subfigure}[b]{.45\textwidth}\centering
        \includegraphics[scale=.85]{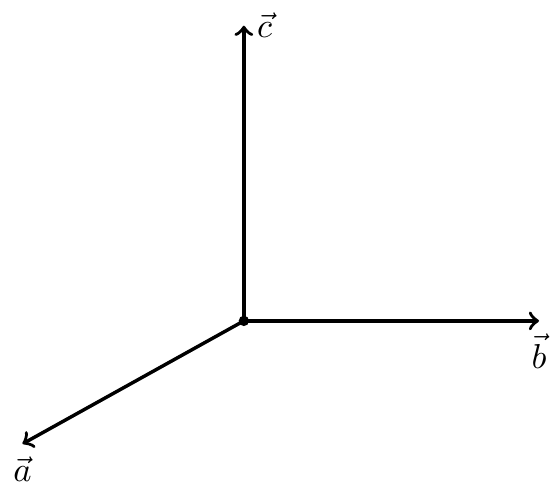}
         \caption{Desni triedar vektora $\vec{a},\,\vec{b},\, \vec{c}$}
          \label{slika52}
   \end{subfigure}\hspace{.5cm}
     \begin{subfigure}[b]{.45\textwidth}\centering
        \includegraphics[scale=.85]{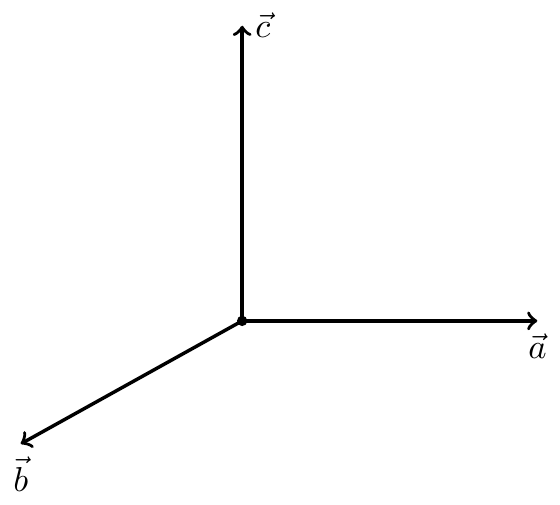}
         \caption{Lijevi triedar vektora $\vec{a},\,\vec{b},\,\vec{c}$}
          \label{slika53}
   \end{subfigure}
        \caption{Desni i lijevi triedar vektora}
     \label{slika46}
  \end{figure}

\begin{definition}[Desni triedar]
  Ka\v ze se da tri nekomplanarna vektora $\vec{a},\,\vec{b},\,\vec{c}$ sa zajedni\v ckim po\v cetkom obrazuju redom desni triedar ako se rotacija vektora $\vec{a}$ prema vektoru $\vec{b},$ najkra\' cim putem, posmatrano sa kraja vektora $\vec{c},$ vr\v si suprotno kretanju kazaljke na satu (Slika \ref{slika52}).
\end{definition}
Na sli\v can na\v cin defini\v se se lijevi triedar koji obrazuju tri nekomplanarna vektora $\vec{a},\,\vec{b},\,\vec{c}$ sa zajedni\v ckim po\v cetkom (Slika \ref{slika53}).
\index{vektori!desni triedar}
 \\
Razlikuju\'ci desni i lijevi triedar, mo\v zemo sada navesti definiciju vektorskog proizvoda.

\begin{definition}[Vektorski proizvod]
   Ako je $\vec{n}_0$ jedini\v cni vektor normalan na ravan koju obrazuju vektori $\vec{x}$ i $\vec{y},$ pri \v cemu $\vec{x},\,\vec{y}$ i $\vec{n}_0$ obrazuju desni triedar, onda se vektor \[|\vec{x}||\vec{y}|\sin\measuredangle(\vec{x},\vec{y})\vec{n}_0\]
   naziva vektorski proizvod vektora $\vec{x}$ i $\vec{y}$ i obilje\v zava sa $\vec{x}\times\vec{y}.$
\end{definition}
\index{vektori! vektorski proizvod}
Dakle, vrijedi
\begin{empheq}[box=\mymath]{equation*}
\vec{x}\times\vec{y}=|\vec{x}||\vec{y}|\sin\measuredangle(\vec{x},\vec{y})\vec{n}_0
\end{empheq}
\newpage
i
\begin{empheq}[box=\mymath]{equation*}
|\vec{x}\times\vec{y}|=|\vec{x}||\vec{y}||\sin\measuredangle(\vec{x},\vec{y})|.
\end{empheq}
\index{vektori! ra\v cunanje vektorskog proizvoda}
Vektorski proizvod dva vektora $\vec{x},\,\vec{y}$ je vektor $\vec{z}$, koji je normalan na ravan koju obrazuju vektori $\vec{x}$ i $\vec{y},$ intenzitet mu je jednak $|\vec{x}||\vec{y}||\sin\measuredangle(\vec{x},\vec{y})|,$ a smjer mu je takav da $\vec{x},\,\vec{y}$ i $\vec{z}$ obrazuju triedar desne orjentacije.

\begin{figure}[!h]\centering
        \includegraphics[scale=.85]{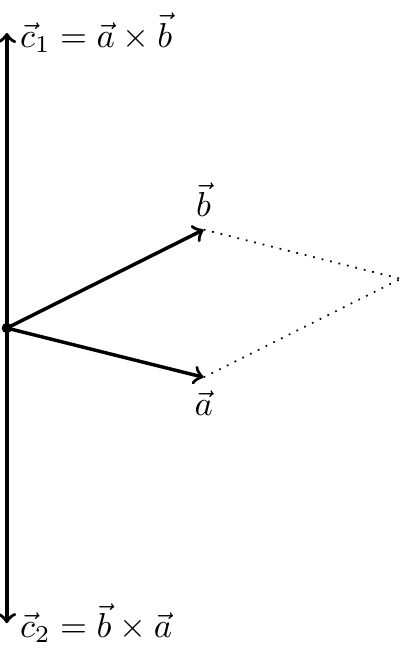}
         \caption{Vektori $\vec{c}_1$ i $\vec{c}_2$}
          \label{slika54}
  \end{figure}

Ako je $\vec{c}_1=\vec{a}\times\vec{b}$ i $\vec{c}_2=\vec{b}\times\vec{a},$ me\dj usobni raspored vektora prikazan je na Slici \ref{slika54}.
Iz same definicije i Slike \ref{slika54} vidimo da je $\vec{c}_1=-\vec{c}_2,$ drugim rije\v cima vektori $\vec{c}_1$ i $\vec{c}_2$ su suprotni, a operacija vektorski proizvod je antikomutativna operacija.

\begin{figure}[!h]\centering
        \includegraphics[scale=.85]{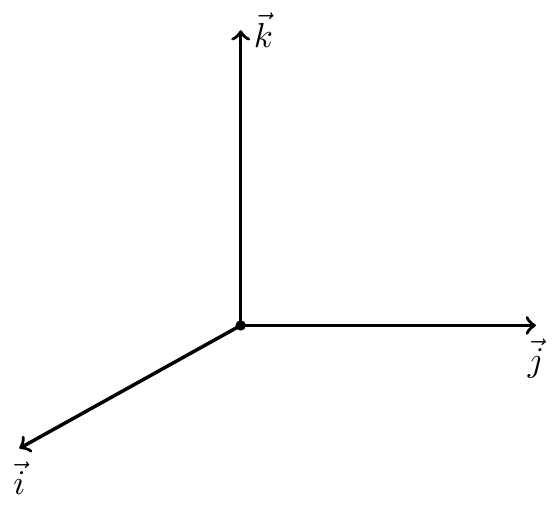}
         \caption{Ort--vektori $\vec{i},\,\vec{j},\,\vec{k}$}
          \label{slika55}
  \end{figure}

Kako ortovi $\vec{i},\,\vec{j},\,\vec{k}$ obrazuju triedar desne orjentacije, kao na Slici \ref{slika55},
to na osnovu definicije vektorskog proizvoda dobijamo
\begin{align*}
 \vec{i}\times\vec{i}&=\vec{0}&& \vec{i}\times\vec{j}=\vec{k}&&\vec{i}\times\vec{k}=-\vec{j}\\
 \vec{j}\times\vec{i}&=-\vec{k}&& \vec{j}\times\vec{j}=\vec{0}&&\vec{j}\times\vec{k}=\vec{i}\\
 \vec{k}\times\vec{i}&=\vec{j}&& \vec{k}\times\vec{j}=-\vec{i}&&\vec{k}\times\vec{k}=\vec{0}.
\end{align*}
\newpage
Vrijedi sljede\' ca teorema.

\index{vektori! osobine vektorskog proizvoda}

\begin{theorem}[Osobine vektorskog proizvoda]\index{teorema! o osobinama vektorskog proizvoda}
  Ako su $\vec{x},\,\vec{y},\,\vec{z}$ proizvoljni vektori i $\alpha$ je proizvoljan skalar, onda je
   \begin{enumerate}[$(1)$]
      \item $\vec{x}\times\vec{y}=-(\vec{y}\times\vec{x});$
      \item $(\alpha\vec{x})\times\vec{y}=\vec{x}\times (\alpha\vec{y})=\alpha(\vec{x}\times\vec{y});$
      \item $\vec{x}\times(\vec{y}+\vec{z})=\vec{x}\times\vec{y}+\vec{x}\times\vec{z}.$
   \end{enumerate}
\end{theorem}

Ako su vektori $\vec{x},\,\vec{y}$ dati u ortonormiranoj bazi, tj. $\vec{x}=(x_1,x_2,x_3)$ i $\vec{y}=(y_1,y_2,y_3)$ tada je
\begin{align*}
  \vec{x}\times\vec{y}&=(x_1\vec{i}+x_2\vec{j}+x_3\vec{k})\times(y_1\vec{i}+y_2\vec{j}+y_3\vec{k})\\
   &=(x_2y_3-x_3y_2)\vec{i}+(x_3y_1-x_1y_3)\vec{j}+(x_1y_2-y_1x_2)\vec{k}.
\end{align*}
Posljednju jednakost mo\v zemo pisati u obliku determinante

\index{vektori! ra\v cunanje vektorskog proizvoda}

\begin{empheq}[box=\mymath]{equation*}
\vec{x}\times\vec{y}=\begin{blockarray}{rrr}\begin{block}{|ccc|}\vec{i}&\vec{j}&\vec{k}\\x_1&x_2&x_3\\y_1&y_2&y_3\\\end{block}\end{blockarray}\,.
\end{empheq}

\paragraph{Primjena vektorskog proizvoda.} Iz definicije vektorskog proizvoda proizilazi da je intenzitet vektorskog proizvoda $\vec{x}$ i $\vec{y}$ jednak 
povr\v sini paraleograma koju formiraju ova dva vektora. Ako ozna\v cimo sa $P$ ovu povr\v sinu, vrijedi
\begin{empheq}[box=\mymath]{equation*}
P=|\vec{x}\times\vec{y}|.
\end{empheq}
 \begin{figure}[!h]\centering
        \includegraphics[scale=.9]{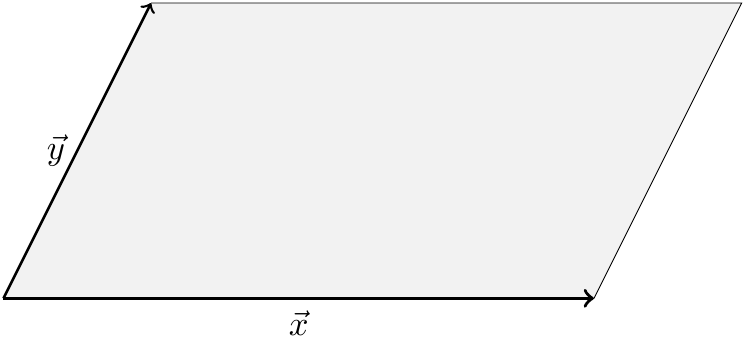}
         \caption{Povr\v sina paralelograma}
          \label{slika57}
  \end{figure}

U  slu\v caju da su vektori $\vec{x}$ i $\vec{y}$ paralelni ili ako imaju isti nosa\v c tada je
\[\vec{x}\times\vec{y}=\vec{0},\]
jer je $\measuredangle(\vec{x},\vec{y})=0$, pa je $\sin \measuredangle(\vec{x},\vec{y})=0$ ($\vec{x}\times\vec{y}=|\vec{x}||\vec{y}|\sin \measuredangle(\vec{x},\vec{y})\vec{n}_0$). Koriste\' ci vektorski proizvod mo\v zemo ispitati kolinearnost vektora, tj. vrijedi
\begin{empheq}[box=\mymath]{equation*}
\vec{x}\times\vec{y}=\vec{0}\Leftrightarrow (\vec{x}\parallel\vec{y},\,\vec{x}\neq \vec{0}\wedge\vec{y}\neq \vec{0}).
\end{empheq}

\begin{example}
Ako je $|\vec{a}|=10,\,|\vec{b}|=5,\,\vec{a}\cdot\vec{b}=12,$ izra\v cunati $|\vec{a}\times\vec{b}|.$\\\\
\noindent Rje\v senje:\\\\
  Vrijedi
  \[\vec{a}\cdot\vec{b}=|\vec{a}||\vec{b}|\cos\varphi \Leftrightarrow \cos\varphi=\frac{\vec{a}\cdot\vec{b}}{|\vec{a}||\vec{b}|}=\frac{12}{50}=\frac{6}{25},\]
  \[|\vec{a}\times\vec{b}|=|\vec{a}||\vec{b}||\sin\varphi|
             =5\cdot 10\sqrt{1-\cos^2\varphi}=5\cdot 10\sqrt{1-\left( \frac{6}{25}\right)^2}=2\sqrt{589}.\]
\end{example}

\begin{figure}[!h]\centering
        \includegraphics[scale=1.]{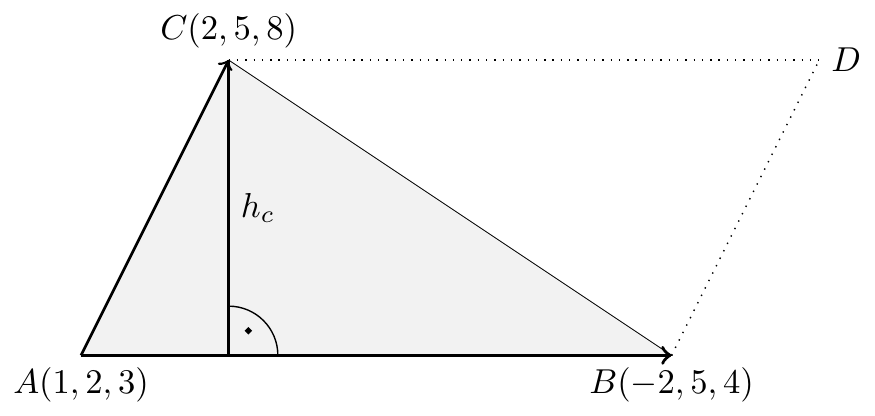}
         \caption{Povr\v sina trougla}
          \label{slika58}
  \end{figure}

\begin{example}
 Izra\v cunati povr\v sinu i visinu $h_C$ trougla $\triangle ABC,$ \v ciji su vrhovi\\
 $A(1,2,3),\,B(-2,5,4),\,C(2,5,8).$\\\\
\noindent Rje\v senje:\\\\
Vidi Sliku \ref{slika58}, povr\v sinu trougla ra\v cunamo
\[P_{\triangle ABC}=\frac{1}{2}P_{\square ABCD}=\frac{1}{2}|\overrightarrow{AB}\times\overrightarrow{AC}|.\]
Vrijedi
\begin{align*}
 \overrightarrow{AB}&=(-3,3,1)\\
 \overrightarrow{AC}&=(1,3,5)\\
 \overrightarrow{AB}\times\overrightarrow{AC}&=\begin{blockarray}{ccc}\begin{block}{|rrr|}\vec{i}&\vec{j}&\vec{k}\\-3&3&1\\1&3&5\\\end{block} \end{blockarray}
                  =12\vec{i}+16\vec{j}-12\vec{k}\\
P_{\triangle ABC}&= \frac{1}{2}|\overrightarrow{AB}\times\overrightarrow{AC}|=\frac{1}{2}\sqrt{12^2+16^2+12^2}=\frac{1}{2}\sqrt{544}=2\sqrt{34}.
\end{align*}
Povr\v sina trougla je $P_{\triangle ABC}=2\sqrt{34}.$

Visinu mo\v zemo izra\v cunati i na drugi na\v cin $P_{\triangle ABC}=\frac{1}{2}|\overrightarrow{AB}|\cdot h_c,$ pa je
\[h_c=\frac{2 P_{\triangle ABC}}{ | \overrightarrow{AB} | }=\frac{4\sqrt{34}}{\sqrt{9+9+1}}=\frac{4\sqrt{34}}{\sqrt{19}}=\frac{4\sqrt{646}}{19},\]
visina je $h_C=\frac{4\sqrt{646}}{19}.$
\end{example}

\begin{example}
 Izra\v cunati projekciju vektora $\vec{a}=(3,-12,4)$ na vektor $\vec{b}=\vec{c}\times\vec{d},$ ako je $\vec{c}=(1,0,-2),\,\vec{d}=(1,3,-4).$\\\\
\noindent Rje\v senje:\\\\
 Izra\v cunajmo prvo vektor $\vec{b}=\vec{c}\times\vec{d}$

 \[\vec{b}=\begin{blockarray}{ccc}\begin{block}{|rrr|}i&j&k\\1&0&-2\\1&3&-4\\\end{block} \end{blockarray}
             =6\vec{i}+2\vec{j}+3\vec{k},\]
sada je projekcija
\[\proj_{\overrightarrow{b}}\vec{a}=\frac{\vec{a}\cdot\vec{b}}{|\vec{b}|}=\frac{(3,-12,4)\cdot(6,2,3)}{\sqrt{6^2+2^2+3^2}}=\frac{6}{7}.\]
\end{example}

\begin{example}
Neka je $|\vec{a}|=5,\,|\vec{b}|=5,\,\measuredangle(\vec{a},\vec{b})=\frac{\pi}{4},$ izra\v cunati povr\v sinu paralelograma konstruisanog nad vektorima $\vec{m}=2\vec{b}-\vec{a}$ i $\vec{n}=3\vec{a}+2\vec{b}.$ \\\\
\noindent Rje\v senje:\\\\
   \begin{align*}
      |\vec{m}\times\vec{n}|&=|(2\vec{b}-\vec{a}) \times(3\vec{a}+2\vec{b})|=|6\vec{b}\times\vec{a}
                +\underbrace{4\vec{b}\times\vec{b}}_{=0}-\underbrace{3\vec{a}\times\vec{a}}_{=0}-2\vec{a}\times\vec{b}|\\
            &=|6\vec{b}\times\vec{a}+2\vec{b}\times\vec{a}|=8|\vec{b}\times\vec{a}|=8|\vec{b}||\vec{a}||\sin\measuredangle(\vec{a},\vec{b})|
                   =8\cdot 5\cdot 5\cdot\frac{\sqrt{2}}{2}=100\sqrt{2}.
   \end{align*}
\end{example}

\begin{example}
    Odrediti jedini\v cni vektor koji je normalan na ravan odre\dj enu   vektorima $\vec{a}=(1,1,0),$ $b=(1,-1,1).$\\\\
\noindent Rje\v senje:
  \begin{align*}
     \vec{n}&=\begin{blockarray}{ccc}\begin{block}{|rrr|}i&j&k\\1&1&0\\1&-1&1\\ \end{block}\end{blockarray}=\vec{i}-\vec{j}-2\vec{k},\\
     \vec{n}_0&=\frac{\vec{n}}{|\vec{n}|}=\frac{\vec{i}-\vec{j}-2\vec{k}}{\sqrt{1+1+4}}=\frac{1}{\sqrt{6}}\vec{i}-\frac{1}{\sqrt{6}}\vec{j}-\frac{2}{\sqrt{6}}\vec{k}
           =\left(\frac{1}{\sqrt{6}},-\frac{1}{\sqrt{6}},-\frac{2}{\sqrt{6}}\right).
  \end{align*}
\end{example}

\section{M\lowercase{je\v soviti proizvod tri vektora}}
\begin{figure}[!h]\centering
        \includegraphics[scale=1.]{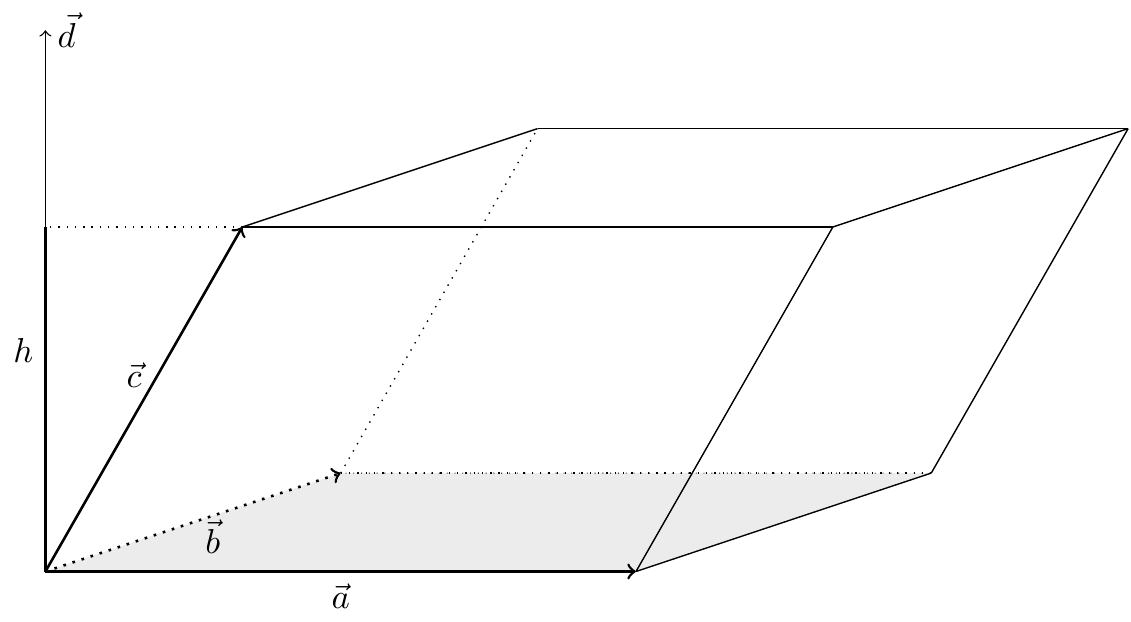}
         \caption{Zapremina paralelopipeda}
          \label{slika59}
  \end{figure}

Neka su data tri nekomplanarna vektora $\vec{a},\,\vec{b},\,\vec{c}$ i vektor $\vec{d}$ normalan na vektore $\vec{a}$ i $\vec{b}.$  Pomno\v zimo  li prvo $\vec{a}$ i $\vec{b}$ vektorski, a zatim $\vec{a}\times\vec{b}$ pomno\v zimo sa $\vec{c}$ skalarno, dobijamo   mje\v soviti proizvod
\begin{empheq}[box=\mymath]{equation*}
(\vec{a}\times\vec{b})\cdot\vec{c},
\end{empheq}\index{vektori! mje\v soviti proizvod tri vektora}
tri vektora $\vec{a},\,\vec{b}$ i $\vec{c}.$ Ova tri vektora konstrui\v su (razapinju) paralelopiped Slika \ref{slika59}. U ovom slu\v caju  $\vec{a},\,\vec{b}$ i $\vec{c},$ obrazuju triedar desne orjentacije. Povr\v sina  baze paralelopipeda (Slika \ref{slika59} paraleogram obojen sivo) jednaka je $|\vec{a}\times\vec{b}|,$ dok je visina $h$ jednaka $h=\proj_{\vec{d}}\vec{c}.$  Pa je sada

\[(\vec{a}\times\vec{b})\cdot\vec{c}=|\vec{a}\times\vec{b}|\vec{d}_0\cdot\vec{c}=|\vec{a}\times\vec{b}|\proj_{\vec{d}_0}\vec{c}=|\vec{a}\times\vec{b}|h=V,\]
gdje je $\vec{d}_0=\frac{\vec{d}}{|\vec{d}|}.$

Sa druge strane, u slu\v caju druga\v cije orjentacije, tj. kada vektori $\vec{a},\,\vec{b},\,\vec{c}$ obrazuju  triedar lijeve orjentacije, projekcija $\proj_{\vec{d}_0}\vec{c}$  bi bila negativna. U ovom slu\v caju,
\[V=-(\vec{a}\times\vec{b})\cdot\vec{c}.\]
Prema tome, zapreminu paralelopipeda ra\v cunamo
\begin{empheq}[box=\mymath]{equation*}
V=\pm(\vec{a}\times\vec{b})\cdot\vec{c},
\end{empheq}
gdje predznak biramo tako da zapremina ne bude negativna.

Ako su vektori $\vec{a},\,\vec{b},\,\vec{c}$ dati u ortonormiranoj bazi i ako vrijedi $\vec{a}=(a_1,a_2,a_3),\,\vec{b}=(b_1,b_2,b_3),\,\vec{c}=(c_1,c_2,c_3),$ tada zapreminu mo\v zemo ra\v cunati ovako:
\begin{empheq}[box=\mymath]{equation*}
V=\pm\begin{blockarray}{ccc}\begin{block}{|rrr|}a_1&a_2&a_3\\b_1&b_2&b_3\\c_1&c_2&c_3\\\end{block}\end{blockarray}\,.
\end{empheq}
\index{vektori! ra\v cunanje mje\v sovitog proizvoda}
\begin{remark}
 Tri nenulta vektora su komplanarna ako je mje\v soviti proizvod tih vektora jednak $0.$
\end{remark}

\ \\
\begin{example}
Odrediti parametar $t$ tako da vektori $\vec{a}=(\ln(t-2),-2,6),\,\vec{b}=(t,-2,5),\,\vec{c}=(0,-1,3)$ budu komplanarni.\\\\
\noindent Rje\v senje:\\\\
 U slu\v caju da su vektori $\vec{a},\,\vec{b},\,\vec{c}$ komplanarni, da le\v ze u istoj ravni, njihov mje\v soviti proizvod bio bi jednak $0.$    Pa iz uslova $(\vec{a}\times\vec{b})\cdot\vec{c}=0$ ra\v cunamo vrijednost parametra $t,$ vrijedi
 \begin{align*}
    (\vec{a}\times\vec{b})\cdot\vec{c}&=\begin{blockarray}{ccc}\begin{block}{|crr|}\ln(t-2)&-2&6\\ t&-2&5\\0&-1&3\\\end{block}\end{blockarray}
       =-6\ln(t-2)-6t+5\ln(t-2)+6t\\&=-\ln(t-2)=0,\\
       &\ln(t-2)=0\Leftrightarrow e^0=t-2\Leftrightarrow t=3.
 \end{align*}
\end{example}

\begin{figure}[!h]\centering
        \includegraphics[scale=.8]{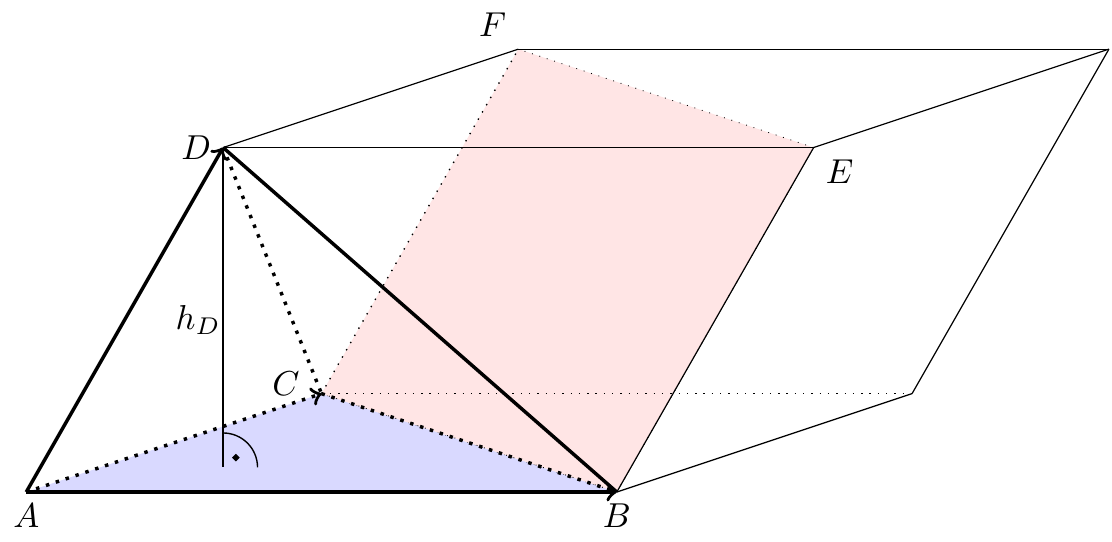}
         \caption{Zapremina tetraedra}
          \label{slika60}
  \end{figure}

\begin{example}
Odrediti zapreminu tetraedra \v ciji su vrhovi $A(3,1,-2),\,B(-4,2,3),\,C(1,5,-1),$\\$D(-5,-1,2).$ Kolika je visina tetraedra ako se za bazu uzme trougao $\triangle ABC$?\\ \ \\
\noindent Rje\v senje:\\\\
Vidi Sliku \ref{slika60}, paralelopiped kojeg razapinju vektori $\overrightarrow{AB},\,\overrightarrow{AC},\,\overrightarrow{AD}$ podijelimo na dva dijela sa ravni, \v ciji je dio predstavljen crvenim (osjen\v cenim) paralelogramom $\square BEFC.$ Sada treba da izra\v cunamo zapreminu tetraedra $V_T$ odre\dj enog vrhovima $A,B,C,D.$  Zapremina ovog tetredra predstavlja $\frac{1}{3}$ zapremine prizme $V_P$ odre\dj ene sa vrhovima $A,B,C,E,F,D.$

Pa je
\begin{align*}
  V_T&=\frac{1}{3}V_P=\frac{1}{3}\cdot\frac{1}{2}V=\frac{1}{6}V\\
  V&=\pm(\overrightarrow{AB}\times\overrightarrow{AC})\cdot\overrightarrow{AD}\\
  \overrightarrow{AB}&=(-7,1,5), \:  \overrightarrow{AC}=(-2,4,1), \: \overrightarrow{AD}=(-8,-2,4),\\
  (\overrightarrow{AB}\times\overrightarrow{AC})\cdot\overrightarrow{AD}&=
      \begin{blockarray}{ccc}\begin{block}{|rrr|}-7&1&5\\-2&4&1\\-8&-2&4\\ \end{block}\end{blockarray}=54,\\
    V_T&=\frac{1}{6}\cdot 54=9.
\end{align*}

Visinu ra\v cunamo na sljede\' ci na\v cin

\begin{align*}
   V_T&=\frac{1}{3}P_{\triangle ABC}\cdot h_D\Leftrightarrow h_D=\frac{3V_T}{P_{\triangle ABC}}\\
   P_{\triangle}&=\frac{1}{2}|\overrightarrow{AB}\times\overrightarrow{AC}|,\:
    \overrightarrow{AB}\times\overrightarrow{AC}=\begin{blockarray}{ccc}\begin{block}{|rrr|}\vec{i}&\vec{j}&\vec{k}\\-7&1&5\\-2&4&1\\\end{block}\end{blockarray}
       =-19\vec{i}-3\vec{j}-26\vec{k}\\
   |\overrightarrow{AB}\times\overrightarrow{AC}|&=\sqrt{(-19)^2+(-3)^2+(-26)^2}=\sqrt{1046}\\
   h_D&=\frac{54}{\sqrt{1046}}=\frac{54\sqrt{1046}}{1046}.
\end{align*}
\end{example}

\newpage

\section{Z\lowercase{adaci}} \index{Zadaci za vje\v zbu! vektori}

\paragraph{Linearna zavisnost vektora. Razlaganje vektora po bazi}

\begin{enumerate}
   \item 	
   \begin{inparaenum}
      \item  Da li su vektori $\vec{a}=4\vec{i}-6\vec{j}+10\vec{k}$ i $\vec{b}=-6\vec{i}+9\vec{j}-15\vec{k}$ kolinearni?	

      \item Odrediti parametre $u$ i $v$, tako da vektori $\vec{a}=(u,1,-2)$ i $\vec{b}=(-1,3,v),$ budu kolinearni.

      \item Ako je $\vec{a}=\left(\tfrac{2}{3},-\tfrac{3}{5},\tfrac{4}{3}\right)$, odrediti $x,z$ vektora $\vec{b}=\left(x,4,z\right)$,
         tako da $\vec{a}$  i $\vec{b}$ budu kolinearni.
   \end{inparaenum}
  \item
    \begin{enumerate}
       \item Odrediti parametar $k$, tako da vektori $\vec{a}=(-1,3,2),\,\vec{b}=(2,k,-4),\,\vec{c}=(k,12,6),$ budu komplanarni. Za
        takvu vrijednost parametra $k$ razlo\v ziti vektor $\vec{a}$ po pravcima vektora $\vec{b}$ i $\vec{c}$.
       \item Ispitati linearnu zavisnost vektora $\vec{a}=(1,-1,2),\,\vec{b}=(4,2,0)$ i $\vec{c}=(1,1,1).$
       \item Ispitati linearnu zavisnost vektora $\vec{l},\vec{m},\vec{n}$, a ako su zavisni razlo\v ziti vektor $\vec{l}$ na $\vec{m}$ i $\vec{n}$,

    \begin{inparaenum}[(i)]
        \item $\vec{l}=\left(2,-1,-1\right),\vec{m}=\left(-1,2,-1\right),\vec{n}=\left(-1,-1,2\right)$;\,

        \item $\vec{l}=\left(1,1,1\right),\vec{m}=\left(0,1,1\right),\vec{n}=\left(-1,0,1\right)$.
    \end{inparaenum}
  \end{enumerate}
\end{enumerate}

 \paragraph{Skalarni proizvod vektora}
  \begin{enumerate}\setcounter{enumi}{2}
  	\item
  	\begin{inparaenum}
    \item Dati su vektori $\vec{a}=(2,1,-2),\,\vec{b}=(-1,4,3).$
          Izra\v cunati $|\vec{a}|,\,|\vec{b}|,\,\vec{a}\cdot\vec{b},\measuredangle(\vec{a},\vec{b}),$ $\proj_{\vec{b}}\vec{a}$ i $\proj_{\vec{a}}\vec{b}.$

    \item Odrediti projekciju vektora $\vec{d}=4\vec{a}-3\vec{b}$ na vektor $\vec{c}$, ako su $\vec{a}=(2,-1,3),\,\vec{b}=(3,-2,1),\,\vec{c}=(2,-4,5).$

    \item Neka je $|a|=3,\,|b|=6,\,\measuredangle(\vec{a},\vec{b})=\frac{\pi}{2}.$ Izra\v cunati $|\vec{a}+\vec{b}|$ i $|2\vec{a}-\vec{b}|.$
    \end{inparaenum}
    \item
    \begin{inparaenum}
        \item Tjemena trougla su $A(2,-1,3),\,B(1,3,-4),\,C(0,2,4).$ Odrediti unutra\v snje\\ uglove trougla i du\v zine stranica trougla.

        \item Izra\v cunati du\v zinu dijagonala paralelograma, ako su mu stranice vektori $\vec{a}=2\vec{m}+4\vec{n},\,\vec{b}=-3\vec{m}+5\vec{n}$
          i $|\vec{m}|=2\sqrt{2},\,|\vec{n}|=3,\,\measuredangle(\vec{m},\vec{n})=\frac{\pi}{6}.$
    \end{inparaenum}

    \item Za koje su vrijednosti parametra $m$ vektori $\vec{a}=(-2,1,m),\,\vec{b}=(1,-2,3)$ ortogonalni.
  \end{enumerate}

 \paragraph{Vektorski proizvod vektora}
   \begin{enumerate}\setcounter{enumi}{5}
   	\item
   	 \begin{enumerate}
      \item Ako je $|\vec{a}|=5,\,|\vec{b}|=12$ i $\vec{a}\cdot \vec{b}=45,$ izra\v cunati $|\vec{a}\times \vec{b}|.$
      \item Ako je $\vec{p}=2\vec{a}+\vec{b},\,\vec{q}=\vec{a}-2\vec{b}$ i $|\vec{a}|=3,\,|\vec{b}|=4$ i $\measuredangle(\vec{a},\vec{b})=\frac{\pi}{6},$
            izra\v cunati $\vec{p}\times \vec{q}.$
     \end{enumerate}
     \item
    \begin{enumerate}
      \item Stranice paralelograma date su vektorima $\vec{p}=-3\vec{a}+2\vec{b},\,\vec{q}=3\vec{a}+4\vec{b}$  i  $|\vec{a}|=2,\,|\vec{b}|=3$
           i $\measuredangle(\vec{a},\vec{b})=\frac{\pi}{4}.$    Izra\v cunati povr\v sinu paralelograma i ugao izmedju dijagonala.
      \item Neka je $\overrightarrow{AB}=3\vec{p}-4\vec{q},\,\overrightarrow{BC}=\vec{p}+5\vec{q}$
            i $|\vec{q}|=|\vec{q}|=1,\,\measuredangle(\vec{p},\vec{q})=\frac{\pi}{4}.$ Izra\v cunati povr\v sinu trougla i visinu $h_c$.
    \end{enumerate}
    \item
    \begin{enumerate}
      \item Izra\v cunati projekciju $\vec{a}=(2,1,-3)$ na vektor $\vec{b}=\vec{c}\times \vec{a},$ ako je $\vec{c}=(1,0,-2)$ i $\vec{b}=1,3,-4).$
      \item Dati su vektori $\vec{a}=(1,1,-1),\vec{b}=(-2,-1,2),\,\vec{c}=(1,-1,2).$
           \begin{enumerate}[(i)]
             \item Razlo\v ziti vektor $\vec{c}$  po pravcima vektora $\vec{a},\vec{b},\,\vec{a}\times \vec{b};$
             \item Odrediti ugao koji obrazuju vektor $\vec{c}$ sa ravni odredjenom vektorima $\vec{a}$ i $\vec{b}$.
          \end{enumerate}
     \end{enumerate}
   \end{enumerate}

\paragraph{Mje\v soviti vektorski proizvod}

   \begin{enumerate}  \setcounter{enumi}{8}
      \item Ispitati da li ta\v cke $A(1,2,-1),\,B(0,1,5),\,C(-1,2,1),\,D(2,13)$ pripaju istoj ravni.
      \item
       \begin{inparaenum}[(i)]
         \item Dokazati da su vektori $\vec{a}=(-1,3,2),\,\vec{b}=(2,-3,4),\,\vec{c}=(-3,12,6)$ komplanarni.

         \item Izra\v cunati zapreminu tetraedra i visinu $h_D,$ \v ciji su vrhovi $A(2,0,0),\,B(0,3,0),$ $C(0,0,6),\,D(2,3,8).$
       \end{inparaenum}
   \end{enumerate}

\paragraph{Razni zadaci}
  \begin{enumerate}  \setcounter{enumi}{10}
    \item Dati su vektori $\vec{a}=-2\vec{p}+4\vec{q},\,\vec{b}=2\vec{p}-2\vec{q},\,\measuredangle(\vec{a},\vec{b})=\frac{\pi}{6}.$
          Izra\v cunati $|\vec{a}|,\,|\vec{b}|,\,\vec{a}\cdot\vec{b},\measuredangle(\vec{a},\vec{b}),\,\proj_{\vec{b}}\vec{a}$ i $\proj_{\vec{a}}\vec{b}.$
    \item Ako je $|\overrightarrow{AB}|=2,\,|\overrightarrow{AC}|=4,\,\measuredangle(\vec{AB},\overrightarrow{AC})=\frac{\pi}{3}$ i ta\v cka M je sredina du\v zi $BC$,
           izra\v cunati $\overrightarrow{AM}$ preko $\overrightarrow{AB},\,\overrightarrow{AC},$ a zatim $|\overrightarrow{AM}|.$
    \item
     \begin{enumerate}
       \item Izra\v cunati ugao izmedju simetrala koordinatnih osa $yOz$ i $xOz$.
       \item Izra\v cunati ugao izmedju vektora $2\vec{a}-4\vec{b}$ i $3\vec{a}+2\vec{b},$ ako je $\vec{a}=(2,-1,3),\,\vec{b}=(3,-2,2).$
     \end{enumerate}

    \item Odrediti projekciju vektora $\vec{a}=(2,4,\sqrt{5})$ na osu koja zaklapa uglove $\alpha=\frac{\pi}{3},\,\beta=\frac{\pi}{6}$ a sa $z$-osom tup ugao.

\item
  \begin{inparaenum}
  	\item U trouglu $\Delta \,ABC$ poznato je $A(2,1,-3),\,\overrightarrow{AB}=(2,-3,5),\,\overrightarrow{BC}=(3,-2,4).$ Odrediti koordinate vrhova $B$ i $C,$ koeficijente
  	vektora $\overrightarrow{AC}$ i vrijednost ugla $\gamma. $
  	
  	\item Izra\v cunati povr\v sinu trougla \v cije su dvije stranice vektori $\vec{a}=4\vec{j}-\vec{k}+3\vec{k}$ i $\vec{b}=(5,-3,7).$
  	
  	\item Dat je trougao \v ciji su vrhovi $A(5,2,-4),\,B(9,-8,-3)$ i $C(16,-6,-11).$ Odrediti povr\v sinu $\triangle\,ABC$ i unutra\v snje uglove trougla.
  \end{inparaenum}

   \item
   \begin{inparaenum}
       \item Odrediti parametar $\lambda$ tako da vektori $\vec{a}=14\vec{j}+\lambda \vec{k}+3\vec{i}$ i $\vec{b}=(\lambda,-2,\lambda)$ budu okomiti.

       \item Odrediti parametar $k$ tako da vektori $\vec{a}=(-1,3,2),\,\vec{b}=(2,k,-4)$ i $\vec{c}=(k,12,6)$ budu komplanarni. Za tako dobijenu vrijednost parametra $k$ razlo\v ziti vektor $\vec{a}$ po pravcima vektora $\vec{b}$ i $\vec{c}$.
   \end{inparaenum}

\item Dati su vektori $\vec{a}=(2,-1,7),\,\vec{b}=(0,-4,3)$ i $\vec{c}=(-1,-2+1).$ Ispitati komplanarnost vektora. Ako nisu komplanarni izra\v cunati visinu tetraedra koju obrazuju vektori $\vec{a},\,\vec{b}$ i $\vec{c}$, ako se za bazu uzmu vektori $\vec{a}$ i $\vec{b}$.

\item Izra\v cunati du\v zinu visine trougla koju obrazuju vektori $\vec{a}=(1,-4,5)$ i $\vec{b}=(0,2,-4)$ i to onu koja odgovara stranici koja je odredjena vektorom $\vec{a}.$

\item
    \begin{inparaenum}
    	\item Izra\v cunati povr\v sinu nad vektorima $\vec{a}=(2,1,-k)$ i $\vec{b}=(-1,1,-7)$ i izra\v cunati ugao izmedju vektora $\vec{a}$ i $\vec{a}+\vec{2b}$.
    	    	
        \item Dati su vektori $\overrightarrow{OA}=(5,2,-1),\,\overrightarrow{OB}=(1,-3,4),\,\overrightarrow{OC}=(-2,1,3),\,\overrightarrow{OD}=(2,6,-2).$  Pokazati da je $\square\,ABCD$ paralelogram i izra\v cunati ugao izmedju dijagonala.

        \item Data su tri tjemena paralelograma $A(2,-1,0),\,B(3,4,1),\,C(-1,0,-1).$ Odrediti \v cetvrto tjeme i ugao izmedju dijagonala.
    \end{inparaenum}

\item Izra\v cunati visinu paralelopipeda kojeg obrazuju vektori $\vec{a}=(1,-4,5),\,\vec{b}=(0,2,-4)$ i $\vec{c}=(2,-2,1).$

\end{enumerate}

\chapter{Ravan i prava}



 Ravan i prava su osnovni objekti u geometriji, a u ovom poglavlju cilj nam je predstaviti ih i analiti\v ckim putem, tj. predstaviti ih nekim formulama. 
 Ovakav pristup znatno olak\v sava i pro\v siruje primjenu ovih geometrijskih objekata. Da bi ostvarili pomenuti cilj koristimo se ve\'c ste\v cenim znanjem iz matematike. Veliki zna\v caj u analiti\v ckom prestavljanju ravni i prave imaju vektori, determinante i sistemi linearnih algebarskih jedna\v cina, koji su obra\dj eni u prethodnim poglavljima.

\section[R\lowercase{avan}]{Ravan}
Ravan mo\v zemo predstaviti na vi\v se na\v cina analiti\v ckim putem. U nastavku su dati oblici jedna\v cine ravni koji se naj\v ce\v s\'ce koriste.

\subsection{Jedna\v cine ravni}
\begin{figure}[!h]\centering
        \includegraphics[scale=1.]{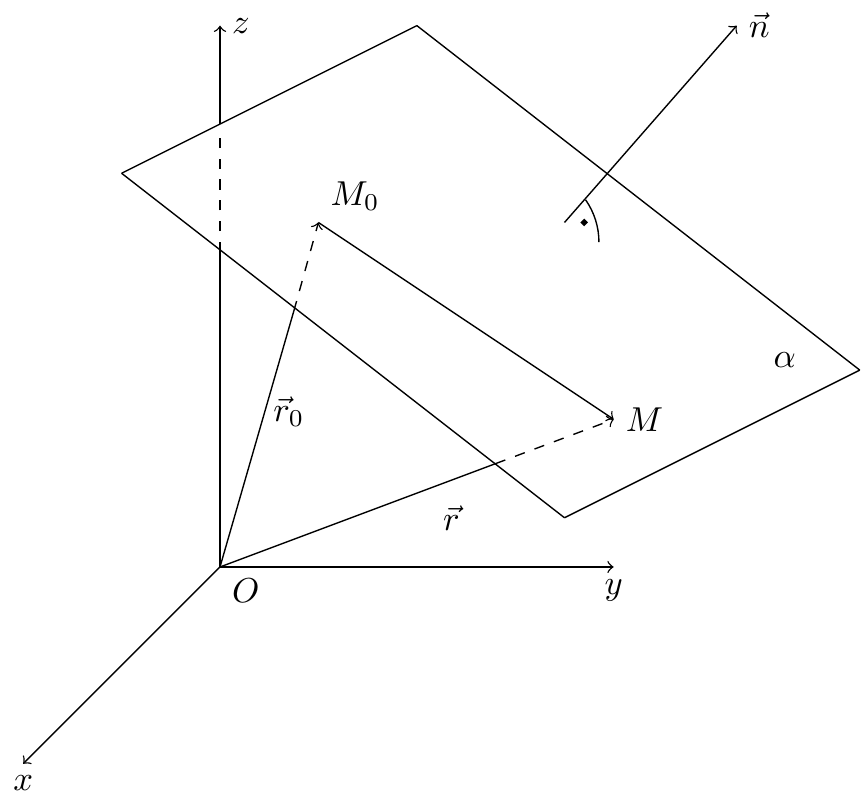}
        \caption{Ravan $\alpha$}
     \label{slika61b}
  \end{figure}

Posmatrajmo Sliku \ref{slika61b}.
U pravouglom Descartesovom koordinatnom sistemu data je ravan $	\alpha,$ odre\dj ena ta\v ckom $M_0(x_0,y_0,z_0)$ i vektorom $\vec{n}=(A,B,C),$ koji je normalan na ravan $\alpha.$ Neka je $M(x,y,z)$ proizvoljna ta\v cka ravni $\alpha$ i neka je $\vec{r}=\overrightarrow{OM}$ njen vektor polo\v zaja (vektor \v ciji je po\v cetak u ta\v cki $O$ a zavr\v setak u ta\v cki $M$). Vektor polo\v zaja ta\v cke $M_0$ je $\vec{r}_0=\overrightarrow{OM_0}.$ Iz trougla $\triangle OMM_0$ vidimo da je $\overrightarrow{M_0M}=\overrightarrow{OM}-\overrightarrow{OM_0}=\vec{r}-\vec{r}_0.$  Kako je vektor $\vec{n}$ normalan na ravan $\alpha$ samim tim normalan je i na svaku pravu i vektor koji le\v ze u toj ravni, te je tako vektor $\vec{n}$ normalan i na vektor $\vec{r}-\vec{r}_0.$   Iz poglavlja sa vektorima poznat je uslov $\vec{x}\perp\vec{y}\Leftrightarrow (\vec{x}\cdot\vec{y}=0,\,\vec{x}\neq \vec{0},\vec{y}\neq \vec{0}).$ Ako sada iskoristimo ovaj uslov, dobijamo

\begin{equation}
    \vec{n}(\vec{r}-\vec{r}_0)=0.
 \label{ravan1}
\end{equation}
Ako iskoristimo oznaku $\vec{n}\cdot\vec{r}_0=a,$ onda \eqref{ravan1} mo\v zemo pisati u obliku

\begin{equation}
   \tcbhighmath[mojstil1]{ \vec{n}\cdot\vec{r}=a.}
 \label{ravan2}
\end{equation}
Jedna\v cina \eqref{ravan2} naziva se op\v sta vektorska jedna\v cina ravni.

\index{ravan! op\v sta vektorska jedna\v cina ravni}

Napi\v simo sada vektore $\vec{n},\vec{r}$ i $\vec{r}_0$ preko koordinata: $ \vec{n}=A\vec{i}+B\vec{j}+C\vec{k},\,\vec{r}= x\vec{i}+y\vec{j}+z\vec{k},\,\vec{r}_0=
x_0\vec{i}+y_0\vec{j}+z_0\vec{k}.$ Sada iz \eqref{ravan1} dobijamo
\[(A\vec{i}+B\vec{j}+C\vec{k})\cdot[(x-x_0)\vec{i}+(y-y_0)\vec{j}+(z-z_0)\vec{k}]=0,\]odnosno
\[\tcbhighmath[mojstil1]{A(x-x_0)+B(y-y_0)+C(z-z_0)=0.}\]
Ova jedna\v cina predstavlja skalarnu jedna\v cinu ravni kroz datu ta\v cku sa poznatim vekto-\\rom normale. U ovoj  jedna\v cini $x,y,z$ su koordinate bilo koje ta\v cke ravni,  $x_0,y_0,z_0$ su koordinate poznate ta\v cke, dok su $A,B,C$ koeficijenti vektora normale. Dalje je
\[Ax+By+Cz-(Ax_0+By_0+Cz_0)=0, \]
pa ozna\v cimo li $-(Ax_0+By_0+Cz_0)=D,$ dobijamo
\begin{equation}
\tcbhighmath[mojstil1]{   Ax+By+Cz+D=0.}
\label{ravan3}
\end{equation}
Jednakost \eqref{ravan3} predstavlja op\v stu skalarnu jedna\v cinu ravni.

\index{ravan! op\v sta skalarna jedna\v cina ravni}

\begin{figure}[!h]\centering
        \includegraphics[scale=1.]{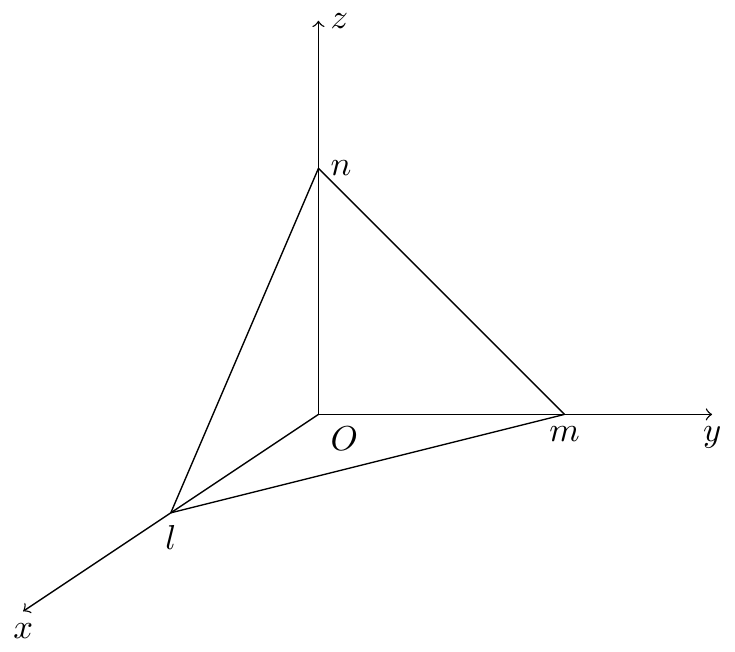}
        \caption{Ravan $\alpha$ sa odsje\v ccima (segmentima) $l,m,n$ }
     \label{slika62}
  \end{figure}

Ako je $A\neq 0,\,B\neq 0,\, C\neq 0,\,D\neq 0,$ onda vrijedi

\begin{align}
        Ax+By+Cz+D=0\Leftrightarrow
       &Ax+By+Cz=-D\Leftrightarrow
       \frac{x}{-\frac{D}{A}}+\frac{y}{-\frac{D}{B}}+\frac{z}{-\frac{D}{C}}=1\Leftrightarrow\nonumber\\
       &\tcbhighmath[mojstil1]{\frac{x}{l}+\frac{y}{m}+\frac{z}{n}=1,}\label{ravan4}
\end{align}
gdje su $-\frac{D}{A}=l,\,-\frac{D}{B}=m,\,-\frac{D}{C}=n.$ Jednakost \eqref{ravan4} predstavlja segmentni oblik jedna\v cine ravni.  Brojevi $l,m,n$ su odje\v cci ravni $\alpha$ na koordinatnim osama $x,y,z,$ respektivno (Slika \ref{slika62}).

\index{ravan! segmenti oblik jedna\v cine ravni}

Jedini\v cni vektor, u oznaci $\vec{a}_0,$ proizvoljnog vektora $\vec{a}$ dobijamo tako \v sto vektor $\vec{a}$ (koji nije nula vektor) podijelimo sa njegovim intenzitetom $|\vec{a}|,$ pa je
\[\tcbhighmath[mojstil1]{\vec{a}_0=\frac{\vec{a}}{|\vec{a}|}.}\]

\begin{figure}[!h]\centering
       \begin{subfigure}[b]{.3\textwidth}\centering
        \includegraphics[scale=.8]{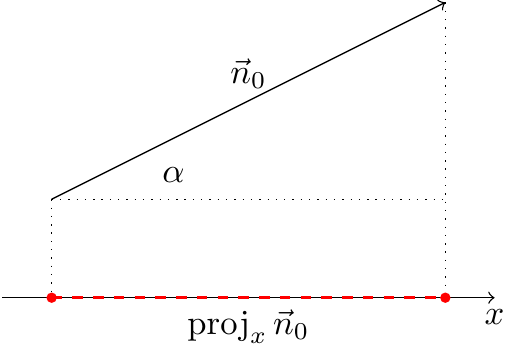}
         \caption{}
          \label{slika61}
   \end{subfigure}
   \begin{subfigure}[b]{.3\textwidth}\centering
        \includegraphics[scale=.8]{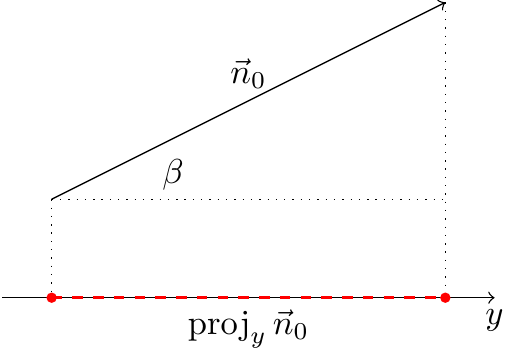}
         \caption{}
      \label{slika36}
     \end{subfigure}
     \begin{subfigure}[b]{.3\textwidth}\centering\hspace{.5cm}
        \includegraphics[scale=.8]{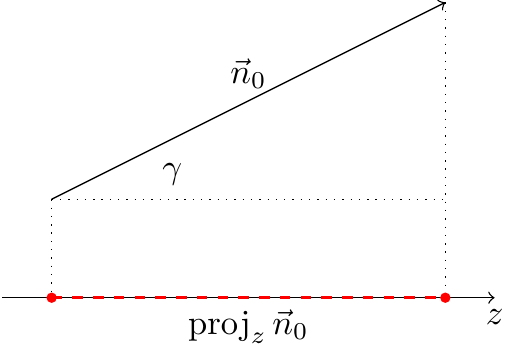}
         \caption{}
      \label{slika37}
     \end{subfigure}
        \caption{Projekcije jedini\v cnog vektora $\vec{n}_0$ na $x,y$ i $z$ osu}
     \label{slika63}
  \end{figure}

Neka je $\vec{n}_0$ jedini\v cni vektor i neka vektor $\vec{n}_0$ zahvata uglove $\alpha,\beta,\gamma$ sa $x,y,z$ osama respektivno.
Tada vrijedi (vidjeti Slike \ref{slika63})
\begin{align*}
    \tcbhighmath[mojstil1]{\proj_{x}\vec{n}_0=\cos\alpha,\:\proj_{y}\vec{n}_0=\cos\beta,\:\proj_{z}\vec{n}_0=\cos\gamma,}\:
\end{align*}
i
\begin{align*}
   \tcbhighmath[mojstil1]{ \vec{n}_0=\cos\alpha\vec{i}+\cos\beta\vec{j}+\cos\gamma\vec{k}.}
\end{align*}

\index{ravan! jedini\v cni vektor normale ravni}


Kako je $\vec{n}_0$ jedini\v cni vektor ($|\vec{n}_0|=1$), vrijedi $|\vec{n}_0|=\sqrt{\cos^2\alpha+\cos^2\beta+\cos^2\gamma}$ i
\[\tcbhighmath[mojstil1]{\cos^2\alpha+\cos^2\beta+\cos^2\gamma=1.}\]

Posmatrajmo Sliku \ref{slika67}.
 \begin{figure}[!h]\centering
       \begin{subfigure}[b]{.5\textwidth}\centering
        \includegraphics[scale=1]{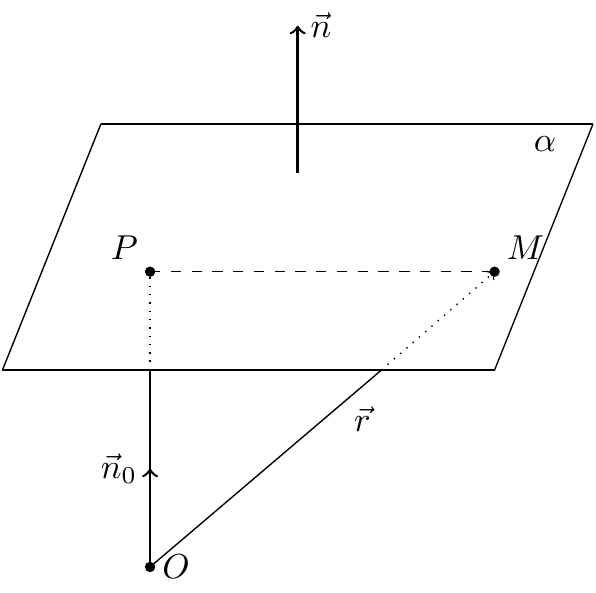}
         \caption{}
          \label{slika61a}
   \end{subfigure}
   \begin{subfigure}[b]{.5\textwidth}\centering
        \includegraphics[scale=1]{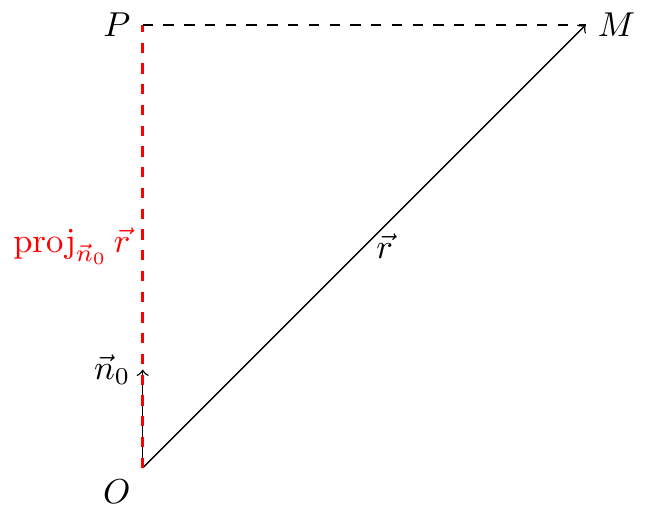}
         \caption{}
      \label{slika67}
     \end{subfigure}
        \caption{Projekcije vektora $\vec{r}$ na jedini\v cni vektor $\vec{n}_0$}
     \label{slika67}
  \end{figure}

Vektoru normale $\vec{n}$ odgovara jedini\v cni vektor $\vec{n}_0.$ Po\v sto radimo sa slobodnim vekto-\\rima, prenesimo vektor $\vec{n}_0$ u koordinatni po\v cetak, kao na Slici \ref{slika61a}.  Du\v zina du\v zi $|OP|,$ tj. rastojanje ravni od koordinatnog po\v cetka, jednaka je apsolutnoj vrijednosti projekcije vektora $\vec{r}$ na osu koja je odre\dj ena vektorom $\vec{n}_0.$ Ova osa normalna je na ravan $\alpha.$   Dakle, vrijedi

\[\proj_{\vec{n}_0}\vec{r}=\pm|OP|=\frac{\vec{n}_0\cdot\vec{r}}{|\vec{n}_0|}=\vec{n}_0\cdot\vec{r}.\]
Ozna\v cimo ovu projekciju sa $p$, pa je sada
\[p=\vec{n}_0\cdot\vec{r}.\]
Ako se sada vratimo u jedna\v cinu ravni \eqref{ravan2}, $\vec{n}\cdot \vec{r}=a$ i podijelimo je sa $|\vec{n}|$ dobijamo
\[\frac{\vec{n}\cdot\vec{r}}{|\vec{n}|}=\frac{a}{|\vec{n}|}\Leftrightarrow \frac{\vec{n}}{|\vec{n}|}\cdot\vec{r}=\frac{a}{|\vec{n}|}
      \Leftrightarrow\vec{n}_0\cdot\vec{r}=\frac{a}{|\vec{n}|}.\]
S obzirom da je $\vec{n}_0\cdot\vec{r}=p,$    to
\begin{equation}
\tcbhighmath[mojstil1]{\vec{n}_0\cdot\vec{r}=p,}
\label{ravan5}
\end{equation}
predstavlja normalnu jedna\v cina ravni u vektorskom obliku.

\index{ravan! normalna jedna\v cina ravni\\ u vektorskom obliku}

Poslije skalarnog mno\v zenja vektora $\vec{r}=x\vec{i}+y\vec{j}+z\vec{k}$ i $\vec{n}_0=\cos\alpha\vec{i}+\cos\beta\vec{j}+\cos\gamma\vec{k}$ u \eqref{ravan5}, dobijamo
\begin{equation*}
\tcbhighmath[mojstil1]{  x\cos\alpha+y\cos\beta+z\cos\gamma-p=0. }
\label{ravan6}
\end{equation*}
Ovo je normalna jedna\v cina ravni u skalarnom obliku.

\index{ravan! normalna jedna\v cina ravni u skalarnom obliku}

\begin{figure}[!h]\centering
        \includegraphics[scale=1]{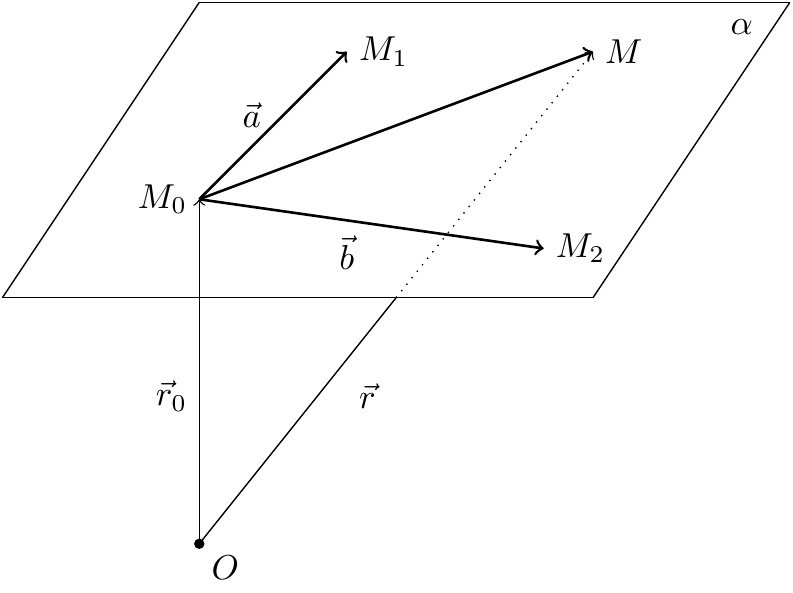}
        \caption{Tri ta\v cke u ravni    }
     \label{slika68}
  \end{figure}

Da bi izveli jo\v s jedan oblik jedna\v cine ravni, posmatrajmo Sliku \ref{slika68}. Neka je data ravan $\alpha$ i tri ta\v cke te ravni i to: $M_0(x_0,y_0,z_0),\,M_1(x_1,y_1,z_1)$ i $M_2(x_2,y_2,z_2)$. Neka vrijedi i $\vec{a}=\overrightarrow{M_0M_1},\,\vec{b}=\overrightarrow{M_0M_2},\,\vec{r}_0=\overrightarrow{OM_0},\,\vec{r}=\overrightarrow{OM}$ i $M\in\alpha.$ Jasno je da je $\overrightarrow{M_0M}=\vec{r}-\vec{r}_0,$  te da su vektori $\vec{a},\,\vec{b}$ i $\vec{r}-\vec{r}_0$ linearno zavisni po\v sto pripadaju istoj ravni $\alpha.$ Zbog toga vektor $\vec{r}-\vec{r}_0$ mo\v zemo predstaviti kao linearnu kombinaciju vektora $\vec{a}$ i $\vec{b}.$ Drugim rije\v cima, postoje realni brojevi $u$ i $v$ (zva\' cemo ih parametri) takvi da je
\[\vec{r}-\vec{r}_0=u\vec{a}+v\vec{b},\]
ili
\begin{equation}
\tcbhighmath[mojstil1]{   \vec{r}=\vec{r}_0+u\vec{a}+v\vec{b.}}
\label{ravan7}
\end{equation}

\index{ravan! parametarska jedna\v cina ravni u vektorskom obliku}

Mijenjaju\' ci koordinate ta\v cke $M$ ($M$ ostaje u ravni $\alpha$) mijenjaju se i vrijednosti parame-\\tara $u$ i $v,$ zato jedna\v cinu \eqref{ravan7} zovemo vektorska parametarska jedna\v cina ravni.

Prelaskom na koordinate, iz \eqref{ravan7} dobijamo
\begin{align}
  x&=x_0+ua_1+vb_1\nonumber\\
  y&=y_0+ua_2+vb_2\label{ravan8}\\
  z&=z_0+ua_3+vb_3,\nonumber
\end{align}
gdje su $\vec{a}=(a_1,a_2,a_3)$ i $\vec{b}=(b_1,b_2,b_3).$ Jedna\v cine \eqref{ravan8} predstavljaju parametarske jedna\v cine ravni.

\index{ravan! parametarske jedna\v cine ravni\\ u skalarnom obliku}

Ako iskoristimo uslov komplanarnosti vektora $\overrightarrow{M_0M},\,\overrightarrow{M_0M_1},\,\overrightarrow{M_0M_2},$ \\ tj.  $\left(\overrightarrow{M_0M}\times\overrightarrow{M_0M_1}\right)\cdot\overrightarrow{M_0M_2}=0,$ dobijamo jedna\v cinu ravni kroz tri ta\v cke
\begin{equation*}
 \tcbhighmath[mojstil1]{\begin{blockarray}{ccc}\begin{block}{|ccc|}x-x_0&y-y_0&z-z_0\\x_1-x_0&y_1-y_0&z_1-z_0\\x_2-x_0&y_2-y_0&z_2-z_0\\ \end{block}\end{blockarray}=0.}
\label{ravan9}
\end{equation*}

\index{ravan! jedna\v cina ravni kroz tri ta\v cke}

\begin{figure}[!h]\centering
        \includegraphics[scale=1]{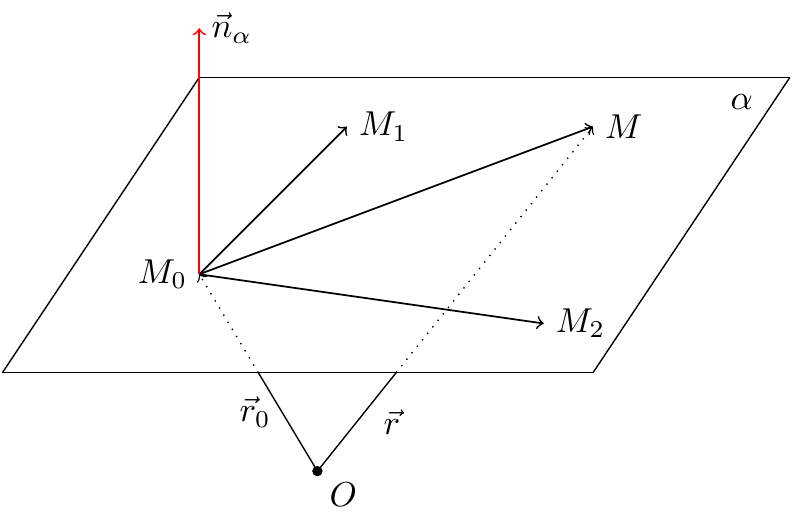}
        \caption{Ravan $\alpha$}
     \label{slika70}
  \end{figure}

\begin{example}
   Napisati sve oblike jedna\v cine ravni (koje su izvedene u prethodnom dijelu ovog poglavlja) koja je odre\dj ena ta\v ckama $M_0(1,1,0),\,M_1(-2,0,4),\,M_2(2,3,-1).$\\\\
\noindent Rje\v senje:\\\\
Vidi Sliku \ref{slika70}, odredimo prvo vektore $\vec{r}_0=\overrightarrow{OM_0},\,\vec{r}=\overrightarrow{OM},\,\overrightarrow{M_0M_1},\,\overrightarrow{M_0M_2}$
\begin{align*}
  \vec{r}_0&=(1,1,0)\\
  \vec{r}&=(x,y,z)\\
  \overrightarrow{M_0M_1}&=(-3,-1,4)\\
  \overrightarrow{M_0M_2}&=(1,2-1),
\end{align*}
sada vektor normale $\vec{n}_{\alpha}$
\begin{align*}
   \vec{n}_{\alpha}&=\overrightarrow{M_0M_2}\times\overrightarrow{M_0M_1}
        =\begin{blockarray}{rrr}\begin{block}{|rrr|}\vec{i}&\vec{j}&\vec{k}\\1&2&-1\\-3&-1&4\\\end{block}\end{blockarray}=(7,-1,5).
\end{align*}
\begin{enumerate}[$\bullet$]
\item Za oblik jedna\v cine ravni $\vec{n}_\alpha\cdot(\vec{r}-\vec{r}_0)=0$
  \[(7\vec{i}-\vec{j}+5\vec{k})\cdot(\vec{r}-\vec{i}-\vec{j})=0;\]
\item za $\vec{n}_\alpha\cdot\vec{r}=a$ je
\[(-7\vec{i}-\vec{j}+5\vec{k})\cdot\vec{r}=6;\]
\item za $A(x-x_0)+B(y-y_0)+C(z-z_0)=0$ je
\[7(x-1)-(y-1)+5(z-0)=0;\]
\item za $Ax+By+Cz+D=0$ je
\[7x-y+5z-6=0;\]
\item za $\frac{x}{l}+\frac{y}{m}+\frac{z}{n}=1$ je
\[\frac{x}{\frac{6}{7}}+\frac{y}{-6}+\frac{z}{\frac{6}{5}}=1;\]

\item za $x\cos\alpha+y\cos\beta+z\cos\gamma-p=0,$ vrijedi
     \begin{align*}
        p&=\frac{\vec{n}_\alpha\cdot\vec{r}_0}{ \left| \vec{n}_\alpha\right|}=\frac{(7,-1,5)(1,1,0)}{\sqrt{49+1+25}}=\frac{6}{\sqrt{75}}\\
        \cos\alpha&=\frac{A}{|\vec{n}_\alpha|}=\frac{7}{\sqrt{75}}\\
        \cos\beta&=\frac{B}{|\vec{n}_\alpha|}=\frac{-1}{\sqrt{75}}\\
        \cos\gamma&=\frac{C}{|\vec{n}_\alpha|}=\frac{5}{\sqrt{75}},
     \end{align*}
pa je jedna\v cina ravni u ovom slu\v caju
\[x\frac{7}{\sqrt{75}}+y\frac{-1}{\sqrt{75}}+z\frac{5}{\sqrt{75}}-\frac{6}{\sqrt{75}}=0;\]

\item za $\vec{r}=\vec{r}_0+u\overrightarrow{M_0M_2}+v\overrightarrow{M_0M_1}$ je
     \[\vec{r}=\vec{i}+\vec{j}+u(\vec{i}+2\vec{j}-\vec{k})+v(-3\vec{i}-\vec{j}+4\vec{k});\]
\item parametarski oblik
             \begin{align*}
               x&=1+u-3v\\
               y&=1+2u-v\\
               z&=-u+4v;
             \end{align*}

\item jedna\v cina ravni kroz tri ta\v cke
       \[\begin{blockarray}{ccc}\begin{block}{|rrr|}x-1&y-1&z-0\\-2-1&0-1&4-0\\2-1&3-1&-1-0\\ \end{block}\end{blockarray}=0.\]
\end{enumerate}
\end{example}

\subsection{Rastojanje ta\v cke od ravni}
\begin{figure}[!h]\centering
        \includegraphics[scale=1]{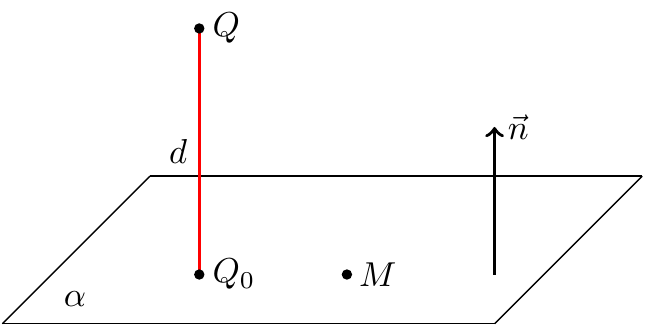}
        \caption{Rastojanje ta\v cke od ravni}
     \label{slika69}
  \end{figure}

Neka je data ravan $\alpha:Ax+By+Cz+D=0$ i neka ta\v cka $M(x_2,y_2,z_2)$ pripada ravni $\alpha,$ tada vrijedi $Ax_2+By_2+Cz_2+D=0,$ tj.  koordinate ta\v cke $M$ zadovoljavaju jedna\v cinu ravni $\alpha.$ Zanima nas kako da izra\v cunamo rastojanje neke ta\v cke, koja ne pripada ravni, od te ravni. Posmatrajmo Sliku \ref{slika69}. Ta\v cka $Q(x_1,y_1,z_1)$ ne pripada ravni $\alpha$. Kako da izra\v cunamo rastojanje $d$ ta\v cke $Q$ od ravni $\alpha?$
Spustimo normalu iz ta\v cke $Q$ na ravan, dobijamo ta\v cku $Q_0(x_0,y_0,z_0).$ Po\v sto ta\v cka $Q_0$ pripada ravni $\alpha$ njene koordinate zadovoljavaju jedna\v cinu ravni
\[Ax_0+By_0+Cz_0+D=0,\]
ali $Q\notin\alpha$, pa je $Ax_1+By_1+Cz_1+D\neq 0.$ Vrijedi $\vec{n}=A\vec{i}+B\vec{j}+C\vec{k},$  dok za vektor \v cija je po\v cetna ta\v cka $Q_0$ a krajnja $Q,\:$ vrijedi  $\overrightarrow{Q_0Q}=(x_1-x_0)\vec{i}+(y_1-y_0)\vec{j}+(z_1-z_0)\vec{k}.$ Sada je
\begin{align*}
 Ax_1+By_1+Cz_1+D&=Ax_1+By_1+Cz_1+D-(\underbrace{Ax_0+By_0+Cz_0+D}_{=0})\\
               &=A(x_1-x_0)+B(y_1-y_0)+C(z_1-z_0)=\vec{n}\cdot\overrightarrow{Q_0Q}\\
               &=|\vec{n}|\left| \overrightarrow{Q_0Q}\right|\cos\measuredangle(\vec{n},\overrightarrow{Q_0Q}),
\end{align*}
a kako su vektori $\vec{n}$ i $\overrightarrow{Q_0Q}$ kolinearni, to mo\v ze biti $\cos\measuredangle(\vec{n},\overrightarrow{Q_0Q})=\pm 1,$ i $\left|\overrightarrow{Q_0Q}\right|=d$  pa je
\[Ax_1+By_1+Cz_1+D=\pm|\vec{n}|d,\]
i zbog \[|\vec{n}|=\sqrt{A^2+B^2+C^2}\]
dobijamo
\begin{equation}
  \tcbhighmath[mojstil1]{d=\left|\frac{Ax_1+By_1+Cz_1+D}{\sqrt{A^2+B^2+C^2}} \right|. }
\label{ravan10}
\end{equation}
Jednakost \eqref{ravan10} predstavlja obrazac za ra\v cunanje rastojanja ta\v cke od ravni.

\index{ravan! rastojanje ta\v cke od ravni}

\begin{example}
Izra\v cunati visinu piramide $h_s$ \v ciji su vrhovi $S(0,6,4),\,A(3,5,3),\,B(-2,11,-5),$\\ $C(1,-1,4).$\\\\
\noindent Rje\v senje:\\\\
 Visinu mo\v zemo izra\v cunati kao rastojanje ta\v cke $S$ od ravni koja je odre\dj ena ta\v ckama $A,B,C.$ Odredimo prvo jedna\v cinu ravni kroz ove tri posljednje ta\v cke. Vrijedi

\begin{align*}
   \begin{blockarray}{ccc}\begin{block}{|ccc|} x-x_0&y-y_0&z-z_0\\x_1-x_0&y_1-y_0&z_1-z_0\\
           x_2-x_0&y_2-y_0&z_2-z_0\\\end{block}\end{blockarray}=
    \begin{blockarray}{ccc}\begin{block}{|ccc|}x-3&y-5&z-3\\-5&6&-8\\-2&-6&1\\ \end{block}
\end{blockarray}=0,
\end{align*}
 pa je $2x-y-2z+5=0,$ dakle $A=2,\,B=-1,\,C=-2,\,D=5.$ Rastojanje je
 \begin{align*}
   d&=\left| \frac{Ax_s+By_s+Cz_s+D}{\sqrt{A^2+B^2+C^2}}\right|=\left|\frac{0-6-8+5}{\sqrt{4+1+4}}\right|=3.
 \end{align*}
\end{example}

\subsection{Ugao izme\dj u dvije ravni}
Neka su date dvije ravni
\begin{align*}
  \alpha:A_1x+B_1y+C_1z+D_1=0\\
  \beta:A_2x+B_2y+C_2z+D_2=0.
\end{align*}

Ako se ravni $\alpha$ i $\beta$ sijeku tada se ugao  izme\dj u ove dvije ravni $\phi,$ defini\v se kao ugao izme\dj u njihovih vektora normala $\vec{n}_{\alpha}$ i $\vec{n}_{\beta}.$ Kako je
\begin{align*}
   \vec{n}_{\alpha}&=(A_1,B_1,C_1),\\
   \vec{n}_{\beta}&=(A_2,B_2,C_2),
\end{align*}
i $\cos\phi=\frac{\vec{n}_{\alpha}\cdot\vec{n}_{\beta}}{|\vec{n}_{\alpha}||\vec{n}_{\beta}|},$ to vrijedi
\begin{equation*}
  \tcbhighmath[mojstil1]{\cos\phi=\frac{A_1A_2+B_1B_2+C_1C_2}{\sqrt{A^2_1+B^2_1+C^2_1}\sqrt{A^2_2+B^2_2+C^2_2}}.}
\label{ravan11}
\end{equation*}

\index{ravan! ugao izme\dj u dvije ravni}

U slu\v caju da su ravni paralelne $\alpha\parallel\beta,$ vektori normala $\vec{n}_{\alpha}$ i $\vec{n}_{\beta}$ su kolinearni, \v sto je ekvivalentno uslovu
\begin{equation*}
     \tcbhighmath[mojstil1]{ \frac{A_1}{A_2}=\frac{B_1}{B_2}=\frac{C_1}{C_2}.}
\label{ravan12}
\end{equation*}

\index{ravan! uslov paralelnosti}

Ako su ravni $\alpha$ i $\beta$ okomite, tj. $\alpha\perp\beta,$ to je ekvivalentno uslovu
\begin{equation*}
    \tcbhighmath[mojstil1]{A_1A_2+B_1B_2+C_1C_2=0.}
\label{ravan13}
\end{equation*}

\index{ravan! uslov okomitosti}

\begin{example}
 Izra\v cunati ugao izme\dj u ravni $\alpha:x+3y-4z+5=0$ i $\beta:2x+2y+2z-7=0.$ \\\\
\noindent Rje\v senje:\\\\
  Vektori normala su $\vec{n}_\alpha=(1,3,-4)$ i $\vec{n}_\beta=(2,2,2).$ Ozna\v cimo ugao izme\dj u ovih ravni $\phi=\measuredangle(\vec{n}_\alpha,\vec{n}_\beta).$ Vrijedi
  \begin{align*}
     \cos\phi&=\frac{\vec{n}_\alpha\cdot\vec{n}_\beta}{|\vec{n}_\alpha||\vec{n}_\beta|}
            =\frac{(1,3,-4)\cdot(2,2,2)}{\sqrt{1^2+3^2+(-4)^2}\, \sqrt{2^2+2^2+2^2}}
            =\frac{2+6-8}{\sqrt{26}\,\sqrt{12}}=0,
  \end{align*}
  pa je $\phi=\frac{\pi}{2}.$
\end{example}

\section[P\lowercase{rava}]{Prava}

\index{prava}

\begin{figure}[!h]\centering
        \includegraphics[scale=.95]{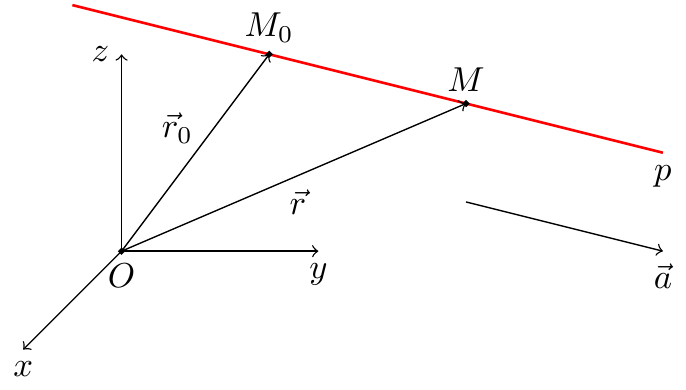}
        \caption{Prava}
     \label{slika71}
  \end{figure}

Neka kroz datu ta\v cku $M_0(x_0,y_0,z_0)$ paralelno sa vektorom $\vec{a}=(l,m,n)$ prolazi prava $p,$ Slika \ref{slika71}. Ova prava je jedinstvena.
Sa slike vidimo da je $\overrightarrow{M_0M}=\vec{r}-\vec{r}_0$ te da je $\overrightarrow{M_0M}\parallel \vec{a},$ a kako se radi o slobodnim vektorima to su $\vec{r}-\vec{r}_0$ i $\vec{a}$ linearno  zavisni, te mo\v zemo $\vec{r}-\vec{r}_0$ izraziti preko $\vec{a}$ i obrnuto, tj.
\[\vec{r}-\vec{r}_0=t\vec{a}\] odsnosno
\begin{equation}
 \tcbhighmath[mojstil1]{ \vec{r}=\vec{r}_0+t\vec{a},\,t\in\mathbb{R}.}
 \label{prava1}
\end{equation}
Jedna\v cina \eqref{prava1} je parametarska vektorska jedna\v cina prave. Vektor $\vec{a}$ je vektor pravca prave (ili karakteristi\v cni vektor prave ili vektor prave) $p.$

\index{prava! parametarska vektorska jedna\v cina}

Kako je $\vec{r}-\vec{r}_0\parallel\vec{a}$ to je $\vec{a}\times(\vec{r}-\vec{r}_0)=\vec{0},$ a ako  ozna\v cimo $\vec{b}=\vec{a}\times\vec{r}_0,$ dobijamo jednakost
\begin{equation}
   \tcbhighmath[mojstil1]{\vec{a}\times\vec{r}=\vec{b},}
 \label{prava2}
\end{equation}
koja predstavlja op\v stu vektorsku jedna\v cinu prave. U posljednoj jedna\v cini prave \eqref{prava2} vektori $\vec{a}$ i $\vec{b}$ su dati, dok je $\vec{r}$ vektor polo\v zaja proizvoljne ta\v cke prave, u ovom slu\v caju ta\v cke $M.$

\index{prava! op\v sta vektorska jedna\v cina prave}

Neka su $x,y,z$ koordinate ta\v cke $M\in p,$ a $x_0,y_0,z_0$ su koordinate ta\v cke $M_0$ iz jednakosti \eqref{prava1} dobijamo
\begin{equation}
\tcbhighmath[mojstil1]{
\begin{aligned}
      x&=x_0+lt\\
      y&=y_0+mt\\
      z&=z_0+nt,\,t\in\mathbb{R}.
\end{aligned}
}
\label{prava3}
\end{equation}
\ \\
Jednakosti \eqref{prava3} predstavljaju parametarske jedna\v cine prave u skalarnom obliku.   Sada iz jedna\v cina datih u  \eqref{prava3} izrazimo parametar $t$

\index{prava! parametarske jedna\v cine prave\\ u skalarnom obliku}

\begin{align*}
   x&=x_0+lt\Leftrightarrow t=\frac{x-x_0}{l}\\
   y&=y_0+mt\Leftrightarrow t=\frac{y-y_0}{m}\\
   z&=z_0+nt\Leftrightarrow t=\frac{z-z_0}{n}
\end{align*}
pa dobijamo
\begin{equation}
  \tcbhighmath[mojstil1]{  \frac{x-x_0}{l}=\frac{y-y_0}{m}=\frac{z-z_0}{n}.}
\label{prava4}
\end{equation}
\noindent Posljednje jednakosti predstavljaju kanonski oblik jedna\v cine prave.
\begin{remark}
 Specijalno. Ako prava prolazi kroz dvije poznate ta\v cke 	$A(x_1,y_1,z_1)$ i $B(x_2,y_2,z_2),$ tada na osnovu \eqref{prava4} i predstavljanja vektora kroz dvije ta\v cke \eqref{vektorDvijeTacke}, dobijamo jedna\v cinu prave kroz dvije ta\v cke  \[ \frac{x-x_1}{x_2-x_1}=\frac{y-y_1}{y_2-y_1}=\frac{z-z_1}{z_2-z_1}.\]
\end{remark}

\index{prava! kanonske jedna\v cine prave}

\begin{figure}[!h]\centering
	\includegraphics[scale=.8]{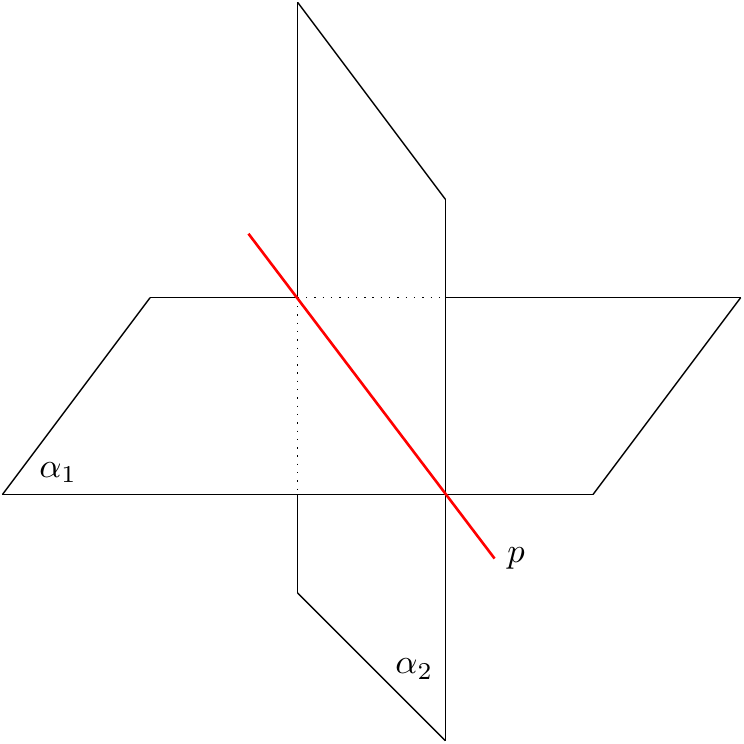}
	\caption{Prava $p$ u presjeku dvije ravni $\alpha_1$ i $\alpha_2$}
	\label{slika72}
\end{figure}
\newpage

Pravu mo\v zemo predstaviti i kao dvije ravni koje se sijeku. Neka su date ravni
\begin{align*}
\alpha_1&:A_1x+B_1y+C_1z+D_1=0\\
\alpha_2&:A_2x+B_2y+C_2z+D_2=0,
\end{align*}
koje se sijeku, tada
\begin{equation}
\tcbhighmath[mojstil1]{ p:\left\{ \begin{array}{c}A_1x+B_1y+C_1z+D_1=0\\A_2x+B_2y+C_2z+D_2=0,\end{array}\right.}
\label{prava5}
\end{equation}
\ \\
predstavlja jedna\v cinu prave Slika \ref{slika72}.

\index{prava! jedna\v cina prave}

\begin{example}
Pravu $p:\begin{cases}2x-y-z-4=0\\2x-3y-2z+7=0 \end{cases}$ napisati u \\
\begin{inparaenum}[$(a)$]
   \item kanonskom;
   \item parametarskom obliku.
 \end{inparaenum}\ \\\\
\noindent Rje\v senje:\\\\
Jedna\v cina prave u kanonskom obliku je
$\frac{x-x_0}{l}=\frac{y-y_0}{m}=\frac{z-z_0}{n}.$ Potrebno je da odredimo vektor pravca $\vec{a}=(l,m,n)$ i ta\v cku $M(x_0,y_0.z_0)$ koja pripada pravoj.   Odredimo prvo vektor pravca $\vec{a}$ prave kao $\vec{a}=\vec{n}_\alpha\times\vec{n}_\beta,$ gdje su $\alpha:2x-y-z-4=0$ i $\beta:2x-3y-2z+7=0.$ Po\v sto je $\vec{n}_\alpha=(2,-1,-1)$ $\vec{n}_\beta=(2,-3,-2),$ vrijedi
  \begin{align*}
    \vec{a}&=\vec{n}_\alpha\times\vec{n}_\beta
    =\begin{blockarray}{ccc}\begin{block}{|ccc|}\vec{i}&\vec{j}&\vec{k}\\2&-1&-1\\2&-3&-2\\\end{block}\end{blockarray}=
    -\vec{i}+2\vec{j}-4\vec{k}.
  \end{align*}
Ta\v cku $M$ dobi\' cemo rje\v savaju\' ci sistem
$\begin{cases}2x-y-z-4=0\\2x-3y-2z+7=0. \end{cases}$ Sistem je neodre\dj en i ima beskona\v cno mnogo rje\v senja. Svako rje\v senje ovog sistema predstavlja jednu ta\v cku prave $p.$ Uvrstimo na primjer $x=0$ i rije\v simo sistem, dobijamo $y=15,$ $z=-19.$ Pa je tra\v zena jedna\v cina prave u kanonskom obliku
\[p:\frac{x-0}{-1}=\frac{y-15}{2}=\frac{z+19}{-4}.\]
Iz posljednje jednakosti je $\frac{x-0}{-1}=\frac{y-15}{2}=\frac{z+19}{-4}=t,$ pa su parametarske jedna\v cina prave
\begin{align*}
  \frac{x-0}{-1}&=t\Leftrightarrow x=-t\\
  \frac{y-15}{2}&=t\Leftrightarrow y=2t+15\\
  \frac{z+19}{-4}&=t\Leftrightarrow z=-4t-19.
\end{align*}
\end{example}

\subsection{Uzajamni odnos dvije prave}

Neka su date dvije prave $p_1$ i $p_2$ u kanonskim oblicima
\begin{align*}
p_1&:\frac{x-x_1}{l_1}=\frac{y-y_1}{m_1}=\frac{z-z_1}{n_1},\\
p_2&:\frac{x-x_2}{l_2}=\frac{y-y_2}{m_2}=\frac{z-z_2}{n_2},
\end{align*}
gdje su $\vec{a}_1=(l_1,m_1,n_1)$ i $\vec{a}_2=(l_2,m_2,n_2)$ vektori pravca pravih $p_1$ i $p_2,$ respektivno i znamo da je $M_1(x_1,y_1,z_1)\in p_1$ i $M_2(x_2,y_2,z_2)\in p_2$.

\begin{figure}[!h]\centering
        \includegraphics[scale=.9]{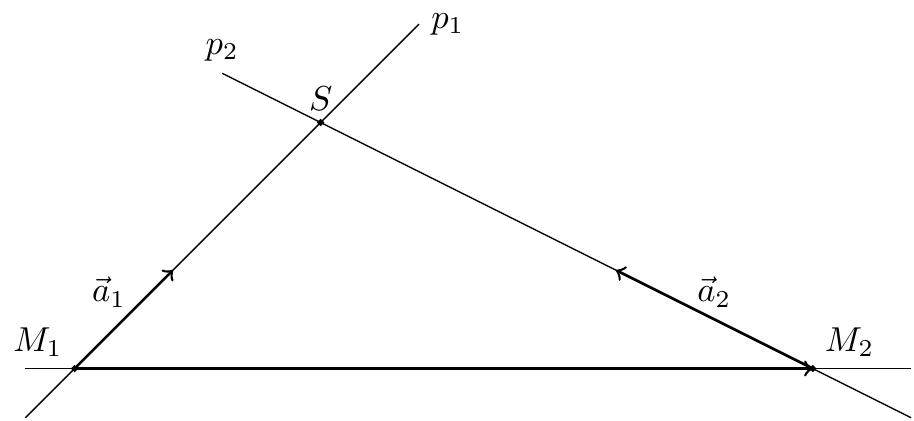}
        \caption{ Prave $p_1$ i $p_2$}
     \label{slika73}
  \end{figure}

\paragraph{Ugao izme\dj u dvije prave.} Pod uglom $\phi$ izme\dj u dvije prave smatramo ugao izme\dj u njihovih vektora pravaca, u na\v sem slu\v caju
$\phi=\measuredangle(\vec{a}_1,\vec{a}_2),$ pa  je
\begin{equation*}
  \tcbhighmath[mojstil1]{ \cos\phi=\frac{\vec{a}_1\cdot\vec{a}_2}{|\vec{a}_1||\vec{a}_2|}=
      \frac{l_1l_2+m_1m_2+n_1n_2}{\sqrt{l^2_1+m^2_1+n^2_1}\sqrt{l^2_2+m^2_2+n^2_2}}.}
\label{prava6}
\end{equation*}
\index{prava! ugao izme\dj u dvije prave}
U slu\v caju $\phi=\frac{\pi}{2},$ vrijedi  $p_1\perp p_2\Leftrightarrow\vec{a}_1\perp\vec{a}_2,$ pa je
\begin{equation*}
 \tcbhighmath[mojstil1]{l_1l_2+m_1m_2+n_1n_2=0.}
\label{prava7}
\end{equation*}
\index{prava! uslov okomitosti}
Za $\phi=0$ sada je $p_1\parallel p_2\Leftrightarrow\vec{a}_1\parallel\vec{a}_2,$ pa je
\begin{equation*}
  \tcbhighmath[mojstil1]{\frac{l_1}{l_2}=\frac{m_1}{m_2}=\frac{n_1}{n_2}.}
\label{prava8}
\end{equation*}
\index{prava! uslov paralelnosti}
\paragraph{Uslov komplanarnosti i uslov da se dvije prave sijeku.} Posmatrajmo Sliku \ref{slika73}. Ako su vektori $\vec{a}_1,\,\vec{a}_2$ i $\overrightarrow{M_1M_2}$ komplanarni,
tada je
\begin{equation*}
\tcbhighmath[mojstil1]{\begin{blockarray}{ccc}\begin{block}{|ccc|}x_2-x_1&y_2-y_1&z_2-z_1\\l_1&m_1&n_1\\l_2&m_2&n_2\\\end{block}\end{blockarray}=0,}
\label{prava9}
\end{equation*}
\newpage
\noindent vrijedi i obrnuto. Sada mogu nastupiti dva slu\v caja, da se prave sijeku ili da su paralelne, tj.

\begin{figure}[!h]\centering
        \includegraphics[scale=1]{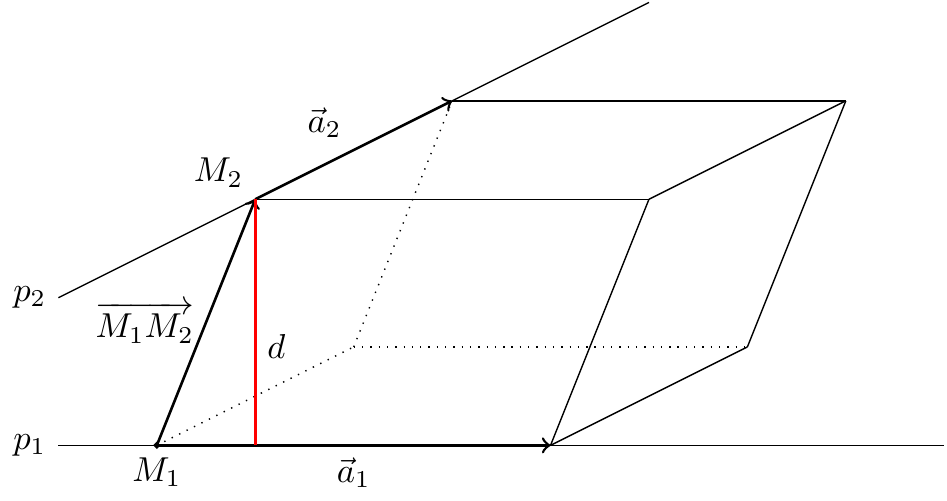}
        \caption{ Mimoilazne prave $p_1$ i $p_2$ }
     \label{slika74a}
  \end{figure}

\begin{enumerate}
\item Ako je ispunjen uslov \eqref{prava9} i ako su elementi druge vrste proporcionalni elementima tre\' ce vrste, tada su prava paralelne $p_1\parallel p_2.$
\item Ako je ispunjen uslov \eqref{prava9}, ali ako elementi druge vrste nisu proporcionalni elementima tre\' ce vrste tada se prave sijeku.
\end{enumerate}
\index{prava! uslov komplanarnosti}
\index{prava! uslov da se dvije prave sijeku}
\paragraph{Rastojanje izme\dj u dvije mimoilazne prave.} Mimoilazne prave su prave koje nisu paralelne i nemaju zajedni\v ckih ta\v caka. Zapreminu paralelopipeda, koji je odre\dj en vekto-\\rima $\vec{a}_1,$ $\vec{a}_2$ i $\overrightarrow{M_1M_2}$ (vidjeti Sliku \ref{slika74a}),  mo\v zemo izra\v cunati na sljede\' ci na\v cin
\[V=\left| (\vec{a}_1\times\vec{a}_2)\cdot\overrightarrow{M_1M_2}\right|\]
kao i
\[V=\left| \vec{a}_1\times\vec{a}_2\right|d,\]

\begin{figure}[!h]\centering
	\includegraphics[scale=1]{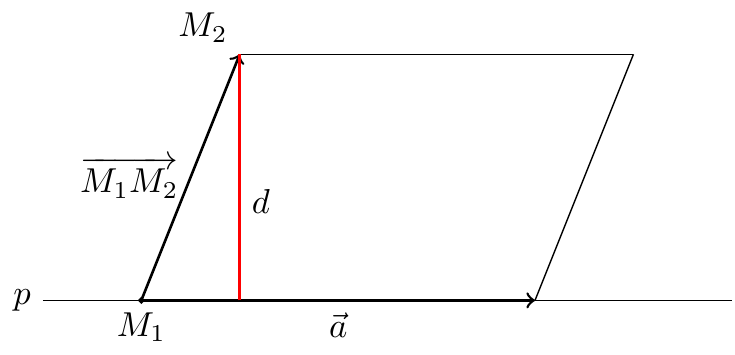}
	\caption{Prava $p$ i  ta\v cka $M_2$ }
	\label{tackaPrava}
\end{figure}

pa je sada
\begin{equation*}
 \tcbhighmath[mojstil1]{d=\frac{\left|(\vec{a}_1\times\vec{a}_2)\cdot\overrightarrow{M_1M_2} \right|}
                  {\left| \vec{a}_1\times\vec{a}_2\right|}.}
\label{prava10}
\end{equation*}
\index{prava! rastojanje izme\dj u dvije mimoilazne prave}
\paragraph{Rastojanje ta\v cke od prave.}  \index{prava! rastojanje ta\v cke od prave} Neka je potrebno odrediti rastojanje ta\v cke, \v cije koordi-\\nate znamo, $M_2(x_2,y_2,z_2)$  od prave $p:\frac{x-x_1}{l}=\frac{y-y_1}{m}=\frac{z-z_1}{n},$ Slika \ref{tackaPrava}. Neka je $M_1(x_1,y_1,z_1)$ fiksirana ta\v cka na pravoj $p.$ Po\v cetak vektora pravca $\vec{a}=(l,m,n),$ prave $p,$  dovedemo u ta\v cku $M_1.$

Vektori $\vec{a}$ i $\overrightarrow{M_1M_2}$ konstrui\v su paralelogram, \v cija je povr\v sina
\[P=\left|  \vec{a} \times \overrightarrow{M_1M_2}\right|.\]
\newpage
\noindent Istu povr\v sinu mo\v zemo izra\v cunati kao proizvod intenziteta vektora $\vec{a}$ i visine $d,$
\[P=|\vec{a}| d.\]
Kombinuju\'ci posljednje dvije jednakosti dobijamo formulu za ra\v cunanje rastojanja ta\v cke od prave
\begin{equation*}
\tcbhighmath[mojstil1]
{d=\frac{ \left|   \vec{a} \times \overrightarrow{M_1M_2}  \right|   }{\left|\vec{a}  \right|}.
}
\label{pravaTacka10}
\end{equation*}

\subsection{Uzajamni odnos prave i ravni}
Neka su dati ravan $\alpha$ i prava $p$ sa
\begin{align*}
  \alpha:Ax+By+Cz+D=0,\\
  p:\frac{x-x_0}{l}=\frac{y-y_0}{m}=\frac{z-z_0}{n},
\end{align*}
gdje je $\vec{n}=(A,B,C)$ vektor normale ravni $\alpha,$ $\vec{a}=(l,m,n)$ vektor pravca prave $p$ i ta\v cka $M_0(x_0,y_0.z_0)\in p.$

 \begin{figure}[!h]\centering
       \begin{subfigure}[b]{.3\textwidth}\centering
        \includegraphics[scale=.65]{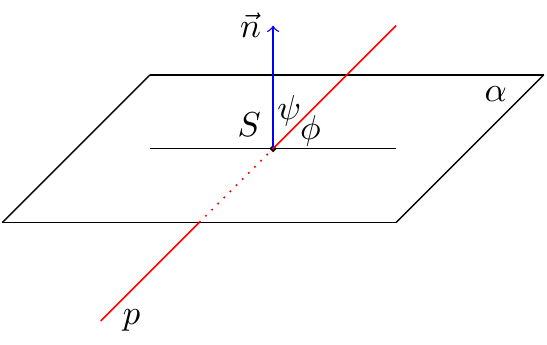}
         \caption{Prodor prave}
          \label{slika75}
   \end{subfigure}
   \begin{subfigure}[b]{.3\textwidth}\centering
        \includegraphics[scale=.65]{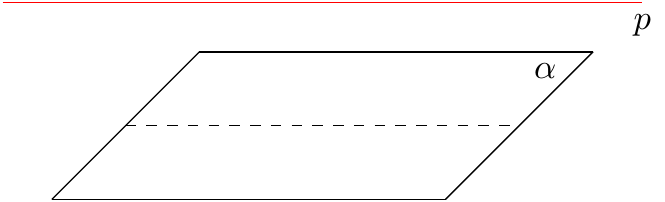}\vspace{1.cm}
         \caption{Prava i ravan su paralelne}
      \label{slika76}
     \end{subfigure}
    \begin{subfigure}[b]{.3\textwidth}\centering\hspace{.5cm}
        \includegraphics[scale=.65]{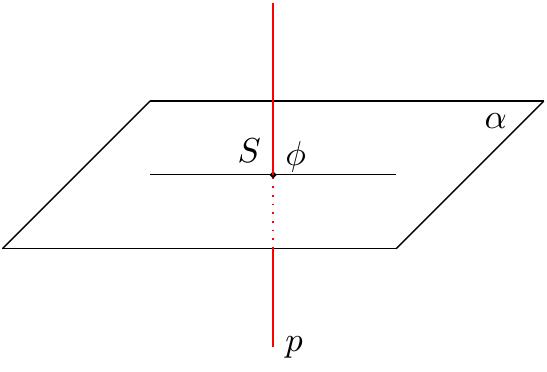}
         \caption{Prava i ravan su okomite $\phi=\frac{\pi}{2}$ }
      \label{slika77}
     \end{subfigure}
       \caption{ Odnos prave i ravni}
     \label{slika74}
  \end{figure}

\index{prava i ravan! uzajamni odnos prave i ravni}

Ugao $\phi$ (u na\v sem slu\v caju Slika \ref{slika75}), prave $p$ prema ravni $\alpha,$ ra\v cunamo kao \[\phi=\frac{\pi}{2}-\psi.\]
\newpage
\noindent Kako je $\cos\psi=\cos\left(\frac{\pi}{2}-\phi\right)=\sin\phi,$ to vrijedi
\begin{equation*}
\tcbhighmath[mojstil1]{ \sin\phi=\frac{Al+Bm+Cn}{\sqrt{A^2+B^2+C^2}\sqrt{l^2+m^2+n^2}}.}
\label{prava11}
\end{equation*}
Ako su prava i ravan paralelne Slika \ref{slika76}, vrijedi
\[p\parallel\alpha\Leftrightarrow\vec{a}\perp\vec{n}\Leftrightarrow\tcbhighmath[mojstil1]{Al+Bm+Cn=0.}\]
Ako su prava i ravan normalne Slika \ref{slika77}, vrijedi
\[ p\perp\alpha\Leftrightarrow\vec{a}\parallel\vec{n}\Leftrightarrow
   \tcbhighmath[mojstil1]{\frac{l}{A}=\frac{m}{B}=\frac{n}{C}.}\]
Da bi prava $p\subset\alpha$ potrebno je i dovoljno da vrijedi $M_0\in\alpha\wedge p\parallel\alpha,$ tj.
\begin{align*}
  Ax_0+By_0+Cz_0+D&=0\\
  Al+Bm+Cn&=0.
\end{align*}
Ta\v cku prodora $S$  prave $p$ kroz ravan $\alpha,$ Slika \ref{slika75}, najlak\v se je odrediti tako \v sto pravu izrazimo preko parametarskih jedna\v cina \eqref{prava3}. Zatim promjenljive $x,y,z$ u jedna\v cini ravni izrazimo preko jedna\v cina \eqref{prava3}. Dobijamo
\[Ax_0+By_0+Cz_0+D+t(Al+Bm+Cn)=0,\]
odakle se ra\v cuna $t.$

\begin{example}
 Odrediti odnos datih pravih i ravni
 \begin{enumerate}[(a)]
 	\item $p:\frac{x-1}{2}=\frac{y}{3}=\frac{z+1}{-1},$ $\alpha:x+y+5z-7=0;$
 	\item $p:\frac{x-2}{3}=\frac{y-1}{-2}=\frac{z-3}{2},$ $\alpha:2x+2y-z-3=0;$
 	\item $p:\frac{x-1}{3}=\frac{y-2}{-2}=\frac{z-3}{1},$ $\alpha:6x-4y+2z+7=0.$ 	
 \end{enumerate}	
\noindent Rje\v senje:
 \begin{enumerate}[(a)]
 	\item U ovom slu\v caju vektor pravca prave je $\vec{a}=(l,m,n)=(2, 3, -1),$ dok je vektor normale ravni $\vec{n}=(A,B,C)=(1,1,5).$  Izra\v cunajmo $\vec{n}\cdot\vec{a},$ vrijedi
 	\[\vec{n}\cdot\vec{a}=Al+Bm+Cn=1\cdot 2+1\cdot 3 +5\cdot (-1)=0,\]
 	dakle vektori $\vec{n}$ i $\vec{a}$ su ortogonalni, pa je ili $\alpha\parallel p$ ili $p\subset \alpha.$  Ta\v cka $M_0(1,0,-1)\in p,$ ispitajmo da li ova ta\v cka pripada i ravni $\alpha.$ Vrijedi
 	\[Ax_0+By_0+Cz_0+D=1\cdot1+1\cdot 0+5\cdot(-1)-7=-11\neq 0, \]
 	ta\v cka $M_0\notin \alpha$ (ta\v cka $M_0$ ne zadovoljava jedna\v cinu ravni), pa zaklju\v cujemo da vrijedi $p\parallel\alpha,$ tj. date prava i ravan su paralelne.
 	\item Sada je $\vec{a}=(3,-2,2)$ i $\vec{n}=(2,2,-1),$ pa vrijedi
 	\[\vec{n}\cdot \vec{a}=2\cdot3+2\cdot (-2)-1\cdot 2=0,\]
 	opet je ili $\alpha\parallel p$ ili $p\subset \alpha.$  	Ispitajmo da li ta\v cka $M_0(2,1,3)$ pripada ravni $\alpha,$ vrijedi
 	\[Ax_0+By_0+Cz_0+D=2\cdot 2+2\cdot 1-1\cdot3-3=0,\]
 	dakle $M_0\in\alpha$ i $p\subset\alpha,$ tj. data prava $p$ le\v zi u datoj ravni $\alpha.$
 	\item Sada je $\vec{a}=(3,-2,1)$ i $\vec{n}=(6,-4,2),$ pa je
 	\[\vec{n}\cdot \vec{a}=6\cdot 3-4\cdot (-2)+2\cdot 1=28\neq 0, \]
 	dakle vektori $\vec{a}$ i $\vec{n}$ nisu ortogonalni pa imamo prodor prave $p$ kroz ravan $\alpha.$ Sada mo\v zemo izra\v cunati koordinate prodorne ta\v cke $P$. Da bi to uradili treba da rije\v simo sistem  $\frac{x-1}{3}=\frac{y-2}{-2}=\frac{z-3}{1} \wedge 6x-4y+2z+7=0.$ Ovo je najlak\v se uraditi na sljede\'ci na\v cin:
 	prvo prevedemo pravu $p$ u parametraski oblik.
 	Vrijedi $\frac{x-1}{3}=\frac{y-2}{-2}=\frac{z-3}{1}=t,$ pa je
 	\begin{align*}
 	\frac{x-1}{3}&=t\Leftrightarrow  x=3t+1 \\
    \frac{y-2}{-2}&=t\Leftrightarrow y=-2t+2 \\
    \frac{z-3}{1}&=t\Leftrightarrow z=t+3.
 	\end{align*}
 	Sada u jedna\v cini ravni $\alpha$ teku\'ce koordinate $x,\, y,\,z$ izrazimo preko parametra $t$ koriste\'ci upravo dobijene parametarske jedna\v cine prave, vrijedi
 	\[6(3t+1)-4(-2t+2)+2(t+3)+7=0\Leftrightarrow t=-\frac{11}{28}.\]
 	Dobijenu vrijednost parametra $t$ uvrstimo u parametarske jedna\v cine prave, dobijamo
 	\begin{align*}
 	x&=3t+1=-\frac{5}{28}\\
 	y&=-2t+2=\frac{78}{28}\\
 	z&=t+3=\frac{73}{28},
 	\end{align*}
 	ili $P(-\frac{5}{28}, \frac{78}{28}, \frac{73}{28}).$\\
 	Ugao $\phi$ izmedju prave i ravni ra\v cunamo iz $\sin\phi=\frac{Al+Bm+Cn}{\sqrt{A^2+B^2+C^2}\sqrt{l^2+m^2+n^2}}.$ Vrijedi
 	\[\sin\phi=\frac{Al+Bm+Cn}{\sqrt{A^2+B^2+C^2}\sqrt{l^2+m^2+n^2}}=\frac{6\cdot 3-4\cdot (-2)+2\cdot 1}{\sqrt{6^2+(-4)^2+2^2} \sqrt{3^2+(-2)^2+1^2}  }=1, \]
 	pa je $\phi=\frac{\pi}{2}.$ Do istog rezultata mogli smo do\'ci na osnovu uslova kolinearnosti dva vektora, vrijedi
 	\[\frac{l}{A}=\frac{m}{B}=\frac{n}{C}=\frac{3}{6}=\frac{-2}{-4}=\frac{1}{2}\quad \left( =\frac{1}{2}\right),\]
 	dakle vektor pravca $\vec{a}$ prave $p$ i vektor normale $\vec{n}$ ravni $\alpha$ su kolinearni, pa je $p\perp \alpha,$ tj. $\phi=\frac{\pi}{2}$ (ili $90^0$). 	
 \end{enumerate}	
\end{example}

\section{Z\lowercase{adaci}}

\index{Zadaci za vje\v zbu! ravan i prava}

\begin{enumerate}
\item Izra\v cunati ugao izmedju ravni $\alpha:-x+3y-z-4=0$ i $\beta:-3x+6z-6=0. $
\item
  \begin{enumerate}
  	\item Odrediti jedna\v cinu ravni $\alpha,$ koja prolazi kroz ta\v cke $A(2,-1,0)$ i $B(3,2-5),$ a normalna je na ravan $\beta:2x-y+3z-7=0.$
    \item Napisati jedna\v cinu ravni koja prolazi kroz ta\v cku $M(-2,3,-1)$ i
      \begin{enumerate}[(i)]
         \item na koordinatnim osama odsijeca jednake odsje\v cke;
         \item prolazi  kroz $y-$osu;
         \item prolazi kroz koordinatni po\v cetak i ta\v cku $A(2,1-5).$
      \end{enumerate}
    \item Izra\v cunati jedna\v cinu ravni u obliku $\vec{r}\cdot \vec{n}=\alpha,$ koja je paralelna sa ravni $\vec{r}\cdot(3\vec{i}-\vec{j}+\vec{k})=-5$
       i prolazi kroz ta\v cku $M(0,1,2).$
    \item Napisati jedna\v cinu ravni koja sadr\v zi ta\v cke $M(2,1,-1),\,N(-1,0,1)$ a okomita je na ravan $2x-y+4z-1=0.$
    \item Napisati jedna\v cinu ravni koja sadr\v zi ta\v cke $M(-1,1,2),\,N(0,2,1)$ a okomita je na ravan  $x-3y+4z-7=0$.
\end{enumerate}

\item Izra\v cunati visinu piramide $h_S$, \v ciji su vrhovi u ta\v ckama $S(1,-2,3),\,A(2,-4,2),$\\ $B(2,3,4),\,C(1,2,3).$

\item Pravu $l:\left\{ \begin{array}{c} -2x-y+3z-4=0\\x+2y-z=1, \end{array}\right.$  napisati u
    \begin{enumerate}
      \item kanonskom obliku;
      \item parametarskom obliku.
    \end{enumerate}

\item
  \begin{enumerate}
  	\item Odrediti jedna\v cinu ravni $\alpha$ koja sadr\v zi pravu $p:\dfrac{2-x}{3}=\dfrac{y-1}{2}=\dfrac{-2z}{3}$ i normalna je ravan $\alpha:3x-4y+2z-1=0.$
    \item Izra\v cunati jedna\v cinu ravni kojoj pripada prava $p:\dfrac{x-2}{1}=\dfrac{y+2}{1}=\dfrac{z}{2},$ a paralelna je
      pravoj $q:\dfrac{x-2}{-1}=\dfrac{y}{0}=\dfrac{z+3}{1}.$
    \item Napisati jedna\v cinu ravni koja sadr\v zi prave $p:\dfrac{x-1}{2}=\dfrac{y+1}{-1}=\dfrac{z}{1}$ i $q:\dfrac{x+1}{-4}=\dfrac{y}{2}=\dfrac{z+1}{-2}.$
    \item Napisati jedna\v cinu normale povu\v cenu iz ta\v cke $M(1,1,2)$ na pravu $p:\dfrac{x-1}{3}=\dfrac{y}{-1}=\dfrac{z+2}{2}.$
    \item Napisati jedna\v cinu ravni koja sadr\v zi pravu $p:\dfrac{x-1}{2}=\dfrac{y}{1}=\dfrac{z+1}{3}$ i prolazi kroz ta\v cku $M(2,1-,0).$
    \item Izra\v cunati jedna\v cinu ravni koja sadr\v zi ta\v cu $M(1,0,1)$\\ i pravu $\left\{ \begin{array}{l}2x+3y-z-5=0\\x-3y+2z+2=0.\end{array} \right.$
   \end{enumerate}

\item Odrediti rastojanje ta\v cke $A(2,1,0)$ od prave   $l:\left\{ \begin{array}{c} -x-y+3z-4=0\\3x+2y-z=1. \end{array}\right.$

\item Data je jedna\v cina prave u vektorskom obliku $\vec{r}\times (2\vec{i}+\vec{j}-\vec{k})=2\vec{i}-\vec{j}+2\vec{k}.$ Odrediti odgovaraju\' ci kanonski
      i parametarski oblik jedna\v cine prave.
\item Odrediti najkra\' ce rastojanje izmedju pravih  $p:\dfrac{x-1}{2}=\dfrac{y+2}{1}=\dfrac{z-5}{-1}$ i $q:\dfrac{x+3}{1}=\dfrac{y-3}{2}=\dfrac{z}{-3}.$

\item Odrediti jedna\v cinu normale povu\v cene iz ta\v cke $P(1,3,2)$ na pravu $l:\dfrac{x-2}{3}=\dfrac{y}{-1}=	\dfrac{z+1}{2}.$
\item Odrediti jedna\v cinu ortogonalne projekcije prave $\dfrac{x-2}{2}=\dfrac{y-1}{4}=\dfrac{z+1}{3}$ na ravan $x+2y+z-3=0.$
\end{enumerate}

\chapter{Dodatak}

\section[P\lowercase{rocentni ra\v cun.} R\lowercase{a\v cun smjese.}]{Procentni ra\v cun. Ra\v cun smjese.}
\pagestyle{fancy}


\subsection{Razmjere i proporcije} \index{razmjera ili omjer}
U raznim prora\v cunima, u prirodnim, in\v zinjerskim  naukama i dr., \v cesto se pojavljuju razmjere, odnosi ili omjeri nekih veli\v cina.
\begin{definition}[Razmjera] Razmjera (ili omjer) dva broja $a$ i $b$ ($a,\,b>0$) je njihov koli\v cnik. Pi\v semo \[a:b\text{ ili } \frac{a}{b},\:a,b> 0.\]
\end{definition}
Ako su dvije razmjere $a:b$ i $c:d$ jednake, dobijamo proporciju.
\begin{definition}[Prosta proporcija] \index{proporcija!prosta}
Jednakost koju formiraju jednake razmjere nazivamo proporcija. Pi\v semo
\[a:b=c:d\text{ ili } \frac{a}{b}=\frac{c}{d}.\]
\end{definition}
\noindent Proporcije imaju osobinu
\begin{empheq}[box=\mymath]{equation*}
a:b=c:d=(a\pm c):(b\pm d).
\end{empheq}

\noindent Ako je vi\v se razmjera jednako, vrijedi
\begin{equation*}\label{prop}
\tcbhighmath[mojstil1]{a:b=c:d=e:f=p:q\Leftrightarrow \frac{a}{b}=\frac{c}{d}=\frac{e}{f}=\frac{p}{q}\Leftrightarrow a:c:e:p=b:d:f:q.}
\end{equation*}
Proporcije u  prethodnom izrazu nazivamo produ\v zene proporcije.\index{proporcija!produ\v zena}

\begin{example}
    Odrediti $x$ iz proporcije $(x+9):6=x:5.$ \\\\
\noindent Rje\v senje:\\\\
  $(x+9):6=x:5\Leftrightarrow 5(x+9)=6x\Leftrightarrow x=45.$
\end{example}

\begin{example}
       Primjenom osobina produ\v zene proporcije odrediti $x,y,z$ i $t,$ ako je\\ $x+y+z+t=198$ i $x:y:z:t=1:2:3:5.$\\\\
\noindent Rje\v senje:\\\\  Iz produ\v zene proporcije vrijedi
          \[\frac{x}{1}=\frac{y}{2}=\frac{z}{3}=\frac{t}{5}=k,\]
          pa je
          \[x=k,\,y=2k,\,z=3k,\,t=5k.\]
          Sada u jednakosti $x+y+z+t=198,$ zamijenimo $x,y,z$ i $t$ preko $k,$ dobijamo
          \[k+2k+3k+5k=198\Leftrightarrow 11k=198\Leftrightarrow k=18,\]
          pa je na kraju
          \[x=18,\,y=36,\,z=54,\,t=90.\]
\end{example}

\begin{example}[Direktna proporcija.] \index{proporcija!direktna}
Tri majice ko\v staju $90\km.$ Koliko majica mo\v zemo kupiti za $240\km$?\\\\
\noindent Rje\v senje:\\\\ Vrijedi $3:90\km=x:240\km\Leftrightarrow x\cdot 90\km=3\cdot240\km\Leftrightarrow x=8.$
\end{example}
%
\begin{example}[Obrnuta proporcija.]\index{proporcija!produ\v zena}
 Jedan vagon vre\' ca krompira istovarila su tri radnika za 12 sati. Za koliko \' ce sati taj vagon istovariti 4 radnika?\\\\
\noindent Rje\v senje: \\\\ Vrijedi $3:4=x\,\text{sati}:12\,\text{sati}\Leftrightarrow 4x=36\Leftrightarrow x=9\,\text{ sati.}$
\end{example}

\subsection{Procentni ra\v cun}
\index{procentni ra\v cun}
Ako ozna\v cimo sa
 \begin{align*}
    G-&\text{ glavnicu}\\
    p-&\text{ procenat}\\
    I-&\text{ iznos,}
 \end{align*}
tada vrijedi proporcija
\begin{equation*}
\tcbhighmath[mojstil1]{ G:100=I:p.}
\label{procenat}
\end{equation*}
Glavnica je ukupan iznos ili koli\v cina nekog dobra, iznos je dio glavnice, dok je procenat stoti dio glavnice. Oznaka za $1$ procenat je
\[1 \text{ procenat} =1\%.\]
Dakle $1\%=\frac{1}{100} \text{ glavnice}.$ Glavnici odgovara $100\%,$ a iznosu $I$ vrijednost $p.$


\begin{example}
  U grupi je bilo 32 studenta, a sljede\' cu godinu upisalo je 30--oro. Koliki je procenat prolaznosti? \\\\
\noindent Rje\v senje:\\\\ Vrijedi $G=32,\,I=30,$ pa je
\[32:100\%=30:p\Leftrightarrow 32p =3000\%\Leftrightarrow p=93.75\%.\]
Prolaznost je $93.75\%.$
\end{example}

\begin{example}
  Cijena nekog proizvoda je smanjena za 10\%, a zatim je pove\' cana za 15\% i sada iznosi 60 KM. Kolika je prvobitna cijena? \\\\
\noindent  Rje\v senje:\\\\ Koristi\' cemo oznaku $C$ za cijenu i to $C_1,\,C_2,\,C_3$ za prvu, drugu i tre\' cu cijenu, respektivno. Vrijedi proporcija
\[C_1:100\%=C_2:(100\%-10\%)\Leftrightarrow C_2=\frac{9}{10}C_1,\]
ili
\[C_2=C_1-\frac{1}{10}C_1=\frac{9}{10}C_1,\:\left(10\% \text{ odgovara } \frac{1}{10}\right).\]
Dalje imamo
\[C_3=C_2+0.15C_2=1.15 C_2=1.15\cdot 0.9 C_1\Leftrightarrow C_1=\frac{C_3}{1.15\cdot 0.9}=\frac{60}{1.15\cdot 0.9}\approx57.97\km.  \]
\end{example}
\begin{remark}
Osim procenta koristi se i promil. Promil je hiljaditi dio od neke cjeline, oznaka za 1 promil je
   \[1\text{ promil }=1\text{\textperthousand}. \]
\end{remark}

\subsection{Ra\v cun smjese}
\index{ra\v cun smjese}
\paragraph{Jednostavni (ili prosti) ra\v cun smjese} Ponekad je potrebno pomije\v sati  proizvode razli\v citih cijena kako bi se dobio proizvod neke zadane cijene. Na primjer, trgovac ima u skladi\v stu kafu od $3\km$ i kafu od $15\km.$ Prva kafa je lo\v sijeg kvaliteta, a druga je preskupa za maloprodaju. On se odlu\v cuje da ih pomije\v sa da bi dobio kafu unaprijed zadane mase, koja ko\v sta $8\km$ po kilogramu, koja ne\' ce biti ni preskupa, a bi\' ce zadovoljavaju\' ceg kvaliteta.

Sa matemati\v cke strane, isti problem je u npr. mije\v sanju dvije kiseline razli\v citih koncentracija da bi se dobila kiselina sa unaprijed zadanom koncetracijom i unaprijed zadane zapremine.  Preciznija formulacija problema glasi:\\
U kojem odnosu i u kojim koli\v cinama treba pomije\v sati ta\v cno dvije veli\v cine (dva sastojka) $x_1$ i $x_2$ koje imaju neko zajedni\v cko svostvo razli\v citih vrijednosti (intenziteta) $s_1$ i $s_2,$ tako da se dobije smjesa ukupne koli\v cine $x$ i \v zeljenog intenziteta $s?$
Problem rje\v savamo rje\v savanjem sistema
 \begin{gather*}
  x_1+x_2=x\\x_1s_1+x_2s_2=xs,
\end{gather*}
sada zbog $x_2=x-x_1,$ dobijamo linearnu jedna\v cinu sa jednom nepoznatom
\[x_1s_1+(x-x_1)s_2=xs,\]
\v cijim rje\v savanjem dobijamo $x_1,$ a $x_2$ iz $x_2=x-x_1.$ Tra\v zeni odnos ra\v cunamo iz \[\frac{x_1}{x_2}.\]



\begin{example}[Prost ra\v cun smjese.] \index{ra\v cun smjese!prost}
Imamo $48\%$ i $78\%$ kiselinu. Koliko je potrebno jedne, a koliko druge nasuti u posudu da bi se dobilo $10$ litara, $60\%$--ne kiseline?\\\\
\noindent Rje\v senje:\\\\
Ozna\v cimo jednu koncentraciju sa $s_1=48\%,$ odgovaraju\' cu koli\v cinu sa $x_1,$ druga koncentracija je $ \:s_2=78\%$ i odgovaraju\' ca koli\v cina je
$x_2.$ Ukupna koli\v cina je $x_1+x_2=x=10\lit,$ dok je odgovaraju\' ca koncetracija poslije mije\v sanja $s=60\%.$\\
Sada vrijedi
\begin{gather*}
  x_1+x_2=10\lit \\ x_1s_1+x_2s_2=xs,
\end{gather*}
kako je $x_2=10\lit-x_1,$ dobijamo
\[x_1\cdot 0.48+(10-x_1)\cdot 0.78=10\cdot 0.6\Leftrightarrow x_1=6\lit .\]
Potrebno je $6\lit,\,48\%$-- kiseline i $4\lit,\,78\%$--kiseline.
\end{example}

\begin{example}[Jednostavni ra\v cun smjese] \index{ra\v cun smjese!jednostavni}
  Sa koliko postotnom kiselinom, treba pomije\v sati $6\lit,\,48\%$--ne kiseline, da bi se dobilo $10\lit,\,60\%$--ne kiseline?  \\\\
\noindent   Rje\v senje:\\\\  Sada vrijedi
 \[4p+6\cdot 0.48=10\cdot 0.6\Leftrightarrow p=0.78,\]
 ili u procentima $78\%.$
\end{example}

\paragraph{Slo\v zeni ra\v cun smjese} Ako imamo vi\v se od dvije veli\v cine koje trebamo pomije\v sati dobijamo sljede\' ci problem: U kojem odnosu i u kojim koli\v cinama treba pomije\v sati neke veli\v cine iste vrste $x_1,x_2,\ldots,x_n,$ koje imaju neko zajedni\v cko svojstvo razli\v citih vrijednosti $s_1,s_2,\ldots,s_n$ da bi dobili smjesu ukupne koli\v cine $x=x_1+x_2+\ldots+x_n$ i \v zeljenog intenziteta $s?$

Ovaj problem rje\v savamo algoritmom \v sema zvijezde koji \' cemo ilustrovati na sljede\' cem primjeru.

\begin{figure}[!h]\centering
 \includegraphics[scale=.8]{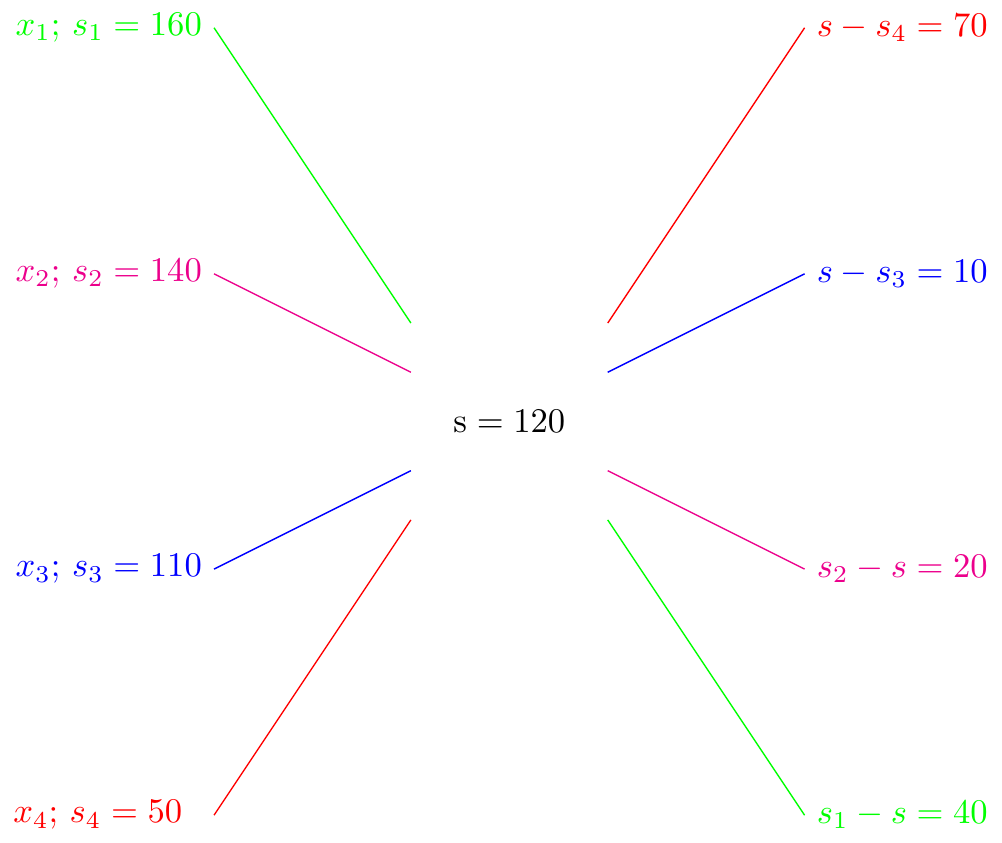}
\caption{\v Sema za dobijanje pro\v sirene proporcije}
\label{slozena}
\end{figure}

\begin{example}[Slo\v zeni ra\v cun smjese] \index{ra\v cun smjese!slo\v zeni}
U skladi\v stu su 4 vrste robe po cijenama $160,\,140,\,110$ i $50\km$ (svojstva). Kako treba izmije\v sati ovu robu po vrstama da bi dobili $560\kg$ ove robe po cijeni od $120\km$ po kilogramu? \\\\
\noindent Rje\v senje:\\\\
U \v semi su une\v sene na lijevoj strani cijene pojedinih roba (vidjeti Sliku \ref{slozena}), tj. svojstva $s_1,s_2,s_3,s_4$. Treba izra\v cunati kolike su mase pojedinih vrsta robe koje je potrebno izmije\v sati (u ovom zadatku ove mase su ozna\v cena sa $x_1,x_2,x_3,x_4$). U sredini je tra\v zena cijena, dok su na desnoj strani tra\v zene vrijednosti iz proporcije (ove vrijednosti trebaju biti pozitivne, tj. uzimaju se po apsolutnoj vrijednosti, zato uvijek oduzimamo od ve\' ce vrijednosti manju).

Iz sheme na Slici \ref{slozena} dobijemo sljede\' cu pro\v sirenu proporciju
\[x_1:x_2:x_3:x_4=70:10:20:40\Leftrightarrow x_1:x_2:x_3:x_4=7:1:2:4.\]
Iz prethodne pro\v sirene proporcije vrijedi
\begin{align*}
  \frac{x_1}{7}&=k\Leftrightarrow x_1=7k\\
  \frac{x_2}{1}&=k\Leftrightarrow x_2=k\\
  \frac{x_3}{2}&=2k\Leftrightarrow x_3=2k\\
  \frac{x_4}{4}&=4k\Leftrightarrow x_4=4k.
\end{align*}
Koristimo sada ukupnu masu robe
\[7k+k+2k+4k=560\Leftrightarrow 14k=560\Leftrightarrow k=40.\]
Tra\v zene vrijednosti dobijamo na sljede\' ci na\v cin:\\\\
Robe \v cija je cijena $s_1=160\km$ potrebno je $x_1=7\cdot 40\kg=280\kg;$\\
Robe \v cija je cijena $s_2=140\km$ potrebno je $x_2=1\cdot 40\kg=40\kg;$\\
Robe \v cija je cijena $s_3=110\km$ potrebno je $x_3=2\cdot 40\kg=80\kg;$\\
Robe \v cija je cijena $s_4=50\km$  potrebno je $x_4=4\cdot 40\kg=160\kg.$\\\\
Provjera
\[280+40+80+160=560,\]
i
\begin{align*}
280\cdot 160+40\cdot 140+80\cdot 110+160\cdot 50&=560\cdot 120\\
67200&=67200.
\end{align*}
\end{example}

\begin{remark}Ovaj zadatak nema jedinstveno rje\v senje.
\end{remark}

\section{Z\lowercase{adaci}} \index{Zadaci za vje\v zbu!ra\v cun smjese}
\paragraph{Proporcije}
\begin{enumerate}
 \item
	\begin{enumerate}
      \item Rije\v siti proporciju
         \begin{inparaenum}[(i)]
            \item $(x-3):4=10:6$;
            \item $x:\left(1+\frac{1}{11}\right)=\left(3+\frac{2}{3}\right):2.$
        \end{inparaenum}
       \item Podijeliti du\v z od $456\m$ na tri jednaka dijela \v cije \' ce du\v zine biti redom proporcionalne brojevima $\frac{2}{3},\,\frac{9}{8}$ i $\frac{7}{12}.$
       \item Unuk, otac i djed imaju zajedno 120 godina. Koliko svaki od njih ima godina, ako su im godine u razmjeri $1:5:9?$
       \item Tri osobe ulo\v zile su u jedan posao ove svote novca: osoba A $1200\km,$ osoba B $9000\km$ i osoba C $15000\km.$ Ako je ukupna zarada od tog posla $210000\km,$
        koji dio zarade \' ce pripasti osobi C, (zarada treba da je proporcionalna ulo\v zenoj svoti novca)?
       \item  Brojevi $a,b$ i $c$ su u razmjeri $2:3:4.$ Ako je njihova aritmeti\v cka sredina 15, koliko iznosi najmanji od njih?
       \item Omjer \v se\' cera i maslaca u kola\v cu je $4:3.$ U kola\v c je stavljeno $0.3\kg$ maslaca. Koliko \' cemo staviti \v se\' cera u kola\v c?
    \end{enumerate}
\end{enumerate}
\paragraph{Procentni ra\v cun}
\begin{enumerate} \setcounter{enumi}{1}
	\item
  \begin{enumerate}
    \item Sanduk vo\' ca te\v zak je $90\kg$,  od \v cega je neupotrebljivo $4\%.$ Kolika je te\v zina upotrebljivog vo\' ca?
    \item Odsutna su 4 studenta, \v sto iznosi $12.5\%$ od ukupnog broja studenata jedne grupe. Koliko studenata ima u toj grupi?
    \item Koliko litara mlije\v cne masti ima u $300\lit$ mlijeka, ako to mlijeko sadr\v zi $2.8\%$ mlije\v cne kiseline?
     \item Cijena cipela je $270\km$. Kolika \' ce cijena biti nakon sni\v zenja od $15\%$?
     \item Poslije prelaska na novo radno mjesto jednom radniku plata je pove\' cana za $20\%.$ Kolika je plata bila ako je to pove\' canje $132\km$?
     \item Cijena knjige sni\v zena je $10\%$, a zatim za $20\%$ i sada ko\v sta $45\km.$ Kolika je bila cijena prije prvog sni\v zenja?
     \item Na ispitu radila su se 3 zadatka. Pri tome $12\%$ studenata nije rije\v silo ni jedan zadatak, $32\%$ studenata rije\v silo je jedan ili dva zadatka,
        dok je $14\%$ studenata rije\v silo sva tri zadatka. Koliko je ukupno studenata radilo ispit?
     \item Poslije 3 pojeftinjenja od po $10\%$  cijena robe iznosi $2187\km.$ Izra\v cunati prvobitnu cijenu robe?
  \end{enumerate}
\end{enumerate}
\paragraph{Ra\v cun smjese}
\begin{enumerate} \setcounter{enumi}{2}
	\item
     \begin{enumerate}
     	\item U tri vre\' ce ima $64.2\kg$ bra\v sna. U prvoj vre\' ci ima $20\%$ manje bra\v sna nego u drugoj,  a u tre\' coj $42.5\%$ od koli\v cine bra\v sna
     	iz prve vre\' ce. Koliko bra\v sna ima u svakoj vre\' ci?
     	\item Na skladi\v stu ima kafe po cijeni od $7.5\km$ i $5.5\km$ po kilogramu. Napraviti $120\kg$   mje\v savine kafe koja \' ce ko\v stati po $6.8\km$ po kilogramu.
     	\item Imamo 4 vrste neke robe po cijeni od $120,\,100,\,70$ i $50\km$. Kako treba pomije\v sati tu robu da dobijemo $400\kg$ robe po cijeni od $80\km ?$
     	  Koliko mogu\' cih rje\v senja postoji?
     \end{enumerate}

  \item
     \begin{enumerate}
     	\item Koliko vode temperature $40^0C$ i vode temperature $25^0C$ treba pomije\v sati da se dobije $90\lit$ vode temperature $30^0C?$
     	\item Koliko treba uzeti sumporne kiseline ja\v cine $52\%$ , a koliko ja\v cine $88\% $ da se dobije mje\v savina od $144\lit$ ja\v cine $72\%?$
     	\item Koliko litara $80\% $ alkohola treba dodati u $1\lit$ vode da se dobije $20\%$--tni alkohol?
     	\item Razbla\v zen je $75\%$ \v spirit sa $12\lit$ vode i dobijen je $51\%$-- postotni \v spirit. Kolika je bila prvobitna koli\v cina \v spirita?
     	\item Koliko vode treba izmije\v sati sa $150\g,\:12\%$--nog rastvora kiseline da smjesa bude $4\%$--postotna?
     	\item Kada se pomije\v sa $60\lit$ ruma od $72\%$ sa $70\lit$ alkohola od $96\%.$ Koliko treba vode dosuti u ovu smjesu da se dobije rum od $46\%$?
     	\item Koliko $75\%$--tne otopine soli treba dodati u $20\lit$ da se dobije $60\%$--tna otopina soli?
   \end{enumerate}

  \item
     \begin{enumerate}
     	\item Komad  bronze mase $7.5\kg$ sadr\v zi $72.\%$ bakra. Kada se ovaj komad stopi sa drugim dobije se $10\kg$ bronze koja sadr\v zi $70\%$ bakra.
     	   Koliko je procenata bakra bilo u drugom komadu bakra?
     	\item Koliko bakra treba pomije\v sati sa $21\g$ \v cistog zlata da se dobije smjesa fino\' ce $0.75?$
     	\item Zlatar izmije\v sa dvije vrste zlata \v cije su fino\' ce $0.75$ i $0.81$ ($f$--fino\' ca zlatne smjese $0\leqslant f\leqslant 1$) i dobije se $50\g$
     	      zlata 	fino\' ce $0.774.$ Po koliko je grama zlata uzeo obje vrste zlatar?
    \end{enumerate}

\end{enumerate}

   \backmatter
   \pagestyle{empty}
   \addcontentsline{toc}{chapter}{Literatura i reference}
   \printbibliography[keyword={knjiga},title={Literatura i reference}]
  \printbibliography[keyword={samir},title={Nau\v cni rezultati -- Samir}]
  \printbibliography[keyword={sanela},title={Nau\v cni rezultati  -- Sanela}]

@article{samir2010scheme,
	author = {Duvnjakovi\'{c},E. and Karasulji\'{c},S. and Oki\v{c}i\'{c},N.}, 
	title = "{Difference Scheme for Semilinear Reaction-Diffusion Problem}",
	booktitle = {14th International Research/Expert Conference Trends in the Development of Machinery 
		and Associated Technology TMT 2010, 7. Mediterranean Cruise }, 
	year = {2010},            
	pages = {793--796},
	url={http://www.tmt.unze.ba/zbornik/TMT2010/199-TMT10-138.pdf}
}

@article{samir2011scheme,
	author = { Duvnjakovi\'c, E. and Karasulji\'c, S.},
	title = {Difference Scheme for Semilinear Reaction-Diffusion Problem on a Mesh of Bakhvalov Type}, 
	journal = {Mathematica Balkanica},
	volume = {25, Fasc. 5},
	pages ={499--504},
	year = {2011},
	url={http://www.mathbalkanica.info/toc/siteCONT25-5.pdf}
}

@article{samir2015uniformlyconvergent,
	author={Duvnjakovi{\'c}, E. and Karasulji{\'c}, S. and Pa{\v s}i{\'c}, V. and Zarin, H.},
	title={A uniformly convergent difference scheme on a modified Shishkin mesh for the singularly perturbed reaction-diffusion boundary value problem},
	year={2015},
	journal={Journal of Modern Methods in Numerical Mathematics},
	volume={6},
	number={1},
	pages={28-43},
	url={http://www.m-sciences.com/index.php?journal=jmmnm&page=issue&op=view&path[]=72},
	doi={10.20454/jmmnm.2015.971}
}

@article{samir2015uniformly,
	author={Karasulji{\'c}, S. and Duvnjakovi{\'c}, E. and  Zarin, H.},
	title={Uniformly convergent difference scheme for a semilinear reaction-diffusion problem},
	journal={Advances in Mathematics: Scientific Journal},
	volume={4},
	number={2},
	year={2015},
	pages={139--159},
	url={http://research-publication.com/wp-content/uploads/2019/03/AMSJ-2015-N2-6.pdf}
}

@article{samir2017construction,
	author={Karasulji{\' c}, S. and Duvnjakovi{\' c}, E. and Pa{\v s}i{\' c}, V. and Barakovi{\' c}, E.},
	title={Construction of a global solution for the one dimensional singularly--perturbed boundary value problem}, 
	journal={Journal of Modern Methods in Numerical Mathematics},
	volume={8},
	number={1--2},
	year={2017},
	doi={10.20454/jmmnm.2017.1275 },
	pages={52-65}
}

@article{samir2018uniformly,
	author={Karasulji{\' c}, S. and Duvnjakovi{\' c}, E. and Memi{\' c}, E.},
	title={Uniformly Convergent Difference Scheme for a Semilinear Reaction-Diffusion Problem on a Shishkin Mesh},
	journal={Advances in Mathematics: Scientific Journal},
	volume={7},
	number={1},
	year={2018},
	pages={23--38},
	url={http://research-publication.com/wp-content/uploads/2019/03/AMSJ-2018-N1-4.pdf}
}

@article{karasuljic2019class,
	title={A class of difference schemes uniformly convergent on a modified Bakhvalov mesh},
	author={Karasulji{\'c}, Samir and Zarin, Helena and Duvnjakovi{\'c}, Enes},
	journal={Journal of Modern Methods in Numerical Mathematics},
	volume={10},
	number={1-2},
	pages={16--35},
	year={2019},
	publisher={ModernScience Publishers},
	doi={10.20454/jmmnm.2019.1513}
}

@article{liseikin2020numerical,
	title={Numerical analysis of grid--clustering rules for problems with power of the first type boundary layers},
	author={Liseikin, V. D. and Karasulji{\'c}, S.},
	journal={Computational Technologies},
	year={2020},
	volume={25},
	number={1},
	pages={49--66},
	doi={10.25743/ICT.2020.25.1.004}
}

@article{karasuljic2020onconstruction,
	title={On construction of a global numerical solution for a semilinear singularly-perturbed reaction 
		diffusion boundary value problem},
	author={Karasulji{\'c}, S. and Ljevakovi{\'c}, H.},
	journal={Mat. bilten},
	volume={44},
	number={2},
	pages={131--148},
	year={2020},
	doi={10.37560/matbil2020131k}
}

@article{liseikin2021oncomprehensive,
	title={On a Comprehensive Grid for Solving Problems Having Exponential or Power-of-First-Type Layers},
	author={Liseikin, V. D. and Karasulji{\'c}, S. and Mukhortov A.V. and Passonen V.I.},
	
	journal={Numerical Geometry, Grid Generation and Scientific Computing, Lecture Notes in Computational Science and Engineering 143},
	publisher={Springer Nature Switzerland AG 2021},
	year={2021},
	doi={10.1007/978-3-030-76798-3_14}
}

@inproceedings{samir2010tmt, 
	author = {Duvnjakovi\'{c},E. and Karasulji\'{c},S. and Oki\v{c}i\'{c},N.}, 
	title = {Difference Scheme for Semilinear Reaction-Diffusion Problem},
	booktitle = {14th International Research/Expert Conference Trends in the Development of Machinery 
		and Associated Technology TMT 2010, 7. Mediterranean Cruise }, 
	year = {2010},            
	pages = {793--796},
	url={ http://www.tmt.unze.ba/zbornik/TMT2010/199-TMT10-138.pdf}
}

@inproceedings{samir2011uniformnly,
	author = {Duvnjakovi\'{c},E. and Karasulji\'{c},S.},
	Title={Uniformly Convergente Difference Scheme for Semilinear Reaction-Diffusion Problem },
	booktitle = {Conference on Appllied and Scietific Computing, Trogir, Croatia },
	pages = {25},
	year = {2011 },
	url={http://applmath11.math.hr/abs_book.pd}
}

@inproceedings{samir2011skoplje,
	author = {Duvnjakovi\'{c},E.  and Karasulji\'{c},S.},
	Title={Uniformly Convergente Difference Scheme for Semilinear Reaction-Diffusion Problem  },
	booktitle = {SEE Doctoral Year Evaluation Workshop, Skopje, Macedonia},
	year = {2011 },
	url={http://users.sch.gr/alexiouth/index.php/workshops/evaluation-ws}
}

@inproceedings{samir2012uniformnly,
	author = {Karasulji\'{c},S. and Duvnjakovi\' {c},E.},
	Title={Difference Scheme for Semilinear Reaction-Diffusion Problem on Layer--adapted Mesh},
	booktitle = {The Seventh Bosnian-Herzegovinian Mathematical Conference, Sarajevo, BiH},
	year = {2012 },
	url={http://www.anubih.ba/Journals/vol.8,no-2,y12/16matskup7.pdf}
}

@inproceedings{samir2012class,
	author = {Duvnjakovi\'{c},E.  and Karasulji\'{c},S.},
	Title={Class of Difference Scheme for Semilinear Reaction-Diffusion Problem on Shishkin Mesh},
	booktitle = {MASSEE International Congress on Mathematics - MICOM 2012, Sarajevo, Bosnia and Herzegovina},
	year = {2012 }
}

@inproceedings{samir2013collocation,
	author = {Duvnjakovi\'{c},E. and Karasulji\'{c},S.},
	Title={Collocation Spline Methods for Semilinear Reaction-Diffusion Problem on Shishkin Mesh},
	booktitle = {IECMSA-2013, Second International Eurasian Conference on Mathematical Sciences and Applications, Sarajevo, Bosnia and Herzegovina},
	year = {2013 }
}

@inproceedings{samir2015construction,
	author = {Karasulji{\'{c}},S. },
	Title={Construction of the Difference Scheme for Semilinear Reaction-Diffusion  Problem on a Bakhvalov Type Mesh},
	booktitle = {The Eighth Bosnian-Herzegovinian Mathematical Conference, Sarajevo, BiH},
	year = {2015 },
	% url={http://www.anubih.ba/Journals/vol.11,no-2,y15/14Resume.pdf},
	keywords={samir}
}

@inproceedings{liseikin2019rules,
	title={On Rules for Grid Clustering in the Zones of Boundary and Interior Layers},
	author={Liseikin, V.D. and Kudryavtsev, A. N. and Paasonen, V. I. and Karasulji{\'c}, S. and Mukhortov, A. V.},
	booktitle={Mathematics and its Applications. International Conference  in honor of the 90th birthday of Sergei K. Godunov},
	year={2019},
	organization={Novosibirsk, Russia}
}

@inproceedings{liseikin2020comprehensiveproceedings,
	title={On a comprehensive grid for solving problems having exponential or power--of--first--type layers},
	author={Liseikin, Vladimir and Karasulji{\'c}, Samir and Mukhortov, Aleksandr and Paasonen, Viktor},
	booktitle={NUMGRID 2020,  Moscow, Russia},
	year={2020},
	organization={Dorodnicyn Computing Center FRC CSC RAS}
}

@book{crnjac1994matematika,
	author={Crnjac, M. and Juki\' c, D. and Scitovski, R.},
	title={Matematika},
	year={1994},
	publisher={Ekonomski fakultet Osijek, R.Hrvatska},
	isbn={953-6073-03-X}
}

@book{dedagic1997uvod,
	author={Dedagi\' c, F.},
	title={Uvod u vi\v su matematiku},
	publisher={Univerzitet u Tuzli, Tuzla, BiH},
	year={1997}
}

@book{dedagic2005matematicka,
	author={Dedagi\'c, F.},
	title={Matemati\v cka analiza, prvi dio},
	publisher={Univerzitet u Tuzli},
	address={Tuzla, BiH},
	isbn={9958-609-39-8},
	year={2005}	
}

@book{jukic1998matematika1,
	author={Juki\' c, D. and Scitovski, R.},    
	title={Matematika 1},
	publisher={Prehrambeno tehnolo\v ski fakultet  i Elektrotehni\v cki fakultet Osijek, R. Hrvatska},
	year={1998},
	isbn={953-6032-18-X}
}

@book{jukic2017matematika1,
	author={Juki\' c, D. and Scitovski, R.},    
	title={Matematika 1},
	publisher={STUDIO HS Internet d.o.o., Osijek, R. Hrvatska},
	year={2017},
	isbn={978-953-6032-18-1}
}

@book{mitrinovic1973linearna,
	author={Mitrinovi\' c, D.S. and  Mihailovi\' c, D. and Vasi\' c, P.M.},
	title={Linearna algebra, polinomi, analiti\v cka geometrija},
	year={1973},
	publisher={Beogradski izdava\v cki zavod, Beograd, SFRJ}
}

@book{uscumlic2002elementi,
	author={U\v s\' cumli\' c, M.P. and Mili\v ci\' c, P.M.},
	title={Elementi vi\v se matematike I},
	publisher={IP "Nauka", Beograd, Srbija} ,
	year={2002}
}

@book{vugdalic2009matematika,
	author={R. Vugdali\'c},
	title={Diferencijalni i integralni ra\v cun realne funkcije jedne realne promjenljive},
	series={Manualia universitas studiorum Tuzlaensis},
	year={2009},
	publisher={Ars grafika},
	address = {Tuzla, BiH},
	isbn={978-9958-9072-8-9}
}

@book{vugdalic2014matematika1,
	author={R. Vugdali\'c},
	title={Matematika 1},
	year={2014},
	publisher={IN SCAN d.o.o.},
	address={Tuzla, BiH},
	isbn={978-9958-894-22-0}
}

@article{sanela2017spectra,
	  title={The spectra of certain nonlinear superposition operators in the spaces of sequences},
	  author={Halilovi\'c, S. and Vugdali\'c, R.},
	  journal={Gulf Journal fo Mathematics},
	  volume={5},
	  number={2},
	  pages={20--30},
	  year={2017}
}

@article{sanela2014rhodius,
	     author={Halilovi\'c, S. and Vugdali\'c, R.},
	     title={The Rhodius Spectra of Some Nonlinear Superposition Operators in the Spaces of Sequences},
	     journal={Adv. Math., Sci. J.},
	     volume={3},
	     number={2},
	     pages={83--96},
	     year={2014}
}

@article{dedagic2012solvatility,
	author={Dedagi\'c, F. and Halilovi\'c, S. and Barakovi\'c, E.},
	title={On the solvability of discrete nonlinear Hammerstein systems in lp, $\sigma$ spaces},
	journal={Mathematica Balkanica},
	volume={26},
	number={3--4},
	pages={325--333},
	year={2012}
}

@article{sanela2014neuberger,
	author={Halilovi\'c, S. and Vugdali\'c,  R.},
	title={The Neuberger Spectra of Nonlinear Superposition Operators in the Spaces of Sequences},
	journal={Journal of the International Mathematical Virtual Institute},
	volume={4},
	%number={4},
	pages={97--119},
	year={2014},
	keywords={sanela}
}

@article{sanela2017general,
	author={Vugdali\'c,  R. and Halilovi\'c, S. },
	title={On general cosine operator function in Banach space},
	journal={Adv.Math., Sci.J.},
	volume={6},
	number={1},
	pages={23--27},
	year={2017}
}

@article{sanela2010nonlinear,
	author={Halilovi\'c, S. },
	title={Nonlinear Spectral Theories and Solvability of Nonlinear Hammerstein Equations},
	journal={Mathematica Balkanica},
	volume={24},
	number={3--4},
	pages={321--328},
	year={2010}
}

@article{sanela2021some,
	author={Halilovi\'c, S. },
	title={Some Spectra of Superposition operators Generated by an Exponential Function},
	journal={Communications in Mathematics and Applications},
	volume={12},
	number={1},
	pages={221--229},
	year={2021}
}

@inproceedings{sanela2021spectra,
	title={Spectra of superposition operators generated by an exponential function},
	author={Halilovi\'c, S.},
	booktitle={9th International Euroasian Conference on Mathematical Sciences and Applications},
	year={2020},
	%organization={Novosibirsk, Russia},
	keywords={sanela}
}

@article{sanela2016formula,
	author={Vugdali\'c, R. and Halilovi\'c, S. },
	title={A formula for $n$-times integrated $Co$ group of operators ($n\in\mathbb{N}$)},
	journal={Advances in Mathematics: Scientific Journal},
	volume={5},
	number={2},
	pages={161--166},
	year={2016}
}

@thesis{sanela2015spektri,
       author={Halilovi\'c, Sanela},
       title={Spektri nelinearnog operatora superpozicije u Banachovim prostorima nizova i primjene},
       year={2015},
       note={Doktorska disertacija}
}

@inproceedings{sanela2015spectra,
	title={The spectra of nonlinear superposition operators},
	author={Halilovi\'c, S. and Vugdali\'c, R.},
	booktitle={The eight Bosnian-Herzegovinian Mathematical Conference, Sarajevo},
	year={2015},
	%organization={Novosibirsk, Russia},
	keywords={sanela}
}

@article{sanela2011spektru,
	author={Halilovi\'c, S. },
	title={O spektru nelinearnih operatora},
	journal={Math. e},
	volume={20},
	%number={2},
	pages={25--36},
	year={2011},
	keywords={sanela}
}

@inproceedings{sanela201ospectral,
	title={Spectral radius of the operator superposition in the spaces of sequences},
	author={Halilovi\'c, S. },
	booktitle={The fifth Bosnian-Herzegovinian Mathematical Conference},
	year={2010},
	%organization={Novosibirsk, Russia},
	keywords={sanela}
}

@thesis{sanela2009spektralne,
	author={Halilovi\'c, Sanela},
	title={Spektralne teorije nelinearnih operatora u Banachovom prostoru},
	year={2009},
	note={Magistarski rad}
}

@inproceedings{sanela2009onsolvability,
	title={On the Solvability of Discrete Nonlinear Hammerstein Systems},
	author={Dedagi\'c, F. and Halilovi\'c, S. },
	booktitle={International Congress on Mathematics MICOM},
	year={2009},
	%organization={Novosibirsk, Russia},
	keywords={sanela}
}

@thesis{sanela2004prostor,
	author={Alihod\v zi\'c, Sanela},
	title={Prostor i reprezentacija linearnih funkcionala},
	year={2004},
	note={Diplomski rad}
}

@article{samir20152d,
         author={Karasulji{\'c}, S. and Duvnjakovi{\'c}, E. and  Zarin, H.},
        title={A uniformly convergent difference scheme on a Shishkin type mesh for the 2D singular perturbation boundary value problem},
         year={2015},
         note={Submitted}
          }

@article{samir2015modifiedbah,
         author={Karasulji{\'c}, S.  and Duvnjakovi{\'c}, E. and Zarin, H.},
         title={A uniformly convergent difference scheme on a modified Bakhvalov mesh for the singular perturbation boundary value problem},
         year={2018},
         note={Submitted}
        }

@article{karasuljic2020construction,
	title={On construction of a global numerical solution for a semilinear singularly--perturbed reaction diffusion boundary value problem},
	author={Karasulji{\'c}, Samir and Ljevakovi{\'c}, Hidajeta},
	journal={Mat. Bilten},
	volume={44},
	number={2},
	pages={131--148},
	year={2020},
	publisher={Union of Mathematicians of Macedonia},
	doi={10.37560/matbil2020131k}
}

@article{duvnjakovic2001variational,
         author={Duvnjakovi{\'c}, E.},
         title={A Variational Difference Scheme for Singularly Perturbed One-dimensional Reaction-diffusion problem},
         journal={ Radovi Matemati\v cki, ANUBiH-Sarajevo, },
        % publisher={ANUBiH-Sarajevo}, 
         volume={10},
         pages={ 109-119},
         year={ 2001} 
}

@article{duvnjakovic2004difference,
         author={Oki{\v c}i{\'c}, N. and Duvnjakovi{\'c}, E.},
         title={Difference scheme for nonlinear singular perturbation problem},
         journal={Zbornik radova PMF},
         volume={1},
         pages={53-60},
         year={2004}
         }

@article{duvnjakovic2005procjena,
         author={Duvnjakovi{\'c}, E.},
         title={Procjena rje{\v s}enja kvazilinearnog elipti{\v c}kog rubnog problema},
         journal={Zbornik radova PMF},
         volume={2},
         pages={54-59},
         year={2005}
          }

@inproceedings{duvnjakovic2008difference,
               author={Duvnjakovi{\' c}, E. and Oki{\v c}i{\'c}, N.},
               title={Difference scheme for semilinear reaction-diffusion problem of elliptic type},
               booktitle={ 12th  International Research/Expert Conference "Trends in The Development of Machinery and Associated Technology" TMT 2008, Istambul, Turkey},
               year={2008},
               pages={1117-1120},
               url={http://www.tmt.unze.ba/zbornik/TMT2008/280-TMT08-104.pdf}
                }

@article{allen1955relaxation,
         author={Allen, D.N. de G. and Southwell, R.V.},    
         title={Relaxation methods applied to determine the motion, in two dimensions, 
                of a  viscous fluid past a fixed cylinder},
         journal={Q J Mechanics Appl Math},
         year={1955},
         volume={8},
         number={2},
         doi={10.1093/qjmam/8.2.129}
       }

@article{axelsson1990monotonicity,
      author={Axelsson, O and Kolotilina, L. },
      title={Monotonicity and discretization error estimates},
      journal={SIAM J. Numer. Anal.},
      number={27},
      volume={6},
      pages={1591--1611},
      year={1990},
      ISSN = {00361429},
      URL = {http://www.jstor.org/stable/2157772}
          }

@article{bahvalov1969,
        title={Towards optimization of methods for solving boundary value problems in 
                the presence of boundary layers},
        journal={Zh. Vychisl. Mat. i Mat. Fiz.},
        author={Bakhvalov, N. S.},
        volume={9},
        pages={841--859},
        year={1969},
        language={(In Russian)},
        doi={10.1016/0041-5553(69)90038-X},
        url={http://mi.mathnet.ru/eng/zvmmf/v9/i4/p841 }
       }

@article{bahvalov1969a,
%         author={\foreignlanguage{russian}{Bahvalov}, \foreignlanguage{russian}{N. S.}},
%         title={\foreignlanguage{russian}{ K optimizacii metodov resheniya kraevyh zadach pri nalichi pogra--nichnogo sloya}},
%         journal={{\foreignlanguage{russian}{Zhurnal Bychislitel}}b{\foreignlanguage{russian}{no{\u i} Matematiki i Matematichesko{\u i} Fiziki}}},       
%         journal={{\foreignlanguage{russian}{Zhurnal Bychislitel\char126 no{\u i} Matematiki i Matematichesko{\u i} Fiziki}}},       
%         volume={ Tom 9, No 4},
%         year={1969},
%         pages={555-568},
%         language={(Russian)}
%         }

@article{beckett2000convergence,
title = "Convergence analysis of finite difference approximations on equidistributed grids to a singularly perturbed boundary value problem ",
journal = "Applied Numerical Mathematics ",
volume = "35",
number = "2",
pages = "87 - 109",
year = "2000",
issn = "0168-9274",
doi = {10.1016/S0168-9274(99)00065-3},
%url = "http://www.sciencedirect.com/science/article/pii/S0168927499000653",
url={http://goo.gl/hmJCfg},
author = "Beckett, G. and Mackenzie, J.A."
}

@article{berger1970asymptotic,
        author={ Berger, M. S. and Fraenkel, L. E.},
        title={On the asymptotic solution of a nonlinear Dirichlet problem},
        journal={J. Math. Mech.},
        volume={ 19},
        year={1970},
        pages={ 553-585}
        }

@article{berger1981analysis,
          title={ An analysis of a uniformly accurate difference method for a singular perturbation problem},
          author={  Berger,A. E. and   Solomon, J. M. and  Ciment, M.},
          journal={ Math. Comp.},
          volume={ 37},
          year={1981},
          pages={ 79-94},
          doi={10.1090/S0025-5718-1981-0616361-0 }
          }

@article{blatov1992galerkin,
         author={Blatov,  I. A.},
         title={On the Galerkin finite-element method for elliptic quasilinear singularly perturbed boundary value
                problems I},
         journal={ Differential Equations},
         volume={ 28},
         number={7},
         year={1992},
         pages={ 931-940},
         note={Russian}
          }

@article{blatov1992agalerkin,
         author={Blatov,  I.A.},
         title={On the Galerkin finite-element method for elliptic quasilinear singularly perturbed boundary value
                problems II},
         journal={ Differential Equations},
         volume={ 28},
         number={10},
         year={1992},
         pages={ 1469-1477},
          note={Russian}
          }

@article{blatov1994galerkin,
         author={Blatov,  I. A.},
         title={On the Galerkin finite-element method for elliptic quasilinear singularly perturbed boundary value
                problems III. Problems with angular boundary layers},
         journal={ Differential Equations},
         volume={ 30},
         number={3},
         year={1994},
         pages={ 432-444},
          note={Russian}
          }

@article{boglaev1984approximate,
       title = "Approximate solution of a non-linear boundary value problem with a small parameter for the highest-order differential ",
       journal={ Zh. Vychisl. Mat. i Mat. Fiz.},
       volume = "24",
       number = "11",
       pages = "1649 -- 1656",
       year = "1984",
       doi = {10.1016/0041-5553(84)90005-3},
       url = "http://www.mathnet.ru/links/3d5582cba89c929035e5d7e541089694/zvmmf4282.pdf",
       author = {Boglaev, I. P.},
       note={In Russian}
       }

@article{boglaev1988numerical,
         author={Boglaev, I. P.},
         title={A numerical method for solving a quasilinear elliptic equation with a small parameter at the highest derivatives},
         journal={ Zh. Vychisl. Mat. Mat. Fiz.},
         year={1988},
         volume={28},
         issue={4},
         pages={492--502},
         doi={ 10.1016/0041-5553(88)90156-5},
         url ={ http://www.mathnet.ru/links/905d71bee117e940745a8489cc67fa55/zvmmf3658.pdf},
         note={In Russian}
         }

@article{boglaev2005uniform,
         author={Boglaev, I.P.},
         title={Uniform convergence of monotone iterative methods for semilinear singularly perturbed problems of elliptic and parabolic types},
         volume={ 20},
         pages={ 86-103},
         year={ 2005},
         journal={Electronic Transactions on Numerical Analysis (ETNA)},
         publisher={Kent State University},
         url={https://eudml.org/doc/125397}
         }

@article{boglaev2006monotone,
  title={Monotone finite difference domain decomposition algorithms and applications to nonlinear singularly perturbed reaction-diffusion problems},
  author={Boglaev, I. and Hardy, M.},
  journal={Advances in Difference Equations},
  volume={2006},
  number={1},
  pages={1-38},
  year={2006},
  publisher={Hindawi Publishing Corporation 410 PARK AVENUE, 15TH FLOOR,\# 287 PMB, NEW YORK, NY 10022 USA},
  URL = {http://www.advancesindifferenceequations.com/content/2006/1/070325},
  doi={10.1155/ADE/2006/70325}
}

@article{desanti1986boundary,
         author={DeSanti, A. J.},
         title={ Boundary and interior layer behavior of solutions of some singularly perturbed 
                 semilinear elliptic boundary value problems},
         journal={ J. Math. Pures Appl.},
         volume={ 65},
         year={1986},
         pages={ 227-262}
        }

@article{desanti1987singularly,
          title={  Singularly perturbed Dirichlet problems with subquadratic nonlinearities},
          author={  DeSanti, A. J.},
          journal={ Trans. Amer. Math. Soc.},
          volume={ 303},
          year={1987},
          pages={ 585-593},
          note={MathSciNet review: 902786},
          url={ http://www.ams.org/journals/tran/1987-303-02/S0002-9947-1987-0902786-4/},
          doi={10.1090/S0002-9947-1987-0902786-4 }
          }

@article{dyke1994niteenth,
author = { van Dyke, M.},
title = {Nineteenth-Century Roots of the Boundary-Layer Idea},
journal = {SIAM Review},
volume = {36},
number = {3},
pages = {415-424},
year = {1994},
doi = {10.1137/1036097},
URL = { http://dx.doi.org/10.1137/1036097},
eprint = {  http://dx.doi.org/10.1137/1036097}
}

@article {elmistikawy1978numerical,
	title = {Numerical method for boundary layers with blowing- The exponential box scheme},
	journal = {AIAA Journal},
	volume = {16},
	number={7},
	year = {1978},
	pages = {749-751},
	author = { El-Mistikawy, T.M.A. and Werle, M.J.},
	doi={10.2514/3.7573 }
           }

@article{farrell1988sufficient,
         author={Farrell, P. A.},           
         title={Sufficient Conditions for the Uniform Convergence of a Difference Scheme for a Singularly Perturbed Turning Point Problem},
         journal ={Siam J. Numer. Anal.},
         volume={25},
         number={3},
         year={1988},
         doi={10.1137/0725038}
         }

@article{fife1973semilinear,
          year={1973},
          journal={Archive for Rational Mechanics and Analysis},
          volume={52},
          number={3},
          doi={10.1007/BF00247733},
          title={Semilinear elliptic boundary value problems with small parameters},
          url={http://dx.doi.org/10.1007/BF00247733},
          publisher={Springer-Verlag},
          author={Fife, P. C.},
          pages={205-232},
          language={English}
          }

@article{fife1974interior,
  author={Fife, P. C. and Greenlee, W.M.},
  title={Interior transition layers for elliptic boundary value problems with a small parameter},
  journal={Russian Mathematical Surveys},
  volume={29},
  number={4},
  pages={103},
  url={http://stacks.iop.org/0036-0279/29/i=4/a=A04},
  year={1974},
  doi={10.1070/RM1974v029n04ABEH001291}
     }

@article{gartland1987uniform,
         author={Gartland, E. C. Jr.},
         title={Uniform High-Order Difference Schemes for a Singularly Perturbed Two-Point Boundary Value Problem},
         journal={Mathematics of Computation},
         volume={48},
         number={78},
         year={1987},
         pages={551-564},
         url={http://www.jstor.org/stable/2007827}
         }

@article{gartland1988graded,
          author={Gartland, E. C. Jr.},
          title={Graded-mesh difference schemes for singularly perturbed two-point boundary value problem},
          journal={Math. Comput.},
           year={1988},
           pages={631--657},
           doi={10.1090/S0025-5718-1988-0935072-1},
           volume={51},
            note={MR935072}
             }

@article{denisov2008cornerboundary,
        author={Denisov,I.V.}, 
        title={Corner boundary layer in nonlinear singularly perturbed elliptic problems},
        journal={Zh. Vychisl. Mat. Fiz.,},
        year={2008},
        volume={48},
        number={1},
        pages={62--79},
        url={http://www.mathnet.ru/links/ba9903a715b401d25da409c4545a4a80/zvmmf195.pdf},
        doi={10.1007/s11470-008-1005-7}
        }

@book{hemker1977numerical,
  title={A numerical study of stiff two-point boundary problems},
  author={Hemker, P.W.},
  isbn={9789061961468},
  lccn={78316930},
  series={Mathematical Centre tracts},
  url={https://books.google.ba/books?id=JjfvAAAAMAAJ},
  year={1977},
  publisher={Mathematisch Centrum}
}

@article{herceg1981some,
          author={Herceg, D. and Vulanovi{\' c}, R.},
          title={Some Finite-Difference Schemes for a Singular Perturbation Problem on a Non-Uniform Mesh},
          journal={Novi Sad J. Math},
          volume={11},
          pages={117-134},
          year={1981}
          }

@article{herceg1982numerical,
          author={Cvetkovi{\' c}, Lj. and Herceg, D.},
          title={On a Numerical Solution of the Boundary Value Problem Using an Optimal Numerical Differentiation},
          journal={Novi Sad J. Math},
          volume={12},
          pages={177-189},
          year={1982},
          url={http://www.dmi.uns.ac.rs/nsjom/Papers/12/NSJOM_12_177_189.pdf}
          }

@article{herceg1982some,
          author={Herceg, D.},
          title={Some Difference Schemes for Two Point Boundary Value Problems},
          journal={Novi Sad J. Math},
          volume={12},
          pages={123-137},
          year={1982},
          url={http://www.dmi.uns.ac.rs/nsjom/Papers/12/NSJOM_12_123_137.pdf}
          }

@article{herceg1987numerical,
          author={Herceg, D. and Petrovi{\' c}, N.},
          title={On Numerical Solution of a Singularly Perturbed Boundary Value Problem II},
          journal={Novi Sad J. Math},
          volume={17},
          number={1},
          pages={163-186},
          year={1987}
          }

@article{herceg1990,
       year={1989},
       journal={Numer. Math.},
       volume={56},
       number={7},
       doi={10.1007/BF01405196},
       title={Uniform fourth order difference scheme for a singular perturbation problem},
       url={http://dx.doi.org/10.1007/BF01405196},
       publisher={Springer-Verlag},
       author={Herceg, D.},
       pages={675--693},
       url = {http://eudml.org/doc/133421}
       }

@article{herceg1990numerical,
          author={Herceg, D.},
          title={On the Numerical Solution of a Singularly Perturbed Nonlocal Problem},
          journal={Novi Sad J. Math},
          volume={20},
          number={2},
          pages={1-10},
          year={1990}
          }

@article{herceg1991,
        author={Herceg,D. and Surla,K.},
        title={Solving a nonlocal singularly perturbed problem by spline in tension},
        journal={Review of research Faculty of Sciences-University of Novi Sad},
        volume={Vol.21, No.2},
        pages={119-132},
        year={1991},
        url={http://www.dmi.uns.ac.rs/nsjom/Papers/21_2/NSJOM_21_2_119_132.pdf}
        }

@article{herceg1993discrete,
          title={On a Discrete Analogue for Boundary Value Problem},
          author={Herceg, D.},
          journal={Novi Sad J. Math},
          volume={23},
          number={2},
          pages={401-410},
          year={1993}
          }

@article{herceg1994numerical,
          author={Herceg, D.},
          title={Numerical Solution of Some Discrete Analogues of Boundary Value Problem},
          journal={Novi Sad J. Math},
          volume={24},
          number={2},
          pages={187-196},
          year={1994}
          }

@article{herceg1994finite,
          author={Herceg, D.},
          title={On a Finite Difference Analogue for a Boundary Value Problem},
          journal={Novi Sad J. Math},
          volume={24},
          number={1},
          pages={381-388},
          year={1994}
          }

@article{herceg1998,
        author={Herceg, D. and Surla,K. and Rapaji\'{c},S.},
        title={Cubic spline difference scheme on a mesh of a Bakhvalov type},
        journal={Novi Sad J. Math.},
        volume={28},
        number={3},
        year={1998},
        pages={41--49},
        url={http://www.dmi.uns.ac.rs/nsjom/Papers/28_3/NSJOM_28_3_041_049.pdf}
        }

@article{herceg1999acceleration,
         author ={Herceg, D. and Mali{\v c} i{\'c}, H.},
         title={On Acceleration of Solving Singularly Perturbed Boundary Value Problem},
         journal={Novi Sad J. Math.},
         volume={29},
         number={1},
         year={1999},
         pages={155-168}
         }

@article{herceg2000,
        author={Herceg,D. and Mali\v{c}i\'c,H. and Liki\'c,I.},
        title={On a finite difference analogue of fourth order for boundary value problem},
        journal={Novi Sad J. Math.},
        volume={30},
        number={1},
        year={2000},
        pages={197--203}
       }

@article{herceg2000themesh,
          author={Herceg, D. and  Mali\v{c}i\'c,H.},
          title={The mesh chopping algorithm for singularly perturbed boundary value problem},
          journal={Novi Sad J.Math},
          volume={30},
          number={1},
          year={2000},
          pages={155--164}
        }

@article{herceg2000finite,
          author={Herceg, D. and Kreji{\' c}, N. and Mali{\v c}i{\'c}, H.},
          title={On a Finite Difference Analogue of a Singular Boundary Value Problems},
          journal={Novi Sad J. Math},
          volume={30},
          number={1},
          pages={23-31},
          year={2000},
          url={http://www.dmi.uns.ac.rs/nsjom/Papers/30_1/NSJOM_30_1_023_031.pdf}
          }

@article{herceg2000mesh,
          author={Herceg, D. and Mali{\v c}i{\' c}, H.},
          title={The Mesh Chopping Algorithm for Singularly Perturbed Boundary Value Problem},
          journal={Novi Sad J. Math},
          volume={30},
          number={1},
          year={2000},
          pages={155-164},
          url={http://www.dmi.uns.ac.rs/nsjom/Papers/30_1/NSJOM_30_1_155_164.pdf}
          }

@article{herceg2001numerical,
          author={Herceg, D. and Surla, K. and Radeka, I. and Mali{\v c}i{\' c}, H.},
          title={Numerical Experiments with Different Schemes for a Singularly Perturbed Problem},
          journal={Novi Sad J. Math},
          volume={31},
          number={1},
          year={2001},
          pages={93-101},
          url={http://www.dmi.uns.ac.rs/nsjom/Papers/31_1/NSJOM_31_1_093_101.pdf}
          }

@article{herceg2001some,
          author={Surla, K. and Herceg, D. and Mali{\v c} i{\' c}, H.},
          title={Some Comparisons of Difference Schemes on Meshes of Shishkin and Bakhvalov Type},
          journal={Novi Sad J. Math},
          volume={31},
          number={1},
          pages={133-140},
          year={2001}
          }

@article{herceg2003,
        title={On numerical solution of semilinear singular perturbation problems by using the Hermite scheme on a new Bakhvalov-type mesh},
        author={Herceg, D. and Miloradovi\'{c}, M.},
        journal={Novi Sad J. Math},
        volume={33},
        number={1},
        pages={145--162},
        year={2003},
        url={http://www.dmi.uns.ac.rs/nsjom/Papers/33_1/nsjom_33_1_145_162.pdf}
        }

@article{herceg2003a,
        title={On a fourth-order finite difference method for nonlinear two-point boundary value problems},
        author={Herceg, D. and Herceg, {Dj}.},
        journal={Novi Sad J. Math},
        volume={33},
        number={2},
        pages={173--180},
        year={2003},
        url={http://www.dmi.uns.ac.rs/nsjom/Papers/33_2/nsjom_33_2_173_180.pdf}
        }

@article{herceg2003high,
          author={Radeka, I. and Herceg, D.},
          title={High-Order Methods for Semilinear Singularly Perturbed Boundary Value Problems},
          journal={Novi Sad J. Math},
          volume={33},
          number={2},
          pages={139-161},
          year={2003}
          }

@article{herceg2004nonequdistant,
          author={Herceg, D. and Herceg, Dj.},
          title={On a Nonequdistant Difference Scheme of Chawla Type},
          journal={Novi Sad J. Math},
          volume={34},
          number={1},
          pages={171-190},
          year={2004},
          url={http://www.dmi.uns.ac.rs/nsjom/Papers/34_1/nsjom_34_1_171_190.pdf}
           }

@article{howes1979singularly,
          author={Howes, F. A.},
          title={ Singularly perturbed semilinear elliptic boundary value problems},
          journal={Communications in Partial Differential Equations},
          volume={ 4},
          issue={1},
          year={1979},
          pages={ 1--39},
          doi={10.1080/03605307908820090}
          }

@article{illin1969differencing,
%        author={Il'in, Arlen Mikha{\"{i}}lovitch and Minachin},
%        journal={Mathematical notes of the Academy of Sciences of the USSR},
%        note={in Russian},
%        year={1969}, 
%        volume={6},
%        number ={2},
%        pages={596-602},
%        title={Differencing scheme for a differential equation with a small parameter affecting the highest derivative}
%           }

@article{jiang2001local,
title = "Local exponentially fitted finite element schemes for singularly perturbed convection--diffusion problems ",
journal = "Journal of Computational and Applied Mathematics ",
volume = "132",
number = "2",
pages = "277 - 293",
year = "2001",
%doi = "10.1016/S0377-0427(00)00329-0",
url = "http://www.sciencedirect.com/science/article/pii/S0377042700003290",
author = "Lishang, J. and Xingye, Y.",
        doi={10.1016/S0377-0427(00)00329-0}
 }

@article{kadalbajoo2010brief,
%title = "A brief survey on numerical methods for solving singularly perturbed problems ",
%journal = "Applied Mathematics and Computation ",
%volume = "217",
%umber = "8",
%pages = "3641 - 3716",
%year = "2010",
%doi = "10.1016/j.amc.2010.09.059",
%url = "http://www.sciencedirect.com/science/article/pii/S0096300310010532",
%author = "Kadalbajoo, M.K. and Gupta, V."
%}

@article{kelley1990semilinear,
title = "Semilinear elliptic singular perturbation problems with nonuniform interior behavior ",
journal = "Journal of Differential Equations ",
volume = "86",
number = "1",
pages = "88 - 101",
year = "1990",
issn = "0022-0396",
doi = "10.1016/0022-0396(90)90042-N",
url = "http://www.sciencedirect.com/science/article/pii/002203969090042N",
author = " Kelley, W. and  Ko, B."
}

@article{kopteva1999uniform,
  title={Uniform convergence with respect to a small parameter of a scheme with central difference on refining grids},
  author={Kopteva, N. V.},
  journal={Computational mathematics and mathematical physics},
  volume={39},
  number={10},
  pages={1594--1610},
  year={1999},
  publisher={Oxford; New York: Pergamon Press, c1992-},
  url={http://mi.mathnet.ru/eng/zvmmf/v39/i10/p1662 }
}

@article{kopteva2001,
  title={Uniform second-order pointwise convergence of a central difference approximation for a quasilinear convection-diffusion problem},
  author={Kopteva, N. and Lin{\ss}, T.},
  journal={J. Comput. Appl. Math.},
  volume={137},
  number={2},
  pages={257--267},
  year={2001},
  publisher={Elsevier},
  doi={10.1016/S0377-0427(01)00353-3}
}

@article{kopteva2001maximum,
  title={Maximum norm a posteriori error estimates for a one-dimensional convection-diffusion problem},
  author={Kopteva, N.},
  journal={SIAM Journal on numerical analysis},
  volume={39},
  number={2},
  pages={423--441},
  year={2001},
  publisher={SIAM},
  DOI={10.1137/S0036142900368642}
}

@article{kopteva2001robust,
  title={A robust adaptive method for a quasi-linear one-dimensional convection-diffusion problem},
  author={Kopteva, N. and Stynes, M.},
  journal={SIAM Journal on Numerical Analysis},
  volume={39},
  number={4},
  pages={1446--1467},
  year={2001},
  publisher={SIAM},
  doi={10.1137/S003614290138471X}
}

@article{kopteva2001uniform,
  title={Uniform pointwise convergence of difference schemes for convection-diffusion problems on layer-adapted meshes},
  author={Kopteva, N.},
  journal={Computing},
  volume={66},
  number={2},
  pages={179--197},
  year={2001},
  publisher={Springer-Verlag},
  doi={10.1007/s006070170034}
}

@article{kopteva2003error,
  title={Error expansion for an upwind scheme applied to a two-dimensional convection-diffusion problem},
  author={Kopteva, N.},
  journal={SIAM Journal on numerical analysis},
  volume={41},
  number={5},
  pages={1851--1869},
  year={2003},
  publisher={SIAM},
  url={http://ftp.demec.ufpr.br/CFD/bibliografia/MER/Kopteva_2003.pdf}%,
  %doi={10.1137/S003614290241074X}
  }

@article{kopteva2004,
  title={Numerical analysis of a singularly perturbed nonlinear reaction--diffusion problem with multiple solutions},
  author={Kopteva, N. and Stynes, M.},
  journal={Appl. Numer. Math.},
  volume={51},
  number={2},
  pages={273--288},
  year={2004},
  publisher={Elsevier},
  doi={10.1016/j.apnum.2004.07.001}
}

@article{kopteva2004accurate,
  title={How accurate is the streamline-diffusion FEM inside characteristic (boundary and interior) layers?},
  author={Kopteva, N.},
  journal={Computer methods in applied mechanics and engineering},
  volume={193},
  number={45},
  pages={4875--4889},
  year={2004},
  publisher={Elsevier}
}

@article{kopteva2005error,
  title={Error analysis of a 2D singularly perturbed semilinear reaction-diffusion problem1},
  author={Kopteva, N.},
  journal={Mathematical Modelling and Analysis},
  pages={227--233},
  year={2005},
  url={http://www.staff.ul.ie/natalia/pdf/Kopteva_trakai.pdf}
}

@article{kopteva2007maximum,
  title={Maximum norm error analysis of a 2d singularly perturbed semilinear reaction-diffusion problem},
  author={Kopteva, N.},
  journal={Mathematics of Computation},
  volume={76},
  number={258},
  pages={631--646},
  year={2007},
  doi={10.1090/S0025-5718-06-01938-7 },
  url={http://www.ams.org/journals/mcom/2007-76-258/S0025-5718-06-01938-7},
  publisher={Washington, DC: National Academy of Sciences-National Research Council,[1960?-}
}

@article{kopteva2007maximum1d,
  title={Maximum norm a posteriori error estimates for a 1D singularly perturbed semilinear reaction--diffusion problem},
  author={Kopteva, N.},
  journal={IMA Journal of numerical analysis},
  volume={27},
  number={3},
  pages={576--592},
  year={2007},
  publisher={Oxford University Press},
  doi={ 10.1093/imanum/drl020}
}

@article{kopteva2008maximum,
  title={Maximum norm a posteriori error estimate for a 2D singularly perturbed semilinear reaction-diffusion problem},
  author={Kopteva, N.},
  journal={SIAM Journal on Numerical Analysis},
  volume={46},
  number={3},
  pages={1602--1618},
  year={2008},
  publisher={SIAM},
  doi={10.1137/060677616}
}

@article{kopteva2009,
  title={A robust overlapping Schwarz method for a singularly perturbed semilinear reaction-diffusion problem with multiple solutions},
  author={Kopteva, N. and Pickett, M. and Purtill, H.},
  journal={Int. J. Numer. Anal. Model},
  volume={6},
  pages={680--695},
  year={2009},
  url={http://www.global-sci.org/ijnam/readabs.php?vol=6&no=4&doc=680&year=2009&ppage=695}
}

@article{kopteva2010singularly,
title = "A singularly perturbed semilinear reaction--diffusion problem in a polygonal domain ",
journal = "Journal of Differential Equations ",
volume = "248",
number = "1",
pages = "184 - 208",
year = "2010",
doi = "10.1016/j.jde.2009.08.020",
url = "http://www.sciencedirect.com/science/article/pii/S0022039609003209",
author = "Kellogg, B. R. and  Kopteva, N."
}

@article{kopteva2010shishkin,
         author="Kopteva, N. and O'Riordan, E.",
         title={Shishkin meshes in the numerical solution of singularly perturbed
                  differential equations},
         journal="Int. J. Numer. Anal. Model. ",
         volume="7",
         number="3",
         pages="393-415",
         year="2010",
         url={http://www.math.ualberta.ca/ijnam/Volume-7-2010/No-3-10/2010-03-01.pdf}
        }

@article{kopteva2011robust,
%  title={A robust grid equidistribution method for a one--dimensional singularly perturbed semilinear reaction--diffusion problem},
%  author={Chadha, N. M. and Kopteva, N.},
%  journal={IMA Journal of numerical analysis},
%  %volume={31},
 % %number={1},
 %%pages={188--211},
%  pages={1--23};
%  %year={2011},
%  year={2005},
%  publisher={Oxford University Press},
%  %doi={10.1093/imanum/drp033},
%  doi={10.1093/imanum/dri000}
%   }

@article{kopteva2011maximum,
  title={Maximum norm a posteriori error estimate for a 3d singularly perturbed semilinear reaction-diffusion problem},
  author={Chadha, N. M. and Kopteva, N.},
  journal={Advances in Computational Mathematics},
  volume={35},
  number={1},
  pages={33--55},
  year={2011},
  publisher={Springer},
  doi={10.1007/s10444-010-9163-2}
}

@article{kopteva2011stabilised,
  title={Stabilised approximation of interior-layer solutions of a singularly perturbed semilinear reaction--diffusion problem},
  author={Kopteva, N. and Stynes, M.},
  journal={Numerische Mathematik},
  volume={119},
  number={4},
  pages={787--810},
  year={2011},
  publisher={Springer},
  doi={10.1007/s00211-011-0395-y}
}

@article{kopteva2012second,
  title={A second-order overlapping Schwarz method for a 2D singularly perturbed semilinear reaction-diffusion problem},
  author={Kopteva, N. and Pickett, M.},
  journal={Mathematics of Computation},
  volume={81},
  number={277},
  pages={81--105},
  year={2012},
  doi={10.1090/S0025-5718-2011-02521-4 }
}

@article{kopteva2014maximum,
  title={Maximum-norm a posteriori error estimates for singularly perturbed elliptic reaction-diffusion problems},
  author={Demlow, A. and Kopteva, N.},
  year={2014}
}

@article{kreiss1986numerical,
           title={ Numerical Methods for Stiff Two-Point Boundary Value Problems},
           author={ Kreiss, H.O. and  Nichols,  N. K. and  Brown, D. L.},
           journal={SIAM Journal on Numerical Analysis},
           volume={ 23},
           number={2},
           year={1986},
           pages={ 325-368},
           publisher={Society for Industrial and Applied Mathematics},
           url={ http://www.jstor.org/stable/2157471}
          }

@article{ladyzhnskaya1965quasi,
         author={Ladyzhenskaya,O. A. and Ural'tseva, N. N.},
         title={Quasi-linear elliptic equations and variational problems in several independent variables},
         journal={ Russ. Math. Surv. },
         volume={16},
         number={1},
         year={1961},
         pages={19--90},
         doi={10.1070/RM1961v016n01AB-\\EH004099},
         url={http://goo.gl/9rajxk}
         %url={http://stacks.iop.org/0036-0279/16/i=1/a=R02}
         }

@article{lee1997fast,
           author={ Lee, J.Y. and  Greengard, L.},
           doi={10.1137/S1064827594272797},
           journal={SIAM J. Sci. Comput.},
           volume={ 18},
           number={2},
           pages={403--429},
           title={A Fast Adaptive Numerical Method for Stiff Two-Point Boundary Value Problems},
           year={1997}
           }

@article{li1999global,
title = "A global uniformly convergent finite element method for a quasi-linear singularly perturbed elliptic problem ",
journal = "Computers and  Mathematics with Applications ",
volume = "38",
number = "5--6",
pages = "197 - 206",
year = "1999",
doi = "10.1016/S0898-1221(99)00226-6",
url = "http://www.sciencedirect.com/science/article/pii/S0898122199002266",
author = "Li, J.  and Navon, M."
}

@article{lins1999anupwind,
title = "An upwind difference scheme on a novel Shishkin-type mesh for a linear convection--diffusion problem ",
journal = "Journal of Computational and Applied Mathematics ",
volume = "110",
number = "1",
pages = "93 - 104",
year = "1999",
doi = "10.1016/S0377-0427(99)00198-3",
url = "http://www.sciencedirect.com/science/article/pii/S0377042799001983",
author = "Lin{\ss}, T."
}

@article{linss2000analysis,
        author={Lin{\ss}, T.},
        title={Analysis of a Galerkin finite element method on a Bakhvalov\\ -Shishkin mesh for a linear 
               convection--diffusion problem},
        journal={IMA J. Numer. Anal.},
        year={2000},
        volume={20},
        number={4},
        pages={ 621-632},
        doi={ 10.1093/imanum/20.4.621}
         }

@article{linss2000uniform,
  title={Uniform pointwise convergence on Shishkin-type meshes for quasi-linear convection-diffusion problems},
  author={Lin{\ss}, T. and Roos, H.G. and Vulanovi{\'c}, R.},
  journal={SIAM J. Numer. Anal.},
  volume={38},
  number={3},
  pages={897--912},
  year={2000},
  publisher={SIAM},
  doi={10.1137/S0036142999355957}
       }

@article{lins2008robust,
author = {Lin\ss, T.},
title = {Robust convergence of a compact fourth-order finite difference scheme for reaction--diffusion problems},
journal = {Numerische Mathematik},
volume = {111},
year = {2008},
pages = {239--249},
issue = {2},
doi = {10.1007/s00211-008-0184-4},
masid = {16020132}
}

@article{linss2001uniform, 
         year={2001},
         journal={Computing},
         volume={66},
         number={1},
         doi={10.1007/s006070170037},
         title={Uniform Pointwise Convergence of Finite Difference Schemes Using Grid Equidistribution},
         url={http://goo.gl/mWpsEo},
         %url={http://dx.doi.org/10.1007/s006070170037},
         publisher={Springer-Verlag},
         author={Lin{\ss}, T.},
         pages={27-39}
         }

@article{linss2002auniformly,
         author = {Wollstein, D. and Lin{\ss}, T. and Roos, H.G.},
         title = {A Uniformly Accurate Finite Volume Discretization for a Convection-Diffusion Problem},
     	year = {2002},
     	Volume={13},
     	pages={1-11},
     	journal={Eletronic Transactions on Numerical Analysis},
     	url={http://citeseerx.ist.psu.edu/viewdoc/summary?doi=10.1.1.399.2736}
     	}

@article{linss2008afinite,
        title={A Finite Difference Method on Layer-Adapted Meshes for an Elliptic Reaction-Diffusion System in Two Dimensions},
        author={ Kellogg, B.R. and Lin{\ss}, T. and Stynes, M.},
        journal={Mathematics of Computation},
        volume={77}, 
        number={264},
        year={2008},
        pages={2085-2096},
        doi={10.1090/S0025-5718-08-02125-X }
}

@article{linss2012approximation,
  title={Approximation of singularly perturbed reaction-diffusion problems by quadratic $C^1$-splines},
  author={Lin{\ss}, T. and Radojev, G. and Zarin, H.},
  journal={Numerical Algorithms},
  volume={61},
  number={1},
  pages={35--55},
  year={2012},
  publisher={Springer US},
  doi={10.1007/s11075-011-9529-7}
         }

@article{lorenz1982stability,
         author={Lorenz, J.},
         title={Stability and monotonicity properties of stiff quasilinear boundary problems},
         journal={Zb.rad. Prir. Mat. Fak. Univ. Novom Sadu, Ser. Mat.},
         volume={ 12},
         year={1982},
         pages={ 151--176},
         note={ MR 85e:34046},
         url={http://www.emis.de/journals/NSJOM/Papers/12/NSJOM_12_151_175.pdf}
         }

@article{mackenzie96uniform,
    author = {Mackenzie, J.A. },
    title = {Uniform Convergence Analysis of an Upwind Finite-Difference Approximation of a Convection-Diffusion Boundary Value Problem on an Adaptive Grid.},
    year = {1996},
    journal={IMA J. Numer. Anal.},
    year={1999},
    volume={19},
    number={2},
    pages={233-249},
    doi={ 10.1093/imanum/19.2.233}
}

@article{malley1977numerical,
            title={ The numerical solution of boundary value problems for stiff differential equations},
            author={ Flaherty, J. E. and  O'Malley, R. E.},
            journal={ Math. Comp. },
            volume={31},
            year={1977},
            pages={ 66-93},
            doi={ 10.1090/S0025-5718-1977-0657396-0 }
            }

@article{miller1976construction,
          title={Construction of a Fem for a Singularly Perturbed Problem in 2 Dimensions},
          author={Miller, J.J.H. },
          publisher={Birkh{\"{a}}user Basel},
          year={1976},
          pages={165-169},
          doi={10.1007/978-3-0348-5328-6_10},
          volume={31},
          journal={Numerische Behandlung von Differentialgleichungen Band 2}        
          }

@article{niijima ,
        author =       { Niijima, K.},
        title =        {A uniformly convergent difference scheme for a
                 semilinear singular perturbation problem},
        journal =      {Numer.Math.},
        volume =       {43},
        number =       {2},
        pages =        {175--198},
        year =         {1984},
        doi={10.1007/BF01390122}
         }

@article{osher1981nonlinear,
        author={ Osher, S.},
        doi={10.1137/0718010},
        journal={ SIAM J. Numer. Anal.},
        volume={ 18},
        number={1},
        pages={ 129--144},
        title={Nonlinear Singular Perturbation Problems and One Sided Difference Schemes},
        year={1981}
         }

@article{prandtl1904uber,
    author = {Prandtl, L.},
    journal = {Verhandl. III. Intern. Math. Kongr. Heidelberg},
    title = {{\"{U}}ber Fl{\"{u}}ssigkeitsbewegung bei sehr kleiner Reibung},
    year = {1904}
}

@article{riordan1986uniform, 
        author = {O'Riordan, E. and  Stynes, M.},
        journal = {Numer. Math.},
        pages = {519-532},
        title = {$L^1$ and $L^{\infty}$ Uniform Convergence of a Difference Scheme for a Semilinear Singular Perturbation Problem.},
        url = {http://eudml.org/doc/133169},
        volume = {50},
        year = {1986/87}
        }

@article{riordan1986uniformly,
          author ={O'Riordan, E. and  Stynes, M.},
          title={A uniformly accurate finite-element method for a singularly perturbed one-dimensional reaction-diffusion problem},
          journal ={Mathematics of Computation},
          volume={47},
          number={176},
          pages={555--570},
          year={1986},
          doi={10.1090/S0025-5718-1986-0856702-7 }
      }

@article{riordan2001parameter,
%         author={Farrell,A., Paul and O'Riordan, Eugene and Miller,J.,H., John and Shishkin, I.,  Grigorii},
%         title={Parameter uniform Fitted Mesh Method for Quasilinear Differential Equations with Boundary Layers},
%         journal={Computational Methods in Applied Mathematics Comput. Methods Appl. Math.},
%         volume={1},
%         number={2},
%         doi={10.2478/cmam-2001-0011},
%         pages={154--172},
%         year={2001}
%         }

@article{riordan2001parameter1,
author={Farrell, P.A. and O'Riordan, E. and Miller,J.J.,H. and Shishkin,G.I.},
title = {PARAMETER-UNIFORM FITTED MESH METHOD FOR QUASILINEAR DIFFERENTIAL EQUATIONS WITH BOUNDARY LAYERS},
 journal={Computational Methods in Applied Mathematics},
 volume={1},
 number={2},
   year={2001},
pages={154--172},
masid = {2320826}
}

@article{roos1990,
   title = {Global uniformly convergent schemes for a singularly perturbed boundary-value problem using patched base spline-functions},
   journal = {Journal of Computational and Applied Mathematics},
   volume = {29},
   number = {1},
   pages = {69--77},
   year = {1990},
   author = {Roos, H.G.},
   doi={10.1016/0377-0427(90)90196-7}
         }

@article{roos1997comparison,
         author = {Roos,H.G. and Skalick{\' y},T.},
         title={ A comparison of the finite element method on Shishkin and Gartland-type meshes for
                convection-diffusion problems},
         journal={QWI Quarterly},
         volume={10},
         number={3 {\&} 4},
         pages={ 277-300},
         year={ 1997},
         url={http://oai.cwi.nl/oai/asset/18320/18320B.pdf}
          }

@article{roos1999sufficient,
 author = {Roos, H.G. and Lin{\ss}, T.},
 title = {Sufficient Conditions for Uniform Convergence on Layer-adapted Grids},
 journal = {Computing},
 volume = {63},
 number = {1},
 pages = {27--45},
 url = {http://dx.doi.org/10.1007/s006070050049},
 doi = {10.1007/s006070050049},
 acmid = {312246},
 publisher = {Springer-Verlag New York, Inc.},
 year={1999}
}

@article{roos2006,
         author = {Roos, H.G.},
        journal = {Applications of Mathematics},
        number = {1},
         pages = {63--72},
       publisher = {Springer},
         title = {Error Estimates for Linear Finite Elements on Bakhvalov-Type Meshes},
        volume = {51},
         year = {2006}
       }

@article{roos2008robust,
        title = {Robust Numerical Methods for Singularly Perturbed Differential Equations: A Survey Covering 2008--2012},         
         author = {Roos, H.G.},
        journal = {ISRN Applied Mathematics},
         publisher = {Hindawi Publishing Corporation},
         volume = {2012},
         doi = {10.5402/2012/379547},
         pages = {1--30},
         year = {2012}
        }

@article{samarskij1980use,
title = "Use of exact difference schemes for estimating the rate of convergence of the method of straight lines ",
journal = "\{USSR\} Computational Mathematics and Mathematical Physics ",
volume = "20",
number = "2",
pages = "102 -- 119",
year = "1980",
issn = "0041-5553",
doi = "10.1016/0041-5553(80)90027-0",
url={http://www.mathnet.ru/php/archive.phtml?wshow=paper&jrnid=zvmmf&paperid=5171&option_lang=eng },
%url = "http://www.sciencedirect.com/science/article/pii/0041555380900270",
author = " Makarov, V.L. and Samarskii, A.A."
}

@article{shishkin1988grid,
  title={Grid approximation of singularly perturbed parabolic equations with internal layers},
  author={Shishkin, G. I.},
  journal={Sov. J. Numer. Anal. M.Russian Journal of Numerical Analysis and Mathematical Modelling},
  volume={3},
  number={5},
  pages={393--408},
  year={1988},
  language={ In Russian},
  doi={ 10.1515/rnam.1988.3.5.393} 
}

@article{stynes1986finite,
         journal={Numerische Mathematik},
         year={1986/87},
         volume ={50},
         number={1},
         pages={1-15},
         title={A finite element method for a singularly perturbed boundary value problem},
         author={Stynes, M. and O'Riordan, E.},
         doi={10.1007/BF01389664}
           
       }

@article{stynes1987adaptive,
   author = {Stynes, M.},
   title = {An Adaptive Uniformly Convergent Numerical Method for a Semilinear Singular Perturbation Problem},
   journal = {SIAM Journal on Numerical Analysis},
   volume = {26},
   number = {2},
   pages = {442-455},
   year = {1989},
   doi = {10.1137/0726025},
   URL = {http://dx.doi.org/10.1137/0726025},
   eprint = {http://dx.doi.org/10.1137/0726025}
       }

@article{stynes1987 ,
        author={Stynes,M. and O'Riordan,E.},
        title={$L^1$ and $L^{\infty}$ uniform convergence of a difference scheme for a semilinear singular perturbation problem},
        journal={Numer. Math.},
        volume={50},
        year={1987},
        pages={519--531},
        doi={10.1007/BF01408573}
        }

@article{stynes1989adaptive,
  title={An adaptive uniformly convergent numerical method for a semilinear singular perturbation problem},
  author={Stynes, M.},
  journal={SIAM journal on numerical analysis},
  volume={26},
  number={2},
  pages={442--455},
  year={1989},
  publisher={SIAM}
}

@article{stynes1991analysis,
        title={ An Analysis of a Singularly Perturbed Two-Point Boundary Value Problem Using Only Finite Element Techniques},
        author={Stynes, M. and O'Riordan, E.},
        journal={Mathematics of Computation},
        volume={ 56},
        number={194},
        year={1991},
        pages={ 663-675},
        publisher={American Mathematical Society},
        url={ http://www.jstor.org/stable/2008400}
         }

@article{stynes1991globally,
         title={A globally uniformly convergent finite element method for a singularly perturbed elliptic problem in two dimensions},
         author={O'Riordan, E. and Stynes, M.},
         journal={ Math. Comp.},
         volume ={57},
         year={1991}, 
         pages={47-62},
         note={MR 1079029},
         doi={10.1090/S0025-5718-1991-1079029-1 }
         }

@article{stynes1991uniformly,
         author={Roos, H.G. and Stynes M.},
         title={A uniformly convergent discretization method for a fourth order singular perturbation problem},
         journal={Bonn. Math. Schrift.},
         volume={228},
         year={1991},
         pages={30-40}   
          }

@book{stynes1993uniformly,
  title={A uniformly convergent method for a singularly perturbed semilinear reaction-diffusion problem with nonunique solutions},
  author={Sun, G. and Stynes, M.},
  year={1993},
  publisher={University College Cork. Department of Mathematics}
}

@article{stynes1996,
         author = {Sun, G. and Stynes, M.},
         journal = {Math. Comput.},
         number ={ 215},
         pages = {1085--1109},
         title = {A uniformly convergent method for a singularly perturbed semilinear reaction-diffusion problem with multiple solutions},
         volume = {65},
         year = {1996},
         url={http://www.jstor.org/stable/2153793}
        }

@article{stynes1995almost,
  title={An almost fourth order uniformly convergent difference scheme for a semilinear singularly perturbed reaction-diffusion problem},
  author={Sun, G. and Stynes, M.},
  journal={Numerische Mathematik},
  volume={70},
  number={4},
  pages={487--500},
  year={1995},
  publisher={Springer-Verlag},
  doi={10.1007/s002110050130}
}

@article{stynes2006,
%         title = {Numerical analysis of singularly perturbed nonlinear reaction-diffusion problems with multiple solutions},
%         journal = {Computers and Mathematics with Applications },
%         volume = {51},
%         number = {5},
%         pages = {857 -- 864},
%         year = {2006},
%        author = {Stynes, Martin and Kopteva, Natalia},
%         doi={10.1016/j.camwa.2006.03.013}
%       }

@article{stynes2003jejune,
         author={Stynes, M.},         
         title={ A Jejune Heuristic Mesh Theorem},
         journal={Computational Methods in Applied Mathematics},
         number={3},
         volume={3},
         pages={488--492},
         doi={10.2478/cmam-2003-0031},
         url={https://www.degruyter.com/downloadpdf/j/cmam.2003.3.issue-3/cmam-2003-0031/cmam-2003-0031.pdf}
         }

@article{stynes2006numerical,
  title={Numerical analysis of singularly perturbed nonlinear reaction-diffusion problems with multiple solutions},
  author={Stynes, M. and Kopteva, N.},
  journal={Computers and Mathematics with Applications},
  volume={51},
  number={5},
  pages={857--864},
  year={2006},
  publisher={Elsevier},
  doi={10.1016/j.camwa.2006.03.013}
}

@article{surla1982,
        author={Surla, K.},
        title = {On convergence of some finite difference schemes for a singular perturbation problem},
        journal ={Univ. u Novom Sadu Zb. Rad, Prirod-Mat. Fak. Ser. Mat},
        volume={12},
        year={1982},
        pages={191--203}
        }

@article{surla1984,
         author={Surla, K.},
         title={Accurate Increase for Some Spline Solutions of Two-Point Boundary Value Problems},
         journal={Novi Sad J. Math.},
         volume={14},
         number={1},
         year={1984},
         pages={51-61}
         }

@article{surla1984on,
         author={ Surla, K.},
         title={On an Error Estimation for Quintic Spline Solutions of Boundary Value Problems},
         journal={Novi Sad J. Math.},
         volume={14},
         number={2},
         year={1984},
         pages={112-115}
         }

@article{surla1986 ,
        author={Surla, K.},
        title={The singularly perturbed spline collocation method for boundary value problems with mixed boundary conditions},
        journal={Univ. u Novom Sadu Zb. Rad, Prirod-Mat. Fak. Ser. Mat},
        volume={16},
        number={2},
        year={1986},
        pages={132--143}
        }

@article{surla1987,
        author={Surla, K.},
        title={A uniformly convergent spline difference scheme for a self-adjoint singular perturbation problem},
        journal={Univ. u Novom Sadu Zb. Rad, Prirod-Mat. Fak. Ser. Mat},
        volume={17},
        number={2},
        year={1987},
        pages={31--38},
        url={http://www.dmi.pmf.uns.ac.rs/nsjom/Papers/17_2/NSJOM_17_2_031_038.pdf}
        }

@article{surla1989spline,
          author={Surla, K. and Jerkovi{\' c}, V.},
          title={Spline Difference Scheme on a Non-Uniform Mesh},
          journal={Novi Sad J. Math.},
          volume={19},
          number={2},
          pages={149-160},
          year={1989}
          }

@article{surla1989exponential, 
        author={Surla, K. and Jerkovi{\' c}, V.}, 
        title={An Exponential Fitted Quadratic Spline Difference Scheme on a Non-Uniform Mesh},
        journal={Novi Sad J. Math.},
        volume={19},
        number={1},
        year={1989},
        pages={1-10}
        }

@article{surla1990,
         author={ Surla, K. and Herceg, D. and Cvetkovi{\' c}, Lj.},
         title={A Family of Exponential Spline Difference Schemes},
         journal={Novi Sad J. Math.},
         volume={20},
         number={1},
         year={1990},
         pages={17-26}
         }

@article{ surla1993,
        author={Surla,K.},
        title={Collocation spline methods in solving boundary value problems},
        journal={Univ. u Novom Sadu Zb. Rad, Prirod-Mat. Fak. Ser. Mat},
        volume={ 23},
        number={2},
        year={1993},
         pages={345--359},
         url={http://www.dmi.uns.ac.rs/nsjom/Papers/23_2/NSJOM_23_2_345_350.pdf}
        }

@article{surla1993analysis,
         author={ Surla, K. and Uzelac, Z.},
         title={An Analysis and Improvement of the El Mistakawy and Werle Scheme},
         journal={Publications de L'Institut Math{\' e}matique},
         volume={54},
         number={68},
         year={1993},
         pages={144-155}
         }

@article{surla1994 ,
        author={Surla,K.},
        title={Some difference schemes derived via finite element method},
        journal={Univ. u Novom Sadu Zb. Rad, Prirod-Mat. Fak. Ser. Mat},
        volume={24},
        number={1},
        year={1994},
        pages={341--346},
        url={http://www.dmi.pmf.uns.ac.rs/nsjom/Papers/24_1/NSJOM_24_1_341_346.pdf}
        }

@article{surla1994a ,
        author={Surla,K.},
        title={On numerical solving singularly perturbed boundary vlaue problems by spline in tension},
        journal={Univ. u Novom Sadu Zb. Rad, Prirod-Mat. Fak. Ser. Mat},
        volume={24},
        number={2},
        year={1994},
        pages={175--186},
        url={http://www.dmi.uns.ac.rs/nsjom/Papers/24_2/NSJOM_24_2_175_186.pdf}
        }

@article{surla1995 ,
        author={Surla,K.},
        title={A spline difference scheme for boundary value problem with small parameter},
        journal={Univ. u Novom Sadu Zb. Rad, Prirod-Mat. Fak. Ser. Mat},
        volume={25},
        number={2},
        year={1995},
         pages={159--168},
         url={http://www.dmi.uns.ac.rs/nsjom/Papers/25_2/NSJOM_25_2_159_168.pdf}
        }

@article{surla1995adifference,
         author={ Surla, K. and Uzelac, Z.},
         title={A Difference Scheme for Boundary Problem with Turning Point},
         journal={Novi Sad J. Math.},
         volume={25},
         number={1},
         year={1995},
         pages={67-73}
         }

@article{surla1995spline,
         author={ Surla, K. and Vukoslav{\v c}evi{\' c}, V.},
         title={A Spline Difference Scheme for Boundary Value Problem with Small Parameter},
         journal={Novi Sad J. Math.},
         volume={25},
         number={2},
         year={1995},
         pages={159-168}
         }

@article{surla1996 ,
        author={Surla,K. and Uzelac,Z.},
        title={A uniformly accurate difference scheme for singular perturbation problem},
        journal={Indian J. pure appl. Math.},
        volume={Vol.27(10)},
        year={1996},
         pages={1005--1016}
        }

@article{surla1996a ,
        author={Vukoslav\v{c}evi\'c,V. and Surla,K.},
        title={Finite element method for solving self adjoint singularly perturbed boundary value problems},
        journal={Mathematica montisnigri},
        volume={Vol. VII},
        year={1996},
         pages={79--86}
       }

@article{surla1996spectral,
         author={ Cvetkovi{\' c}, Lj. and Surla, K.},
         title={Spectral Analysis in Connection with Iterative Solution of Convection-Diffusion Equation},
         journal={Novi Sad J. Math.},
         volume={26},
         number={1},
         year={1996},
         pages={129-134}
         }

@article{surla1998,
        author={Surla,K.},
        title={Exponential functions as boundary layer functions},
        journal={Novi Sad J. Math.},
        volume={28},
        number={3},
        pages={159--176},
        year={1998},
        url={http://www.dmi.uns.ac.rs/nsjom/Papers/28_3/NSJOM_28_3_159_176.pdf}
        }

@article{surla1998some,
         author={Vukoslav{\v c}evi{\' c}, V. and Surla, K. and Rapaji{\' c}, S.},
         title={Some Uniformly Convergent Schemes on Shishkin Mesh},
         journal={Novi Sad J. Math.},
         volume={28},
         number={2},
         year={1998},
         pages={33-42}
          }

@article{surla1998higher,
         author={ Uzelac, Z. and Surla, K.},
         title={A Higher Order Approximation to a Singular Perturbation Problem},
         journal={Novi Sad J. Math.},
         volume={28},
         number={2},
         year={1998},
         pages={19-31}
          }

@article{surla1998thecubic,
         author={ Surla, K.},
         title={The Cubic Spline Difference Scheme on Non-Uniform Mesh},
         journal={Novi Sad J. Math.},
         volume={28},
         number={2},
         year={1998},
         pages={9-17}
         }

@article{surla2000collocation,
         author={Uzelac, Z. and  Surla, K. and Pavlovi{\' c}, Lj.},
         title={On Collocation Methods for Singular Perturbation Problems},
         journal={Novi Sad J. Math.},
         volume={30},
         number={3},
         year={2000},
         pages={173-183} ,
         url={http://www.dmi.uns.ac.rs/nsjom/Papers/30_3/NSJOM_30_3_173_183.pdf}
          }

@article{surla2001collocationmethod,
         author={Surla, K. and Teofanov, Lj. and Uzelac Z.},
         title={On Collocation Methods for Singular Perturbation Problems of Convection-Diffusion Type},
         journal={Novi Sad J. Math.},
         volume={31},
         number={1},
         year={2001},
         pages={125-132} 
          }

@article{surla2001stability,
         author={ Surla, K. and Vukoslav{\v c}evi{\'c}, V.},
         title={On Stability of Some Difference Schemes for Parabolic Diffusion Equations},
         journal={Novi Sad J. Math.},
         volume={31},
         number={2},
         year={2001},
         pages={27-34}
         }

@article{surla2003uniformly,
         author = {Uzelac, Z. and Surla, K.}, 
         title = {A uniformly accurate collocation method for a singularly perturbed problem},
         journal = {Novi Sad J. Math.},
         volume = { 33},
         number={ 1},
         year = {2003},
         pages = {133--143},
         url={http://www.dmi.uns.ac.rs/nsjom/Papers/33_1/nsjom_33_1_133_143.pdf}
         }

@article{surla2003,
         author={ Surla, K. and Uzelac, Z.},
         title={On Stability of Spline Difference Scheme for Reaction-Diffusion Time-Dependent Singularly Perturbed Problem},
         journal={Novi Sad J. Math.},
         volume={33},
         number={2},
         year={2003},
         pages={89-94},
         url={http://www.dmi.uns.ac.rs/nsjom/Papers/33_2/nsjom_33_2_089_094.pdf}
         }

@article{surla2005structure,
         author={ Surla, K. and Teofanov, Lj. and Uzelac, Z.},
         title={The Structure of Spline Collocation Matrix Singularly Pertubation Problems with Two Small Parameters},
         journal={Novi Sad J. Math.},
         volume={35},
         number={1},
         year={2005},
         pages={41-48}
         }

@article{surla2007,
         author={ Surla, K. and Teofanov, Lj. and Uzelac, Z.},
         title={Spline Difference Scheme and Minimum Principle for a Reaction-Diffusion Problem},
         journal={Novi Sad J. Math.},
         volume={37},
         number={2},
         year={2007},
         pages={249-258}
         }

@article{vulanovic1982,
        author={Vulanovi{\'c}, R.},
        title={An exponential fitted scheme on a non-uniform mesh},
        journal={Univ. u Novom Sadu Zb. Rad, Prirod-Mat. Fak. Ser. Mat},
        volume={12},
        year={1982},
        pages={205-215}
       }

@article{vulanovic1983numerical,
         author={Vulanovi{\' c}, R.},
         title={On a Numerical Solution of a Type of Singularly Perturbed Boundary Value Problem by Using a Special Discretization Mesh},
         journal={Novi Sad J. Math.},
         volume={13},
         year={1983},
         pages={187-201},
         url={http://www.dmi.uns.ac.rs/nsjom/Papers/13/NSJOM_13_187_201.pdf}
         }

@article{vulanovic1986extrapolation,
  title={On the extrapolation for a singularly perturbed boundary value problem},
  author={Vulanovi{\'c}, R. and Herceg, D. and Petrovi{\'c}, N.},
  journal={Computing},
  volume={36},
  number={1-2},
  pages={69--79},
  year={1986},
  publisher={Springer-Verlag},
  doi={10.1007/BF02238193}
        }

@article{vulanovic1983,
        author={Vulanovi\'c, R.},
        title={On a numerical solution of a type of singularly perturbed boundary value problem by using a special discretization mesh},
        journal={Univ. u Novom Sadu Zb. Rad, Prirod-Mat. Fak. Ser. Mat},
        volume={13},
        year={1983},
        pages={187--201},
        url={http://www.dmi.pmf.uns.ac.rs/nsjom/Papers/13/NSJOM_13_187_201.pdf}
       }

@article{vulanovic1987numerical,
         author={Vulanovi{\' c}, R.},
         title={On Numerical Solution of a Turning Point Problem},
         publisher={Institute for Numerical and Computational Analysis, Dublin, Ireland},
         volume={8},
         year={1987},
         pages={1-17}
         }

@article{vulanovic1988second,
         author={Vulanovi{\' c}, R.},
         title={A Second Order Uniform Numerical Method for a Turning Point Problem},
         journal={Novi Sad J. Math.},
         volume={18},
         number={1},
         year={1988},
         pages={17-30}
         }

@article{vulanovic1989,
        author={Vulanovi\'c, R.},
        title={Mesh generation methods for numerical solution of quasilinear singular perturbation problems},
        journal={Univ. u Novom Sadu Zb. Rad, Prirod-Mat. Fak. Ser. Mat},
        volume={19},
        number={2},
        year={1989},
         pages={171--193},
         url={http://www.emis.ams.org/journals/NSJOM/Papers/19_2/NSJOM_19_2_171_193.pdf}
       }

@article{vulanovic1989numerical,
         author={Vulanovi{\' c}, R.},
         title={A Numerical Solution of the Singular Perturbation Problem Arising from a Weakly Coupled System},
         journal={Novi Sad J. Math.},
         volume={19},
         number={1},
         year={1989},
         pages={25-32}
         }

@article{vulanovic1989onnumerical,
         author={Vulanovi{\' c}, R.},
         title={On Numerical Solution of a Turning Point Problem},
         journal={Novi Sad J. Math.},
         volume={19},
         number={1},
         year={1989},
         pages={11-24}
         }

@article{vulanovic1991non,
  title={Non-equidistant finite difference methods for elliptic singular perturbation problems},
  author={Vulanovi{\'c}, R.},
  journal={Computational Methods for Boundary and Interior Layers in Several Dimensions, Bool Press, Dublin},
  pages={203--223},
  year={1991}
}

@article{vulanovic1991second,
  title={A second order numerical method for non-linear singular perturbation problems without turning points},
  author={Vulanovi{\'c}, R.},
  journal={USSR Comp. Math. Math+.},
  volume={31},
  number={4},
  pages={522--532},
  year={1991}
}

@article{vulanovic1993,
        author={Vulanovi\'c, R.},
        title={On numerical solution of semilinear singular perturbation problems by using the Hermite scheme},
        journal={Univ. u Novom Sadu Zb. Rad, Prirod-Mat. Fak. Ser. Mat},
        volume={23},
        number={2},
        year={1993},
        pages={363--379},
        url={http://www.dmi.uns.ac.rs/nsjom/Papers/23_2/NSJOM_23_2_363_379.pdf}
       }

@article{vulanovic1993numerical,
         author={Vulanovi{\' c}, R.},
         title={On Numerical Methods for Quasilinear Singular Perturbation Problems without Turning Points},
         journal={Novi Sad J. Math.},
         volume={23},
         number={1},
         year={1993},
         pages={361-370}
         }

@article{vulanovic1994mesh,
         author={Vulanovi{\' c}, R. and Ovcin, Z.},
         title={On Mesh Generation for Singular Perturbation Problems},
         journal={Novi Sad J. Math.},
         volume={24},
         number={1},
         year={1994},
         pages={331-340}
         }

@article{vulanovic1997fourth,
	affiliation = {Kent State University},
	author = {Vulanovi{\' c}, R.},
	copyright = {Kluwer Academic Publishers},
	doi = {10.1023/A:1019187013584},
	journal = {Numerical Algorithms},
	language = {eng},
	number = {2},
	pages = {117-128},
	title = {Fourth order algorithms for a semilinear singular perturbation problem},
	volume = {16},
	year = {1997}
	      }

@article{vulanovic2001special,
         author={Vulanovi{\' c}, R.},
         title={Special Meshes and Higher-Order Schemes for Singularly Perturbed Boundary Value Problems},
         journal={Novi Sad J. Math.},
         volume={31},
         number={1},
         year={2001},
         pages={1-7},
         url={http://www.dmi.uns.ac.rs/nsjom/Papers/31_1/NSJOM_31_1_001_008.pdf}
         }

@article{vulanovic2001uniform,
         author={Lin{\ss}, T. and Vulanovi{\' c}, R.},
         title={Uniform Methods for Semilinear Problems with an Attractive Boundary Turning Point},
         journal={Novi Sad J. Math.},
         volume={31},
         number={2},
         year={2001},
         pages={99-114},
         url={http://www.dmi.uns.ac.rs/nsjom/Papers/31_2/nsjom_31_2_099_114.pdf}
         }

@article{vulanovic2001higher,
 author = {Vulanovi\'{c}, R.},
 title = {A Higher-order Scheme for Quasilinear Boundary Value Problems with Two Small Parameters},
 journal = {Computing},
 volume = {67},
 number = {4},
 pages = {287--303},
 url = {http://dx.doi.org/10.1007/s006070170002},
 doi = {10.1007/s006070170002},
 acmid = {566399},
 publisher = {Springer-Verlag New York, Inc.},
 year={2001}
}

@article{vulanovic2001priori,
  title={A priori meshes for singularly perturbed quasilinear two-point boundary value problems},
  author={Vulanovi{\'c}, R.},
  journal={IMA Journal of Numerical Analysis},
  volume={21},
  number={1},
  pages={349--366},
  year={2001},
  publisher={Oxford University Press},
  doi={10.1093/imanum/21.1.349}
}

@article{vulanovic2004 ,
        author={Vulanovi\'c, R.},
        title={An almost sixth-order finite-difference method for semilinear singular perturbation problems},
        journal={Computational methods in applied  mathematics},
        volume={4},
        number={3},
        year={2004},
        pages={368--383},
        doi={10.2478/cmam-2004-0020}
       }

@article{zarin2012numerical,
         author={Zarin, H. and Gordi{\' c}, S.},
         title={Numerical Solving of Singularly Perturbed Boundary Value Problems with Discontinuities},
         journal={Novi Sad J. Math.},
         volume={42},
         number={1},
         year={2012},
         pages={131-145}
         }

@article{qiu1991analysis,
title = "Analysis of difference approximations to a singularly perturbed two-point boundary value problem on an adaptively generated grid ",
journal = "Journal of Computational and Applied Mathematics ",
volume = "101",
number = "1--2",
pages = "1 - 25",
year = "1999",
doi = {10.1016/S0377-0427(98)00136-8},
url = "http://www.sciencedirect.com/science/article/pii/S0377042798001368",
author = "Qiu, Y. and Sloan, D.M."
}

@article{qiu2000numerical,
title = "Numerical solution of a singularly perturbed two-point boundary value problem using equidistribution: analysis of convergence ",
journal = "Journal of Computational and Applied Mathematics ",
volume = "116",
number = "1",
pages = "121 - 143",
year = "2000",
doi = {10.1016/S0377-0427(99)00315-5}, 
url = "http://www.sciencedirect.com/science/article/pii/S0377042799003155",
author = "Qiu, Y. and Sloan, D.M. and Tang, T."
}

@article{kadalbajoo2003singularly,
title = "Singularly perturbed problems in partial differential equations: a survey ",
journal = "Applied Mathematics and Computation ",
volume = "134",
number = "2--3",
pages = "371 - 429",
year = "2003",
issn = "0096-3003",
doi = "10.1016/S0096-3003(01)00291-0",
url = "http://www.sciencedirect.com/science/article/pii/S0096300301002910",
author = "K. Kadalbajoo, M.K. and Patidar, K.C."
}

@article{kumar2007recent,
 author = {Kumar, M. and Singh, P. and Mishra, H. K.},
 title = {A Recent Survey on Computational Techniques for Solving Singularly Perturbed Boundary Value Problems},
 journal = {Int. J. Comput. Math.},
 issue_date = {October 2007},
 volume = {84},
 number = {10},
 month = oct,
 year = {2007},
 issn = {0020-7160},
 pages = {1439--1463},
 numpages = {25},
 url = {http://dx.doi.org/10.1080/00207160701295712},
 doi = {10.1080/00207160701295712},
 acmid = {1393364},
 publisher = {Taylor \& Francis, Inc.},
 address = {Bristol, PA, USA}
}

@article{kumar2011recent,
        title={A Recent Development of Computer Methods for Solving Singularly Perturbed Boundary Value Problems},
        journal={International Journal of Differential Equations},
        year={2011},
        author={Kumar, M. and Parul}, 
        doi={10.1155/2011/404276},
        url={ttp://dx.doi.org/10.1155/2011/404276},
        volume={2011},
        pages={32}
         }

@article{kumar2014singular,
  title={Singular Perturbation Problems in Nonlinear Elliptic Partial Differential Equations: A Survey},
  author={Kumar, M. and Singh, A. K.},
  journal={International Journal of Nonlinear Science},
  volume={17},
  number={3},
  pages={195--214},
  year={2014},
  %url={http://www.internonlinearscience.org/upload/papers/IJNS%20Vol17%20No%203%20Papers%20%201%20%20Singular%20Perturbation%20Problems%20in%20Nonlinear%20Elliptic%20Partial%20Differential.pdf}
  url={http://goo.gl/fXXiQr}
  }

@article{roos2012robust,
         author={Roos, H.G.},
         title={Robust Numerical Methods for Singularly Perturbed Differential Equations: A Survey Covering 2008--2012},
         journal={ ISRN Applied Mathematics},
         volume={ 2012},
         pages={1- 30},
         year={ 2012},
         doi={10.5402/2012/379547}
         }

@article{shukla2011recent,
  title={A recent development of numerical methods for solving convection-diffusion problems},
  author={Shukla, A. and Singh, A.K. and Singh, P.},
  journal={Applied Mathematics},
  volume={1},
  number={1},
  pages={1--12},
  year={2011},
  publisher={Scientific \& Academic Publishing},
  doi={10.5923/j.am.20110101.01}
  }

@book{ascher1988,
        author =  {Ascher,U.M. and Mattheij,R.M.M. and Russell, R.D.}, 
        title = {Numerical Solution of Boundary Value Problems for Ordinary Differential Equations},
        publisher  = {Prentice Hall,Englewood Cliffs, New Yersey, USA}, 
        year= {1988}
          }

@book{atkinson2014,
       author={Atkinson,K.E},
       title={An Introduction to Numerical Analysis},
       publisher={Wiley \& Sons},
       year={2014},
       isbn={9788126518500}
       }

@book{axelsson2001,
         AUTHOR = {Axelsson,O. and  Barker,V.A.},
         TITLE = { Finite Element Solution of Boundary Value Problems: Theory and Computation},
         PUBLISHER = {SIAM, Philadelphia, USA},
         YEAR= { 2001}
         }

@book{burden2001,
     author={Burden,R.L. and Faires,J.D.},
     title={Numerical Analysis},
     publisher={Brooks/Cole-Thompson Learning, Pacific Grove, USA},
     year={2001}
      }

@book{bahvalov1977,
     author={Bakhvalov,N.S},
     title={Numerical Methods: Analysis, Algebra, Ordinary Differential Equations},
     publisher={MIR},
     year={1977},
     isbn={9780714712079}
     }

@book{bohl1981finite,
         author={Bohl, E.},
         title={Finite Modelle gew\"{o}hnlicher Randwertaufgaben},
         year={1981},
         publisher={Vieweg+Teubner Verlag, Stuttgart},
         doi={10.1007/978-3-322-93092-7},
         isbn={9783322930927}
         }

@book{celic2008,
     author={\v Celi\' c,M.},
     title={Numeri\v cka matematika},
     publisher={Glas srpski, Banjaluka, BiH},
     year={2008},
     isbn={9789993837671}
     }

@book{cheney2004,
      author={Cheney,W. and Kincaid,D.},
      title={Numerical Mathematics and Computing},
      publisher={Brooks/Cole-Thompson Learning, Belmont, USA},
      year={2004}
      }

@BOOK{chag1984,
          author={Chag,K.W. and Howes,F.A.},
          title={Nonlinear Singular Perturbation Phenomena: Theory and Application},
          publisher={Springer-Verlag, Ney York Berlin Heidelberg Tokyo },
          year={1984}
            }

@book{doolan1980uniform,
   title = "Uniform numerical methods for problems with initial and boundary layers",
   author = "Doolan, E.P.  and Miller, J. J. H. and Schilders, W. H. A.",
   publisher = "Boole Press, Dublin",
%   address = "Dublin (P.O. Box 5, 51 Sandycove Rd., Dun Laoire, Co. Dublin)",
%   url = "http://opac.inria.fr/record=b1086535",
%   isbn = "0-906783-02-X",
%   note = "Includes index",
   year = {1980},
%   doi={10.2307/2007353 },
    isbn={9780906783023}
   }

@book{evans2000numerical,
      author={Evans, G. and Blackledge, J. and Yardley, P.},
      title={Numerical Methods for Partial Differential Equations},
      publisher={Springer-Verlag London Limited},
      year={2000},            
      isbn={9788184894592}
      }

@book{farrell2000,
%       author={Farrell,P.A. and Hegarty,A.F. and Miller, J.J.H. and  O'Riordan,E. and Shishkin,G.I.},
%       publisher={Chapman{\&}Hall/CRC, Boca Raton, Florida, USA},
%       title={Robust Computational Techniques for Boundary Layers},
%      year={2000}
%      }

@book{griffel2002,
       author={Griffel,D.H.},
       title={Applied Functional Analysis},
       publisher={Dover Publications, Inc., Mineola, New York, USA},
       year={2002}
         }

@book{herceg1988,
      author={Stojakovi\'{c},Z. and Herceg,D.},
      title={Numeri\v{c}ke metode linearne algebre},
      publisher={IRO \textquotedblleft Gradjevinska knjiga\textquotedblright, Beograd, SFRJ},
      year={1988}
      }

@book{ilin1992matching,
   title = "Matching of asymptotic expansions of solutions of boundary value problems",
   author = "Il'in, A. M. and Minachin, V. and Ivanov, S.",
   series = "Translations of mathematical monographs",
   publisher = "Providence, R.I. American Mathematical Society",
   url = "http://opac.inria.fr/record=b1133280",
   year = 1992
   }

@book{issacson1966,
     title={Analysis of numerical methods,}, 
     author={Issacson, E. and Keller, H.B.},
     publisher={John Wiley \& Sons, New York},
     year={1966}
     }

@book{iserless1996,
       author={Iserles,A.},
       title={A First Course in the Numerical Analysis of Differential Equations},
       publisher={Cambridge University Press, Cambridge, UK},
       year={1996}
       }

@book{johnson2009, 
      author={Johnson,C.},
      title={Numerical Solution of Partial Differential Equations by the Finite Element Methods},
      publisher={Dover Publications, Inc., Mineola, New York, USA},
      year={2009}
      }

@book{keller1992,
      author={Keller,H.B.},
      title={Numerical Methods for Two-Point Boundary-Value Problems},
      publisher={Dover Publications, Inc., New York, USA},
      year={1992},
      isbn={0-486-66925-4}
       }

@book{kovacevic1998uvod,
      author={Kova\v cevi\' c, I.M. and Ralevi\' c, N.M. and Va\v s, L.V.},
      title={Uvod u matemati\v cku analizu},
      publisher={Fakultet tehni\v ckih nauka u Novom Sadu, Srbija},
      year={1998}
      }

@book{kubichek2008,
       author={Kub\'{i}\v{c}ek,M. and Hlav\'{a}\v{c}ek,V.},
       title={Numerical Solution of Nonlinear Boundary Value Problems  with Applications},
       publisher={Dover Publications, Inc., Mineola, New York, USA},
       year={2008}
       }

@book{kress1998numerical,
     author={Kress, R.},
     title={Numerical analysis},
     publisher={Springer},
     year={1998},
     isbn={9788184896022}
     }

@book{linss2010,
  title={Layer-adapted meshes for reaction-convection-diffusion problems},
  author={Lin{\ss}, T.},
  year={2010},
  publisher={Springer-Verlag Berlin Heidelberg},
  doi={10.1007/978-3-642-05134-0}
}

@book{liseikin1989,
    author       = {Lisejkin, V.D. and Petrenko, V.E.},
    title        = {The adaptive-invariant method for the numerical solution of problems with boundary and 
                     interior layers. (Adaptivno-invariantnyj metod chislennogo resheniya zadach 
                     s pogranichnymi i vnutrennimi sloyami.)},  
    year         = {1989},
    pages        = {259 p.},
    publisher    = {Novosibirsk: Vychislitel'nyj\\ Tsentr SO AN SSSR}
  }

@book{liseikin2001,
      author={Liseikin,V. D.},
      title={Layer Resolving Grids and Transformations for Singular Perturbation Problems},
      publisher={VSP BV, AH Zeist, The Netherlands},
      year={2001},
      isbn={9067643467}
      }

@book{miller1991computational,
       author={Miller, J.J.H.},
       title={Computational Methods for Boundary and Interior Layers in Several Dimensions},
       publisher = {Boole Press, Dublin},
       year={1991},
       isbn={97809067-\\83924}
       }

@book{miller1996,
      author={Miller,J.J.H. and O'Riordan,E. and Shishkin,G. I.},
      publisher={World Scientific Publishing\\ Co.Pte.Ltd.Singapore},
      title={Fitted Numerical Methods for Singular Perturbation Problems},
      year={1996},
      isbn={9810224621}
      }

@book{miller2000robust,
     author={Farrell,P.A. and Hegarty, A. F. and Miller, J.J.H. and O'Riordan, E. and Shishkin, G.I.},
     title={Robust Computational Techniques for Boundary Layers},
     publisher={CHAPMAN {\&} HALL/CRC Boca Raton},
     year={2000},
     isbn={97815848-81926}
      }

@book{melenk2002,
      author={Melenk,J.M.},
      title={hp-Finite Element Methods for Singular Perturbations},
      publisher={Springer-Verlag Berlin Heidelberg},
      year={2002},
      isbn={9783540442011}
      }

@book{morton2014numerical,
      author={Morton, K. W. and Mayers D. F.}, 
      title={Numerical Solution of Partial Differential Equations},
      publisher={Cambridge University Press},
      year={2014},
      isbn={9781107447462},
      url={https://drive.google.com/open?id=1YOfvpWDhstiOuviKPfMXMDthzJ-rc33u}
      }

@book{protter1999,
        author={Protter, M.H. and Weinberger,H.F.},
        title={Maximum Principles in Differential Equations},
        publisher={Springer-Verlag New York, Inc},
        year={1999},
        doi={10.1007/978-1-4612-5282-5},
        isbn={978-1-4612-9769-7}
         }

@book{powers1999boundary,
      author={Powers, D. L.},
      title={Boundary value problems},
      publisher={Harcourt/Academic Press},
      year={1999},
      isbn={9780125637343}
         }

@book{radic1978algebra,
      author={Radi\' c, M.},
      title={Algebra},
      publisher={\v Skolska knjiga, Zagreb, SFRJ},
      year={1978}
      }

@book{elementary1993rosen,
     author={Rosen, K.H.},
     title={Elementary number theory and its applications},
     publisher={Addison--Wesley publishing company, Reading, Massachusetts, USA}, 
     year={1993},
     isbn={0201578891}
     }

@book{sewell2005numerical,
      author={Sewell, G.},
      title={The Numerical Solution of Ordinary and Partial Differential Equations},
      publisher={John Wiley \& Sons, Inc., Hoboken, New Jersey, USA},
      year={2005},
      isbn={9780471735809},
      url={https://drive.google.com/open?id=1HJ30z0YxnBqhjIjFKcTu5JWF9xKa7nEb}
}

@book{shihskin2009,
       author={Shishkin,G. I. and Shishkina,L. P.},       
       title={Difference Methods for Singular Perturbation Problems},
       publisher={Chapman {\&} Hall/CRC, Boca Raton, USA},
       year={2009},
       doi={10.1201/9780203492413},
       isbn={9781584884590}
        }

@book{smith2000,
         author={Smith,D.R.}, 
         title={Singular-perturbation Theory-An Introduction with Applications},
         publisher={Cambridge University Press, Cambridge, UK},
         year={1985}
         }

@book{stynes2008robust,
        author={Roos, H.G. and Stynes, M. and Tobiska, L.},
        title={Robust Numerical Methods for Singularly Perturbed Differential Equations},
        publisher={Springer-Verlag Berlin Heidelberg},
        year={2008},
        isbn={9783540344667},
        doi={10.1007/978-3-540-34467-4}
      }

@book{stynes2007numerical,
      author={Grossmann, C. and Roos, H.G. and Stynes, M.},
      title={Numerical Treatment of Partial Differential Equations},
      publisher={Springer-Verlag Berlin Heidelberg},   
      year={2007},
      isbn={9783540715825},
      doi={10.1007/978-3-540-71584-9},
      url={https://drive.google.com/open?id=1yw87XKslNv9Mcfk8C6AJWHoDcofSiWtc}
      }

@book{stoer1993introduction,
      author={Stoer, J. and Bulirsch, R.},
      title={Introduction to Numerical Analysis},
      publisher={Springer-Verlag, Inc.},
      year={1993},
      isbn={038797878X}
      }

@book{tosic2004,
      title={Uvod u numeri\v cku analizu},
      author={To\v si\' c, D.\DJ.},
      publisher={Akademska misao, Beograd},
      year={2004},
      isbn={8674661548}
      }

@book{olga1968,
      author={Ladyzhenskaya,O. A. and Ural'tseva, N. N.},
      title={Linear and Quasi-Linear Elliptic Equations},
      publisher={Academic Press, New York},
      year={1968},
      isbn={978008095-\\5544}
      }

@book{omalley2013singular,
 author = {O'Malley, R. E.},
 title = {Singular Perturbation Methods for Ordinary Differential Equations},
 year = {2013},
 isbn = {9781461269687},
 publisher = {Springer Publishing Company, Incorporated}
}

@book{ortega2000,
        author={Ortega,J. M. and Rheinboldt,W. C.},
        title={Iterative Solution of Nonlinear Equations in Several Variables},
        publisher={SIAM, Philadelphia, USA},
        year={2000},
        isbn={9780898714616}
         }

@book{radunovic2004numericke,
      author={Radunovi\' c, D.},
      title={Numeri\v cke metode},
      publisher={Akademska misao, Beograd},
      year={2004},
      isbn={8674661319}
     }

@book{rogina2003numericka,
      author={Rogina, M., and Singer, S. and Singer, S.},
      title={Numeri\v cka analiza},
      year={2003},
      publisher={Sveu\v cili\v ste u Zagrebu, PMF - Matemati\v cki odjel}
     }

@book{scitovski2015,
      author={Scitovski, R.},
      title={Numeri\v cka matematika},
      publisher={Sveu\v cili\v ste Josipa Jurja Stossmayera u Osijeku},
      year={2015},
      isbn={9789536931798},
      url={http://www.mathos.unios.hr/images/homepages/scitowsk/Num-2015.pdf}     
     }

@book{stakgold1998green,
       author={Stakgold, I.},
       title={Green's functions and boundary value problem},
       publisher={A Wiley-Interscience Publication, JOHN WILEY \& SONS, INC},
       year={1998},
       isbn={9780471610229}
         }

@book{vasileva1995boundary,
   title = "The boundary function method for singular perturbation problems",
   author = "Vasil'eva, A. B. and Butuzov, V. F. and Kalachev, L. V.",
   series = "SIAM studies in applied mathematics",
   volume = "14",
   publisher = "Society for Industrial and Applied Mathematics",
   address = "Philadelphia",
   url = "http://dx.doi.org/10.1137/1.9781611970784",
   isbn = "0-89871-333-1",
   year = {1995}
}

@book{vysnik1957regular,
%        
%        title={Regular degeneration and boundary layer for linear differential equations with small parameter},
%        publisher={Uspehi Mat. Nauk, 12},
%        year={1957},
%        note={(in Russian)}
%         }

@inproceedings{boglaev2005parallel,
author = "Hardy, M. and Boglaev, I.",
title = "Parallel implementation of a monotone domain decomposition algorithm for nonlinear reaction-diffusion problems",
journal = "ANZIAM J.",
volume = "46",
pages = "C290--C303",
year = 2005,
booktitle = "Proc. of 12th Computational Techniques and Applications Conference CTAC-2004",
%doi={10.0000/anziamj.v46i0.960},
url={http://journal.austms.org.au/ojs/index.php/ANZIAMJ/article/view/960}
}

@inproceedings{duvnjakovic1999class,
              author={Duvnjakovi{\'c}, E.},
              title={ A class of difference schemes for singular perturbation  problem},
              booktitle={   Proceedings of the 7.International Conference on Operation  Research, Croatian OR Society, Osijek},
              year={ 1999}, 
              pages={ 197-208}
                 }

@inproceedings{emelianov1981difference, 
               author={Emel'ianov, K. V.},
               title={Difference Schemes for Singularly Perturbed Boundary Value Problems},
               booktitle={BAIL IV; Proceedings of the Fourth International Conference on Boundary and 
                           Interior Layers - Computational and Asymptotic Methods, Novosibirsk, USSR},
               publisher={Boole Press, Dublin},               
               year={1986},           
               pages={51-60}
               }

@inproceedings{kopteva2006pointwise,
  title={Pointwise error estimates for 2nd singularly perturbed semilinear reaction-diffusion problems},
  author={Kopteva, N.},
  booktitle={Finite Difference Methods: Theory and Applications. Proceedings of the 4th International Conference, Lozenetz, Bulgaria},
  pages={105--114},
  year={2006},
  url={ http://www.staff.ul.ie/natalia/pdf/FDM_blg_06.pdf}
       }

@inproceedings{lisejkin1884numerical,
        author={Lisejkin, D, V. and Janenko, N.N.},
        title={On the numerical solution of equations with interior and exterior boundary layers on 
               a nonuniform mesh},
        booktitle={BAIL III; Proceedings of the 3rd International Conference on Boundary and Interior Layers},
        publisher={Dublin, Ireland},
        pages={68-80},
        year={1984}
         }

@inproceedings{riordan1986analysis,
         author = {O'Riordan, E. and  Stynes, M.},
         title = {Analysis of difference schemes for singularly perturbed differential equations 
                  using a discrete Green's function},
         booktitle={BAIL IV; Proceedings of the Fourth International Conference on Boundary and Interior Layers 
                          - Computational and Asymptotic Methods, Novosibirsk, USSR},
         publisher = {Boole Press, Dublin},
         year = {1986},
         pages = {157-168}
         }

@thesis{farrell1983uniformly,
         author={Farrell, P.A.},
         title={Uniformly  convergent difference schemes for singularly perturbed turning and non--turning point problem,},
         institution={Trinity Colege Dublin},
         year={1983},
         note={Ph.D.Dissertation}
         }

@thesis{dannunzio1986,
         author={D'Annunzio,C.M.},
         title={Numerical analysis of a singular perturbation problem with multiple solution},
         institution={University of Maryland at College Park},
         year={1986},         
         note={Ph.D. Dissertation, (unpublished)}
         }

@thesis{vulanovic1986doktorat,
        author={Vulanovi\'{c}, R.},
        title={Mesh construction for discretization of singularly perturbed boundary value problems},
        publisher={Univerzitet u Novom Sadu},
        year={1986},
        note={Doctoral dissertation},
        url={http://elibrary.matf.bg.ac.rs/bitstream/handle/123456789/81/phdReljaVulanovic.pdf}
         }
  \setcounter{section}{6}
   \pagestyle{empty}

   \newpage 
   \addcontentsline{toc}{chapter}{Indeks pojmova}
   \printindex
\end{document}